# PARALLEL OPTIMIZATION

## Theory, Algorithms, and Applications

YAIR CENSOR
Department of Mathematics and Computer Science
University of Haifa

STAVROS A. ZENIOS
Department of Public and Business Administration
University of Cyprus

**This 30 years old book is posted for free download by the generosity of Oxford University Press that granted permission to allow access to future generations of young scientists. A printed copy can be purchased from the publisher. The authors received for this book "The 1999 INFORMS Computing Society Prize for Research Excellence in the Interface Between Operations Research and Computer Science" from the Institute for Operations Research and the Management Sciences (INFORMS). For a book review by Janos D. Pinter go to: Journal of Global Optimization, 16: 107-108, 2000.**

NUMERICAL MATHEMATICS AND SCIENTIFIC COMPUTATION

*Series Editors*

G. H. Golub

K. W. Morton

R. Jeltsch

W. A. Light

E. Süli

NUMERICAL MATHEMATICS AND SCIENTIFIC COMPUTATION

**The Algebraic Eigenvalue Problem*
J. H. Wilkinson
1965

**Direct Methods for Sparse Matrices*
I. S. Duff, A. M. Erisman, and J. K. Reid
1987

**Computer Solution of Linear Problems*
J. L. Nazareth
1988

**Curve and Surface Fittings with Splines*
P. Dierckx
1993

**Numerical Solution of Sturm-Liouville Problems*
J. D. Pryce
1993

*Parallel and Sequential Methods for Ordinary Differential Equations*
K. Burrage
1995

*Parallel Optimization: Theory, Algorithms, and Applications*
Y. Censor and S. A. Zenios
1997

Monographs marked with an asterisk (*) appeared in the series "Monographs on Numerical Analysis," which has been folded into, and is continued by, the current series.

# PARALLEL OPTIMIZATION

## Theory, Algorithms, and Applications

YAIR CENSOR
Department of Mathematics and Computer Science
University of Haifa

STAVROS A. ZENIOS
Department of Public and Business Administration
University of Cyprus

*To Erga, Aviv, Nitzan, and Keren — Y.C.*

*To Christiana, Efy, and Elena — S.A.Z.*

# Foreword

This book is a must for anyone interested in entering the fascinating new world of *parallel optimization* using parallel processors—computers capable of executing a multitude of complex operations in a nanosecond.

The authors are among the pioneers of this fascinating new world and they tell us about new applications they have explored, which algorithms appear to work best, how various parallel processors differ in their design, and the comparative results they obtained using different types of algorithms on different types of parallel processors.

According to an old adage, the whole can sometimes be much more than the sum of its parts. I am thoroughly in agreement with the authors' belief in the added value of bringing together *applications*, *mathematical algorithms* and *parallel computing techniques.* This is exactly what they found true in their own research and report on in the book.

Many years ago, I, too, experienced the thrill of combining three diverse disciplines: the *application* (in my case linear programs), the *solution algorithm* (the simplex method), and the then *new tool* (the serial computer). The union of the three made possible the optimization of many real-world problems. Parallel processors are the new generation of computers and they have the power to tackle a variety of applications, such as those which require solution in real time, or have model parameters that are not known with certainty, or have a vast number of variables and constraints.

Image restoration tomography, radiation therapy, finance, industrial planning, transportation and economics are the sources for many of the interesting practical problems used by the authors to test the methodology.

George B. Dantzig
Stanford University, 1996

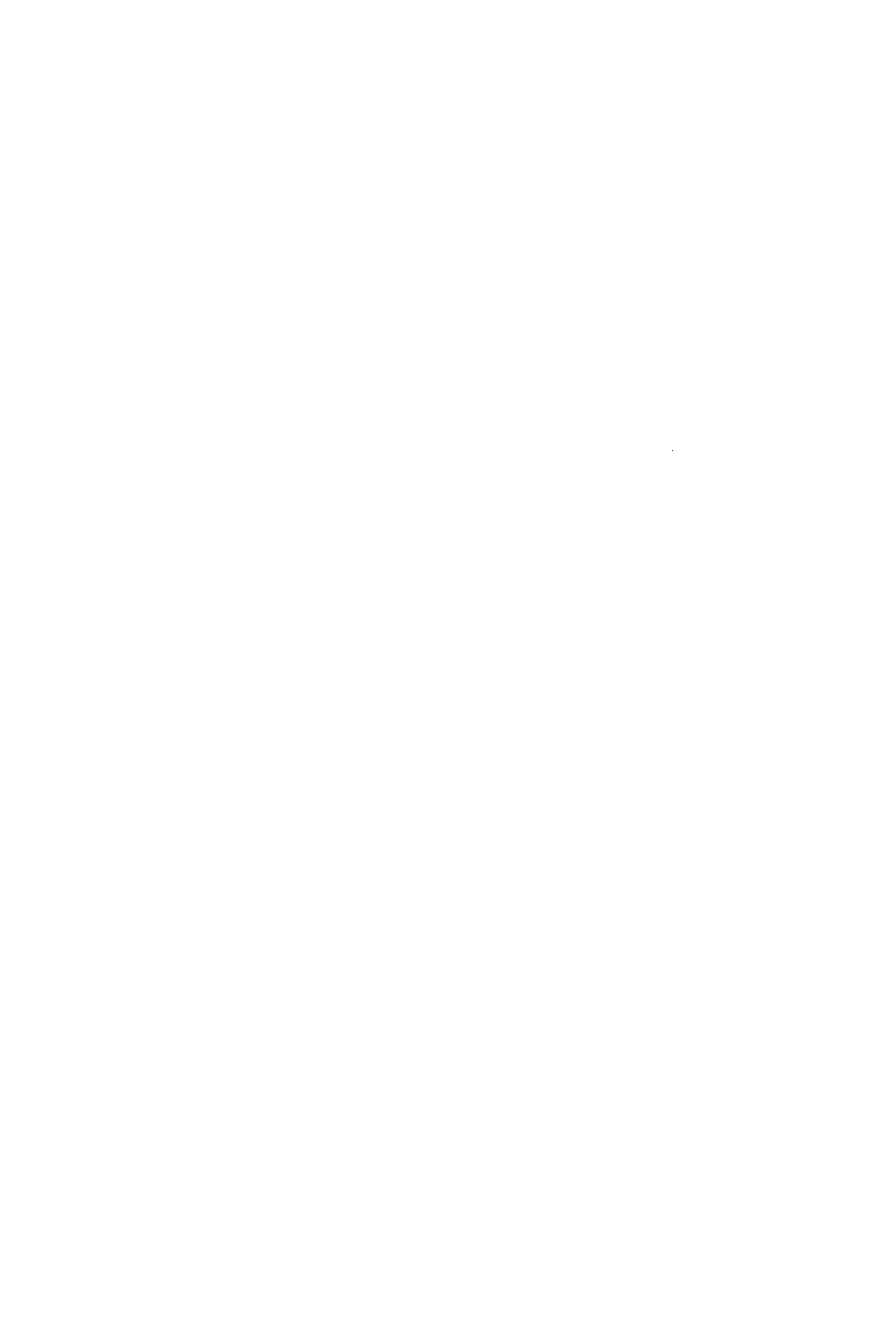

# Preface

> As the sun eclipses the stars by his brilliancy so the one of knowledge will eclipse the fame of the assemblies of the people if he proposes algebraic problems, and still more if he solves them.
>
> Brahmagupta, 650 AD

Problems of mathematical optimization are encountered in diverse areas of the exact sciences, the natural sciences, the social sciences and engineering. Many of them are rooted in real-world applications. Developments in the vast field of optimization are, to a great extent, motivated by these applications and over the years have drawn from both mathematics and computer science. Mathematics creates the foundation for the design and analysis of optimization algorithms. Computer science provides the tools for the design of data structures and for the translation of the mathematical algorithms into numerical procedures that are implementable on a computer. The efficient and robust implementation of an optimization algorithm becomes crucial when one is faced with the solution of large-scale, real-world applications.

The introduction of parallel computer architectures has led to technological innovations which are benefiting every area of scientific computing where large-scale problems are attacked. In this book we give an introduction to methods of *parallel optimization*, by introducing parallel computing ideas and techniques into both optimization theory and into numerical algorithms for large-scale optimization problems. We also examine significant and broad areas of application where the problems are particularly suitable for solution on parallel machines, and where substantial progress has been made in recent years with the application of parallel optimization algorithms.

Some mathematical algorithms that are recognized today as parallel algorithms date back to the 1920s; and in the late 1970s some attempts were made to solve optimization problems with parallel computers using the Illiac IV array processor at the University of Illinois. However, it was not until the early 1980s that concentrated and systematic efforts were instigated by several researchers in the field of parallel optimization. Some of the contributions made over the last two decades have matured to the point where a coherent theoretical framework has been developed, extensive numerical experiments have been carried out, and large-scale problems from diverse areas of application have been solved successfully.

This book gives a comprehensive account of these developments. The coverage is unavoidably not exhaustive, since parallel computing technology has influenced recent developments in all areas of optimization. A series of books could be written on the many applications of parallel computing, including its use in linear programming, large-scale constrained optimization, unconstrained optimization, global optimization and combinatorial optimization (see Section 1.5 for references). This book focuses on parallel optimization methods for large-scale constrained optimization problems and structured linear programs, and as such it provides a comprehensive chart of part of the vast intersection between parallel computing and optimization. We set out to describe a domain where parallel computing is showing a great impact—precisely because of the large-scale nature of the applications—and where many of the recent research developments have occurred. Even within this domain we do not claim that the material presented here about theory, parallel algorithms, and applications is exhaustive. However, related developments that are not treated in the book are discussed in extensive *Notes and References* sections at the end of each chapter.

What then, has determined our selection of theory, algorithms, and applications? Primarily, we have focused on methods where substantial computational experiences have been accumulated over the years, and where, we feel, sufficient integration has been achieved between the theory, the algorithms, and the applications.

Quite often the implementation of an algorithm changes one's perspective on the important features of the algorithm, and such accumulated experience that we have acquired through our own work in the field is reflected in our treatment. In addition, the intricacies of exploiting the problem structure are only fully revealed during an implementation. Finally, it is only with computational experiments that we can have full confidence in the efficiency and robustness of an algorithm. The material presented in this book leads to *implementable* parallel algorithms that have undergone scrutiny on a variety of parallel architectures. And our choice of topics is broad enough so that readers can get a comprehensive view of the landscape of parallel optimization methods.

While not all currently known parallel algorithms are discussed, the book introduces several, from three broad families of algorithms for constrained optimization. Those are defined later in the book as *iterative projection algorithms*, *model decomposition algorithms*, and *interior point algorithms*. When correctly viewed these algorithms satisfy the design characteristics of *good* parallel algorithms.

We begin with the fundamentals of parallel computers: what they are, how to assess their performance, how to design and implement parallel algorithms. This core knowledge on parallel computers is then linked with the theoretical algorithms. The mathematical algorithms and parallel computing techniques are then combined to bear on the solution of several important applications: image reconstruction from projections, matrix balancing, network optimization, nonlinear programming for planning under uncertainty, and financial planning. We also address implementation issues and study results of recent experiments that highlight the efficiency of the developed algorithms when implemented on suitable parallel computer platforms.

We believe that the value of bringing together applications, mathematical algorithms, and parallel computing techniques extends beyond the successful solution of the specific problems at hand. Namely, it familiarizes the reader with the complete process from the modeling of a problem through the design of solution algorithms while it introduces the art and science of parallel computations. It is not possible to study these three disciplinary efforts—modeling, mathematics of algorithms, and parallel computing—in isolation. The successful solution of real-world problems in scientific computing is the result of coordinated efforts across all three fronts. We hope that this book will help the reader develop such a broad perspective and, thus, follow Brahmagupta's admonishment.

To keep the book to a reasonable length we had to make some decisions on priority of topics and on prerequisite knowledge. Many important and relevant topics have been left out or are mentioned only casually. These include questions on rate of convergence, computational complexity, stopping criteria, behavior of the algorithms in inconsistent cases, and so on.

We assume that the reader has been systematically exposed to differential and integral calculus, linear algebra, convex analysis and optimization theory. Sections 10.2, 10.3 and Chapter 13 assume familiarity with notions from probability theory.

Finally, although we have compiled an extensive bibliography, we may have missed relevant references or erred in crediting work. We are grateful to readers who bring to our attention such omissions so that we can correct them in the future.

## Organization of the Book

The material in this book is organized in three parts: Theory, Algorithms, and Applications. There is an introductory chapter (Chapter 1) which outlines the fundamental topics on parallel computing.

**Part I: Theory** develops the theory of generalized distances and generalized projections (Chapter 2) and the theory for their use in solving linear programming problems via proximal minimization (Chapter 3). The theory of penalty and barrier methods and augmented Lagrangians is discussed and explained in Chapter 4. This material provides the theoretical foundation upon which the algorithms in Chapters 5 to 8 are developed.

**Part II: Algorithms** develops iterative projection algorithms, model decomposition algorithms, and interior point algorithms. Chapter 5 discusses iterative algorithms for the solution of convex feasibility problems, using the theory of generalized projections. Similarly, Chapter 6 uses the theory of generalized distances and generalized projections to develop algorithms for linearly constrained optimization problems. Chapter 7 develops model decomposition algorithms, based on the theory of penalty methods and augmented Lagrangians. Chapter 8 introduces interior point algorithms for linear and quadratic programming, and explains ways in which the structure of some large-scale optimization problems can be exploited for parallel computations by these algorithms.

**Part III: Applications** discusses applications from several diverse real-world domains where the parallel algorithms are suitable. Each chapter contains a description of the real-world application; it develops one or more mathematical models for each problem; and discusses solution algorithms from one or more of the algorithm classes of Part II. Chapter 9 discusses problems of matrix estimation. Chapter 10 examines problems of image reconstruction from projections. Chapter 11 focuses on the problem of radiation therapy treatment planning. Chapter 12 features problems in transportation and the multicommodity flow problem. Chapter 13 looks at problems of planning under uncertainty using stochastic programming and robust optimization models.

Finally, two chapters are devoted to the parallel implementation and testing of the algorithms. Implementations are discussed in Chapter 14, and while this issue is an important one, it is implicitly linked closely to the computer architecture. To the extent possible our discussion applies to a whole class of machines rather than to a specific hardware model. Chapter 15 summarizes numerical experiments with several of the algorithms that demonstrate their effectiveness for solving large-scale problems when implemented on parallel machines.

Figure 0.1 illustrates the interdependencies among the chapters. The sequence of chapters indicated in this diagram must be followed to appreciate fully a line of development from the theoretical foundations to the algorithms and their applications. While we emphasize the importance of studying the continuum of theory-algorithms-applications, the chapters are written in such a way that they can be used as a reference, and need not be studied sequentially. Readers interested only in applications and mathematical models may read the relevant chapters from Part III, although to fully appreciate the solution algorithms they must study the earlier chapters as well. But even then, a solution algorithm can be understood without referring to the relevant chapters on theory from Part I, unless the reader wishes to understand the proof of convergence as well. The book can, therefore, be used either as a textbook or as a reference book.

## Suggested course outlines

The book is organized in a way that allows it to be used as a text for graduate courses in large-scale optimization, parallel computing, or large-scale mathematical modeling. There are three different avenues that an instructor may follow in teaching this material, especially considering that the whole book cannot be covered in the usual timeframe of a one-semester course. Regardless of the avenue an instructor decides to follow, Chapter 1 gives a general introduction to the material and should be covered first.

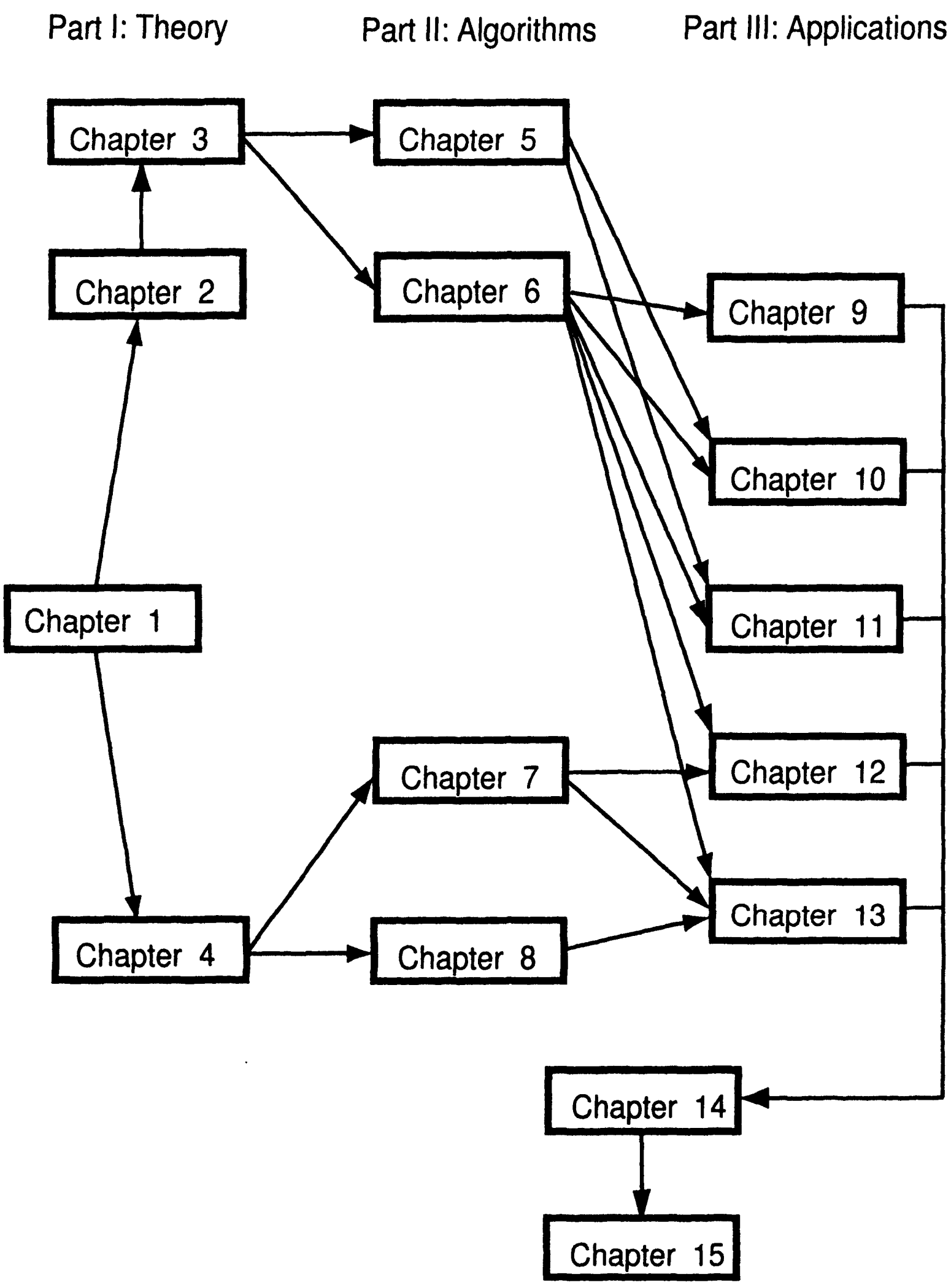

**Figure 0.1** Organization of chapters in the book.

One approach is to teach a course on theory and algorithms for constrained optimization and feasibility problems. Such a course will cover the theory part of the book, Chapters 2, 3 and 4, followed by the algorithms in Chapters 5, 6, 7 and 8. Any one of the Chapters 3, 7 or 8 could be omitted without loss of continuity, but a balanced treatment of different families of optimization algorithms should include Chapters 6, 7 and 8. This course focuses on the theoretical aspects of parallel optimization, and references to parallel computations can be cursory.

A second approach is to teach a course on numerical methods for large-scale structured optimization problems. Such a course will focus on the algorithm section, without prior introduction to the theory that is essential for establishing convergence. Here the emphasis is on developing the students' understanding of the structure of an algorithm, and assuming the underlying theory. Such a course will cover the material from Chapters 6, 7, and 8, which present general algorithms, and also selected sections from Chapters 9, 10, 12, and 13, which illustrate the development of algorithms for specific problems. In this course the exploitation of special structure is a key issue, and references to parallel computing become crucial. This course could also discuss the implementation of algorithms on parallel machines, with coverage of the material in Chapter 14.

Yet a third approach is to focus on applications of optimization, and present—in cookbook fashion—implementable algorithms for the solution of real-world instances of large-scale problems. In this course material from the chapters on applications, with references to the corresponding chapters in Part II where specific algorithms have been developed for the applications at hand will be taught. This course will start with the iterative optimization algorithms of Chapter 6 and move on to applications in Chapters 9, 10, 12, and 13. As in the previous course, the exploitation of special structure is a key issue and references to parallel computing become crucial. This course could also discuss the implementation of algorithms on parallel machines in Chapter 14.

The book can also be used for a course on parallel computing. Several of the algorithms are simple enough so that students with little background in optimization can readily understand them, and such algorithms can be the focus for implementation exercises. Furthermore, the structures of the underlying models vary from the very simple dense matrix (e.g., the dense transportation problems of Sections 12.2.1 and 12.4.1), to sparse and structured matrices and graph problems (e.g., the multicommodity transportation problems and the stochastic networks of Sections 12.2.2, 12.4.3, and 13.8, respectively). Hence, the material can be used to introduce students to the art of implementing algorithms on parallel machines. Such a course will provide some motivation by introducing applications from Chapters 9, 10, 12, and 13. For each application an algorithm can be introduced (specific implementable algorithms are found in the applications chapters), and references made to the implementation techniques of Chapter 14. The course should follow the sequence of topics as discussed in Chapter 14, but before teaching any material from this chapter the corresponding application chapter should first be introduced. Finally, Chapter 15 can be used as a reference for students who wish to test the efficiency of their implementation, or the performance of their parallel machine.

## Acknowledgments

Some of the material presented in this book is based on our own published work, and we wish to express our appreciation to past and present collaborators from whom we have learned a great deal. Particular thanks for many useful discussions and support, and for reading various parts of the book, are extended to Marty Altschuler, Dimitri Bertsekas, Dan Butnariu, Charlie Byrne, Alvaro De Pierro, Jitka Dupačova, Jonathan Eckstein, Tommy Elfving, Gabor Herman, Dan Gordon, Alfredo Iusem, Elizabeth Jessup, Arnold Lent, Robert Lewitt, Olvi Mangasarian, Jill Mesirov, Bob Meyer, John Mulvey, Soren Nielsen, Mustafa Pinar, Simeon Reich, Uri Rothblum, Michael Schneider, Jay Udupa, Paul Tseng, and Dafeng Yang.

A draft version was read by Dimitri Bertsekas, Jonathan Eckstein, Tommy Elfving, Michael Ferris, and Gabor Herman. We thank them for their constructive comments. Any remaining errors or imperfections are, of course, our sole responsibility.

Part of the work of Yair Censor in this area was undertaken in collaboration with the Medical Image Processing Group (MIPG) at the Department of Radiology, University of Pennsylvania. The support and encouragement of Gabor Herman for this continued collaboration is gratefully acknowledged. The work of Stavros A. Zenios was undertaken while he was on the faculty at the Department of Operations and Information Management, The Wharton School, University of Pennsylvania, and while on leave at the Operations Research Center of the Sloan School, Massachussetts Institute of Technology, and with Thinking Machines Corporation, Cambridge, Mass. Substantial portions of Stavros Zenios' work on this book were completed during his visits to the University of Bergamo, Italy, and the support of Marida Bertocchi in making these visits possible is gratefully acknowledged. We express our appreciation to these organizations as well as to our current institutions, the University of Haifa, Haifa, Israel and the University of Cyprus, Nicosia, Cyprus, for creating an environment wherein international collaborations are fostered and where long-term undertakings, such as the writing of this book, are encouraged and supported.

This text grew out of lecture notes that we prepared and delivered at the *19th Brazilian Mathematical Colloquium*, which took place in Rio de Janeiro in the summer of 1993. We thank the organizers of that Colloquium and Professor Jacob Palis, the Director of the Instituto de Mathemática Pura e Aplicada (IMPA) where the Colloquium was held, for giving us the opportunity to present our work at this forum, and thereby giving us the impetus to complete the project of writing this book.

This work was supported by grants from the National Institutes of Health (HL–28438), the Air Force Office of Scientific Research (AFOSR-91–0168) and the National Science Foundation (CCR–9104042 and SES–9100216) in the USA, and the National Research Council (CNR), Italy.

Computing resources were made available through the North-East Parallel Architectures Center (NPAC) at Syracuse University, the Army High-Performance Computing Research Center (AHPCRC) at the University of Minnesota, Thinking Machines Corporation, Cray Research Inc., and Argonne National Laboratory. We also thank Ms. Danielle Friedlander and Ms. Giuseppina La Mantia for their work in processing parts of the manuscript in LaTeX. Finally, we express our appreciation to Senior Editor Bill Zobrist and Editorial Assistant Krysia Bebick from the Oxford University Press office in New York for their help.

Haifa, Nicosia, and Philadelphia<br>December 1996

Yair Censor<br>Stavros A. Zenios

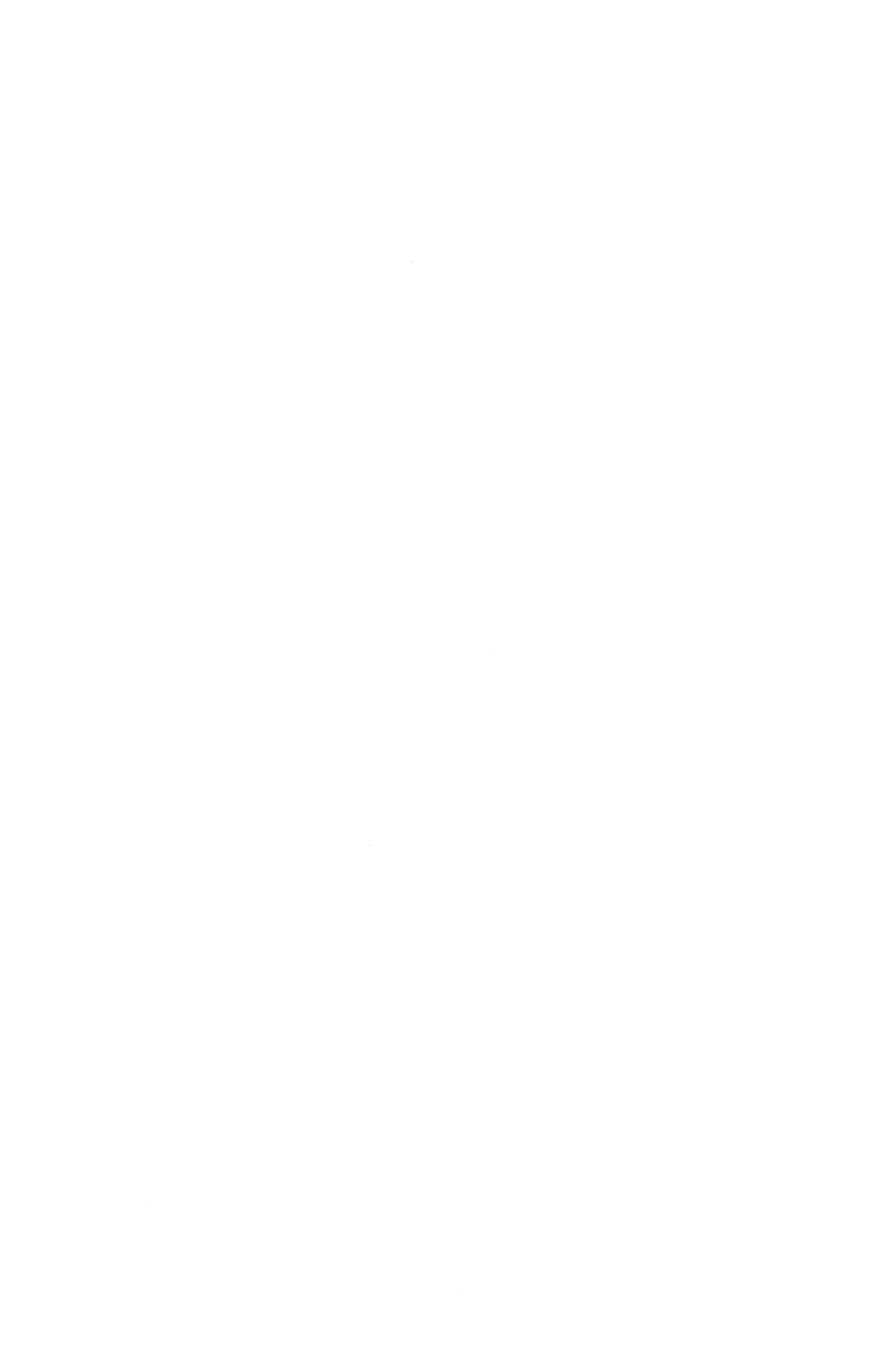

# Contents

# Glossary of Symbols

|  |  |
|---|---|
| $\doteq$ | the equality defines. |
| $\approx$ | approximate equality. |
| $x \leftarrow b$ | assign value $b$ to the variable $x$. |
| $\exists$ | there exists. |
| $\blacksquare$ | end of proof. |
| $\emptyset$ | empty set. |
| $\overline{C}$ | closure of a set $C$. |
| $\text{bd}\, C$ | boundary of a set $C$. |
| $\text{int}\, C$ | interior of a set $C$. |
| $\prod_{i=1}^{m} Q_i$ | product set $Q_1 \times Q_2 \times \ldots \times Q_m$ of the sets $Q_i$. |
| $K \setminus L$ | set difference containing all elements in the set $K$ and not in $L$. |
| $K \oplus L$ | the direct sum of sets $K$ and $L$, i.e., $K \oplus L = \{a + b \mid a \in K,\ b \in$ |
| $\text{conv}(K)$ | convex hull of set $K$. |
| $\text{card}(K)$ | cardinality of set $K$. |
| $\subseteq$ | set inclusion. |
| $\subset$ | strict set inclusion. |
| $\not\subseteq$ | set non-inclusion. |
| $\cap$ | set intersection. |
| $\cup$ | set union. |
| $\mathbb{R}^n$ | $n$-dimensional Euclidean space. |
| $\mathbb{R}^n_+$ | nonnegative orthant of $\mathbb{R}^n$. |
| $\mathbb{R}^n_{++}$ | positive orthant of $\mathbb{R}^n$. |
| $\mathbb{R}^{m \times n}$ | space of all real $m \times n$ matrices. |
| $N_0$ | set of all nonnegative integers. |
| $\nabla f(x)$ | gradient of $f$ at $x$. |
| $df/dx,\ f'$ | derivative of a function $f$. |
| $f'(u; v)$ | directional derivative of $f$ at $u$ in the direction $v$. |
| $\nabla_x h(x, y)$ | gradient of the function $h$ with respect to the vector variable $x$. |
| $\nabla^2 f(x)$ | Hessian matrix of $f$ at $x$. |
| $\partial f(x)$ | subdifferential set, i.e., set of all subgradients of the function $f$ at |
| $\text{epi}(f)$ | epigraph of a function $f$. |
| $g^+(x)$ | function obtained from $g(x)$ by $g^+(x) \doteq \max(0, g(x))$. |
| $\mathcal{R}f$ | the Radon transform of $f$. |
| ent | Shannon's entropy function. |
| exp | exponential function. |
| log | natural logarithm function. |
| $e$ | basis of the natural logarithm function. |
| Min (min) | minimum, the minimal value. |
| Max (max) | maximum, the maximal value. |
| inf | infimum. |
| argmin | point where a minimal value of a function is attained. |
| argmax | point where a maximal value of a function is attained. |
| liminf | the smallest accumulation point of a sequence. |
| $\mathrm{mid}(a, b, c)$ | the median of the three real numbers $a, b$, and $c$. |
| $\langle \cdot, \cdot \rangle$ | inner product in $\mathbb{R}^n$. |
| $\|\cdot\|$, $\|\cdot\|_2$ | Euclidean norm in $\mathbb{R}^n$. |
| $\|\cdot\|_1$ | $\ell_1$-norm in $\mathbb{R}^n$. |
| $\mid t \mid$ | the absolute value of the real number $t$. |
| $d(x, Q)$ | Euclidean distance between the point $x$ and the set $Q$ in $\mathbb{R}^n$. |
| $B(x, \rho)$ | ball with radius $\rho$ and center at $x$, i.e., the set defined by $B(x, \rho) \doteq \{y \in \mathbb{R}^n \mid \|x - y\| \leq \rho\}$. |
| $H(a, b)$ | hyperplane $H = \{x \in \mathbb{R}^n \mid \langle a, x \rangle = b\}$ for $a \in \mathbb{R}^n$ and $b \in \mathbb{R}$. |
| $\mathcal{B}(S)$ | family of Bregman functions with zone $S$. |
| $D_f(x, y)$ | generalized distance between $x$ and $y$. |
| $P_\Omega(y)$, $P_\Omega^f(y)$ | generalized projection of $y$ onto a set $\Omega$ associated with $f$. |
| $P_{i,\lambda}$ | relaxed projection onto a set $Q_i$. |
| $P_{w,\lambda}$ | relaxed convex combination of projections. |
| $(r, t)$ | open interval, on the real line, between $r$ and $t$. |
| $[r, t]$ | closed interval, on the real line, between $r$ and $t$. |
| $A^{\mathrm{T}}$, $x^{\mathrm{T}}$ | transposed matrix or vector. |
| $A^{-1}$ | inverse matrix. |
| $A^\dagger$ | generalized inverse of the matrix $A$. |
| $\det A$ | determinant of the matrix $A$. |
| $\mathcal{R}(A)$ | range of the matrix $A$. |
| $\mathcal{N}(S)$ | null space of an operator $S$. |
| $I_n$ | $n \times n$ identity matrix. |
| $I$ | (used to denote a matrix) is a conformable identity matrix. |
| $I$ | (used to denote a set) is the set $I \doteq \{1, 2, \ldots, m\}$. |
| $\mathbf{1} \in \mathbb{R}^n$ | vector with all components equal to 1. |
| $e^l \in \mathbb{R}^n$ | $l$th basis vector in $\mathbb{R}^n$, the $l$th column of identity matrix $I_n$. |
| $\mathrm{diag}\,(x)$, $\mathrm{diag}\,(x_1, x_2, \ldots, x_n)$ | $n \times n$ diagonal matrix with the $j$th diagonal entry being the $j$th component of the vector $x$. |
| $\mathrm{diag}\,(S_1, S_2, \ldots, S_n)$ | block-diagonal matrix with submatrices $S_l$, $l = 1, 2, \ldots, n$. |
| $A_{i\cdot}$ | row $i$ of the matrix $A$. |
| $A_{\cdot j}$ | column $j$ of the matrix $A$. |
| $a^i$ | column $i$ of the matrix $A^{\mathrm{T}}$. |

# PARALLEL OPTIMIZATION

## CHAPTER 1

# Introduction

> Capable of solving scientific problems so complex that all previously known methods of solution were considered impractical, an electronic robot, known as Eniac—Electronic Numerical Integrator and Computer—has been announced by the War Department. It is able to compute 1,000 times faster than the most advanced general-purpose calculating machine, and solves in hours problems which would take years on a mechanical machine. Containing nearly 18,000 vacuum tubes, the 30-ton Eniac occupies a room 30 by 50 feet.
>
> *Scientific American*, June 1946

Computing power improved by a factor of one million in the period 1955–1990, and it is expected to improve by that factor again just within the next decade. How can this accelerated improvement be sustained? The answer is found in the parallel computer architectures: multiple, semi-autonomous processors coordinating for the solution of a single problem. The late 1980s have witnessed the advent of *massive parallelism* whereby the number of processors that are brought to bear on a single problem could be several thousands. This rapid technological development is transforming significantly the exact sciences, the social sciences and all branches of engineering. The American Academy of Arts and Sciences devoted the winter 1992 issue of its journal, *Daedalus*, to these developments, referring to what was termed a "New Era of Computation."

As parallel machines became available, scientists from several disciplines realized that parallelism was a natural way to view their applications. For example, weather forecasters realized that they could split their models naturally by geographical regions and run forecasting models separately for each region. Only boundary information must be communicated across adjacent regions. Air-traffic controllers can track the flow of traffic on high-altitude jet routes using multiple processors, allocating one processor to each control sector or to each route. Information among processors must be exchanged only as aircraft change routes or move from the jurisdiction of one control sector to the next.

What factors motivate these developments? The first designs of parallel computers were motivated by the limitations of the serial, von Neumann architecture. The performance of computers that process scalars serially is limited by the speed with which signals can be communicated across different components of the machine. The speed of light provides an upper bound on the speed at which serial computations can be executed, and parallel designs appear as the only way to overcome this barrier.

Following the early developments of the parallel computing industry in the middle 1970s, it was realized that parallel machines could also be cost-effective. There are economies of scale in utilizing multiple processors around common memory, and the continued progress in very large-scale integrated circuits provides additional justification for packing multiple processors in the same system. Parallelism is an inevitable step in the development process of computing technology, which major breakthroughs in serial processing technology do not change. Instead, they accelerate it as more efficient processors are coupled together to create parallel machines. For instance, the earlier hypercube systems by Intel Corporation were based on hardware found in personal computers. Their latest models coordinate a thousand processors that are typically found in workstations.

Parallelism is improving by a quantum leap the size and complexity of models that can be represented on and solved by a computer. But these improvements are not sustained by the architecture of the machines alone. Equally important is the development of appropriate mathematical algorithms and the proper decomposition of the problem at hand, in order to utilize effectively the parallel architectures. In this book we study methods of parallel computing for optimization problems. We take foremost an algorithmic approach, developing the theory underlying a broad class of iterative algorithms. Specific algorithms are then presented, their mathematical convergence is established, and their suitability for parallel computations is discussed.

The design of an algorithm is, in turn, strongly influenced by the structure of the application. Therefore, we discuss several applications drawn from such diverse areas as operations research, financial modeling, transportation, and image reconstruction from projections in diagnostic medicine. Other areas where similar models and solution techniques are used are also mentioned. Equipped with the knowledge of algorithms and applications, we also embark on implementations of the algorithms on parallel architectures. Alternative approaches to the parallel implementation of an algorithm are discussed, and benchmark results on various computers are summarized. Those results highlight the effectiveness of the algorithms and their parallel implementations.

Several of the algorithms and their underlying theories have a long history that dates back to the 1920s. At that time the issue of computations—that is, implementation of the algorithm as a computer program—was irrelevant. (The calculations involved in an algorithm were in those days carried out by humans, or by humans operating mechanical calculators. The issues underlying the efficient execution of these calculations are dif-ferent from the issues of computation, as we understand them today.) Here we present the algorithms not just with a view toward computations, but more importantly, toward parallel computations. It is when viewed within the context of parallel computing that these algorithms are understood to possess some important features. They are well-suited for implementation on a wide range of parallel architectures, with as few as two processors or as many as millions. It is precisely our aim to show that the theory presented in this book produces implementable algorithms whose efficiency is independent of the specific parallel machine. The implementations we discuss and the associated computational results, are, to some extent, machine dependent. But the underlying mathematical algorithms can be implemented on a wide variety of parallel machines. The theory underlying parallel optimization algorithms and the algorithms themselves (and, of course, the applications) are independent of the rapidly changing landscape of parallel computer architectures.

The remainder of this chapter introduces the basics of parallel computing. Our treatment is not rigorous, in that we do not fully specify a theoretical model of an abstract parallel machine, nor do we completely classify parallel architectures. Instead, we provide a general classification of computer architectures that is suitable for illustrating, later in the book, the suitability of algorithms for parallel computations. We first introduce the basics of parallel computing in Section 1.1. Then, in Section 1.2, we discuss the major issues that need to be addressed in designing algorithms for parallel machines. Section 1.3 provides a general classification of the parallel optimization algorithms. In Section 1.4 we discuss ways for measuring the performance of parallel algorithms implemented on parallel machines. These measures are used in subsequent chapters to establish the efficiency of the algorithms. Notes and references to the literature are given in Section 1.5.

## 1.1 Parallel Computers

Parallelism in computer systems is not a recent concept. ENIAC—the first electronic digital computer built at the University of Pennsylvania between 1943 and 1946—was designed with multiple functional units for adding, multiplying, and so forth. The primary motivation behind this design was to deliver the computing power that was required by applications but was not possible with the electronic technology of that time. The shift from vacuum tubes to transistors, integrated circuits, and very large-scale integrated (VLSI) circuits rendered parallel designs obsolete and uniprocessor systems were predominant throughout the late 1960s.

A milestone in the evolution of parallel computers was the Illiac IV project at the University of Illinois in the 1970s. The array architecture of the Illiac—i.e., its ability to operate efficiently on arrays of data—prompted studies on the design of suitable algorithms for scientific computing on this machine. A study on the use of array processing for linear programming was carried out by Pfefferkorn and Tomlin (1976). The Illiac never went past the stage of the research project, however, and only one machine was ever built.

Another milestone in the evolution of parallel computers was the introduction of the CRAY 1 in 1976, by Cray Research, Inc. The term *supercomputer* was coined to indicate the fastest computer available at any given point in time. The vector architecture of the CRAY—i.e., the presence of vector registers and vector functional units—introduced the notion of vectorization in scientific computing. Designing or restructuring of numerical algorithms to exploit the vector architecture became a critical issue. Vectorization of an application can range from simple modifications of the implementation, with the use of computational kernels that are streamlined for the machine architecture to fundamental changes in data structures and the design of algorithms that are rich in vector operations.

Since the mid-seventies parallel computers have been evolving rapidly in the level of performance they can deliver, the size of available memory, and the number of parallel processors that can be applied to a single task. The Connection Machine CM–2 by Thinking Machines Corporation, for example, can be configured with up to 65,536 very simple processing elements. The CRAY Y-MP/4 has four powerful vector processors. To understand parallel architectures we start with a general taxonomy.

### 1.1.1 Taxonomy of parallel architectures

Several alternative parallel architectures were developed throughout the 1980s and 1990s, and today there is no single widely accepted model for parallel computation. A classification of computer architectures was proposed by Flynn (1972) and is used to distinguish between alternative parallel architectures. Flynn proposed the following four classes, based on the interaction between instruction and data streams of the processor(s):

*SISD – Single Instruction stream Single Data stream.* Computers in this class execute a single instruction on a single piece of data before moving on to the execution of the next instruction on a different piece of data. The traditional von Neumann computer is a scalar uniprocessor that falls in this category.

*SIMD – Single Instruction stream Multiple Data stream.* Computers with this architecture can execute a single instruction simultaneously on multiple data. This type of computation is possible only if the operations of an algorithm are identical over a set of data and if the data can be mapped on multiple processors for concurrent execution. Examples of SIMD systems are the Connection Machine CM–2 from Thinking Machines Corporation, and the DAP (distributed array processor) from Active Memory Technologies.

*MISD – Multiple Instruction stream Single Data stream.* Computers with this architecture execute multiple instructions concurrently on a single piece of data. This form of parallelism has not received much attention in practice. It appears in the taxonomy for completeness.

*MIMD – Multiple Instruction stream Multiple Data stream.* Computers in this class can execute multiple instructions concurrently on multiple pieces of data. Multiple instructions are generated by code modules that may execute independently from each other. Each module may operate either on a subset of the data of the problem or on a copy of all the problem data, or it may access all the data of the problem, together with the other modules, in a way that avoids read or write conflicts.

Whenever multiple data streams are used, i.e., in the MIMD and SIMD systems, another type of classification is needed to describe the organization of memory. Memory organization is classified either as shared, or distributed. In a shared memory system the multiple data streams are stored in memory that is common to all processors. In a distributed memory system each processor has its own local memory, and data are transferred among processors, as needed by the code, using a communication network. The characteristics of the communication network are important for the design of parallel machines, but they are not essential for the implementation of a parallel algorithm. Communication networks are treated only in passing in this book, but references on the topic are given in Section 1.5.

Multiprocessor systems are also characterized by the number of available processors. Small-scale parallel systems have up to 16 processors, medium-scale systems up to 128, and large-scale systems up to 1024. Systems with more than 1024 processors are considered massively parallel. Finally, multiprocessor systems are also characterized as *coarse-grain* or *fine-grain*. In coarse-grain systems each processor is very powerful, typically of the kind found in contemporary workstations, with several megabytes of memory. Fine-grain systems typically use simple processing elements with a few kilobytes of local memory each. For example, an Intel iPSC/860 system with 1024 processors is a coarse-grain parallel machine. The Connection Machine CM–2 with up to 64K (here 1K = 1024) simple processing elements is a fine-grain, massively parallel machine.

A type of computer that deserves special classification is the *vector computer*, which has special vector registers for storing arrays and special vector instructions that operate on vectors. While vector computers are a simple, special case of SIMD machines, they are usually considered as a different class because vector registers and vector instructions appear in the processors of many parallel MIMD systems as well. Also, the development of algorithms or software for a vector computer (as, for example, the CRAY Y-MP/4 or the IBM 3090–600/VF) poses different problems than the design of algorithms for a system with multiple processors that operate synchronously on multiple data (as, for example, the Connection Machine CM-2 or the DAP).

### 1.1.2 Unifying concepts of parallel computing

The landscape of parallel computer architectures is quite complex. From the perspective of the applications programmer, however, some unifying concepts have emerged that make it possible to program parallel machines without relying excessively on the particular architecture. The peak performance of a machine can usually be achieved only with careful implementation of an algorithm that exploit the architecture. This is especially true in the implementation of data communication algorithms whose efficiency depends critically on the effective utilization of the communication network. However, performance close to the peak can be achieved by an implementation that is not tied to a particular architecture. We discuss in this section some general, unifying concepts of parallel computing. The description of parallel architectures at the level treated here suffices for the discussion of the parallel implementations of algorithms in the rest of the book.

#### Single Program Multiple Data (SPMD) Architectures

The distinction between MIMD (multiple instructions, multiple data) and SIMD (single instruction, multiple data) parallel architectures gives way to the unified notion of *Single Program Multiple Data* (SPMD). An SPMD architecture permits the operation of a single program on multiple sets of problem data. Each data set is operated upon by a different processor. The requirement that the program executes synchronously—as in the case of an SIMD architecture—is now relaxed.

The repeated use of a linear programming solver, for example, to solve multiple related instances of a linear program—e.g., in the analysis of multiple scenarios of prices for portfolio optimization models—can be implemented naturally on an SPMD machine. Whether all solvers execute exactly the same steps synchronously or not is irrelevant to the programmer. An SIMD computer would require each solver to execute exactly the same sequence of steps and could be inefficient if each solver should execute different steps, while an MIMD computer would be more efficient in executing this application. However, the analyst who views this as an application on an SPMD architecture need not be concerned whether the computer is, in Flynn's taxonomy, SIMD or MIMD.

#### Parallel Virtual Memory

The distinction between distributed memory and shared memory is not crucial to the applications programmer with the introduction of the notion of *parallel virtual memory.* This idea originated with software systems built to harness the computing power of distributed networks of workstations. Such systems (see Section 1.5 for references) enable users to link together multiple, possibly heterogeneous, computers from a distributed network and apply them to the execution of a single program.

In a distributed memory system each processor has its own local memory and information among processors is exchanged using messages passed along some communication network. In a shared memory system all processors have access to the same memory and information is exchanged by reading or writing data to this common memory. Care should be taken to guarantee the integrity of data when multiple processors need to write to the same memory location.

Using the notion of parallel virtual memory we can view a distributed architecture as a system with shared virtual memory. The algorithm designer has to deal with the problem of data integrity to ensure correctness of the results, but the fact that memory is distributed is of no direct consequence. The developer of the distributed memory environment has to keep track of global objects, such as common variables. This can be achieved by distributing all global objects to local memories. An alternative would be to distribute the addresses of all global objects to local memories and let each processor fetch the data only when it is used by the application. These two approaches trade off efficiency in space with efficiency in time. No matter which way this issue is resolved by the developer of the architecture, the applications programmer is not concerned with the distinction between shared and distributed memory. Indeed, software systems like Express, LINDA or PVM (see Section 1.5 for references) have made it possible to port applications between shared memory machines, distributed memory machines, and distributed heterogeneous networks of workstations with little degradation of performance and no reprogramming effort. The success of these systems led to the development of a message-passing interface standard, the MPI. Figure 1.1 exemplifies these points by showing the performance of a Monte Carlo simulation model for the pricing of financial instruments on a variety of computers with different parallel architectures.

### Heterogeneous Distributed Computing

During the early stages of the development of parallel processing there was a quest for a "winner" parallel computer architecture. With the emergence of the unifying computing paradigms like SPMD and parallel virtual memory, the notion of *heterogeneous distributed computing* began to dominate as a concept for introducing parallelism. Instead of using a particular parallel machine, users may link together heterogeneous networks of workstations to create cost-effective forms of parallel processing. Fiber optics facilitate such links, providing data-transfer rates that were until recently found only on tightly coupled systems. Software systems like Express, LINDA or PVM allow users to decouple their problems in a parallel virtual memory environment. Multiple workstations can execute tasks from this environment, complete them at their own pace depending on the workload and performance of each workstation, and return the results.

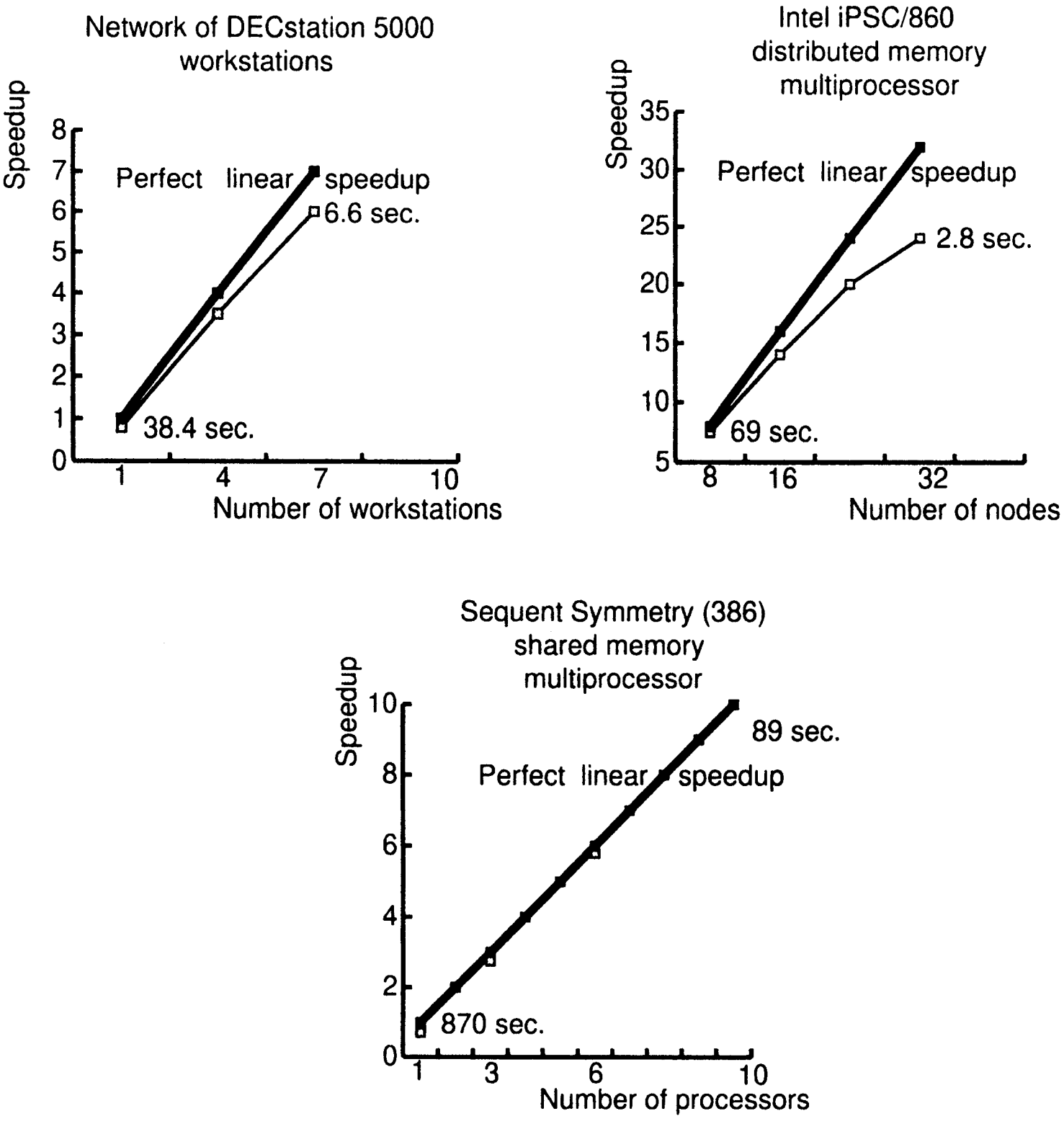

**Figure 1.1** The performance of a Monte Carlo simulation model for pricing financial instruments improves almost linearly on a variety of computer platforms with different parallel architectures, with either shared memory or distributed memory.

There is no need to assume that all the servers on the network are identical workstations or even just workstations. Some of the servers could be more advanced parallel architectures. For example, an array processor attached to the network could be used to execute linear algebra calculations, while workstations could be executing the less homogeneous operations.

### 1.1.3 Control and data parallelism

There are two important paradigms for parallel programming that we wish to discuss: *control parallelism* and *data parallelism*. These are models for programming applications on parallel machines and should not be confused with the models of parallel computer architectures discussed in the previous section. These paradigms facilitate the implementation of parallel algorithms on a broad range of parallel architectures. The architecture of the target machine—the interaction between instruction and data streams, memory organization, and communication network—are not explicitly taken into account by a programmer who uses these paradigms to implement an algorithm. Of course, depending on the complexity of the code and the amount of coordination required among processors, the resulting implementation may not be as efficient as it would have been if it were implemented using a paradigm closely linked to the hardware architecture. As parallel languages are continuously enhanced, and compilers become more sophisticated, such losses in efficiency become less noticeable.

### Control Parallelism

In control parallelism the instructions of a program are partitioned into sets of instructions that are independent of each other and therefore can be executed concurrently. A well-known type of control parallelism is *pipelining.* Consider, for example, the simple operation, commonly called SAXPY,

$$Y \leftarrow \alpha X + Y,$$

where $X = (X(j))_{j=1}^{n}$ and $Y = (Y(j))_{j=1}^{n}$ are $n$-dimensional vectors and $\alpha$ is a scalar. On a scalar uniprocessor computer this operation will be completed in $2n$ steps (assuming that the execution of an addition or a multiplication takes the same number of machine cycles). On a computer with two functional units (an adder and a multiplier), it is possible to pipeline the $\alpha X$ multiplication with the $+Y$ operation. As soon as the first product result $\alpha X(1)$ becomes available it can be added to $Y(1)$, while the multiplier simultaneously computes $\alpha X(2)$. Both the adder and the multiplier operate concurrently. If the adder and multiplier operate on scalars, then pipelining the two operations will complete the SAXPY in just $n + 1$ steps. The maximum improvement in performance (defined as *relative speedup*, see Definition 1.4.3 in Section 1.4) from the control parallel execution of SAXPY is by a factor of $2n/(n + 1) \approx 2$.

In general, the speedup achieved through control parallelism is limited by the number of operations of the algorithm that must be executed serially. If $s$ denotes the fraction of the code that is executed serially, then the speedup from parallelizing the remaining $p = 1 - s$ fraction using $P$ processors is bounded by Amdahl's law:

$$\text{Relative Speedup} \leq \frac{1}{s + p/P}.$$

The speedup that can be "achieved" with an infinite number of processors is bounded by $1/s$.

### Data Parallelism

In data parallelism, instructions are performed on many data elements simultaneously by many processors. For example, with $n$ processors the SAXPY operation can be completed in two steps: one step executes the product $\alpha X$ for all components of $X$ simultaneously and the next step executes $+Y$ on all components of $Y$. The speedup with the data-parallel execution of the SAXPY operation is proportional to the size of the problem, i.e., to $n$.

The speedup achieved through data parallelism is restricted by the number of data elements that can be handled concurrently, while the performance of control parallelism is restricted by the serial fraction of the code. Data parallelism is, in general, favored over control parallelism as a parallel programming paradigm. The reason for this is that as problems get bigger the amount of data increases, while the complexity of the algorithm and the code remains unchanged. Hence, more processors can be utilized effectively for large-scale problems with a data-parallel paradigm.

Data parallelism creates *tasks* (see Section 1.2) by partitioning the problem data and distributing them to multiple processors. The actual number of processors on the target machine is not critical, since, if necessary, multiple tasks can be scheduled for execution on the same processor. We may assume that each data partition is distributed to a different *virtual processor* (VP), and that a physical processor may emulate multiple virtual processors (see Section 14.1). The mapping of virtual to physical processors depends on the communication patterns needed by the algorithm for the exchange of data, and on the communication network of the target machine. For example, two VPs that exchange data frequently should be mapped to two physical processors that are on adjacent nodes of the communication network, i.e., nodes of the network that are connected by a communication link. Two VPs whose data elements depend on each other, e.g., through an assignment operation such as $Y \leftarrow \alpha X$, should be mapped to the same physical processor. The mapping of virtual to physical processors is done by the compiler. Modern data-parallel programming languages usually contain primitives that allow the programmer to specify mappings that are efficient.

Early applications of data parallelism were done on SIMD architectures, such as the Connection Machines CM–1 and CM–2, the AMT DAP, and Masspar MP-1. These architectures facilitated the coordination of thousands of processors through a global control unit that broadcasts identical instructions to all processors. Data parallelism has recently been realized on MIMD architectures, such as the Connection Machine CM–5, so that each processor decodes its own instructions and executes a program, asynchronously, on subsets of the problem data. Synchronization only takes place when processors need to exchange data. But the need to exchange data is determined concisely by the compiler: it is needed every time an operation uses data that are mapped on different physical processors. An MIMD architecture that executes a data-parallel programming paradigm is a single program multiple data (SPMD) machine (see Section 1.1.2).

## 1.2 How Does Parallelism Affect Computing?

Mathematical theory and computational techniques always went hand-in-hand in the field of optimization. This has been especially true for large-scale applications. For small-scale problems the correctness of a theoretical algorithm is usually sufficient to guarantee that the underlying problem can be solved successfully. For large-scale problems theoretical correctness does not suffice. Can the algorithm be translated into a computer code that does not exceed the memory limits of available computers? Will the algorithm work when there are (small) numerical errors in the execution of its steps, or will the errors be amplified? Does the algorithm require an inordinate number of operations in order to arrive at the proximity of a solution, and hence will it exhaust the user's computer budget and patience?

These issues are important whether the algorithm is implemented on a uniprocessor or on a multiprocessor parallel computer. However, additional considerations become important when a parallel computer is used for the implementation. For example, an algorithm that requires a large number of operations to reach a solution could be more efficient than an algorithm that requires fewer operations, if the steps of the former can be executed concurrently. Hence, it is important to discuss the ways in which parallelism affects computations. In particular, the following issues need to be addressed when designing an algorithm for parallel implementation.

*Task partitioning:* This is the problem of partitioning the operations of the algorithm, or the data of the problem, into independent tasks that can be mapped on multiple processors. The operations of the algorithm can be partitioned into sets that are independent of each other and can proceed concurrently. The problem data can be partitioned into blocks without interdependencies, so that multiple blocks can be processed in parallel. The partitions of the operations of the algorithm, or the blocks of independent data, are called *tasks*. Task partitioning appears, at first glance, to be a problem for the programmer of an algorithm, not for the mathematician who is designing the algorithm. But this is not quite so. The problem of task partitioning is easier to solve for algorithms designed to have independent operations or algorithms that use only small subsets of the problem data at each step. Hence, the mathematician can assist the applications programmer in addressing the problem of task partitioning by designing suitable algorithms.

*Task scheduling:* The issue here is how to assign the tasks to one or more processors for simultaneous execution. This appears, once more, to be a problem for the programmer and not for the mathematician who is designing the algorithm. However, since parallel computers come in a variety of architectures, this problem can not be left only to the programmer. At the one extreme of parallel architectures we have coarse-grain, small-scale machines such as the CRAY Y-MP/4 and the IBM 3090-600. At the other extreme we have fine-grain, massively parallel systems such as the Connection Machine CM–2, the MassPar, and the DAP. Scheduling of tasks to processors, for this variety of architectures, can be done satisfactorily when the underlying theoretical algorithm is flexible. If the operations of the algorithm can be structured into as many independent tasks as the number of available processors then the scheduling problem can be resolved by the programmer. The role of the algorithm designer is therefore important in this respect.

*Task synchronization:* Here we address the question of specifying an order of execution of the tasks, and the instances during execution of the algorithm where information must be exchanged among the tasks, in order to ensure the correct progress of the iterations according to the algorithm. This problem can again be viewed, from a narrow perspective, as one of computer implementation. But algorithms whose convergence is guaranteed without excessive requirements for synchronization are likely to be more efficient when implemented in parallel.

We see that both the structure of the mathematical algorithm and the structure of the underlying application play a crucial role in parallel computations. It may be possible to decompose the problem into domains that can be solved independently. For example, in image reconstruction from projections we may consider the decomposition of a discretized cross section of the image into subsections. If the subsections are properly chosen the image could be reconstructed from these smaller sections.

This is the background on which this book is developed. We argue that it is not enough to just know what a parallel computer is. And for a mathematician, it is not enough to "develop parallel algorithms." In order to successfully attack a significant real-world problem we have to study carefully the modeling process, the appropriate parallel computing architecture, and develop algorithms that are mathematically sound and suitable for both.

We stress, however, that the design of an algorithm should not be linked to the specifics of a parallel computer architecture. Such an approach would render the algorithm obsolete once the architecture changes. Given also the diversity of existing parallel architectures, this approach would find limited applicability. A basic premise of this book is that algorithms should be flexible enough to facilitate their implementation on a range of architectures. Indeed, we present algorithms for several optimization problems

that—when viewed from the proper perspective—are flexible enough to be implementable on a broad spectrum of parallel machines. Examining the general structure of the iterative algorithms, we understand their suitability for parallel implementation, as well as their flexibility for adaptation to different architectures. A general classification of parallel algorithms is given next.

## 1.3 A Classification of Parallel Algorithms

A *parallel algorithm* is any algorithm that can be implemented on a parallel machine in a way that results in clock-time savings when compared to its implementation, for the solution of the same problem, on a serial computer. This property of an algorithm, i.e., of being parallel or not, stems, independently, from either its mathematical structure, or from the features of the problem to which it is applied, or from both. Even an algorithm that is mathematically sequential can become, for certain decomposable problems, a parallel algorithm. We look first at several mathematical structures of algorithms and discuss their parallelism (Section 1.3.1), and then turn to algorithms that become parallel due to problem structure (Section 1.3.2).

We deal with algorithms drawn from three broad algorithmic classes: (i) *iterative projection algorithms*, (ii) *model decomposition algorithms*, and (iii) *interior point algorithms*. These classes of algorithms exploit parallelism in fundamentally different ways. Algorithms in the first class do so mostly due to their mathematical structure, while algorithms in the second and third classes do so for problems with special structure. Iterative projection methods are the most flexible in exploiting parallelism since they do so by suitable reorganization of their operations. The reorganization of operations can be accomplished in many different ways, and a classification of the alternative approaches is given in Section 1.3.1. Model decomposition and interior point algorithms can exploit parallelism by taking advantage of the structure of the model of the problem that is being solved. Section 1.3.2 gives a general preview of the potential parallelism of these two classes of algorithms.

### 1.3.1 Parallelism due to algorithm structure: Iterative projection algorithms

We consider applications that can be modeled by a system of linear or nonlinear equations or inequalities, i.e., a system of the form

$$f_i(x) \star 0, \quad i = 1, 2, \ldots, m, \tag{1.1}$$

where $x \in \mathbb{R}^n$ ($\mathbb{R}^n$ is the $n$-dimensional Euclidean space), $f_i : \mathbb{R}^n \rightarrow \mathbb{R}$, for all $i$, are real-valued functions, and $\star$ stands for equality signs, inequality signs, or a mixture of such. This generic description allows us to classify a broad set of iterative procedures and discuss their suitability for parallel computations. The mathematical classification of algorithms we give next depends on how each algorithm makes use of *row-action* or *block-action* iterations. This classification is purely logical and independent of the architecture of the machine on which the algorithm will be implemented or of the structure of the problem to be solved.

Solving a *feasibility problem* (1.1) may be attempted if the system is *feasible*, i.e., if the set $\{x \in \mathbb{R}^n \mid f_i(x) \star 0, \;\; i = 1, 2, \ldots, \; m\}$ is not empty, or even if the system is infeasible. The mathematical model may also be set up as an *optimization problem* with an objective function $f_0 : \mathbb{R}^n \to \mathbb{R}$ imposed and (1.1) serving as constraints. The optimization problem can be written as:

$$\text{Optimize} \quad f_0(x) \tag{1.2}$$

$$\text{s.t.} \qquad f_i(x) \star 0, \qquad i = 1, 2, \ldots, m. \tag{1.3}$$

Here *optimize* stands for either minimize or maximize. Additional constraints of the form $x \in Q$ can also be imposed, where $Q \subseteq \mathbb{R}^n$ is a given subset describing, for example, box constraints. Specific examples of such problems are described, as they arise in several fields of applications, in the next chapters and references to others are given.

Special-purpose iterative algorithms designed to solve these problems may employ iterations that use in each iterative step a single row of the constraints system (1.1) or a group of rows. A *row-action iteration* has the functional form

$$x^{\nu+1} \doteq R_{i(\nu)}(x^{\nu}, f_{i(\nu)}), \tag{1.4}$$

where $\nu$ is the iteration index, and $i(\nu)$ is the *control index*, $1 \leq i(\nu) \leq m$, specifying the row that is acted upon by the algorithmic operator $R_{i(\nu)}$. The algorithmic operator generates, in some specified manner, the new iterate $x^{\nu+1}$ from the current iterate $x^{\nu}$ and from the information contained in $f_{i(\nu)}$. $R_{i(\nu)}$ may depend on additional parameters that vary from iteration to iteration, such as relaxation parameters, weights, tolerances, etc.

The system (1.1) may be decomposed into $M$ groups of constraints (called *blocks*) by choosing integers $\{m_t\}_{t=0}^{M}$ such that

$$0 = m_0 < m_1 < \cdots < m_{M-1} < m_M = m, \tag{1.5}$$

and defining for each $t$, $1 \leq t \leq M$, the subset

$$I_t \doteq \{m_{t-1} + 1, \;\; m_{t-1} + 2, \ldots, m_t\}. \tag{1.6}$$

This yields a partition of the set:

$$I \doteq \{1, 2, \ldots, \; m\} = I_1 \cup I_2 \cup \cdots \cup I_M. \tag{1.7}$$

A *block-action iteration* has then the functional form

$$x^{\nu+1} \doteq B_{t(\nu)}(x^\nu, \{f_i\}_{i\in I_{t(\nu)}}), \tag{1.8}$$

where $t(\nu)$ is the control index, $1 \le t(\nu) \le M$, specifying the block which is used when the algorithmic operator $B_{t(\nu)}$ generates $x^{\nu+1}$ from $x^\nu$ and from the information contained in all rows of (1.1) whose indices belong to $I_{t(\nu)}$. Again, additional parameters may be included in each $B_{t(\nu)}$, such as those specified for the row-action operator $R_{i(\nu)}$. The special-purpose iterative algorithms that we consider in our classification scheme may address any problem of the form (1.1) or (1.2)–(1.3) and may be classified as having one of the following four basic structures.

*Sequential Algorithms.* For this class of algorithms we define a control sequence $\{i(\nu)\}_{\nu=0}^{\infty}$ and the algorithm performs, in a strictly sequential manner, row-action iterations according to (1.4), from an appropriate initial point until a stopping rule applies.

*Simultaneous Algorithms.* Algorithms in this class first execute simultaneously row-action iterations on all rows

$$x^{\nu+1,i} \doteq R_i(x^\nu, f_i), \qquad i = 1, 2, \ldots, m, \tag{1.9}$$

using the same current iterate $x^\nu$. The next iterate $x^{\nu+1}$ is then generated from all intermediate ones $x^{\nu+1,i}$ by an additional operation

$$x^{\nu+1} \doteq S(\{x^{\nu+1,i}\}_{i=1}^{m}). \tag{1.10}$$

Here $S$ and $R_i$ are algorithmic operators, the $R_i$'s being all of the row-action type.

*Sequential Block-Iterative Algorithms.* Here the system (1.1) is decomposed into fixed blocks according to (1.7), and a control sequence $\{t(\nu)\}_{\nu=0}^{\infty}$ over the set $\{1, 2, \ldots, M\}$ is defined. The algorithm performs sequentially, according to the control sequence, block iterations of the form (1.8).

*Simultaneous Block-Iterative Algorithms.* In this class, block iterations are first performed using the same current iterate $x^\nu$, on all blocks simultaneously

$$x^{\nu+1,t} \doteq B_t(x^\nu, \{f_i\}_{i\in I_t}), \qquad t = 1, 2, \ldots, M. \tag{1.11}$$

The next iterate $x^{\nu+1}$ is then generated from the intermediate ones $x^{\nu+1,t}$ by

$$x^{\nu+1} \doteq S(\{x^{\nu+1,t}\}_{t=1}^{M}). \tag{1.12}$$

Here $S$ and $B_t$ are algorithmic operators, the $B_t$'s being block-iterative (1.8).

An important generalization of block-iterative algorithms refers to *variable block-iterative* algorithms, wherein block iterations are performed according to formula (1.8). The blocks $I_{t(\nu)}$ however, are not determined a priori according to (1.7) and kept fixed, but the algorithm employs an infinite sequence $\{I_{t(\nu)}\}_{\nu=0}^{\infty}$ of nonempty blocks $I_{t(\nu)} \subseteq I$. This means that as the iterations proceed, both the sizes of the blocks and the assignment of constraints to blocks may change in a much more general manner than the scheme with fixed blocks. Of course, convergence theory of such an algorithm could still impose certain restrictions on how the sets $I_{t(\nu)}$ may be constructed.

Although not formulated as such before, many special-purpose iterative algorithms that were proposed and implemented over the years can be identified according to this classification. As we proceed through the book these general algorithmic structures will become clearer, with the presentation of specific algorithms.

### Parallel Computing with Iterative Algorithms

We examine two approaches for introducing parallelism in both the *sequential* and the *simultaneous* algorithms described above. We address in particular the problem of task partitioning, by identifying operations of the algorithm that can be executed concurrently.

A sequential row-action algorithm uses information contained in a single row $f_{i(\nu)}(x) \star 0$ and the current iterate $x^{\nu}$ in order to generate the next iterate $x^{\nu+1}$. The control sequence $\{i(\nu)\}_{\nu=0}^{\infty}$ specifies the order in which rows are selected. Is it possible to select rows in such a way that two (or more) rows can be operated upon simultaneously without altering the mathematical structure of the algorithm? This can be achieved if at every iteration only different components of the updated vector are changed. If the algorithmic operators $R_{i(\nu)}$ and $R_{i(\nu+1)}$, which are vector-to-vector mappings, use different components of $x^{\nu}$ and $x^{\nu+1}$ to operate upon and leave all other components unchanged, then the operations can proceed concurrently. Identifying row-indices $i(\nu)$ so that the operators can be applied concurrently on such rows depends on the structure of the family $\{f_i(x)\}_{i=1}^{m}$. We will see that several important problems have a structure that allows us to identify such row-indices. The parallel algorithm is, with this approach, mathematically identical to the serial algorithm.

The parallelism in the simultaneous algorithms is obvious. The iterations (1.9) can be executed concurrently using up to $m$ parallel processors. The step (1.10) is a synchronization step where the processors must cooperate, and exchange information contained in the $m$ vectors $\{x^{\nu+1,i}\}_{i=1}^{m}$ in order to compute the next iterate $x^{\nu+1}$.

The block-iterative algorithms—sequential or simultaneous—lead to parallel computing identical to the two ways outlined above. However, the introduction of blocks permits more flexibility in handling the task schedul-ing problem. The blocks can be chosen so as to ensure that all processors receive tasks of the same computational difficulty. Hence, all processors complete their work in roughly the same amount of time. Delays while processors wait for each other are minimized. The specification of block-sizes will typically depend on the computer architecture. For fine-grain, massively parallel systems it is preferable to create a very large number of small blocks. Fewer, but larger, blocks are needed for coarse-grain parallel computers.

### 1.3.2 Parallelism due to problem structure: Model decomposition and interior point algorithms

We consider now the solution of *structured* optimization problems. We introduce first the following definition.

**Definition 1.3.1 (Block Separability)** *A function $F$ is block-separable if it can be written as the sum of functions of subsets of its variables. Each subset is called a block.*

Consider, for example, the vector $\boldsymbol{x} \in \mathbb{R}^{nK}$, which is the concatenation of $K$ subvectors, each of dimension $n$,

$$\boldsymbol{x} = \left((x^1)^{\mathrm{T}}, (x^2)^{\mathrm{T}}, \ldots, (x^K)^{\mathrm{T}}\right)^{\mathrm{T}},$$

where $x^k \in \mathbb{R}^n$, for all $k = 1, 2, \ldots, K$. (Boldface letters denote vectors in the product space $\mathbb{R}^{nK}$). Then $F : \mathbb{R}^{nK} \to \mathbb{R}$ is block-separable if it can be written as $F(\boldsymbol{x}) \doteq \sum_{k=1}^{K} f_k(x^k)$, where $f_k : \mathbb{R}^n \to \mathbb{R}$.

We consider now the following constrained optimization problem:

**Problem** $[\mathcal{P}]$

$$\text{Optimize} \quad F(\boldsymbol{x}) \doteq \sum_{k=1}^{K} f_k(x^k) \tag{1.13}$$

$$\text{s.t.} \quad x^k \in X_k \subseteq \mathbb{R}^n, \qquad \text{for all } k = 1, 2, \ldots, K, \tag{1.14}$$

$$\boldsymbol{x} \in \Omega \subseteq \mathbb{R}^{nK}. \tag{1.15}$$

The block-separability of the objective function and the structure of the constraint sets of this problem can be exploited for parallel computing. Two general algorithms that exploit this structure are described below.

#### Model Decomposition Algorithms

If the product set $X_1 \times X_2 \times \cdots \times X_K$ of problem $[\mathcal{P}]$ is a subset of $\Omega$, then the problem can be solved by simply ignoring the constraints $\boldsymbol{x} \in \Omega$. A solution can be obtained by solving $K$ independent subproblems in each of the $x^k$ vector variables. These subproblems can be solved in parallel.

In general, however, the *complicating constraints* $\boldsymbol{x} \in \Omega$ cannot be ignored. A model decomposition applies a *modifier* to problem $[\mathcal{P}]$ to obtain a problem $[\mathcal{P}']$ in which the complicating constraints are not explicitly present. It then employs a suitable algorithm to solve $[\mathcal{P}']$. Since the modified problem does not have complicating constraints it can be decomposed into $K$ problems, each one involving only the $x^k$ vector and the constraint $x^k \in X_k$, for $k = 1, 2, \ldots, K$. These problems are not only smaller and simpler than the original problem $[\mathcal{P}]$, but they are also independent of each other, and can be solved in parallel. If the solution to the modified problem $[\mathcal{P}']$ is "sufficiently close" to a solution of the original problem $[\mathcal{P}]$ then the process terminates. Otherwise the current solution is used to construct a new, modified problem, and the process repeats until some stopping criterion is satisfied (see Figure 1.2). A wide variety of algorithms fall under this broad framework, and specific instances are discussed in Chapter 7.

Of particular interest is the application of model decomposition algorithms to problems with linear constraints, with a constraint matrix that is block-angular or dual block-angular. For example, the constraint matrix of the classic *multicommodity network flow* problem (Section 12.2.2) is of the following, block-angular form

$$A = \begin{pmatrix} A_1 & & & \\ & A_2 & & \\ \vdots & & \ddots & \\ & & & A_K \\ I & I & \ldots & I \end{pmatrix}. \tag{1.16}$$

Here $\{A_k\}$, $k = 1, 2, \ldots, K$, are $m \times n$ network flow constraints matrices, and $I$ denotes the $n \times n$ identity matrix. The constraint set of this problem can be put in the form of problem $[\mathcal{P}]$ by defining $X$ as the product set $X \doteq \prod_{k=1}^{K} X_k$, where $X_k = \{x^k \mid A_k x^k = b^k,\ 0 \leq x^k \leq u^k\}$, where $b^k$ denotes the right-hand-side of the equality constraints and $u^k$ are vectors of upper bounds on the variables, and by defining $\Omega \doteq \{\boldsymbol{x} \mid \boldsymbol{x} \in X,\ E\boldsymbol{x} \leq U\}$, where $E = [I \mid I \mid \ldots \mid I]$. A model decomposition algorithm can now be applied to solve the multicommodity flow problem (see Section 12.5).

The constraint matrix of linear programming formulations of two-stage stochastic programming problems (Section 13.3.3) is of the dual block-angular form

$$A = \begin{pmatrix} A_0 & & & \\ T_1 & W_1 & & \\ \vdots & & \ddots & \\ T_N & & & W_N \end{pmatrix}. \tag{1.17}$$

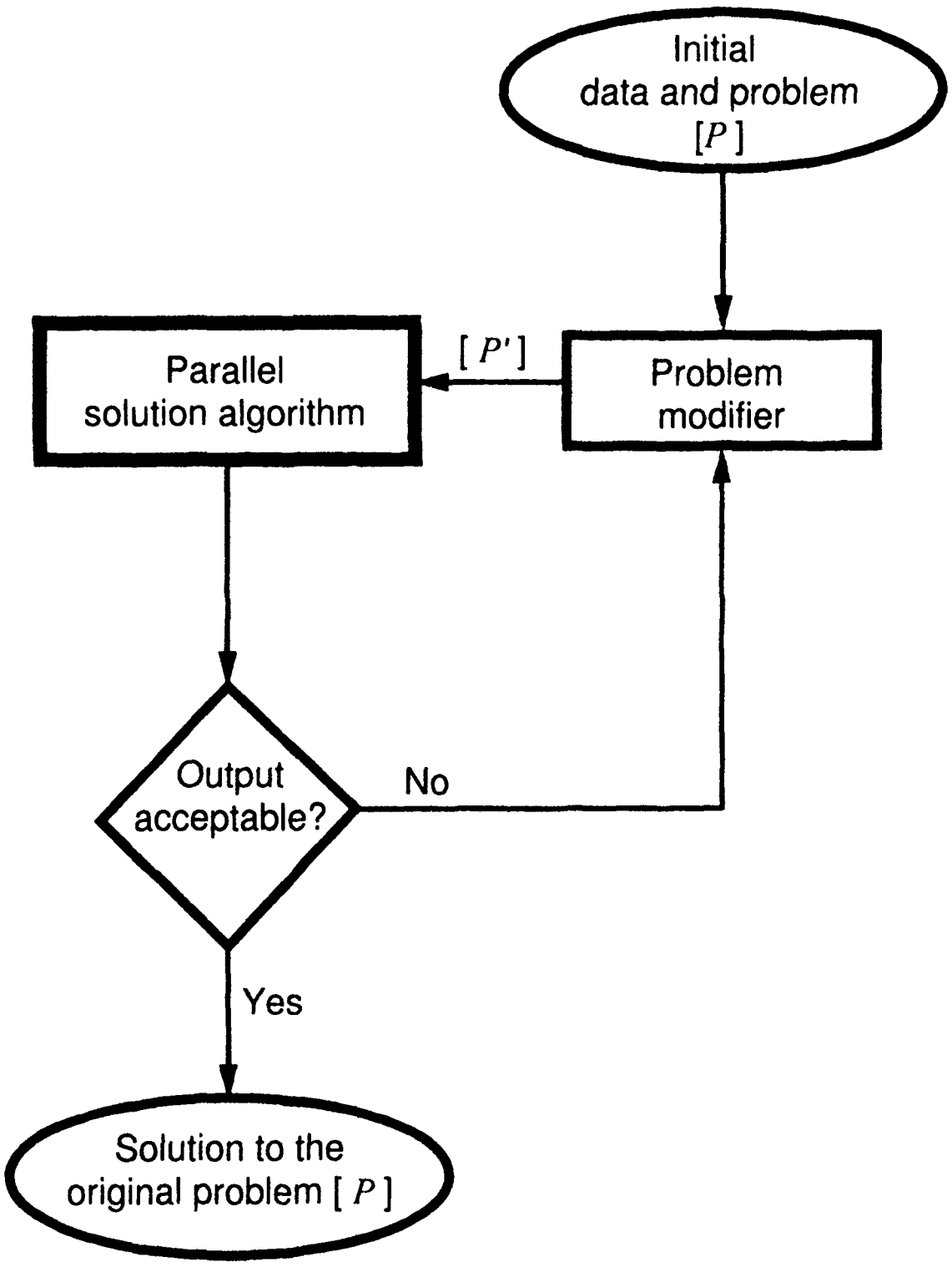

**Figure 1.2** Structure of model decomposition algorithms: Modifications $[\mathcal{P}']$ of the original problem $[\mathcal{P}]$ are solved repeatedly using a parallel solution algorithm. The solution of each modified problem is used to construct the next modification until an acceptable solution to the original problem is reached.

Here $A_0$ is an $m_0 \times n_0$ constraint matrix, and $\{T_l\}$, $l = 1, 2, \ldots, N$ are $m_l \times n_0$, and $\{W_l\}$, $l = 1, 2, \ldots, N$ are $m_l \times n_l$, constraint matrices. The dual formulation of this problem, which is formulated using the transpose $A^{\mathrm{T}}$ of the constraint matrix, can be put in the form of $[\mathcal{P}]$ using definitions for $X$ and $\Omega$ similar to those introduced above for the multicommodity flow problem.

### Interior Point Algorithms

Consider now special cases of problem $[\mathcal{P}]$, when the objective function is linear or separable quadratic, and the constraints sets $X$ and $\Omega$ are defined by linear equalities and inequalities. Then we can solve the problem directly using a general-purpose interior point algorithm for linear or quadratic programs (see Chapter 8).

The major computational effort of these algorithms appears in the factorization and solution of a system of equations, for the unknown vector $\Delta y$, of the form

$$(A\Theta A^{\mathrm{T}})\Delta y = \psi, \tag{1.18}$$

where $A$ is the constraint matrix of the problem, $\psi$ is a given vector, $\Theta$ is a given diagonal matrix, and $\Delta y$ denotes a dual step direction of the interior point algorithm.

When the constraint matrix $A$ is of the structure shown in (1.16) or (1.17) it is possible to design special-purpose procedures that solve the system of equations (1.18) by solving subsystems involving the smaller matrices $A_1, \ldots, A_K$, or $A_0, W_1, \ldots, W_N, T_1, \ldots, T_N$ (see Section 8.3). These subsystems are independent of each other and can be solved in parallel. Hence, while interior point algorithms are not parallel algorithms per se they can exploit parallel architectures when executing the numerical linear algebra calculations for problems with special structures.

## 1.4 Measuring the Performance of Parallel Algorithms

How does one evaluate the performance of an algorithm executing on a parallel computer? There is no single simple answer to this question. Instead, several measures of performance are usually employed. Depending on the objective of the experimental study, different measures may be appropriate. In this section we define several such measures.

**Definition 1.4.1 (FLOPS)** *The number of floating-point operations per second executed by an implementation of an algorithm is called FLOPS.*

Observed FLOPS rates of a computer program are considered the yardstick for the evaluation of algorithms on supercomputers. The prefixes M, G and T are used to denote MFLOPS (Mega-FLOPS or $10^6$ FLOPS), GFLOPS (Giga-FLOPS or $10^9$ FLOPS), and TFLOPS (Tera-FLOPS or $10^{12}$ FLOPS), respectively.

**Definition 1.4.2 (Speedup)** *Speedup is the ratio of the solution time for a problem with the fastest available serial code executed on a single processor, to the solution time with the parallel algorithm executed on multiple processors. The sequential algorithm is executed on one of the processors of the parallel system.*

**Definition 1.4.3 (Relative Speedup)** *Relative speedup is the ratio of the solution time for a problem with a parallel algorithm executed on a single processor of the parallel computer, to the solution time with the same algorithm when executed on multiple processors.*

**Definition 1.4.4 (Absolute Speedup)** *Absolute speedup is the ratio of the solution time for a problem with the fastest available serial code executed on* *the fastest available serial computer, to the solution time with the parallel algorithm executed on multiple processors of the parallel computer.*

**Definition 1.4.5 (Efficiency)** *Efficiency is the ratio of relative speedup to the number of processors.*

Efficiency provides a way to measure the performance of an algorithm independently of the number of processors on the particular computer. *Linear speedup* is observed when a parallel algorithm on $P$ processors runs $P$ times faster than on a single processor. Linear speedup corresponds to efficiency of 1. *Sublinear speedup*, or efficiency less than 1, is achieved when the improvement in performance is less than $P$. Sublinear speedup is common due to the presence of sequential segments of code, delays for processor synchronization, overhead in spawning independent tasks, and so on. *Superlinear speedup*, with efficiency exceeding 1, is unusual. It indicates that the parallel algorithm follows a different, more efficient solution path than the sequential algorithm. It is often possible in such situations to improve the performance of the sequential algorithm based on insights gained from the parallel algorithm.

**Definition 1.4.6 (Amdahl's law)** *Amdahl's law gives an upper bound on the relative speedup achieved on a system with P processors:*

$$\textit{Relative Speedup} \leq \frac{1}{s + p/P}, \tag{1.19}$$

*where s is the fraction of the execution time of the application performed serially, and p is the fraction of execution time that is performed in parallel ($s+p = 1$). An upper bound on relative speedup, "achieved" with an infinite number of processors, is $1/s$.*

It has been observed that in practice, 10 percent of the execution time of typical applications cannot be parallelized (this is the time spent in data input/output, initialization, serial bottlenecks of the algorithm). Amdahl's law then leads us to the conclusion that a speedup of ten is the best one could expect from the use of parallel machines. This is a very limited view of parallelism. Parallel computers are not used exclusively to solve existing applications faster; rather multiple processors are often used to solve in constant time much bigger instances of the same problem. This leads to the need for the following definitions.

**Definition 1.4.7 (Scalability)** *Scalability is the ability of a parallel algorithm to maintain a constant level of efficiency as the number of processors increases, by solving problems whose size increases in proportion to the number of processors.*

When an application is scaled in size to fit the larger number of processors, the serial part usually remains unchanged. It is the parallel part that scales up, and linear speedup can be expected. Consider a problem with a fraction $s$ of serial execution time and a parallelizable part that consumes a fraction $p$ of the execution time. When a $P$ processors system becomes available the application would be scaled. If the parallel part would scale linearly with the number of available procesors then the execution time for the larger application would be scaled by $s + pP$.

**Definition 1.4.8 (Scaled Speedup)** *A modified Amdahl's law that takes into account the scaling of the problem to fit in the parallel machine gives the scaled speedup by:*

$$\textit{Scaled Speedup} = s + pP, \tag{1.20}$$

*where $s$ is the fraction of the execution time of the application performed serially, and $p$ is the fraction of execution time that is performed in parallel $(s + p = 1)$.*

All measures presented above provide useful information on the performance of a parallel algorithm, but each provides only a partial measure of performance. High relative speedup, for example, with an inherently slow algorithm is of little use when an alternative algorithm converges much faster to a solution when executing on a uniprocessor. Similarly, a high FLOPS rate with an algorithm that performs excessive calculations is not as meaningful as a reduced FLOPS rate of an algorithm that requires fewer operations.

## 1.5 Notes and References

Textbook treatments are available for several aspects of parallelism, including parallel computer architectures, parallel languages, parallel programming, and the use of parallel programming in several areas of application. We mention first some general references. For treatment of hardware aspects see, for example, Almasi and Gottlieb (1994). For parallel programming paradigms see Leighton (1992), Lester (1993), and Kumar et al. (1994). The development of programs for parallel machines is discussed by Carriero and Gelernter (1992). For alternative viewpoints on parallel computing in general, see Buzbee and Sharp (1985), Deng, Glimm and Sharp (1992), and the special issue of the journal *Daedalus* (1992). A collection of references on parallel computing has been compiled by the Association for Computing Machinery (1991).

For perspectives on parallel computing for optimization problems, see Kindervater, Lenstra, and Rinnooy Kan (1989), and Zenios (1994b). Advanced treatment of numerical methods for parallel and distributed computing, including extensive coverage of optimization problems, is given in Bertsekas and Tsitsiklis (1989). Special issues of journals focusing on parallel optimization were edited by Meyer and Zenios (1988), Mangasarian and Meyer (1988, 1991, 1994), and Rosen (1990). Pardalos, Phillips, and Rosen (1992) treat some aspects of mathematical programming with respect to their potential for parallel implementations. Extensive lists of references on parallel computing for optimization are collected in Schnabel (1985) and Zenios (1989).

Ortega (1988) provides an introduction to the parallel and vector solution of linear systems. Methods of parallel computing for global optimization problems are collected in the volume by Pardalos (1992), see also Pardalos (1986). A discussion of parallel methods for unconstrained optimization and for the solution of nonlinear equations is found in Schnabel (1995)

**1.1** Extensive treatment of parallel computer architectures can be found in Hwang and Briggs (1984) and Hwang (1993). Kumar et al. (1994) also contains material on hardware aspects. A historical note on the development of parallel architectures is given by Desrochers (1987).

**1.1.1** The taxonomy of parallel architectures based on the interaction of instructions and data streams is due to Flynn (1972). Vector processing is discussed in Desrochers (1987) and Hwang and Briggs (1984). The classification of parallel computers based on the organization of memory was given by Bell (1992).

The use of array processors for linear programming was first discussed by Pfefferkorn and Tomlin (1976), and the use of vector computing for the solution of large-scale nonlinear optimization problems by Zenios and Mulvey (1986a, 1988b).

**1.1.2** SPMD was motivated by the computing paradigm of SIMD architectures, as introduced by Hillis (1985) for the Connection Machine (see also Blelloch 1990). The software systems that motivated the ideas of parallel virtual memory, and facilitated the use of distributed networks of heterogeneous workstations, are LINDA from Scientific Computing Associates, see Carriero and Gelernter (1992), Express from Parasoft Corporation (1990), and PVM of Beguelin et al. (1991). The message-passing interface standard MPI is described in a report published by the University of Tennessee (1995). The graphs of Figure 1.1 are from Cagan, Carriero and Zenios (1993), and were obtained using LINDA.

**1.1.3** In his PhD thesis Hillis (1985) suggested a machine for implementing data parallel computing. His thesis provided the foundation for the development of a commercial massively parallel machine, the Connection Machine CM–1, and later the CM–2. Hillis and Steele (1986) discuss data parallel algorithms.

**1.3** An extensive list of references for iterative algorithms is given in Sections 5.11 and 6.10; for model decomposition algorithms in Section 7.3; and for interior point algorithms in Section 8.4. Extensive discussion on the use of block-separability for the solution of large-scale optimization problems can be found in Griewank and Toint (1982). Block-separability was exploited, in the context of parallel computing, by Lescrenier (1988), Lescrenier and Toint (1988), and Zenios and Pinar (1992). Problems with block- or dual block-angular structure appear in many applications domains, such as those discussed in Chapters 12 and 13. However, it might be possible to detect a block-angular structure even in more general optimization problems using the techniques developed by Ferris and Horn (1994).

**1.4** Performance metrics for parallel algorithms are discussed in papers by Jordan (1987) and Dritz and Boyle (1987). Barr and Hickman (1993) review performance metrics and develop guidelines for reporting computational experiences with parallel algorithms. See also Section 15.1. A good discussion of scalability is given in Kumar et al. (1994). Amdahl's law was proposed in Amdahl (1967), and its modified variant incorporating the effects of scaling is given in Gustafson (1988), who attributes it to E. Barsis from Sandia National Laboratory, USA.

# Part I

# THEORY

The most practical thing in the world is a good theory.

Attributed to H. von Helmholtz

CHAPTER 2

# Generalized Distances and Generalized Projections

Many of the algorithms that we study in this book employ projections onto convex sets in order to achieve either feasibility or optimization. A projection of a given point onto a convex set is defined as another point which has two properties. First, it belongs to the set onto which the projection operation is performed and, second, it renders a minimal value to the distance between the given point and any point of the set. If the Euclidean distance $\| x - y \|$ is used in this context then the projection is called a metric projection. But it turns out to be very useful to employ additional means to measure the "distance" between the vectors $x$ and $y$. The main purpose of this chapter is to present a theory of generalized distances and generalized projections, which include the classic metric projections as a special case. These generalized distances and projections were first proposed by Bregman in 1967, from whose work the specific, more focused version given here was subsequently developed.

A *generalized distance* is a real-valued nonnegative function of two vector variables $D(x, y)$ defined in a specific manner, which may be used to "measure the distance" between $x$ and $y$ in some generalized sense. When defining generalized distances it is customary to allow nonsymmetry, i.e., $D(x, y) \neq D(y, x)$ and not insist on the triangle inequality that a traditional distance function must obey. Because of the lack of symmetry such distances are sometimes referred to as *directed distances.* Under certain conditions one can use a generalized distance to define projections onto convex sets, by minimizing $D(x, y)$ over all $x$ in the set for some fixed $y$ whose projection onto the set is sought. In this chapter we present a family of generalized distances and their associated *generalized projections* that serve as basic tools in the subsequent developments of parallel algorithms in optimization theory.

In Section 2.1 we define Bregman functions and associate with them generalized distances. Generalized projections are also defined and examples are given. Of particular interest are generalized projections onto hyperplanes which are characterized and discussed in Section 2.2. In Section 2.3 a characterization is given of Bregman functions, whose zone is the whole space. In Section 2.4 we present an important result which can be viewed as a generalized triangle inequality for generalized distances (Theorem 2.4.1), and from which we derive the well-known Kolmogorov criterion (Theorem 2.4.2). Section 2.5 gives a short discussion of Csiszár's $\varphi$-divergences, which are also generalized distances, but different from the ones we study. Notes and references to the literature are given in Section 2.6.

## 2.1 Bregman Functions and Generalized Projections

Let $S$ be a nonempty, open, convex set, such that $\overline{S} \subseteq \Lambda$, i.e., its closure $\overline{S}$ is contained in $\Lambda$, where $\Lambda$ is the domain of a function $f : \Lambda \subseteq \mathbb{R}^n \to \mathbb{R}$. Assume that $f(x)$ has continuous first partial derivatives at every $x \in S$, and denote by $\nabla f(x)$ its gradient at $x$. From $f(x)$ construct the function $D_f(x,y)$, where $D_f : \overline{S} \times S \subseteq \mathbb{R}^{2n} \to \mathbb{R}$, by

$$D_f(x,y) \doteq f(x) - f(y) - \langle \nabla f(y), x - y \rangle. \tag{2.1}$$

This function is called *generalized distance function* or *D-function.* $D_f(x,y)$ may be interpreted as the difference $f(x) - h(x)$ where $h(z)$ represents the hyperplane $H$, which is tangent to the epigraph of $f$ at the point $(y, f(y))$ in $\mathbb{R}^{n+1}$ (Figure 2.1 illustrates the situation geometrically). This can be seen through elementary convex analysis. The epigraph of $f$, epi$(f)$, is the set epi$(f) \doteq \{(x,u) \mid x \in \mathbb{R}^n,\ u \in \mathbb{R}, f(x) \leq u\}$ of points on and above the graph of $f$ in $\mathbb{R}^{n+1}$. When $f$ is a convex function, epi$(f)$ is a convex set in $\mathbb{R}^{n+1}$. The normal vector $N$ to the set epi $(f)$ at the point $(y, f(y))$ on the graph of $f$ is given by $N \doteq (\nabla f(y), -1)$.

A nonvertical supporting hyperplane (see Definition 5.5.2) $H$ to the set epi$(f)$ at the point $(y, f(y))$ is given by

$$H = \{(z, h(z)) \mid h(z) = f(y) + \langle \nabla f(y), z - y \rangle\} \tag{2.2}$$

and we see that $D_f(x,y) = f(x) - h(x)$.

We adopt the following notation for the *partial level sets* of $D_f(x,y)$. For any $\alpha \in \mathbb{R}$,

$$\begin{aligned} L_1^f(y,\alpha) &\doteq \{x \in \overline{S} \mid D_f(x,y) \leq \alpha\}, \\ L_2^f(x,\alpha) &\doteq \{y \in S \mid D_f(x,y) \leq \alpha\}. \end{aligned} \tag{2.3}$$

We define below a family of functions that have certain common properties. The reason for imposing conditions $(i)$–$(vi)$ in the following definition under one heading is that these are precisely the conditions which are needed to ensure the applicability of the iterative methods of Chapters 5 and 6 for linearly constrained, equality, inequality and interval convex programming problems, when the function $D_f(x,y)$ has the form (2.1).

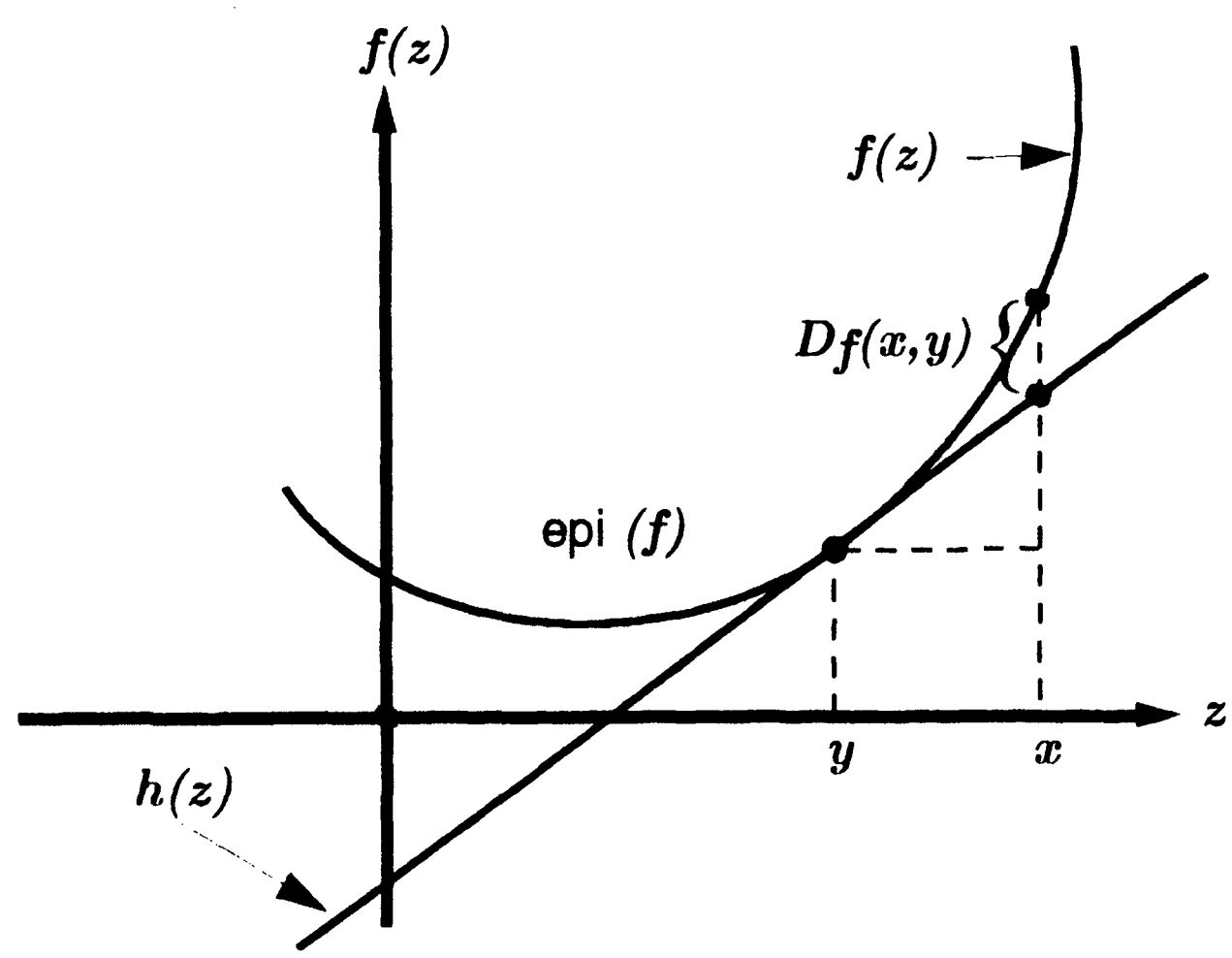

**Figure 2.1** Geometric interpretation of a generalized distance $D_f(x, y)$.

**Definition 2.1.1 (Bregman Functions)** *A function* $f : \Lambda \subseteq \mathbb{R}^n \rightarrow \mathbb{R}$ *is called a Bregman function if there exists a nonempty, open, convex set* $S$*, such that* $\overline{S} \subseteq \Lambda$ *and the following conditions hold:*

(*i*) $f(x)$ *has continuous first partial derivatives at every* $x \in S$,

(*ii*) $f(x)$ *is strictly convex on* $\overline{S}$,

(*iii*) $f(x)$ *is continuous on* $\overline{S}$,

(*iv*) *for every* $\alpha \in \mathbb{R}$, *the partial level sets* $L_1^f(y, \alpha)$ *and* $L_2^f(x, \alpha)$ *are bounded, for every* $y \in S$ *and for every* $x \in \overline{S}$, *respectively,*

(*v*) *If* $y^\nu \in S$, *for all* $\nu \geq 0$, *and* $\lim_{\nu \rightarrow \infty} y^\nu = y^*$ *then* $\lim_{\nu \rightarrow \infty} D_f(y^*, y^\nu) = 0$,

(*vi*) *If* $y^\nu \in S$ *and* $x^\nu \in \overline{S}$, *for all* $\nu \geq 0$, *and if* $\lim_{\nu \rightarrow \infty} D_f(x^\nu, y^\nu) = 0$, *and* $\lim_{\nu \rightarrow \infty} y^\nu = y^*$, *and* $\{x^\nu\}$ *is bounded, then* $\lim_{\nu \rightarrow \infty} x^\nu = y^*$.

We denote the family of Bregman functions by $\mathcal{B}(S)$, refer to the set $S$ as the *zone* of the function $f$, and write $f \in \mathcal{B}(S)$. We caution the reader that this definition of Bregman functions is not universal. Since it was first coined modifications were made to the list of conditions. Sometimes additional conditions are needed to achieve a certain goal, such as those found when dealing with proximal minimization (see Chapter 3). In other cases the conditions are modified to extend the family of Bregman functions so that more functions can be included (see Chapter 6). Recent research on Bregman functions in which other modifications of our definitions are used is mentioned in Section 2.6.

**Lemma 2.1.1** *For every* $f \in \mathcal{B}(S)$, *we have* $D_f(x, y) \geq 0$, *for all* $x \in \overline{S}$

*and all $y \in S$, and $D_f(x,y) = 0$ if and only if $x = y$.*

**Proof** This is a basic result in convex analysis which says that the assertion of this lemma is equivalent to the strict convexity of $f$ (Definition 2.1.1($ii$)), when $f$ is differentiable (Definition 2.1.1($i$)). See, e.g., Theorem 3.3.3 in Bazaraa and Shetty (1979). ■

We use now the generalized distance function (2.1) to define generalized projections.

**Definition 2.1.2 (Generalized Projections)** *Given $\Omega \subseteq \mathbb{R}^n$, $f \in \mathcal{B}(S)$, and $y \in S$, a point $x^* \in \Omega \cap \overline{S}$ for which*

$$\min_{z \in \Omega \cap \overline{S}} D_f(z,y) = D_f(x^*,y) \tag{2.4}$$

*(i.e., the minimum exists and the equality holds) is denoted by $P_\Omega(y)$ and is called a generalized projection (or simply, a projection) of the point $y$ onto the set $\Omega$.*

The next lemma guarantees the existence and uniqueness of generalized projections.

**Lemma 2.1.2** *If $f \in \mathcal{B}(S)$, then for any closed convex set $\Omega \subseteq \mathbb{R}^n$, such that $\Omega \cap \overline{S} \neq \emptyset$, and for any $y \in S$, there exists a unique generalized projection $x^* \doteq P_\Omega(y)$.*

**Proof** For any $w \in \Omega \cap \overline{S}$, the set

$$B \doteq \{z \in \overline{S} \mid D_f(z,y) \leq D_f(w,y)\} \tag{2.5}$$

is bounded (because of the boundedness of the partial level sets, Definition 2.1.1($iv$)) and closed (because $D_f(z,y)$ is continuous in $z$ on $\overline{S}$, by Definition 2.1.1($iii$)). Therefore, the nonempty (since $w \in B \cap \Omega$) set $T \doteq (\Omega \cap \overline{S}) \cap B$ is bounded. Since the intersection of closed sets is closed, $T$ is also closed, hence compact. Consequently, $D_f(z,y)$, as a continuous function of $z$, takes its infimum on the compact set $T$ at a point $x^* \in T$. For every $z \in \Omega \cap \overline{S}$ which lies outside $B$, we have $D_f(w,y) < D_f(z,y)$. Hence, $x^*$ satisfies (2.4) and is a generalized projection of $y$ onto $\Omega$.

To show uniqueness, suppose that there are two points $u, v \in \Omega \cap \overline{S}$, such that $u \neq v$ and $D_f(u,y) = D_f(v,y) = \min_{z \in \Omega \cap \overline{S}} D_f(z,y)$, for some $y \in S$. Then, $\frac{1}{2}(u+v) \in \Omega \cap \overline{S}$, because of the convexity of $\Omega \cap \overline{S}$. Since $f$ is strictly convex we have

$$\begin{aligned} D_f(\tfrac{1}{2}(u+v),y) &= f(\tfrac{1}{2}(u+v)) - f(y) - \langle \nabla f(y), \tfrac{1}{2}(u+v) - y\rangle \\ &< \tfrac{1}{2}f(u) + \tfrac{1}{2}f(v) - \tfrac{1}{2}f(y) - \tfrac{1}{2}f(y) \\ &\quad -\tfrac{1}{2}\langle \nabla f(y), u-y\rangle - \tfrac{1}{2}\langle \nabla f(y), v-y\rangle \\ &= \tfrac{1}{2}(D_f(u,y) + D_f(v,y)) = \min\nolimits_{z \in \Omega \cap \overline{S}} D_f(z,y), \end{aligned}$$

which is a contradiction because the generalized distance between $\frac{1}{2}(u+v)$ and $y$ cannot be strictly smaller than $\min_{z\in\Omega\cap\overline{S}} D_f(z,y)$. ■

We give now some examples. The first two deal with the most popular Bregman functions.

**Example 2.1.1** *The function*

$$f(x) = \frac{1}{2}\,\|x\|^2 \tag{2.6}$$

*is a Bregman function, with* $\Lambda = S = \overline{S} = \mathbb{R}^n$, *and* $D_f(x,y)$ *is then, according to (2.1),*

$$D_f(x,y) = \frac{1}{2}\,\|x-y\|^2\,. \tag{2.7}$$

This half-squared Euclidean norm of the difference has the symmetry property in $x$ and $y$, not present in other generalized distances. The generalized projection of any point $y$ onto any closed, convex set $\Omega$ is, in this case, the usual orthogonal (metric) projection.

**Example 2.1.2** *The "$x\log x$" entropy function,* $\operatorname{ent} x$, *maps the nonnegative orthant* $\mathbb{R}^n_+$ *into* $\mathbb{R}$ *according to*

$$\operatorname{ent} x \doteq -\sum_{j=1}^{n} x_j \log x_j. \tag{2.8}$$

*Here "log" denotes the natural logarithms and, by definition,* $0\log 0 \doteq 0$.

This function, known as Shannon's entropy, has its origins in information theory. It is easy to verify that for $f(x) = -\operatorname{ent} x$ we have

$$D_f(x,y) = \sum_{j=1}^{n} x_j \left(\log\left(\frac{x_j}{y_j}\right) - 1\right) + \sum_{j=1}^{n} y_j, \tag{2.9}$$

which is the well-known *Kullback-Leibler cross entropy function* from statistics.

The next lemma shows that the negative of $\operatorname{ent} x$ is indeed a Bregman function and that the associated $D_f$-function is given by (2.9).

**Lemma 2.1.3** *The function* $-\operatorname{ent} x$ *is a Bregman function with* $\Lambda = \mathbb{R}^n_+$ *and zone* $S_e$ *defined by*

$$S_e \doteq \{x \in \mathbb{R}^n \mid x_j > 0,\ \textit{for all } j = 1,2,\ldots,n\}. \tag{2.10}$$

**Proof** It is easy to show that conditions $(i)$ and $(ii)$ of Definition 2.1.1 are satisfied. Condition $(iii)$ holds, due to the convention $0\log 0 = 0$. For a given, fixed, $y$ let any component of $x$ go to $+\infty$. Then $D_f(x,y) \to +\infty$ as well and, therefore, we must conclude that for any $\alpha \in \mathbb{R}$, $L_1^f(y,\alpha)$ is bounded. A similar argument shows that $L_2^f(x,\alpha)$ is bounded, proving $(iv)$. Condition $(v)$ also follows from (2.9).

Assume now that the premises of condition $(vi)$ are satisfied. It is sufficient to show that any convergent subsequence of $\{x^\nu\}$ converges to $y^*$. Consider a general term of (2.9), namely $t : \mathbb{R}_+ \times (\mathbb{R}_+ \setminus \{0\}) \to \mathbb{R}$ defined by

$$t(x,y) \doteq x\left(\log\left(\frac{x}{y}\right) - 1\right) + y, \tag{2.11}$$

for $x \geq 0$ and $y > 0$. For any fixed $y$, $t(x,y) \geq 0$, for all $x \geq 0$, and

$$t(x,y) = 0 \text{ if and only if } x = y. \tag{2.12}$$

Now consider a convergent subsequence $\{x^{\nu_s}\}$ of $\{x^\nu\}$ and assume that $\lim_{s\to\infty} x^{\nu_s} = \overline{x}$. Then

$$\lim_{s\to\infty} D_f(x^{\nu_s}, y^{\nu_s}) = \lim_{s\to\infty} \sum_{j=1}^n t(x_j^{\nu_s}, y_j^{\nu_s}) = 0. \tag{2.13}$$

From the nonnegativity of $t(x,y)$ it follows that, for each $j$,

$$\lim_{s\to\infty} t(x_j^{\nu_s}, y_j^{\nu_s}) = 0. \tag{2.14}$$

Noting that $x_j^{\nu_s} \to \overline{x}_j$ and $y_j^{\nu_s} \to y_j^*$, we consider the following two cases. If $y_j^* > 0$, then (2.12) and (2.14) imply that $y_j^* = \overline{x}_j$. If $y_j^* = 0$, then (2.11) and (2.14) imply that $\overline{x}_j = 0$, and again $y_j^* = \overline{x}_j$. Hence $y^* = \overline{x}$. ■

**Example 2.1.3** *Denote the Cartesian product of $n$ intervals $(a,b)$ by $(a,b)^n$, let $S = (-1,1)^n$, and define*

$$f(x) \doteq -\sum_{j=1}^n \sqrt{1 - x_j^2}. \tag{2.15}$$

*Then it is straightforward to confirm that $f$ meets conditions $(i)$–$(vi)$ of Definition 2.1.1. By suitable scaling and translation of the argument, $f$ can be transformed so that its zone is the Cartesian product of a collection of $n$ arbitrary open bounded intervals.*

**Example 2.1.4** *If $f$ is a Bregman function, then $f(x) + \langle c, x\rangle + \beta$, for any $c \in \mathbb{R}^n$ and $\beta \in \mathbb{R}$, also satisfy the conditions of Definition 2.1.1 and produces exactly the same $D_f$-function, and therefore is also a Bregman function.*

## 2.2 Generalized Projections onto Hyperplanes

A key role in the iterative projection methods for linear feasibility problems and for linearly constrained optimization problems discussed in Chapters 5 and 6 is played by generalized projections onto hyperplanes. Therefore we take a closer look at these projections now. A *hyperplane* is a set of the form

$$H = \{x \in \mathbb{R}^n \mid \langle a, x \rangle = b\}, \tag{2.16}$$

where $a \in \mathbb{R}^n$ and $b \in \mathbb{R}$ are given. To facilitate our presentation we introduce the following definition.

**Definition 2.2.1 (Zone Consistency)** *(i) A function $f \in \mathcal{B}(S)$ with zone $S$ is said to have the zone consistency property with respect to the hyperplane $H$ if, for every $y \in S$, we have $P_H(y) \in S$. That is, the generalized projection of any point $y \in S$ onto the hyperplane $H$ remains in $S$.*
*(ii) $f \in \mathcal{B}(S)$ with zone $S$ is strongly zone consistent with respect to the hyperplane $H$ if it is zone consistent with respect to $H$, and with respect to every other hyperplane $H'$ which is parallel to $H$ and lies between $y$ and $H$.*

The following lemma characterizes generalized projections onto hyperplanes.

**Lemma 2.2.1** *Let $f \in \mathcal{B}(S)$, $H = \{x \mid \langle a, x \rangle = b\}$, and assume that $f$ is zone consistent with respect to $H$. For any given $y \in S$, the system*

$$\nabla f(x^*) = \nabla f(y) + \lambda a, \tag{2.17}$$

$$\langle a, x^* \rangle = b, \tag{2.18}$$

*determines uniquely the point $x^*$, which is the generalized projection of $y$ onto $H$. For a fixed representation of $H$ (i.e., for fixed $a \in \mathbb{R}^n$ and $b \in \mathbb{R}$), the system also determines uniquely the real number $\lambda$.*

**Proof** Consider the constrained minimization problem (2.4) with $\Omega = H$. Due to the zone consistency assumption, we look for a minimum within the open set $S$, and therefore we can write the Lagrangian of this problem as

$$L(z, y, \lambda) = f(z) - f(y) - \langle \nabla f(y), z - y \rangle - \lambda(\langle a, z \rangle - b), \tag{2.19}$$

and the necessary conditions for a point to be the generalized projection of $y$ onto $H$ are then

$$\nabla_z L(z, y, \lambda) = 0, \quad \nabla_\lambda L(z, y, \lambda) = 0, \tag{2.20}$$

from which (2.17)–(2.18) follow. To show the uniqueness of $x^*$ and of the parameter $\lambda$, for given $a$ and $b$, assume that $x^{**} \in \mathbb{R}^n$ and the parameter $\mu$ also solve (2.17)–(2.18), i.e.,

$$\nabla f(x^{**}) = \nabla f(y) + \mu a, \tag{2.21}$$

$$\langle a, x^{**} \rangle = b, \tag{2.22}$$

and that $x^* \neq x^{**}$. Then, multiplying (2.17) through by $(x^{**} - x^*)$ and using (2.18) and (2.22), we get

$$\langle \nabla f(x^*), x^{**} - x^* \rangle = \langle \nabla f(y), x^{**} - x^* \rangle. \tag{2.23}$$

Since $x^* \neq x^{**}$ we know, from (2.1) and Lemma 2.1.1, that

$$f(x^{**}) - f(x^*) > \langle \nabla f(x^*), x^{**} - x^* \rangle, \tag{2.24}$$

thus, combining (2.23) and (2.24),

$$f(x^{**}) - f(x^*) > \langle \nabla f(y), x^{**} - x^* \rangle. \tag{2.25}$$

Using similar arguments, after multiplying (2.21) through by $(x^* - x^{**})$, we obtain

$$f(x^*) - f(x^{**}) > \langle \nabla f(y), x^* - x^{**} \rangle. \tag{2.26}$$

The contradiction $0 > 0$ is then obtained by adding up (2.25) and (2.26). Consequently $x^* = x^{**}$, which also implies that $\lambda = \mu$. ■

For some fixed $a \in \mathbb{R}^n$ and $b \in \mathbb{R}$ representing the hyperplane $H$, the $\lambda$ obtained from the system (2.17)–(2.18) is called the *generalized projection parameter associated with the generalized projection of $y$ onto $H$*. It is also important to note, and easy to check that, for a prespecified $\lambda$ and for given $a$ and $y$, equation (2.17) uniquely determines $x^*$.

The next result is a statement about the signs of the parameters associated with generalized projections onto hyperplanes. It is fairly obvious in the case of orthogonal projections but needs verification for the generalized projections used here.

**Lemma 2.2.2** *Let $H = \{x \mid \langle a, x \rangle = b\}$ and $y \in S$. For any $f \in \mathcal{B}(S)$ which is zone consistent with respect to $H$, the parameter $\lambda$ associated with the generalized projection of $y$ onto some particular representation of $H$ (i.e., $a \neq 0$ and $b$ given), satisfies*

$$\lambda(b - \langle a, y \rangle) > 0, \quad \text{if } y \notin H, \tag{2.27}$$

$$\lambda = 0, \quad \text{if } y \in H. \tag{2.28}$$

**Proof** From (2.17)–(2.18) we get, if $x^* \neq y$,

$$\langle \nabla f(x^*) - \nabla f(y), x^* - y \rangle = \lambda \langle a, x^* - y \rangle.$$

But $\langle \nabla f(x^*) - \nabla f(y), x^* - y\rangle = D_f(x^*, y) + D_f(y, x^*)$, and also $\langle a, x^* - y\rangle = b - \langle a, y\rangle$. Therefore, if $y \notin H$ then

$$\lambda(b - \langle a, y\rangle) = D_f(x^*, y) + D_f(y, x^*), \tag{2.29}$$

and the result follows from Lemma 2.1.1. If $\langle a, y\rangle = b$, then $x^* = y$ is the unique generalized projection of $y$ onto $H$, thus, $\lambda = 0$. ■

**Corollary 2.2.1** *Let $H = \{x \mid \langle a, x\rangle = b\}$ and $f \in \mathcal{B}(S)$ which is zone consistent with respect to $H$. If $y^1 \in S$ and $y^2 \in S$ are in the two opposing half-spaces defined by $H$, i.e., $\langle a, y^1\rangle < b$ and $\langle a, y^2\rangle > b$, then the parameters $\lambda_1$ and $\lambda_2$ associated with the generalized projections of $y^1$ and $y^2$, respectively, onto $H$ have inverse signs.*

**Proof** This is an immediate consequence of Lemma 2.2.2. ■

Consider orthogonal projections that are performed successively onto two given parallel hyperplanes starting from some initial point. In this case the final point is the same point that we would reach by orthogonally projecting the initial point directly onto the second hyperplane. The same thing happens with generalized projections, as the next lemma shows.

**Lemma 2.2.3** *Let $H_1$ and $H_2$ be two parallel hyperplanes in $\mathbb{R}^n$ with representations $H_1 = \{x \mid \langle a, x\rangle = b_1\}$, $H_2 = \{x \mid \langle a, x\rangle = b_2\}$, and let $y \in S$. Then, for any $f \in \mathcal{B}(S)$ which is zone consistent with respect to both hyperplanes,*

$$P_{H_2}(P_{H_1}(y)) = P_{H_2}(y), \tag{2.30}$$

*and the associated projection parameters obey the equation*

$$\lambda_1 + \hat{\lambda}_2 = \lambda_2, \tag{2.31}$$

*where $\lambda_i$, $i = 1, 2$, is the parameter associated with projecting $y$ onto $H_i$, and $\hat{\lambda}_2$ is the parameter associated with projecting $P_{H_1}(y)$ onto $H_2$.*

**Proof** From Lemma 2.2.1 we have

$$\nabla f(x^1) = \nabla f(y) + \lambda_1 a, \quad \langle a, x^1\rangle = b_1, \tag{2.32}$$

$$\nabla f(x^2) = \nabla f(y) + \lambda_2 a, \quad \langle a, x^2\rangle = b_2, \tag{2.33}$$

$$\nabla f(\hat{x}^2) = \nabla f(x^1) + \hat{\lambda}_2 a, \quad \langle a, \hat{x}^2\rangle = b_2, \tag{2.34}$$

describing the projections of $y$ onto $H_1$, of $y$ onto $H_2$, and of $x^1 = P_{H_1}(y)$ onto $H_2$, respectively. Now, substitute (2.32) into (2.34) to obtain

$$\nabla f(\hat{x}^2) = \nabla f(y) + (\lambda_1 + \hat{\lambda}_2)a, \quad \langle a, \hat{x}^2\rangle = b_2. \tag{2.35}$$

From the uniqueness of the projection of $y$ onto $H_2$ it follows that $x^2 = \hat{x}^2$ and that $\lambda_1 + \hat{\lambda}_2 = \lambda_2$ . ■

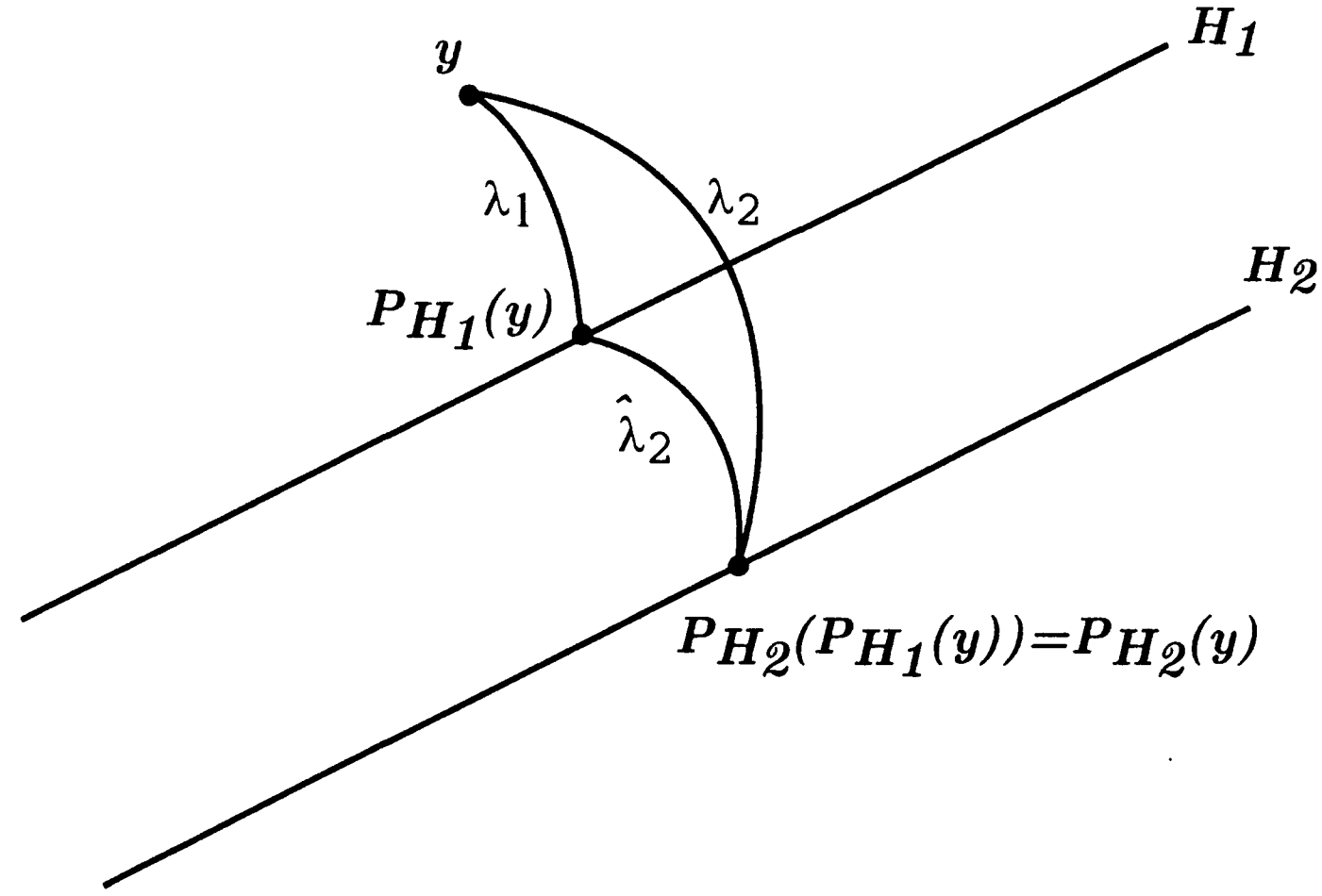

**Figure 2.2** Geometric interpretation of Lemma 2.2.3.

We cannot really describe geometrically the path of generalized projections onto parallel hyperplanes for all possible Bregman functions. But if the arcs in Figure 2.2 are interpreted as trajectories of projected points then this figure illustrates Lemma 2.2.3 fairly well.

**Lemma 2.2.4** *If the conditions of Lemma 2.2.3 hold, with hyperplanes* $H_1 = \{x \mid \langle a,x\rangle = \gamma\}$ *and* $H_2 = \{x \mid \langle a,x\rangle = \delta\}$, *then* $\Gamma \leq \Delta$, *if and only if* $\gamma \leq \delta$, *where* $\Gamma$ *and* $\Delta$ *are the parameters associated with the projections of* $y$ *onto* $H_1$ *and* $H_2$, *respectively.*

**Proof** Observe that

$$\nabla f(z^1) = \nabla f(y) + \Gamma a, \quad \langle a, z^1\rangle = \gamma, \tag{2.36}$$

$$\nabla f(z^2) = \nabla f(y) + \Delta a, \quad \langle a, z^2\rangle = \delta. \tag{2.37}$$

Subtracting (while assuming $\gamma \neq \delta$) and then taking the inner product with $(z^2 - z^1)$ we obtain

$$\langle \nabla f(z^2) - \nabla f(z^1), z^2 - z^1\rangle = (\Delta - \Gamma)(\delta - \gamma). \tag{2.38}$$

The left-hand side equals $D_f(z^1, z^2) + D_f(z^2, z^1)$, which is nonnegative by Lemma 2.1.1. Therefore, $\delta - \gamma \geq 0$, implies $\Delta - \Gamma \geq 0$, and vice versa. ■

## 2.3 Bregman Functions on the Whole Space

The six conditions of Definition 2.1.1 may be difficult to verify in some specific instances. Also, not all conditions are directly formulated for the function $f$; some are related to $f$ through the level sets or through properties of the $D_f$-function. Therefore it is useful to provide conditions that are easier to verify and would guarantee that $f \in \mathcal{B}(S)$. In this subsection we present such conditions for functions defined on the whole space $\mathbb{R}^n$. These are sufficient but not necessary conditions.

**Theorem 2.3.1** *Let $f : \mathbb{R}^n \to \mathbb{R}$ be a function that satisfies the following properties:*

*(I) $f$ is twice continuously differentiable and its Hessian matrix $\nabla^2 f(x)$ is positive definite, for all $x \in \mathbb{R}^n$.*

*(II)* $\displaystyle \lim_{\|x\| \to \infty} \frac{f(x)}{\| x \|} = +\infty.$

*Then $f \in \mathcal{B}(\mathbb{R}^n)$.*

It is worthwhile to note that under the differentiability assumptions in $(I)$, strongly convex functions (i.e., functions that satisfy the condition $\langle \nabla f(x) - \nabla f(y), x - y \rangle \geq \alpha \| x - y \|^2$, for all $x, y \in \mathbb{R}^n$, for some positive constant $\alpha$) are exactly the functions with positive definite Hessians. Also, a continuous function $g(x)$ that is defined on all of $\mathbb{R}^n$ is called *coercive* if $\lim_{\|x\| \to \infty} g(x) = +\infty$; thus property $(II)$ above and property $(IV)$ below are coerciveness conditions.

The proof of Theorem 2.3.1 requires four lemmas. Before studying them, observe that in Definition 2.1.1 conditions $(v)$ and $(vi)$ hold trivially in the interior of the zone, as a consequence of condition $(i)$, and need to be checked only on the boundary of $S$. So, when considering functions defined on all $\mathbb{R}^n$, those conditions are trivially satisfied. Also, in this case, condition $(ii)$ of Definition 2.1.1 is a consequence of $(I)$, which also implies conditions $(i)$ and $(iii)$ of Definition 2.1.1. The only difficulty lies in verifying condition $(iv)$.

Consider a function $g : \mathbb{R}^n \setminus \{0\} \to \mathbb{R}$ with the following properties:

$(III)$ $g$ is twice continuously differentiable,

$(IV)$ $\lim_{\|x\| \to \infty} g(x) = +\infty$, and

$(V)$ the function $h : \mathbb{R}^n \to \mathbb{R}$ defined by $h(x) \doteq g(x) \cdot \| x \|$, if $x \neq 0$ and $h(0) \doteq 0$, has a positive definite Hessian $\nabla^2 h(x)$, for all $x \neq 0$.

For $x \neq 0$, define the auxiliary function $\varphi_x : \mathbb{R}_{++} \to \mathbb{R}$ as

$$\varphi_x(\lambda) \doteq \langle x, \nabla g(\lambda x) \rangle, \tag{2.39}$$

where $\mathbb{R}_{++}$ is the set of positive real numbers. Property $(III)$ guarantees that $\varphi_x$ is continuously differentiable. The sequence of four lemmas that leads to the proof of Theorem 2.3.1 now follows.

**Lemma 2.3.1** *For any $x \neq 0$, $\varphi_x(1) \leq 0$ implies that $\varphi_x(\lambda) \leq 0$, for all $0 < \lambda \leq 1$.*

**Proof** We argue by negation. Assume that there exists a $\lambda_0$, $0 < \lambda_0 < 1$, such that $\varphi_x(\lambda_0) > 0$. Let $\lambda_1 \doteq \inf\{\lambda \mid \lambda > \lambda_0, \ \varphi_x(\lambda) = 0\}$, so that

$$\varphi_x(\lambda) > 0, \text{ for } \lambda_0 \leq \lambda < \lambda_1. \tag{2.40}$$

Since $\varphi_x$ is continuous, $\lambda_1$ is well-defined and $\varphi_x(\lambda_1) = 0$. Also, the derivative of $\varphi_x(\lambda)$ (with respect to $\lambda$) at the point $\lambda_1$ satisfies $\varphi'_x(\lambda_1) \leq 0$, because otherwise, $\varphi_x(\lambda) < \varphi_x(\lambda_1)$, for $\lambda_1 - \epsilon < \lambda < \lambda_1$, for some positive $\epsilon$, in contradiction to (2.40). Therefore, we have

$$0 \geq 2\lambda_1\varphi_x(\lambda_1) + \lambda_1^2\varphi'_x(\lambda_1) = 2\lambda_1\langle x, \nabla g(\lambda_1 x)\rangle + \lambda_1^2\langle x, \nabla^2 g(\lambda_1 x)x\rangle$$

$$= 2\langle y, \nabla g(y)\rangle + \langle y, \nabla^2 g(y)y\rangle, \quad \text{with } y \doteq \lambda_1 x, \tag{2.41}$$

where $\nabla^2 g$ is the Hessian matrix of $g$. Furthermore, property $(V)$ above postulates that the Hessian of $h(y)$ is positive definite, i.e., $0 < \langle y, \nabla^2 h(y)y\rangle$ for all $y \neq 0$. Differentiation yields

$$0 < \langle y, \nabla^2 h(y)y\rangle = 2\| y \|\langle y, \nabla g(y)\rangle + \| y \|\langle y, \nabla^2 g(y)y\rangle - \| y \| g(y),$$

thus, $g(y) < 2\langle y, \nabla g(y)\rangle + \langle y, \nabla^2 g(y)y\rangle$, for all $y \neq 0$. Choosing $y = \lambda_1 x$ for an $x$ with a sufficiently large $\| x \|$ guarantees, in the light of property $(IV)$, that $g(y) > 0$, in contradiction to (2.41). ■

**Lemma 2.3.2** *For every $x \neq 0$, there exists a $\mu > 1$ such that $\varphi_x(\mu) > 0$.*

**Proof** Otherwise, $0 \geq \varphi_x(\mu) = \langle x, \nabla g(\mu x)\rangle$, for all $\mu > 1$. Then, by the mean-value theorem, $g(\mu x) \leq g(x)$ would hold for all $\mu > 1$, in contradiction to property $(IV)$. ■

**Lemma 2.3.3** *The following set is bounded:*

$$V \doteq \{x \in \mathbb{R}^n \setminus \{0\} \mid \langle x, \nabla g(x)\rangle \leq 0\} = \{x \in \mathbb{R}^n \setminus \{0\} \mid \varphi_x(1) \leq 0\}.$$

**Proof** We argue by negation. Suppose that there exists a sequence $\{x^\nu\}_{\nu=0}^\infty$ such that $\{x^\nu\}_{\nu=0}^\infty \subseteq V$, and $\lim_{\nu\to\infty} \| x^\nu \| = +\infty$. Let $x^*$ be an accumulation point of $\{x^\nu / \| x^\nu \|\}$ and take any $\mu \geq 1$. Given any $\epsilon > 0$, take $\nu$ such that $\| (x^\nu / \| x^\nu \|) - x^* \| < \epsilon/\mu$ and $\| x^\nu \| > \mu$. Let $\lambda_\nu \doteq \mu / \| x^\nu \| < 1$. Then,

$$\| \lambda_\nu x^\nu - \mu x^* \| < \epsilon. \tag{2.42}$$

Since $x^\nu \in V$ and $\lambda_\nu > 0$, Lemma 2.3.1 guarantees that $\lambda_\nu x^\nu \in V$. $\mu \geq 1$ implies that $\| \lambda_\nu x^\nu \| \geq 1$. Because of continuity of $\varphi_x$, the set $U \doteq V \cap \{x \in \mathbb{R}^n \mid \| x \| \geq 1\}$ is closed, and, therefore, (2.42) implies that $\mu x^* \in U \subset V$. So, $\varphi_{x^*}(\mu) \leq 0$, for all $\mu \geq 1$, in contradiction to Lemma 2.3.2. ■

**Lemma 2.3.4** *If $f$ has properties $(I)$ and $(II)$ of Theorem 2.3.1, then*

(*i*) $\lim_{\|x\|\to\infty}(f(x) - \langle a, x\rangle) = +\infty$, *for all* $a \in \mathbb{R}^n$; *and*
(*ii*) $\lim_{\|x\|\to\infty}(\langle x - a, \nabla f(x)\rangle - f(x)) = +\infty$, *for all* $a \in \mathbb{R}^n$.

**Proof** (*i*) Take any $a \in \mathbb{R}^n$; then from property $(II)$, it follows that there exists a finite $M$ such that $f(x)/\|x\| \geq 2\|a\|$, for all $x$ with $\|x\| > M$. This means that $\|a\|\|x\| \leq f(x) - \|a\|\|x\| \leq f(x) - \langle a, x\rangle$, and the left-hand side tends to infinity as $\|x\| \to \infty$, if $a \neq 0$. If $a = 0$, (*i*) follows directly from property $(II)$. (*ii*) For an arbitrary $\rho > 0$, consider the function $\overline{f}(x) \doteq f(x + a) + \rho$. Clearly, $\overline{f}$ satisfies properties $(I)$ and $(II)$. So, $g(x) \doteq \overline{f}(x)/\|x\|$ satisfies properties $(III)$, $(IV)$, $(V)$, and, from Lemma 2.3.3, the set $V$ is bounded. Take $M$ such that, for all $y$ with $\|y\| > M$, $\langle y, \nabla g(y)\rangle \geq 0$. Then,

$$0 \leq \langle y, \nabla g(y)\rangle = (1/\|y\|)(\langle y, \nabla\overline{f}(y)\rangle - \overline{f}(y)),$$

which implies that $0 \leq \langle y, \nabla\overline{f}(y)\rangle - \overline{f}(y)$, for all $y$ with $\|y\| > M$.

Letting $x \doteq y + a$, we obtain that, for all $x$ with $\|x\| > M + \|a\|$,

$$0 \leq \langle x - a, \nabla f(x)\rangle - f(x) - \rho.$$

This implies that $\rho \leq \langle x - a, \nabla f(x)\rangle - f(x)$, and since $\rho$ is arbitrary, (*ii*) holds. ■

**Proof (of Theorem 2.3.1)** As noted before, only condition (*iv*) in Definition 2.1.1 remains to be checked. We have

$$L_1^f(y, \alpha) = \{x \mid f(x) - \langle\nabla f(y), x\rangle \leq \alpha + f(y) - \langle\nabla f(y), y\rangle\}.$$

Applying Lemma 2.3.4(*i*) with $a \doteq \nabla f(y)$, yields

$$\lim_{\|x\|\to\infty}(f(x) - \langle\nabla f(y), x\rangle) = +\infty,$$

thus $L_1^f(y, \alpha)$ is bounded, for all $y \in \mathbb{R}^n$. For

$$L_2^f(x, \alpha) = \{y \mid \langle y - x, \nabla f(y)\rangle - f(y) \leq \alpha - f(x)\},$$

apply Lemma 2.3.4(*ii*) with $a \doteq x$, and conclude that

$$\lim_{\|y\|\to\infty}(\langle y - x, \nabla f(y)\rangle - f(y)) = +\infty.$$

So, $L_2^f(y, \alpha)$ is bounded, for all $x \in \mathbb{R}^n$. ■

Property $(II)$ is not a necessary condition for a function $f$ to belong to $\mathcal{B}$, even for twice continuously differentiable functions defined in all $\mathbb{R}^n$, as the following example shows. Let

$$f(x) \doteq \begin{cases} \frac{1}{2}(x^2 - 4x + 3), & \text{if } x \leq 1, \\ -\log x, & \text{if } x \geq 1. \end{cases}$$

It is straightforward to verify that $f$ is a twice continuously differentiable Bregman function but $\lim_{x \to +\infty}(f(x)/x) = 0$.

## 2.4 Characterization of Generalized Projections

We give now a characterization of the generalized projections introduced in Section 2.2. The characterization depends on the following basic result, which we will use repeatedly in the next chapters.

**Theorem 2.4.1** *Let $f \in \mathcal{B}(S)$ and let $\Omega \subseteq \mathbb{R}^n$ be a closed convex set such that $\Omega \cap \overline{S} \neq \emptyset$. Assume that $y \in S$, implies $P_\Omega(y) \in S$. Let $z \in \Omega \cap \overline{S}$, then, for any $y \in S$, the following inequality holds*

$$D_f(P_\Omega(y), y) \leq D_f(z, y) - D_f(z, P_\Omega(y)). \tag{2.43}$$

**Proof** Define the function

$$G(u) \doteq D_f(u, y) - D_f(u, P_\Omega(y)). \tag{2.44}$$

Expanding $G(u)$ according to the definition of $D_f$ shows that it is an affine function of the form $\langle u, a \rangle + b$, where the vector $a$ and the real $b$ are independent of $u$, so that $G(u)$ is convex. For any $\lambda$, $0 \leq \lambda \leq 1$, we denote $u_\lambda \doteq \lambda z + (1 - \lambda)P_\Omega(y)$ and obtain, due to the convexity of G($u$), that

$$\begin{aligned} &D_f(u_\lambda, y) - D_f(u_\lambda, P_\Omega(y)) \\ &\quad \leq \lambda(D_f(z, y) - D_f(z, P_\Omega(y))) + (1 - \lambda)D_f(P_\Omega(y), y). \end{aligned} \tag{2.45}$$

This leads, for $\lambda > 0$, to

$$\begin{aligned} &D_f(z, y) - D_f(z, P_\Omega(y)) - D_f(P_\Omega(y), y) \\ &\quad \geq (1/\lambda)(D_f(u_\lambda, y) - D_f(P_\Omega(y), y)) - (1/\lambda)D_f(u_\lambda, P_\Omega(y)). \end{aligned} \tag{2.46}$$

The first term on the right-hand-side of (2.46) is nonnegative because of (2.4). The second term tends to zero as $\lambda \to 0$. To see this we note that a straightforward calculation, using (2.1) and Definition 2.1.1($i$), shows that

$$\nabla_x D_f(x, P_\Omega(y)) \mid_{x = P_\Omega(y)} = 0. \tag{2.47}$$

The well-known relationship between the directional derivative $F'(u; v)$, at the point $u$ in the direction $v$, of a differentiable function $F$ and its gradient $\nabla F$ that we need next is

$$F'(u;v) = \langle \nabla F(u), v\rangle, \quad \text{for all } v \in \mathbb{R}^n. \tag{2.48}$$

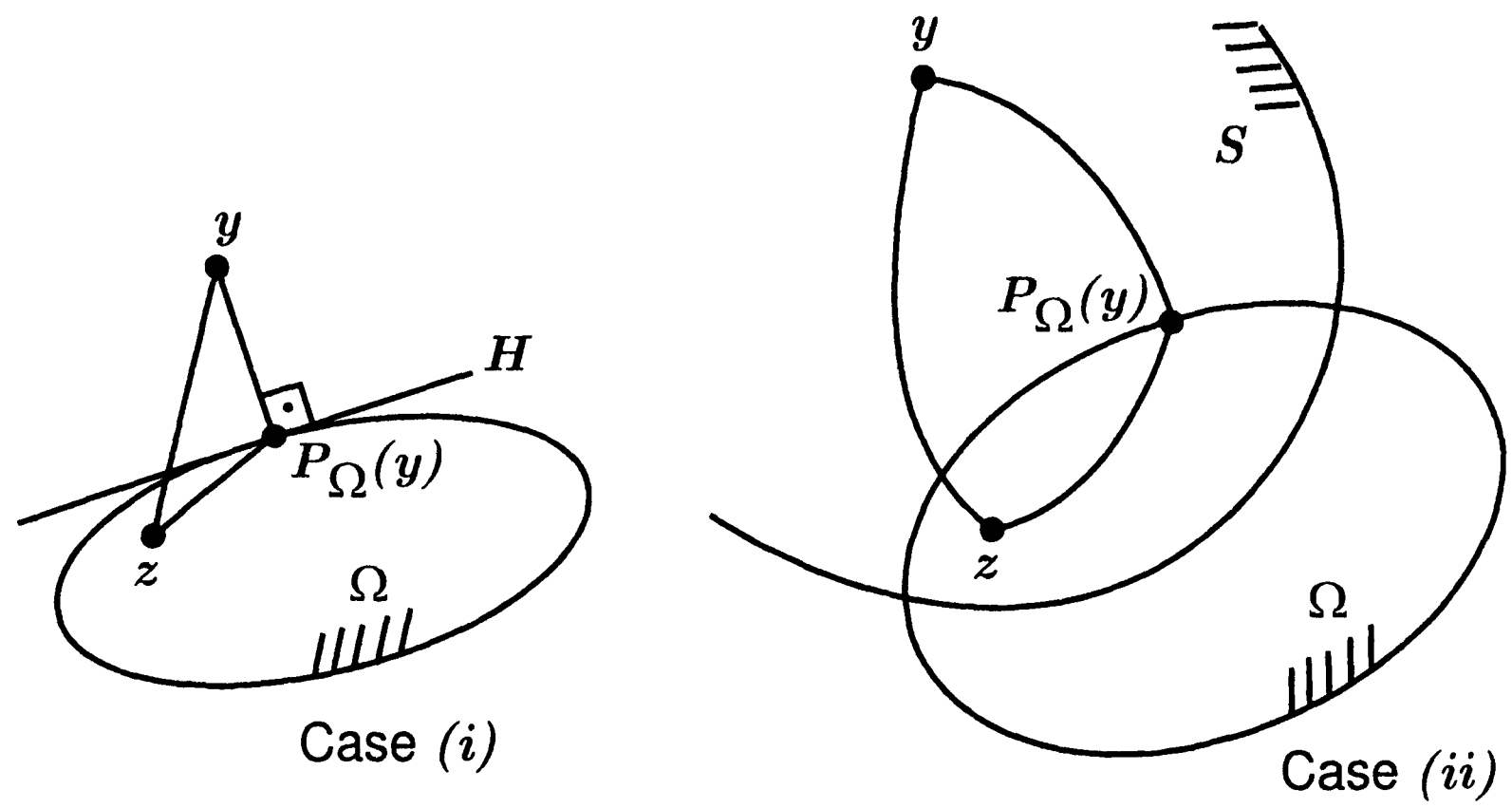

**Figure 2.3** Geometric description for Theorem 2.4.1: (*i*) the case of orthogonal projections; (*ii*) the case of generalized projections.

Using this relationship together with (2.47) we obtain

$$\begin{aligned}
&\lim_{\lambda\to 0^+} \frac{D_f(u_\lambda, P_\Omega(y))}{\lambda} \\
&= \lim_{\lambda\to 0^+} \frac{D_f(P_\Omega(y) + \lambda(z - P_\Omega(y)), P_\Omega(y)) - D_f(P_\Omega(y), P_\Omega(y))}{\lambda} \\
&= \langle \nabla_x D_f(x, P_\Omega(y))|_{x=P_\Omega(y)}, z - P_\Omega(y)\rangle = 0.
\end{aligned} \tag{2.49}$$

■

Figure 2.3 illustrates Theorem 2.4.1. The classic case with orthogonal projections appears in Figure 2.3(*i*) and it means that a hyperplane $H$ through the point $P_\Omega(y)$, and which is perpendicular to $y - P_\Omega(y)$, supports the convex set $\Omega$ at the point $P_\Omega(y)$. That is, the set $\Omega$ lies entirely on "one side" of $H$. Figure 2.3(*ii*) describes the situation for generalized projections, but which are not necessarily orthogonal.

From this result the following characterization of generalized projections is obtained.

**Theorem 2.4.2** *Under the assumptions of Theorem 2.4.1, for any $x \in S$, $x^* \in \Omega \cap \overline{S}$ is $P_\Omega(x)$ if and only if*

$$\langle u - x^*, \nabla f(x) - \nabla f(x^*)\rangle \le 0, \quad \textit{for all}\, u \in \Omega \cap \overline{S}. \tag{2.50}$$

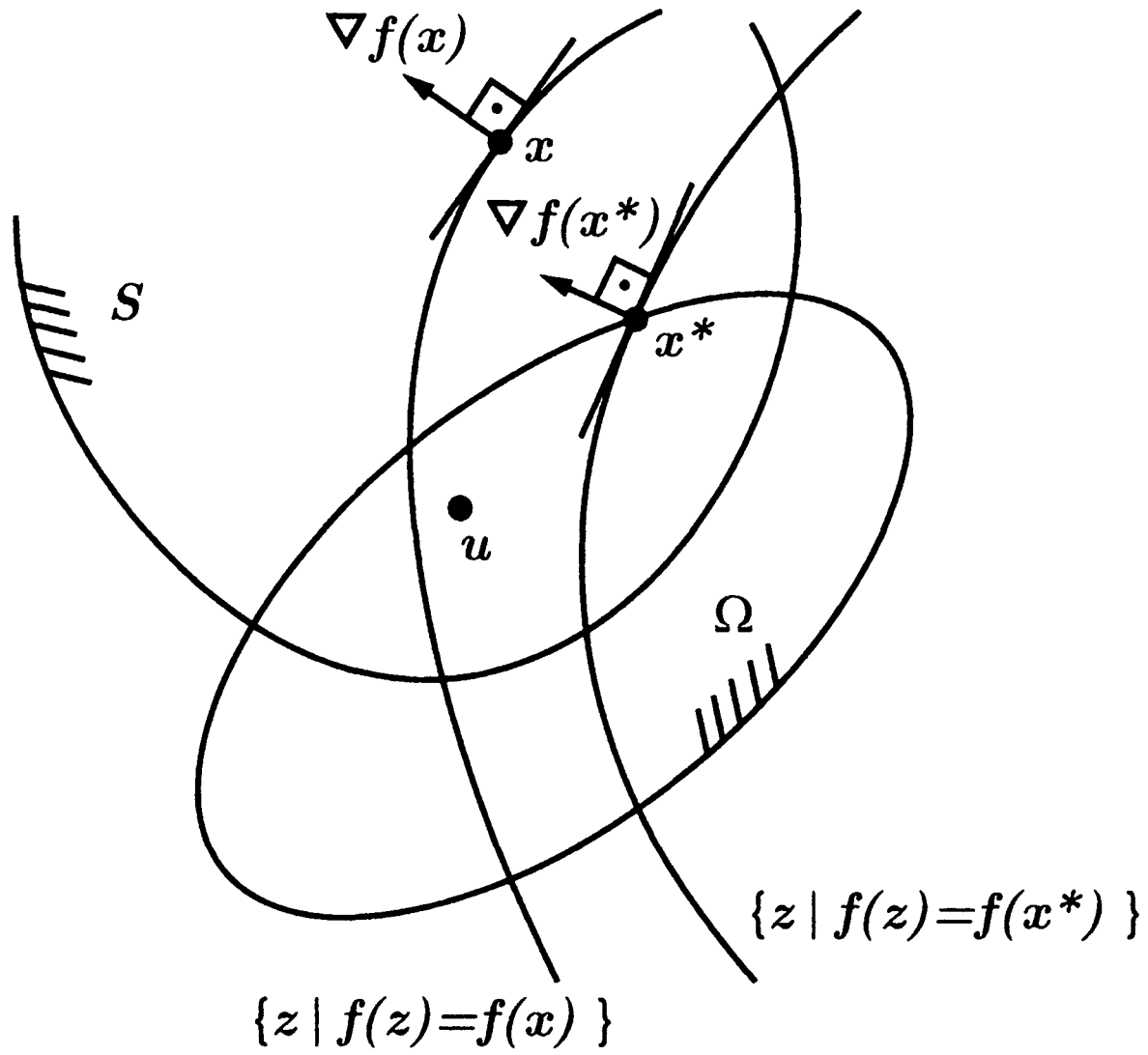

**Figure 2.4** Kolmogorov's Criterion: Geometric description of Theorem 2.4.2 for a set $\Omega$ and $f \in \mathcal{B}(S)$.

**Proof** From the definition of $D_f$, we verify that

$$D_f(v,x) = D_f(u,x) - D_f(u,v) + \langle \nabla f(x) - \nabla f(v), u - v \rangle. \tag{2.51}$$

If $v = x^* = P_\Omega(x)$, then (2.50) follows from (2.43) and (2.51). In the other direction, if (2.50) holds, then (2.51) reduces, because of the nonnegativity of $D_f$, to the inequality

$$D_f(x^*,x) \leq D_f(u,x), \text{ for all } u \in \Omega \cap \overline{S}, \tag{2.52}$$

which means that $x^* = P_\Omega(x)$. ■

For the special case $f(x) = \frac{1}{2} \parallel x \parallel^2$, with $S = \mathbb{R}^n$, Theorem 2.4.2 reduces to a well-known result on orthogonal projections, sometimes referred to as *Kolmogorov's criterion*. In this case it says that

$$\langle u - P_\Omega(x), x - P_\Omega(x) \rangle \leq 0, \text{ for all } u \in \Omega, \tag{2.53}$$

which means that the angle between $(x - P_\Omega(x))$ and $(u - P_\Omega(x))$ is obtuse, for any $u \in \Omega$. In the notation of Figure 2.3($i$) we have that $(y - P_\Omega(y))$ makes an obtuse angle with $(z - P_\Omega(y))$ for any $z \in \Omega$. Figure 2.4 illustrates Theorem 2.4.2 for a general Bregman function $f \in \mathcal{B}(S)$.

## 2.5 Csiszár $\varphi$-divergences

Generalized distances exist that are defined differently from the function $D_f(x, y)$ in (2.1). One important family of such distances is the family of $\varphi$-divergences.

**Definition 2.5.1 ($\varphi$-divergence)** *Let $\varphi : \mathbb{R}_{++} \to \mathbb{R}$ be twice continuously differentiable and strictly convex on $\mathbb{R}_{++}$ with $\varphi(1) = \varphi'(1) = 0$ and $\varphi''(1) > 0$, and $\lim_{t \to 0^+} \varphi'(t) = -\infty$. For such a function $\varphi$, the function $d_\varphi : \mathbb{R}^n_{++} \times \mathbb{R}^n_{++} \to \mathbb{R}$ defined by*

$$d_\varphi(x, y) \doteq \sum_{j=1}^{n} y_j \varphi\left(\frac{x_j}{y_j}\right), \tag{2.54}$$

*is called a $\varphi$-divergence.*

Some basic properties of $\varphi$-divergences are given in the next proposition.

**Proposition 2.5.1** *Let $d_\varphi(x, y)$ be a $\varphi$-divergence. Then,*

*(i) $d_\varphi(x, y) \geq 0$, for all $x, y \in \mathbb{R}_{++}$.*

*(ii) $d_\varphi(x, y) = 0$ if and only if $x = y$.*

*(iii) The level sets of $d_\varphi(\cdot, y)$ are bounded, for all $y \in \mathbb{R}_{++}$.*

*(iv) The level sets of $d_\varphi(x, \cdot)$ are bounded, for all $x \in \mathbb{R}_{++}$.*

*(v) $d_\varphi(x, y)$ is jointly convex in $(x, y)$ and strictly convex in $x$.*

*(vi) $\lim_{\nu \to \infty} d_\varphi(y^*, y^\nu) = 0$ if and only if $\lim_{\nu \to \infty} y^\nu = y^*$.*

**Proof** The strict convexity of $\varphi$ and its zero global minimal value at $t = 1$ imply its nonnegativity, which together with positivity of $y$ imply $(i)$. The fact that $\varphi$ has zero value only at $t = 1$ leads to $(ii)$. The level set $\{x \in \mathbb{R}_{++} \mid d_\varphi(x, y) \leq \alpha\}$ must be bounded, for any $y \in \mathbb{R}_{++}$ and any real $\alpha$, because if we let any component $x_j$ go to $+\infty$ then so would $d_\varphi(x, y)$, violating its bound $\alpha$. Thus, $(iii)$ is true, and a similar argument proves $(iv)$. To prove $(v)$ one calculates the following derivatives:

$$\frac{\partial}{\partial x_j} d_\varphi(x, y) = \varphi'\left(\frac{x_j}{y_j}\right), \quad \frac{\partial^2}{\partial x_j^2} d_\varphi(x, y) = \frac{1}{y_j} \varphi''\left(\frac{x_j}{y_j}\right), \quad \frac{\partial^2}{\partial x_i \partial x_j} d_\varphi(x, y) = 0,$$

$$\frac{\partial}{\partial y_j} d_\varphi(x, y) = \varphi\left(\frac{x_j}{y_j}\right) - \left(\frac{x_j}{y_j}\right) \varphi'\left(\frac{x_j}{y_j}\right), \quad \frac{\partial^2}{\partial y_j^2} d_\varphi(x, y) = \frac{x_j^2}{y_j^3} \varphi''\left(\frac{x_j}{y_j}\right),$$

$$\frac{\partial^2}{\partial y_i \partial y_j} d_\varphi(x, y) = 0.$$

The Hessian matrix of $d_\varphi(x, y)$ with respect to $x$ is diagonal with positive diagonal elements, thus positive definite, which proves strict convexity of $d_\varphi(x,y)$ in $x$. By looking at the Hessian of $d_\varphi(x,y)$ with respect to $(x,y)$, a similar argument guarantees joint convexity of $d_\varphi(x,y)$ in $(x,y)$. Continuity and $\varphi(1)=0$ imply $(vi)$. ■

**Example 2.5.1** $\varphi_1(t) \doteq t\log t - t + 1$ *satisfies the properties listed in Definition 2.5.1, and in this case*

$$d_{\varphi_1}(x,y) = \sum_{j=1}^{n}(x_j \log\left(\frac{x_j}{y_j}\right) + y_j - x_j), \tag{2.55}$$

*is the Kullback-Leibler cross entropy function (2.9), which can be extended to* $\mathbb{R}^n_+ \times \mathbb{R}^n_{++}$.

This is the only known case in which a member of the family of $\varphi$-divergences coincides with a member of the family of Bregman distances $D_f(x,y)$, up to additive linear terms and multiplicative constants.

**Example 2.5.2** $\varphi_2(t) \doteq t - \log t - 1$ *also satisfies the required properties and yields* $d_{\varphi_2}(x,y) = d_{\varphi_1}(y,x)$. *This is the* $\varphi_1$*-divergence of the previous example with exchanged variables* $x$ *and* $y$.

**Example 2.5.3** $\varphi_3(t) \doteq (\sqrt{t}-1)^2$ *gives rise to* $d_{\varphi_3}(x,y) = (\sqrt{x_j} - \sqrt{y_j})^2$.

The squared Euclidean distance $\frac{1}{2} \parallel x - y \parallel^2$, which is a generalized distance for $f(x) = \frac{1}{2} \parallel x \parallel^2$ (see Example 2.1.1) cannot be derived from a $\varphi$-divergence $d_\varphi(x,y)$ for any $\varphi$. However, a certain relationship exists between $D_f(x,y)$ and $d_\varphi(x,y)$, which is formulated next.

**Proposition 2.5.2** *Let* $f(x) = \sum_{j=1}^{n}\varphi(x_j)$ *where* $\varphi(t)$ *has the properties listed in Definition 2.5.1 on* $\mathbb{R}$ *(not just on* $\mathbb{R}_{++}$*) and define the two-dimensional function*

$$B_\varphi(s,t) \doteq \varphi(s) - \varphi(t) - (s-t)\varphi'(t). \tag{2.56}$$

*Then*

$$D_f(x,y) = \sum_{j=1}^{n} B_\varphi(x_j,y_j), \tag{2.57}$$

*and*

$$d_\varphi(x,y) = \sum_{j=1}^{n} y_j B_\varphi\left(\frac{x_j}{y_j},1\right), \tag{2.58}$$

*for all* $x,y \in \mathbb{R}^n$ *and* $x,y \in \mathbb{R}^n_+$, *respectively.*

**Proof** Follows immediately from the definitions. ■

A further logical link connecting generalized distances and $\varphi-$divergences is the fact that it is possible to obtain these two families by an axiomatic approach to statistical inference. Reference to the literature on this axiomatic approach is given in Section 2.6 and we do not give it further consideration here. In the rest of the book we concentrate on Bregman distances. Some results on generalized distances have their $\varphi$-divergence analogs, while for others the question of extending them to $\varphi$-divergences is still open.

## 2.6 Notes and References

For many years generalized distances of various kinds have been playing an important role in minimum distance estimation in statistical analysis. See, e.g., Titterington, Smith, and Makov (1985, Chapter 4) or Liese and Vajda (1987), or the review paper by Teboulle (1992a). But it was Bregman (1967b) who introduced the family of the $D_f(x, y)$ generalized distances as a tool for developing optimization techniques. His work was continued by Censor and Lent (1981) and then by De Pierro and Iusem (1986) and others. It is interesting to note that Hestenes (1975, Chapter 5) also sensed the importance of $D_f(x, y)$, which he called $E$-function.

**2.1** The definition of Bregman functions (Definition 2.1.1) is based on Bregman (1967b). It first appeared in this form in the paper by Censor and Lent (1981) where Lemma 2.1.2 and the other results presented here appear. Examples of additional Bregman functions and distances can be found in the papers of Teboulle (1992b), Iusem, Svaiter and Teboulle (1994), Eckstein (1992), and Iusem (1995).
Additional conditions on Bregman functions appear in Iusem (1995). Extensions of the class of Bregman functions can be found in Bauschke and Borwein (1995) who gain new insights by studying what they call Bregman/Legendre functions. Kiwiel (1994) extends the class of Bregman functions to generalized Bregman functions (which he calls $B$-functions), and uses this more general class in his subsequent study (1995a).

Since only conditions $(i)$–$(ii)$ of Definition 2.1.1 are used in Lemma 2.1.1, it can be seen that the requirement on the boundedness of the sets $L_2^f(x, \alpha)$, but not of $L_1^f(y, \alpha)$ is redundant. See Kiwiel (1994).

**2.2** This material was developed by Censor and Lent (1981). The property embodied in Lemma 2.2.3 is termed *parallel transitivity* in Csiszár (1991). Additional properties of $D_f$ functions are related to other axioms presented in Csiszár's paper.

**2.3** This treatment of Bregman functions with $S = \mathbb{R}^n$ was derived by De Pierro and Iusem (1986). However, as noted by Professor Iusem, (in private communication) they erroneously assumed only strict convexity in properties $(I)$ and $(V)$ instead of positive definitness of the Hessian matrix. For further information on such functions and on strongly convex functions see, e.g., Polyak (1987), Hiriart-Urruty and Lemaréchal (1993) or Bertsekas and Tsitsiklis (1989). Part $(i)$ of Lemma 2.3.4 can also be found in Ortega and Rheinboldt (1970, E4.3-9). For a definition of coercive function see, e.g., Peressini, Sullivan, and Uhl (1988).

**2.4** Theorem 2.4.1 is a specialization of Bregman (1967b, Lemma 1) for the case of $D_f$ functions. Kolmogorov's criterion, i.e., Theorem 2.4.2 for $f(x) = \frac{1}{2} \parallel x \parallel^2$ with $S = \mathbb{R}^n$, can be found, e.g., in Kinderlehrer and Stampacchia (1980) and in Zarantonello (1971). Theorem 2.4.2 appears in Censor and Reich (1996) and has been proved in a different way by applying optimality conditions to the generalized distance minimization involved in the projection; see Teboulle (1992b).

**2.5** The literature on Csiszár's $\varphi$-divergences is very large. The concept was introduced by Csiszár (1963, 1967), and independently by Ali and Silvey (1966). Further examples of $d_\varphi(x,y)$ for various $\varphi$ can be found, e.g., in Ben-Tal and Teboulle (1987), in Kapur and Kesavan (1992), and in Iusem (1994). Proposition 2.5.2 appears in Teboulle (1992b). An axiomatic approach to statistical inference which leads to the families $D_f(x,y)$ and $d_\varphi(x,y)$ is given in Csiszár (1991). See also the related work of Jones and Trutzer (1990), Jones and Byrne (1990), Shore and Johnson (1981), Snickars and Weibull (1977), and Csiszár (1977).

CHAPTER 3

# Proximal Minimization with $D$-Functions

In the next chapters we will see how generalized distances, and projections derived from them, give rise to useful iterative algorithms. These include the multiprojection algorithm, which is particularly suitable for the feasibility problem (see Chapter 5) and the various algorithms for optimization of Bregman functions (see Chapter 6). The use of generalized distances does not stop there. In this chapter we study the use of generalized distances for the *regularization* of optimization problems, and as a result we obtain a class of optimization algorithms called *proximal minimization algorithms*. These algorithms originally used Euclidean distances, and we show here how generalized distances can be incorporated. Later we examine how computational results shed light on the practical advantages of this approach (see Chapter 15).

The proximal minimization algorithm is designed to solve the problem

$$\text{Minimize} \quad F(x) \tag{3.1}$$

$$\text{s.t.} \quad x \in X, \tag{3.2}$$

where $F : \mathbb{R}^n \to \mathbb{R}$ is a given convex function and $X \subseteq \mathbb{R}^n$ is a nonempty closed convex subset of $\mathbb{R}^n$. The approach is based on converting (3.1)–(3.2) into a sequence of optimization problems by adding terms to $F(x)$. These added terms measure the distance between the variable vector $x$ and the current iterate $x^\nu$, either in the Euclidean sense, or in a generalized sense, according to some $D_f$-function, where $f \in \mathcal{B}(S)$ is a Bregman function (as defined and studied in Chapter 2).

The proximal minimization algorithm can be viewed within the general approach known as regularization. Suppose that the problem (3.1)–(3.2) is ill-conditioned in some sense. By this we mean simply that it lacks some desirable property, such as strict convexity, for example. In such a situation it is sometimes possible to modify the original objective function $F(x)$ by adding a term $\gamma G(x)$ such that the perturbed function $F(x) + \gamma G(x)$ has the desired property, which $F(x)$ lacks. The quantity $\gamma > 0$ is called the *regularization parameter*. For the regularization method to work we usually look for conditions that will guarantee that as $\gamma \to 0^+$ (i.e., the value of $\gamma$ approaches zero from the right) the solution of the perturbed problems will converge to a solution of the original problem. That is, if $x(\gamma) \doteq \operatorname{argmin}(F(x) + \gamma G(x) \mid x \in X)$, for $\gamma > 0$, then $\lim_{\gamma \to 0^+} x(\gamma) = x^*$, where $x^*$ is a solution of problem (3.1)–(3.2). For example, if one uses a quadratic regularizing function $G(x) \doteq \frac{1}{2} \parallel x - a \parallel^2$, with some fixed $a \in \mathbb{R}^n$, and a sequence $\{\gamma_\nu\}_{\nu=0}^{\infty}$, such that $\lim_{\nu\to\infty} \gamma_\nu = 0$, of positive regularizing parameters, then the iteration

$$x^{\nu+1} = \operatorname*{argmin}_{x \in X} \left( F(x) + \frac{\gamma_\nu}{2} \parallel x - a \parallel^2 \right), \tag{3.3}$$

is obtained. A difficulty associated with such an iteration is that, as $\gamma_\nu$ approaches zero, the effect of the additive quadratic term diminishes to a point where—practically, though not theoretically—the objective function of (3.3) again lacks the desired property whose absence we are trying to remedy in the first place. A possible way around this difficulty is to refrain from forcing the sequence $\{\gamma_\nu\}_{\nu=0}^{\infty}$ to go to zero, and at the same time replace $a$ in the quadratic $G(x)$ by the current iterate $x^\nu$. This results in a scheme

$$x^{\nu+1} = \operatorname*{argmin}_{x \in X} \left( F(x) + \frac{\gamma_\nu}{2} \parallel x - x^\nu \parallel^2 \right), \tag{3.4}$$

where the sequence $\{\gamma_\nu\}_{\nu=0}^{\infty}$ is bounded away from zero. In this chapter, we show that not only does such a scheme converge, but that it is possible to replace the quadratic additive term by any $D_f$-function, i.e., by a term $D_f(x, x^\nu)$ for any Bregman function $f \in \mathcal{B}(S)$. The interest in regularization methods for parallel optimization is motivated by the fact that many parallel algorithms are available for solving optimizations problems with $D_f$-functions (as, for example, the algorithms in Chapter 6). Therefore, using the regularization approach enables us to take a problem of the form (3.1)–(3.2), for which parallel algorithms may not be available, and solve it by solving a sequence of regularized subproblems using a parallel algorithm. The mathematical foundation of this general approach is the subject of this chapter. Its practical viability is later demonstrated through applications to various real-world problems (Chapter 15).

## 3.1 The Proximal Minimization Algorithm

We assume a given or constructed sequence $\{c_\nu\}$ of positive numbers for all $\nu \in N_0$, $N_0 \doteq \{0, 1, 2, 3, \ldots\}$ with

$$\liminf_{\nu \to \infty} c_\nu = c > 0. \tag{3.5}$$

The original quadratic proximal minimization method is described as follows.

**Algorithm 3.1.1 The Quadratic Proximal Minimization Algorithm**

**Step 0:** (Initialization.) $y^0 \in \mathbb{R}^n$ is arbitrary.
**Step 1:** (Iterative step.)

$$x^{\nu+1} = \underset{x \in X}{\operatorname{argmin}} \left( F(x) + \frac{1}{2c_\nu} \| x - x^\nu \|^2 \right). \tag{3.6}$$

The *proximal minimization algorithm with $D_f$-functions,* henceforth abbreviated PMD, is defined as follows: given are a Bregman function $f$ with zone $S$, and a positive sequence $\{c_\nu\}$ for which (3.5) holds.

**Algorithm 3.1.2 The Proximal Minimization Algorithm with $D$-functions (PMD)**

**Step 0:** (Initialization.) $x^0 \in S$ is arbitrary.
**Step 1:** (Iterative step.)

$$x^{\nu+1} = \underset{x \in X \cap \overline{S}}{\operatorname{argmin}} \left( F(x) + \frac{1}{c_\nu} D_f(x, x^\nu) \right). \tag{3.7}$$

In order for this algorithm to be well-defined we make the following assumption.

**Assumption 3.1.1** *The PMD algorithm (Algorithm 3.1.2) generates a sequence $\{x^\nu\}$ such that $x^\nu \in S$, for every $\nu = 0, 1, 2, \ldots$ .*

This assumption is needed, because $D_f$ is defined on $\overline{S} \times S$. It actually tells us that, given $F$ and $X$ of (3.1)–(3.2), we are free to choose only $f$ and $S$ such that Assumption 3.1.1 would hold. If $X \subseteq S$, then Assumption 3.1.1 trivially holds, which is true for the quadratic case $f(x) = \frac{1}{2} \| x \|^2$, where $S = \mathbb{R}^n$. We show later that it holds also for the entropy case.

## 3.2 Convergence Analysis of the PMD Algorithm

In our analysis we first establish the existence and uniqueness of the minimum of $(F(x) + \frac{1}{c} D_f(x, y))$. Denote by $X^*$ the solution set of problem (3.1)–(3.2), i.e.,

$$X^* \doteq \{x^* \in X \mid F(x^*) \leq F(x), \text{ for all } x \in X\}. \tag{3.8}$$

**Proposition 3.2.1** *Let $f \in \mathcal{B}(S)$. For every $y \in S$ and $c > 0$, the minimum of $F(x) + 1/c\, D_f(x, y)$ over $X \cap \overline{S}$ is attained at a unique point, denoted by $x_f(y, c)$, provided that $F(x)$ is bounded below over $X$.*

**Proof** For all $c > 0$ and $y \in S$, the level sets

$$\{x \in X \cap \overline{S} \mid F(x) + \frac{1}{c} D_f(x, y) \leq \alpha\}, \qquad \alpha \in \mathbb{R}, \tag{3.9}$$

are bounded. This is true because otherwise, for some $c > 0$ and $y \in \mathbb{R}^n$, there would exist an unbounded sequence $\{x^\nu\} \subseteq X \cap \overline{S}$ for which

$$D_f(x^\nu, y) \leq c\,(\alpha - L), \tag{3.10}$$

where $L$ is the lower bound for $F(x)$ over $X$. But this would contradict $(iv)$ of Definition 2.1.1; thus, the level sets (3.9) must be bounded. This allows us to equivalently search for the minimum of $F(x) + 1/c\,D_f(x, y)$ over a compact subset of $X \cap \overline{S}$ instead of over $X \cap \overline{S}$. The Weierstrass theorem, which says that a real continuous function on a compact set in $\mathbb{R}^n$ attains its minimum on the set, then implies that the above-mentioned minimum is attained.

Here we use the fact that since $F$ takes only finite values and is convex on all of $\mathbb{R}^n$ then it is necessarily continuous (see, e.g., Rockafellar 1970, Corollary 10.1.1). The strict convexity of $D_f(x, y)$ with respect to $x$, for fixed $y$, stemming from the strict convexity of the function $f$, ensures the uniqueness. ■

**Proposition 3.2.2** *If $D_f(x, y)$ is jointly convex with respect to both $x$ and $y$, i.e., as a function on $\mathbb{R}^{2n}$, then the function $\Phi_c : S \to \mathbb{R}$, defined by*

$$\Phi_c(y) \doteq \min_{x \in X \cap \overline{S}} \left( F(x) + \frac{1}{c} D_f(x, y) \right), \tag{3.11}$$

*is convex over $S$.*

**Proof** Let $y^1, y^2 \in S$ and let $0 \leq \alpha \leq 1$. Denote $x_f^i = x_f(y^i, c)$, for $i = 1, 2$, according to the notation introduced in Proposition 3.2.1. Using the convexity of $F$ and joint convexity of $D_f(x, y)$, we have

$$\begin{aligned}
&\alpha \Phi_c(y^1) + (1 - \alpha)\Phi_c(y^2) \\
&= \alpha(F(x_f^1) + \frac{1}{c} D_f(x_f^1, y^1)) + (1 - \alpha)(F(x_f^2) + \frac{1}{c} D_f(x_f^2, y^2)) \\
&\geq F(\alpha x_f^1 + (1 - \alpha)x_f^2) + \frac{1}{c} D_f(\alpha x_f^1 + (1 - \alpha)x_f^2, \alpha y^1 + (1 - \alpha)y^2) \\
&\geq \min_{x \in X \cap \overline{S}} (F(x) + \frac{1}{c} D_f(x, \alpha y^1 + (1 - \alpha)y^2)) \\
&= \Phi_c(\alpha y^1 + (1 - \alpha)y^2),
\end{aligned}$$

and the result follows. ■

**Proposition 3.2.3** *Let $f \in \mathcal{B}(S)$ be twice continuously differentiable, and let $D_f(x,y)$ be jointly covex with respect to both $x$ and $y$. Then the function $\Phi_c(y)$ is continuously differentiable on $S$ and its gradient is given by*

$$\nabla\Phi_c(y) = [\nabla^2 f(y)]^{\mathrm{T}} \left(\frac{y - x_f(y,c)}{c}\right), \tag{3.12}$$

*where $\nabla^2 f(y)$ stands for the Hessian matrix of $f$ at $y$.*

**Proof** Consider any $y \in S$, $d \in \mathbb{R}^n$, and $\alpha > 0$ such that $y + \alpha d \in S$. Using the directional derivative of $\Phi_c$ at the point $y$ in the direction $d$, $\Phi'_c(y;d)$, we have

$$\begin{aligned} & F(x_f(y,c)) + \frac{1}{c} D_f(x_f(y,c), y + \alpha d) \\ & \geq \Phi_c(y + \alpha d) \geq \Phi_c(y) + \alpha\Phi'_c(y;d) \\ & = F(x_f(y,c)) + \frac{1}{c} D_f(x_f(y,c), y) + \alpha\Phi'_c(y;d), \end{aligned} \tag{3.13}$$

where the second inequality in (3.13) follows from the convexity of $\Phi_c$ (Proposition 3.2.2). Therefore, using (2.1) we get from (3.13)

$$\begin{aligned} & \frac{1}{c} \left(f(y) - f(y + \alpha d) + \langle \nabla f(y) - \nabla f(y + \alpha d), x_f(y,c) - y\rangle\right. \\ & \left. + \langle \nabla f(y + \alpha d), \alpha d \rangle\right) \geq \alpha\Phi'_c(y;d). \end{aligned} \tag{3.14}$$

Since

$$\lim_{\alpha \to 0} \left( (f(y) - f(y + \alpha d))/\alpha + \frac{1}{\alpha} \langle \nabla f(y + \alpha d), \alpha d \rangle \right) = 0, \tag{3.15}$$

$$\lim_{\alpha \to 0} \frac{\nabla f(y) - \nabla f(y + \alpha d)}{\alpha} = -\nabla^2 f(y) d, \tag{3.16}$$

we obtain from (3.14), by dividing by $\alpha$ and letting $\alpha \to 0$,

$$\langle \nabla^2 f(y) d, (y - x_f(y,c))/c \rangle \geq \Phi'_c(y;d), \quad \text{for all } d \in \mathbb{R}^n. \tag{3.17}$$

Replacing $d$ by $-d$ in (3.17), we get

$$-\langle \nabla^2 f(y) d, (y - x_f(y,c))/c \rangle \geq \Phi'_c(y;-d) \geq -\Phi'_c(y;d), \tag{3.18}$$

where the second inequality is a standard relation for directional derivatives of convex functions.

The relations (3.17)–(3.18) imply that

$$\Phi'_c(y;d) = \langle \nabla^2 f(y) d, (y - x_f(y,c))/c \rangle, \quad \text{for all } d \in \mathbb{R}^n, \tag{3.19}$$

or, equivalently, that $\Phi_c$ is differentiable and that its gradient is given by (3.12). Since $\Phi_c$ is convex (see Proposition 3.2.2) it follows that its gradient is continuous. ■

Define $S^* \doteq \{y^* \in S \mid \Phi_c(y^*) \leq \Phi_c(y), \text{ for all } y \in S\}$. The next proposition gives a relationship between $S^*$, which is the set of minimizers of $\Phi_c(y)$, the zone $S$ of the Bregman function $f$, and the solution set $X^*$. For a function $f$ with zone $S = \mathbb{R}^n$, we get from this relationship, as a special case, that $X^* = S^*$.

**Proposition 3.2.4** *Let the Hessian $\nabla^2 f(z)$ be nonsingular, for all $z \in S^*$. Then,*

$$X^* \cap S = S^*. \tag{3.20}$$

**Proof** The function $F(x) + \frac{1}{c} D_f(x,y)$ takes, for $x = y$, the value $F(y)$ because of Lemma 2.1.1. It follows that

$$\Phi_c(y) \leq F(y), \text{ for all } y \in X \cap S. \tag{3.21}$$

If $z^* \in X^* \cap S$ then (3.21) holds and we have, for all $y \in S$,

$$\Phi_c(z^*) \leq F(z^*) \leq F(x_f(y,c)) \leq F(x_f(y,c)) + \frac{1}{c} D_f(x_f(y,c),y) = \Phi_c(y),$$

because $D_f(x,y) \geq 0$. Thus, $z^*$ minimizes $\Phi_c(y)$ over $S$, i.e., $z^* \in S^*$. Conversely, if $z^* \in S^*$ we have, from (3.12),

$$c\nabla\Phi_c(z^*) = \nabla^2 f(z^*)^{\mathrm{T}}(z^* - x_f(z^*,c)) = 0, \tag{3.22}$$

which implies that $z^* = x_f(z^*,c) \in X \cap S$. Using again (3.21), we have

$$F(z^*) = \Phi_c(z^*) \leq \Phi_c(y) \leq F(y), \text{ for all } y \in X \cap S, \tag{3.23}$$

and therefore $z^* \in X \cap S$. ■

We can now present the convergence proof of the PMD algorithm.

**Theorem 3.2.1** *Let $f \in \mathcal{B}(S)$; let Assumption 3.1.1 hold; and assume that $X^* \cap \overline{S} \neq \emptyset$. Then any sequence $\{x^\nu\}$ generated by the PMD algorithm (Algorithm 3.1.2), where $c_\nu > 0$ and $\liminf_{\nu\to\infty} c_\nu = c > 0$, converges to an element of $X^*$.*

**Proof** The proof consists of three steps:
**Step 1**. The sequence $\{x^\nu\}$ is bounded.
**Step 2**. All the accumulation points of $\{x^\nu\}$ belong to $X^*$.
**Step 3**. There is a unique limit point.
We now proceed with the proof of each step.

**Step 1.** We have, for all $x \in X \cap \overline{S}$,

$$F(x^{\nu+1}) + \frac{1}{c_\nu} D_f(x^{\nu+1}, x^\nu) \leq F(x) + \frac{1}{c_\nu} D_f(x, x^\nu), \tag{3.24}$$

from which follows that, for all $x \in X \cap \overline{S}$ with $F(x) \leq F(x^{\nu+1})$, it is true that

$$D_f(x^{\nu+1}, x^\nu) \leq D_f(x, x^\nu). \tag{3.25}$$

Therefore, $x^{\nu+1}$ is the unique generalized projection of $x^\nu$ onto the convex set

$$\Omega \doteq \{x \in X \mid F(x) \leq F(x^{\nu+1})\}. \tag{3.26}$$

Using Theorem 2.4.1 and the fact that $X^* \subseteq \Omega$, we have

$$0 \leq D_f(x^{\nu+1}, x^\nu) \leq D_f(x^*, x^\nu) - D_f(x^*, x^{\nu+1}), \tag{3.27}$$

for every $x^* \in X^* \cap \overline{S}$. Thus, for all $\nu \geq 0$,

$$D_f(x^*, x^{\nu+1}) \leq D_f(x^*, x^\nu), \quad \text{for any } x^* \in X^* \cap \overline{S}. \tag{3.28}$$

A sequence $\{x^\nu\}$ that satisfies this last inequality is called *$D_f$-Fejér-monotone with respect to the set $X^* \cap \overline{S}$* (see also Definition 5.3.1). This implies that $\{x^\nu\}$ is bounded because it means that $x^\nu \in L_2^f(x^*, \alpha)$, for all $\nu \geq 0$, with $\alpha \doteq D_f(x^*, x^0)$, and condition (*iv*) of Definition 2.1.1 applies.

**Step 2.** Let $\{x^\nu\}_{\nu \in T}$, where $T \subseteq N_0$, be a subsequence converging to a point $x^\infty \in X \cap \overline{S}$. Recall that, by Lemma 2.1.1 and (3.28), the sequence $\{D_f(x^*, x^\nu)\}_{\nu \in N_0}$ is nonnegative and nonincreasing. It follows that $\lim_{\nu \to \infty} D_f(x^*, x^\nu)$ exists for any $x^* \in X^* \cap \overline{S}$. In view of (3.27), $\lim_{\nu \to \infty} D_f(x^{\nu+1}, x^\nu) = 0$, thus, also

$$\lim_{\substack{\nu \to \infty \\ \nu \in T}} D_f(x^{\nu+1}, x^\nu) = 0, \tag{3.29}$$

which by condition (*vi*) of Definition 2.1.1 implies that $\{x^{\nu+1}\}_{\nu \in T}$ also converges to $x^\infty$.

Next, observe that (3.24) remains true with $x = x^\nu$. Thus,

$$F(x^{\nu+1}) + \frac{1}{c_\nu} D_f(x^{\nu+1}, x^\nu) \leq F(x^\nu), \qquad \text{for all } \nu \geq 0, \tag{3.30}$$

because $D_f(x^\nu, x^\nu) = 0$. Therefore $F(x^\nu) - F(x^{\nu+1}) \geq (1/c_\nu)\, D_f(x^{\nu+1}, x^\nu) \geq 0$, for all $\nu \geq 0$, i.e., $\{F(x^\nu)\}_{\nu \in N_0}$ is nonincreasing, so that $\{F(x^\nu)\}_{\nu \in T}$ must converge to $F(x^\infty)$.

Let $x^* \in X^* \cap \overline{S}$ and $\alpha \in (0,1)$, and set $x \doteq \alpha x^* + (1-\alpha)x^{\nu+1}$ in (3.24). From the convexity of $F(x)$, we get

$$\begin{aligned} & F(x^{\nu+1}) + \frac{1}{c_\nu} D_f(x^{\nu+1}, x^\nu) \\ & \leq F(\alpha x^* + (1-\alpha)x^{\nu+1}) + \frac{1}{c_\nu} D_f(\alpha x^* + (1-\alpha)x^{\nu+1}, x^\nu) \\ & \leq \alpha F(x^*) + (1-\alpha)F(x^{\nu+1}) + \frac{1}{c_\nu} D_f(\alpha x^* + (1-\alpha)x^{\nu+1}, x^\nu), \end{aligned} \tag{3.31}$$

which, by the definition of the $D_f$ function, can be rewritten as

$$\begin{aligned} & \alpha c_\nu (F(x^{\nu+1}) - F(x^*)) \leq D_f(\alpha x^* + (1-\alpha)x^{\nu+1}, x^\nu) + (1/c_\nu) D_f(x^{\nu+1}, x^\nu) \\ & = f(x^{\nu+1} - \alpha(x^* - x^{\nu+1})) - f(x^{\nu+1}) - \alpha \langle \nabla f(x^\nu), x^* - x^{\nu+1} \rangle. \end{aligned} \tag{3.32}$$

Dividing by $\alpha$, taking the limit as $\alpha \to 0^+$, and denoting by $f'(\cdot\,;\cdot)$ the directional derivative, we get

$$\begin{aligned} & c_\nu (F(x^{\nu+1}) - F(x^*)) \\ & \leq f'(x^{\nu+1}; x^* - x^{\nu+1}) - \langle \nabla f(x^\nu), x^* - x^{\nu+1} \rangle \\ & = \langle \nabla f(x^{\nu+1}) - \nabla f(x^\nu), x^* - x^{\nu+1} \rangle \\ & = D_f(x^*, x^\nu) - D_f(x^*, x^{\nu+1}) - D_f(x^{\nu+1}, x^\nu). \end{aligned} \tag{3.33}$$

The optimality of $x^*$ and the fact that, for some $\epsilon > 0$ we have $c_\nu \geq \epsilon > 0$, for all $\nu \geq 0$, guarantee the nonnegativity of the left-hand side of (3.33) for all $\nu \geq 0$.

From the existence of $\lim_{\nu \to \infty} D_f(x^*, x^\nu)$, for any $x^* \in X^* \cap \overline{S}$, we obtain

$$\lim_{\substack{\nu \to \infty \\ \nu \in T}} (D_f(x^*, x^\nu) - D_f(x^*, x^{\nu+1})) = 0, \tag{3.34}$$

and the third term in the right-hand side of (3.33) also tends to zero by (3.29). Thus, $0 = \lim_{\substack{\nu \to \infty \\ \nu \in T}} (F(x^{\nu+1}) - F(x^*)) = F(x^\infty) - F(x^*)$, and so $x^\infty \in X^*$.

**Step 3.** Let $\{x^\nu\}_{\nu \in T_1}$ and $\{x^\nu\}_{\nu \in T_2}$ be two convergent subsequences of $\{x^\nu\}_{\nu=0}^\infty$, generated by the PMD algorithm, i.e.,

$$\lim_{\substack{\nu \to \infty \\ \nu \in T_1}} x^\nu = x^* \in X^* \cap \overline{S}, \tag{3.35}$$

$$\lim_{\substack{\nu \to \infty \\ \nu \in T_2}} x^\nu = x^{**} \in X^* \cap \overline{S}, \tag{3.36}$$

for some $T_1 \subseteq N_0$ and $T_2 \subseteq N_0$. Then, by Definition 2.1.1($v$) we have

$$\lim_{\substack{\nu\to\infty\\ \nu\in T_1}} D_f(x^*,x^\nu)=0,\quad \lim_{\substack{\nu\to\infty\\ \nu\in T_2}} D_f(x^{**},x^\nu)=0. \tag{3.37}$$

$D_f$-Fejér monotonicity then implies that $D_f(x^*,x^{\nu+1})\le D_f(x^*,x^\nu)$ for all $\nu\ge 0$, thus $\lim_{\nu\to\infty} D_f(x^*,x^\nu)$ exists, and because of (3.37) is equal to zero. This implies that $\lim_{\substack{\nu\to\infty\\ \nu\in T_2}} D_f(x^*,x^\nu)=0$, and Definition 2.1.1$(vi)$ yields $x^*=x^{**}$. ■

## 3.3 Special Cases: Quadratic and Entropic PMD

Choosing the Bregman function $f(x)=\frac{1}{2}\,\|x\|^2$ with $\Lambda=S=\overline{S}=\mathbb{R}^n$ gives $D_f(x,y)=\frac{1}{2}\,\|x-y\|^2$ and immediately returns the PMD algorithm to its original form with quadratic additive terms.

For the negative Shannon's entropy $f(x)=-\operatorname{ent}x$ (Example 2.1.2) we have

$$D_f(x,y)=\sum_{j=1}^n x_j\left(\log\left(\frac{x_j}{y_j}\right)-1\right)+\sum_{j=1}^n y_j, \tag{3.38}$$

and the iterative step of the entropy-type PMD algorithm can be obtained from (3.7). When problem (3.1)–(3.2) is linear programming, i.e., when $F(x)=\langle b,x\rangle$ for some given $b=(b_j)\in\mathbb{R}^n$, and $x\in X$ are linear constraints, then (3.7) becomes

$$x^{\nu+1}=\operatorname*{argmin}_{x\in X\cap\mathbb{R}^n_+}\left\{\sum_{j=1}^n x_jb_j+\frac{1}{c_\nu}\sum_{j=1}^n x_j\left(\log\left(\frac{x_j}{x_j^\nu}\right)-1\right)+\frac{1}{c_\nu}\sum_{j=1}^n x_j^\nu\right\}.$$

This is essentially a linearly constrained entropy optimization problem, obtained by subsuming all linear terms into the entropy functional. For such problems, the iterative algorithms of Chapter 6, which lend themselves efficiently to parallel computation, are applicable. The advantages, in practice, of an entropic or other PMD algorithm over the original quadratic proximal minimization algorithm depend largely on the specific form of the original problem (3.1)–(3.2), and on the availability of efficient special-purpose algorithms for performing the step (3.7). These practical questions will be further examined later (see Chapter 15).

## 3.4 Notes and References

The proximal point algorithm is essentially the method of successive approximations applied to the problem of finding fixed points of some suitable nonexpansive mappings in Hilbert spaces; see, e.g., the recent survey paper of Lemaire (1989). Martinet (1970, 1972) seems to have been the first who applied it to convex optimization and Rockafellar (1976a, 1976b) further developed the algorithm and extensively analyzed it. The regularization approach toward the proximal minimization algorithm appears in various publications; see, e.g., Polyak (1987, Chapter 6) or the survey article by Iusem (1995).

**3.1–3.2** The historical development of the proximal point algorithm for solving the problem $0 \in T(z)$ for an arbitrary maximal monotone operator $T$ and its specialization for $T = \partial F$, the subdifferential of $F$, explain why quadratic additive terms in Algorithm 3.1.1 were originally mandatory; see Rockafellar (1976b), and for some recent work in this field see Mahey et al. (1992). The quadratic proximal minimization algorithm was widely studied in recent years. In addition to the sources mentioned above we note some recent papers by Güler (1991), Eckstein and Bertsekas (1992) and Wright (1990). See also Chapter XV.4 in Hiriart-Urruty and Lemaréchal (1993) and Chapter 3.4.3 in Bertsekas and Tsitsiklis (1989). The idea of replacing quadratic additive terms by nonquadratic ones already exists with respect to other algorithms; see, for example, Bertsekas (1982, Chapter 5) and Polyak (1992). Ben-Tal (1985) and Ben-Tal and Teboulle (1987) used $\varphi$-divergences (see Section 2.5) as penalty functions in stochastic programming, and Ben-Tal, Charnes and Teboulle (1989) used them to develop entropic means. Eggermont (1990) studied an algorithm of the same form as Algorithm 3.1.1 with the Kullback-Leibler distance (3.38) replacing the quadratic term and with $X = \mathbb{R}^n_+$. This was done within a framework of multiplicative iterative algorithms; see De Pierro (1991) for a review of such methods.

The PMD algorithm (Algorithm 3.1.2) was proposed by Censor and Zenios (1992), who established its convergence. Eckstein (1992, 1996) extended the PMD algorithm to the setting of maximal monotone operators and to approximate iterations, respectively. Further results related to Eckstein's work appear in Censor and Reich (1996). The ergodic convergence of the PMD algorithm in the setting of monotone operators was investigated by Güler (1994). Ferris (1991) showed finite termination of the proximal point algorithm under an additional assumption that the objective function has a weak sharp minimum. See also Ferris (1988) and Burke and Ferris (1993) for related work on weak sharp minima in mathematical programming. Chen and Teboulle (1993) supplied an alternative convergence proof for the PMD algorithm, and Teboulle (1992b) discusses entropic proximal maps, that use either a $D_f(x, y)$ generalized distance or a $d_\varphi(x, y)$ Csiszár $\varphi$-divergence instead of quadratic terms to construct generalized augmented Lagrangians. Iusem, Teboulle, and co-workers have extensively studied a companion algorithm to the PMD algorithm, in which the original quadratic term in (3.6) is replaced by $\varphi$-divergences rather than by $D_f(x, y)$ generalized distances; see Iusem, Svaiter, and Teboulle (1994). Iusem and Teboulle (1993, 1995) study convergence rates of both the PMD algorithm and their entropic proximal method with $\varphi$-divergences (EPMC). Ibaraki, Fukushima, and Ibaraki (1992) develop a primal-dual version of the original proximal point algorithm for linearly constrained convex programming problems. Auslender and Haddou (1995) and Burachik (1995) extend the PMD and EPMC algorithms to variational inequalities.

In a recent report Censor, Iusem, and Zenios (1995) extended the use of Bregman functions to the variational inequality problem on convex sets. Kiwiel (1995a) studied proximal minimization methods with generalized Bregman functions that were introduced in Kiwiel (1994). Finally, we mention the recent work of Iusem, Svaiter, and da Cruz Neto (1995), who study the relationship of the PMD algorithm with central paths defined by arbitrary barriers and with Cauchy trajectories in Riemannian manifolds. All these concepts are important in connection with variational inequalities problems. Bertsekas and Tseng (1994) proposed and analyzed the so-called partial proximal minimization algorithms.

The Weierstrass theorem mentioned in the proof of Proposition 3.2.1 can be found in any text on analysis, e.g., Rudin (1953). The definition and properties of directional derivatives used in the proof of Proposition 3.2.3 can be found, e.g., in Rockafellar (1970).

**3.3** The idea of replacing a linear programming problem with a sequence of entropy problems was heuristically suggested in an unpublished report by Eriksson (1985). He also discusses a specific strategy for choosing the parameters $\{c_\nu\}$ and a solution algorithm, but gives no convergence analysis.

CHAPTER 4

# Penalty Methods, Barrier Methods and Augmented Lagrangians

---

We now turn to study some general optimization methods upon which we can later base the design of model decomposition algorithms (see Chapter 1). These general methods for constrained nonlinear optimization are the penalty methods that employ penalty or barrier functions giving rise to so-called exterior penalty methods or interior penalty methods, respectively. These methods are closely related to the primal-dual Lagrangian algorithmic scheme and to the method of augmented Lagrangians (sometimes referred to as the multiplier method), which we also discuss here.

Such algorithms are not only considered most efficacious in solving large-scale constrained optimization problems, but are also suitable *problem modifiers* in our *model decomposition algorithms* (Chapter 7). Barrier-function methods are also helpful when studying the primal-dual path following algorithm for linear programming (see Chapter 8). This algorithm is developed using a logarithmic barrier function. Finally, there is a close relationship between the method of augmented Lagrangians and the proximal minimization algorithms (see Chapter 3), which should not be overlooked.

Section 4.1 examines penalty functions and their use in penalty methods, wherein the optimal solution of the problem is approached from outside the feasible set (exterior penalty). In Section 4.2 we describe barrier functions and their use in approaching the optimal solution from within the feasible set. Section 4.3 explains the basics of duality theory and the primal-dual algorithmic framework. Section 4.4 gives an account of augmented Lagrangian methods and how they are related to the proximal minimization approach developed earlier (see Chapter 3). Notes and References are given in Section 4.5.

## 4.1 Penalty Methods

Penalty methods for constrained nonlinear optimization are based on the following idea: a penalty function is defined that imposes a penalty for constraint violations by raising the value of the objective function that has to be minimized. The penalty is greater for points that are farther away from the feasible set and is equal to zero at feasible points that satisfy all constraints of the problem. The optimal solution to the constrained problem is obtained as the limit of solutions of a sequence of unconstrained penalized problems.

To elaborate this idea we consider the following optimization problem

$$\text{Minimize} \quad F(x) \tag{4.1}$$

$$\text{s.t.} \quad x \in \Omega, \tag{4.2}$$

where $F:\mathbb{R}^n \rightarrow \mathbb{R}$ is a continuous function and $\Omega \subseteq \mathbb{R}^n$ is the feasible set.

**Definition 4.1.1 (Penalty Functions)** *Given a nonempty set $\Omega \subseteq \mathbb{R}^n$, a function $p : \mathbb{R}^n \rightarrow \mathbb{R}$ is called an (exterior) penalty function with respect to $\Omega$ if the following hold: (i) $p(x)$ is continuous, (ii) $p(x) \geq 0$, for all $x \in \mathbb{R}^n$, and (iii) $p(x) = 0$ if and only if $x \in \Omega$.*

**Example 4.1.1** *When $\Omega \doteq \{x \in \mathbb{R}^n \mid g_i(x) \leq 0,\ i = 1, 2, \ldots, m\}$ describes the feasible set by inequality constraints, where $g_i : \mathbb{R}^n \rightarrow \mathbb{R}$, a commonly used penalty function is the quadratic penalty function:*

$$p(x) \doteq \frac{1}{2} \parallel g^+(x) \parallel^2 = \frac{1}{2} \sum_{i=1}^{m} (\max(0, g_i(x)))^2 . \tag{4.3}$$

Using a penalty function we can associate with (4.1)–(4.2) an unconstrained penalized problem

$$\underset{x \in \mathbb{R}^n}{\text{Minimize}} \ (F(x) + cp(x)), \tag{4.4}$$

with a positive penalty parameter $c > 0$. If the penalized problem (4.4) yields the exact solution of the original problem (4.1)–(4.2) for a finite value of the parameter $c$, then we call $p(x)$ an *exact penalty function* for the original problem. If this is not the case then the penalty method for solving (4.1)–(4.2) works as follows.

### Algorithm 4.1.1 The Penalty Algorithmic Scheme

**Step 0:** (Initialization.) Let $\{c_\nu\}_{\nu=0}^{\infty}$ be a monotone sequence of penalty parameters such that $c_\nu \geq 0$ and $c_{\nu+1} > c_\nu$, for all $\nu \geq 0$, and $\lim_{\nu \rightarrow \infty} c_\nu = +\infty$.

**Step 1:** (Iterative step.) For every $\nu \geq 0$, solve the unconstrained penalized optimization problem and set

$$x^\nu = \underset{x \in \mathbb{R}^n}{\text{argmin}} \ (F(x) + c_\nu p(x)). \tag{4.5}$$

The following convergence result gives conditions, on the original problem, under which any accumulation point of a sequence $\{x^\nu\}_{\nu=0}^\infty$, generated by Algorithm 4.1.1, is a solution of the original problem.

**Theorem 4.1.1** *Let $p(x)$ be a penalty function for problem (4.1)–(4.2) in which $F(x)$ is a continuous function, $\Omega$ is a closed set and one of the following conditions holds:*

*(i)* $\lim_{\|x\|\to\infty} F(x) = +\infty$,

*(ii)* $\Omega$ *is bounded and* $\lim_{\|x\|\to\infty} p(x) = +\infty$.

*Then any sequence $\{x^\nu\}_{\nu=0}^\infty$, generated by the penalty algorithmic scheme (Algorithm 4.1.1) has at least one accumulation point, every such point is an optimal solution of problem (4.1)–(4.2), and* $\lim_{\nu\to\infty} p(x^\nu) = 0$.

**Proof** The Weierstrass theorem, which says that a real continuous function on a compact set in $\mathbb{R}^n$ obtains its minimum on the set, guarantees the existence of an optimal solution $x^* \in \Omega$ of the original problem. This means that $F(x) \geq F(x^*)$, for all $x \in \Omega$. This is obvious if condition $(ii)$ holds because then $\Omega$ is closed and bounded. If condition $(i)$ holds then we take any $\overline{x} \in \mathbb{R}^n$ and deduce that there exists $M > 0$, such that $\| x \| \geq M$ implies $F(x) \geq F(\overline{x})$. Thus, solving the problem (4.1)–(4.2) reduces to the minimization of $F(x)$ over the set $\Omega \cap \{x \in \mathbb{R}^n \mid \| x \| \leq M\}$, which is compact, and the Weierstrass theorem applies. Next we see that from the nonnegativity of the penalty function and the penalty parameters, from the strict monotone increase of the latter, and from (4.5) we have, for all $\nu \geq 0$,

$$F(x^{\nu+1}) + c_{\nu+1}p(x^{\nu+1}) > F(x^{\nu+1}) + c_\nu p(x^{\nu+1}) \geq F(x^\nu) + c_\nu p(x^\nu). \quad (4.6)$$

Also, from the nature of $x^\nu$ and $x^{\nu+1}$ as minimizers of the penalized problems at the respective iterations, we have, for all $\nu \geq 0$,

$$F(x^\nu) + c_\nu p(x^\nu) \leq F(x^{\nu+1}) + c_\nu p(x^{\nu+1}).$$

The last two inequalities yield $(c_{\nu+1} - c_\nu)p(x^{\nu+1}) \leq (c_{\nu+1} - c_\nu)p(x^\nu)$, which shows that, for all $\nu \geq 0$,

$$p(x^{\nu+1}) \leq p(x^\nu), \quad (4.7)$$

because $\{c_\nu\}_{\nu=0}^\infty$ is monotonically increasing. The optimality of $x^* \in \Omega$ makes $p(x^*) = 0$ and it follows that, for all $\nu \geq 0$,

$$F(x^\nu) \leq F(x^\nu) + c_\nu p(x^\nu) \leq F(x^*) + c_\nu p(x^*) = F(x^*). \quad (4.8)$$

Any sequence $\{x^\nu\}_{\nu=0}^\infty$, generated by the penalty method, must be bounded because otherwise it would contradict (4.8) if condition $(i)$ holds or it would contradict $p(x^\nu) \leq p(x^0)$, for all $\nu \geq 0$, which follows from (4.7) if condition $(ii)$ holds. This proves the existence of an accumulation point $\hat{x}$, and we let $\{x^\nu\}_{\nu \in K}$, $K \subseteq N_0 \doteq \{0, 1, 2, \ldots\}$ be a subsequence converging to $\hat{x}$. From continuity of $F$ and from (4.8) we obtain $\lim_{\nu \in K} F(x^\nu) = F(\hat{x}) \leq F(x^*)$. The sequence $\{F(x^\nu) + c_\nu p(x^\nu)\}_{\nu=0}^{\infty}$ is, by (4.6), monotonically increasing and, by (4.8), bounded from above; thus

$$\lim_{\nu \to \infty} (F(x^\nu) + c_\nu p(x^\nu)) \doteq q^* \leq F(x^*),$$

which shows that, for $\nu \in K$, $\lim_{\nu \in K} c_\nu p(x^\nu) = q^* - F(\hat{x})$. Therefore, since the sequence $\{c_\nu\}$ is unbounded, we must have, by (4.7), $\lim_{\nu \in K} p(x^\nu) = 0$, which ensures that $p(\hat{x}) = 0$, by continuity of $p(x)$. This means that $\hat{x} \in \Omega$ and since we showed above that $F(\hat{x}) \leq F(x^*)$, it follows that $F(\hat{x}) = F(x^*)$ and thus $\hat{x}$ is an optimal solution. ■

This theorem shows that the penalty algorithmic scheme is an *exterior* penalty method because the iterates $\{x^\nu\}_{\nu=0}^{\infty}$ have a decreasing level of infeasibility—if we measure the level of infeasibility by the values of the penalty function, see (4.7); and the function values $\{F(x^\nu)\}_{\nu=0}^{\infty}$ approach the minimal value $F(x^*)$ from below. In contrast to this behavior, the barrier methods (see below) approach an optimal point from the interior of the feasible set.

## 4.2 Barrier Methods

Barrier methods for the solution of constrained nonlinear optimization problems are similar to penalty methods in that they replace the original problem by a sequence of unconstrained problems whose objective functions are augmented modifications of the original $F(x)$. The main difference, however, lies in the nature of the function added to $F(x)$. Here this function builds a barrier along the boundary of the feasible set $\Omega$, which prevents the iterates produced by the algorithm from approaching the boundary while searching for a minimum inside $\Omega$. We must assume here that the set $\Omega$ in problem (4.1)–(4.2) is such that its interior is not empty and each point on the boundary is the limit of a sequence of points of its interior int $\Omega$. Such sets are referred to as *robust.*

**Definition 4.2.1 (Barrier Functions)** *Given a robust set $\Omega \subseteq \mathbb{R}^n$, a function $q : \mathbb{R}^n \to \mathbb{R}$ is called a barrier function with respect to $\Omega$ (also called an interior penalty function) if the following hold: (i) $q(x)$ is continuous on* int $\Omega$, *(ii) $q(x) \geq 0$, for all $x \in$* int $\Omega$, *and (iii)* $\lim q(x) = +\infty$, *as $x$ tends to any boundary point of $\Omega$.*

**Example 4.2.1** *When $\Omega$ is as in Example 4.1.1 and we assume that the functions $g_i$ are continuous on $\mathbb{R}^n$, for all $i = 1, 2, \ldots, m$, and that the set*

*$\Omega$ is robust and that* $\text{int}\ \Omega \doteq \{x \in \mathbb{R}^n \mid g_i(x) < 0,\ i = 1, 2, \ldots, m\}$, *then the function* $q(x) \doteq -\sum_{i=1}^{m} \frac{1}{g_i(x)}$ *is a barrier function for* $\Omega$.

This time we associate with the original constrained problem a modified constrained problem with an objective function augmented by a barrier function, namely,

$$\text{Minimize} \quad F(x) + cq(x) \tag{4.9}$$

$$\text{s.t.} \quad x \in \text{int}\ \Omega, \tag{4.10}$$

with a positive *barrier parameter* $c > 0$. Although this is a constrained problem, it can still be solved by a search technique for unconstrained problems because if such a search is initialized at a point in $\text{int}\,\Omega$, then all points obtained by a search technique for the solution of the problem $\min_{x \in \mathbb{R}^n}(F(x) + cq(x))$ will stay in int $\Omega$. The barrier method consists of solving a sequence of such modified problems with an ever-decreasing sequence $c_\nu \to 0$ of barrier parameters, as described by the following algorithmic scheme.

**Algorithm 4.2.1 The Barrier Algorithmic Scheme**

**Step 0:** (Initialization.) Let $\{c_\nu\}_{\nu=0}^{\infty}$ be a monotone sequence of barrier parameters such that $c_\nu > 0$ and $c_{\nu+1} < c_\nu$, for all $\nu \geq 0$, and $\lim_{\nu \to \infty} c_\nu = 0$.

**Step 1:** (Iterative step.) For every $\nu \geq 0$, solve the modified problem and set

$$x^\nu = \operatorname*{argmin}_{x \in \text{int}\ \Omega} (F(x) + c_\nu q(x)). \tag{4.11}$$

The next theorem shows that any accumulation point of a sequence $\{x^\nu\}_{\nu=0}^{\infty}$ generated by this algorithm is a solution of the original problem.

**Theorem 4.2.1** *Let $q(x)$ be a barrier function for problem (4.1)–(4.2) with a nonempty, closed and robust feasible set $\Omega$, and with $F(x)$ a continuous function, and assume that one of the following conditions holds:*

*(i)* $\lim_{\|x\| \to \infty} F(x) = +\infty$,

*(ii) $\Omega$ is bounded.*

*Then any sequence $\{x^\nu\}_{\nu=0}^{\infty}$ generated by the barrier algorithmic scheme (Algorithm 4.2.1) has at least one accumulation point, every such point is an optimal solution of problem (4.1)–(4.2), and* $\lim_{\nu \to \infty} c_\nu q(x^\nu) = 0$.

**Proof** Using the same argument as in the proof of Theorem 4.1.1, the existence of $x^* \in \Omega$, which is optimal for the original problem, is guaranteed. Now we can write, for every $\nu \geq 0$,

$$F(x^*) \leq F(x^\nu) \leq F(x^\nu) + c_\nu q(x^\nu), \tag{4.12}$$

and use the continuity of $F$ and the robustness of $\Omega$ to conclude that, for every $\epsilon > 0$, there exists a point $\tilde{x} \in \text{int } \Omega$ such that $F(\tilde{x}) \leq F(x^*) + \epsilon$. Therefore, for every $\nu \geq 0$, we get

$$F(x^*) + \epsilon + c_\nu q(\tilde{x}) \geq F(\tilde{x}) + c_\nu q(\tilde{x}) \geq F(x^\nu) + c_\nu q(x^\nu),$$

which shows that $\lim_{\nu \to \infty} F(x^\nu) + c_\nu q(x^\nu) \leq F(x^*) + \epsilon$. Combining this with the limiting behavior of (4.12) we get

$$\lim_{\nu \to \infty} F(x^\nu) + c_\nu q(x^\nu) = F(x^*) = \lim_{\nu \to \infty} F(x^\nu) \tag{4.13}$$

and $\lim_{\nu \to \infty} c_\nu q(x^\nu) = 0$. The sequence $\{x^\nu\}_{\nu=0}^\infty$ must be bounded because otherwise it would contradict (4.13) if condition $(i)$ holds or, if condition $(ii)$ holds, this follows from (4.11) since the interior of a bounded set is also bounded. Let $\hat{x}$ be an accumulation point (guaranteed to exist since $\{x^\nu\}_{\nu=0}^\infty$ is bounded) and assume that $\lim_{\nu \in K} x^\nu = \hat{x}$, for some infinite subset $K \subseteq N_0$. Then continuity of $F$ yields $F(\hat{x}) = F(x^*)$ and, since $\Omega$ is closed, $\hat{x} \in \Omega$. Thus, every accumulation point of $\{x^\nu\}_{\nu=0}^\infty$ is a minimum of the original problem. ■

It is easy to show, following arguments similar to those in the proof of Theorem 4.1.1, that in the barrier algorithmic scheme we get, for all $\nu \geq 0$,

$$\begin{aligned} F(x^{\nu+1}) + c_{\nu+1} q(x^{\nu+1}) &< F(x^\nu) + c_\nu q(x^\nu), \\ q(x^\nu) &\leq q(x^{\nu+1}), \\ F(x^{\nu+1}) &\leq F(x^\nu), \end{aligned}$$

showing that the values of the original objective function monotonically decrease during the iteration of the barrier method.

**Example 4.2.2 (The Logarithmic Barrier Function)** *For a set $\Omega$ as described in Example 4.2.1, the function $q(x) \doteq -\sum_{i=1}^m \log g_i(x)$ is a barrier function. Observe the relation of this with Burg's entropy function defined in (6.149); see also Section 6.9.2.*

## 4.3 The Primal-Dual Algorithmic Scheme

Duality, Lagrange multipliers, and primal-dual algorithms for constrained optimization play a prominent role in optimization theory. In this section we explain their basic ideas in order to use them in conjunction with penalty function methods. The combination of primal-dual algorithms with penalty function methods gives rise to the method of augmented Lagrangians, the latter of which is used in the development of model decompositions (see

Chapters 7 and 8). For simplicity of presentation we deal here only with equality constrained problems of the form

$$\text{Minimize} \quad F(x) \tag{4.14}$$

$$\text{s.t.} \quad h_i(x) = 0, \qquad i = 1, 2, \ldots, m, \tag{4.15}$$

$$x \in \Omega, \tag{4.16}$$

where $F : \mathbb{R}^n \to \mathbb{R}$, $h_i : \mathbb{R}^n \to \mathbb{R}$, for all $i = 1, 2, \ldots, m$ are given functions and $\Omega \subseteq \mathbb{R}^n$ is a given subset. The *Lagrangian* of this problem is defined as

$$L(x, \pi) \doteq F(x) + \sum_{i=1}^{m} \pi_i h_i(x), \tag{4.17}$$

for every $x \in \mathbb{R}^n$ and every $\pi = (\pi_i) \in \mathbb{R}^m$. The *dual function* associated with this problem is

$$g(\pi) \doteq \min_{x \in \Omega} L(x, \pi), \tag{4.18}$$

and the unconstrained optimization problem

$$\text{Maximize} \quad g(\pi) \tag{4.19}$$

$$\text{s.t.} \quad \pi \in \mathbb{R}^m, \tag{4.20}$$

is the dual problem to (4.14)–(4.16). The known *local duality theorem* relates the solutions $x^*$ and $\pi^*$, of the primal (original) problem and its dual problem, respectively, for the case when $\Omega = \mathbb{R}^n$, and we present it here without proof.

The following terms are used in the theorem. A point $x^*$ is a *local extremum* (minimum or maximum) point of an optimization problem if it is a *feasible point*, i.e., fulfills the constraints of the problem, and there exists a neighborhood of the point where all values of the *objective function* $F(x)$ are not smaller than (or, for maximization problems, not greater than) $F(x^*)$. A *Lagrange multiplier vector* $\pi^*$ for problem (4.14)–(4.15) is a vector that satisfies the first-order necessary local optimality conditions at a local extremum point $x^*$, i.e.,

$$\nabla F(x^*) + \sum_{i=1}^{m} \pi_i^* \nabla h_i(x^*) = 0. \tag{4.21}$$

A point $x^*$ satisfying the constraints $h_i(x^*) = 0$, $i = 1, 2, \ldots, m$ is said to be a *regular point* of these constraints if the gradient vectors $\{\nabla h_i(x^*)\}_{i=1}^m$ are linearly independent. The Hessian of a real-valued, twice continuously differentiable function $u : \mathbb{R}^n \to \mathbb{R}$, denoted by $\nabla^2 u(x)$, is an $n \times n$ matrix whose $(i,j)$th entry is $(\nabla^2 u(x))_{i,j} = \dfrac{\partial^2 u(x)}{\partial x_i \partial x_j}$ and, since $\dfrac{\partial^2 u}{\partial x_i \partial x_j} = \dfrac{\partial^2 u}{\partial x_j \partial x_i}$, it is symmetric. The Hessian of the Lagrangian (4.17), with respect to $x$, at the point $x^*$ is

$$L^* \doteq \nabla^2 L(x^*) = \nabla^2 F(x^*) + \sum_{i=1}^{m} \pi_i \nabla^2 h_i(x^*), \tag{4.22}$$

where $\nabla^2 F$ and $\nabla^2 h_i$ are the Hessian matrices of the respective functions. The local duality theorem requires that this matrix be positive definite, i.e., that $\langle x, L^* x\rangle > 0$, for every $x \in \mathbb{R}^n$, such that $x \neq 0$.

**Theorem 4.3.1 (Local Duality Theorem)** *Suppose that $x^*$ is a local minimum point of the optimization problem (4.14)–(4.15) with a minimal value $F(x^*) = r_*$ and a Lagrange multiplier vector $\pi^*$. Suppose also that $x^*$ is a regular point of the constraints (4.15) and that the corresponding Hessian of the Lagrangian $L^* = \nabla^2 L(x^*)$ is positive definite. Then the dual problem (4.19)–(4.20) has a local maximum at $\pi^*$ with a maximal value $g(\pi^*) = r_*$, and $x^*$ as the corresponding point to $\pi^*$ in (4.18).*

A similar result can be obtained for problems having inequality constraints in addition to the equality constraints. If appropriate convexity assumptions are added, then the local extrema (mentioned in the theorem) can be replaced by global extrema. Finally, it is not necessary to include all the constraints of a problem in the definition of the dual functional $g(\pi)$. Local duality can be defined with respect to any subset of functional constraints; this is called *partial duality.* For example, the constraints $h_i(x) = 0$ might be separated into two subsets, easy and hard constraints. The hard ones can be dualized, i.e., removed from the constraint set and incorporated into the Lagrangian, while the easy ones remain as constraints.

The primal-dual approach for constrained optimization problems aims, in view of this duality theorem, at solving the dual unconstrained problem (4.19)–(4.20). This is done by an iterative scheme which alternates between the minimization of the Lagrangian in (4.18) and the application of a steepest ascent iteration to the dual problem.

#### Algorithm 4.3.1 The Primal-Dual Algorithmic Scheme

**Step 0:** (Initialization.) $\pi^0 \in \mathbb{R}^m$ is arbitrary.

**Step 1:** (Iterative step.) Given $\pi^\nu$, for some $\nu \geq 0$,

(i) Solve the minimization problem

$$x^\nu = \operatorname*{argmin}_{x \in \Omega} L(x, \pi^\nu). \tag{4.23}$$

(ii) Do a steepest ascent step to calculate $\pi^{\nu+1}$ via

$$\pi_i^{\nu+1} = \pi_i^{\nu} + c_\nu h_i(x^\nu), \quad i = 1, 2, \ldots, m, \tag{4.24}$$

where $\{c_\nu\}_{\nu=0}^{\infty}$ is a nondecreasing sequence of positive numbers, called *stepsizes*, and return to the beginning of the iterative step.

The general scheme of Algorithm 4.3.1 can give rise to various specific algorithms, depending on the particular algorithm employed in (4.23) and the choice of stepsizes in (4.24). Although a detailed study of general primal-dual algorithms is beyond the scope of this book, we will examine in detail the convergence of certain specific primal-dual algorithms (see Chapter 6, in particular, the discussion in Section 6.3).

From the structure of Algorithm 4.3.1, we can assess some of its drawbacks. First, the positive definiteness of $L^*$—referred to as *local convexity*—is necessary for Theorem 4.3.1. Without it the dual problem (4.19)–(4.20) is not well-defined and the iteration (4.24) is not meaningful. Second, slow convergence of the iterates $\{\pi^\nu\}_{\nu=0}^{\infty}$ in (4.24) usually necessitates many repeated returns to the minimization in (4.23), thus limiting the usefulness of the whole algorithm to problems whose structure makes these minimizations efficient. This happens in several important areas of applications when the original problem is *separable.* This means (see Definition 1.3.1) that the vector $x \in \mathbb{R}^m$ can be partitioned into $q$ subvectors $x^l \in \mathbb{R}^{n_l}, l = 1, 2, \ldots, q$, such that $\sum_{l=1}^{q} n_l = n$, and the objective function and the constraints separate into sums of functions in the following form:

$$\text{Minimize} \quad \sum_{l=1}^{q} F_l(x^l) \tag{4.25}$$

$$\text{s.t.} \quad \sum_{l=1}^{q} \overline{h}_l\,(x^l) = 0. \tag{4.26}$$

In this formulation, $\overline{h}_l : \mathbb{R}^{n_l} \to \mathbb{R}^m$ and the equation (4.26) represent the original equality constraints $h_i(x) = 0$, $i = 1, 2, \ldots, m$. When inequality constraints are present in the original problem, i.e., $g_i(x) \leq 0$, $i = 1, 2, \ldots, p$, then in the separable formulation they take the form $\sum_{l=1}^{q} \overline{g}_l(x_l) \leq 0$, where $\overline{g}_l : \mathbb{R}^{n_l} \to \mathbb{R}^p$. For such separable problems the minimization involved in (4.23) decomposes into small subproblems. This class of problems is later treated in detail (see Chapter 7).

## 4.4 Augmented Lagrangian Methods

Augmented Lagrangian methods, also known as *multiplier methods,* are a class of very effective general methods for constrained nonlinear optimization. They combine the basic idea of penalty function methods (see Section 4.1) with the primal-dual algorithmic approach, based on classic Lagrange duality theory (see above). Penalty methods suffer from slow convergence and numerical instabilities that result from the need to increase the penalty parameters $c_\nu$ in equation (4.5) in order to achieve feasibility. Classical Lagrangian methods of the primal-dual type have the drawbacks explained in the previous section. But augmented Lagrangian methods, when properly employed and under mild technical conditions, can be made to behave better than both methods.

We look again at the constrained nonlinear optimization problem (4.14)–(4.16), assuming now that $F : \mathbb{R}^n \longrightarrow \mathbb{R}$ is a convex function and that $\Omega \subseteq \mathbb{R}^n$ is a nonempty polyhedral set, i.e., it can be expressed as the intersection of finitely many half-spaces. The basic idea of augmented Lagrangian methods is to replace first this original problem by an equivalent problem

$$\text{Minimize} \quad F(x) + \frac{1}{c} f(ch(x)) \tag{4.27}$$

$$\text{s.t.} \quad h(x) = 0, \tag{4.28}$$

$$x \in \Omega, \tag{4.29}$$

where $h : \mathbb{R}^n \longrightarrow \mathbb{R}^m$ is the vector of functions $h = (h_i)$, thus (4.28) represents the same equality constraints as (4.15). For this problem to be equivalent to the original problem, the function $f$ must have the property that $f(0) = 0$ and the scalar parameter $c$ is positive. The Lagrangian of this equivalent problem, called the *augmented Lagrangian*, has the form

$$L_c(x, \pi) \doteq F(x) + \langle \pi, h(x) \rangle + \frac{1}{c} f(ch(x)), \tag{4.30}$$

and the dual problem is

$$\text{Maximize} \quad g_c(\pi) \tag{4.31}$$

$$\text{s.t.} \quad \pi \in \mathbb{R}^m, \tag{4.32}$$

where the dual function is

$$g_c(\pi) \doteq \min_{x \in \Omega} L_c(x, \pi). \tag{4.33}$$

Thus, keeping the original constraints (4.15) and (4.28) in the equivalent problem, in addition to having them built into (4.27), really means that a penalty term of the form $(1/c)f(ch(x))$ has been added to the Lagrangian rather than to the original objective function $F(x)$. To generate an algorithmic scheme an iterative process is used, which alternates between the minimization of the augmented Lagrangian (4.30) and the application of an ascent iteration to the dual problem.

The method was originally proposed with the function $f(z) = \frac{1}{2} \parallel z \parallel^2$ and was later extended to functions of the form $f(z) = \sum_{i=1}^{m} \phi(z_i)$, where $\phi : \mathbb{R} \to \mathbb{R}$ belongs to the class of penalty functions $P_E$ defined as follows.

**Definition 4.4.1 (The Class $P_E$ of Penalty Functions)** *The function $\phi : \mathbb{R} \to \mathbb{R}$ belongs to the class of penalty functions $P_E$ if it has the following properties:*

*(i) $\phi$ is continuously differentiable and strictly convex on $\mathbb{R}$.*

*(ii) $\phi(0) = 0$ and $\frac{d\phi}{dt}(0) = 0$.*

*(iii) $\lim_{t \to -\infty} \frac{d\phi}{dt}(t) = -\infty$ and $\lim_{t \to +\infty} \frac{d\phi}{dt}(t) = +\infty$.*

Examples of functions in this class are:

$\phi(t) \doteq \frac{1}{2}t^2$, this is the classic case with $f(z) = \frac{1}{2} \parallel z \parallel^2$ .
$\phi(t) \doteq \rho^{-1}|t|^\rho$, for $\rho > 1$.
$\phi(t) \doteq \rho^{-1}|t|^\rho + \frac{1}{2}t^2$, for $\rho > 1$.
$\phi(t) \doteq \cosh(t) - 1$, where $\cosh(t)$ is the hyperbolic cosine function.

The next proposition expresses the gradient of the dual function in terms of the constraint functions. We will assume that $x^*$ is a local solution of problem (4.27)–(4.29) and a regular point of the constraints (4.28). Denote by $x(\pi)$ the unique solution of (4.33) in the neighborhood of $x^*$.

**Proposition 4.4.1** *If all involved functions $F(x)$, $h_i(x)$, $i = 1, 2, \ldots, m$, and $x(\pi)$ are continuously differentiable, then*

$$\nabla g_c(\pi) = h(x(\pi)). \tag{4.34}$$

**Proof** From (4.33) and the definition of $x(\pi)$ we have

$$g_c(\pi) = F(x(\pi)) + \langle \pi, h(x(\pi)) \rangle + \frac{1}{c} f(ch(x(\pi))). \tag{4.35}$$

From this we obtain

$$\frac{\partial g_c(\pi)}{\partial \pi_i} = \sum_{j=1}^{n} \frac{\partial F}{\partial x_j} \frac{\partial x_j}{\partial \pi_i} + h_i(x(\pi)) + \sum_{s=1}^{m} \pi_s \left( \sum_{j=1}^{n} \frac{\partial h_s}{\partial x_j} \frac{\partial x_j}{\partial \pi_i} \right)$$
$$+ \frac{1}{c} \sum_{s=1}^{m} \frac{\partial f}{\partial z_s} c \left( \sum_{t=1}^{n} \frac{\partial h_s}{\partial x_t} \frac{\partial x_t}{\partial \pi_i} \right). \tag{4.36}$$

Denoting by $\nabla h(x)$ the $n \times m$ matrix whose $i$th column consists of the gradient $\nabla h_i(x)$, and by $\nabla x(\pi)$ the $m \times n$ matrix whose $j$th column consists of the gradient $\nabla x_j(\pi)$, we may rewrite this as

$$\nabla g_c(\pi) = \nabla x(\pi)\left(\nabla F(x(\pi)) + \nabla h(x(\pi))\pi + \frac{1}{c}\nabla h(x(\pi))c\nabla f(h(x(\pi)))\right) \\ +h(x(\pi)). \tag{4.37}$$

The vector which multiplies the matrix $\nabla x(\pi)$ is nothing but the gradient with respect to $x$ of the augmented Lagrangian $\nabla_x L_c(x(\pi), \pi)$, which is zero by the definition of $x(\pi)$; thus we obtain (4.34). ■

Let $\phi \in P_E$ be any function of the class of penalty functions given in Definition 4.4.1. In order to formulate the augmented Lagrangian algorithmic scheme for this class of functions we need to define, for each $i = 1, 2, \ldots, m$, the scalar quantities

$$\beta_i(\pi, c) \doteq \int_0^1 \frac{d^2\phi}{dt^2}(\alpha c h_i(x(\pi)))d\alpha, \tag{4.38}$$

where the function under the integral sign is the second derivative of $\phi$ at a point $t = \alpha c h_i(x(\pi))$. Let $B(\pi, c) \doteq \text{diag}\ (\beta_1(\pi, c), \beta_2(\pi, c), \ldots, \beta_m(\pi, c))$ be the diagonal $m \times m$ matrix with the $\beta_i$'s on its diagonal. If $\phi(t) = \frac{1}{2}t^2$ we get $\beta_i(\pi, c) = 1$ for all $i$, and $B(\pi, c)$ is the identity matrix. In this case, the augmented Lagrangian algorithmic scheme, given next, becomes the classic *quadratic method of multipliers.*

### Algorithm 4.4.1 The Augmented Lagrangian Algorithmic Scheme

**Step 0:** (Initialization.) $\pi^0 \in \mathbb{R}^m$ is arbitrary.

**Step 1:** (Iterative step.) Given $\pi^\nu$, for some $\nu \geq 0$,

(i) Solve the minimization problem,

$$x^\nu = \operatorname*{argmin}_{x \in \Omega} L_{c_\nu}(x, \pi^\nu). \tag{4.39}$$

(ii) Do an ascent step to calculate $\pi^{\nu+1}$ via

$$\pi^{\nu+1} = \pi^\nu + c_\nu B(\pi^\nu, c_\nu)\nabla g_{c_\nu}(\pi^\nu), \tag{4.40}$$

where $\{c_\nu\}_{\nu=0}^\infty$ is an appropriately chosen sequence of positive parameters, and return to the beginning of the iterative step.

This is quite a general scheme which can be made specific according to the method chosen for minimization of the augmented Lagrangian in (4.39) and the choice of the sequence $\{c_\nu\}_{\nu=0}^\infty$. When calculating the entries of $B(\pi^\nu, c_\nu)$ according to (4.38) we should remember that by definition, $x(\pi^\nu) = x^\nu$. The iteration formula (4.40) is a deflected-gradient type iteration whose behavior depends on $\{c_\nu\}_{\nu=0}^\infty$ and on the matrix $B(\pi^\nu, c_\nu)$, which deflects the direction of the gradient $\nabla g_{c_\nu}(\pi^\nu)$. For functions $\phi \in P_E$, the ascent nature of (4.40) can be guaranteed. This formula can be simplified in the following way. By using a first-order Taylor expansion formula for the function $\frac{d\phi}{dt}$ around the point $t = 0$, and using the fact that $\phi(0) = 0$, we get

$$\frac{d\phi}{dt}(ch_i(x(\pi)) = \int_0^1 \frac{d^2\phi}{dt^2}(\alpha ch_i(x(\pi)))ch_i(x(\pi))d\alpha. \tag{4.41}$$

This means that $\frac{d\phi}{dt}(ch_i(x(\pi)) = \beta_i(\pi)ch_i(x(\pi))$, and so the iteration (4.40) is actually identical with

$$\pi_i^{\nu+1} = \pi_i^\nu + \frac{d\phi}{dt}(c_\nu h_i(x^\nu)). \tag{4.42}$$

Finally, it is worth noting that with constant parameters $c_\nu = c > 0$, for all $\nu \geq 0$, the sequence of matrices $\{B(\pi^\nu, c)\}_{\nu=0}^\infty$ tends to the identity matrix, as $x(\pi^\nu) \rightarrow x^*$, assuming that $\frac{d^2\phi}{dt^2}(0) = 1$, thus making (4.40) a fixed stepsize, steepest ascent iteration, as $\pi^\nu$ tends to $\pi^*$—the optimal solution of the dual problem.

We wish now to analyze the convergence of Algorithm 4.4.1 when applied to the original problem (4.27)–(4.29) for the case of linear equality constraints, i.e., $h(x) = Ax - b$, for some given $m \times n$ matrix $A$ and given vector $b \in \mathbb{R}^m$. It turns out that this algorithm is then closely related to the proximal minimization algorithm (Algorithm 3.1.2) that we studied earlier (see Chapter 3). Specifically, the augmented Lagrangian algorithm is then equivalent to the proximal minimization algorithm applied to the dual problem (4.31)–(4.32). Let us assume here that the minimum in (4.39) is always uniquely attained, and introduce the new vector variable $z \doteq Ax - b$. Then

$$\min_{x\in\Omega} L_{c_\nu}(x, \pi^\nu) = \min_{x\in\Omega}\left(F(x) + \langle \pi^\nu, Ax - b\rangle + \frac{1}{c_\nu} f(c_\nu(Ax - b))\right), \tag{4.43}$$

and we may write

$$\begin{pmatrix} x^\nu \\ z^\nu \end{pmatrix} = \underset{z=Ax-b,\ x\in\Omega,\ z\in\mathbb{R}^m}{\operatorname{argmin}} \left(F(x) + \langle \pi^\nu, z\rangle + \frac{1}{c_\nu} f(c_\nu z)\right), \tag{4.44}$$

where $z^\nu = Ax^\nu - b$. This constrained optimization problem has a dual function

$$\hat{g}(\mu) = \min_{x\in\Omega,\ z\in\mathbb{R}^m}\left(F(x) + \langle \pi^\nu, z\rangle + \frac{1}{c_\nu} f(c_\nu z) + \langle \mu, Ax - b - z\rangle\right), \tag{4.45}$$

and it decomposes into two minimizations:

$$\hat{g}(\mu) = \min_{x\in\Omega} \theta_1(x,\mu) + \min_{z\in\mathbb{R}^m} \theta_2(z,\mu), \tag{4.46}$$

where

$$\theta_1(x,\mu) \doteq F(x) + \langle \mu, Ax-b\rangle, \tag{4.47}$$

and

$$\theta_2(z,\mu) \doteq \langle \pi^\nu - \mu, z\rangle + \frac{1}{c_\nu} f(c_\nu z). \tag{4.48}$$

Applying a duality theorem under the appropriate conditions (such as Theorem 4.3.1) we are guaranteed that an optimal dual solution $\mu^*$ exists, i.e., $\mu^* = \text{argmax}_{\mu\in\mathbb{R}^m} \hat{g}(\mu)$. The first right-hand side summand of (4.46),

$$\min_{x\in\Omega} \theta_1(x,\mu) = g(\mu) = \min_{x\in\Omega} \left(F(x) + \langle \mu, Ax-b\rangle\right), \tag{4.49}$$

is precisely the dual function (4.18) of the original problem (4.27)–(4.29). The unconstrained minimum in the second summand of (4.46), when $\mu = \mu^*$, is attained at $z^\nu$, for which

$$\mu_i^* - \pi_i^\nu = \frac{d\phi}{dt}(c_\nu z_i^\nu), \quad \text{for } i = 1,2,\ldots,m, \tag{4.50}$$

as follows from solving the system $\dfrac{\partial\theta_2}{\partial z_i}(z,\mu^*) = 0$, $i = 1,2,\ldots,m$. From (4.34) with $h(x) = Ax-b$ along with (4.42) and (4.50) we obtain then that $\mu^* = \pi^{\nu+1}$. With similar calculations we know that the minimum of

$$\min_{z\in\mathbb{R}^m} \theta_2(z,\mu) = \min_{z\in\mathbb{R}^m} \left( \langle \pi^\nu - \mu, z\rangle + \frac{1}{c_\nu}\sum_{i=1}^m \phi(c_\nu z_i)\right) \tag{4.51}$$

is attained at a point $z^*$ for which

$$\mu_i - \pi_i^\nu = \frac{d\phi}{dt}(c_\nu z_i^*), \quad i = 1,2,\ldots,m. \tag{4.52}$$

In order to obtain the minimal value of (4.51), i.e., calculate $\theta_2(z^*,\mu)$, we must extract $z^*$ from (4.52) and substitute it into (4.51). This cannot be done unless we specify $\phi(t)$; taking, for example, $\phi(t) \doteq \frac{1}{2}t^2$ we obtain

$$\begin{aligned}
\theta_2(z^*,\mu) &= \langle \pi^\nu - \mu, z^*\rangle + \frac{1}{c_\nu}\sum_{i=1}^m \frac{1}{2} c_\nu^2 (z_i^*)^2 \\
&= \langle \pi^\nu - \mu, \frac{\mu - \pi^\nu}{c_\nu}\rangle + \frac{c_\nu}{2}\sum_{i=1}^m \left(\frac{\mu_i - \pi_i^\nu}{c_\nu}\right)^2 \\
&= -(1/2c_\nu) \parallel \mu - \pi^\nu \parallel^2 .
\end{aligned} \tag{4.53}$$

Combining this with (4.46)–(4.49) and with the fact that $\mu^* = \pi^{\nu+1}$, we

conclude that

$$\mu^* = \operatorname*{argmax}_{\mu \in \mathbb{R}^m} \hat{g}(\mu) = \pi^{\nu+1} = \operatorname*{argmax}_{\mu \in \mathbb{R}^m} (g(\mu) - \frac{1}{2c_\nu} \| \mu - \pi^\nu \|^2), \quad (4.54)$$

which is precisely the quadratic proximal minimization algorithm (i.e., Algorithm 3.1.2 with $f(x) = \frac{1}{2} \| x \|^2$) applied to the dual problem (4.19)–(4.20) of the original problem. Whether or not this equivalence is true—or under what conditions it could be true—for general Bregman functions, is still an open problem. Another specific instance in which a similar equivalence holds is between the *exponential method of multipliers* and the PMD algorithm (Algorithm 3.1.2) with an entropic distance function (see Section 4.5 for references).

## 4.5 Notes and References

Penalty and barrier methods and augmented Lagrangians are important classes of methods for constrained nonlinear optimization; see, for example, Gill et al. (1989) for a survey, and Luenberger (1984) for a textbook treatment.

**4.1** A thorough understanding of the properties of penalty functions is important for the study of augmented Lagrangian methods. Consult the survey paper by Bertsekas (1976) or the general remarks in the introduction to his book (1982). Fiacco and McCormick (1968) pioneered the use of quadratic penalty functions in sequential unconstrained minimization; see also Polak (1971). The suggestions to study constrained problems via related unconstrained problems seem to originate in the work of Courant (1943) and Frisch (1955). Penalty methods are discussed in many advanced books in the field, see, e.g., Fletcher (1987) or Minoux (1986), whose presentations provide the foundation for our discussion. A recent unified treatment and review of exact penalty functions and methods appear in Burke (1991).

**4.2** Frisch (1955) proposed the logarithmic barrier function. Carroll (1961) proposed the function presented in Example 4.2.1. Most books mentioned above contain material on barrier function methods; see also the discussion in Gill, Murray, and Wright (1981) which is geared toward practical applications of the methods.

**4.3** We limited the scope of our discussion of this central subject in constrained optimization to the basics of sequential Lagrangian minimization. For further details on duality theory in convex analysis see Rockafellar (1970) or Hiriart-Urruty and Lemaréchal (1993). Luenberger (1984) or Bazaraa and Shetty (1979) can be used as introductions. For a review see Rockafellar (1976c).

**4.4** Augmented Lagrangian methods (also referred to as multiplier methods) originated with the works of Hestenes (1969) and Powell (1969);

see also Hestenes (1975). Extended treatments of this subject can be found in almost every modern advanced book in the field. An extensive treatise is Bertsekas (1982) where one can also find a thorough treatment of nonquadratic penalty functions and an extensive list of references to the contributions made to this field. The exponential multiplier method, originally proposed by Kort and Bertsekas (1972), has been recently analyzed by Tseng and Bertsekas (1993) who showed its equivalence with the entropic version of the PMD algorithm (Algorithm 3.1.2). Teboulle (1992b) also constructed generalized augmented Lagrangians by using either Bregman distances or Csiszár distances (see Chapter 2).

# Part II

# ALGORITHMS

The road to wisdom?—Well it's plain and simple to express:
Err
and err
and err again
but less
and less
and less.

Piet Hein, *The Road to Wisdom*, 1966

CHAPTER 5

# Iterative Methods for Convex Feasibility Problems

---

The *convex feasibility problem* (CFP) is how to find a point in the nonempty intersection $Q \doteq \cap_{i=1}^{m} Q_i \neq \emptyset$ of a finite family of closed convex sets $Q_i \subseteq \mathbb{R}^n$, $i \in I \doteq \{1, 2, \ldots, m\}$ in the $n$-dimensional Euclidean space $\mathbb{R}^n$. This fundamental problem has many applications in mathematics as well as in other fields. When the sets are given in the form

$$Q_i = \{x \in \mathbb{R}^n \mid g_i(x) \leq 0, \ g_i \ \text{ is a convex function}\},$$

we are faced with the problem of solving a system of inequalities with convex functions, of which the linear case is an important and special case.

In this chapter we present several iterative methods for the convex feasibility problem, aiming at a description of both the methods and some of the apparent connections between them. We start with the *method of successive orthogonal projections* (SOP) of Gubin, Polyak, and Raik, also known in recent literature on image recovery as the *method of projections onto convex sets* (POCS). As the names suggest, this method performs orthogonal projections successively onto the individual convex sets $Q_i$ (see Example 2.1.1); therefore, it is a sequential row-action algorithm, according to the classification in Section 1.3.

The other algorithms that we discuss represent two main directions of development. One direction aims at the replacement of the orthogonal projections onto the individual convex sets by other operations that are potentially easier to execute. The second direction aims at creating algorithms that employ either orthogonal projections or other alternative operations in a nonsequential manner. In Section 5.1 we look at control sequences and relaxation parameters. The method of successive orthogonal projections for solving the convex feasibility problem is presented in Section 5.2, and the method of cyclic subgradient projections is presented and studied in Section 5.3. In Section 5.4 several additional methods are briefly discussed with emphasis on their relationship to the cyclic subgradient projections method. The $(\delta, \eta)$-Algorithm, studied in Section 5.5, replaces orthogonal projections onto the convex sets of the convex feasibility problem with projections onto separating hyperplanes. In Section 5.6 the block-iterative version (with variable blocks) of the method of successive orthogonal projections is presented and its convergence analyzed. A block-iterative version of the $(\delta, \eta)$-Algorithm is given in Section 5.7. Successive generalized projections for the convex feasibility problem are discussed in Section 5.8 and a multiprojection simultaneous algorithm, which allows the use of different types of projections within a single application of the algorithm, is the subject of Section 5.9. Section 5.10 introduces a special algorithm, which employs an automatic relaxation strategy for the solution of the linear interval feasibility problem. Notes and references in Section 5.11 conclude this chapter.

## 5.1 Preliminaries: Control Sequences and Relaxation Parameters

Many of the methods studied in this book obey a specific control sequence and employ relaxation parameters. A *control sequence* $\{i(\nu)\}_{\nu=0}^{\infty}$ is a sequence of indices according to which individual sets $Q_i$—or, in other cases, blocks that are groups of sets—may be chosen for the execution of an iterative algorithm. Here are some important control sequences.

**Definition 5.1.1 (Control sequences)**

1. **Cyclic control.** *A control sequence is cyclic if* $i(\nu) = \nu \bmod m + 1$, *where* $m$ *is the total number of sets in the convex feasibility problem.*
2. **Almost cyclic control.** $\{i(\nu)\}_{\nu=0}^{\infty}$ *is almost cyclic on* $I \doteq \{1, 2, \ldots, m\}$ *if* $i(\nu) \in I$, *for all* $\nu \geq 0$, *and there exists an integer* $C \geq m$ *(called the almost cyclicality constant) such that, for all* $\nu \geq 0$, $I \subseteq \{i(\nu + 1), i(\nu + 2), \ldots, i(\nu + C)\}$.
3. **Repetitive control.** *A control sequence* $\{i(\nu)\}_{\nu=0}^{\infty}$ *is called repetitive on* $\{1, 2, \ldots, m\}$ *if* $i(\nu) \in I$ *for all* $\nu \geq 0$, *and for every* $l \in I$ *and every* $\nu \geq 0$ *there exists an index* $\mu > \nu$ *such that* $i(\mu) = l$.
4. **Remotest set control.** *This is obtained by determining* $i(\nu)$ *such that*

$$d(x^\nu, Q_{i(\nu)}) = \max\{d(x^\nu, Q_i) \mid i \in I\}, \tag{5.1}$$

*where* $d(x^\nu, Q_i)$ *is the Euclidean distance between the* $\nu$*th iterate* $x^\nu$ *and the set* $Q_i$, *i.e.,* $d(x^\nu, Q_i) \doteq \min\{\| x^\nu - y \| \mid y \in Q_i\}$.

5. **Approximately remotest set control.** *Denote the maximal distance by* $\max\{d(x^\nu, Q_i) \mid i \in I\} \doteq \theta(x^\nu)$, *and choose the indices* $i(\nu)$ *to satisfy the condition*

$$\lim_{\nu \to \infty} d(x^\nu, Q_{i(\nu)}) = 0 \Rightarrow \lim_{\nu \to \infty} \theta(x^\nu) = 0. \tag{5.2}$$

6. **Most violated constraint control.** *This control, closely related to the remotest set control, is obtained by determining which constraint is most violated by the iterate* $x^\nu$. *If* $Q_i \doteq \{x \in \mathbb{R}^n \mid g_i(x) \leq 0\}$, *for all* $i \in I$, *are the sets under consideration, then one has to determine at each step*
$$g_{i(\nu)}(x^\nu) = \max\{g_i(x^\nu) \mid i \in I\}, \tag{5.3}$$
*and take* $i(\nu)$ *as the control index if* $g_{i(\nu)}(x^\nu) > 0$.

The almost cyclic control means that every sequence of $C$ consecutive indices must include all set indices. Almost cyclic controls are less restrictive than cyclic controls and therefore add an important option as to how the application of a method to a particular problem may be carried out. As will be seen later, this flexibility is important in creating specialized versions of an algorithm and in devising parallel implementations. A cyclic control is almost cyclic with $C = m$.

The repetitive control means that, no matter how far we proceed with the algorithm, every index must appear again (i.e., be repeated).

Obviously, every remotest set control is an approximately remotest set control; but it is also true that every cyclic control is an approximately remotest set control, and these facts make the approximately remotest set control an important theoretical tool.

Several of the methods discussed in this and other chapters employ a sequence $\{\lambda_\nu\}_{\nu=0}^{\infty}$ of *relaxation parameters.* Loosely speaking, these parameters overdo or underdo the move prescribed in an iterative step. Relaxation parameters add an extra degree of freedom to the way a method might actually be implemented and have important consequences for the performance of the method in practice. The effects of relaxation parameters on the iterations—and technical conditions on them that guarantee convergence of an algorithm—are discussed for each method in the sections that follow.

## 5.2 The Method of Successive Orthogonal Projections

We start with a well-known iterative algorithm to solve the convex feasibility problem: the method of *successive orthogonal projections* (SOP) of Gubin, Polyak, and Raik. Starting from an arbitrary point the method generates a sequence $\{x^\nu\}_{\nu=0}^{\infty}$ that converges to a point in $Q$ by performing successive orthogonal projections onto the individual convex sets $Q_i$. The projections can also be relaxed via relaxation parameters. Figure 5.1 demonstrates the algorithm in the case that all relaxation parameters are set to unity.

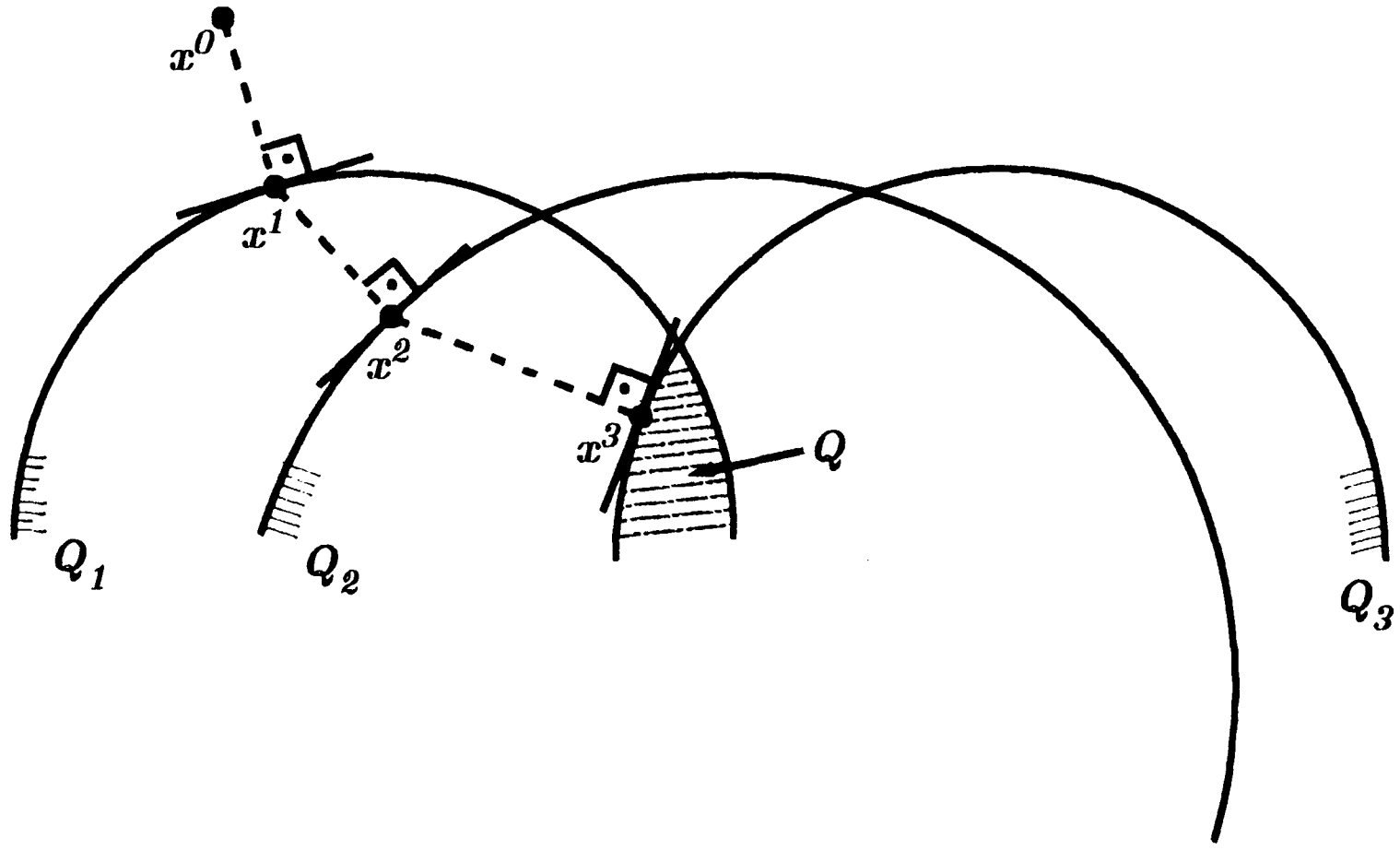

**Figure 5.1** The method of successive orthogonal projection.

### Algorithm 5.2.1 The Method of Successive Orthogonal Projections (SOP)

**Step 0:** (Initialization.) $x^0 \in \mathbb{R}^n$ is arbitrary.

**Step 1:** (Iterative step.) For all $\nu \geq 0$,

$$x^{\nu+1} = x^\nu + \lambda_\nu (P_{Q_{i(\nu)}}(x^\nu) - x^\nu), \tag{5.4}$$

where $P_{Q_{i(\nu)}}(x^\nu)$ stands for the orthogonal projection of $x^\nu$ onto the set $Q_{i(\nu)}$; see (5.5).

**Relaxation parameters:** $\{\lambda_\nu\}_{\nu=0}^\infty$ are confined to the interval $\epsilon_1 \leq \lambda_\nu \leq 2 - \epsilon_2$ for all $\nu \geq 0$ with some arbitrarily small $\epsilon_1, \epsilon_2 > 0$.

**Control:** The sequence $\{i(\nu)\}_{\nu=0}^\infty$ is almost cyclic on $I$.

Figure 5.2 shows the effect of relaxation on the iterative step of Algorithm 5.2.1. If $\lambda_\nu = 1$ then $x^{\nu+1} = P_{Q_{i(\nu)}}(x^\nu)$ and the step is unrelaxed. If $\eta \leq \lambda_\nu < 1$ then $x^{\nu+1}$ is an underrelaxed projection of $x^\nu$ onto $Q_{i(\nu)}$, and if $1 < \lambda_\nu \leq 2 - \eta$ then $x^{\nu+1}$ is an overrelaxed projection of $x^\nu$ onto $Q_{i(\nu)}$.

The SOP method is particularly useful when the projections onto the individual sets are easily calculated, such as when the sets are hyperplanes or half-spaces. In general however, application of this method requires, at each iterative step, the solution of a subsidiary minimization problem associated with the projection onto the selected set, namely,

$$P_{Q_{i(\nu)}}(x^\nu) = \operatorname*{argmin}_{y \in Q_{i(\nu)}} \| x^\nu - y \|, \tag{5.5}$$

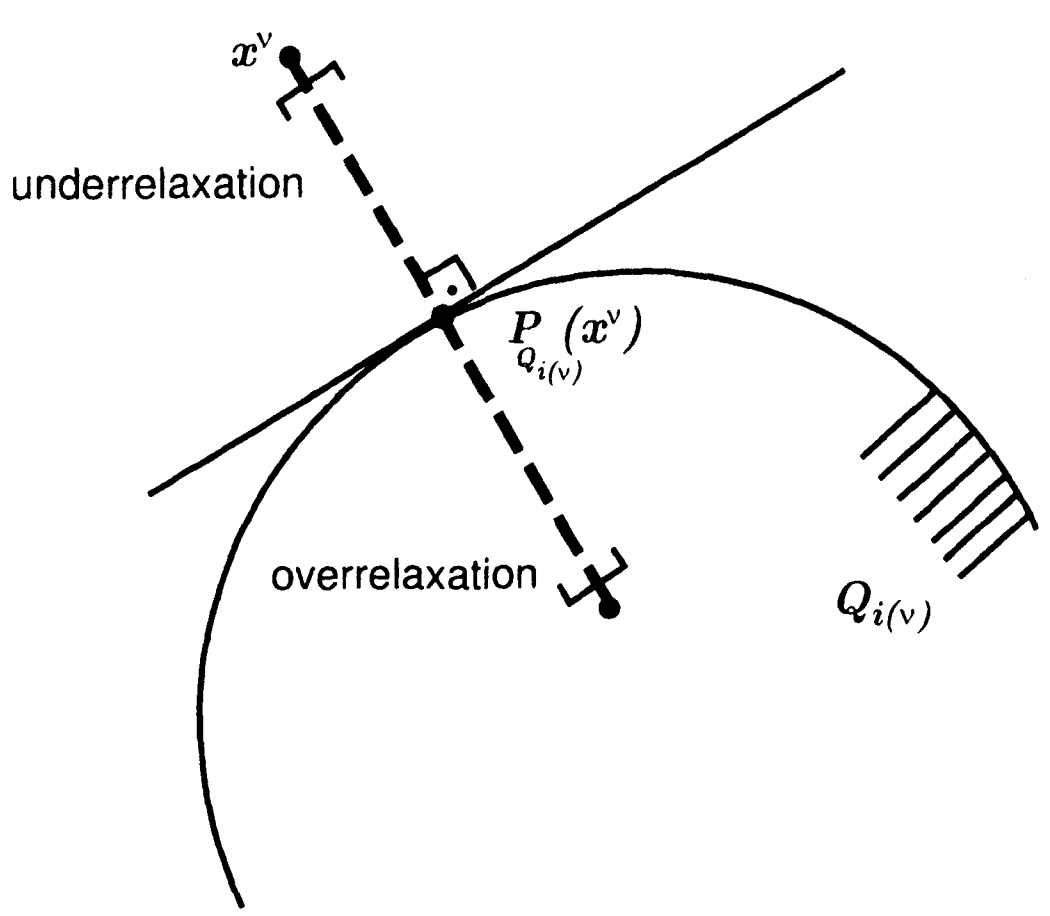

**Figure 5.2** The effect of relaxation on an iterative step of Algorithm 5.2.1.

where $\|\cdot\|$ stands for the Euclidean norm. An instance of simple projections is the linear case, discussed later (see section 5.4.4). There, $Q_i$ are half-spaces in $\mathbb{R}^n$ and the SOP method coincides with another well-known algorithm, the relaxation method of Agmon, and Motzkin and Schoenberg. We will not give a proof of the convergence of the SOP method (Algorithm 5.2.1) since it follows as a special case from the convergence proofs of more general schemes that we study later.

For the general convex feasibility problem with convex sets of the form $Q_i = \{x \in \mathbb{R}^n \mid g_i(x) \leq 0\}$, with $g_i : \mathbb{R}^n \to \mathbb{R}$ convex differentiable functions for all $i \in I$, the iterative step of the SOP method, in fact, requires that a move be made in a direction determined by $\nabla g_{i(\nu)}(x^{\nu+1})$, the gradient of $g_{i(\nu)}$ calculated at the next, and not yet known, iterate $x^{\nu+1}$. This happens because it is precisely the direction necessary for a move toward the projection point $P_{Q_{i(\nu)}}(x^\nu)$. A method that circumvents this difficulty is the cyclic subgradient projections method, which we study next.

## 5.3 The Cyclic Subgradient Projections Method

The cyclic subgradient projections method (CSP) differs from the SOP method in that it requires at each iterative step a move in the direction of the gradient (or subgradient) calculated at the current iterate, which is already known, i.e., $\nabla g_{i(\nu)}(x^\nu)$ and not $\nabla g_{i(\nu)}(x^{\nu+1})$. In this way some of the computational difficulties associated with the SOP method are circumvented because the need to solve (5.5) is replaced by gradient (or subgradient) calculations at an available point, i.e., $\nabla g_{i(\nu)}(x^\nu)$. There is further discussion on this point in Section 5.4.1.

In order to describe the CSP method we recall that a vector $t \in \mathbb{R}^n$ is said to be a *subgradient* of a convex function $g$ at a point $y$, if

$$\langle t, x - y \rangle \leq g(x) - g(y), \text{ for every } x \in \mathbb{R}^n. \tag{5.6}$$

This inequality is referred to as the *subgradient inequality.* If $g$ is differentiable at $y$ then its gradient $\nabla g(y)$ is the unique subgradient of $g$ at $y$. A convex function always has a subgradient and the set of all subgradients of $g$ at $y$ is denoted by $\partial g(y)$ (see Section 5.11 for a reference on these).

**Algorithm 5.3.1 The Method of Cyclic Subgradient Projections (CSP)**

**Step 0:** (Initialization.) $x^0 \in \mathbb{R}^n$ is arbitrary.
**Step 1:** (Iterative step.)

$$x^{\nu+1} = \begin{cases} x^\nu - \lambda_\nu \dfrac{g_{i(\nu)}(x^\nu)}{\| t^\nu \|^2} t^\nu, & \text{if } g_{i(\nu)}(x^\nu) > 0, \\ x^\nu, & \text{if } g_{i(\nu)}(x^\nu) \leq 0, \end{cases} \tag{5.7}$$

where $t^\nu \in \partial g_{i(\nu)}(x^\nu)$ is a subgradient of $g_{i(\nu)}$ at the point $x^\nu$.

**Relaxation Parameters:** $\{\lambda_\nu\}_{\nu=0}^{\infty}$ are confined to the interval $\epsilon_1 \leq \lambda_\nu \leq 2 - \epsilon_2$, for all $\nu \geq 0$, with some, arbitrarily small, $\epsilon_1, \epsilon_2 > 0$.

**Control:** The sequence $\{i(\nu)\}_{\nu=0}^{\infty}$ is almost cyclic on $I$.

Observe that if $t^\nu = 0$, then $g_{i(\nu)}$ takes its minimal value at $x^\nu$, implying by the nonemptiness of $Q$ that $g_{i(\nu)}(x^\nu) \leq 0$, so that $x^{\nu+1} = x^\nu$.

The relationship of the CSP method to other iterative methods for solving the convex feasibility problem and to the relaxation method for solving linear inequalities will be clarified later. We turn now to the convergence theorem for the CSP method, using the following definition.

**Definition 5.3.1 (Fejér-monotonicity)** *A sequence* $\{x^\nu\}_{\nu=0}^{\infty}$ *is Fejér-monotone with respect to some fixed set* $Q \subseteq \mathbb{R}^n$ *if, for every* $x \in Q$,

$$\| x^{\nu+1} - x \| \leq \| x^\nu - x \|, \textit{ for all } \nu \geq 0. \tag{5.8}$$

It is easy to see that every Fejér-monotone sequence is bounded because, for a fixed $x \in Q$, $\| x^{\nu+1} - x \| \leq \| x^\nu - x \| \leq \ldots \leq \| x^0 - x \|$, and applying the triangle inequality gives $\| x^{\nu+1} \| - \| x \| \leq \| x^{\nu+1} - x \| \leq \| x^0 - x \|$; thus, $\| x^{\nu+1} \| \leq \| x^0 - x \| + \| x \| = M$, for all $\nu \geq 0$, and $M$ is constant since $x$ is given.

**Theorem 5.3.1** *Let the following assumptions hold:*

(*i*) *the functions* $g_i(x)$ *are continuous and convex on* $\mathbb{R}^n$, *for all* $i \in I$,

*(ii)* $Q \doteq \cap_{i\in I} Q_i \neq \emptyset$, *and*

*(iii) for some* $\hat{x} \in Q$ *there exists a constant* $K = K(\hat{x})$ *such that* $\| t \| \leq K$, *for all subgradients* $t \in \partial g_i(x)$, *for all* $i \in I$, *and for all* $x \in \mathbb{R}^n$ *for which* $\| x - \hat{x} \| \leq \| x^0 - \hat{x} \|$. *This assumption is called uniform boundedness of the subgradients.*

*Then any sequence* $\{x^\nu\}$ *produced by Algorithm 5.3.1 converges to a solution of the convex feasibility problem, i.e.,* $x^\nu \to x^* \in Q$, *as* $\nu \to \infty$.

**Proof** The proof consists of the following five steps:
Step 1. $\{x^\nu\}$ is a Fejér-monotone sequence with respect to $Q$.
Step 2. $\lim_{\nu\to\infty} g_{i(\nu)}(x^\nu) = 0$.
Step 3. $\lim_{\nu\to\infty} \| x^{\nu+1} - x^\nu \| = 0$.
Step 4. $\lim_{\nu\to\infty} g_i(x^\nu) = 0$, for every *fixed* $i \in I$ .
Step 5. $\lim_{\nu\to\infty} x^\nu = x^* \in Q$.
We now proceed with the proof of each step.

**Step 1.** Take some $x \in Q$. Abbreviating here $i \doteq i(\nu)$ and denoting $\alpha_\nu \doteq \lambda_\nu \dfrac{g_i(x^\nu)}{\| t^\nu \|^2}$, we have, for an active step—i.e., when a change is made to the current iterate when $g_i(x^\nu) > 0$—in (5.7)

$$\begin{aligned} \| x^{\nu+1} - x \|^2 &= \| x^\nu - \alpha_\nu t^\nu - x \|^2 \\ &\leq \| x^\nu - x \|^2 + \alpha_\nu^2 \| t^\nu \|^2 - 2\alpha_\nu \langle t^\nu, x^\nu - x\rangle. \end{aligned} \tag{5.9}$$

From the subgradient inequality (5.6) and by $g_i(x) \leq 0$, (since $x \in Q$) we obtain

$$\begin{aligned} \| x^{\nu+1} - x \|^2 &\leq \| x^\nu - x \|^2 + \alpha_\nu^2 \| t^\nu \|^2 - 2\alpha_\nu g_i(x^\nu) \\ &= \| x^\nu - x \|^2 + (\lambda_\nu^2 - 2\lambda_\nu) \frac{(g_i(x^\nu))^2}{\| t^\nu \|^2}. \end{aligned} \tag{5.10}$$

From the fact that $\epsilon_1 \leq \lambda_\nu \leq 2 - \epsilon_2$, for all $\nu \geq 0$, we get

$$\| x^{\nu+1} - x \|^2 \leq \| x^\nu - x \|^2 - \epsilon_1\epsilon_2 \frac{(g_i(x^\nu))^2}{\| t^\nu \|^2}, \tag{5.11}$$

from which Fejér-monotonicity follows.

**Step 2.** For $x \in Q$ the sequence $\{\| x^\nu - x \|\}$ is monotonically decreasing and bounded below, therefore, $\lim_{\nu\to\infty} \| x^\nu - x \| = d$, for example. This implies at once, via (5.11), that

$$\lim_{\nu\to\infty} \frac{(g_{i(\nu)}(x^\nu))^2}{\| t^\nu \|^2} = 0. \tag{5.12}$$

Now we use the assumption on uniform boundedness of the subgradients. Let $\hat{x} \in Q$ be the point whose existence is assumed in the theorem.

If we denote

$$T_{\hat{x}} \doteq \{x \in \mathbb{R}^n \mid \| x - \hat{x} \| \leq \| x^0 - \hat{x} \|\}, \tag{5.13}$$

then $x^\nu \in T_{\hat{x}}$ for every $\nu \geq 0$, because $\hat{x} \in Q$, $x^\nu \in \mathbb{R}^n$, and by repeated application of Fejér-monotonicity. Assumption (*iii*) of the Theorem now ensures that $\| t^\nu \| \leq K$, and combining this with (5.12) we obtain

$$\lim_{\nu \to \infty} g_{i(\nu)}(x^\nu) = 0. \tag{5.14}$$

**Step 3.** By substitution from (5.7) and use of (5.12), we get

$$\| x^{\nu+1} - x^\nu \|^2 \leq \lambda_\nu^2 \frac{(g_{i(\nu)}(x^\nu))^2}{\| t^\nu \|^2} \tag{5.15}$$

and that the right-hand side of it tends to zero as $\nu \to \infty$. Note that this implies also that

$$\lim_{\nu \to \infty} \| x^{\nu+j} - x^\nu \| = 0, \tag{5.16}$$

for every integer $j$.

**Step 4.** Let $i \in I$ be a *fixed* index. Then

$$| g_i(x^\nu) | \leq | g_i(x^\nu) - g_i(x^\mu) | + | g_i(x^\mu) |, \tag{5.17}$$

where $\mu$ is chosen to be the integer closest to $\nu$ such that, for some $r$

$$\mu = i + rm. \tag{5.18}$$

For all $\nu \geq 0$, we then have $| \nu - \mu | \leq C$, with $C$ the constant of almost cyclic control (see Definition 5.1.1).

The set $T_{\hat{x}}$ defined by (5.13) is compact; therefore $g_i(x)$ is uniformly continuous on it. Thus (5.16) implies that $\lim_{\nu \to \infty} | g_i(x^\nu) - g_i(x^\mu) | = 0$. The choice (5.18) of $\mu$ means that for the $i$ that we fixed above $i = i(\mu)$; therefore (5.14) implies that $\lim_{\nu \to \infty} | g_i(x^\mu) | = 0$. Thus, (5.17) gives the required result that

$$\lim_{\nu \to \infty} g_i(x^\nu) = 0, \text{ for every fixed } i \in I. \tag{5.19}$$

**Step 5.** Fejér-monotonicity of $\{x^\nu\}$, proven in Step 1, implies boundedness. Therefore, $x^\nu$ must have a convergent subsequence, i.e.,

$$\lim_{s \to \infty} x^{\nu_s} = x^*. \tag{5.20}$$

From (5.19) and the continuity of $g_i$ we know that $x^* \in Q$. In Step 2 we showed that $\lim_{\nu \to \infty} \| x^\nu - x^* \| = d$, but now we have the additional information that $\lim_{s \to \infty} \| x^{\nu_s} - x^* \| = 0$; Thus, $\lim_{\nu \to \infty} \| x^\nu - x^* \| = 0$, and the proof is complete. ∎

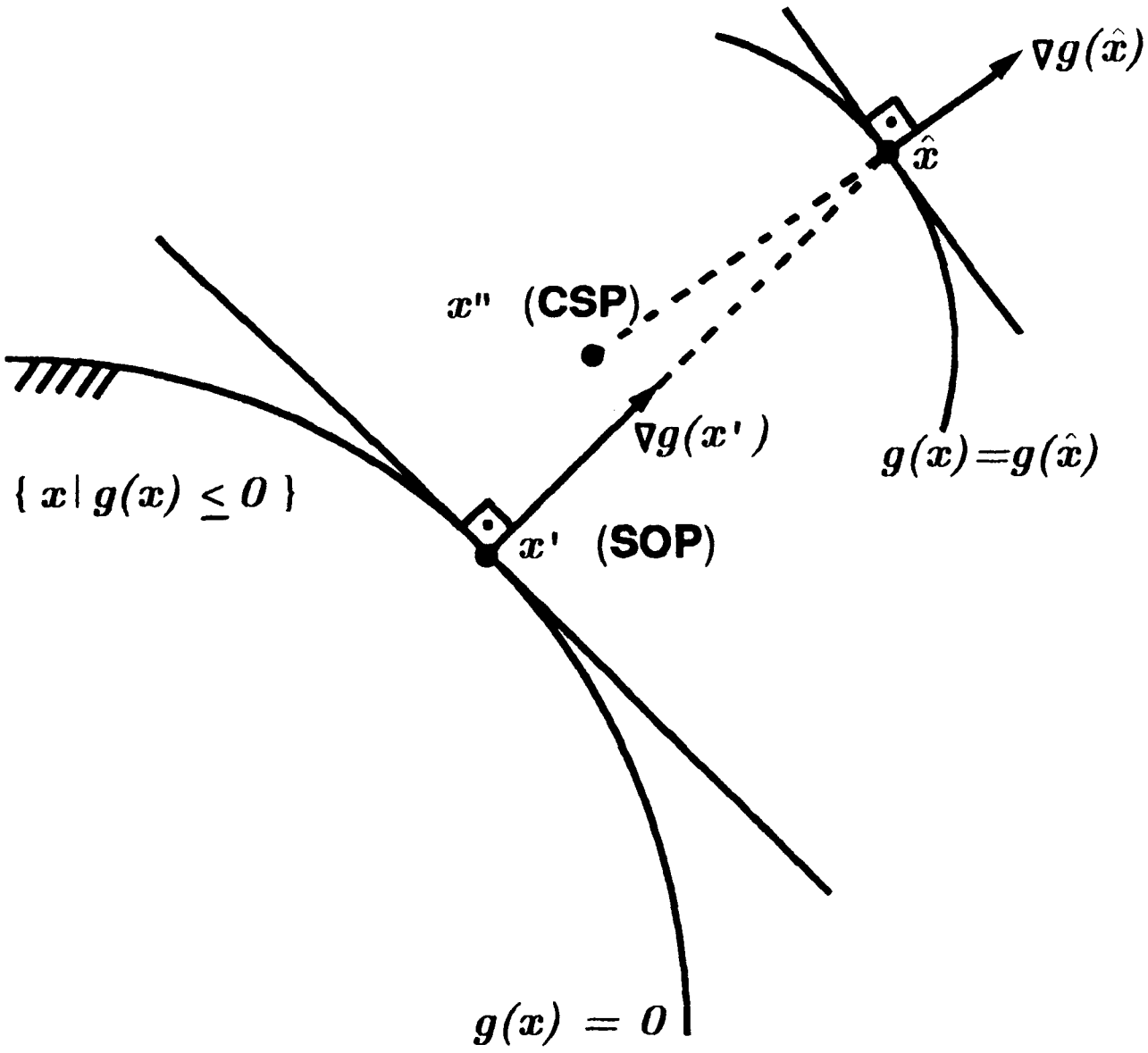

**Figure 5.3** Cyclic Subgradient Projections (CSP) and Successive Orthogonal Projections (SOP) – Comparison of (unrelaxed) moves in an iterative step.

## 5.4 The Relationship of CSP to Other Methods

In this section we briefly discuss several additional methods for solving the convex feasibility problem and describe their relationships to the CSP method.

### 5.4.1 The method of successive orthogonal projections

The SOP method (Algorithm 5.2.1) is particularly useful when the projections onto the individual sets are easily calculated. In general however, application of the SOP method would require at each step the solution of an auxiliary minimization problem associated with the projection onto the current set.

Figure 5.3 illustrates the basic difference between SOP and CSP. Performing an orthogonal projection of a point $\hat{x}$ onto a convex set defined by $\{x \in \mathbb{R}^n \mid g(x) \leq 0\}$ amounts to making a move in the direction of the negative gradient of $g$ at the point $x'$, which is the orthogonal projection of $\hat{x}$ onto the set. But the point $x'$ is not yet known, so the direction of the move cannot be calculated via gradient calculation, but rather through the minimization of $\| \hat{x} - x \|$, over all $x$ for which $g(x) \leq 0$, thus minimizing the distance between $\hat{x}$ and the set.

In contrast to this, in order to find the next iterate $x''$, the CSP method makes a move from $\hat{x}$ in the direction of the negative gradient of $g$ at the point $\hat{x}$ itself. Such a move might be less efficient than an orthogonal projection in terms of the speed of the initial convergence of the process, but it enables us to replace the distance minimization problems by simpler gradient or subgradient calculations.

### 5.4.2 A remotest set controlled subgradient projections method

A control sequence $\{i(\nu)\}_{\nu=0}^{\infty}$ for the application of an iterative method to the convex feasibility problem (see Definition 5.1.1) is a *remotest set control* if $i(\nu)$ is determined at each step by

$$g_{i(\nu)}(x^{\nu}) = \max\{g_i(x^{\nu}) \mid 1 \leq i \leq m\}. \tag{5.21}$$

A subgradient projections method for the convex feasibility problem with the remotest set control has the same form as the CSP method (Algorithm 5.3.1) except that the control is replaced by the remotest set control. It is, however, in the introduction of the cyclic (or almost cyclic) control that the advantage of the CSP method lies. Having a cyclic method eliminates the need to determine the remotest set, which in turn requires the performance of orthogonal projections onto all sets at each iterative step. In a real-world problem, which might involve a large number of functions, this flexibility is of great practical significance.

### 5.4.3 The scheme of Oettli

The method of Oettli is an algorithmic scheme for the solution of the convex feasibility problem with sets of the form $Q_i = \{x \in \mathbb{R}^n \mid g_i(x) \leq 0\}$, $i \in I \doteq \{1, 2, \ldots, m\}$. Again $Q = \cap_{i \in I} Q_i \neq \emptyset$. To describe this method we need the notion of a monotonic norm. A norm $p(\cdot)$ on $\mathbb{R}^m$ is *monotonic* if, for any $x, y \in \mathbb{R}^m$, for which the absolute values of their components obey $\mid x_i \mid \leq \mid y_i \mid$, for all $i \in I$, we have $p(x) \leq p(y)$. It is easy to see that this is equivalent to the condition that $p(x) = p(\mid x \mid)$, for all $x \in \mathbb{R}^m$, where $\mid x \mid$ denotes the vector obtained by taking the absolute values of all components of $x$, i.e., $\mid x \mid = (\mid x_j \mid)_{j=1}^n$. In particular, every $\ell_p$-norm, $\| x \|_p \doteq (\sum_{j=1}^n \mid x_j \mid^p)^{1/p}$, for $1 \leq p \leq \infty$, is monotonic. An example of a nonmonotonic norm in $\mathbb{R}^2$ is $p(x) \doteq \max(2 \mid x_1 + x_2 \mid, \mid x_1 - x_2 \mid)$.

**Algorithm 5.4.1 Oettli's Scheme**

**Step 0:** (Initialization.) $x^0 \in \mathbb{R}^n$ is arbitrary.
**Step 1:** (Iterative step.)

$$x^{\nu+1} = x^{\nu} - \lambda_{\nu} \frac{\varphi(x^{\nu})}{\| t^{\nu} \|^2} t^{\nu}, \quad \text{if } x^{\nu} \notin Q. \tag{5.22}$$

The function $\varphi : \mathbb{R}^n \to \mathbb{R}$ is defined by

$$\varphi(x) \doteq p(g_1^+(x), g_2^+(x), \ldots, g_m^+(x)), \tag{5.23}$$

where $g_i^+(x) \doteq \max(0, g_i(x))$, and $p$: $\mathbb{R}^m \to \mathbb{R}$ is a monotonic norm. $t^\nu \in \partial\varphi(x^\nu)$ is any subgradient of $\varphi$ at $x^\nu$.

**Relaxation parameters:** $\lambda_\nu \in (0,2)$ are such that $\sum_{\nu=0}^{\infty} \lambda_\nu(2-\lambda_\nu) = +\infty$.

**Control:** The control of the method is determined by the choice of the monotonic norm $p$.

Different choices of the monotonic norm $p$ in Oettli's scheme give rise to various concrete algorithms. Choosing $p$ to be the $\ell_\infty$-norm in $\mathbb{R}^m$, $\| x \|_\infty \doteq \max_{1 \le j \le m} | x_j |$ yields the remotest set controlled method of subgradient projections for solving the convex feasibility problem (see Section 5.4.2). However, we are unaware of any monotonic norm $p$ that, upon substitution into Oettli's scheme, yields the CSP method (Algorithm 5.3.1), i.e., with cyclic control.

### 5.4.4 The linear feasibility problem: solving linear inequalities

When all the functions $g_i$ are linear, i.e., $g_i(x) \doteq \langle a^i, x\rangle - b_i$ for all indices $i \in I \doteq \{1, 2, \ldots, m\}$, and $a^i \in \mathbb{R}^n$ and $b_i \in \mathbb{R}$ are given, we have the problem of solving a system of linear inequalities. Both the cyclically controlled SOP method (Algorithm 5.2.1) and the CSP method (Algorithm 5.3.1) then coincide and reduce to the cyclically controlled relaxation method of Agmon, Motzkin, and Schoenberg (AMS).

**Algorithm 5.4.2 The Relaxation Method of Agmon, Motzkin, and Schoenberg (AMS)**

**Step 0:** (Initialization.) $x^0 \in \mathbb{R}^n$ is arbitrary.

**Step 1:** (Iterative step.)

$$x^{\nu+1} = x^\nu + c_\nu a^{i(\nu)}, \tag{5.24}$$

where $c_\nu \doteq \min\left(0, \lambda_\nu \dfrac{b_{i(\nu)} - \langle a^{i(\nu)}, x^{(\nu)}\rangle}{\| a^{i(\nu)} \|^2}\right)$.

**Relaxation parameters:** $\lambda_\nu$ are such that $\epsilon_1 \le \lambda_\nu \le 2 - \epsilon_2$, for all $\nu \ge 0$, with some, arbitrarily small, $\epsilon_1, \epsilon_2 > 0$.

**Control:** The sequence $\{i(\nu)\}_{\nu=0}^{\infty}$ is almost cyclic on $I$.

This algorithm proceeds in the following manner. At any iterative step, the iterate $x^\nu$ and the closed half-space determined by the $i(\nu)$th inequality $Q_{i(\nu)} \doteq \{x \in \mathbb{R}^n \mid \langle a^{i(\nu)}, x\rangle \le b_{i(\nu)}\}$ are at hand. If $x^\nu \in Q_{i(\nu)}$, then $x^{\nu+1} = x^\nu$, but if the current iterate $x^\nu$ violates the $i(\nu)$th inequality, i.e., $x^\nu \notin Q_{i(\nu)}$ then $x^{\nu+1}$ lies on the line through $x^\nu$, perpendicular to $H_{i(\nu)}$, the bounding hyperplane of $Q_{i(\nu)}$. For unity relaxation, $x^{\nu+1}$ is the orthogonal projection of $x^\nu$ onto $Q_{i(\nu)}$ (see Figure 5.4). The relaxation option actually allows the next iterate $x^{\nu+1}$ to be a point inside the line segment connecting $x^\nu$ and its orthogonal reflection with respect to $H_{i(\nu)}$.

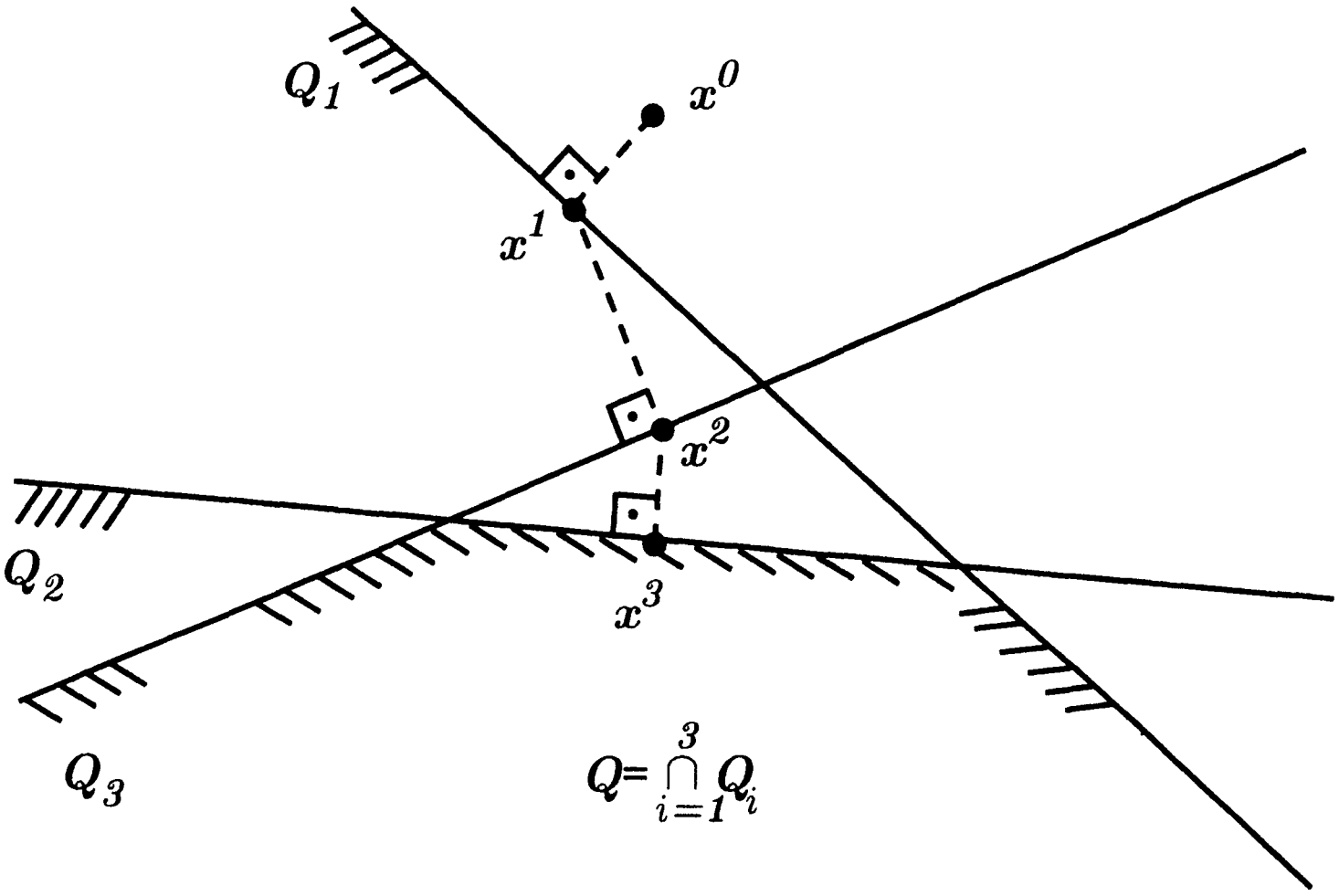

**Figure 5.4** The relaxation method of Agmon, Motzkin, and Schoenberg (AMS) for solving a system of linear inequalities, illustrated for $n = 2$ and $m = 3$.

### 5.4.5 Kaczmarz's algorithm for systems of linear equations and its nonlinear extension

Kaczmarz's algorithm for the solution of the system $\langle a^i, x\rangle - b_i = 0$, for $i \in I$ where $a^i \in \mathbb{R}^n$ and $b_i \in \mathbb{R}$ are given, may be derived from the AMS relaxation method (Algorithm 5.4.2) by replacing the equality signs with pairs of opposite inequality signs, or derived directly from the SOP method (Algorithm 5.2.1). Kaczmarz's algorithm was independently discovered in the field of image reconstruction from projections (see Chapter 10) where it was called Algebraic Reconstruction Technique (ART).

**Algorithm 5.4.3 Kaczmarz's Algorithm, Algebraic Reconstruction Technique (ART)**

**Step 0:** (Initialization.) $x^0 \in \mathbb{R}^n$ is arbitrary.

**Step 1:** (Iterative step.)

$$x^{\nu+1} = x^\nu + \lambda_\nu \frac{b_{i(\nu)} - \langle a^{i(\nu)}, x^\nu\rangle}{\| a^{i(\nu)} \|^2} a^{i(\nu)}, \tag{5.25}$$

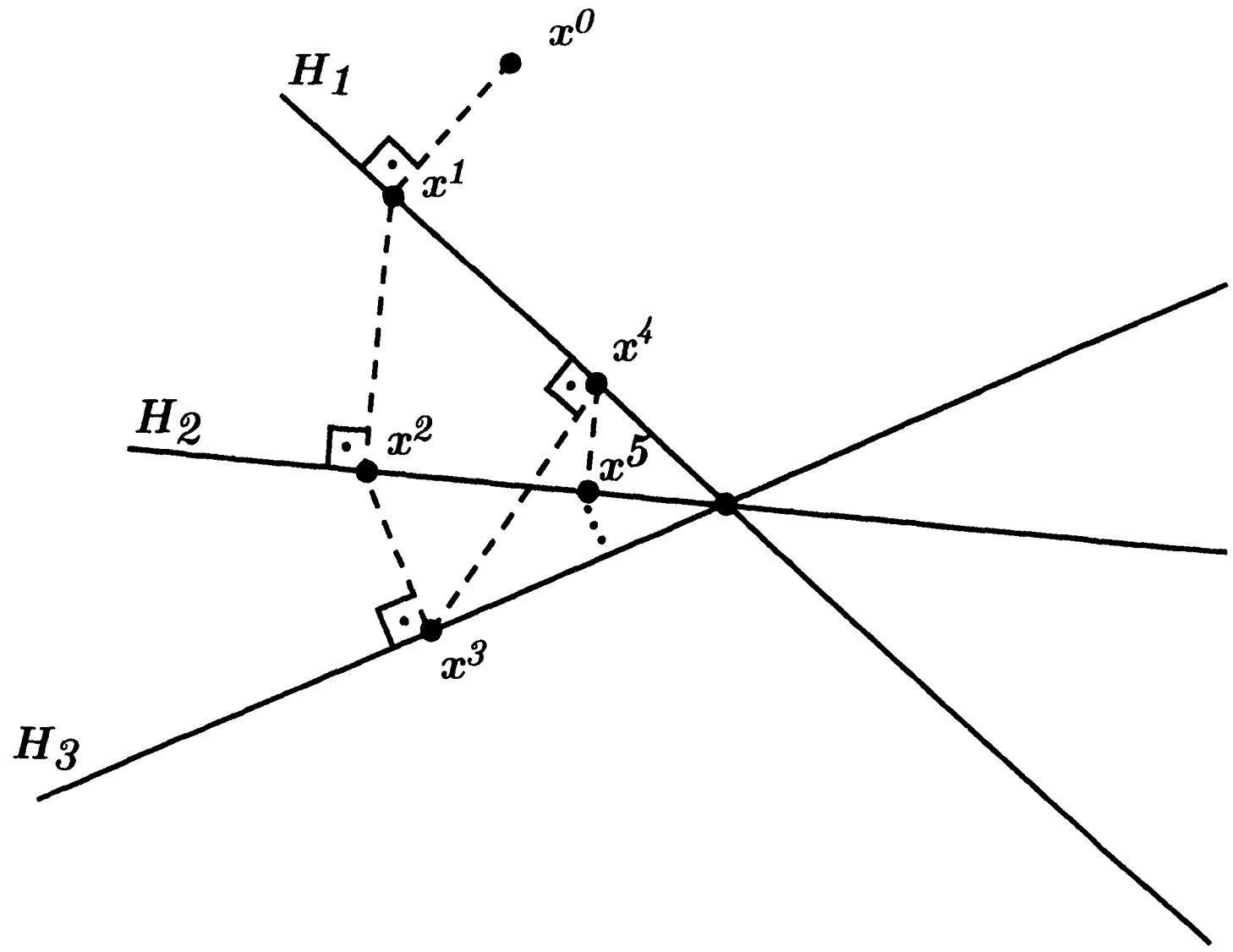

**Figure 5.5** Kaczmarz' method (with unity relaxation) for solving a system of linear equations, illustrated for $n = 2$ and $m = 3$.

**Relaxation parameters:** For all $\nu \geq 0$, $\lambda_\nu$ are confined to $\epsilon_1 \leq \lambda_\nu \leq 2 - \epsilon_2$, for any, arbitrarily small, $\epsilon_1, \epsilon_2 > 0$.

**Control:** The sequence $\{i(\nu)\}_{\nu=0}^{\infty}$ is almost cyclic on $I$.

We can now interpret the algorithm according to (5.25). Given the current iterate $x^\nu$ and the hyperplane $H_{i(\nu)} = \{x \in \mathbb{R}^n \mid \langle a^{i(\nu)}, x \rangle = b_{i(\nu)}\}$, determined by the $i(\nu)$th equation, the new iterate $x^{\nu+1}$ lies on the line through $x^\nu$, which is perpendicular to $H_{i(\nu)}$. For unity relaxation, $\lambda_\nu = 1$, $x^{\nu+1}$ is the orthogonal projection of $x^\nu$ onto $H_{i(\nu)}$ (see Figure 5.5). The relaxation option again allows the next iterate $x^{\nu+1}$ to be a point inside the line segment connecting $x^\nu$ and its orthogonal reflection with respect to $H_{i(\nu)}$.

Kaczmarz's algorithm was extended to the nonlinear case for various controls, including the version for a cyclic control given next. Consider the system of nonlinear equations $f_i(x) = 0$, $i \in I = \{1, 2, \ldots, m\}$. For each $i \in I$, $f_i : \mathbb{R}^n \to \mathbb{R}$ is a (generally) nonlinear function whose first partial derivatives exist and are continuous for all $x$. $\nabla f_{i(\nu)}(x^\nu)$, the gradient of $f_{i(\nu)}$ at $x^\nu$, replaces the normal vector $a^{i(\nu)}$ in (5.25).

### Algorithm 5.4.4 Nonlinear Extension of Kaczmarz's Algorithm

**Step 0:** (Initialization.) $x^0 \in \mathbb{R}^n$ is arbitrary.

**Step 1:** (Iterative step.)

$$x^{\nu+1} = x^{\nu} - \frac{f_{i(\nu)}(x^{\nu})}{\| \nabla f_{i(\nu)}(x^{\nu}) \|^2} \nabla f_{i(\nu)}(x^{\nu}). \tag{5.26}$$

**Relaxation parameters:** Unity, i.e., no relaxation is incorporated.

**Control:** The sequence $\{i(\nu)\}_{\nu=0}^{\infty}$ is cyclic over $I$.

This method converges locally to a solution of the system. The conditions for such convergence include second-order differentiability but no convexity assumptions. Relaxation parameters are not incorporated (see Section 5.11 for references).

The close relationship of Algorithm 5.4.4 to the CSP method (Algorithm 5.3.1) is obvious, but there seems to be no obvious way to derive this algorithm from the CSP method. Assuming convexity of all $f_i$'s and replacing the equality signs in the system $f_i(x) = 0$, $i \in I$, with pairs of inequalities leads to inequalities of two types: convex inequalities, i.e., $f_i(x) \leq 0$, with convex $f_i$'s, which represent convex sets in $\mathbb{R}^n$; and *complementary-convex* inequalities, i.e., $-f_i(x) \leq 0$, $f_i$ convex, which represent sets in $\mathbb{R}^n$ that are complements of convex sets. Therefore, the CSP method as presented in Algorithm 5.3.1 cannot be applied.

## 5.5 The $(\delta, \eta)$-Algorithm

To study the next class of algorithms we must first define the concepts, *separating hyperplane* and *supporting hyperplane.*

**Definition 5.5.1 (Separating Hyperplane)** *Let $C_1$ and $C_2$ be two nonempty sets in $\mathbb{R}^n$. A hyperplane $H \doteq \{x \in \mathbb{R}^n \mid \langle a, x \rangle = b\}$ separates $C_1$ and $C_2$ if $\langle a, x \rangle \geq b$ for all $x \in C_1$, and $\langle a, x \rangle \leq b$ for all $x \in C_2$. If in addition, $C_1 \cup C_2 \not\subseteq H$, then $H$ is said to properly separate the two sets.*

**Definition 5.5.2 (Supporting Hyperplane)** *Let $C$ be a nonempty set in $\mathbb{R}^n$ and let $\overline{x} \in \operatorname{bd} C$, i.e., $\overline{x}$ belongs to the boundary $\operatorname{bd} C$ of $C$. A hyperplane $H \doteq \{x \in \mathbb{R}^n \mid \langle a, x - \overline{x} \rangle = 0\}$ supports $C$ at $\overline{x}$ if either $\langle a, x - \overline{x} \rangle \geq 0$ for all $x \in C$, or else $\langle a, x - \overline{x} \rangle \leq 0$ for all $x \in C$. If in addition, $C \not\subseteq H$, then $H$ is a proper supporting hyperplane of $C$ at $\overline{x}$.*

An orthogonal projection of a point $x$ onto a set $C$ can be viewed as an orthogonal projection of $x$ onto the particular hyperplane $H$, which separates $x$ from $C$ and supports $C$ at $x'$, the closest point to $x$ in $C$ (see Figure 5.6). But, of course, when executing such an orthogonal projection neither the point $x'$ nor the separating and supporting hyperplane $H$ are available. In view of the simplicity of an orthogonal projection onto a hyperplane, it is natural to ask whether one could use other separating supporting hyperplanes such as $\hat{H}$ in Figure 5.6, instead of that particular hyperplane $H$ through the closest point to $x$. Aside from theoretical interest, this approach leads to algorithms useful in practice, provided that the computational efforts of finding such other hyperplanes compete favorably with the work involved in performing orthogonal projections directly onto the given sets.

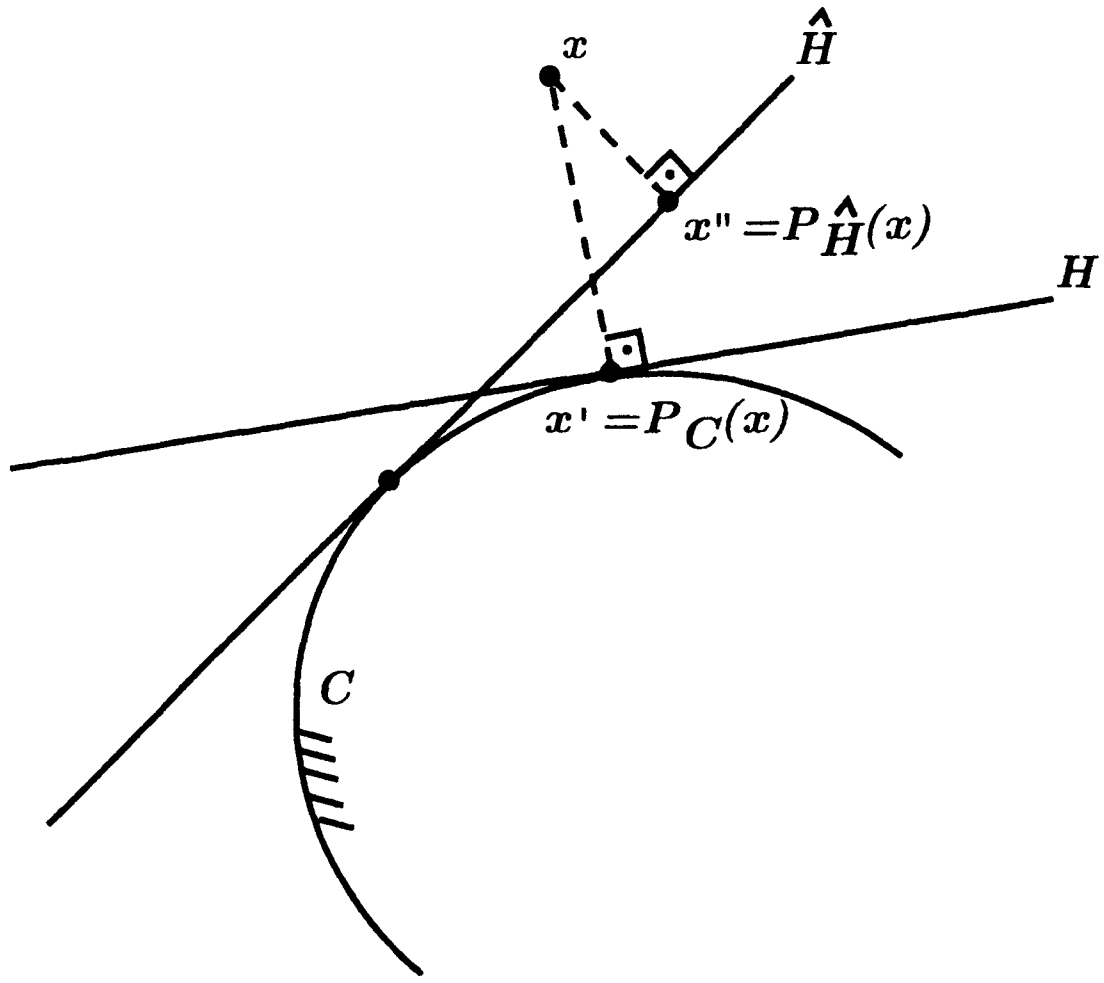

**Figure 5.6** Can we replace the orthogonal projection $x'$ of $x$ onto $C$ by $x''$ obtained by projecting $x$ onto another separating supporting hyperplane?

With these thoughts in mind, a general framework for the design of algorithms for the convex feasibility problem is given now by the so-called $(\delta, \eta)$-Algorithm. The description of the $(\delta, \eta)$-Algorithm involves the construction of a set $A_C(x, \lambda)$, used to describe the iterative step, and the definition of a control sequence through which individual sets are chosen. With a given point $x \in \mathbb{R}^n$ and a given closed and convex set $C \subseteq \mathbb{R}^n$, we associate a set $A_C(x, \lambda)$ in the following manner. Choose some $\delta$, $0 \leq \delta \leq 1$, and some $\eta$, $0 \leq \eta \leq 1$, and let $B \doteq B(x, \delta d(x, C))$ be the ball centered at $x$ with radius $\delta d(x, C)$, where $d(x, C)$, the Euclidean distance between $x$ and $C$, is defined by $d(x, C) \doteq \min\{\| x - z \| \mid z \in C\}$. For $x \notin C$ denote by $\mathcal{H}_{x,C}$ the family of all hyperplanes that separate $B$ from $C$. With this notation define:

$$A_C(x, \lambda) \doteq \begin{cases} \{x\}, & \text{if } x \in C, \\ \{x + \lambda(P_H(x) - x) \mid H \in \mathcal{H}_{x,C},\ \eta \leq \lambda \leq 2 - \eta\}, & \text{if } x \notin C. \end{cases} \tag{5.27}$$

This means that for a point $x$ outside $C$, the set $A_C(x, \lambda)$ consists of all, possibly relaxed by $\lambda$, projections of $x$ onto *any* hyperplane that separates the set $C$ from a ball with radius $\delta d(x, C)$ centered at $x$.

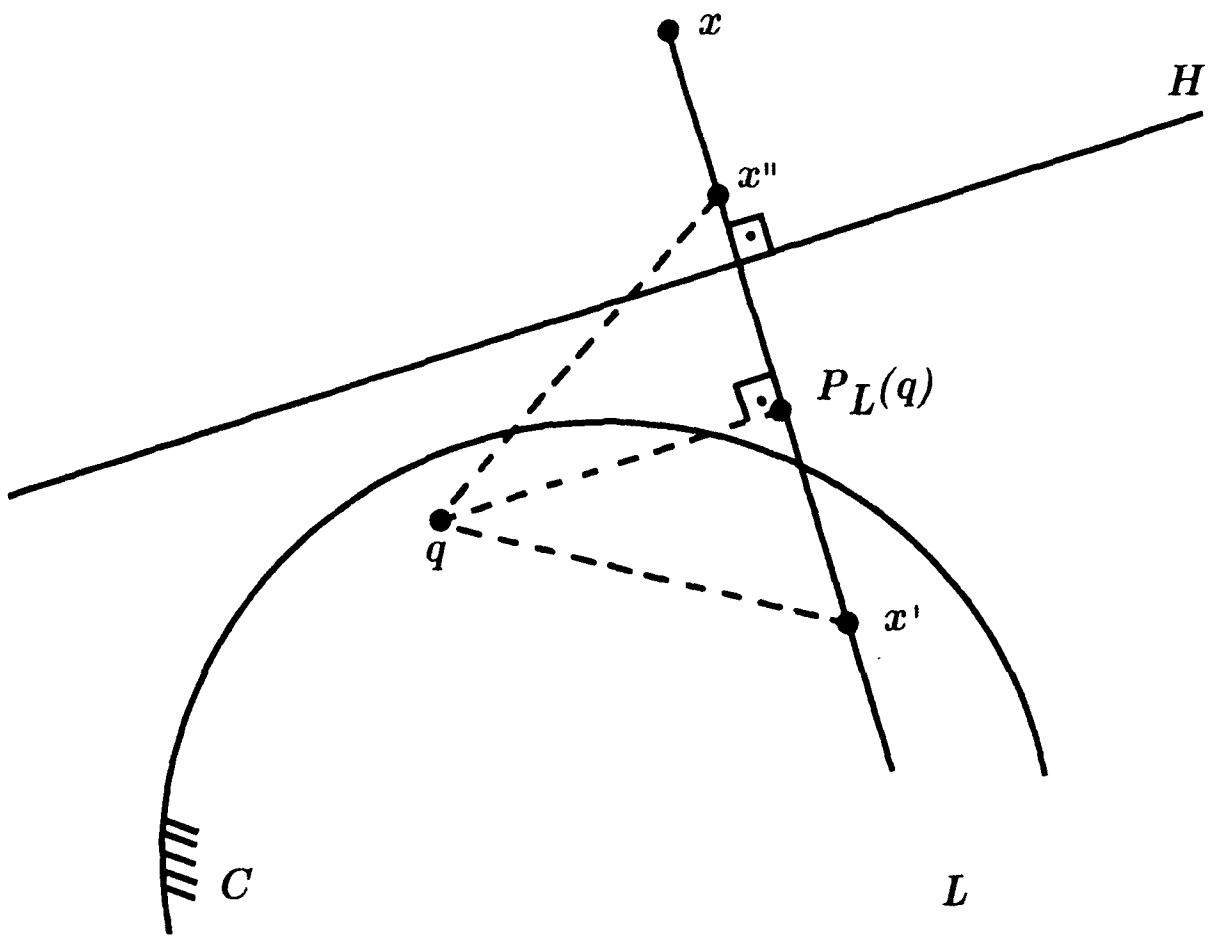

**Figure 5.7** Proof of Lemma 5.5.1.

The following results will be used in the convergence analysis that follows.

**Lemma 5.5.1** *Let $C \subseteq \mathbb{R}^n$ be a convex set. If $q \in C$ and $x' \in A_C(x, \lambda)$, then*

$$\| x' - q \| \leq \| x - q \| . \tag{5.28}$$

**Proof** Let $L$ denote the line through $x$ and $x'$ and set

$$x'' \doteq \begin{cases} x', & \text{if } \| x' - x \| \leq \| P_L(q) - x \|, \\ 2P_L(q) - x', & \text{otherwise.} \end{cases} \tag{5.29}$$

Basic geometric considerations show that $\| x'' - q \| = \| x' - q \|$. Using Theorem 2.4.1, as applied to the Bregman function $f(x) = \frac{1}{2} \| x \|^2$, with $\Omega$ being the half-space whose bounding hyperplane is perpendicular to the line $L$ and passes through the point $x''$, $y = x$, $P_\Omega(y) = x''$ and $z = q$, yields $0 \leq \| x'' - x \|^2 \leq \| q - x \|^2 - \| q - x'' \|^2$ and the required result follows. ■

Figure 5.7 illustrates the situation. The next lemma refines this result.

**Lemma 5.5.2** *Let $q$ be a point in a convex set $C \subseteq \mathbb{R}^n$; let $V$ be a bounded set in $\mathbb{R}^n$; and let $\delta$ and $\eta$, used in the definition of $A_C(x, \lambda)$, be such that $0 < \delta \leq 1$ and $0 \leq \eta \leq 1$. Then there exists a constant $\alpha > 0$ such that for every $x \in V \setminus C$ and for every $x' \in A_C(x, \lambda)$*

$$\| x' - q \| \leq \| x - q \| \, (1 - \alpha d^2(x, C)). \tag{5.30}$$

**Proof** Consult Figure 5.7, remembering that $q$ is fixed in $C$. As in the previous lemma we have $\| x - q \|^2 \geq \| x'' - q \|^2 + \| x - x'' \|^2$ because the angle $qx''x$ is obtuse. Therefore,

$$\begin{aligned} \| x' - q \|^2 = \| x'' - q \|^2 &\leq \| x - q \|^2 - \| x - x'' \|^2 \\ &= \| x - q \|^2 \left(1 - \frac{\| x - x'' \|^2}{\| x - q \|^2}\right) \\ &\leq \| x - q \|^2 \left(1 - \frac{\| x - x'' \|^2}{\| x - q \|^2} + \frac{1}{4}\frac{\| x - x'' \|^4}{\| x - q \|^4}\right) \\ &= \| x - q \|^2 \left(1 - \frac{1}{2}\frac{\| x - x'' \|^2}{\| x - q \|^2}\right)^2 . \end{aligned} \tag{5.31}$$

But $\| x - x'' \| \geq \eta\delta d(x, C)$, and since $V$ is bounded, there exists an $r$ such that for every $x \in V$, $\| x - q \| \leq r$. Therefore,

$$\| x' - q \| \leq \| x - q \| \left(1 - \frac{\eta^2\delta^2 d^2(x, C)}{2r^2}\right),$$

and (5.30) holds with $\alpha \doteq \eta^2\delta^2/2r^2$. ■

Observe that if $x \in C$ in the last lemma, then (5.30) reduces to (5.28). We present now the algorithmic scheme called the $(\delta, \eta)$-Algorithm.

**Algorithm 5.5.1 The $(\delta, \eta)$-Algorithm**

**Step 0:** (Initialization.) $x^0 \in \mathbb{R}^n$ is arbitrary.
**Step 1:** (Iterative step.)

$$x^{\nu+1} \in A_{Q_{i(\nu)}}(x^\nu, \lambda_\nu), \tag{5.32}$$

**Control:** $\{i(\nu)\}_{\nu=0}^{\infty}$ is repetitive (cf. Definition 5.1.1) on $I$.

Observe that for $\delta = 1$ and $\eta > 0$, the $(\delta, \eta)$-Algorithm coincides with the SOP method (Algorithm 5.2.1) because then all the hyperplanes in the family $\mathcal{H}_{x,Q}$ are also supporting to $Q$ at the point closest to $x$. Therefore, the proof of convergence of the $(\delta, \eta)$-Algorithm given next, for positive $\delta$ and $\eta$, also proves the convergence of the SOP method.

**Theorem 5.5.1** *A sequence $\{x^\nu\}_{\nu=0}^{\infty}$ generated by a $(\delta, \eta)$-algorithm with $0 < \delta \leq 1$ and $0 < \eta \leq 1$ converges to a point in the nonempty intersection set $Q = \cap_{i=1}^m Q_i$.*

In the proof we will need the following lemmas.

**Lemma 5.5.3** *A sequence $\{x^\nu\}_{\nu=0}^{\infty}$ generated by a $(\delta, \eta)$-Algorithm, with $0 \leq \delta \leq 1$, $0 \leq \eta \leq 1$, is Fejér-monotone with respect to the solution set $Q \doteq \cap_{i\in I} Q_i$ of the convex feasibility problem, i.e.,*

$$\| x^{\nu+1} - s \| \leq \| x^{\nu} - s \|, \tag{5.33}$$

*for every $s \in Q$ and all $\nu \geq 0$.*

**Proof** Follows by repeated use of Lemma 5.5.1. ■

**Lemma 5.5.4** *Given any $s \in Q$ and any $u \notin Q$ there exists an $r > 0$ such that if $u \notin Q_i$ then, for every $x \in B(u,r)$ any $x^i$ such that $x^i \in A_{Q_i}(x,\lambda)$, has the property*

$$\| x^i - s \| < \| u - s \| . \tag{5.34}$$

**Proof** Let $V$ be a ball centered at $u$. By Lemma 5.5.2 there exist positive numbers $\alpha_1, \alpha_2, \ldots, \alpha_m$, such that for any $i = 1, 2, \ldots, m$, $x \in V \setminus Q_i$ and $x^i \in A_{Q_i}(x,\lambda)$ implies that

$$\| x^i - s \| \leq \| x - s \| \, (1 - \alpha_i d^2(x, Q_i)).$$

Let $\alpha \doteq \min\{\alpha_i \mid 1 \leq i \leq m\}$, then

$$\| x^i - s \| \leq \| x - s \| \, (1 - \alpha d^2(x, Q_i)). \tag{5.35}$$

Choose $r' > 0$ such that if $u$ is outside any $Q_i$ then a ball with radius $r'$ around it will not meet with $Q_i$, i.e., $u \notin Q_i$ implies $B(u,r') \cap Q_i = \emptyset$. Then, for $x \in B(u,r')$ and $u \notin Q_i$,

$$d(x, Q_i) \geq d(u, Q_i) - r'. \tag{5.36}$$

If we denote $\overline{I}(u) \doteq \{i \mid u \notin Q_i\}$ and $c \doteq \min\{(d(u,Q_i) - r') \mid i \in \overline{I}(u)\}$ then $\| x^i - s \| \leq \| x - s \| \, (1 - \alpha c^2)$. Now, for any $r$ such that $0 < r \leq r'$, and $x \in B(u,r)$

$$\| x^i - s \| \leq \| (x - u) + (u - s) \| \, (1 - \alpha c^2) \leq (r + \| u - s \|)(1 - \alpha c^2),$$

and therefore, since $1 - \alpha c^2 < 1$, (5.34) holds if $r$ is sufficiently small. ■

**Proof** **(of Theorem 5.5.1)** By Lemma 5.5.3, the sequence $\{x^\nu\}_{\nu=0}^{\infty}$ is Fejér-monotone with respect to $Q$. Since every Fejér-monotone sequence is bounded, so is $\{x^\nu\}_{\nu=0}^{\infty}$ and thus, has an accumulation point $u$. This means that for any arbitrary $\epsilon > 0$, and for every $\nu$, there exists a $\mu > \nu$ such that

$$\| x^{\mu} - u \| \leq \epsilon. \tag{5.37}$$

Let $s$ be any point fixed in $Q$. Since by Lemma 5.5.3, $\| x^{\mu} - s \| \leq \| x^{\nu} - s \|$ for $\mu > \nu$, it follows from (5.37) that

$$\| u - s \| \leq \| x^{\nu} - s \| \quad \text{for all } \nu \geq 0. \tag{5.38}$$

If $u \in Q$, then by Lemma 5.5.3 $\lim_{\nu \to \infty} \| x^{\nu} - u \|$ exists, and from (5.37) this limit is zero.

To complete the proof we show that $u$ must be in $Q$. For $r$ as in Lemma 5.5.4, consider $x^\nu$ such that $\| x^\nu - u \| \leq r$. Suppose that $u \notin Q_{i(\nu)}$, then by Lemma 5.5.4 $\| x^{\nu+1} - s \| < \| u - s \|$, contradicting (5.38). Thus, $u \in Q_{i(\nu)}$, and by Lemma 5.5.1 $\| x^{\nu+1} - u \| \leq \| x^\nu - u \| \leq r$, and repetition of the argument yields $u \in Q_{i(\nu+1)}$.

Similarly, $u \in Q_{i(\nu)}$ for all $\mu \geq \nu$ and since the control $\{i(\nu)\}_{\nu=0}^\infty$ is repetitive, $u \in Q_i$ for all $i = 1, 2, \ldots, m$, implying that $u \in Q$. ■

Observe that if the interior of $Q$, int $Q \neq \emptyset$, then Lemmas 5.5.2 and 5.5.4 are not needed, because then the accumulation point $u$ has to be a limit point of $\{x^\nu\}_{\nu=0}^\infty$. To see this, suppose that $u^1$ and $u^2$ are two accumulation points of $\{x^\nu\}_{\nu=0}^\infty$. Then, for every $s \in Q$ we would have $\| u^1 - s \| = \| u^2 - s \|$ by Lemma 5.5.3, and since $Q$ has a nonempty interior, $u^1 = u^2$. Finally, the fact that the left-hand side of

$$\| x^{\nu+1} - x^\nu \| \geq \delta\eta d(x^\nu, Q_{i(\nu)}) \text{ for all } \nu \geq 0, \tag{5.39}$$

tends to zero as $\nu$ tends to infinity, and the facts that $\delta$ and $\eta$ are positive, and $\{i(\nu)\}_{\nu=0}^\infty$ is repetitive, imply that $d(u, Q_i) = 0$, for all $i = 1, 2, \ldots, m$, so that $u \in Q$.

The algorithm described next represents a particular realization of the $(\delta, \eta)$-Algorithm discussed above. This algorithm requires that an interior point in each individual set $Q_i$ be available beforehand. In exchange for this requirement the algorithm eliminates the need to project onto the sets $Q_i$ directly, as in the SOP method. As before, bd $C$ denotes the boundary of the set $C$.

**Algorithm 5.5.2 A Realization of the $(\delta, \eta)$-Algorithm Using Interior Points**

**Step 0:** Find a set of interior points $y^i \in Q_i$, for all $i = 1, 2, \ldots, m$.
**Step 1:** Choose $x^0 \in \mathbb{R}^n$ arbitrarily.
**Step 2:** Given $x^\nu$, pick $Q_{i(\nu)}$ according to the repetitive control sequence.
**Step 3:** If $x^\nu \in Q_{i(\nu)}$ then set $x^{\nu+1} = x^\nu$ and return to Step 2.
**Step 4:** If $x^\nu \notin Q_{i(\nu)}$ then construct the line $L_\nu$ through the points $x^\nu$ and $y^{i(\nu)}$.
**Step 5:** Denote by $z^\nu$ the point closest to $x^\nu$ in the set $L_\nu \cap Q_{i(\nu)}$.
**Step 6:** Construct a hyperplane $H_\nu$ separating $x^\nu$ from $Q_{i(\nu)}$ and supporting $Q_{i(\nu)}$ at $z^\nu$.
**Step 7:** Project $x^\nu$ orthogonally onto $H_\nu$ to get the point $P_{H_\nu}(x^\nu)$.
**Step 8:** Choose a relaxation parameter $\lambda_\nu$, such that $0 < \eta \leq \lambda_\nu \leq 2 - \eta$, and compute the next iterate by

$$x^{\nu+1} = x^\nu + \lambda_\nu (P_{H_\nu}(x^\nu) - x^\nu), \tag{5.40}$$

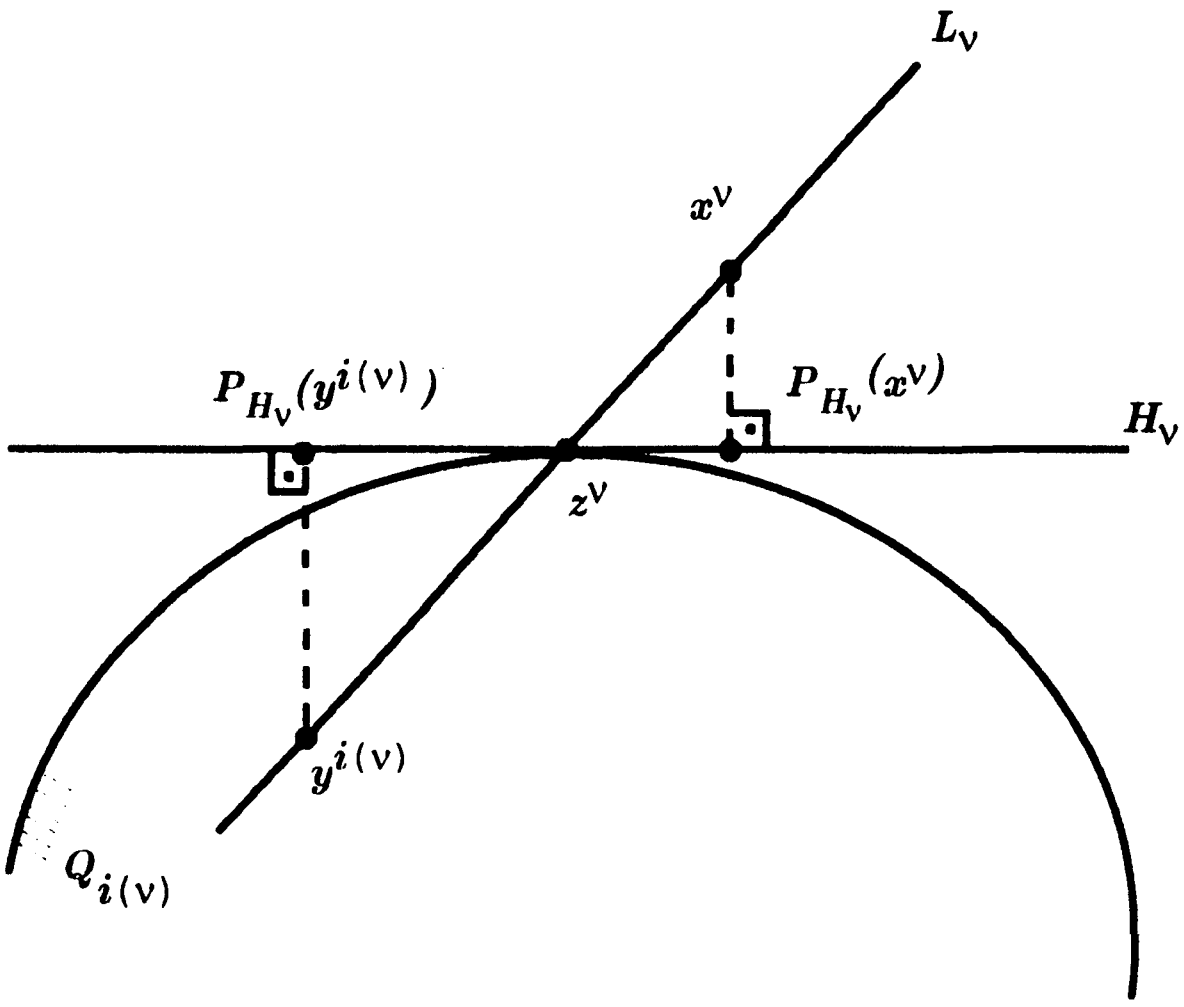

**Figure 5.8** Geometric description for the proof of Theorem 5.5.2.

and go back to Step 2.

We prove that Algorithm 5.5.2 converges to a solution of the convex feasibility problem by showing that it is an appropriate $(\delta, \eta)$-Algorithm.

**Theorem 5.5.2** *A sequence $\{x^\nu\}_{\nu=0}^\infty$ generated by Algorithm 5.5.2 converges to a solution $x^* \in Q$ of the convex feasibility problem.*

**Proof** The Algorithm 5.5.2 is obviously a $(\delta, \eta)$-Algorithm with $\eta > 0$ and $0 \leq \delta \leq 1$ (see Figure 5.8). To invoke Theorem 5.5.1 we have to show that for such an algorithm $\delta > 0$ always holds.

By Lemma 5.5.3 the sequence $\{x^\nu\}_{\nu=0}^\infty$ is Fejér-monotone with respect to $Q$ and thus bounded. Figure 5.8 helps to visualize the rest of the proof. Suppose $x^\nu \notin Q_{i(\nu)}$, then

$$\| P_{H_\nu}(x^\nu) - x^\nu \| = \frac{\| x^\nu - z^\nu \| \| y^{i(\nu)} - P_{H_\nu}(y^{i(\nu)}) \|}{\| y^{i(\nu)} - z^\nu \|}. \tag{5.41}$$

Also,

$$\| x^\nu - z^\nu \| \geq d(x^\nu, Q_{i(\nu)}), \tag{5.42}$$

and

$$\| y^{i(\nu)} - P_{H_\nu}(y^{i(\nu)}) \| \geq d(y^{i(\nu)}, \text{bd}\, Q_{i(\nu)}) \geq d > 0, \tag{5.43}$$

where

$$d \doteq \min\{d(y^i, \text{bd}\, Q_i) \mid 1 \leq i \leq m\}. \tag{5.44}$$

Since $\{x^\nu\}_{\nu=0}^\infty$ is bounded, there exists a positive $T$ such that

$$\| y^{i(\nu)} - z^\nu \| \leq T, \text{ for all } \nu \geq 0. \tag{5.45}$$

Combining these inequalities with (5.41) validates that

$$\| P_{H_\nu}(x^\nu) - x^\nu \| \geq \frac{d}{T} d(x^\nu, Q_{i(\nu)}), \text{ for all } \nu \geq 0, \tag{5.46}$$

proving that the algorithm is a $(\delta, \eta)$-Algorithm with $\delta \doteq d/T > 0$. ∎

We conclude our study of $(\delta, \eta)$-Algorithms by writing down a special version of the last algorithm for the case when the sets $Q_i$ are represented via convex functions, i.e., $Q_i \doteq \{x \in \mathbb{R}^n \mid g_i(x) \leq 0\}$, where $g_i : \mathbb{R}^n \to \mathbb{R}$ are convex for all $i = 1, 2, \ldots, m$, and where the control sequence $\{i(\nu)\}$ is cyclic.

**Algorithm 5.5.3 $(\delta, \eta)$-Algorithm for Convex Inequalities**

**Step 0:** Find points $y^1, y^2, \ldots, y^m$, such that $g_i(y^i) < 0$ for all $i = 1, \ldots, m$.
**Step 1:** Choose $x^0 \in \mathbb{R}^n$ arbitrarily.
**Step 2:** Given $x^\nu$ pick $Q_{i(\nu)}$ using a cyclic control $i(\nu) = \nu \bmod m + 1$.
**Step 3:** If $g_{i(\nu)}(x^\nu) \leq 0$, set $x^{\nu+1} = x^\nu$ and return to Step 2.
**Step 4:** If $g_{i(\nu)}(x^\nu) > 0$, choose some $\theta$ such that $0 \leq \theta \leq 1$ and define $z(\theta) \doteq \theta y^{i(\nu)} + (1 - \theta)x^\nu$.
**Step 5:** Solve the single (nonlinear) equation

$$g_{i(\nu)}(z(\theta)) = 0, \tag{5.47}$$

denote by $\theta_\nu$ the smallest value of $\theta$ for which $z(\theta)$ solves (5.47), and set $z^\nu \doteq z(\theta_\nu)$.

**Step 6:** Calculate a subgradient

$$t^\nu \in \partial g_{i(\nu)}(z^\nu) \tag{5.48}$$

(if $g_{i(\nu)}$ is differentiable at $z^\nu$ then $t^\nu = \nabla g_{i(\nu)}(z^\nu)$) and construct a supporting separating hyperplane $H_\nu$ given by

$$H_\nu = \{x \in \mathbb{R}^n \mid \langle t^\nu, x \rangle = \langle t^\nu, z^\nu \rangle\}. \tag{5.49}$$

**Step 7:** Compute the next iterate $x^{\nu+1}$ according to

$$x^{\nu+1} = x^\nu + \lambda_\nu \frac{\langle t^\nu, z^\nu \rangle - \langle t^\nu, x^\nu \rangle}{\| t^\nu \|^2} t^\nu, \tag{5.50}$$

where $\lambda_\nu$ is a user-chosen relaxation parameter, $0 < \eta \leq \lambda_\nu \leq 2 - \eta$, and return to Step 2.

Observe that the subgradient $t^\nu$ in (5.48) gives a direction that is normal to $Q_{i(\nu)}$ at the point $z^\nu$. Equations (5.49) and (5.50) are well defined because $t^\nu \neq 0$. For if $t^\nu = 0$ then $z^\nu$ is a minimum point of $g_{i(\nu)}$, but then $g_{i(\nu)}(y^{i(\nu)}) < 0$ and $g_{i(\nu)}(z^\nu) = 0$ would constitute a contradiction.

## 5.6 The Block-Iterative Projections Algorithm

The method of successive orthogonal projections SOP (Algorithm 5.2.1), as well as its generalizations in the cyclic subgradients projections CSP method (Algorithm 5.3.1) and in the $(\delta, \eta)$-Algorithm (Algorithm 5.5.1), are all sequential algorithms as defined in Section 1.3. In contrast to the row-action nature of these algorithms, a *fully simultaneous* iterative projections method is obtained by projecting the current iterate $x^\nu$ onto each and every set $Q_i$ , $i = 1, 2, \ldots, m$, and then taking the next iterate $x^{\nu+1}$ to be a convex combination of all the projections $P_{Q_i}(x^\nu)$. The projections onto multiple sets can be performed in parallel.

The Block-Iterative Projections (BIP) method, which we study next, is an algorithmic scheme that encompasses both the sequential SOP and the fully simultaneous ideas. But the BIP method is even more general because it allows the processing of blocks (i.e., groups of sets $Q_i$) that need not be fixed in advance, but which could change dynamically throughout the iterations. The number of blocks, their sizes, and the assignment of sets $Q_i$ to blocks may all vary, provided that the weights attached to the sets fulfill a simple technical condition, given in Theorem 5.6.1 below.

Let $I \doteq \{1, 2, \ldots, m\}$ and let $\{Q_i \mid i \in I\}$ be a finite family of closed convex sets. The intersection $Q \doteq \cap_{i \in I} Q_i$ is again assumed to be nonempty. We introduce *weight vectors*, $w = (w(i))_{i=1}^m$ such that $w(i) \geq 0$, for all $i \in I$, and $\sum_{i \in I} w(i) = 1$.

A sequence $\{w^\nu\}_{\nu=0}^\infty$ of weight vectors is called *fair* if, for every $i \in I$, there exist infinitely many values of $\nu$ for which $w^\nu(i) > 0$. Given a weight vector $w$ we define the convex combination $P_w(x) \doteq \sum_{i \in I} w(i) P_i(x)$, where $P_i(x)$ is the orthogonal projection of $x$ onto the set $Q_i$. The general scheme for block-iterative projections is as follows.

### Algorithm 5.6.1 The Block-Iterative Projections Algorithm (BIP)

**Step 0:** (Initialization.) $x^0 \in \mathbb{R}^n$ is arbitrary.

**Step 1:** (Iterative step.)

$$x^{\nu+1} = x^\nu + \lambda_\nu \left( P_{w^\nu}(x^\nu) - x^\nu \right), \tag{5.51}$$

where $\{w^\nu\}_{\nu=0}^\infty$ is a fair sequence of weight vectors.

**Relaxation parameters:** $\{\lambda_\nu\}_{\nu=0}^\infty$ is a sequence of user-determined relaxation parameters such that for all $\nu \geq 0$, $\epsilon_1 \leq \lambda_\nu \leq 2 - \epsilon_2$, for any arbitrarily small $\epsilon_1, \epsilon_2 \geq 0$.

The special case when the weight vectors are given by $w^\nu \doteq e^{i(\nu)}$, with $e^q \in \mathbb{R}^n$ being the $q$th standard basis vector, having one in its $q$th coordinate and zeros elsewhere, gives rise to a sequential row-action method. With this choice of weight vectors the BIP algorithm coincides with the SOP method (Algorithm 5.2.1), and $\{i(\nu)\}$ is the control sequence of the algorithm. For example, a cyclic control sequence dictates $i(\nu) = \nu \bmod m + 1$.

At the other extreme, choosing any sequence of weight vectors $\{w^\nu\}$ with $w^\nu(i) \neq 0$ for all $\nu \geq 0$ and all $i \in I$, leads to a fully simultaneous algorithm in which *all* sets $\{Q_i\}_{i=1}^m$ are being acted upon in every iterative step. This approach was first suggested for the linear case by Cimmino (see Section 5.11), and therefore fully simultaneous algorithms are often said to be of the *Cimmino-type.*

A block-iterative version with fixed blocks is obtained from Algorithm 5.6.1 by partitioning the indices of $I$ as $I = I_1 \cup I_2 \cup \cdots \cup I_M$ into $M$ blocks and using weight vectors of the form $w^\nu = \Sigma_{i \in I_{t(\nu)}} w^\nu(i) e^i$, where $\{t(\nu)\}$ is a control sequence over the set $\{1, 2, \ldots, M\}$ of block indices. In this case, if one considers the linear inequalities feasibility problem with

$$Q_i \doteq \{x \in \mathbb{R}^n \mid \langle a^i, x \rangle \leq b_i\}, \tag{5.52}$$

for every $i \in I$, where $a^i \in \mathbb{R}^n$, $b_i \in \mathbb{R}$, then a fixed-block version of the AMS relaxation method (Algorithm 5.4.2) of the following form is obtained.

**Algorithm 5.6.2 The Block-AMS Algorithm**

**Step 0:** (Initialization.) $x^0 \in \mathbb{R}^n$ is arbitrary.
**Step 1:** (Iterative step.)

$$x^{\nu+1} = x^\nu + \lambda_\nu \left( \sum_{i \in I_{t(\nu)}} w^\nu(i) c_i(x^\nu) a^i \right), \tag{5.53}$$

where $\{t(\nu)\}$ is a cyclic (or almost cyclic) control sequence on $\{1, 2, \ldots, M\}$, and $c_i(x^\nu)$ are defined, for $i = 1, 2, \ldots, m$, by

$$c_i(x^\nu) \doteq \min \left( 0, \frac{b_i - \langle a^i, x^\nu \rangle}{\| a^i \|^2} \right). \tag{5.54}$$

**Relaxation parameters:** For all $\nu \geq 0$, $\lambda_\nu$ is confined to the interval $\epsilon_1 \leq \lambda_\nu \leq 2 - \epsilon_2$ for any arbitrarily small $\epsilon_1, \epsilon_2 > 0$.

The generality of the definition of a fair sequence of weight vectors also permits variable block sizes and variable assignments of inequalities into the blocks that can be used. Figure 5.9 demonstrates an iterative step of the block-AMS algorithm for a block of three half-spaces $Q_i$ of the form (5.52). The next iterate $x^{\nu+1}$ is a convex combination, weighted by the components $w^\nu(i)$ of the current weight vector $w^\nu$ of the orthogonal projections of the current iterate $x^\nu$ onto each of the half-spaces. Choice of a relaxation parameter $\lambda_\nu \neq 1$ would place $x^{\nu+1}$ somewhere along the line through $x^\nu$ and the point marked $x^{\nu+1}(\lambda_\nu = 1)$ in the figure.

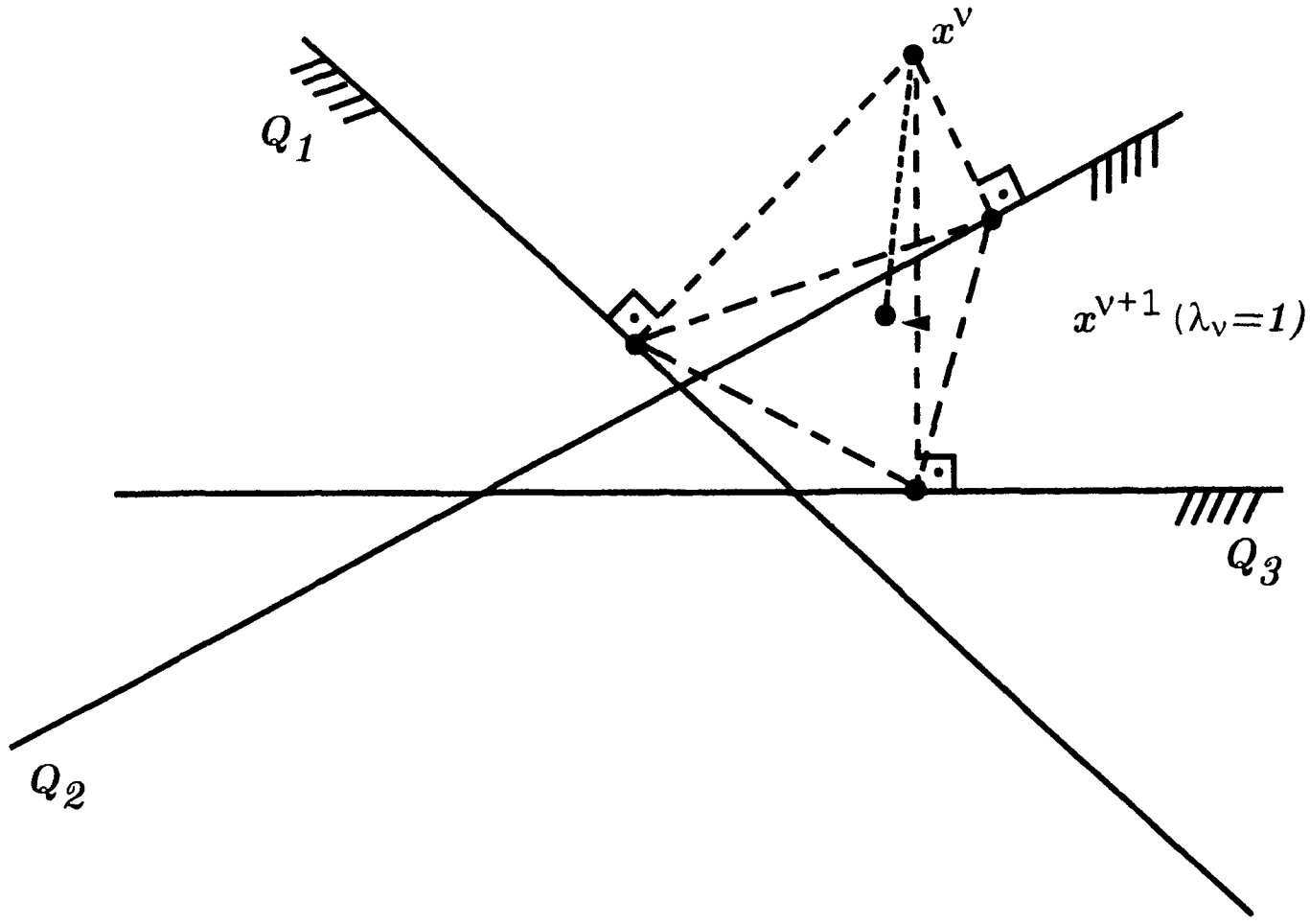

**Figure 5.9** The iterative step of the block-AMS algorithm (Algorithm 5.6.2) for a block of three inequalities and unity relaxation, $\lambda_\nu = 1$.

### 5.6.1 Convergence of the BIP algorithm

First we explain the notation that will be used. For any $J \subseteq I \doteq \{1, \ldots, m\}$ and any weight vector $w$, define $w(J) \doteq \sum_{j \in J} w(j)$. For any set $B \subseteq \mathbb{R}^n$ denote by $J(B)$ the set of indices of those sets $Q_i$ which do not meet $B$, i.e., $J(B) \doteq \{i \in I \mid B \cap Q_i = \emptyset\}$. For a singleton $B = \{x\}$ we write $J(\{x\}) = J(x) = \{i \in I \mid x \notin Q_i\}$. $B(x, \rho) \doteq \{y \in \mathbb{R}^n \mid \| x - y \| \leq \rho\}$ is the ball with radius $\rho$ centered at $x \in \mathbb{R}^n$. Finally, we denote a relaxed projection of a point $x \in \mathbb{R}^n$ onto $Q_i$ by $x + \lambda(P_i(x) - x) \doteq P_{i,\lambda}(x)$ and a relaxed convex combination of projections by $x + \lambda(P_w(x) - x) \doteq P_{w,\lambda}(x)$, where $\lambda \in \mathbb{R}$ is a relaxation parameter and $w$ is a weight vector.

**Proposition 5.6.1** *Given any $i$, $1 \leq i \leq m$, if $x \in \mathbb{R}^n$ then for every $y \in Q_i$ and every $\lambda$ such that $\epsilon_1 \leq \lambda \leq 2 - \epsilon_2$, with some arbitrarily small fixed $\epsilon_1, \epsilon_2 > 0$,*

$$\| P_{i,\lambda}(x) - y \| \leq \| x - y \| . \tag{5.55}$$

*Moreover, if $x \notin Q_i$, then the inequality is strict.*

**Proof** The following holds

$$\begin{aligned}\| P_{i,\lambda}(x) - y \|^2 &= \| x - y \|^2 + \lambda^2 \| P_i(x) - x \|^2 + 2\lambda\langle x - y, P_i(x) - x\rangle \\ &= \| x - y \|^2 + \lambda(\lambda - 2) \| P_i(x) - x \|^2 + 2\lambda\langle P_i(x) - x, P_i(x) - y\rangle \\ &\le \| x - y \|^2 - \epsilon_1\epsilon_2 \| P_i(x) - x \|^2 \le \| x - y \|^2 . \end{aligned} \tag{5.56}$$

This follows from the fact that the vector $x - P_i(x)$ is normal to $Q_i$ at $P_i(x)$, which provides that $\langle P_i(x) - x, P_i(x) - y\rangle \le 0$ for every $y \in Q_i$. This last inequality is Kolmogorov's criterion; use Theorem 2.4.2 and apply (2.50) for the case $f(x) = \frac{1}{2} \| x \|^2$. The last inequality of (5.56) is strict if $x \notin Q_i$. ■

**Proposition 5.6.2** *Let $q \in Q$ and let $\lambda$ be such that $\epsilon_1 \le \lambda \le 2 - \epsilon_2$, with some arbitrarily small fixed $\epsilon_1, \epsilon_2 > 0$, and let $w$ be any weight vector. Then, for every $x \in \mathbb{R}^n$,*

$$\| P_{w,\lambda}(x) - q \| \le \| x - q \| . \tag{5.57}$$

**Proof** Obviously, $P_{w,\lambda}(x) = \sum_{i\in I} w(i)P_{i,\lambda}(x)$, and repeated use of Proposition 5.6.1 with $y = q$ shows that $P_{i,\lambda}(x) \in B = B(q, \| x - q \|)$. The convex combination must, therefore, also be in $B$. ■

**Proposition 5.6.3** *Let $u \in \mathbb{R}^n$, $J = J(u)$, let $\epsilon_1 \le \lambda \le 2 - \epsilon_2$, with some fixed $\epsilon_1, \epsilon_2 > 0$, and let $w$ be any weight vector. Then, for every $r > 0$, there exists a real nonnegative $\gamma$ such that if $\| x \| \le r$ then*

$$\| P_{w,\lambda}(x) - u \| \le \| x - u \| + \gamma w(J). \tag{5.58}$$

**Proof** Define

$$\gamma_1 \doteq \max\{\| P_{j,\lambda}(x) - u \| \mid \| x \| \le r,\ j \in J,\ \epsilon_1 \le \lambda \le 2 - \epsilon_2\},$$

and $\gamma \doteq \max\{\gamma_1, r + \| u \|\}$. Then, using Proposition 5.6.1 we obtain $\| P_{j,\lambda}(x) - u \| \le \| x - u \|$ whenever $j \notin J$. This implies that, based on the definitions of $\gamma_1$ and $\gamma$,

$$\begin{aligned}\| P_{w,\lambda}(x) - u \| &= \| \sum_{j\in J} w(j)\,(P_{j,\lambda}(x) - u) + \sum_{j\notin J} w(j)\,(P_{j,\lambda}(x) - u) \| \\ &\le w(J)\gamma_1 + (1 - w(J)) \| x - u \| \\ &= \| x - u \| + w(J)(\gamma_1 - \| x - u \|) \\ &\le \| x - u \| + \gamma w(J), \end{aligned} \tag{5.59}$$

which completes the proof. ■

**Proposition 5.6.4** *Let $q \in Q$, and let $B \subseteq \mathbb{R}^n$ be any compact set, $J \doteq J(B)$. Let $\lambda$ be such that $\epsilon_1 \le \lambda \le 2 - \epsilon_2$ with some fixed $\epsilon_1, \epsilon_2 > 0$, and* *let $w$ be any weight vector. Then there exists an $\alpha > 0$ such that, for every $x \in B$,*

$$\| P_{w,\lambda}(x) - q \| \leq \| x - q \| - \alpha w(J). \tag{5.60}$$

**Proof** Define

$$\alpha \doteq \min\{(\| x - q \| - \| P_{j,\lambda}(x) - q \|) \mid x \in B, j \in J, \epsilon_1 \leq \lambda \leq 2 - \epsilon_2\}.$$

From Proposition 5.6.1 and the compactness of $B$, it follows that $\alpha > 0$. From the definition of $\alpha$ we have that $\| P_{j,\lambda}(x) - q \| \leq \| x - q \| - \alpha$, for every $x \in B$ and every $j \in J$. Using this and the fact that if $j \notin J$ then $\| P_{j,\lambda}(x) - q \| \leq \| x - q \|$, we get

$$\| P_{w,\lambda}(x) - q \| = \| \sum_{j \in J} w(j)\,(P_{j,\lambda}(x) - q) + \sum_{j \notin J} w(j)\,(P_{j,\lambda}(x) - q) \|$$

$$\leq w(J)(\| x - q \| - \alpha) + (1 - w(J))\, \| x - q \| = \| x - q \| - \alpha w(J),$$

and the proof is completed. ∎

To formulate the convergence theorem we use some additional notation. For a given fair sequence $\{w^\nu\}$ of weight vectors we define the index set $F(\{w^\nu\}) \doteq \{i \in I \mid \sum_{\nu=0}^{\infty} w^\nu(i) = +\infty\}$. This condition guarantees that for each set $Q_i$ with $i \in F(\{w^\nu\})$ the weights $w^\nu(i)$ attached to the $i$th set $Q_i$, even if they decrease as $\nu$ goes to infinity, do not go to zero too quickly. We define the set $\hat{Q} \doteq \cap_{i \in F(\{w^\nu\})} Q_i$ with the convention that if $F(\{w^\nu\}) = \emptyset$ then $\hat{Q} = \mathbb{R}^n$.

**Theorem 5.6.1** *If $Q \neq \emptyset$, if $\{w^\nu\}$ is any fair sequence of weight vectors, and if $\{\lambda_\nu\}$ is any sequence of relaxation parameters for which $\epsilon_1 \leq \lambda_\nu \leq 2 - \epsilon_2$, for all $\nu \geq 0$, where $\epsilon_1, \epsilon_2 > 0$ are arbitrarily small but fixed, then any sequence $\{x^\nu\}$ generated by the BIP algorithm (Algorithm 5.6.1) converges to a point $x^* \in \hat{Q}$.*

**Proof** Any sequence $\{x^\nu\}$ generated by the BIP Algorithm is Fejér-monotone with respect to $Q$, i.e., $\| x^{\nu+1} - q \| \leq \| x^\nu - q \|$ for every $\nu \geq 0$ and any $q \in Q$. This follows from Proposition 5.6.2 and, in turn, implies that $\{x^\nu\}$ is bounded. Next we show that $\{x^\nu\}$ is convergent. Assuming the contrary, it has two or more distinct accumulation points. Let $u$ be one of them, let $v$ be the one closest to $u$ (if there are several such $v$, pick any one), and let $r \doteq \| u - v \|$. We first show that $u \in Q$. The sequence $\{\| x^\nu - q \|\}$ is monotonically decreasing and bounded from below. Since $u$ is an accumulation point of $\{x^\nu\}$ it follows that $\lim_{\nu \to \infty} \| x^\nu - q \| = \| u - q \|$, and that, for all $\nu = 0, 1, 2, \ldots$,

$$\| x^\nu - q \| \geq \| u - q \| . \tag{5.61}$$

Suppose that $u \notin Q$. Choose $\rho > 0$ such that $\rho < r/2$ and such that the ball $B = B(u,\rho)$ satisfies $B \cap Q_j = \emptyset$ for every $j \in J(u)$. Let $J = J(B)$, and let $\gamma$ and $\alpha$ be as in Propositions 5.6.3 and 5.6.4, respectively. Define $\tau \doteq \rho\alpha/(\gamma+\alpha)$ and choose an index $\nu$ such that $\| x^\nu - u \| < \tau$. Since $v$ is also an accumulation point and $\rho < r/2$, there exists an index $\mu > \nu$ such that $x^\mu \notin B$. Choose the first such $\mu$. Then, by Proposition 5.6.4 and the triangle inequality,

$$\| x^\mu - q \| \leq \| x^\nu - q \| - \alpha \sum_{t=\nu}^{\mu-1} w^t(J) < \| u - q \| + \tau - \alpha \sum_{t=\nu}^{\mu-1} w^t(J). \quad (5.62)$$

Using (5.61) and (5.62) we get

$$\sum_{t=\nu}^{\mu-1} w^t(J) < \frac{\tau}{\alpha}. \quad (5.63)$$

On the other hand, by Proposition 5.6.3

$$\| x^\mu - u \| \leq \| x^\nu - u \| + \gamma \sum_{t=\nu}^{\mu-1} w^t(J), \quad (5.64)$$

which, together with (5.63) yields

$$\| x^\mu - u \| < \tau + \frac{\gamma\tau}{\alpha} = \rho, \quad (5.65)$$

contradicting $x^\mu \notin B$, and hence showing that $u \in Q$. Using Fejér-monotonicity with regard to $u$, which amounts to the monotonic decrease of the sequence $\{\| x^\nu - u \|\}$, the conclusion that the sequence $\{x^\nu\}$ converges to $u$ follows. Having secured convergence of $x^\nu$ to $u$ still does not tell us about the nature of $u$ because $u \in Q$ has only been shown under the (false) assumption of several distinct accumulation points. Hence the limit $x^*$ of any sequence $\{x^\nu\}$ generated by the BIP algorithm still needs to be characterized, and we do so by showing that it belongs to $\hat{Q}$.

Assume that $x^* \notin Q_{i_0}$, for some $i_0 \in F(\{w^\nu\})$, i.e., $\sum_{\nu=0}^{\infty} w^\nu(i_0) = +\infty$. Choose a ball $B_1$ centered at $x^*$ such that $B_1 \cap Q_{i_0} = \emptyset$. Let $\nu$ be such that $x^\mu \in B_1$ whenever $\mu \geq \nu$. By Proposition 5.6.4 there exists $\alpha > 0$ for which

$$\| x^\mu - q \| \leq \| x^\nu - q \| - \alpha \sum_{t=\nu}^{\mu-1} w^t(J(B_1)), \quad (5.66)$$

for every $\mu \geq \nu$. But since $i_0 \in J(B_1)$ and $i_0 \in F(\{w^\nu\})$, we have

$$\lim_{\mu \to \infty} \sum_{t=\nu}^{\mu-1} w^t(J(B_1)) = +\infty,$$

implying that $\lim_{\mu \to \infty} \| x^\mu - q \| = -\infty$, which is impossible. ∎

The condition $\sum_{\nu=0}^{\infty} w^\nu(i) = +\infty$ is quite mild and, as mentioned before, it prevents the weights $w^\nu(i)$ attached to any particular set $Q_i$ from diminishing too quickly as iterations proceed. In practical cases, such as the row-action case, or cases where fixed weights are attached to sets $Q_i$, the condition holds for all $i \in I$ and then $\hat{Q} = Q$.

## 5.7 The Block-Iterative $(\delta, \eta)$-Algorithm

The $(\delta, \eta)$-Algorithm (Algorithm 5.5.1) replaces, as we have seen, the orthogonal projection of the current iterate $x^\nu$ onto the selected set $Q_{i(\nu)}$ by a projection of $x^\nu$ onto an appropriately constructed hyperplane $H_{i(\nu)}$. The BIP method uses orthogonal projections onto the sets but allows their grouping into blocks. Can these two different extensions of the SOP method be combined? In other words, if we perform steps of the $(\delta, \eta)$-Algorithm simultaneously onto the sets of a block $\{Q_i\}_{i \in I_t}$, for some $I_t \subseteq I$, and then take the (possibly relaxed) convex combination of the resulting points as the next iterate, will the process converge to a solution of the convex feasibility problem? The answer is affirmative and we can formulate the block-iterative-$(\delta, \eta)$-Algorithm as follows.

### Algorithm 5.7.1 The Block-Iterative $(\delta, \eta)$-Algorithm

**Step 0:** (Initialization.) $x^0 \in \mathbb{R}^n$ is arbitrary.

**Step 1:** (Iterative step.) Given $x^\nu$, choose a weight vector $w^\nu$ from a fair sequence and, for all $i$ for which $w^\nu(i) > 0$, pick vectors $\pi^i(x^\nu) \in A_{Q_i}(x^\nu, \lambda_{\nu,i})$ and compute the next iterate

$$x^{\nu+1} = x^\nu + \theta_\nu \left( \sum_{i=1}^{m} w^\nu(i) \pi^i(x^\nu) - x^\nu \right). \tag{5.67}$$

Here the sets $A_{Q_i}(x^\nu, \lambda_{\nu,i})$ are constructed according to (5.27) for some arbitrary fixed $0 < \eta \leq 1$.

**Relaxation parameters:** In addition to the "internal" relaxation parameters $\lambda_{\nu,i}$, which are determined when the choice of $\pi^i(x^\nu)$ is made, the "outer" relaxation parameters $\{\theta_\nu\}$ are confined to the interval $\epsilon_1 \leq \theta_\nu \leq 2 - \epsilon_2$, for some fixed arbitrary $\epsilon_1, \epsilon_2 > 0$.

The following convergence theorem can be proven by combining the proofs of Theorems 5.5.1 and 5.6.1.

**Theorem 5.7.1** *If $Q = \cap_{i \in I} Q_i \neq \emptyset$, if $\{w^\nu\}$ is any fair sequence of weight vectors, and if $0 < \delta \leq 1$ and $0 < \eta \leq 1$, then any sequence $\{x^\nu\}_{\nu=0}^\infty$ generated by Algorithm 5.7.1 converges to an $x^* \in \hat{Q}$.*

## 5.8 The Method of Successive Generalized Projections

In the previous sections of this chapter we used the method of successive orthogonal projection (SOP) as a starting point. From this starting point we developed and studied algorithms that allowed us to replace the orthogonal projection involved in each iterative step of SOP with another, hopefully less demanding, operation. The CSP algorithm and the $(\delta, \eta)$-Algorithm are examples. The development of block-iterations in BIP and in the block-iterative $(\delta, \eta)$-Algorithm extended the scope of the overall structure of algorithms. In particular, we went from the inherently sequential row-action algorithms to fully simultaneous or to block-iterations, which employ in every iterative step groups $\{Q_i\}_{i \in I_t}$, $I_t \subseteq I$, of the constraint sets of the convex feasibility problem.

Now we go back to our starting point, the SOP method. Since the orthogonal projection is but one case of the family of generalized projections that we studied in Chapter 2, it is quite natural to inquire whether we could replace the orthogonal projections in the SOP method with generalized projections and retain the overall convergence? The affirmative answer to this question is given by constructing a sequential algorithm, which asymptotically solves the convex feasibility problem (CFP) by performing successive generalized projections onto the individual sets $Q_i$, starting from an arbitrary point. We present this algorithm now.

Let $f \in \mathcal{B}(S)$ be some given Bregman function with zone $S$ (see Chapter 2), and let a family of closed convex sets $Q_i \subseteq \mathbb{R}^n, i \in I$, be given.

**Algorithm 5.8.1 The Method of Successive Generalized Projections (SGP)**

**Step 0:** (Initialization.) $x^0 \in S$ is arbitrary.
**Step 1:** (Iterative step.)

$$x^{\nu+1} = P_{Q_{i(\nu)}}(x^\nu), \ \nu = 0, 1, 2, \ldots, \tag{5.68}$$

where $P_{Q_{i(\nu)}}(x^\nu)$ now stands for a generalized projection of the point $x^\nu$ onto the set $Q_{i(\nu)}$, as given by Definition 2.1.2.

**Relaxation parameters:** No relaxation parameters are allowed.

**Control:** The sequence $\{i(\nu)\}_{\nu=0}^\infty$ is cyclic on $I$.

For this algorithm to be well defined we need the following assumption:

**Assumption 5.8.1** *The function $f \in \mathcal{B}(S)$ is zone consistent with respect to each $Q_i, i = 1, 2, \ldots, m$, (see Definition 2.2.1).*

**Theorem 5.8.1** *Let $f \in \mathcal{B}(S)$ be a given Bregman function with zone $S$. Let a CFP be given by a finite family of closed convex sets $Q_i$, $i \in I$, such that $Q \doteq \cap_{i \in I} Q_i$ and $Q \cap \overline{S} \neq \emptyset$, and assume that Assumption 5.8.1 holds. Under these conditions, any sequence $\{x^\nu\}_{\nu=0}^\infty$ generated by Algorithm 5.8.1 converges to a point $x^* \in Q$.*

**Proof** The proof consists of the following three steps:

**Step 1.** $\{x^\nu\}_{\nu=0}^\infty$ is bounded, hence contains at least one accumulation point.

**Step 2.** Every accumulation point of $\{x^\nu\}_{\nu=0}^\infty$ belongs to $Q$.

**Step 3.** $\{x^\nu\}_{\nu=0}^\infty$ has a unique limit point.

We now proceed to prove each of these steps.

**Step 1.** For any $z \in Q \cap \overline{S}$ we have, according to Theorem 2.4.1 and Algorithm 5.8.1, under Assumption 5.8.1,

$$D_f(x^{\nu+1}, x^\nu) \leq D_f(z, x^\nu) - D_f(z, x^{\nu+1}), \tag{5.69}$$

which yields, using Lemma 2.1.1, the inequality,

$$D_f(z, x^{\nu+1}) \leq D_f(z, x^\nu), \text{ for all } z \in Q \cap \overline{S}, \text{ for all } \nu \geq 0. \tag{5.70}$$

This property is called *$D_f$-Fejér-monotonicity of the sequence $\{x^\nu\}_{\nu=0}^\infty$ with respect to the set $Q$*, because it extends to $D_f$-functions the property of Fejér-monotonicity (see Definition 5.3.1). It implies that $\{x^\nu\}_{\nu=0}^\infty$ is bounded because repeated use of (5.70) gives

$$x^\nu \in L_2^f(z, \alpha), \text{ for all } \nu \geq 0, \tag{5.71}$$

with $\alpha \doteq D_f(z, x^0)$ and condition $(iv)$ of Definition 2.1.1 applies.

**Step 2.** Let $\{x^\nu\}_{\nu \in K}$ with $K \subseteq N_0 \doteq \{0, 1, 2, 3, \ldots\}$ be a subsequence converging to $x^*$. Let $\{x^\nu\}_{\nu \in M}$, $M \subseteq K$ be a subsequence of the former subsequence such that all its elements belong to one set, for example $Q_1$. From the sequences $\{x^{\nu+i-1}\}_{\nu \in M}$ for $i = 2, 3, \ldots, m$ one can extract convergent subsequences, so we assume, without loss of generality, that

$$\lim_{\nu \to \infty,\, \nu \in M} x^{\nu+i-1} = x^{*i}, \quad i = 1, 2, \ldots, m, \tag{5.72}$$

and that $x^{*1} = x^*$. Since $\{x^{\nu+i-1}\}_{\nu \in M} \subseteq Q_i$ and the sets $Q_i$ are closed, it follows that

$$x^{*i} \in Q_i, \text{ for all } i \in I. \tag{5.73}$$

By (5.70) and the nonnegativity of $D_f$, $\lim_{\nu \to \infty} D_f(z, x^\nu)$ exists for any $z \in Q \cap \overline{S}$. In view of (5.69) $\lim_{\nu \to \infty} D_f(x^{\nu+1}, x^\nu) = 0$, thus also

$$\lim_{\nu \to \infty,\ \nu \in M} D_f(x^{\nu+1}, x^\nu) = 0. \tag{5.74}$$

Condition $(vi)$ of Definition 2.1.1 can be used repeatedly to show that from (5.72) and (5.73),

$$x^* = x^{*1} = x^{*2} = \cdots = x^{*m}, \tag{5.75}$$

which by (5.73) shows that $x^* \in Q$.

**Step 3.** Let

$$\lim_{\nu \to \infty,\ \nu \in K_1} x^\nu = x^* \in Q \cap \overline{S}, \tag{5.76}$$

$$\lim_{\nu \to \infty,\ \nu \in K_2} x^\nu = x^{**} \in Q \cap \overline{S}, \tag{5.77}$$

for some $K_1 \subseteq N_0$ and $K_2 \subseteq N_0$, be two convergent subsequences of $\{x^\nu\}_{\nu=0}^\infty$, generated by Algorithm 5.8.1. Then, by Definition 2.1.1$(v)$, we have

$$\lim_{\nu \to \infty,\ \nu \in K_1} D_f(x^*, x^\nu) = 0, \tag{5.78}$$

and also

$$\lim_{\nu \to \infty,\ \nu \in K_2} D_f(x^{**}, x^\nu) = 0. \tag{5.79}$$

But $D_f$-Fejér-monotonicity with respect to $Q$ implies that

$$D_f(x^*, x^{\nu+1}) \leq D_f(x^*, x^\nu), \text{ for all } \nu \geq 0, \tag{5.80}$$

thus $\lim_{\nu \to \infty} D_f(x^*, x^\nu)$ exists and because of (5.78) is equal to zero. This implies that

$$\lim_{\nu \to \infty,\ \nu \in K_2} D_f(x^*, x^\nu) = 0, \tag{5.81}$$

and Definition 2.1.1$(vi)$ then yields $x^* = x^{**}$. ■

The SGP method generalizes the SOP method from orthogonal projections to generalized projections, but only when unity relaxation, i.e., $\lambda_\nu = 1$, for all $\nu \geq 0$, is used in Algorithm 5.2.1. No relaxation parameters have yet been built into Algorithm 5.8.1.

At this point it might seem unclear to the reader why one would prefer the SGP method over the SOP method. It is difficult enough, in general, to perform orthogonal projections at each iterative step, so that there seems to be no reason to replace those by nonorthogonal generalized projections. It turns out, however, that performing successive generalized projections leads to some very useful, special-purpose optimization methods which will be discussed in Chapter 6. Another important consequence of the SGP algorithm is the possibility of using different types of generalized projections within a single algorithm. This is the subject matter of the next section.

## 5.9 The Multiprojections Algorithm

Generalized distances give rise to generalized projections onto convex sets. An important question is whether or not one can use within the same projection algorithm different types of such generalized projections. This question has practical consequences in the area of signal detection and image recovery in situations that can be formulated mathematically as a convex feasibility problem. Using a product space formalism, we study here the convergence of a multiprojection algorithm. Our algorithm is fully simultaneous, i.e., it uses in each iterative step all sets of the convex feasibility problem. Different multiprojection algorithms can be derived from our algorithmic scheme by a judicious choice of the Bregman functions that govern the process.

### 5.9.1 The product space setup

We now introduce the product space setup and then discuss Bregman functions and generalized $D$-functions in this setup.

Let $\mathbb{R}^n$ be the Euclidean space with scalar product, norm, and distance denoted, respectively by $\langle\cdot,\cdot\rangle$, $\|\cdot\|$ and $d(\cdot,\cdot)$. Let the product space $\boldsymbol{V} \doteq \mathbb{R}^{nm} = \mathbb{R}^n \times \mathbb{R}^n \times \cdots \times \mathbb{R}^n$ be such that any vector $\boldsymbol{x} \in \boldsymbol{V}$ is $\boldsymbol{x} = (x^1, x^2, \ldots, x^m)$, where $x^i \in \mathbb{R}^n$, $i \in I \doteq \{1, 2, \ldots, m\}$. The scalar product in $\boldsymbol{V}$ is denoted and defined by

$$\langle\langle \boldsymbol{x}, \boldsymbol{y} \rangle\rangle \doteq \sum_{i=1}^{m} \langle x^i, y^i \rangle, \quad \text{for all } \boldsymbol{x}, \boldsymbol{y} \in \boldsymbol{V}, \tag{5.82}$$

and induces the norm $|||\cdot|||$ and distance $\boldsymbol{d}((\cdot,\cdot))$ in $\boldsymbol{V}$. Note the typographical convention of denoting all quantities related to the product space with boldface letters.

Given a finite family of closed convex sets $Q_i \subseteq \mathbb{R}^n, i \in I$, let $Q \doteq \cap_{i\in I} Q_i$ and assume that $Q \neq \emptyset$. The convex feasibility problem (CFP) in $\mathbb{R}^n$ is to find an $x^* \in Q$. Next, let us define in $\boldsymbol{V}$ the product set $\boldsymbol{Q} \doteq Q_1 \times \cdots \times Q_m = \prod_{i\in I} Q_i$ and the set

$$\boldsymbol{\Delta} \doteq \{\boldsymbol{b} \in \boldsymbol{V} \mid \boldsymbol{b} = (b, b, \ldots, b),\ b \in \mathbb{R}^n\}. \tag{5.83}$$

Note that $\boldsymbol{Q}$ is a closed convex subset of $\boldsymbol{V}$ and $\boldsymbol{\Delta}$ is a linear subspace of $\boldsymbol{V}$. The set $\boldsymbol{\Delta}$ is called the *diagonal set* and the mapping $\boldsymbol{\delta} : \mathbb{R}^n \to \boldsymbol{\Delta}$ defined by

$$\boldsymbol{\delta}(b) \doteq (b, b, \ldots, b) = \boldsymbol{b} \in \boldsymbol{\Delta}, \tag{5.84}$$

is called the *canonical mapping.*

In this setup it is obvious that the CFP in $\mathbb{R}^n$ is equivalent to finding a point in the nonempty intersection $\boldsymbol{Q} \cap \boldsymbol{\Delta}$ in $\boldsymbol{V}$, i.e.,

$$x^* \in Q \text{ if and only if } \boldsymbol{\delta}(x^*) \in \boldsymbol{Q} \cap \boldsymbol{\Delta}. \tag{5.85}$$

Let us now introduce Bregman functions. Let $f_i : \mathbb{R}^n \to \mathbb{R}$, $i \in I$, be a finite family of Bregman functions with zones $S_i \subseteq \mathbb{R}^n$, $i \in I$, respectively, i.e., $f_i \in \mathcal{B}(S_i)$, $i \in I$.

**Definition 5.9.1** *Given $f_i \in \mathcal{B}(S_i)$, $i \in I$, and a positive weight vector $w = (w_i) \in \mathbb{R}^m$, i.e., $w_i > 0$ for all $i \in I$, and $\sum_{i \in I} w_i = 1$, define the function $\boldsymbol{F}_w : \boldsymbol{V} \to \mathbb{R}$ by*

$$\boldsymbol{F}_w(\boldsymbol{v}) \doteq \sum_{i=1}^{m} w_i f_i(v^i), \tag{5.86}$$

*and the product set*

$$\boldsymbol{S} \doteq \prod_{i=1}^{m} S_i. \tag{5.87}$$

The $\boldsymbol{D}$-function associated with $\boldsymbol{F}_w$ can be expressed as follows.

**Lemma 5.9.1** *Given $f_i \in \mathcal{B}(S_i)$, $i \in I$, and any positive weight vector $w \in \mathbb{R}^m$, the following hold:*

$$\boldsymbol{D}_{\boldsymbol{F}}(\boldsymbol{x}, \boldsymbol{y}) = \sum_{i=1}^{m} w_i D_{f_i}(x^i, y^i), \tag{5.88}$$

*where $\boldsymbol{D}_{\boldsymbol{F}} : \overline{\boldsymbol{S}} \times \boldsymbol{S} \subseteq \boldsymbol{V}^2 \to \mathbb{R}$ is the D-function associated with $\boldsymbol{F}$ in $\boldsymbol{V}$ and the subscript $w$ has been, and will henceforth be, omitted from $\boldsymbol{F}_w$ for the sake of brevity, and*

$$\boldsymbol{F} \in \mathcal{B}(\boldsymbol{S}). \tag{5.89}$$

**Proof** By definition,

$$\boldsymbol{D}_{\boldsymbol{F}}(\boldsymbol{x}, \boldsymbol{y}) = \boldsymbol{F}(\boldsymbol{x}) - \boldsymbol{F}(\boldsymbol{y}) - \langle\langle \nabla \boldsymbol{F}(\boldsymbol{y}), \boldsymbol{x} - \boldsymbol{y} \rangle\rangle, \tag{5.90}$$

and the result (5.88) follows from (5.86), (5.82), and from use of (2.1) for each $f_i$.

Conditions $(i)$–$(iii)$ of Definition 2.1.1 are easy to verify. Condition $(iv)$ can be handled by negation. To treat in detail condition $(vi)$ of Definition 2.1.1 (condition $(v)$ can be treated in a similar way) note that, by assumption condition $(vi)$ holds for each $f_i$ for all $i \in I$, that is if $\lim_{\nu\to\infty} D_{f_i}(x^{i,\nu}, y^{i,\nu}) = 0$, $\lim_{\nu\to\infty} y^{i,\nu} = y^{i,*} \in \overline{S_i}$ and $\{x^{i,\nu}\}$ is bounded, then $\lim_{\nu\to\infty} x^{i,\nu} = y^{i,*}$ for every $i \in I$. Now, $\lim_{\nu\to\infty} \boldsymbol{D}_{\boldsymbol{F}}(\boldsymbol{x}^\nu, \boldsymbol{y}^\nu) = 0$ implies $\lim_{\nu\to\infty} D_{f_i}(x^{i,\nu}, y^{i,\nu}) = 0$, for each $i \in I$, by (5.88) and Lemma 2.1.1 and using $w > 0$. Also, $\lim_{\nu\to\infty} \boldsymbol{y}^\nu = \boldsymbol{y}^* \in \overline{\boldsymbol{S}}$ if and only if $\lim_{\nu\to\infty} y^{i,\nu} = y^{i,*} \in \overline{S_i}$ for all $i \in I$. This is so because of the equivalence between convergence in the product space and convergence in each component space, i.e., $\lim_{\nu\to\infty} |||\boldsymbol{y}^* - \boldsymbol{y}^\nu||| = 0$ if and only if $\lim_{\nu\to\infty} ||y^{i,*} - y^{i,\nu}|| = 0$ for all $i \in I$. Finally, boundedness of $\{\boldsymbol{x}^\nu\}$ in $\boldsymbol{V}$ is equivalent to the boundedness of $\{x^{i,\nu}\}$, for all $i \in I$, in $\mathbb{R}^n$, and the result (5.89) follows. ■

The function $\boldsymbol{F}_w$ defined by (5.86) will be called the Bregman function in $\boldsymbol{V}$ induced by the family $\{f_i\}_{i \in I}$ and the positive weight vector $w \in \mathbb{R}^m$.

### 5.9.2 Generalized projections in the product space

Next we characterize generalized projections onto $\boldsymbol{Q}$ and $\boldsymbol{\Delta}$. In what follows $\boldsymbol{P}_{\boldsymbol{Q}}^{\boldsymbol{F}}(\boldsymbol{b})$ denotes the Bregman projection of a point $\boldsymbol{b} \in \boldsymbol{V}$ onto a closed convex set $\boldsymbol{Q}$, where $\boldsymbol{F}$ is defined in (5.86).

**Lemma 5.9.2** *Given a finite family of Bregman functions $f_i \in \mathcal{B}(S_i)$ with zones $S_i$, $i \in I$, any positive weight vector $w \in \mathbb{R}^m$, and a finite family of closed convex subsets $Q_i \subseteq \mathbb{R}^n$ such that $Q_i \cap \overline{S_i} \neq \emptyset$, $i \in I$, then for any $\boldsymbol{x} \in \boldsymbol{V} \cap \boldsymbol{S}$,*

$$\boldsymbol{P}_{\boldsymbol{Q}}^{\boldsymbol{F}}(\boldsymbol{x}) = (P_{Q_1}^{f_1}(x^1), P_{Q_2}^{f_2}(x^2), \ldots, P_{Q_m}^{f_m}(x^m)). \tag{5.91}$$

**Proof** Let $\boldsymbol{P}_{\boldsymbol{Q}}^{\boldsymbol{F}}(\boldsymbol{x}) = \boldsymbol{v}$ , then we must solve

$$\boldsymbol{v} = \operatorname*{argmin}_{\boldsymbol{z} \in \boldsymbol{Q} \cap \overline{\boldsymbol{S}}} \boldsymbol{D}_{\boldsymbol{F}}(\boldsymbol{z}, \boldsymbol{x}). \tag{5.92}$$

By Lemma 5.9.1, this amounts to solving

$$\min\{\sum_{i=1}^{m} w_i D_{f_i}(z^i, x^i) \mid z^i \in Q_i \cap \overline{S}_i, \ i \in I\}, \tag{5.93}$$

but $D_{f_i}(z^i, x^i) \geq 0$ and $w_i > 0$ for all $i \in I$ means that this minimum is obtained when each $D_{f_i}$ reaches its minimum, i.e., when

$$v^i = \operatorname*{argmin}_{z^i \in Q_i \cap \overline{S}_i,} D_{f_i}(z^i, x^i) = P_{Q_i}^{f_i}(x^i), \ \text{ for all } i \in I, \tag{5.94}$$

and the proof is complete. ■

Given any family of Bregman functions $f_i \in \mathcal{B}(S_i)$, $i \in I$, a positive weight vector $w \in \mathbb{R}^m$ and some fixed $\boldsymbol{v} \in \boldsymbol{S}$, we define

$$\varphi(h) \doteq \sum_{i=1}^{m} w_i D_{f_i}(h, v^i). \tag{5.95}$$

Also, denote $S \doteq \cap_{i=1}^m S_i$ and assume that $S \neq \emptyset$.

**Assumption 5.9.1** *The Bregman functions $\{f_i\}_{i\in I}$ are such that for any $\boldsymbol{v} \in \boldsymbol{V}$ and any positive $w \in \mathbb{R}^m$:*

$$\operatorname{argmin}\{\varphi(h) \mid h \in \overline{S}\} \subseteq S. \tag{5.96}$$

This assumption, guaranteeing that the minimum of $\varphi(h)$ over $\overline{S}$ is attained at a point where $\nabla\varphi = 0$, is needed in the next lemma.

**Lemma 5.9.3** *If $S \doteq \cap_{i=1}^m S_i \neq \emptyset$ and Assumption 5.9.1 holds, then for any $\boldsymbol{v} \in \boldsymbol{V}$ and positive $w \in \mathbb{R}^m$,*

$$\boldsymbol{P}_{\boldsymbol{\Delta}}^{\boldsymbol{F}}(\boldsymbol{v}) = \boldsymbol{\delta}(b), \tag{5.97}$$

*where $b$ is the unique solution of*

$$\sum_{i=1}^{m} w_i \nabla f_i(b) = \sum_{i=1}^{m} w_i \nabla f_i(v^i). \tag{5.98}$$

**Proof** The left-hand side of (5.97) must be of the form $\boldsymbol{\delta}(b)$, for some $b \in \mathbb{R}^n$, where $\boldsymbol{\delta}(b)$ is the canonical mapping. By definition,

$$\boldsymbol{\delta}(b) = \operatorname*{argmin}_{\boldsymbol{\delta}(h)\in\overline{\boldsymbol{S}}} \boldsymbol{D}_{\boldsymbol{F}}(\boldsymbol{\delta}(h), \boldsymbol{v}). \tag{5.99}$$

The solution of the problem

$$\min_{h\in\overline{S}} \sum_{i=1}^{m} w_i D_{f_i}(h, v^i) \tag{5.100}$$

is attained at a point where $\nabla\varphi(h) = 0$ with

$$\varphi(h) = \sum_{i=1}^{m} w_i(f_i(h) - f_i(v^i) - \langle\nabla f_i(v^i), h - v^i\rangle), \tag{5.101}$$

and the result follows. ■

### 5.9.3 The simultaneous multiprojections algorithm and the split feasibility problem

We now use Algorithm 5.8.1 for finding a point in $\boldsymbol{Q} \cap \boldsymbol{\Delta}$. Let $f_i \in \mathcal{B}(S_i)$, for $i \in I$, be a finite family of Bregman functions and $w \in \mathbb{R}^m$ a given positive weight vector. Construct the induced Bregman function $\boldsymbol{F}$ on the product space and assume that $\boldsymbol{F}$ is zone consistent with respect to $\boldsymbol{\Delta}$ and $\boldsymbol{Q}$, i.e.,

$$\boldsymbol{y} \in \boldsymbol{S} \Rightarrow \boldsymbol{P}_{\boldsymbol{\Delta}}^{\boldsymbol{F}}(\boldsymbol{y}) \in \boldsymbol{S}, \tag{5.102}$$

$$\boldsymbol{y} \in \boldsymbol{S} \Rightarrow \boldsymbol{P}_{\boldsymbol{Q}}^{\boldsymbol{F}}(\boldsymbol{y}) \in \boldsymbol{S}. \tag{5.103}$$

Applying Algorithm 5.8.1 with $\boldsymbol{F} \in \mathcal{B}(\boldsymbol{S})$ to the two sets $\boldsymbol{Q}$ and $\boldsymbol{\Delta}$ in $\boldsymbol{V}$ with a cyclic control sequence, i.e., by alternating between the two sets, yields the following algorithm.

**Algorithm 5.9.1 Successive Generalized Projections (SGP) for Two Sets in the Product Space**

**Step 0:** (Initialization.) $\boldsymbol{b}^0 \in \boldsymbol{\Delta} \cap \boldsymbol{S}$ is arbitrary.
**Step 1:** (Iterative step.)

$$\boldsymbol{b}^{\nu+1} = \boldsymbol{P}_{\boldsymbol{\Delta}}^{\boldsymbol{F}}(\boldsymbol{P}_{\boldsymbol{Q}}^{\boldsymbol{F}}(\boldsymbol{b}^{\nu})), \ \nu = 0, 1, 2, \ldots \ . \tag{5.104}$$

Note that by (5.102) and (5.103) $\boldsymbol{b}^{\nu+1} \in \boldsymbol{S}$, so that $b^{\nu} \in S = \cap_{i=1}^{m} S_i$. Taking into account Lemmas 5.9.2 and 5.9.3 and that $\boldsymbol{b}^{\nu} = \boldsymbol{\delta}(b^{\nu})$ with $b^{\nu} \in \mathbb{R}^n$, we obtain from Algorithm 5.9.1 the following algorithm in $\mathbb{R}^n$.

**Algorithm 5.9.2 The Simultaneous Multiprojections Algorithm**

**Step 0:** (Initialization.) $b^0 \in S$ is arbitrary and $w \in \mathbb{R}^m$ is an arbitrary positive weight vector.
**Step 1:** (Iterative step.) Given $b^{\nu}$ calculate the projections

$$v^{i,\nu+1} = P_{Q_i}^{f_i}(b^{\nu}), \ i = 1, 2, \ldots, m, \tag{5.105}$$

then solve the equation

$$\sum_{i=1}^{m} w_i \nabla f_i(b^{\nu+1}) = \sum_{i=1}^{m} w_i \nabla f_i(v^{i,\nu+1}), \tag{5.106}$$

to find $b^{\nu+1}$. Set $b^{\nu} \leftarrow b^{\nu+1}$ and repeat the iterative step.

The next theorem establishes the convergence of any sequence generated by Algorithm 5.9.2.

**Theorem 5.9.1** *Let $f_i \in \mathcal{B}(S_i)$, $i \in I$, be a family of Bregman functions and let a CFP in $\mathbb{R}^n$ be given by a finite family of closed convex sets $Q_i \subseteq \mathbb{R}^n$, $i \in I$, such that $Q \doteq \cap_{i \in I} Q_i$ and $Q \cap \overline{S} \neq \emptyset$. Assume that for each $i$, $f_i$ is zone consistent with respect to $Q_i$, and that Assumption 5.9.1 holds. Under these conditions any sequence $\{b^{\nu}\}_{\nu=0}^{\infty}$ generated by Algorithm 5.9.2 converges to a point $b^* \in Q \cap \overline{S}$.*

**Proof** Assumption 5.9.1 and the zone consistency of each $f_i$ with respect to $Q_i$ guarantee (5.102) and (5.103), respectively. Likewise, $Q \cap \overline{S} \neq \emptyset$ yields $\boldsymbol{Q} \cap \boldsymbol{\Delta} \cap \overline{\boldsymbol{S}} \neq \emptyset$. Thus the convergence of iterates produced by (5.104)

$$\lim_{\nu \to \infty} \boldsymbol{b}^\nu = \boldsymbol{b}^* \in \boldsymbol{Q} \cap \boldsymbol{\Delta} \tag{5.107}$$

follows from Theorem 5.8.1 and the desired result follows from (5.85). ■

By choosing $f_i(x) = \frac{1}{2} \parallel x \parallel^2$ with zones $S_i = \mathbb{R}^n$ for all $i \in I$, we have $\nabla f_i(x) = x$ and $P_{Q_i}^{f_i}(b^\nu)$ becomes the orthogonal projection of $b^\nu$ onto $Q_i$. In this case Algorithm 5.9.2 is precisely the fully simultaneous projections method of Cimmino. It is identical to the special case of the BIP method (Algorithm 5.6.1) where $\lambda_\nu = 1$ for all $\nu \geq 0$, i.e., unity relaxation, and all weight vectors of BIP are $w^\nu = w$, for all $\nu \geq 0$. Other choices of Bregman functions $f_i \in \mathcal{B}(S_i)$ give rise to new simultaneous algorithms derived from Algorithm 5.9.2.

An important problem that Algorithm 5.9.2 helps to solve is the following *split* feasibility problem,

$$\text{find } x \in C \subseteq \mathbb{R}^n \text{ such that } Ax \in Q \subseteq \mathbb{R}^n, \tag{5.108}$$

where $C$ and $Q$ are closed convex sets and $A : \mathbb{R}^n \to \mathbb{R}^n$ is a linear bijection (i.e., one-to-one and onto mapping) represented by a full-rank matrix also denoted by $A$. Such problems arise in signal processing, specifically in phase retrieval and other image restoration or recovery problems.

To consider this problem let,

$$A(C) \doteq \{y \in \mathbb{R}^n \mid y = Ax, \ x \in C\}. \tag{5.109}$$

This set is easily seen to be convex. Based on the fact that $A$ is a bijection it also follows that $A(C)$ is closed (see, e.g., Taylor 1967, Theorem 4.2-E). Assuming that $Q \cap A(C) \neq \emptyset$, we can immediately apply the SOP method (Algorithm 5.2.1) to the two sets $Q$ and $A(C)$. For unity relaxation, the iterative step would have the form

$$x^{\nu+1} = A^{-1} P_Q(P_{A(C)}(Ax^\nu)), \quad \nu \geq 0, \tag{5.110}$$

where $A^{-1}$ is the inverse matrix of $A$. The orthogonal projection $P_{A(C)}(Ax^\nu)$ can be replaced by the minimization problem $\min_{x \in C} \parallel x - x^\nu \parallel_{A^T A}$, where $\parallel x \parallel_{A^T A} \doteq \langle x, A^T A x\rangle^{1/2}$. Therefore,

$$P_{A(C)}(Ax^\nu) = A\hat{P}_C(x^\nu), \tag{5.111}$$

where $\hat{P}$ denotes projection with respect to the distance defined by the norm $\parallel \cdot \parallel_{A^T A}$. Such oblique projections are special cases of the generalized projections arising from generalized distances $D_f$. We may thus rewrite (5.110) as

$$x^{\nu+1} = A^{-1}P_Q(A\hat{P}_C(x^\nu)), \quad \nu \geq 0, \tag{5.112}$$

and convergence would still be guaranteed from the convergence of the SOP method.

However, there are practical difficulties related to the nature of the set $C$ in some applications such as image recovery. For example, the set $A(C)$ is not simple to characterize, thus to project onto; and it is relatively simple to project onto $C$ orthogonally but performing the nonorthogonal projection $\hat{P}_C$ is not simple. The first difficulty is resolved by replacing projection onto $A(C)$ with projection onto $C$. But to resolve the second one would need to replace $\hat{P}_C$ with the orthogonal projection $P_C$. If this is done in (5.112) it would result in $P_{A(C)}$ being replaced by a nonorthogonal projection such as $\overline{P}_{A(C)}$ in (5.110). So, all in all, we see that the situation calls for executing within the same iterative process two different projections. While there is no sequential algorithm to accommodate this need, the simultaneous multiprojections algorithm (Algorithm 5.9.2) can, and actually has been used to, solve such problems.

## 5.10 Automatic Relaxation for Linear Interval Feasibility Problems

We return now to the linear case and consider the problem of solving iteratively large and possibly sparse systems of interval linear inequalities of the form

$$\alpha_i \leq \langle a^i, x\rangle \leq \beta_i, \quad i \in I \doteq \{1, 2, \ldots, m\}, \tag{5.113}$$

where $a^i \in \mathbb{R}^n$ is given, for all $i$, $\alpha = (\alpha_i) \in \mathbb{R}^m$, and $\beta = (\beta_i) \in \mathbb{R}^m$ are given too. Assuming that the system is feasible, an $x^* \in \mathbb{R}^n$ which solves (5.113) is required.

Geometrically the system represents $m$ nonempty hyperslabs in $\mathbb{R}^n$, each of which is the nonempty intersection of a pair of half-spaces. If we are willing to ignore the slab-structure of the problem it could be addressed as a system of $2m$ linear one-sided inequalities and solved by the relaxation method of Agmon, Motzkin, and Schoenberg (AMS, Algorithm 5.4.2). The method we study next however, takes advantage of the interval structure of the problem by handling in every iterative step a pair of inequalities. The method also realizes a certain relaxation principle in an automatic manner. External relaxation parameters are available on top of the built-in relaxation principle.

As will be seen later, in order to prove convergence the linear problem (5.113) is first represented as a nonlinear problem. This can be done simply and, once we overcome the barrier of transforming a linear problem to a nonlinear one, the proof follows by applying the method of cyclic subgradient projections (CSP, Algorithm 5.3.1).

For every hyperslab of the system (5.113) denote by

$$\overline{H}_i \doteq \{x \in \mathbb{R}^n \mid \langle a^i, x\rangle = \beta_i\} \tag{5.114}$$

and by

$$\underline{H}_i \doteq \{x \in \mathbb{R}^n \mid \langle a^i, x\rangle = \alpha_i\} \tag{5.115}$$

its bounding hyperplanes. The *median hyperplane* will be

$$H_i \doteq \{x \in \mathbb{R}^n \mid \langle a^i, x\rangle = \frac{1}{2}(\alpha_i + \beta_i)\}, \tag{5.116}$$

and the *half-width* $\theta_i$ of the hyperslab is

$$\theta_i = \frac{\beta_i - \alpha_i}{2 \parallel a^i \parallel}. \tag{5.117}$$

The signed distance of a point $u \in \mathbb{R}^n$ from the $i$th median hyperplane $H_i$ is given by

$$\hat{d}(u, H_i) = \frac{\langle a^i, u\rangle - \frac{1}{2}(\alpha_i + \beta_i)}{\parallel a^i \parallel} \tag{5.118}$$

and we will denote

$$\hat{d}_i = d(x^\nu, H_i). \tag{5.119}$$

The absolute values of these quantities measure the distance to the hyperplane $H_i$ but the quantities themselves are positive or negative depending on which side of the hyperplane the point is located.

The Automatic Relaxation Method (ARM) is as follows.

**Algorithm 5.10.1 Automatic Relaxation Method (ARM)**

**Step 0:** (Initialization.) $x^0 \in \mathbb{R}^n$ is arbitrary.
**Step 1:** (Iterative step.) Given $x^\nu$ calculate

$$x^{\nu+1} = \begin{cases} x^\nu, & \text{if } |\hat{d}_i| \leq \theta_i, \\ x^\nu - \dfrac{\lambda_\nu}{2}\left(\dfrac{\hat{d}_i^2 - \theta_i^2}{\hat{d}_i}\right)\dfrac{a^i}{\parallel a^i \parallel}, & \text{otherwise.} \end{cases} \tag{5.120}$$

**Relaxation parameters:** $\{\lambda_\nu\}_{\nu=0}^\infty$ are confined, for all $\nu \geq 0$, to

$$\epsilon_1 \leq \lambda_\nu \leq 2 - \epsilon_2, \qquad \text{for some } \epsilon_1, \epsilon_2 > 0. \tag{5.121}$$

**Control:** The sequence $\{i(\nu)\}_{\nu=0}^\infty$ is almost cyclic on $I$.

Before proving convergence of the algorithm we discuss its structure. As can be seen from (5.120) this algorithm is a row-action method. This means that only the original data $a^i$, $\alpha_i$, $\beta_i$, $i = 1, 2, \ldots, m$, to which no changes are made, are used throughout the iterations. No matrix operations are involved and in a single iterative step access is required to only one row of the system (5.113). Also, the method works with slabs rather than with the individual bounding hyperplanes, and only $x^\nu$, the immediate predecessor of $x^{\nu+1}$, is needed at any iterative step. These features, together with its simplicity, facilitate the application of this algorithm to very large and sparse systems.

To understand the automatic relaxation built into the algorithm, assume for a moment that $\lambda_\nu = 2$, for some $\nu \geq 0$. Rewrite the active part, i.e., the second line of (5.120), with $i \equiv i(\nu)$, as

$$x^{\nu+1} = x^\nu - \hat{d}_i \frac{a_i}{\| a^i \|} + \frac{\theta_i^2}{\hat{d}_i} \frac{a_i}{\| a^i \|}, \tag{5.122}$$

and consider sequentially the two additive changes made in this step to $x^\nu$. The point $\hat{x}^\nu$ obtained from $x^\nu$ by

$$\hat{x}^\nu = x^\nu - \hat{d}_i \frac{a_i}{\| a^i \|} \tag{5.123}$$

is the orthogonal projection of $x^\nu$ onto the median hyperplane $H_{i(\nu)}$. Correcting $\hat{x}^\nu$ to obtain $x^{\nu+1}$ according to the rest of (5.122) represents a move along the line through $x^\nu$ and $\hat{x}^\nu$ back toward $x^\nu$. The extent of this move depends on the distance of $x^\nu$ from $H_{i(\nu)}$ relative to the half-width $\theta_{i(\nu)}$ of the present hyperslab. If $x^\nu$ is far from the slab then the *move back* correction is small, whereas for an $x^\nu$ near the slab the total outcome of (5.122) will lie closer to the nearest hyperplane $\overline{H}_i$ or $\underline{H}_i$. At any rate, the move back never results in $x^{\nu+1}$ outside the $i(\nu)$th slab. It has been argued (in studies examining use of the relaxation method for solving systems of inequalities) that a good relaxation strategy will start with a relaxation parameter close to one, and when the sequence $\{x^\nu\}$ approaches the feasible set the relaxation parameter should be gradually increased to two.

As can be seen from the explanation above, Algorithm 5.10.1 realizes such a philosophy. The relaxation parameters $\lambda_\nu$ in (5.120) may be used to further control the actual implementation, by consistently underrelaxing the whole process, for example. Caution should be exercised, however, so that the basic idea of automatic relaxation of the algorithm is not negated through their choice. Another algorithm with the capability to enforce a similar relaxation strategy is described next. It was originally developed to solve a fully discretized model of the problem of image reconstruction from projections (see Chapter 10).

The modeling of a real-world problem often leads to a system of linear

equations

$$\langle a^i, x\rangle = b_i, \quad i = 1, 2, \ldots, m. \tag{5.124}$$

As mentioned earlier, such a system might be underdetermined due to lack of information or because it is self-contradictory (inconsistent). Measurement inaccuracies, noise corruption of data, and discretization errors in the model might all lead to the idea that rather than aiming at an exact algebraic solution of (5.124) (which might not even exist), it would be advisable to prescribe tolerances $\epsilon_i > 0$ and replace (5.124) with a system of interval linear inequalities such as (5.113) with

$$\beta_i = b_i + \epsilon_i \quad \text{and} \quad \alpha_i = b_i - \epsilon_i. \tag{5.125}$$

(See Chapter 13 for a discussion of different modeling approaches for dealing with incomplete, erroneous, or noisy data.) A special-purpose relaxation method for such a problem is the following algorithm, called Algebraic Reconstruction Technique 3 (ART3).

### Algorithm 5.10.2 Algebraic Reconstruction Technique 3 (ART3)

**Step 0:** (Initialization.) $x^0 \in \mathbb{R}^n$ is arbitrary.
**Step 1:** (Iterative step.) Given $x^\nu$ calculate

$$x^{\nu+1} = x^\nu - s_\nu \frac{a^i}{\| a^i \|}, \tag{5.126}$$

where

$$s_\nu = \begin{cases} 0, & \text{if } |\hat{d}_i| \le \theta_i, \\ 2(\hat{d}_i - t_i\theta_i), & \text{if } \theta_i < |\hat{d}_i| < 2\theta_i, \\ \hat{d}_i, & \text{if } 2\theta_i \le |\hat{d}_i|, \end{cases} \tag{5.127}$$

and $t_i \doteq \text{sign } \hat{d}_i$ is the sign of $\hat{d}_i$, i.e., $+1$ or $-1$.

**Relaxation parameters:** None.

**Control:** Cyclic, i.e., $i(\nu) = \nu (\text{mod } m) + 1$.

Let us explain the geometry of these iterations. In addition to the basic hyperslab with half-width $\theta_i$ around the median hyperplane $H_i$, we use an enveloping hyperslab. The half-width of the enveloping hyperslab is $2\theta_i$ (see Figure 5.10) and it serves to enforce in ART3 the following relaxation strategy. When the current iterate $x^\nu$ is far, i.e., $2\theta_i \le |\hat{d}_i|$, which means that it is outside the enveloping hyperslab, then $x^{\nu+1}$ is obtained by projecting $x^\nu$ right onto the median hyperplane $H_i$. When $x^\nu$ is near, i.e., $\theta_i < |\hat{d}_i| < 2\theta_i$ which means that it lies within the enveloping hyperslab but outside the basic hyperslab, then $x^{\nu+1}$ is the reflection of $x^\nu$ with respect to the bounding hyperplane of the basic hyperslab. Finally, if $x^\nu$ is inside the basic hyperslab, i.e., $|\hat{d}_i| \leq \theta_i$, then no move is made and $x^{\nu+1} = x^\nu$. This is a strategy similar to the one used by ARM (Algorithm 5.10.1) but it is enforced in a stepwise form in ART3 rather than in the automatic and gradual manner of ARM.

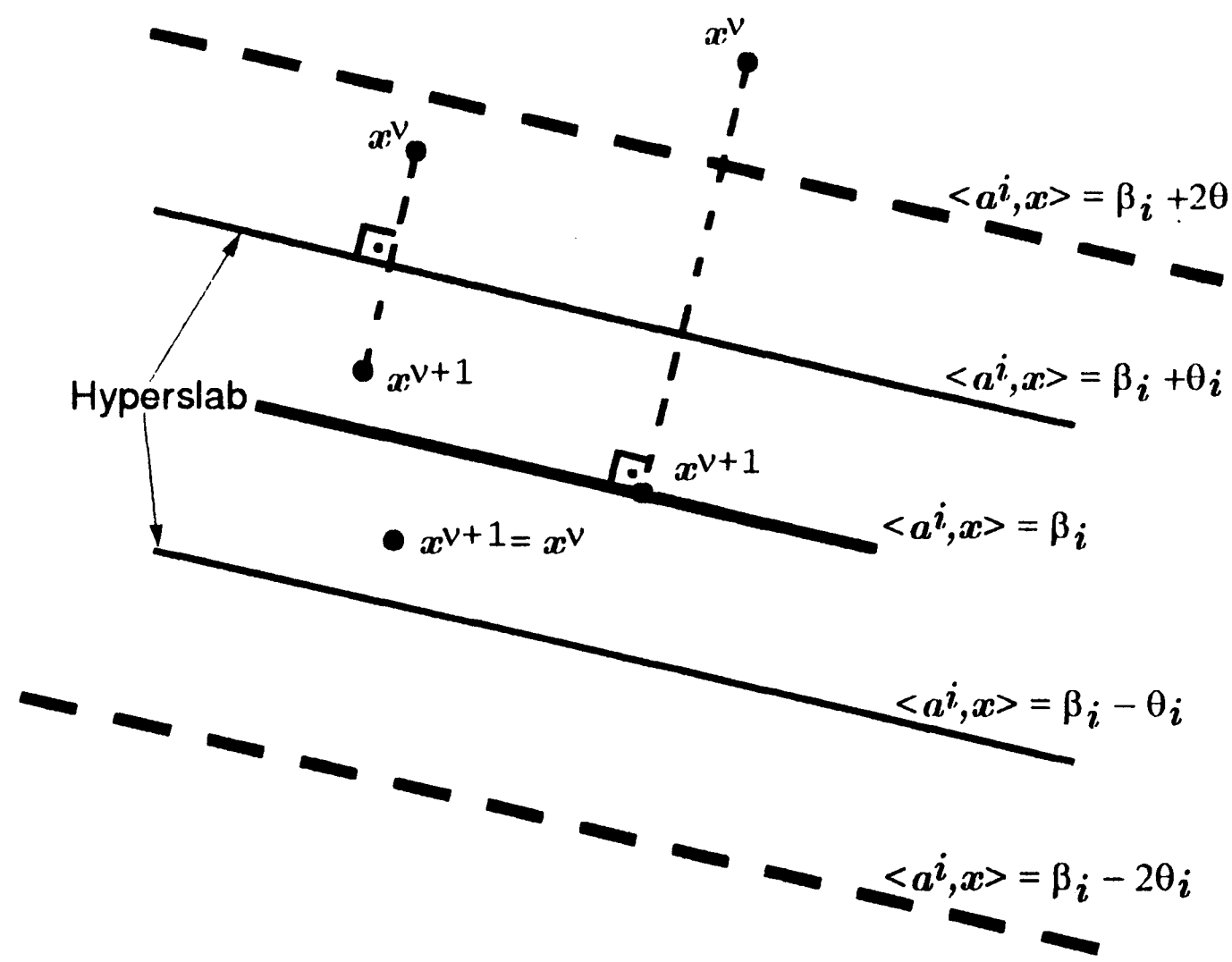

**Figure 5.10** The ART3 algorithm for linear intervals: all possible cases.

It is interesting to compare the actions of one iterative step in relation to some of the methods studied so far. We do this in Figure 5.11 where we describe the quantity $\langle a, x^{\nu+1}\rangle - \langle a, x^\nu\rangle$, as a function of $\langle a, x^\nu\rangle$. The particular hyperplane under consideration is defined by $\langle a, x\rangle = \beta$ and the half-width of the slab is $\theta$. The four methods compared here are: the method of Kaczmarz (ART), Algorithm 5.4.3; the relaxation method of Agmon, Motzkin, and Schoenberg (AMS), Algorithm 5.4.2; the algebraic reconstruction technique (ART3), Algorithm 5.10.2; and the automatic relaxation method (ARM), Algorithm 5.10.1. As seen in the figure, ART3 implements the special relaxation strategy in a stepwise manner whereas ARM does so gradually. This figure shows the behavior of the four algorithms by considering the dependence of $\langle a, x^{\nu+1}\rangle - \langle a, x^\nu\rangle$ on $\langle a, x^\nu\rangle$, for a given $\nu$. The line representing ART goes through the point $(\beta, 0)$ and its slope is minus one. This means that if $\langle a, x^\nu\rangle = \beta$ then $\langle a, x^{\nu+1}\rangle - \langle a, x^\nu\rangle = 0$, i.e., $x^{\nu+1} = x^\nu$. For any point $x^\nu$ which is not on the hyperplane $\langle a, x\rangle = \beta$ we have $\langle a, x^\nu\rangle \neq \beta$ and from the graph for ART we then get

$$\frac{\langle a, x^{\nu+1}\rangle - \langle a, x^\nu\rangle}{\beta - \langle a, x^\nu\rangle} = 1, \tag{5.128}$$

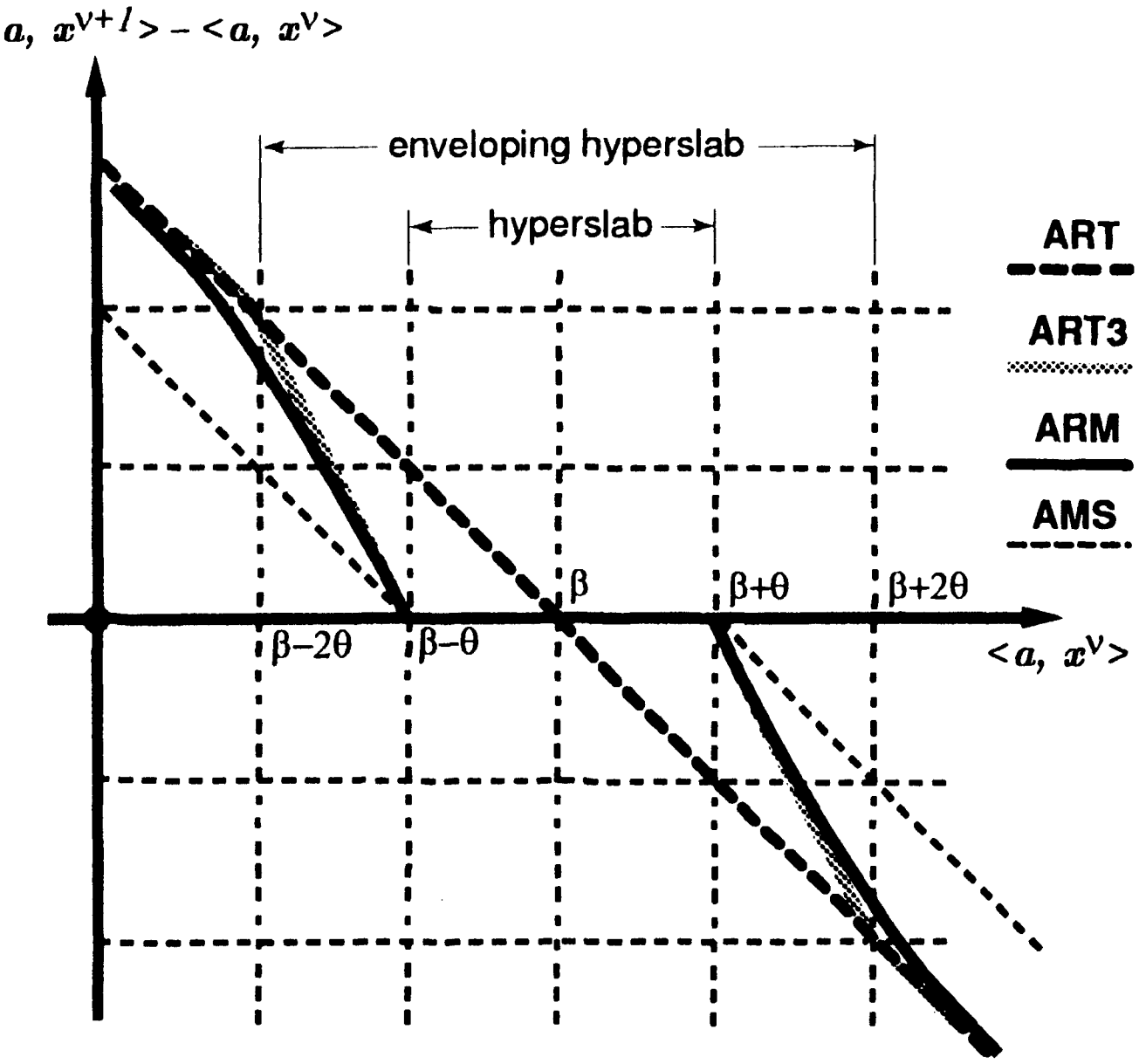

**Figure 5.11** A graphical comparison of four relaxation methods: Kaczmarz's method (ART), the relaxation method of Agmon, Motzkin, and Schoenberg (AMS), the algebraic reconstruction technique (ART3), and the automatic relaxation method (ARM).

or equivalently, $\langle a, x^{\nu+1}\rangle = \beta$. The line representing AMS is horizontal at $\langle a, x^{\nu+1}\rangle - \langle a, x^{\nu}\rangle = 0$, for all $\beta - \theta \leq \langle a, x^{\nu}\rangle \leq \beta + \theta$. This shows that, for the interval feasibility problem, if $x^{\nu}$ is within the hyperslab then $x^{\nu+1} = x^{\nu}$. For $\langle a, x^{\nu}\rangle$ outside the interval $[\beta - \theta, \beta + \theta]$ the slope of the line is again minus one showing that points $x^{\nu}$ which are outside the hyperslab get projected orthogonally onto the closest bounding hyperplane of the hyperslab.

ART3 coincides with AMS when $\langle a, x^{\nu}\rangle$ is within the interval $[\beta - \theta, \beta + \theta]$. The slope of its line which is minus one when $\langle a, x^{\nu}\rangle \leq \beta - 2\theta$ or $\beta + 2\theta \leq \langle a, x^{\nu}\rangle$, shows that when $x^{\nu}$ is outside the enveloping hyperplane then $x^{\nu+1}$ is the projection of $x^{\nu}$ onto $\langle a, x\rangle = \beta$. In the enveloping hyperplane the slope of the line of ART3 is precisely minus two, showing that if $x^{\nu}$ lies there then it is reflected, with respect to the nearest bounding hyperplane of the hyperslab, to obtain $x^{\nu+1}$. Finally, the graph of ARM shows how the relaxation strategy is achieved in a continuous manner for all $x^{\nu}$ that are outside the hyperslab.

We now turn to analyzing the convergence of ARM (Algorithm 5.10.1).

The proof of convergence is based on first describing the problem as a convex feasibility problem and then applying to it the CSP algorithm for solving convex inequalities.

**Theorem 5.10.1** *Any sequence* $\{x^\nu\}_{\nu=0}^{\infty}$ *produced by Algorithm 5.10.1 converges to a solution of the interval linear inequalities problem (5.113) provided that the problem is feasible.*

**Proof** For every $i = 1, 2, \ldots, m$, construct the function $f_i : \mathbb{R}^n \to \mathbb{R}$ by

$$\begin{aligned} f_i(x) &= (\langle a^i, x\rangle - \beta_i)(\langle a^i, x\rangle - \alpha_i) \\ &= \langle a^i, x\rangle^2 - (\alpha_i + \beta_i)\langle a^i, x\rangle + \alpha_i\beta_i. \end{aligned} \tag{5.129}$$

It is clear that, as long as $\alpha_i \leq \beta_i$,

$$\{x \in \mathbb{R}^n \mid \alpha_i \leq \langle a^i, x\rangle \leq \beta_i\} = \{x \in \mathbb{R}^n \mid f_i(x) \leq 0\}. \tag{5.130}$$

It follows that, regardless of the values of $a^i$, $\alpha_i$, and $\beta_i$, $f_i(x)$ is always a convex function on $\mathbb{R}^n$. The gradient of $f_i(x)$ is

$$\nabla f_i(x) = (2\langle a^i, x\rangle - (\alpha_i + \beta_i))a^i. \tag{5.131}$$

Thus, problem (5.113) is identical to the convex feasibility problem:

> Find an $x \in \mathbb{R}^n$ such that $f_i(x) \leq 0$, $i = 1, 2, \ldots, m$, where $f_i(x)$ are given by (5.129).

Applying the CSP method (Algorithm 5.3.1) to this convex feasibility problem we obtain the iterative step, in case $x^{\nu+1} \neq x^\nu$,

$$x^{\nu+1} = x^\nu - \lambda_\nu \frac{(\langle a^i, x^\nu\rangle^2 - (\alpha_i + \beta_i)\langle a^i, x^\nu\rangle + \alpha_i\beta_i)}{(2\langle a^i, x^\nu\rangle - (\alpha_i + \beta_i)) \parallel a^i \parallel^2} a^i, \tag{5.132}$$

which is the same as (5.120). The conditions of Theorem 5.3.1 hold and therefore the desired convergence follows. ∎

## 5.11 Notes and References

Iterative methods for the convex feasibility problem have been utilized for many years and successfully applied in various areas of mathematics and the physical and social sciences. Their study and use continue with vigor. Recent reviews and comprehensive unifying studies include Combettes (1993), Kiwiel (1995b), Cegielski (1993), and the systematic discussions of the topic in Hilbert space of Bauschke and Borwein (1996), and Bauschke, Borwein, and Lewis (1994). Earlier reviews of some of the methods discussed in this chapter can be found in Censor (1981, 1984). Much material on projections onto convex sets in Hilbert space can be found in Zarantonello (1971). Some unifying studies on projection and gradient methods for the convex feasibility problem are: Ottavy (1988), Butnariu and Mehrez (1992), Combettes and Trussell (1990), Barakat and Newsam (1985), and Combettes (1996b). It is beyond the scope of this section to describe in detail the vast literature on this topic. The references listed in the above-mentioned papers represent a fair cross section of the research done to date. Consult Section 10.7 for additional notes and references related to methods discussed here.

**5.1** Some of the control sequences mentioned here appear under different names in the literature. An early reference to the use of the almost cyclic control is in Herman and Lent (1978). The proof that every cyclic control is an approximately remotest set control can be found in Gubin, Polyak, and Raik (1967, Lemma 4). Another interesting control sequence not mentioned here is the *threshold control* of McCormick (1977). Repetitive control is found in the literature as random, chaotic, unrestricted, free-ordering, or free-steering.

**5.2** The method of successive orthogonal projection (SOP) is presented by Gubin, Polyak, and Raik (1967). Von Neumann (1950), whose article is a reprint of mimeographed lecture notes first distributed in 1933, treated alternating orthogonal projections between two closed linear subspaces. The von Neumann convergence theorem was extended to more than two subspaces by Halperin (1962) in Hilbert space, and by Bruck and Reich (1977) for linear and nonlinear projections in uniformly convex Banach spaces. The SOP method became popular in recent literature on image recovery, where it is also known as the POCS method. An early application of POCS to signal processing appeared in Lent and Tuy (1981). It was further studied by Youla (1978, 1987), Youla and Webb (1982), Sezan and Stark (1982, 1987), and Levi and Stark (1984). A self-contained analysis of the convergence of the SOP algorithm appears in Gubin, Polyak, and Raik (1967). They also treat the inconsistent case when $Q = \cap_{i \in I} Q_i = \emptyset$. Weak convergence results in Hilbert space under repetitive controls can be found, e.g., in Amemiya and Ando (1965), Dye, Khamsi, and Reich (1991), and Dye and Reich (1992).

**5.3** Material on subgradients of convex functions can be found in texts on convex analysis; see Rockafellar (1970) or Hiriart-Urruty and Lemaréchal (1993). The CSP method seems to have been first published by Erëmin (1966) and independently by Raik (1967). Censor and Lent (1982) rediscovered it. Censor et al. (1979) proposed (though without proving convergence) and experimented with an extension of the CSP method to the reverse-convex feasibility problem, the latter being a feasibility problem with a family of sets $\{Q_i\}$ some of which are convex and the rest are complements of convex sets. See also in this regard the remarks made in Section 5.4 on Algorithm 5.4.4. This non-linear extension of Kaczmarz's algorithm is due to McCormick (1977). An indispensable source for the closely related subgradient method for nondifferentiable minimization is found in Shor (1985). De Pierro and Iusem (1988) proposed a modified CSP method which, under mild conditions, converges finitely. Iusem and Moledo (1986) devised a simultaneous version of the CSP method.

**5.4** The remotest set controlled method of subgradient projections for the convex feasibility problem was derived by Polyak (1969). Oettli's algorithmic scheme can be found in Oettli (1975). Similar schemes were investigated in the papers by Erëmin (1965, 1966, 1970), and Erëmin and Mazurov (1967). See also Blum and Oettli (1975, Chapter 6). The linear feasibility problem of finding a solution of a system of linear inequalities or of solving a system of linear equations was dealt with separately in many papers. The AMS method was initially studied by Agmon (1954) and by Motzkin and Schoenberg (1954). Goffin (1977, 1980, 1981, 1982) studied the AMS method extensively regarding its rate of convergence, certain condition numbers, computational complexity, and its relationship to other methods, esp. the subgradient method for nondifferentiable optimization. Mandel (1984) studies further the cyclic AMS method and shows how bounds on convergence can be improved if strong underrelaxation, i.e., very small relaxation parameters, is used. Telgen (1982) studies finite convergence of the AMS method.

Kaczmarz's algorithm appeared in his original paper (1937). A textbook account of the algorithm, which also discusses earlier and classical methods such as Southwell's method, Gauss-Seidel's method, and Successive Over-Relaxation (SOR) appears in Gastinel (1970). Consult also Householder (1975) for a unified general treatment of projection methods for the solution of linear systems. The rediscovery of this method, under the name ART (Algebraic Reconstruction Technique), by Gordon, Bender, and Herman (1970) made Kaczmarz's algorithm a cornerstone method in the field of image reconstruction from projections (see Chapter 10). Tanabe (1971), investigated Kaczmarz's method for both singular and nonsingular systems and later put the method within the framework of linear iterative processes (1974). The method was also studied by Marti (1979), Trummer (1981, 1983, 1984), Ansorge (1984), and many others. We provide several further notes and references on Kaczmarz's method in Chapter 10 where ART is discussed. A nonlinear extension of the method appears in McCormick (1977); see also the related work of Martínez and De Sampaio (1986). García-Palomares (1995) studies an Armijo-Newton-like procedure that locates a feasible point of a nonempty system of nonlinear, not necessarily convex, inequalities in a finite number of steps.

Spingarn (1985) proposes a primal-dual algorithm for solving linear inequalities that is similar to the AMS method, but improves its behavior through additional dual variables. His algorithm is based on the more general method of partial inverses of monotone operators that he developed earlier (1983).

**5.5** Separation and support of convex sets are basic notions in convex analysis. Further information can be found in Rockafellar (1970), Avriel (1976), Bazaraa and Shetty (1979), and Luenberger (1969).

The $(\delta, \eta)$-Algorithm and its realizations were first proposed by Aharoni, Berman, and Censor (1983). García-Palomares (1993) studies mathematically and experimentally, in a multiprocessor environment, projected aggregation methods (PAM). These methods are shown there to be, under appropriate conditions, particular instances of the $(\delta, \eta)$-Algorithm.

**5.6** The BIP algorithm was proposed by Aharoni and Censor (1989) in the Euclidean space setting. It was further studied in Hilbert space by Butnariu and Censor (1990, 1994). When the BIP algorithm uses in every iterative step all sets $Q_i$, $i \in I$, then it is a fully simultaneous algorithm structured after the simultaneous algorithm of Cimmino (1938) for solving a system of linear equations; see, e.g., Gastinel (1970). This fully simultaneous case of BIP was studied earlier by Auslender (1976) and by Iusem and De Pierro (1986). A treatment of such algorithms in Banach spaces is presented by Reich (1983). For linear inequalities, the fully simultaneous case of the block-AMS algorithm (Algorithm 5.6.2) was studied by Censor and Elfving (1982). Yang and Murty (1992) construct at each iterative step a surrogate constraint onto which the current iterate is projected. Kiwiel (1995b) also approaches algorithmic structures from such a viewpoint.

**5.7** The block-iterative $(\delta, \eta)$-Algorithm, which combines the features of BIP and the $(\delta, \eta)$-Algorithm was proposed by Flåm and Zowe (1990) where a complete analysis is given. The unified treatment in Bauschke and Borwein (1996) also focuses on this algorithm.

**5.8** The method of successive generalized projections (SGP) originates in the pioneering paper by Bregman (1967b). The result has been extended from the cyclic to any repetitive control, in the more general setting of paracontractions, by Censor and Reich (1996). The convex feasibility problem can be viewed as a special instance of the problem of finding a common fixed point of a finite family of nonexpansive operators; see, e.g., Tseng (1992). Much research has been done on this in various settings; see, for example, the studies of Elsner, Koltracht, and Neumann (1992), Censor and Reich (1996), Flåm (1995), and references therein. Alber and Butnariu (1996) studied the SGP method in reflexive Banach spaces. Butnariu and Flåm (1995) and Butnariu (1995) formulated the stochastic convex feasibility problem and proposed and studied the closely related *expected-projections method* in Hilbert space. Butnariu, Censor, and Reich (1996) proposed an algorithm of iterative averaging of entropic projections for the stochastic convex feasibility problem.

**5.9** The material of this section originates in Censor and Elfving (1994). The product space setup was first used by Pierra (1984). He used the equivalence (5.85) to derive Cimmino's fully simultaneous method (with equal weights) from the SOP method. For the special case $f_i(x) = \frac{1}{2} \parallel x \parallel^2$ with $S = \mathbb{R}^n$, for all $i \in I$, Lemmas 5.9.2 and 5.9.3 were proved by Pierra. The formulation of the split feasibility problem is motivated by the works of Kotzer, Rosen, and Shamir (1995) and of Kotzer, Cohen and Shamir (1995). In a recent sequence of papers, Crombez (1991, 1992, 1993a, 1993b, 1995a, 1995b, 1995c) studies the problem of image recovery, in a Hilbert space setting, via convex combinations of (orthogonal) projections. The Hilbertian convex feasibility problem in both the consistent and inconsistent case has been studied by Combettes (1994, 1995, 1996a). See also related work by Sloboda (1978).

**5.10** Goffin's relaxation principle, upon which the algorithms in this section are based, appeared in his PhD dissertation (1971). Herman (1975) first implemented Goffin's principle in the ART3 algorithm. The automatic relaxation method (ARM) was proposed by Censor (1985). Figure 5.11 was suggested in a private communication by Dr. Arnold Lent from AT&T Bell Laboratories in Holmdel, New Jersey (USA).

CHAPTER 6

# Iterative Algorithms for Linearly Constrained Optimization Problems

In this chapter we discuss mainly row-action methods. Their distinctive feature is that they are iterative procedures which, without making any changes to the original constraints matrix $A$, use the rows of $A$ one at a time. Such methods appear in several areas of applications, and have demonstrated effectiveness for problems with large and sparse matrices that do not have any detectable or usable sparsity pattern, apart from a high degree of sparsity.

We present row-action methods for handling huge and sparse systems (with numbers of unknowns, equations, or inequalities in the range of hundreds of thousands), and examine linear and nonlinear systems of equalities or inequalities.

In describing many different row-action methods together, in a coherent presentation, we have two aims in mind: To facilitate cross-fertilization between various fields of applications where row-action methods are—or could be—used, and to pave the way for better insight into their mathematical nature. The row-action methods for linearly constrained optimization are treated in a unified manner within the framework of generalized Bregman distances and generalized projections (see Chapter 2).

In Section 6.1 we present several approaches and solution concepts as well as the special environment that may affect the user's choice of approach. Row-action methods are defined in Section 6.2. In Section 6.3 we present Bregman's algorithm for linear equality constrained optimization of Bregman functions. The motivation for the construction of this algorithm is given and its convergence to the desired solution is proven. The special modification of this algorithm for interval constrained problems is presented and analyzed in Section 6.4. The specialization of these algorithms to the Bregman function $f(x) = \frac{1}{2} \parallel x \parallel^2$ yields several row-action algorithms for norm minimization, which are discussed in Section 6.5. When the negative Shannon entropy is used, i.e., $f(x) = \sum_{j=1}^{n} x_j \log x_j$, then we obtain a variety of row-action entropy optimization algorithms, as discussed in Section 6.6. The block-MART algorithm for entropy optimization over linear equality constraints is treated in detail in Section 6.7. Section 6.8 brings two additional features into the algorithms studied in Sections 6.3 and 6.4: the introduction of underrelaxation parameters, and the extension of the applicable family of objective functions to Bregman functions with singularities on the boundary of their zone. The so-called hybrid algorithm in Section 6.9 combines the overall structure of Bregman's algorithm with a method of deriving a closed-form formula that replaces the projection step, which otherwise needs to be calculated iteratively. Specific realizations of Shannon's, Burg's and Rényi's entropies are given. Notes and References are given in Section 6.10.

## 6.1 The Problem, Solution Concepts, and the Special Environment

### 6.1.1 The problem

In many fields of applications, including those described in later chapters of this book, the modeling of a physical or mathematical problem leads to a system of linear equations

$$\langle a^i, x\rangle = b_i, \qquad i \in I, \tag{6.1}$$

where $I \doteq \{1,2,\ldots,m\}$, $a^i \in \mathbb{R}^n$, $a^i \neq 0$, $x \in \mathbb{R}^n$, $b_i \in \mathbb{R}$. We also write the system (6.1) as $Ax = b$ where $A$ denotes the $m \times n$ matrix whose $i$th row is $(a^i)^{\mathrm{T}}$.

This system is sometimes underdetermined due to lack of information; often it is greatly overdetermined in which case it may be inconsistent. When modeling real-world problems, we have to be prepared to deal with inconsistent and possibly ill-conditioned systems. Barring these difficulties we might have reason to believe that the exact algebraic solution of the system is less desirable, in terms of the original problem for which the system was set up, than another, differently defined "solution." This belief may be supported by evidence of measurement inaccuracy, noise corruption of data, errors introduced during discretization of the model, etc. In such cases, as well as in the situations described above, it is useful to turn to a different solution concept rather than aim at the exact algebraic solution of the system.

### 6.1.2 Approaches and solution concepts

We should try to use the information contained in the system (6.1) in a way that reflects our limited confidence in the equations, and there are several different ways of doing so. Each approach described below leads to a different solution concept, and the choice of solution concept for a particular application belongs to the user and should be made with reference to the specific problem at hand.

### The Feasibility Approach

Here one seeks a point $x$ that lies within a specific vicinity of all hyperplanes defined by the equations (6.1), by prescribing some positive tolerances $\epsilon_i$, for all $i \in I$, and aiming at any solution $x \in \mathbb{R}^n$ that lies in the intersection of all hyperslabs

$$b_i - \epsilon_i \leq \langle a^i, x \rangle \leq b_i + \epsilon_i, \qquad i \in I. \tag{6.2}$$

Tolerances $\epsilon_i$ have to be chosen so that the feasible region is neither empty nor too large. Inequalities of the form (6.2) are called *interval inequalities* and the associated feasibility problem is a linear interval feasibility problem. If we use a solution method that does not take advantage of the fact that the inequalities come in pairs, we regard the system (6.2) simply as a system of one-sided linear inequalities with twice as many inequalities, obtained by multiplying one of the two sets of inequalities by minus one. The linear interval feasibility problem can be handled by methods discussed in Chapter 5 (see, in particular, Section 5.10).

### The Optimization Approach

Here we specify an objective function $f : \mathbb{R}^n \to \mathbb{R}$, according to which a particular element from the feasible region described by (6.2) is selected as the solution to the original real-world system. The resulting optimization problem

$$\text{Minimize} \quad f(x), \tag{6.3}$$

$$\text{s.t.} \qquad b_i - \epsilon_i \leq \langle a^i, x \rangle \leq b_i + \epsilon_i, \qquad i \in I, \tag{6.4}$$

is an interval constrained problem and one would prefer a solution method that benefits from the slab-structure of the constraints, as opposed to a one-sided linearly constrained optimization technique. If the problem does not have a unique solution, a secondary optimization criterion is sometimes used and another objective function $g(x)$ is optimized over the set of solutions of (6.3)–(6.4).

In both the feasibility approach and the optimization approach additional information about the real-world problem may result in additional inequalities that further restrict the feasible set. A common case is that of box constraints $l_j \leq x_j \leq u_j$, $j = 1, 2, \ldots, n$, which may reflect some prior information that the desired solution $x = (x_j)$ must be within the known bounds $l = (l_j)$ and $u = (u_j)$. Although these inequalities can be regarded as interval inequalities (with $x_j = \langle e^j, x \rangle$, where $e^j$ is the $j$th standard basis vector), a good solution method should handle them in a simpler way.

### The Regularization Approach

Here we consider the problem

$$\text{“solve”} \quad Ax = b \tag{6.5}$$

$$\text{s.t.} \quad x \in Q, \tag{6.6}$$

where $Ax = b$ is the system (6.1) and $Q$ represents additional constraints. This problem is treated by replacing it with

$$\text{Minimize} \quad f(x) + \mu g(Ax - b) \tag{6.7}$$

$$\text{s.t.} \quad x \in Q, \tag{6.8}$$

where $f : \mathbb{R}^n \rightarrow \mathbb{R}$ and $g : \mathbb{R}^m \rightarrow \mathbb{R}$ are some suitably chosen convex functions and $\mu$ is a user-specified parameter reflecting the relative importance attached to each summand. Observe that we speak here about *regularizing* the system of equations (6.1), which gives us the flexibility to choose both functions $f(x)$ and $g(x)$. The least squares regularization, where $f(x)$ and $g(x)$ are both of the form $\| x \|^2$ is by far the most common approach. Application of row-action methods to (6.7)–(6.8) in this case is possible but other choices of $f$ and $g$ can also be used, such as the *entropy regularization approach* where one or both $f(x)$ and $g(x)$ is of the form $\sum_j x_j \log x_j$.

### The Distance Reduction Approach

Here the goal is to reduce the distance between the data $b$ and the range of the matrix $A$, which is defined as $\mathcal{R}(A) \doteq \{y \in \mathbb{R}^n \mid y = Ax, \text{ for some } x \in \mathbb{R}^n\}$. An element $Ax$ of $\mathcal{R}(A)$ is sometimes referred to as the *pseudodata* of the problem. A generalized distance $D_f(x, y)$ is constructed from some chosen Bregman function $f \in \mathcal{B}(S)$, with zone $S$ (see Section 2.1), and the model problem (6.5)–(6.6) is replaced by the optimization problem

$$\text{Minimize} \quad D_f(b, Ax) \tag{6.9}$$

$$\text{s.t.} \quad x \in Q, \tag{6.10}$$

$$Ax \in S. \tag{6.11}$$

For $f(x) = \frac{1}{2} \| x \|^2$, with $S = \mathbb{R}^n$ and $Q = \mathbb{R}^n$, this approach leads to the well-known least squares solution method. If one chooses the negative Shannon's entropy as $f$ (see Example 2.1.2) then we have a Kullback-Leibler distance reduction problem to which the Expectation Maximization (EM) algorithm can be applied (see Section 10.2.1). If the problem is feasible, i.e., $b \in \mathcal{R}(A) \cap Q$, then a distance reduction algorithm will generate an exact solution of (6.5)–(6.6). Otherwise, it should find an $x^* \in Q$ that is closest to $\mathcal{R}(A)$ according to the specified generalized distance.

### 6.1.3 The special computational environment

An important factor of the mathematical problem involved in any of the approaches described above is the environment within which the problem is addressed. The environment is distinguished by the following properties.

*Dimensionality:* The system (6.1) is huge, e.g., $n \geq 10^5$, and $m$ even greater.

*Sparsity:* The matrix $A$ of the system (6.1) is sparse. For example, in a fully discretized model in image reconstruction from projections less than one percent of the entries are nonzero (see Chapter 10).

*Lack of structure:* One fails to recognize any structure in the distribution of nonzero entries throughout $A$. Alternatively, one might sometimes be able to detect some structure but be unable to take any computational advantage of it.

*Time restriction:* In some applications there is an inherent time restriction on the solution process. That is, a solution which is acceptable in terms of the original real-world problem is required within minutes of data collection. Image reconstruction from projections for medical diagnosis and planning under uncertainty for portfolio management are examples where solutions are needed with some urgency.

*Computation power restriction:* Sometimes solution methods are explicitly required to perform efficiently on special-purpose machines with low memory or other specifications. Image reconstruction equipment and parallel computers are examples of machines with such a special environment.

A prime example arises from the practical problem of computerized tomography (see Chapter 10), where some combination of all five properties usually describes the environment for the problem in any of the mentioned approaches. Problems arising in financial planning, especially in portfolio optimization, have a similar environment except that they exhibit a usable structure (see Chapters 8 and 13). Situations where some of these environmental properties hold also arise in various other fields.

In such an environment the use of general-purpose techniques is usually impractical and solution techniques have to be chosen to cope with the environment in a computationally feasible way.

## 6.2 Row-Action Methods

We consider the concept of row-action (RA) methods as a general framework within which various concrete algorithms fall. As will be seen below, the framework is broad enough to include various kinds of algorithms designed for a wide range of problems. Moreover, there are algorithms which are not formulated as RA methods, but lend themselves to a row-action implementation. One such example, which is discussed later, is Hildreth's algorithm (Algorithm 6.5.2). The matrix mentioned in the next definition is either the matrix $A$, if linear constraints are present, or the Jacobian matrix of first partial derivatives in the case of nonlinear constraints.

**Definition 6.2.1 (Row-Action Method)** *A row-action method is a sequential algorithm (in the classification of Section 1.3) with the following properties:*

*(i)* *no changes are made to the original matrix,*

*(ii)* *no operations are performed on the matrix as a whole,*

*(iii)* *at each iterative step access is required to only one row of the matrix,*

*(iv)* *at each iterative step when $x^{\nu+1}$ is calculated, the only iterate needed is the immediate predecessor $x^{\nu}$.*

We are particularly interested, for efficiency reasons, in RA methods that are further characterized by the fact that they present modest arithmetical demands. This requirement is not precisely defined and might well be considered a universal goal in numerical analysis, but we stress it here to make the point that in the special environment described above (particularly dimensionality), performing slightly more complex operations might quickly reduce the practicality of a method. Some of the methods described below comply less to this demand than others. For example, Hildreth's algorithm requires the storage and calculation of a sequence of dual vectors; Bregman's algorithm requires the execution of a generalized projection at each iterative step.

Finally, we consider another practical requirement that makes RA methods favorable. The term *row generation capability* refers to a situation where it is possible to avoid storing the matrix of system (6.1) explicitly, and instead to have both the nonzero entries of the $i$th row and their addresses generated from the experimental data each time anew. Such situations arise, for example, in the reconstruction of medical images from experimental data in computerized tomography. Property *(iii)* in Definition 6.2.1 makes RA methods favorable in such cases. In the absence of row generation capability one would, of course, use any of the existing efficient methods for storing a sparse matrix.

Row-action methods usually employ a control sequence (see Definition 5.1.1), the most important controls being the cyclic, the almost cyclic, the repetitive, and the remotest set controls. Intuition tells us that a control such as the remotest set control might lead to a faster convergence of the RA method to which it is applied. On the other hand, in the special environment we are in, any control other than cyclic or almost cyclic may violate the additional requirement of arithmetical simplicity, and thus make the RA method impractical. For this reason we emphasize below those two control sequences, although some RA methods can accommodate other controls as well.

## 6.3 Bregman's Algorithm for Inequality Constrained Problems

We begin the study of row-action methods by examining closely Bregman's algorithm for the linearly constrained minimization of a Bregman function $f \in \mathcal{B}(S)$ (see Definition 2.1.1). Consider the problem

$$\text{Minimize} \quad f(x) \tag{6.12}$$

$$\text{s.t.} \quad \langle a^i, x\rangle \le b_i, \qquad i \in I, \tag{6.13}$$

$$x \in \overline{S}. \tag{6.14}$$

Let $H_i \doteq \{x \mid \langle a^i, x\rangle = b_i\}$ and $Q_i \doteq \{x \mid \langle a^i, x\rangle \le b_i\}$; denote also $H = \cap_{i=1}^m H_i$, $Q \doteq \cap_{i=1}^m Q_i$, and assume that $Q \cap \overline{S} \neq \emptyset$. $A = (a_{ij})$ is the $m \times n$ matrix whose $i$th row is $(a^i)^{\mathrm{T}}$, and $b = (b_i) \in \mathbb{R}^m$. Assume that all $a^i \neq 0$. Assume that $f \in \mathcal{B}(S)$ is strongly zone consistent with respect to every $H_i$ (see Definition 2.2.1). Define the following sets:

$$Z \doteq \{x \in S \mid \exists \ z \in \mathbb{R}^m, \text{ such that } \nabla f(x) = -A^{\mathrm{T}} z\}, \tag{6.15}$$

$$Z_0 \doteq \{x \in S \mid \exists \ z \in \mathbb{R}^m_+, \text{ such that } \nabla f(x) = -A^{\mathrm{T}} z\}, \tag{6.16}$$

which are assumed nonempty.

Bregman's algorithm takes at each iteration a hyperplane $H_i$, projects onto it the current iterate according to the generalized distance constructed from the objective function $f(x)$, and computes the next iterate from this projection. The sequence of indices of the selected hyperplanes that governs the application of the method to a problem is the control sequence. The original Bregman algorithm dealt only with cyclic control. The less restrictive, almost cyclic control established here is necessary for the development of the interval programming method in the next section.

Before turning to the algorithm itself and its analysis, we motivate its structure. Here, we show how to develop a primal-dual algorithmic scheme for the following linearly constrained optimization problem

$$\text{Minimize} \quad f(x) \tag{6.17}$$

$$\text{s.t.} \quad Ax = b, \tag{6.18}$$

$$x \ge 0, \tag{6.19}$$

where $A$ is a given real $m \times n$ matrix, $b = (b_i) \in \mathbb{R}^m$. The reader should also refer to the treatment of primal-dual methods presented in Section 4.4. The Lagrangian of (6.17)–(6.19) with respect to the equality constraints is

$$L(x, \pi) = f(x) + \langle \pi,\ Ax - b \rangle, \tag{6.20}$$

where $\pi \in \mathbb{R}^m$ is the dual vector of Lagrange multipliers. The dual function $g : \mathbb{R}^m \to \mathbb{R}$ is

$$g(\pi) = \min_{x \in \mathbb{R}^n_+} L(x, \pi), \tag{6.21}$$

where $\mathbb{R}^n_+$ is the nonnegative orthant of $\mathbb{R}^n$. A necessary condition for minimization of (6.21) is $\nabla_x L(x, \pi) = 0$, which yields the Karush-Kuhn-Tucker conditions for our problem

$$\nabla f(x) = -A^{\mathrm{T}} \pi. \tag{6.22}$$

An iterative primal-dual algorithm may be derived by iterating both the primal $\{x^\nu\}_{\nu=0}^\infty$ and the dual $\{\pi^\nu\}_{\nu=0}^\infty$ vectors such that (6.22) holds for each iteration, i.e.,

$$\nabla f(x^\nu) = -A^{\mathrm{T}} \pi^\nu, \quad \nu \geq 0. \tag{6.23}$$

To do so the corrections to the dual vectors have to be prescribed. The Lagrange duality theorem states that

$$\min_{x \in \mathbb{R}^n_+ \cap H} f(x) = \sup_{\pi \in \mathbb{R}^m_+} g(\pi). \tag{6.24}$$

Therefore, dual corrections should at least entail dual ascent, i.e., guarantee that the sequence $\{g(\pi^\nu)\}_{\nu=0}^\infty$ is increasing. Dual correction can be represented by

$$\pi^{\nu+1} = \pi^\nu + v^\nu, \tag{6.25}$$

where $v^\nu \in \mathbb{R}^m$ is a *dual correction vector.* The algorithmic scheme obtained from (6.23) and (6.25) is

$$\nabla f(x^{\nu+1}) = -A^{\mathrm{T}} \pi^{\nu+1} = -A^{\mathrm{T}} \pi^\nu - A^{\mathrm{T}} v^\nu = \nabla f(x^\nu) - A^{\mathrm{T}} v^\nu. \tag{6.26}$$

If a decision is made to change only one component of $\pi^\nu$, i.e., the $i(\nu)$th component, at each iteration, then

$$v^\nu = \theta_\nu e^{i(\nu)}, \tag{6.27}$$

where $e^i$ is the $i$th standard basis vector in $\mathbb{R}^m$, and we get

$$\pi_i^{\nu+1} = \begin{cases} \pi_i^\nu, & i \neq i(\nu), \\ \pi_i^\nu + \theta_\nu, & i = i(\nu), \end{cases} \tag{6.28}$$

resulting in a Gauss-Seidel type algorithmic scheme.

The algorithm designer must also choose the size of the correction term $\theta_\nu$ and the index $i(\nu)$, specifying which component of $\pi^\nu$ will be changed at the $\nu$th iteration. Regardless of the actual $\theta_\nu$ and $i(\nu)$, the decision (6.27) leads to

$$\nabla f(x^{\nu+1}) = \nabla f(x^\nu) - \theta_\nu a^{i(\nu)}, \tag{6.29}$$

which is a row-action scheme in the sense of Definition 6.2.1. If $\theta_\nu$ is calculated so that $x^{\nu+1}$ satisfies the additional requirement

$$\langle x^{\nu+1}, a^{i(\nu)} \rangle = b_{i(\nu)}, \tag{6.30}$$

that is, the next iterate $x^{\nu+1}$ lies on the $i(\nu)$th hyperplane, and if $\{i(\nu)\}$ is an almost cyclic control sequence, then we obtain Bregman's method for convex programming with linear equality constraints. The solution $x^{\nu+1}$ of the system (6.29)–(6.30) is the generalized projection of the current primal vector $x^\nu$ onto the hyperplane $H_{i(\nu)}$ and the resulting $\theta_\nu$ is the generalized projection parameter (see Lemma 2.2.1).

The algorithm for linear inequality constraints also calculates this parameter $\theta_\nu$. However, before proceeding, it compares it with the $i(\nu)$th component of the current dual vector $\pi^\nu$ and uses the smaller of these two in the iteration formula (6.32). The algorithm is described below.

**Algorithm 6.3.1 Bregman's Algorithm for Linear Inequalities**

**Step 0:** (Initialization.) $x^0 \in Z_0$ is arbitrary, and $z^0$ is such that

$$\nabla f(x^0) = -A^{\mathrm{T}} z^0. \tag{6.31}$$

**Step 1:** (Iterative step.) Given $x^\nu$ and $z^\nu$, calculate $x^{\nu+1}$ and $z^{\nu+1}$ from

$$\begin{aligned} \nabla f(x^{\nu+1}) &= \nabla f(x^\nu) + c_\nu a^{i(\nu)}, \\ z^{\nu+1} &= z^\nu - c_\nu e^{i(\nu)}, \\ c_\nu &\doteq \min(z^\nu_{i(\nu)}, \theta_\nu), \end{aligned} \tag{6.32}$$

where $\theta_\nu$ is the parameter associated with the generalized projection of $x^\nu$ onto $H_{i(\nu)}$. We assume throughout that the representation of every hyperplane is fixed during the whole iteration process, so that the values of $\theta_\nu$ are well defined.

**Control:** The sequence $\{i(\nu)\}_{\nu=0}^{\infty}$ is almost cyclic on the index set $I$.

Recall that $\theta_\nu$, the parameter associated with the generalized projection $\hat{x}^{\nu+1}$ of $x^\nu$ onto $H_{i(\nu)}$, is obtained by solving the system

$$\begin{aligned} \nabla f(\hat{x}^{\nu+1}) &= \nabla f(x^\nu) + \theta_\nu a^{i(\nu)}, \\ \langle a^{i(\nu)}, \hat{x}^{\nu+1} \rangle &= b_{i(\nu)}, \end{aligned} \tag{6.33}$$

(see Lemma 2.2.1). The next lemma gives insight into the behavior of Algorithm 6.3.1 and justifies the geometric interpretation given later.

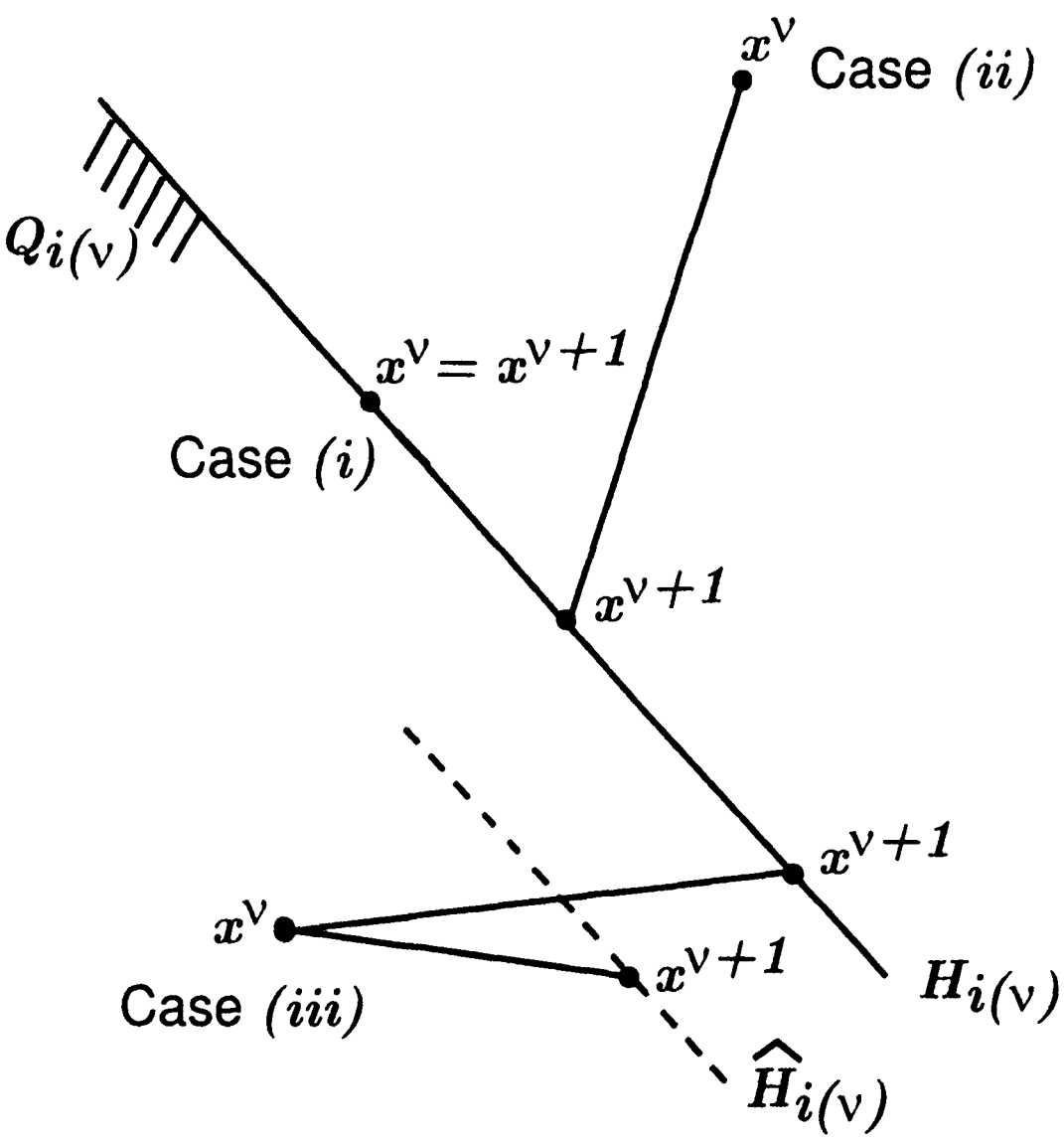

**Figure 6.1** Geometric interpretation of all possible cases in an iterative step of Bregman's algorithm (Algorithm 6.3.1).

**Lemma 6.3.1** *For $\nu = 0, 1, 2, \ldots$, all the primal iterates produced by Algorithm 6.3.1 satisfy $x^{\nu+1} \in Q_{i(\nu)}$. Moreover, if $x^\nu \notin \operatorname{int} Q_{i(\nu)}$, then $x^{\nu+1} \in H_{i(\nu)}$.*

**Proof** Let $x$ and $y$ be any two points in $S$. Then, the following identity is easily derived from (2.1)

$$D_f(x, y) + D_f(y, x) = \langle \nabla f(x) - \nabla f(y), x - y \rangle. \tag{6.34}$$

Let $\hat{x}^{\nu+1}$ be the projection of $x^\nu$ onto $H_{i(\nu)}$, and assume that (6.32) is solvable for $x^{\nu+1}$. Then, replacing $x$ with $x^{\nu+1}$ and replacing $y$ with $\hat{x}^{\nu+1}$ in (6.34), we get, by Lemma 2.1.1,

$$(c_\nu - \theta_\nu)\langle a^{i(\nu)}, x^{\nu+1} - \hat{x}^{\nu+1} \rangle \geq 0. \tag{6.35}$$

From (6.32) we also obtain $c_\nu \leq \theta_\nu$, and by definition $\langle a^{i(\nu)}, \hat{x}^{\nu+1} \rangle = b_{i(\nu)}$. Therefore,

$$\langle a^{i(\nu)}, x^{\nu+1} \rangle \leq b_{i(\nu)} \tag{6.36}$$

follows and proves the first part of the lemma.

To prove the second part, we look ahead to the proof of (6.37) where it is shown that $z^\nu \geq 0$ for all $\nu$. We assume that $x^\nu \notin \operatorname{int} Q_{i(\nu)}$, so that $b_{i(\nu)} - \langle a^{i(\nu)}, x^\nu \rangle \leq 0$, which by Lemma 2.2.2 implies that $\theta_\nu \leq 0$. But then, from (6.32) $c_\nu = \theta_\nu$ and $x^{\nu+1}$ is the generalized projection of $x^\nu$ onto $H_{i(\nu)}$. ■

For the special case $f(x) = \frac{1}{2} \| x \|^2$ (Example 2.1.1) and with the cyclic control, Algorithm 6.3.1 is precisely Hildreth's quadratic programming method (Algorithm 6.5.2). In this quadratic programming case, the iterative step (6.32) takes a very simple form. Generally however (6.32) calls for the solution of a system of $n+1$ nonlinear equations to determine $\theta_\nu$, at each iterative step. The effectiveness of the method chosen to solve this subproblem could determine the overall performance of the algorithm.

Figure 6.1 illustrates the algorithm. Three possible cases may occur in an iterative step when the $i(\nu)$th constraint is considered: if $x^\nu \in H_{i(\nu)}$, then $x^{\nu+1} = x^\nu$; if $x^\nu \notin Q_{i(\nu)}$, then $x^{\nu+1}$ is the projection of $x^\nu$ onto $H_{i(\nu)}$; if $x^\nu \in Q_{i(\nu)}$ then $x^{\nu+1}$ is either the projection of $x^\nu$ onto $H_{i(\nu)}$ or the projection of $x^\nu$ onto another hyperplane, say, $\hat{H}_{i(\nu)}$, which is parallel to $H_{i(\nu)}$ and lies entirely within $Q_{i(\nu)}$. All these possibilities are covered by (6.32) and are represented in the figure.

We now present a convergence theorem for Algorithm 6.3.1 as applied to the standard problem (6.12)–(6.14).

**Theorem 6.3.1** *Assume the following:*

*(i)* $f \in \mathcal{B}(S)$,

*(ii)* $f$ *is strongly zone consistent with respect to each* $H_i$, $i \in I$,

*(iii)* $\{i(\nu)\}_{\nu=0}^{\infty}$ *is almost cyclic with constant of almost cyclicality* $C$,

*(iv)* $Z_0 \neq \emptyset$.

*Then, any sequence* $\{x^\nu\}$ *produced by Algorithm 6.3.1 converges to the point* $x^*$*, which is the solution of (6.12)–(6.14).*

**Proof** The proof consists of the following steps.

**Step 1.** $x^\nu \in Z_0$, for all $\nu \geq 0$.

**Step 2.** Define the Lagrangian $L(x, z)$ of (6.12)–(6.14), and show that the limit $\lim_{\nu\to\infty} L(x^\nu, z^\nu)$ exists.

**Step 3.** $\lim_{\nu\to\infty} D_f(x^{\nu+1}, x^\nu) = 0$.

**Step 4.** $\{x^\nu\}$ is a bounded sequence.

**Step 5.** Any accumulation point of $\{x^\nu\}$ is feasible for (6.12)–(6.14).

**Step 6.** $x^*$ is the solution of (6.12)–(6.14).

We now proceed to prove each step.

**Step 1.** We show that

$$x^\nu \in Z_0, \quad \text{for all } \nu \geq 0. \tag{6.37}$$

The proof proceeds by induction. Observe that $x^0 \in Z_0$ according to the initialization of the algorithm. Assume that $x^\mu \in Z_0$, for all $\mu \leq \nu$; let $i(\nu)$ be the next index and abbreviate $i \equiv i(\nu)$. Then, from (6.32) and using $A^{\mathrm{T}}e^i = a^i$, we get

$$\nabla f(x^{\nu+1}) = -A^{\mathrm{T}}z^\nu + c_\nu a^i = -A^{\mathrm{T}}(z^\nu - c_\nu e^i) = -A^{\mathrm{T}}z^{\nu+1}, \tag{6.38}$$

showing that $x^{\nu+1} \in Z$. But $c_\nu \leq z_i^\nu$ by (6.32), and therefore, $z^\nu \geq 0$ which implies that $z^{\nu+1} \geq 0$, concluding the proof of (6.37).

**Step 2.** Define the Lagrangian of (6.12)–(6.14) as

$$L(x,z) \doteq f(x) + \langle z, Ax - b\rangle. \tag{6.39}$$

The values of the Lagrangian at successive iterates form an increasing sequence, i.e.,

$$d_\nu \doteq L(x^{\nu+1}, z^{\nu+1}) - L(x^\nu, z^\nu) \geq 0, \quad \text{for all } \nu \geq 0. \tag{6.40}$$

This can be shown as follows. From (6.39), (6.40), (2.1), and from

$$\langle z^\nu, Ax^\nu\rangle = \langle A^{\mathrm{T}}z^\nu, x^\nu\rangle = -\langle \nabla f(x^\nu), x^\nu\rangle, \tag{6.41}$$

which holds for all $\nu$ because of (6.38), we may write

$$d_\nu = D_f(x^{\nu+1}, x^\nu) - \langle z^{\nu+1} - z^\nu, b\rangle + \langle \nabla f(x^\nu) - \nabla f(x^{\nu+1}), x^{\nu+1}\rangle, \tag{6.42}$$

which yields, by (6.32),

$$d_\nu = D_f(x^{\nu+1}, x^\nu) + c_\nu(b_i - \langle a^i, x^{\nu+1}\rangle). \tag{6.43}$$

The conclusion that $d_\nu \geq 0$, for all $\nu$, follows by looking at this expression. From Lemma 2.1.1, $D_f(x^{\nu+1}, x^\nu) \geq 0$. For the second summand we know from Lemma 6.3.1 that $x^{\nu+1} \in Q_i$, $i \equiv i(\nu)$, so that $b_i - \langle a^i, x^{\nu+1}\rangle \geq 0$. By (6.37) we know that $z_i^\nu \geq 0$; thus $c_\nu$ can be negative only if it is equal to $\theta_\nu$. But by Lemma 2.2.2 we see that in this case $x^\nu \notin Q_i$, hence (see Lemma 6.3.1), $x^{\nu+1} \in H_i$, implying that $b_i - \langle a^i, x^{\nu+1}\rangle = 0$, thus proving that, in any case, the second summand in (6.43) is also nonnegative.

We now prove existence of $\lim_{\nu\to\infty} L(x^\nu, z^\nu)$. In view of (6.40) we have only to show that $\{L(x^\nu, z^\nu)\}$ is bounded from above, for all $\nu \geq 0$. To see this, take any $z \in Q \cap \overline{S}$ and verify that

$$\langle A^{\mathrm{T}}z^\nu, z - x^\nu\rangle = \langle z^\nu, Az\rangle - \langle z^\nu, Ax^\nu\rangle \leq \langle z^\nu, b - Ax^\nu\rangle, \tag{6.44}$$

because $Az \leq b$, and by (6.37) $z^\nu \geq 0$, for all $\nu \geq 0$. From (2.1), (6.37), and (6.44), we get

$$D_f(z, x^\nu) \leq f(z) - f(x^\nu) + \langle z^\nu, b - Ax^\nu\rangle = f(z) - L(x^\nu, z^\nu), \tag{6.45}$$

from which we infer that, for any $z \in Q \cap \overline{S}$,

$$L(z^\nu, z^\nu) \leq f(z) - D_f(z, x^\nu) \leq f(z), \tag{6.46}$$

because of Lemma 2.1.1. Therefore $\{L(x^\nu, z^\nu)\}$ is bounded from above.

**Step 3.** The convergence to zero of the differences of consecutive iterates—often a cornerstone in proving convergence of iterative methods—can now be deduced without difficulty. Here the difference is measured in terms of the generalized distance $D_f$ associated with $f(x)$, i.e.,

$$\lim_{\nu\to\infty} D_f(x^{\nu+1}, x^\nu) = 0. \tag{6.47}$$

From (6.43) we take $0 \leq D_f(x^{\nu+1}, x^\nu) \leq d_\nu$ and use the existence of $\lim_{\nu\to\infty} L(x^\nu, z^\nu)$ to see that $\lim_{\nu\to\infty} d_\nu = 0$, thus validating (6.47).

**Step 4.** Now we show that the sequence $\{x^\nu\}$ is bounded. From (6.40) and (6.45) we see that $D_f(z, x^\nu) \leq f(z) - L(x^0, z^0) \doteq \alpha$. Therefore, from the boundedness of the partial level set $L_2^f(z, \alpha)$ (see Definition 2.1.1), we obtain the desired result.

**Step 5.** We prove the feasibility of any accumulation point $x^*$ of $\{x^\nu\}$. Let $\lim_{j\to\infty} x^{\nu_j} = x^*$. By repeated application of Definition 2.1.1$(vi)$, together with (6.47) and the boundedness of $\{x^\nu\}$, it follows that

$$\lim_{j\to\infty} x^{\nu_j+t} = x^*, \quad \text{for every } t \in \{0, 1, 2, \ldots, C\}, \tag{6.48}$$

where C is the constant of almost cyclicality.

Consider the semi-infinite array that has $C+1$ rows, with the elements of the sequence $\{x^{\nu_j+t}\}_{j=1}^\infty$ in the $t$th row. We show that, for each $i \in I$, some row of the array contains an infinite number of elements of $Q_i$. Indeed, for each $i \in I$ at least one element in each column of the array must belong to $Q_i$, otherwise almost cyclicality is violated. It follows that there must be some row in the array, infinitely many elements of which belong to $Q_i$, because if all rows had only a finite number of elements belonging to $Q_i$, there would be only a finite number of elements belonging to $Q_i$ in the entire array. Thus one can extract subsequences from the rows of the array, which belong to each of the $Q_i$'s, all of which converge to the same $x^*$. Therefore, $x^* \in Q \cap \overline{S}$ where $Q = \cap_{i=1}^m Q_i$, because $Q_i$ is closed for every $i \in I$ and because of the strong zone consistency assumption.

**Step 6.** Here we show that if $x^*$ is an accumulation point of $\{x^\nu\}$, then it is a solution of (6.12)–(6.14). We then take advantage of the strict convexity of $f(x)$ to conclude the proof. Denote by $I_1$ and $I_2$ the sets of indices of inactive and active constraints at $x^*$, respectively, i.e.,

$$\begin{aligned} I_1 &\doteq \{i \in I \mid \langle a^i, x^* \rangle < b_i\}, \\ I_2 &\doteq \{i \in I \mid \langle a^i, x^* \rangle = b_i\}. \end{aligned} \tag{6.49}$$

In our proof of Step 6 we need the following proposition.

**Proposition 6.3.1** *If $x^* \in Q \cap \overline{S}$ is an accumulation point of $\{x^\nu\}$ produced by Algorithm 6.3.1, then there exists a subsequence $\{x^{\nu_j}\}_{j=0}^{\infty}$ that converges to $x^*$ and has the property that $z_i^{\nu_j} = 0$, for all $i \in I_1$ and all $j \geq 0$.*

We will postpone the proof of this proposition until we have completed the proof of Theorem 6.3.1, and we now proceed with the proof of Step 6. Using the subsequence of Proposition 6.3.1, we may write

$$\langle z^{\nu_j}, Ax^{\nu_j} - b\rangle = \langle A^T z^{\nu_j}, x^{\nu_j} - x^*\rangle, \tag{6.50}$$

because $z_i^{\nu_j} = 0$ for $i \in I_1$, while $\langle a^i, x^*\rangle = b_i$ for $i \in I_2$, and $I = I_1 \cup I_2$.

From (6.38) we have that $A^T z^{\nu_j} = -\nabla f(x^{\nu_j})$, so that

$$\begin{aligned}\langle z^{\nu_j}, Ax^{\nu_j} - b\rangle &= -\langle \nabla f(x^{\nu_j}), x^{\nu_j} - x^*\rangle \\ &= D_f(x^*, x^{\nu_j}) - f(x^*) + f(x^{\nu_j}),\end{aligned} \tag{6.51}$$

which tends to zero, as $j \to \infty$ because of Definition 2.1.1. Therefore using again the continuity of $f$ on $\overline{S}$, we have

$$\lim_{j\to\infty} L(x^{\nu_j}, z^{\nu_j}) = \lim_{j\to\infty} (f(x^{\nu_j}) + \langle z^{\nu_j}, Ax^{\nu_j} - b\rangle) = f(x^*). \tag{6.52}$$

From (6.46) it now follows that $f(x^*) \leq f(z)$, for all $z \in Q \cap \overline{S}$, which proves that $x^*$ is a solution for (6.12)–(6.14). Because of the strict convexity of $f$, $x^*$ must be unique, which proves that $\lim_{\nu\to\infty} x^\nu = x^*$, with $x^*$ being the solution of (6.12)–(6.14). ■

We now prepare for the proof of Proposition 6.3.1, which is constructed by an elaboration on the argument in Step 5 of Theorem 6.3.1. Denote by $\alpha$ the cardinality of $I_1$ in (6.49). We arbitrarily order $I_1$ as $\{i_1, i_2, \ldots, i_\alpha\}$, and for $\beta \leq \alpha$ we write $I_{1,\beta} \doteq \{i_1, i_2, \ldots, i_\beta\}$. In this notation $I_1 = I_{1,\alpha}$. We proceed as follows. First, we show the existence of a subsequence $\{x^{\nu_i}\} \to x^*$, with $z_l^{\nu_i} = 0$, for $l \in I_{1,1}$. From this subsequence a new subsequence is extracted—$\{x^{\nu_s}\}$ say—that also converges to $x^*$ and for which $z_l^{\nu_s} = 0$, for all $l \in I_{1,2}$. The argument is extended inductively to all of $I_1$. We use the nested sets notation to construct subsequences. Let $N_0$ denote the ordered set of nonnegative integers, and let $K$ and $L$ be ordered subsets such that $L \subseteq K \subseteq N_0$. Then, $\{x^\nu\}_{\nu\in K}$ is a subsequence of $\{x^\nu\}_{\nu=0}^{\infty}$, $\{x^\nu\}_{\nu\in L}$ is a subsubsequence, and so on. The construction of the desired subsequence will be affected by the following lemmas.

**Lemma 6.3.2** *Let $\{x^\nu\}_{\nu\in K} \to x^*$, $K \subseteq N_0$. Let $l \in I_1$, and let $i(\nu)$ be the control index. Then, the event that $i(\nu) = l$ and $c_\nu = \theta_\nu$ in the notation of (6.32), can occur only finitely many times.*

**Proof** Suppose the contrary. Let $K' \subseteq K$ be an infinite subset of indices on which the event does occur. Then, from (6.32), $x^{\nu+1} \in H_l$, for all $\nu \in K'$, where $H_l \doteq \{x \mid \langle a^l, x\rangle = b_l\}$. But, from (6.48), $\{x^{\nu+1}\}_{\nu\in K'}$ converges to $x^*$, and the closedness of $H_l$ would force $x^* \in H_l$, contradicting $l \in I_1$. ■

**Lemma 6.3.3** *Let $\{x^\nu\}_{\nu\in K} \to x^*$, $K \subseteq N_0$, and let $l \in I_1$. Then, for some $t$ there exists a convergent subsequence of the form $\{x^{\nu+t}\}_{\nu\in L}$, $L \subseteq K$, and $0 \le t \le C+1$, where $C$ is the constant of almost cyclicality of the control sequence, with the property that $z_l^{\nu+t} = 0$, for all $\nu \in L$.*

**Proof** Construct the array described in Step 5 of the proof of Theorem 6.3.1. By the argument of Step 5 there is at least one index $s$ with the property $\{x^{\nu+s}\}_{\nu\in K'} \to x^*$, with $i(\nu+s) = l$ and $K' \subseteq K$, $K'$ an infinite subset. From Lemma 6.3.2 and (6.32) $c_{\nu+s} < \theta_{\nu+s}$, for $\nu \in L$, where $L$ is the ordered subset of all elements of $K$ from which the indices where equality occurs have been deleted; thus,

$$L \doteq K' \backslash \{\text{ finite number of events where equality occurs }\}.$$

Then, from (6.32) $c_{\nu+s} = z_l^{\nu+s}$ and $z_l^{\nu+s+1} = 0$. But (6.48) shows that $\{x^{\nu+s+1}\}_{\nu\in L} \to x^*$, so the proof is complete. ■

**Lemma 6.3.4** *Let $\{x^\nu\}_{\nu\in K} \to x^*$, $K \subseteq N_0$. Suppose that $z_l^\nu = 0$, for all $\nu \in K$ and all $l \in I_{1,\beta}$. Then, there exists a subsequence $\{x^\nu\}_{\nu\in L}$, $L \subseteq N_0$, with $z_l^\nu = 0$, for all $\nu \in L$ and all $l \in I_{1,\beta+1}$.*

**Proof** The index set $L$, whose existence is asserted, will be the direct sum $K'' \oplus t \doteq \{\mu = \nu + t \mid \nu \in K''\}$, where $K''$ differs from $K$ by only a finite number of elements, and where $t$ is an integer, $0 \le t \le C+1$. From $K$, we first remove any index $\nu$ which, for some integer $s$, $0 \le s \le C$, has the properties $i(\nu+s) \in I_{1,\beta}$ and $c_{\nu+s} = \theta_{\nu+s}$, and call the residual set $K'$. From Lemma 6.3.2, $K'$ is infinite.

Next, we create the semi-infinite array whose $s$th row is $\{x^{\nu+s}\}_{\nu\in K'}$, for $s = 0, 1, \ldots, C$. Since we have removed all indices $\nu$ for which $z_l^\nu, l \in I_{1,\beta}$, might change, it follows that, for all $l \in I_{1,\beta}$ and all $0 \le s \le C+1$, we have $z_l^{\nu+s} = 0$, for all $\nu \in K'$.

Now, use Lemma 6.3.3, taking $l = i(\beta+1)$, to generate a subsequence $\{x^\nu\}_{\nu\in K''}$, $K'' \subset K'$, with $z_i^\nu = 0$ for all $i \in I_{1,\beta+1}$ and all $\nu \in K''$. ■

**Proof** (**of Proposition 6.3.1**) Use Lemma 6.3.3 to start the construction, and then use Lemma 6.3.4 inductively. ■

This brings us to the end of the proof of Theorem 6.3.1. Strong zone consistency is required in this theorem to ensure that the iterates stay within the zone $S$, regardless of the values of $c_\nu$ in (6.32). In the next section, we use this result to establish a variant of the method, which is specifically designed to handle interval constraints in an efficient manner.

## 6.4 Algorithm for Interval-Constrained Problems

Consider the convex programming problem with interval constraints

$$\text{Minimize} \quad f(y) \tag{6.53}$$

$$\text{s.t.} \quad \gamma_j \le \langle \phi^j, y\rangle \le \delta_j, \qquad j \in J \doteq \{1, 2, \ldots, p\}, \tag{6.54}$$

$$y \in \overline{S}, \tag{6.55}$$

where, for all $j \in J$, $\phi^j \in \mathbb{R}^n$ are given vectors, $\gamma_j$ and $\delta_j$ are given real numbers that define the intervals, and $f \in \mathcal{B}(S)$ is a Bregman function.

This problem could easily be transformed to the form of (6.12)–(6.14) by multiplying one side of the interval inequalities by −1. Then Algorithm 6.3.1 could be applied directly. But doing so would result in a problem with $2p$ one-sided inequality constraints, forcing us to deal with as many dual variables. This doubling of the number of dual variables should be avoided when working with a large number of inequalities.

The method for interval convex programming attaches only one dual variable to each interval constraint, and requires only half as many iterative steps as would be needed if Algorithm 6.3.1 were applied directly to the transformed problem. This interval convex programming method is obtained by choosing a particular strategy for the application of Algorithm 6.3.1. The almost cyclic control, introduced into Algorithm 6.3.1 in the previous section, is an indispensable tool for validating the process. Before stating the algorithm, we introduce some definitions and notations.

$$H(\phi^j, \beta_j) \doteq \{y \in \mathbb{R}^n \mid \langle \phi^j, y\rangle = \beta_j\} \tag{6.56}$$

denotes a hyperplane whose representation is assumed *fixed* and determined by its normal vector $\phi^j \in \mathbb{R}^n$ and the scalar $\beta_j$. The constraints of (6.54) will be described by the half-spaces

$$Q_{j^+} \doteq \{y \mid \langle \phi^j, y\rangle \le \delta_j\}, \tag{6.57}$$

$$Q_{j^-} \doteq \{y \mid \gamma_j \le \langle \phi^j, y\rangle\}, \tag{6.58}$$

for all $j \in J$, so that $Q_j = Q_{j^+} \cap Q_{j^-}$ is the hyperslab defined by the $j$th interval constraint. Again denote $Q = \cap_{j=1}^p Q_j$, and assume that $\phi^j \neq 0$, for all $j \in J$. Define

$$U \doteq \{y \in S \mid \exists\, u \in \mathbb{R}^p, \text{ such that } \nabla f(y) = -\Phi^T u\}, \tag{6.59}$$

where $\Phi$ stands for the $p \times n$ matrix whose $j$th row is $(\phi^j)^T$, and assume $U \neq \emptyset$ throughout.

**Algorithm 6.4.1 Row-Action Algorithm for Interval Convex Programming**

**Step 0:** (Initialization.) $y^0 \in U$ is arbitrary, and $u^0$ is such that

$$\nabla f(y^0) = -\Phi^{\mathrm{T}} u^0. \tag{6.60}$$

**Step 1:** (Iterative step.) Given $y^\nu$ and $u^\nu$, calculate $y^{\nu+1}$ and $u^{\nu+1}$ from

$$\begin{aligned} \nabla f(y^{\nu+1}) &= \nabla f(y^\nu) + d_\nu \phi^{j(\nu)}, \\ u^{\nu+1} &= u^\nu - d_\nu e^{j(\nu)}, \\ d_\nu &\doteq \text{mid } (u^\nu_{j(\nu)}, \Delta_\nu, \Gamma_\nu), \end{aligned} \tag{6.61}$$

where $\Delta_\nu$ and $\Gamma_\nu$ are the parameters associated with the generalized projection of $y^\nu$ onto $H(\phi^{j(\nu)}, \delta_{j(\nu)})$ and $H(\phi^{j(\nu)}, \gamma_{j(\nu)})$, respectively, and $\text{mid}(a, b, c)$ denotes the median of the three real numbers $a, b, c$.

**Control:** The sequence $\{j(\nu)\}_{\nu=0}^{\infty}$ is almost cyclic on $J$.

The main idea behind this algorithm can be appreciated from the following remarks. Define

$$a^{j^+} \doteq \phi^j, \qquad b_{j^+} \doteq \delta_j, \tag{6.62}$$

$$a^{j^-} \doteq -\phi^j, \qquad b_{j^-} \doteq -\gamma_j, \tag{6.63}$$

for all $j \in J$. Then, the problem

$$\text{Minimize} \quad f(y) \tag{6.64}$$

$$\text{s.t.} \quad \langle a^{j^+}, y\rangle \le b_{j^+}, \quad \text{for all } j \in J, \tag{6.65}$$

$$\langle a^{j^-}, y\rangle \le b_{j^-}, \quad \text{for all } j \in J, \tag{6.66}$$

$$y \in \overline{S}, \tag{6.67}$$

has the same solution as (6.53)–(6.55). Associate dual variables $z_j^+$ and $z_j^-$ with the constraints of (6.65) and (6.66), respectively, and interlace the vectors $z^+ = (z_j^+)$ and $z^- = (z_j^-)$ obtained in this way into a single dual vector $z$, having $2p$ components, defined by

$$z^{\mathrm{T}} \doteq (z_1^+, z_1^-, z_2^+, z_2^-, \ldots, z_j^+, z_j^-, \ldots, z_p^+, z_p^-). \tag{6.68}$$

The single dual variable used by Algorithm 6.4.1 for the $j(\nu)$th interval constraint of (6.53)–(6.55) is the $j(\nu)$th component of the vector

$$u \doteq z^+ - z^-. \tag{6.69}$$

The proof of Algorithm 6.4.1 consists in applying Algorithm 6.3.1 to the problem (6.64)–(6.67) in such a way that pairs of inequalities originating from interval constraints in (6.53)–(6.55) are taken up in an almost cyclic order, but the decision as to whether to take first $H(\phi^{j(\nu)}, \delta_{j(\nu)})$ and then $H(\phi^{j(\nu)}, \gamma_{j(\nu)})$, or vice versa, varies from pair to pair, depending on the sign of the $j(\nu)$th component of the current vector $u^\nu$ of dual variables. We note that such a strategy would not be permissible without the extra feature of almost cyclicality established in Theorem 6.3.1. It is precisely the flexibility of almost cyclic control that allows us to obtain Algorithm 6.4.1 from Algorithm 6.3.1. The next theorem establishes the convergence of Algorithm 6.4.1.

**Theorem 6.4.1** *Let $f \in \mathcal{B}(S)$ be strongly zone consistent with respect to the hyperplanes $H(\phi^j, \delta_j)$ and $H(\phi^j, \gamma_j)$, for all $j \in J$. Assume that $Q \cap \overline{S} \neq \emptyset$ and that $\phi^j \neq 0$, for all $j \in J$. Assume also that $\{j(\nu)\}_{\nu=0}^{\infty}$ is almost cyclic on $J$ and that $U \neq \emptyset$. Then, any sequence $\{y^\nu\}_{\nu=0}^{\infty}$ generated by Algorithm 6.4.1 converges to a solution of (6.53)–(6.55).*

**Proof** The initialization step of Algorithm 6.3.1 for (6.64)–(6.67) allows us to pick an arbitrary $y^0 \in \hat{Z}_0$, where $\hat{Z}_0$ is defined as in (6.16), with $\hat{A}$ instead of $A$, and $\hat{A}$ being the $2p \times n$ matrix whose odd and even rows are formed from the vectors $a^{j^+}$ and $a^{j^-}$, respectively. The initial dual vector $z^0 \in \mathbb{R}_+^{2p}$ is such that $\nabla f(y^0) = -\hat{A}^T z^0$, and we impose the condition that $z_j^{+0} \cdot z_j^{-0} = 0$, for all $j \in J$.

We take the vectors $z^{+0} = (z_j^{+0})_{j=1}^p$ and $z^{-0} = (z_j^{-0})_{j=1}^p$, such that $z_j^{+0} \doteq \max\,(0, u_j^0)$ and $z_j^{-0} \doteq \max\,(0, -u_j^0)$. It is not difficult to see that this initialization coincides with (6.60) and that $u^0 = z^{+0} - z^{-0}$.

Now we consider the iterative step with given $y^\nu$ and $u^\nu$. The next index $j(\nu)$ is selected from $J$, thereby picking a pair of constraints from (6.64)–(6.67) that correspond to one interval constraint in (6.53)–(6.55). For simplicity we abbreviate here $j \equiv j(\nu)$, since this index does not change within a given step. Now, the sign of $u_j^\nu$ has to be checked. Two cases are possible:

If $u_j^\nu \geq 0$, then apply (6.32) first to the constraint $j^+$ and then to the constraint $j^-$.

If $u_j^\nu < 0$, then apply (6.32) first to the constraint $j^-$ and then to the constraint $j^+$.

We work out the details only for the first case, since the second can be similarly handled. Applying (6.32) to the $j^+$th constraint of (6.64)–(6.67), we have

$$\begin{aligned} \nabla f(y^{\nu+1/2}) &= \nabla f(y^\nu) + c' a^{j^+}, \\ z^{\nu+1/2} &= z^\nu - c' e^{j^+}, \\ c' &= \min((z_j^+)^\nu, \Delta_\nu), \end{aligned} \tag{6.70}$$

where $\Delta_\nu$ is the parameter associated with the projection of $y^\nu$ onto the hyperplane $H(\phi^j, \delta_j)$, and $y^{\nu+1/2}$ and $z^{\nu+1/2}$ denote the primal and dual vectors obtained from $y^\nu$ and $z^\nu$, respectively. Then the other side of the $j$th interval constraint is taken up to produce $y^{\nu+1}$ and $z^{\nu+1}$ from $y^{\nu+1/2}$ and $z^{\nu+1/2}$, respectively,

$$\begin{aligned} \nabla f(y^{\nu+1}) &= \nabla f(y^{\nu+1/2}) + c'' a^{j^-}, \\ z^{\nu+1} &= z^{\nu+1/2} - c'' e^{j^-}, \\ c'' &= \min((z_j^-)^{\nu+1/2}, (\Gamma_-)_{\nu+1/2}), \end{aligned} \tag{6.71}$$

where $(\Gamma_-)_{\nu+1/2}$ is the parameter associated with the projection of $y^{\nu+1/2}$ onto $H(-\phi^j, -\gamma_j)$. Combining (6.70) and (6.71) into a single calculation when going from the $\nu$th to the $(\nu+1)$th iterate, leads to

$$\nabla f(y^{\nu+1}) = \nabla f(y^\nu) + (c' - c'')\phi^j, \tag{6.72}$$

while the changes made to dual variables, which were first

$$(z_j^+)^{\nu+1/2} = (z_j^+)^\nu - c', \quad (z_j^-)^{\nu+1/2} = (z_j^-)^\nu, \tag{6.73}$$

and then

$$(z_j^-)^{\nu+1} = (z_j^-)^{\nu+1/2} - c'', \quad (z_j^+)^{\nu+1} = (z_j^+)^{\nu+1/2}, \tag{6.74}$$

add up to

$$(z_j^+)^{\nu+1} - (z_j^-)^{\nu+1} = (z_j^+)^\nu - (z_j^-)^\nu - (c' - c''), \tag{6.75}$$

which, in view of (6.69), gives

$$u^{\nu+1} = u^\nu - (c' - c'')e^j. \tag{6.76}$$

Next, we show that $d \doteq c' - c''$ coincides with the $d_\nu$ dictated by (6.61). Lemma 2.2.3 plays here an important role by providing us with vital information about the parameters associated with the resulting projections. When applying Lemma 2.2.3 extra care must be taken that the correct representations of the hyperplanes are used. Let $H_1$ of Lemma 2.2.3 be either $H(\phi^j, \delta_j)$ or $H'$, depending on whether $c'$ takes the value $\Delta_\nu$ or $(z_j^+)^\nu$, respectively. Let $H_2$ of Lemma 2.2.3 be $H(\phi^j, \gamma_j)$. The parameter associated with the projection of $y^{\nu+1/2}$ onto $H(\phi^j, \gamma_j)$ will be denoted by $\Gamma_{\nu+1/2}$, and it has the value (see (6.71))

$$\Gamma_{\nu+1/2} = -(\Gamma_-)_{\nu+1/2}. \tag{6.77}$$

Also, call $\Gamma_\nu$ the parameter associated with the projection of $y^\nu$ onto $H(\phi^j, \gamma_j)$. With the correspondence $\lambda_1 = c'$, $\lambda_2 = \Gamma_\nu$, $\hat{\lambda}_2 = \Gamma_{\nu+1/2}$, we get, from Lemma 2.2.3, that

$$c' + \Gamma_{\nu+1/2} = \Gamma_\nu, \tag{6.78}$$

which, by (6.77) is the same as

$$(\Gamma_-)_{\nu+1/2} = c' - \Gamma_\nu. \tag{6.79}$$

Figure 6.2, which describes a special case, helps clarify the argument. $H(-\phi^j, -\gamma_j)$ and $H(\phi^j, \gamma_j)$ are two different representations of the same hyperplane. In Algorithm 6.4.1, the first representation is employed, while the application of Lemma 2.2.3 requires that the second representation be considered. From (6.71), (6.73), and (6.79) we get

$$c'' = \min((z_j^-)^\nu, c' - \Gamma_\nu), \tag{6.80}$$

and then

$$\begin{aligned} d &= c' - \min((z_j^-)^\nu, c' - \Gamma_\nu) \\ &= \max(c' - (z_j^-)^\nu, \Gamma_\nu) \\ &= \max(\Gamma_\nu,\ \min(u_j^\nu, \Delta_\nu - (z_j^-)^\nu)). \end{aligned} \tag{6.81}$$

From Lemma 2.2.4 $\Gamma_\nu \leq \Delta_\nu$, and below we show that $(z_j^-)^\nu = 0$. Therefore,

$$\begin{aligned} d &= \max(\Gamma_\nu, \min(u_j^\nu, \Delta_\nu)) \\ &= \mathrm{mid}(u_j^\nu, \Delta_\nu, \Gamma_\nu). \end{aligned} \tag{6.82}$$

To see that $z_j^{-\nu} = 0$ in the first case, we show that

$$(z_j^+)^\nu \geq 0,\ \ (z_j^-)^\nu \geq 0, \text{ and } (z_j^+)^\nu (z_j^-)^\nu = 0, \ \ \text{for all } j \in J. \tag{6.83}$$

The nonnegativity assertion in (6.83) follows from the initialization of the algorithm by an argument similar to that in Step 1 of the proof of Theorem 6.3.1. The rest is proved via induction on $\nu$, by showing that $(z_j^+)^\nu \cdot (z_j^-)^\nu = 0$ implies $(z_j^+)^{\nu+1} \cdot (z_j^-)^{\nu+1} = 0$.

Since we are in the first case, $u_j^\nu \geq 0$, so that

$$(z_j^-)^\nu = 0, \tag{6.84}$$

by the inductive hypothesis. Expressions for $(z_j^+)^{\nu+1}$ and $(z_j^-)^{\nu+1}$ in terms of $(z_j^+)^\nu$ and $(z_j^-)^\nu$ are given in (6.73) and (6.74) so that we obtain the equations $(z_j^-)^{\nu+1} = (z_j^-)^\nu - c'' = -c''$, and $(z_j^+)^{\nu+1} = (z_j^+)^\nu - c'$.

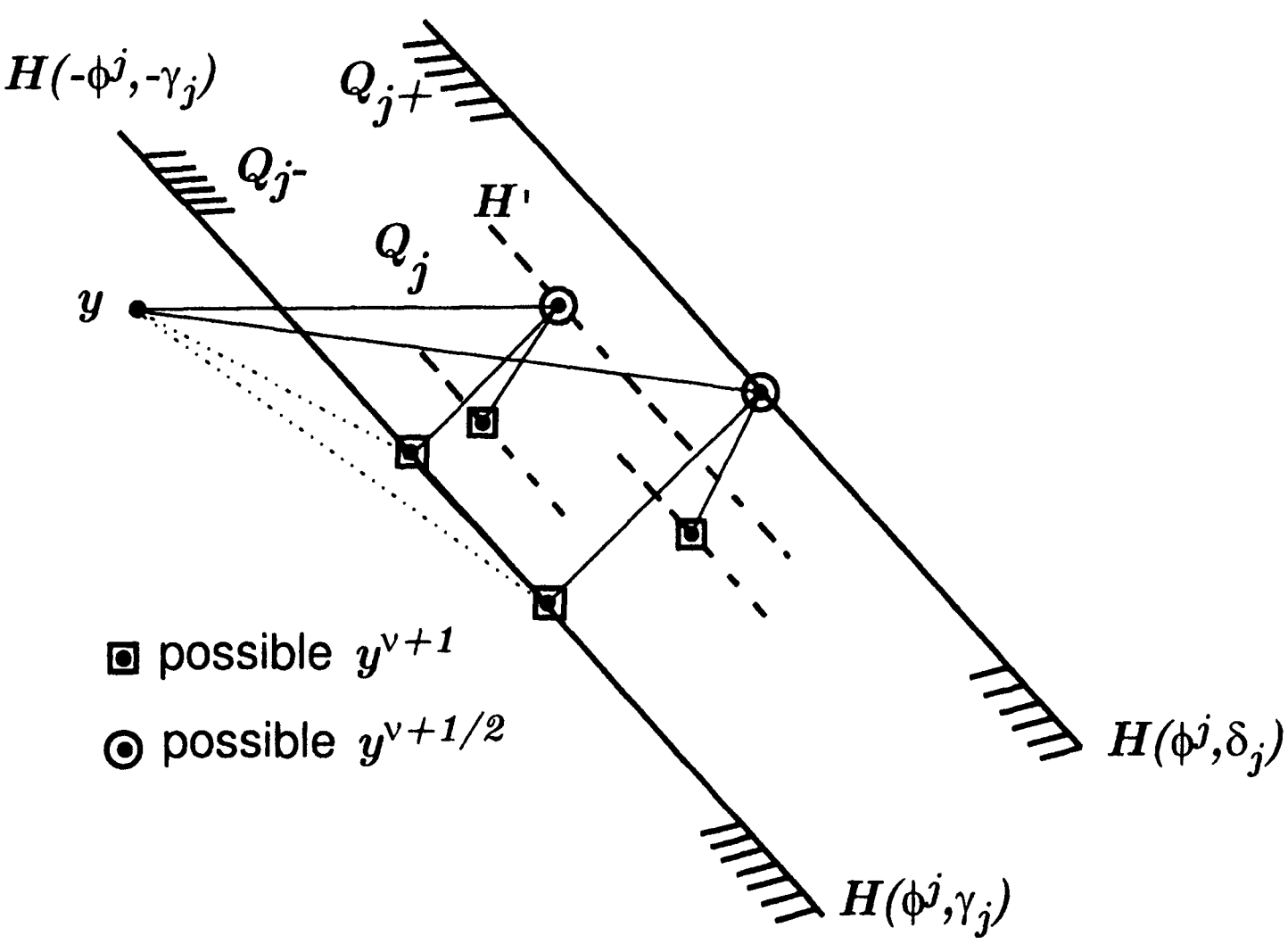

**Figure 6.2** Geometric interpretation of an iterative step in Algorithm 6.4.1.

If $(z_j^+)^\nu \leq \Delta_\nu$, it follows from the definition of $c'$ (see (6.70)), that $(z_j^+)^{\nu+1} = 0$, so we assume that $\Delta_\nu < (z_j^+)^\nu$, implying $c' = \Delta_\nu$. Then, using the information in (6.71), (6.73), (6.79), and (6.84), we get the equation $c'' = \min(0, c' - \Gamma_\nu) = \min(0, \Delta_\nu - \Gamma_\nu)$. But Lemma 2.2.4 shows that $\Delta_\nu \leq \Gamma_\nu$, so $c'' = 0$, completing the inductive step and the proof of this theorem. ■

## 6.5 Row-Action Algorithms for Norm Minimization

Let us now specialize the iterative methods presented in Sections 6.3 and 6.4 for the Bregman function $f(x) = \frac{1}{2} \parallel x \parallel^2$ (see Example 2.1.1). By doing so we obtain, as special cases of the general iterative algorithms (Algorithm 6.3.1 and Algorithm 6.4.1), some well-known special purpose algorithms for norm-minimization. In particular, for linear equality constraints we obtain the algorithm of Kaczmarz; for linear inequality constraints we obtain the algorithm of Hildreth; and for linear interval constraints we obtain the ART4 algorithm of Herman and Lent.

### 6.5.1 The algorithm of Kaczmarz

We have already encountered Kaczmarz's algorithm (Algorithm 5.4.3) in Chapter 5 as an algorithm for solving a linear feasibility problem, i.e., a system of linear equations. Once it is initialized anywhere in the range of the matrix $A^T$, the same algorithm becomes a norm-minimization method.

Apart from the difference in the initialization, the geometric interpretation of the algorithm (see Figure 5.5) remains the same as for Algorithm 5.4.3.

**Algorithm 6.5.1 Kaczmarz's Algorithm for Norm Minimization**

**Step 0:** (Initialization.) $x^0 \in \mathcal{R}(A^T)$, the range of $A^T$.

**Step 1:** (Iterative step.)

$$x^{\nu+1} = x^\nu + \lambda_\nu \frac{b_{i(\nu)} - \langle a^{i(\nu)}, x^\nu \rangle}{\| a^{i(\nu)} \|^2} a^{i(\nu)}, \tag{6.85}$$

**Relaxation parameters:** For all $\nu \geq 0$, $\epsilon_1 \leq \lambda_\nu \leq 2 - \epsilon_2$, for some arbitrarily small but fixed $\epsilon_1, \epsilon_2 > 0$.

**Control:** The sequence $\{i(\nu)\}$ is almost cyclic on $I \doteq \{1, 2, \dots, m\}$.

In the unrelaxed case, i.e., when $\lambda_\nu = 1$ for all $\nu \geq 0$, the convergence of this algorithm can be obtained from the study of Bregman's method (see Section 6.3). With $f(x) = \frac{1}{2} \| x \|^2$, $S = \overline{S} = \mathbb{R}^n$, and only equality constraints, the problem (6.12)–(6.14) becomes,

$$\text{Minimize} \quad \frac{1}{2} \| x \|^2 \tag{6.86}$$

$$\text{s.t.} \quad \langle a^i, x \rangle = b_i \,, \qquad i \in I. \tag{6.87}$$

Kaczmarz's algorithm for norm minimization is then obtained from Algorithm 6.3.1. Lemma 6.3.1 shows that since only equality constraints are present, $x^{\nu+1} \in H_{i(\nu)}$ always holds and $c_\nu = \theta_\nu$ for all $\nu \geq 0$. Thus the dual iterates $z^\nu$ need not be updated at all. This disappearance of the dual variables from the actual iterative step is common to all equality constrained algorithms derived from Algorithm 6.3.1. For the underrelaxed case, i.e., $\epsilon_1 \leq \lambda_\nu \leq 1$ for all $\nu \geq 0$, the convergence can be obtained from the relaxed Bregman method (see Section 6.8). The general case, with the full range of relaxation parameters, cannot be obtained from the algorithmic schemes discussed here, but must be treated separately (see Section 6.10).

### 6.5.2 The algorithm of Hildreth

The algorithm of Hildreth is a special-purpose row-action method obtained from Algorithm 6.3.1 when the Bregman function $f(x) = \frac{1}{2} \| x \|^2$ with zone $S = \overline{S} = \mathbb{R}^n$ is used. Therefore, it is specifically designed to solve the linear inequalities constrained problem

$$\text{Minimize} \quad \frac{1}{2} \| x \|^2 \tag{6.88}$$

$$\text{s.t.} \quad \langle a^i, x \rangle \leq b_i, \qquad i \in I. \tag{6.89}$$

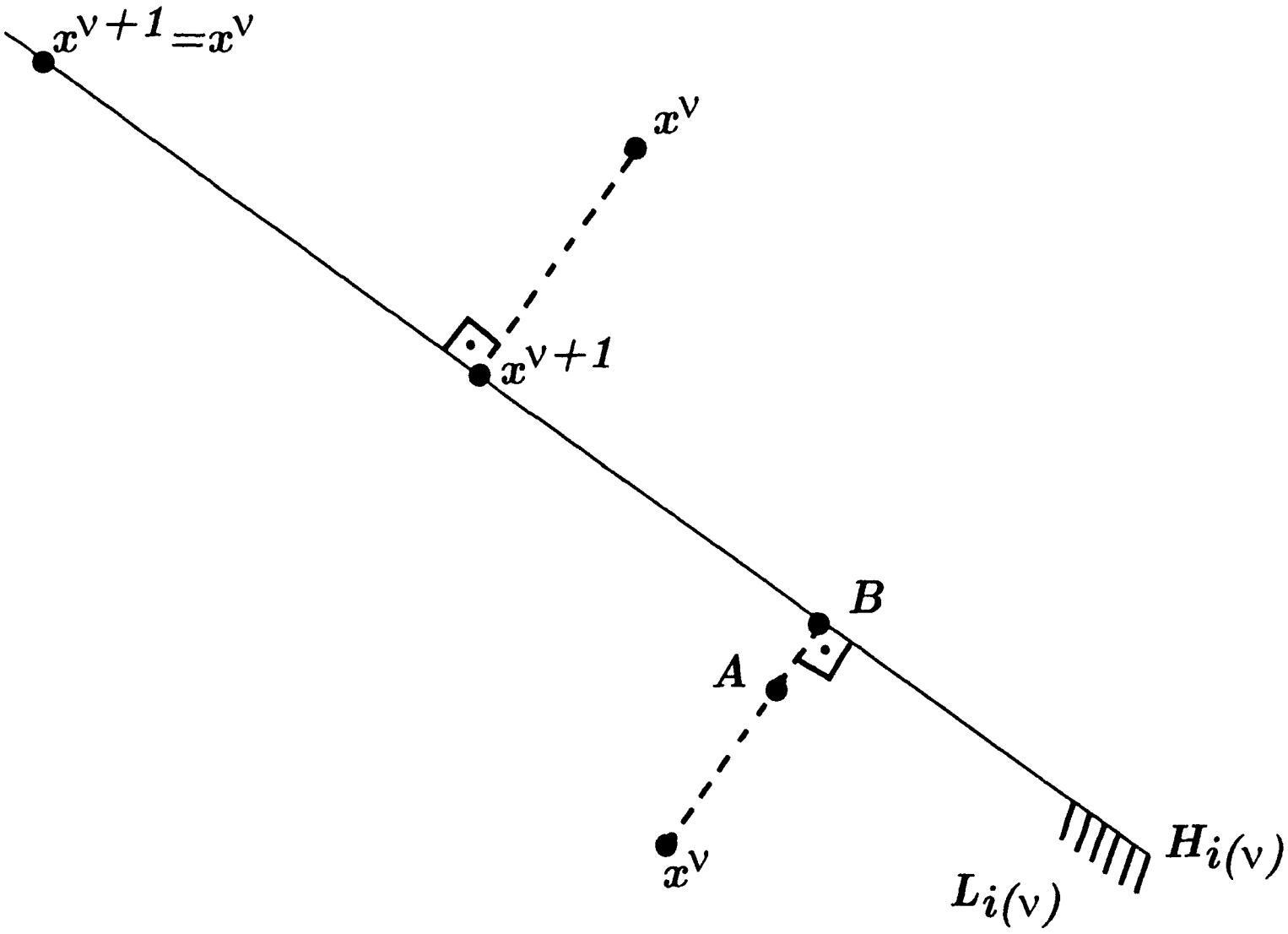

**Figure 6.3** Hildreth's algorithm (Algorithm 6.5.2); all possible cases of an iterative step.

It is a primal-dual algorithm that iteratively updates primal $\{x^\nu\}$ and dual $\{z^\nu\}$ iterates in such a manner that only one component of the dual vector is actually changed in a single iterative step. The size of the change is such that it guarantees dual ascent.

### Algorithm 6.5.2 Hildreth's Algorithm

**Step 0:** (Initialization.) $z^0 \in \mathbb{R}^m_+$ is arbitrary and $x^0 = -A^T z^0$.

**Step 1:** (Iterative step.)

$$x^{\nu+1} = x^\nu + c_\nu a^{i(\nu)}, \tag{6.90}$$

$$z^{\nu+1} = z^\nu - c_\nu e^{i(\nu)}, \tag{6.91}$$

with

$$c_\nu \doteq \min\left( z^\nu_{i(\nu)},\ \lambda_\nu \frac{b_{i(\nu)} - \langle a^{i(\nu)}, x^\nu \rangle}{\| a^{i(\nu)} \|^2} \right). \tag{6.92}$$

**Relaxation parameters:** For all $\nu \geq 0$, $\epsilon_1 \leq \lambda_\nu \leq 2 - \epsilon_2$, for some arbitrarily small but fixed $\epsilon_1, \epsilon_2 > 0$.

**Control:** The sequence $\{i(\nu)\}$ is almost cyclic on $I$.

Again, for unity relaxation where $\lambda_\nu = 1$ for all $\nu \geq 0$, the convergence follows directly from Algorithm 6.3.1. The underrelaxed case can be obtained from the extension of Section 6.8. The convergence for general relaxation parameters in the interval $\epsilon_1 \leq \lambda_\nu \leq 2 - \epsilon_2$ has been treated separately, and thoroughly analyzed in the literature (see Section 6.10).

A geometric interpretation of the algorithm is as follows. Assume that $x^\nu$, $z^\nu$, and the closed half-space $L_{i(\nu)} \doteq \{x \in \mathbb{R}^n \mid \langle a^{i(\nu)}, x\rangle \leq b_{i(\nu)}\}$, determined by the $i(\nu)$th inequality of the problem (6.88)–(6.89) are given. If $x^\nu \notin L_{i(\nu)}$, then $x^{\nu+1}$ is the (possibly relaxed) orthogonal projection of $x^\nu$ onto $L_{i(\nu)}$. If $x^\nu$ belongs to the bounding hyperplane $H_{i(\nu)}$ then $x^{\nu+1} = x^\nu$. Finally, if $x^\nu \in \operatorname{int} L_{i(\nu)}$, i.e., if $\langle a^{i(\nu)}, x^\nu\rangle < b_{i(\nu)}$, then a move perpendicular to the bounding hyperplane $H_{i(\nu)}$ is made. In this case either $c_\nu = z^\nu_{i(\nu)}$ or $x^{\nu+1}$ is the orthogonal projection of $x^\nu$ onto $H_{i(\nu)}$. All these possibilities are depicted in Figure 6.3.

A fully simultaneous version of Hildreth's row-action algorithm has also been proposed and studied (see Section 6.10). In that algorithm a convex combination of individual Hildreth steps with respect to all half-spaces is taken as the next iterate $x^{\nu+1}$, exactly as in the framework of the simultaneous algorithms classification in Section 1.3.

### 6.5.3 ART4—An algorithm for norm minimization over linear intervals

The problem

$$\text{Minimize} \quad \frac{1}{2} \parallel x \parallel^2 \tag{6.93}$$

$$\text{s.t.} \quad \gamma_i \leq \langle a^i, x\rangle \leq \delta_i, \qquad i \in I, \tag{6.94}$$

is of the form (6.53)–(6.55) and Algorithm 6.4.1 applies. The desire to solve a problem with such interval constraints comes from the optimization approach described in Section 6.1.2, where an inconsistent system of equality constraints $\langle a^i, x\rangle = b_i$, $i \in I$ is replaced, by defining $\gamma_i \doteq b_i - \epsilon_i$ and $\delta_i \doteq b_i + \epsilon_i$, for all $i \in I$, by intervals. A practical difficulty is to choose the *tolerances* $\epsilon_i$ as small as possible while ensuring that the system of interval constraints is feasible. A problem of this sort arose in the field of image reconstruction from projections, and thus motivated the development of the ART4 (Algebraic Reconstruction Technique 4) algorithm. It is an extension of Hildreth's algorithm (Algorithm 6.5.2) designed to efficiently handle the interval constraint. ART4 is retrievable from Algorithm 6.4.1 by taking $f(x) = \frac{1}{2} \|x\|^2$ .

**Algorithm 6.5.3 ART4**

**Step 0:** (Initialization.) $x^0 \in \mathbb{R}^n$ and $z^0 \in \mathbb{R}^m$ are such that $x^0 = -A^T z^0$.
**Step 1:** (Iterative step.)

$$x^{\nu+1} = x^\nu + c_\nu a^{i(\nu)}, \tag{6.95}$$

$$z^{\nu+1} = z^\nu - c_\nu e^{i(\nu)}, \tag{6.96}$$

with

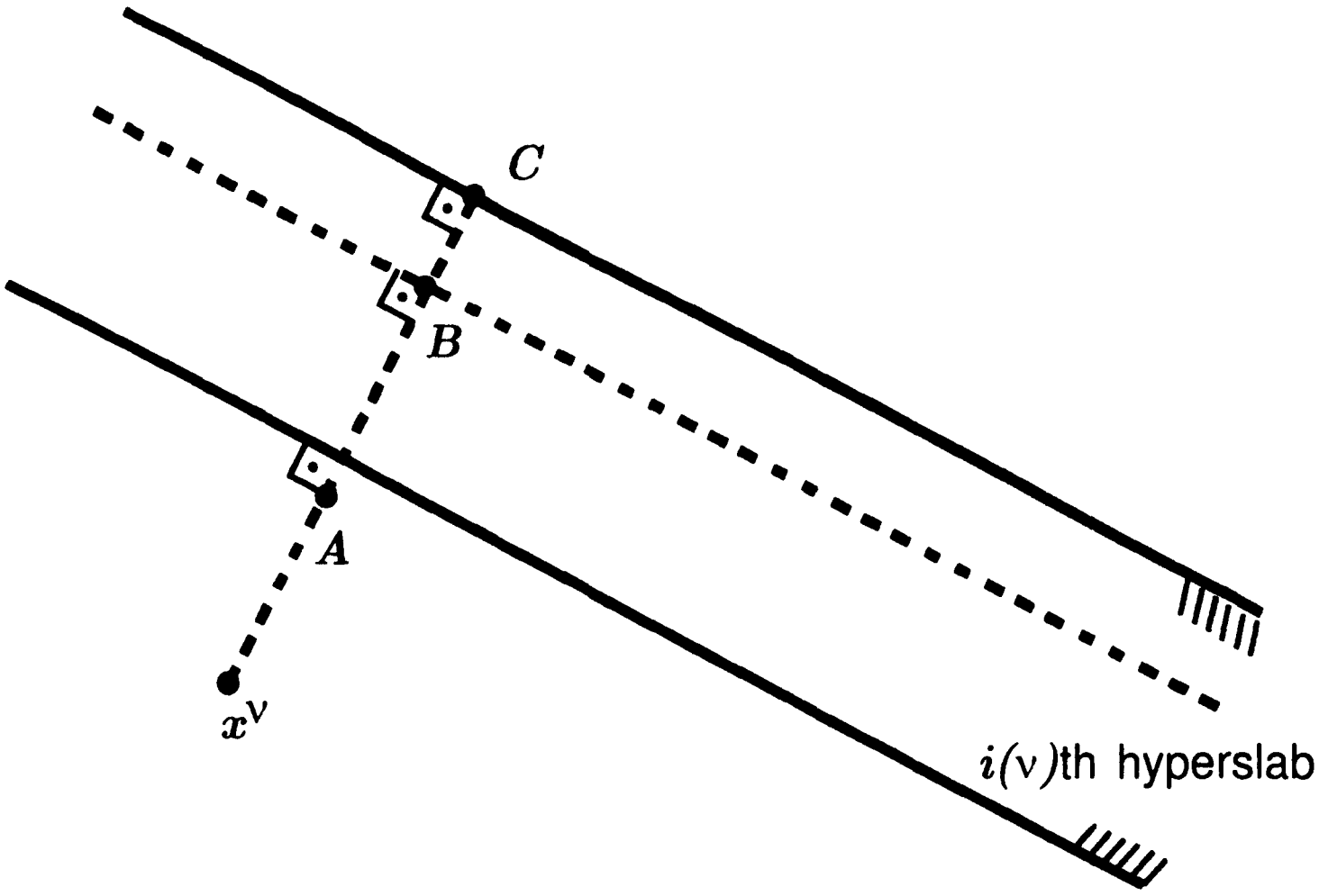

**Figure 6.4** ART4 is a special development of Hildreth's algorithm. In a typical iterative step the next iterate $x^{\nu+1}$ will be one of the points $A, B$ or $C$.

$$c_\nu \doteq \text{mid}\left( z^\nu_{i(\nu)}, \frac{\delta_{i(\nu)} - \langle a^{i(\nu)}, x^\nu \rangle}{\| a^{i(\nu)} \|^2}, \frac{\gamma_{i(\nu)} - \langle a^{i(\nu)}, x^\nu \rangle}{\| a^{i(\nu)} \|^2} \right), \tag{6.97}$$

where mid$(a, b, c)$ stands, as mentioned earlier, for the median of the three real numbers $a, b$ and $c$.

**Relaxation parameters:** Unity, i.e., $\lambda_\nu = 1$, for all $\nu \geq 0$.

**Control:** The sequence $\{i(\nu)\}$ is almost cyclic on $I$.

In a typical iterative step the $i(\nu)$th hyperslab is employed. Given $x^\nu$, the next iterate $x^{\nu+1}$ will be inside the slab or on one of its bounding hyperplanes, according to the value of $c_\nu$. This value of $c_\nu$ is determined by the *mid* operation and depends on the distances of $x^\nu$ from the two bounding hyperplanes and on the value of the $i(\nu)$th component of the $\nu$th dual vector $z^\nu$; see Figure 6.4 where these possibilities are marked by the points $A, B$ or $C$.

If additional box constraints are present in the problem then they can be handled in a simple manner. Consider the following problem.

$$\text{Minimize} \quad \frac{1}{2} \| x \|^2 \tag{6.98}$$

$$\text{s.t.} \quad \gamma_i \leq \langle a^i, x \rangle \leq \delta_i, \qquad i \in I, \tag{6.99}$$

$$l_j \leq x_j \leq u_j, \qquad 1 \leq j \leq n, \tag{6.100}$$

where $l = (l_j)$ and $u = (u_j)$ are some given vectors that determine the box constraints. Instead of replacing $x_j$ by $\langle e^j, x\rangle$ and treating the box constraints as interval constraints, the following algorithm, called ART2, modifies the iterates by resetting variables to the nearest bound whenever the bound is violated.

**Algorithm 6.5.4 ART2**

**Step 0:** (Initialization.) $\hat{x}^0 \in \mathbb{R}^n$ and $z^0 \in \mathbb{R}^m$ are such that $\hat{x}^0 = -A^{\mathrm{T}} z^0$.

**Step 1:** (Iterative step.)

$$\hat{x}^{\nu+1} = \hat{x}^{\nu} + c_\nu a^{i(\nu)}, \tag{6.101}$$

$$z^{\nu+1} = z^{\nu} - c_\nu e^{i(\nu)}, \tag{6.102}$$

with

$$c_\nu \doteq \mathrm{mid}\left(z^\nu_{i(\nu)}, \frac{\delta_{i(\nu)} - \langle a^{i(\nu)}, x^\nu\rangle}{\| a^{i(\nu)} \|^2}, \frac{\gamma_{i(\nu)} - \langle a^{i(\nu)}, x^\nu\rangle}{\| a^{i(\nu)} \|^2}\right), \tag{6.103}$$

and $x^\nu$ denotes the vector whose $t$th component is

$$x_t^\nu \doteq \begin{cases} l_t, & \text{if } \hat{x}_t^\nu < l_t, \\ \hat{x}_t^\nu, & \text{if } l_t \le \hat{x}_t^\nu \le u_t, \\ u_t, & \text{if } u_t < \hat{x}_t^\nu. \end{cases} \tag{6.104}$$

**Relaxation parameters:** Unity, i.e., $\lambda_\nu = 1$, for all $\nu \ge 0$.

**Control:** The sequence $\{i(\nu)\}$ is almost cyclic on $I$.

The convergence of any sequence $\{x^\nu\}$ generated by ART2 to a solution of problem (6.98)–(6.100) follows from Theorem 6.4.1.

## 6.6 Row-Action Algorithms for Shannon's Entropy Optimization

It is possible to specialize the methods of Sections 6.3 and 6.4 for the Bregman function $f(x) = \sum_{j=1}^n x_j \log x_j$ (see Example 2.1.2). By simply calculating gradients, the first equation of (6.32) takes the form

$$x_j^{\nu+1} = x_j^\nu \exp\left(c_\nu a_j^{i(\nu)}\right), \quad j = 1, 2, \ldots, n. \tag{6.105}$$

We look at problems of the form

$$\text{Minimize } \sum_{j=1}^n x_j \log x_j \tag{6.106}$$

$$\text{s.t. } x \in Q \cap \overline{S}, \tag{6.107}$$

where $\overline{S} = \mathbb{R}^n_+$ and $Q$ is one of the following constraints sets:

$$Q_1 \doteq \{x \in \mathbb{R}^n \mid Ax = b\}, \tag{6.108}$$
$$Q_2 \doteq \{x \in \mathbb{R}^n \mid Ax \le b\}, \tag{6.109}$$
$$Q_3 \doteq \{x \in \mathbb{R}^n \mid c \le Ax \le b\}. \tag{6.110}$$

Bregman's algorithm for the solution of (6.106)–(6.107) with $Q = Q_1$ is a row-action iterative procedure which performs successive generalized projections onto the individual hyperplanes of (6.108). Applying Algorithm 6.3.1 in this case yields the following algorithm.

**Algorithm 6.6.1 Bregman's Row-Action Method for Linear Equality Constrained Entropy Maximization**

**Step 0:** (Initialization.) $x^0 \in \mathbb{R}^n_{++}$ is such that for an arbitrary $z^0 \in \mathbb{R}^m_+$

$$x_j^0 = \exp((-A^{\mathrm{T}} z^0)_j - 1), \quad j = 1, 2, \ldots, n. \tag{6.111}$$

**Step 1:** (Iterative step.) Given $x^\nu$ choose a control index $i(\nu)$ and solve the system

$$x_j^{\nu+1} = x_j^\nu \ \exp\,(c_\nu a_j^{i(\nu)}), \quad j = 1, 2, \ldots, n, \tag{6.112}$$
$$\langle x^{\nu+1}, a^{i(\nu)} \rangle = b_{i(\nu)}. \tag{6.113}$$

**Control:** The sequence $\{i(\nu)\}$ is almost cyclic on $I \doteq \{1, 2, \ldots, m\}$.

The system (6.112)–(6.113) represents an *entropy projection* of $x^\nu$ onto the hyperplane $H_{i(\nu)} \doteq \{x \in \mathbb{R}^n \mid \langle a^{i(\nu)}, x \rangle = b_{i(\nu)}\}$ resulting in the next iterate $x^{\nu+1}$. It is a system of $n+1$ equations, of which the first $n$ are not linear in the $n+1$ unknown $x_j^{\nu+1}$, $j = 1, 2, \ldots, n$, and $c_\nu$. It must be solved by some iterative method such as the Newton-Raphson method, or the like. This creates a gap between the theoretical algorithm and its practical implementation. An alternative algorithm for the solution of the same problem with $Q = Q_1$, is MART (Multiplicative Algebraic Reconstruction Technique), which employs a closed-form formula for the iterative updates instead of requiring an iterative inner loop for solving a system such as (6.112)–(6.113).

**Algorithm 6.6.2 Multiplicative Algebraic Reconstruction Technique (MART)**

**Step 0:** (Initialization.) $u^0 \in \mathbb{R}^m$ is arbitrary, and $x^0 \in \mathbb{R}^n$ is defined by

$$1 + \log x_j^0 = (-A^{\mathrm{T}} u^0)_j, \quad j = 1, 2, \ldots, n. \tag{6.114}$$

**Step 1:** (Iterative step.)

$$x_j^{\nu+1} = x_j^{\nu} \left( \frac{b_{i(\nu)}}{\langle a^{i(\nu)}, x^{\nu} \rangle} \right)^{\lambda_\nu a_j^{i(\nu)}}, \quad j = 1, 2, \ldots, n. \tag{6.115}$$

**Relaxation parameters:** $\{\lambda_\nu\}_{\nu=0}^{\infty}$ is a sequence of underrelaxation parameters such that $0 < \epsilon \leq \lambda_\nu \leq 1$ for all $\nu \geq 0$, with some arbitrarily small but fixed $\epsilon > 0$.

**Control:** The sequence $\{i(\nu)\}_{\nu=0}^{\infty}$ is almost cyclic on $I$.

In order to prove convergence of MART, the following assumptions need to be made.

**Assumption 6.6.1** *Feasibility:* $Q_1 \cap \mathbb{R}_+^n \neq \emptyset$.

**Assumption 6.6.2** *Signs:* $a_j^i \geq 0$, $b_i > 0$, *and* $a_j^i \neq 0$, *for all* $i \in I$ *and all* $j = 1, 2, \ldots, n$.

**Assumption 6.6.3** *Normalization:* $Ax = b$ *is scaled so that, for all* $j = 1, 2, \ldots, n$, *and all* $i \in I$, $a_j^i \leq 1$.

Assumption 6.6.2 is not too restrictive because it can usually be made to hold in practice. For example, in image reconstruction from projections (see Chapter 10) $a_j^i \geq 0$ holds by the nature of the problem. $b_i > 0$ means that any equation with zero right-hand side should be removed from the constraints at the start, a very reasonable thing to do when reconstructing an image from its projections. The normalization assumption is also easy to satisfy.

The convergence of a sequence $\{x^\nu\}$ produced by MART to the solution of the linearly constrained entropy optimization problem can be derived indirectly by studying the relationship between MART and Algorithm 6.3.1. Such a study leads to the development of hybrid algorithms presented in Section 6.9. However, a direct proof of convergence can also be given (see Section 6.10) to validate the following result. It also follows from the study of the block-MART algorithm in the next section.

**Theorem 6.6.1** *If Assumptions 6.6.1–6.6.3 hold, then any sequence* $\{x^\nu\}$ *generated by Algorithm 6.6.2 converges to the unique minimizer of problem (6.106)–(6.107) with* $Q = Q_1$ *as in (6.108).*

## 6.7 Block-Iterative MART Algorithm

We study now a block-iterative version of MART (Algorithm 6.6.2) for the solution of problem (6.106)–(6.107) with linear equality constraints. In order to address fixed sets of constraints we assign equations of $Ax = b$ into blocks. This is done by partitioning the index set $I \doteq \{1, 2, \ldots, m\}$ of all row indices of the matrix $A$. Let $M$ be the number of blocks and choose integers $\{m_t\}_{t=0}^{M}$ so that

$$0 = m_0 < m_1 < m_2 < \cdots < m_{M-1} < m_M = m. \tag{6.116}$$

For each $t$, $t = 1, 2, \ldots, M$, the subset $I_t \subseteq I$ defined by

$$I_t = \{m_{t-1} + 1, m_{t-1} + 2, \ldots, m_t\} \tag{6.117}$$

is the index set of the $t$th block. To each block $I_t$, $t = 1, 2, \ldots, M$, a fixed system of weights is assigned by defining

$$0 < w_i^t \leq 1, \text{ for all } i \in I_t, \quad \text{such that } \sum_{i \in I_t} w_i^t = 1. \tag{6.118}$$

If a block contains a single equation, i.e., $I_t$ contains a single index $i$ for some $t$, then $w_i^t = 1$. The block-iterative MART algorithm is formulated as follows.

**Algorithm 6.7.1 Block-Iterative MART**

**Step 0:** (Initialization.) $\pi^0 \in \mathbb{R}_+^m$, and $x^0 \in \mathbb{R}_{++}^n$, so that

$$1 + \log x_j^0 = -(A^T \pi^0)_j, \quad j = 1, 2, \ldots, n. \tag{6.119}$$

**Step 1:** (Iterative step.)

$$x_j^{\nu+1} = \frac{1}{e} \exp\left((-A^T \pi^{\nu+1})_j\right), \quad j = 1, 2, \ldots, n, \tag{6.120}$$

$$\pi_i^{\nu+1} = \begin{cases} \pi_i^\nu, & \text{if } i \notin I_{t(\nu)}, \\ \pi_i^\nu - w_i^{t(\nu)} d_i^\nu, & \text{if } i \in I_{t(\nu)}, \end{cases} \tag{6.121}$$

where, for all $i \in I_{t(\nu)}$,

$$d_i^\nu \doteq \log \frac{b_i}{\langle a^i, x^\nu \rangle}. \tag{6.122}$$

**Control:** The sequence $\{t(\nu)\}_{\nu=0}^\infty$ is almost cyclic on $\{1, 2, \ldots, M\}$.

According to (6.121) dual correction is of the form (6.25) with dual correction vectors that change only components of the dual vector whose indices belong to the $t(\nu)$th block, i.e.,

$$v_i^\nu = \begin{cases} 0, & \text{if } i \notin I_{t(\nu)} \\ -w_i^{t(\nu)} d_i^\nu, & \text{if } i \in I_{t(\nu)}. \end{cases} \tag{6.123}$$

The correction terms $d_i^\nu$, given by (6.122), are different from terms that would have resulted from imposing requirements of the form (6.30). We proceed to prove convergence of Algorithm 6.7.1. It is easy to verify that (6.23) holds. The iterative step of Algorithm 6.7.1 can also be represented by a product as

$$x_j^{\nu+1} = x_j^\nu \prod_{i\in I_{t(\nu)}} \exp\left(w_i^{t(\nu)} d_i^\nu a_j^i\right), \quad j = 1,2,\ldots,n. \tag{6.124}$$

This can be verified by substituting (6.121) into (6.120) which yields

$$\begin{aligned} x_j^{\nu+1} &= \frac{1}{e}\exp\left(-\sum_{i\notin I_{t(\nu)}} a_j^i \pi_i^\nu - \sum_{i\in I_{t(\nu)}} \left(\pi_i^\nu - w_i^{t(\nu)} d_i^\nu\right) a_j^i\right) \\ &= x_j^\nu \exp\left(\sum_{i\in I_{t(\nu)}} w_i^{t(\nu)} d_i^\nu a_j^i\right). \end{aligned} \tag{6.125}$$

Observe also that if $x^0 > 0$ (componentwise) then (6.120) guarantees that

$$x^\nu > 0, \quad \text{for all } \nu \geq 0. \tag{6.126}$$

Next, we use (6.20) to define $L_\nu \doteq L(x^\nu, \pi^\nu)$.

**Proposition 6.7.1** *If Assumptions 6.6.1–6.6.3 hold, then for any $\{x^\nu\}$ and $\{\pi^\nu\}$ generated by Algorithm 6.7.1 $\{L_\nu\}$ is a monotonically increasing sequence.*

**Proof** Using (6.20) and (6.23) we see that

$$L_\nu = f(x^\nu) - \langle \nabla f(x^\nu), x^\nu\rangle - \langle \pi^\nu, b\rangle = -\sum_{j=1}^n x_j^\nu - \langle \pi^\nu, b\rangle. \tag{6.127}$$

Therefore, from (6.121) and (6.124) we get,

$$\begin{aligned} S_\nu \doteq L_{\nu+1} - L_\nu &= \sum_{j=1}^n \left(x_j^\nu - x_j^{\nu+1}\right) - \langle \pi^{\nu+1} - \pi^\nu, b\rangle \\ &= \sum_{j=1}^n x_j^\nu \left(1 - \prod_{i\in I_{t(\nu)}} \exp\left(w_i^{t(\nu)} d_i^\nu a_j^i\right)\right) + \sum_{i\in I_{t(\nu)}} b_i w_i^{t(\nu)} d_i^\nu. \end{aligned} \tag{6.128}$$

To handle (6.128) some inequalities are needed. It is well known from elementary calculus that, for $0 \leq \alpha_i \leq 1$ and any $0 \leq \lambda_i$,

$$1 - \lambda_i^{\alpha_i} \geq \alpha_i(1 - \lambda_i). \tag{6.129}$$

Also, if $\{z_i\}_{i=1}^n$ and $\{r_i\}_{i=1}^n$ are sequences of positive real numbers such that $\sum_{i=1}^n r_i = 1$, then the arithmetic-geometric mean inequality says that

$$\prod_{i=1}^{n} z_i^{r_i} \le \sum_{i=1}^{n} r_i z_i. \tag{6.130}$$

Now, using (6.129) with $\lambda_i \doteq \dfrac{b_i}{\langle a^i, x^\nu \rangle}$, which are positive by Assumption 6.6.2 and by (6.126), and with $\alpha_i \doteq a_j^i$ for a fixed $j$, let $i \in I_{t(\nu)}$, multiply each inequality (6.129) by the weights of (6.118), and sum over $i$. This leads to,

$$1 - \sum_{i \in I_{t(\nu)}} w_i^{t(\nu)} \left( \frac{b_i}{\langle a^i, x^\nu \rangle} \right)^{a_j^i} \ge \sum_{i \in I_{t(\nu)}} w_i^{t(\nu)} a_j^i \left( 1 - \frac{b_i}{\langle a^i, x^\nu \rangle} \right). \tag{6.131}$$

Using (6.130), (6.131), and (6.128), we get

$$\begin{aligned} S_\nu &\ge \sum_{j=1}^{n} x_j^\nu \sum_{i \in I_{t(\nu)}} w_i^{t(\nu)} a_j^i \left( 1 - \frac{b_i}{\langle a^i, x^\nu \rangle} \right) + \sum_{i \in I_{t(\nu)}} b_i w_i^{t(\nu)} \log \left( \frac{b_i}{\langle a^i, x^\nu \rangle} \right) \\ &= \sum_{i \in I_{t(\nu)}} w_i^{t(\nu)} b_i \left( y_i^\nu - 1 - \log\, y_i^\nu \right), \end{aligned} \tag{6.132}$$

where

$$y_i^\nu \doteq \frac{\langle a^i, x^\nu \rangle}{b_i}. \tag{6.133}$$

According to Assumption 6.6.2 and by (6.126) $y_i^\nu > 0$, and thus the expression $y_i^\nu - 1 - \log\, y_i^\nu$ is nonnegative, convex in $y_i^\nu$, and gets its minimal value zero at $y_i^\nu = 1$. Consequently,

$$S_\nu \ge 0, \quad \text{for all } \nu \ge 0, \tag{6.134}$$

which proves the proposition. ∎

**Proposition 6.7.2** *Under Assumptions 6.6.1–6.6.3, the sequence $\{L_\nu\}$ converges for any $\{x^\nu\}$ and $\{\pi^\nu\}$ generated by Algorithm 6.7.1.*

**Proof** In view of Proposition 6.7.1 it is only necessary to show that $\{L_\nu\}$ is bounded from above. Let $y \in Q_1 \cap \mathbb{R}_+^n$ be some fixed feasible point, whose existence is guaranteed by Assumption 6.6.1, and consider the expression

$$D_f(y, x^\nu) = f(y) - f(x^\nu) - \langle \nabla f(x^\nu), y - x^\nu \rangle,$$

with $f(x) = -\text{ent}\, x$, as defined in Example 2.1.2. Convexity and differentiability of $f(x)$ over int $\mathbb{R}_+^n$ ensure that $D_f(y, x^\nu) \ge 0$, for all $\nu \ge 0$ (see Lemma 2.1.1). From (6.20), (6.23) and Assumption 6.6.1, we obtain

$$
\begin{aligned}
D_f(y, x^\nu) &= f(y) - f(x^\nu) + \langle A^{\mathrm{T}} \pi^\nu, y - x^\nu \rangle \\
&= f(y) - f(x^\nu) - \langle \pi^\nu, Ax^\nu - b \rangle \\
&= f(y) - L(x^\nu, \pi^\nu),
\end{aligned}
\tag{6.135}
$$

from which

$$L_\nu \leq f(y), \quad \text{for all } \nu \geq 0, \tag{6.136}$$

follows because of the nonnegativity of $D_f(y, x^\nu)$. ■

**Proposition 6.7.3** *Under Assumptions 6.6.1–6.6.3, any* $\{x^\nu\}_{\nu=0}^\infty$ *generated by Algorithm 6.7.1 is bounded.*

**Proof** The set $\Omega \doteq \{x > 0 \mid D_f(y, x) \leq \alpha\}$ is bounded, for any fixed $y \geq 0$ and any fixed real $\alpha$, because $f(x) = -\text{ent}\, x$ is a Bregman function and has bounded partial level sets (see Lemma 2.1.3). From (6.135) and Proposition 6.7.1, $D_f(y, x^\nu) = f(y) - L_\nu \leq f(y) - L_0 \doteq \alpha$, which puts any sequence $\{x^\nu\}_{\nu=0}^\infty$ generated by Algorithm 6.7.1 in the set $\Omega$. ■

**Proposition 6.7.4** *Under Assumptions 6.6.1–6.6.3, if* $x^*$ *is an accumulation point of a sequence* $\{x^\nu\}$ *generated by Algorithm 6.7.1, then we have* $x^* \in Q_1 \cap \mathbb{R}^n_+$.

**Proof** In view of (6.126) it is only necessary to show that $\langle a^i, x^* \rangle = b_i$, for all $i \in I$. Let

$$\lim_{\ell \to \infty} x^{\nu(\ell)} = x^*, \tag{6.137}$$

and assume first that $i$ is not fixed but rather that it is any element of $I_{t(\nu(\ell))}$. Then, by (6.128), (6.134), and Proposition 6.7.2

$$\lim_{\ell \to \infty} S_{\nu(\ell)} = 0, \tag{6.138}$$

which, by (6.132) and (6.134) leads to

$$\lim_{\ell \to \infty} \frac{\langle a^i, x^{\nu(\ell)} \rangle}{b_i} = 1, \tag{6.139}$$

for any $i \in I_{t(\nu(\ell))}$. Multiplying (6.139) by $b_i$ and using (6.137), we get that $\langle a^i, x^* \rangle = b_i$, for any $i \in I_{t(\nu(\ell))}$. In the case $M = 1$, $I_t = I$, i.e., when all equations are lumped into a single block, the proof would have been completed at this point. In general it is necessary to show that, for any $r \in I$, $\lim_{\ell \to \infty} \langle a^r, x^{\nu(\ell)} \rangle = b_r$. To this end we shift indices in a particular manner. For an arbitrary $r \in I$, we can write

$$\langle a^r, x^{\nu(\ell)} \rangle - b_r = \sum_{j=1}^n a_j^r x_j^s h_j\left(\nu(\ell), s\right) - b_r, \tag{6.140}$$

where each *shift function* $h_j$ is defined by the product

$$h_j(\nu, s) \doteq \prod_{z=s+1}^{\nu} \left( \frac{x_j^z}{x_j^{z-1}} \right), \tag{6.141}$$

for every $j = 1, 2, \ldots, n$. This is well defined because of (6.126). The index $s$, to be used in (6.140), is chosen in the following way: there exists an index $t$, $1 \leq t \leq M$, such that $r \in I_t$; and for every iteration index $\nu = \nu(\ell)$ there exists an earlier index $s = s(\nu)$ such that

$$s \leq \nu, \quad \nu - s \leq C, \quad \text{and} \quad t = t(s), \tag{6.142}$$

i.e., such that $r \in I_{t(s)}$, where $C$ is the constant of almost cyclicality of the control sequence of Algorithm 6.7.1. Using Assumptions 6.6.1–6.6.3, (6.118), and (6.125) the shift functions may be estimated by

$$\begin{aligned} h_j(\nu, s) &= \prod_{z=s+1}^{\nu} \exp\left( \sum_{i \in I_{t(z-1)}} (w_i^{t(z-1)} d_i^{z-1} a_j^i) \right) \\ &\leq \exp\left( \sum_{z=s+1}^{\nu} \sum_{i \in I_{t(z-1)}} d_i^{z-1} \right). \end{aligned} \tag{6.143}$$

Since the sums on the right-hand side of (6.143) range over no more than $C$ and $m$ indices, respectively, we may further write

$$h_j(\nu, s) \leq \exp\ \left( mC \cdot \max\left\{ d_i^{z-1} \mid i \in I_{t(z-1)},\ z \geq s+1 \right\} \right). \tag{6.144}$$

By (6.122) and (6.139) we have that $\lim_{q \to \infty} d_i^q = 0$ for all $i \in I_{t(q)}$, i.e., for all $i$ in the "right" block. Thus,

$$\lim_{\nu \to \infty} h_j(\nu, s) = 1, \quad \text{for all } \ j = 1, 2, \ldots, n. \tag{6.145}$$

From all the above we now get,

$$\begin{aligned} &\mid \langle a^r, x^{\nu(\ell)} \rangle - b_r \mid \\ &\leq \mid \langle a^r, x^s \rangle - b_r \mid + \langle a^r, x^s \rangle \cdot \max_{1 \leq j \leq n} \mid h_j(\nu(\ell), s) - 1 \mid, \end{aligned} \tag{6.146}$$

and as $\ell \to \infty$, so does $s$, and therefore the right-hand side of (6.146) tends to zero by (6.139) and (6.145). ■

Our analysis culminates in the following convergence theorem.

**Theorem 6.7.1** *Under Assumptions 6.6.1–6.6.3 any sequence* $\{x^\nu\}$ *generated by Algorithm 6.7.1, with blocks and weights defined as in (6.117)–*

*(6.118), converges to the unique solution of the problem (6.106)–(6.107) with $Q = Q_1$ as defined by (6.108).*

**Proof** From (6.23), from the continuity of the gradient of $f$, and from (6.137) we deduce that,

$$\begin{aligned}\lim_{\ell\to\infty}\langle \pi^{\nu(\ell)}, Ax^{\nu(\ell)} - b\rangle &= \lim_{\ell\to\infty}\langle -\nabla f(x^{\nu(\ell)}), x^{\nu(\ell)} - x^*\rangle \\ &= -\langle \nabla f(x^*), x^* - x^*\rangle = 0. \qquad (6.147)\end{aligned}$$

By Proposition 6.7.2, $\{L_\nu\}$ converges, and in view of (6.147) its limit is

$$\lim_{\nu\to\infty} L_\nu = \lim_{\ell\to\infty} L_{\nu(\ell)} = \lim_{\ell\to\infty}(f(x^{\nu(\ell)}) + \langle \pi^{\nu(\ell)}, Ax^{\nu(\ell)} - b\rangle) = f(x^*), \tag{6.148}$$

implying by (6.136) that $f(x^*) \leq f(y)$, for any $y \in Q_1 \cap \mathbb{R}^n_+$. This proves that $x^*$ is an optimal solution for the problem. Strict convexity of the function $f(x) = -\operatorname{ent} x$ implies uniqueness of $x^*$, thus it is the limit of the whole sequence $\{x^\nu\}$. ■

## 6.8 Underrelaxation Parameters and Extension of the Family of Bregman Functions

Relaxation parameters allow us to relax any projection operation. When they are incorporated into an iterative algorithm they add to its implementation the extra flexibility of *overdoing* or *underdoing* the iterative step. We have seen earlier (Chapter 5) how this works in algorithms that employ orthogonal projections for the convex feasibility problem. Namely, if $y = P_Q(x)$ is the orthogonal projection of $x$ onto the set $Q$, then $y(\lambda) = x + \lambda(P_Q(x) - x)$ is the relaxed projection with relaxation parameter $\lambda$.

Our first question is how to introduce similar relaxation parameters when a generalized projection onto a hyperplane or a half-space is involved. Another question that we wish to raise here, which is completely independent of the first, is related to the family of Bregman functions $\mathcal{B}(S)$. Consider the following function, called *Burg's entropy*, $B(x)$. It maps the positive orthant $\mathbb{R}^n_{++} \doteq \{x \in \mathbb{R}^n \mid x_j > 0,\ j = 1, 2, \ldots, n\}$ into $\mathbb{R}$ according to

$$B(x) \doteq \sum_{j=1}^{n} \log x_j. \tag{6.149}$$

The function $f(x) = -B(x)$ is not a Bregman function with zone $S = \mathbb{R}^n_{++}$ because it becomes singular on the boundary $\operatorname{bd} S = \mathbb{R}^n_+ \backslash \mathbb{R}^n_{++}$, i.e., when $x_j$ tends to zero for even only one $j$, then the function $-B(x)$ tends to $-\infty$, demonstrating an essential discontinuity.

Is it possible to circumvent this difficulty and apply Bregman's algorithm to this function $f(x) = -B(x)$? If so, can the theory be extended to handle any function with singularities of this sort on the boundary of $\overline{S}$?

In this section we will give some affirmative and constructive answers to these questions. The key to introducing a relaxation parameter is to define an underrelaxed generalized projection of a point $x$ onto a given hyperplane $H$ as a generalized unrelaxed projection of $x$ onto some hyperplane that is parallel to $H$ and lies between $x$ and $H$. In order to get around the singularity problem we define the *singularity set* $T$ of the Bregman function $f \in \mathcal{B}(S)$ and require that the partial level sets of the associated generalized distance $D_f$ be both bounded and *bounded away* from the singularity set $T$ of $f$.

Although the introduction of underrelaxation parameters and the extension of the family of Bregman functions are independent issues, we study them here together. First, we define what it means to say that one set is bounded away from another set.

**Definition 6.8.1 (Bounded Away Sets)** *Given a set $V \subseteq \mathbb{R}^n$ and a set $W \subseteq \mathbb{R}^n$, we say that $W$ is bounded away from $V$ if there exists an open convex set $U$ with closure $\overline{U}$, such that $W \subseteq U$ and $\overline{U} \cap V = \emptyset$.*

In order to define an extended family of Bregman functions we need a set $\widetilde{S}$ that lies between the zone $S$ and its closure $\overline{S}$. Let S be a nonempty, open, convex set and $\widetilde{S}$ a convex set such that $S \subseteq \widetilde{S} \subseteq \overline{S}$ and $\widetilde{S} \subseteq \Lambda$, where $\Lambda$ is the domain of a function $f : \Lambda \subseteq \mathbb{R}^n \to \mathbb{R}$. Assume that $f(x)$ has continuous first partial derivatives at any $x \in S$. For any two sets $A$ and $B$ we denote set subtraction by $A \backslash B \doteq \{x \mid x \in A \text{ and } x \notin B\}$. In the following definition the function $D_f(x, y)$ is as in (2.1) and the partial level sets $L_1^f(y, \alpha)$ and $L_2^f(x, \alpha)$ are as in (2.3).

**Definition 6.8.2 (Extended Bregman Functions)** *We call a function $f : \Lambda \subseteq \mathbb{R}^n \to \mathbb{R}$ an extended Bregman function with zone $S$, and use the notation $f \in \mathcal{EB}(S)$ if there exists a nonempty, open, convex set $S$ and a convex set $\widetilde{S}$ such that $S \subseteq \widetilde{S} \subseteq \overline{S}$, $\widetilde{S} \subseteq \Lambda$ and the following conditions hold:*

(*i*) *f(x) has continuous first partial derivatives at every $x \in S$,*

(*ii*) *f(x) is strictly convex on $\widetilde{S}$,*

(*iii*) *f(x) is continuous on $\widetilde{S}$,*

(*iv*) *For every $\alpha \in \mathbb{R}$, the partial level sets $L_1^f(y, \alpha)$ and $L_2^f(x, \alpha)$ are (a) bounded, and (b) bounded away from $T \doteq \overline{S} \backslash \widetilde{S}$, for every $y \in S$, for every $x \in \widetilde{S}$, respectively,*

(*v*) *If $y^\nu \in S$, for $\nu \geq 0$, and $\lim_{\nu \to \infty} y^\nu = y^* \in \widetilde{S}$ then $\lim_{\nu \to \infty} D_f(y^*, y^\nu) = 0$,*

(*vi*) *If $y^\nu \in S$ and $x^\nu \in \widetilde{S}$, for all $\nu \geq 0$, and if $\lim_{\nu \to \infty} D_f(x^\nu, y^\nu) = 0$, and*

$\lim_{\nu \to \infty} y^\nu = y^* \in \widetilde{S}$, *and* $\{x^\nu\}$ *is bounded, then* $\lim_{\nu \to \infty} x^\nu = y^*$.

This definition extends the previous definition of Bregman functions (Definition 2.1.1) through the introduction of the set $\widetilde{S}$ and condition $(iv)(b)$. Every Bregman function is also an extended Bregman function with $\widetilde{S} = \overline{S}$, but not vice versa. The extension enables us to apply Bregman's method for convex programming to functions that have essential discontinuities at points of the boundary of their zone $S$. The set $T \doteq \overline{S} \backslash \widetilde{S}$, which is taken then as the set of all such singular points and which we call the *singularity set* of $f$, may be the whole boundary of $S$ (i.e., $\widetilde{S} = S$), part of the boundary, or empty, in which case $\widetilde{S} = \overline{S}$. Condition $(iv)(b)$ in Definition 6.8.2 guarantees that all sequences of interest will not approach any point of $T$, thereby allowing the extension of Bregman functions with such singularities. In analogy with Definition 2.1.2 we now define the following:

**Definition 6.8.3 (Extended Projections)** *Given a function* $f \in \mathcal{EB}(S)$ *and a set* $\Omega \subseteq \mathbb{R}^n$ *such that* $\Omega \cap \widetilde{S} \neq \emptyset$, *let* $P_\Omega : S \to \widetilde{S}$ *be defined by*

$$P_\Omega(y) \doteq \operatorname*{argmin}_{x \in \Omega \cap \widetilde{S}} D_f(x, y), \tag{6.150}$$

*and call* $P_\Omega(y)$ *the extended projection of* $y$ *onto* $\Omega$.

**Lemma 6.8.1** *If* $\Omega \subseteq \mathbb{R}^n$ *is a closed and convex set then* $P_\Omega(y)$, *as defined in Definition 6.8.3, exists and is unique.*

**Proof** The proof is similar to that of Lemma 2.1.2. Given $y \in S$, the set $\Omega \cap L_1^f(y, \alpha)$ is convex and compact for $\alpha = D_f(w, y)$, where $w$ is any point that belongs to $\Omega \cap \widetilde{S}$. Convexity is true because $D_f(x, y)$, as a function of $x$, is a sum of a strictly convex function and an affine function; thus it is strictly convex, implying that its level sets $L_1^f(y, \alpha)$ are convex. The set $\Omega$ is convex by definition. Regarding compactness, observe that $L_1^f(y, \alpha)$ may be rewritten as

$$L_1^f(y, \alpha) = \{x \in \overline{S} \cap \overline{U} \mid D_f(x, y) \leq D_f(w, y)\}, \tag{6.151}$$

where $U$ is the set whose existence is guaranteed by $(iv)(b)$ of Definition 6.8.2. $D_f(x, y)$ is a continuous function of $x$ on $\widetilde{S}$, while we have $\overline{S} \cap \overline{U} = \widetilde{S} \cap \overline{U} \subseteq \widetilde{S}$. These facts imply that $L_1^f(y, \alpha)$ is closed and thus, by $(iv)(a)$ of Definition 6.8.2, compact. Therefore, there exists a unique $x^*$ that minimizes $D_f(x, y)$ over $\Omega \cap L_1^f(y, \alpha)$, i.e.,

$$D_f(x^*, y) \leq D_f(w, y) = \alpha, \tag{6.152}$$

for every $w \in \Omega \cap L_1^f(y, \alpha)$. By the definition of $L_1^f(y, \alpha)$, we have that $\alpha < D_f(x, y)$ for every $x \in \Omega \cap (\widetilde{S} \backslash L_1^f(y, \alpha))$, yielding that $x^*$ minimizes

$D_f(x,y)$ over $\Omega \cap \widetilde{S}$. ■

**Definition 6.8.4 (Zone Consistency of Extended Bregman Functions)**

(i) *A function $f \in \mathcal{EB}(S)$ is said to be zone consistent with respect to the hyperplane $H$ if, for every $y \in S$ we have $P_H(y) \in S$.*

(ii) *A function $f \in \mathcal{EB}(S)$ is said to be strongly zone consistent with respect to the hyperplane $H$ if it is zone consistent with respect to $H$ as well as with respect to every hyperplane $H'$, parallel to $H$, which lies between $y$ and $H$.*

This extends Definition 2.2.1. Observe, however, that $P_H(y)$ always belongs to $\widetilde{S}$, so that zone consistency actually means only that the extended projection will not be on the boundary of $S$.

Given a hyperplane

$$H = \{x \in \mathbb{R}^n \mid \langle a, x\rangle = b\}, \tag{6.153}$$

and a function $f \in \mathcal{EB}(S)$, the first order optimality conditions for the problem

$$x^* = \operatorname*{argmin}_{x \in H \cap \widetilde{S}} D_f(x,y), \tag{6.154}$$

where $y \in S$, may be written in the form (see Lemma 2.2.1):

$$\nabla f(x^*) = \nabla f(y) + sa, \tag{6.155}$$

$$\langle a, x^* \rangle = b, \tag{6.156}$$

$$x^* \in \widetilde{S}, \tag{6.157}$$

where $s$ is the Lagrange multiplier associated with equation (6.156). These conditions are both necessary and sufficient because of the strict convexity of $D_f(x,y)$ as a function of $x$. It follows from Lemma 6.8.1 that, given $a \in \mathbb{R}^n$ and $b \in \mathbb{R}$, there exists a unique pair $(x^*, s^*)$ that satisfies (6.155)–(6.157). $x^* = P_H(y)$ and we denote $s^* = \pi_H(y)$ and call it the *projection parameter* associated with the extended projection of $y$ onto $H$. We mention that there may be other solutions $(x, s)$ to (6.155)–(6.156) that fail to satisfy (6.157). Such cases will surface later in our discussion.

In order to present the underrelaxed Bregman algorithm for minimizing a function $f \in \mathcal{EB}(S)$ over a set of linear inequalities, we consider the problem

$$\text{Minimize} \quad f(x) \tag{6.158}$$

$$\text{s.t.} \quad Ax \le b, \tag{6.159}$$

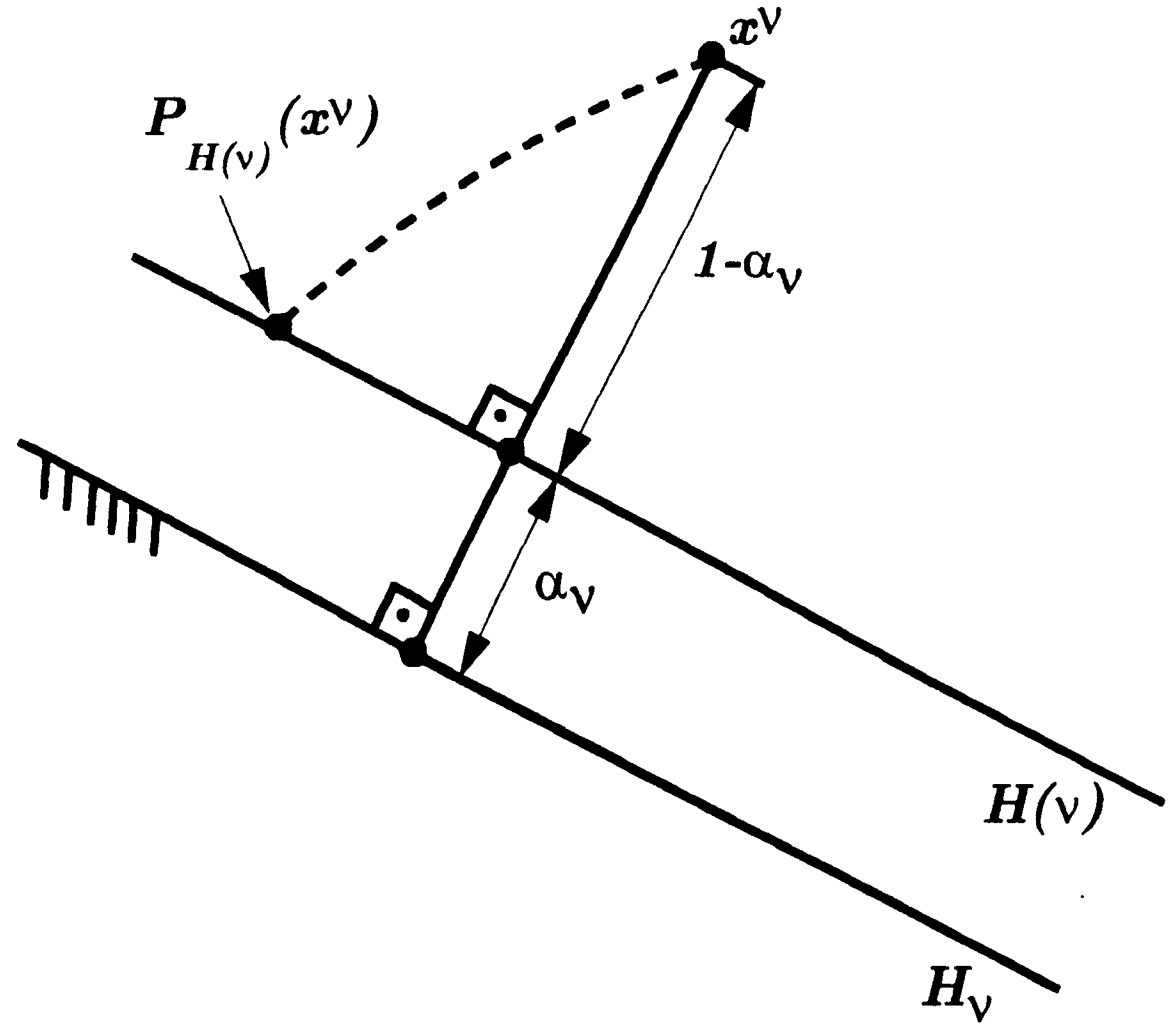

**Figure 6.5** The extended projection onto the relaxed hyperplane $H(\nu)$ is considered the relaxed projection of $x^\nu$ onto the original hyperplane $H_\nu$.

$$x \in \widetilde{S}, \tag{6.160}$$

where $f \in \mathcal{EB}(S)$, $\widetilde{S}$ is as in Definition 6.8.2, and the real $m \times n$ matrix $A$ and vector $b \in \mathbb{R}^m$ are given. We denote $I \doteq \{1, 2, \ldots, m\}$ and

$$Q \doteq \{x \in \mathbb{R}^n \mid Ax \leq b\} = \{x \in \mathbb{R}^n \mid \langle a^i, x\rangle \leq b_i,\ i \in I\}, \tag{6.161}$$

where $a^i \in \mathbb{R}^n$ is the $i$th column of $A^T$, the transpose of $A$, and assume that $a^i \neq 0$, for all $i \in I$, and that $\widetilde{S} \cap Q \neq \emptyset$. Let $\epsilon > 0$ be arbitrarily small but fixed and let $\{\alpha_\nu\}_{\nu=0}^\infty$ be a sequence of relaxation parameters, such that

$$\epsilon \leq \alpha_\nu \leq 1, \qquad \text{for all } \nu \geq 0. \tag{6.162}$$

Underrelaxation is accomplished in the following manner. Recall that the relaxed orthogonal projection of a point $x^\nu$ onto a hyperplane

$$H_\nu \doteq \{x \in \mathbb{R}^n \mid \langle a^{i(\nu)}, x\rangle = b_{i(\nu)}\}, \tag{6.163}$$

is the point $x^{\nu+1}$ given by

$$x^{\nu+1} = x^\nu + \lambda_\nu \frac{b_{i(\nu)} - \langle a^{i(\nu)}, x^\nu\rangle}{\| a^{i(\nu)} \|^2} a^{i(\nu)}. \tag{6.164}$$

For every possible value of the relaxation parameter $\lambda_\nu$ within the interval $0 \leq \lambda_\nu \leq 2$, the point $x^{\nu+1}$ will be located along the line through $x^\nu$ and its orthogonal projection $P_{H_\nu}(x^\nu)$ onto $H_\nu$

$$P_{H_\nu}(x^\nu) = x^\nu + \frac{b_{i(\nu)} - \langle a^{i(\nu)}, x^\nu \rangle}{\| a^{i(\nu)} \|^2} a^{i(\nu)}. \tag{6.165}$$

The underrelaxation, with respect to the hyperplane $H_\nu$, in the case of a generalized projection is obtained by constructing a *relaxed hyperplane*

$$H(\nu) \doteq \{x \in \mathbb{R}^n \mid \langle a^{i(\nu)}, x \rangle = \alpha_\nu b_{i(\nu)} + (1 - \alpha_\nu)\langle a^{i(\nu)}, x^\nu \rangle\}. \tag{6.166}$$

This hyperplane is parallel to $H_\nu$ and, for $0 \leq \alpha_\nu \leq 1$ lies between the point $x^\nu$ and $H_\nu$. The relaxed extended projection of $x^\nu$ onto $H_\nu$ with relaxation parameter $\alpha_\nu$ is then defined to be the unrelaxed extended projection of $x^\nu$ onto $H(\nu)$ (see Figure 6.5). The complete formulation of the algorithm is given next.

**Algorithm 6.8.1 Underrelaxed Algorithm for Extended Bregman Functions Over Linear Inequality Constraints**

**Step 0:** (Initialization.) $x^0 \in Z_0$ is arbitrary, where $Z_0$ is defined as in (6.16), and $z^0$ is such that $\nabla f(x^0) = -A^T z^0$.

**Step 1:** (Iterative step.) Given $\nu, x^\nu$ and $z^\nu$, take an $\alpha_\nu$ such that (6.162) holds and pick a control index $i(\nu)$ and calculate $x^{\nu+1}$ and $z^{\nu+1}$ from

$$\nabla f(x^{\nu+1}) = \nabla f(x^\nu) + c_\nu a^{i(\nu)}, \tag{6.167}$$
$$z^{\nu+1} = z^\nu - c_\nu e^{i(\nu)}, \tag{6.168}$$

where

$$c_\nu \doteq \min(z^\nu_{i(\nu)}, \beta_\nu), \tag{6.169}$$
$$\beta_\nu = \pi_{H(\nu)}(x^\nu), \tag{6.170}$$

in which the hyperplane $H(\nu)$ is given by (6.166).

**Control:** The sequence $\{i(\nu)\}_{\nu=0}^\infty$ is almost cyclic on $I$.

Note that if $c_\nu = \beta_\nu$ then $x^{\nu+1} = P_{H(\nu)}(x^\nu)$. If $c_\nu \neq \beta_\nu$ then $x^{\nu+1}$ is also well defined since, as shown in Proposition 6.8.2 below, in that case $x^{\nu+1}$ is the extended projection of $x^\nu$ onto a hyperplane parallel to $H(\nu)$ that lies between it and $x^\nu$. Strong zone consistency, which is assumed, guarantees that $x^{\nu+1} \in S$ in every case.

The convergence analysis of Algorithm 6.8.1 is similar to the analysis of Algorithm 6.3.1. The replacement of $\overline{S}$ by $\widetilde{S}$ and the introduction of underrelaxation calls for certain modifications pointed out below. The assumptions used here to guarantee convergence are as follows.

**Assumption 6.8.1**

(*i*) $f \in \mathcal{EB}(S)$,

*(ii)* $f$ *is strongly zone consistent with respect to each* $H_i$, $i \in I$ *where* $H_i = \{x \in \mathbb{R}^n \mid \langle a^i, x\rangle = b_i\}$, *are the bounding hyperplanes of the inequality constraints in (6.159),*

*(iii) The sequence* $\{i(\nu)\}_{\nu=0}^{\infty}$ *is almost cyclic on* $I$, *and*

*(iv)* $S \cap Q \neq \emptyset$.

The complementarity properties of the dual iterates are given by the next proposition.

**Proposition 6.8.1** *If Assumption 6.8.1 holds then the iterates* $\{x^\nu\}$ *and* $\{z^\nu\}$ *generated by Algorithm 6.8.1 have, for all* $\nu \geq 0$, *the following properties: (i)* $z^\nu \geq 0$, *and (ii)* $\nabla f(x^\nu) = -A^T z^\nu$.

**Proof** Same as in the proof of Step 1 of Theorem 6.3.1. ■

These properties are used to show that $x^{\nu+1}$ of Algorithm 6.8.1 is indeed an extended projection of $x^\nu$ onto a hyperplane parallel to $H(\nu)$.

**Proposition 6.8.2** *If Assumption 6.8.1 holds then* $x^{\nu+1} = P_{\widetilde{H}(\nu)}(x^\nu)$ *and* $c_\nu = \pi_{\widetilde{H}(\nu)}(x^\nu)$, *where*

$$\widetilde{H}(\nu) \doteq \{x \in \mathbb{R}^n \mid \langle a^{i(\nu)}, x\rangle = \gamma_\nu b_{i(\nu)} + (1-\gamma_\nu)\langle a^{i(\nu)}, x^\nu\rangle\}, \tag{6.171}$$

*for some* $\gamma_\nu$ *such that* $0 \leq \gamma_\nu \leq \alpha_\nu$.

**Proof** Note that $x^{\nu+1}$ is the extended projection of $x^\nu$ on the hyperplane $\{x \in \mathbb{R}^n \mid \langle a^{i(\nu)}, x\rangle = \langle a^{i(\nu)}, x^{\nu+1}\rangle\}$, which is parallel to $H_\nu$ and passes through $x^{\nu+1}$. It remains to show that the right-hand side $\langle a^{i(\nu)}, x^{\nu+1}\rangle$, has the desired form as in (6.171). It is clear that if $c_\nu = \beta_\nu$, then we have $\widetilde{H}(\nu) = H(\nu)$, and we may take $\gamma_\nu = \alpha_\nu$. On the other hand, if $c_\nu \neq \beta_\nu$, then by (6.169) and Proposition 6.8.1

$$0 \leq z^\nu_{i(\nu)} = c_\nu < \beta_\nu. \tag{6.172}$$

Now,

$$\begin{aligned} 0 \leq D_f(x^{\nu+1}, x^\nu) + D_f(x^\nu, x^{\nu+1}) &= \langle \nabla f(x^{\nu+1}) - \nabla f(x^\nu), x^{\nu+1} - x^\nu\rangle \\ &= c_\nu\langle a^{i(\nu)}, x^{\nu+1} - x^\nu\rangle. \end{aligned} \tag{6.173}$$

If $c_\nu = 0$, then $x^{\nu+1} = x^\nu$ by (6.173) and Lemma 2.1.1 and we may take $\gamma_\nu = 0$. From (6.172) the only remaining case is $0 < c_\nu < \beta_\nu$, i.e., $0 < \pi_{\widetilde{H}(\nu)}(x^\nu) < \pi_{H(\nu)}(x^\nu)$. By Lemma 2.2.4,

$$\langle a^{i(\nu)}, x^\nu\rangle < \langle a^{i(\nu)}, x^{\nu+1}\rangle < \alpha_\nu b_{i(\nu)} + (1-\alpha_\nu)\langle a^{i(\nu)}, x^\nu\rangle$$

and there is a $\gamma_\nu \in (0, \alpha_\nu)$ as desired. ■

The Lagrangian of the minimization problem (6.158)–(6.160) is given by $L(x,z) = f(x) + \langle z, Ax - b\rangle$ and the following holds.

**Proposition 6.8.3** *For any $\{x^\nu\}$ and $\{z^\nu\}$ produced by Algorithm 6.8.1, the sequence $L(x^\nu, z^\nu)$ is nondecreasing.*

**Proof** Defining $d_\nu \doteq L(x^{\nu+1}, z^{\nu+1}) - L(x^\nu, z^\nu)$, we first show that

$$D_f(x^{\nu+1}, x^\nu) + (1-\gamma_\nu) D_f(x^\nu, x^{\nu+1}) = \gamma_\nu d_\nu. \tag{6.174}$$

To prove (6.174) we note, using the iterative step formulae of Algorithm 6.8.1, that

$$\begin{aligned}
\gamma_\nu d_\nu &= \gamma_\nu(\langle A^{\mathrm{T}}(z^{\nu+1} - z^\nu), x^{\nu+1}\rangle - \langle z^{\nu+1} - z^\nu, b\rangle + D_f(x^{\nu+1}, x^\nu)) \\
&= \gamma_\nu(\langle \nabla f(x^\nu) - \nabla f(x^{\nu+1}), x^{\nu+1}\rangle - \langle z^{\nu+1} - z^\nu,\ b\rangle + D_f(x^{\nu+1}, x^\nu)) \\
&= \gamma_\nu(c_\nu(b_{i(\nu)} - \langle a^{i(\nu)}, x^{\nu+1}\rangle) + D_f(x^{\nu+1}, x^\nu)).
\end{aligned} \tag{6.175}$$

By Proposition 6.8.2, $x^{\nu+1} = P_{\widetilde{H}(\nu)}(x^\nu)$, so we must have $\langle a^{i(\nu)}, x^{\nu+1}\rangle = \gamma_\nu b_{i(\nu)} + (1-\gamma_\nu)\langle a^{i(\nu)}, x^\nu\rangle$, which implies

$$\gamma_\nu \langle a^{i(\nu)}, x^{\nu+1}\rangle = \gamma_\nu b_{i(\nu)} + (1-\gamma_\nu)\langle a^{i(\nu)}, x^\nu - x^{\nu+1}\rangle.$$

Thus,

$$\gamma_\nu(b_{i(\nu)} - \langle a^{i(\nu)}, x^{\nu+1}\rangle) = (1-\gamma_\nu)\langle a^{i(\nu)}, x^{\nu+1} - x^\nu\rangle. \tag{6.176}$$

Substituting (6.176) into (6.175), we get

$$\begin{aligned}
\gamma_\nu d_\nu &= (1-\gamma_\nu) c_\nu \langle a^{i(\nu)}, x^{\nu+1} - x^\nu\rangle + \gamma_\nu D_f(x^{\nu+1}, x^\nu) \\
&= (1-\gamma_\nu)\langle \nabla f(x^{\nu+1}) - \nabla f(x^\nu), x^{\nu+1} - x^\nu\rangle + \gamma_\nu D_f(x^{\nu+1}, x^\nu) \\
&= (1-\gamma_\nu)(D_f(x^{\nu+1}, x^\nu) + D_f(x^\nu, x^{\nu+1})) + \gamma_\nu D_f(x^{\nu+1}, x^\nu) \\
&= (1-\gamma_\nu) D_f(x^\nu, x^{\nu+1}) + D_f(x^{\nu+1}, x^\nu).
\end{aligned} \tag{6.177}$$

Now we conclude the proof by showing that $d_\nu \geq 0$, for all $\nu \geq 0$. If $\gamma_\nu = 0$, then from Proposition 6.8.2 $x^{\nu+1} = x^\nu$ so that $c_\nu = 0$ and $z^{\nu+1} = z^\nu$. Then $d_\nu = L(x^{\nu+1}, z^{\nu+1}) - L(x^\nu, z^\nu) = 0$. Otherwise, from Proposition 6.8.2, $0 < \gamma_\nu \leq \alpha_\nu \leq 1$ and the result follows from (6.174) and Lemma 2.1.1. ■

Next we show that $\{L(x^\nu, z^\nu)\}_{\nu=0}^\infty$ is also a bounded sequence.

**Proposition 6.8.4** *For any $\{x^\nu\}$ and $\{z^\nu\}$ produced by Algorithm 6.8.1,*

$$L(x^\nu, z^\nu) \leq f(z) - D_f(z, x^\nu) \leq f(z), \text{ for any } z \in Q \cap \widetilde{S}. \tag{6.178}$$

**Proof** Similar to the proof of (6.46) where we had $z \in Q \cap \overline{S}$. ■

Propositions 6.8.3 and 6.8.4 show that $\{L(x^\nu, z^\nu)\}_{\nu=0}^\infty$ is a convergent sequence. Moreover, the next result holds for the sequence $\{x^\nu\}$.

**Proposition 6.8.5** *Under the same assumptions, any sequence $\{x^\nu\}$, produced by Algorithm 6.8.1, is bounded and bounded away from the set $T = \overline{S}\backslash\widetilde{S}$.*

**Proof** Using Propositions 6.8.3 and 6.8.4 recursively we get, for any $z \in Q \cap \widetilde{S}$, that $D_f(z, x^\nu) \leq f(z) - L(x^0, z^0) \doteq \tau$, so that $x^\nu \in L_2^f(z, \tau)$, which completes the proof by Definition 6.8.2($iv$). ■

Note that the fact that the relaxation parameters are bounded away from zero, i.e., $0 < \epsilon \leq \alpha_\nu$, was not used up to this point, and that Proposition 6.8.5 holds even if they satisfy only $0 < \alpha_\nu \leq 1$. In the reminder of the proof, however, the boundedness away from zero of the relaxation parameters is required.

**Proposition 6.8.6** *All accumulation points of $\{x^\nu\}$ belong to $Q \cap \widetilde{S}$.*

**Proof** Any accumulation point $x^*$ of $\{x^\nu\}$ must be in $\widetilde{S}$ because of Proposition 6.8.5. In order to show that $x^* \in Q$ we first take a convergent subsequence $\{x^{\nu_k}\}_{k=0}^\infty$ of the sequence $\{x^\nu\}_{\nu=0}^\infty$, i.e., $\lim_{k\to\infty} x^{\nu_k} = x^*$. We also take a fixed integer $t$ and a sequence $\{l_k\}_{k=0}^\infty$ with $l_k \in \{1, 2, \ldots, t\}$ for all $k \geq 0$. Under these conditions, we claim that

$$\lim_{k\to\infty} x^{\nu_k + l_k} = x^*. \tag{6.179}$$

To see this, consider first the $t$ sequences $\{x^{\nu_k + l}\}_{k=0}^\infty$ with $1 \leq l \leq t$. Since they are subsequences of $\{x^\nu\}_{\nu=0}^\infty$ they are all bounded. From Lemma 2.1.1 and Proposition 6.8.3, since $\gamma_\nu \leq \alpha_\nu \leq 1$ for all $\nu$ we have

$$0 \leq D_f(x^{\nu+1}, x^\nu) \leq \gamma_\nu d_\nu \leq d_\nu, \tag{6.180}$$

and we know that $\lim_{\nu\to\infty} d_\nu = 0$ because $\{L(x^\nu, z^\nu)\}_{\nu=0}^\infty$ is convergent. This implies that

$$\lim_{k\to\infty} D_f(x^{\nu_k+s+1}, x^{\nu_k+s}) = 0, \tag{6.181}$$

for all $s \in \{1, 2, \ldots, t\}$. Applying recursively Definition 6.8.2($vi$), we conclude that for all $l$, $0 \leq l \leq t$,

$$\lim_{k\to\infty} x^{\nu_k + l} = x^*. \tag{6.182}$$

Interlacing these $t + 1$ sequences we can form the sequence

$$\{x^{\nu_1}, x^{\nu_1+1}, \ldots, x^{\nu_1+t}, x^{\nu_2}, x^{\nu_2+1}, \ldots, x^{\nu_2+t}, \ldots, x^{\nu_k}, x^{\nu_k+1}, \ldots, x^{\nu_k+t}, \ldots\},$$

which converges to $x^*$ and of which the sequence of (6.179) is a subsequence, thus validating (6.179). We use this to show that $x^* \in Q$. Take any $p \in$ $\{1, 2, \ldots, m\}$ and $l_k \in \{1, 2, \ldots, C\}$, where $m$ is the number of inequalities in (6.159) and $C$ is the constant of almost cyclicality of the control sequence $\{i(\nu)\}_{\nu=0}^{\infty}$, such that

$$i(\nu_k + l_k) = p, \tag{6.183}$$

and use (6.179) with $t = C$. Then, $\lim_{k\to\infty} x^{\nu_k+l_k} = x^*$ holds and we may take a subsequence $\{x^{s_k}\}_{k=0}^{\infty}$ of $\{x^{\nu_k+l_k}\}_{k=0}^{\infty}$ such that

$$\lim_{k\to\infty} \gamma_{s_k} = \gamma, \quad \text{and} \quad \lim_{k\to\infty} \alpha_{s_k} = \alpha \geq \epsilon, \tag{6.184}$$

would hold. By (6.183) we have,

$$\langle a^p, x^{s_k+1}\rangle = \gamma_{s_k} b_p + (1 - \gamma_{s_k})\langle a^p, x^{s_k}\rangle.$$

Taking limits as $k \to \infty$, we obtain $\langle a^p, x^*\rangle = \gamma b_p + (1-\gamma)\langle a^p, x^*\rangle$, which implies $\gamma(\langle a^p, x^*\rangle - b_p) = 0$. If $\gamma \neq 0$ then

$$\langle a^p, x^*\rangle = b_p. \tag{6.185}$$

If $\gamma = 0$, we have $\gamma_{s_k} \neq \alpha_{s_k}$, for large enough $k$, because of (6.184). Thus,

$$\langle a^p, x^{s_k+1}\rangle < \alpha_{s_k} b_p + (1 - \alpha_{s_k})\langle a^p, x^{s_k}\rangle, \tag{6.186}$$

because of Proposition 6.8.2. Taking now limits as $k \to \infty$ leads to

$$\langle a^p, x^*\rangle \leq \alpha b_p + (1-\alpha)\langle a^p, x^*\rangle, \tag{6.187}$$

which implies $0 \leq \alpha(b_p - \langle a^p, x^*\rangle)$ yielding

$$\langle a^p, x^*\rangle \leq b_p. \tag{6.188}$$

From (6.185) and (6.188) we infer that $x^*$ must satisfy the $p$th constraint, but since $p$ is arbitrary we conclude that $x^* \in Q$. ■

This leads to the final conclusion on convergence, stated in the following theorem.

**Theorem 6.8.1** *Under Assumption 6.8.1 any sequence $\{x^\nu\}$ generated by Algorithm 6.8.1 converges to an $x^* \in \widetilde{S}$ which is a solution of (6.158)–(6.160).*

**Proof** From Proposition 6.8.5 we deduce that there is a convergent subsequence $\lim_{k\to\infty} x^{\nu_k} = x^*$ of $\{x^\nu\}_{\nu=0}^{\infty}$, and from Proposition 6.8.6 we know that $x^* \in Q \cap \widetilde{S}$. Using the notation of (6.49), it is possible to show that

$$z_p^{\nu_k+C+1} = 0, \text{ for all } p \in I_1, \tag{6.189}$$

where $C$ is the constant of almost cyclicality.

To see this, let

$$\rho \doteq (\epsilon/5) \cdot \min_{p \in I_1}\{(b_p - \langle a^p, x^* \rangle)/ \parallel a^p \parallel\} > 0. \tag{6.190}$$

By (6.179), with $t = C + 1$, there exists an index $k_0$ such that

$$\parallel x^{\nu_k + l} - x^* \parallel < \rho, \text{ for all } l \in \{0, 1, \ldots, C+1\}, \text{ and } k \geq k_0. \tag{6.191}$$

Given $p \in I_1$, define

$$l_k \doteq \max_{0 \leq l \leq C}\{l \mid i(\nu_k + l) = p\}, \tag{6.192}$$

the existence of which is guaranteed by almost cyclicality. Let $s_k = l_k + \nu_k$ and assume that $c_{s_k} = \beta_{s_k}$ in (6.169). Then

$$\langle a^p, x^{s_k+1} \rangle = \alpha_{s_k} b_p + (1 - \alpha_{s_k})\langle a^p, x^{s_k} \rangle, \tag{6.193}$$

from which we get that, for all $k \geq k_0$,

$$\begin{aligned} \alpha_{s_k}(b_p - \langle a^p, x^{s_k} \rangle) &< \langle a^p, x^{s_k+1} - x^{s_k} \rangle + \alpha_{s_k}\langle a^p, x^{s_k} - x^* \rangle \\ &\leq \parallel a^p \parallel (\parallel x^{s_k+1} - x^{s_k} \parallel + \alpha_{s_k} \parallel x^{s_k} - x^* \parallel) \\ &\leq \parallel a^p \parallel (\parallel x^{s_k+1} - x^* \parallel + (1 + \alpha_{s_k}) \parallel x^{s_k} - x^* \parallel) \\ &\leq \parallel a^p \parallel (2 + \alpha_{s_k})\rho \leq 4\rho \parallel a^p \parallel . \end{aligned} \tag{6.194}$$

In this way the contradiction $4\rho > \epsilon(b_p - \langle a^p, x^* \rangle)/ \parallel a^p \parallel \geq 5\rho$ is reached, which enforces that, contrary to what we tried to assume before, $c_{s_k} \neq \beta_{s_k}$ in (6.169) implying $c_{s_k} = z_p^{s_k}$ , and thus, $z_p^{s_k+1} = 0$. By the definition of $l_k$, the index $p$ is not used in iteration $\nu_k + l$, for $l_k < l \leq C$, so $z_p^{\nu_k+l}$ remains unaffected. We conclude that $z_p^{\nu_k+C+1} = 0$ for all $p \in I_1$, as we claimed in (6.189). This, together with (6.179) guarantees that

$$\lim_{k \to \infty} x^{\nu_k} = x^* \text{ and } z_p^{\nu_k} = 0, \text{ for } p \in I_1, \text{ and all } k \geq 0, \tag{6.195}$$

(compare with Proposition 6.3.1). From (6.195) we infer next that $x^*$ is a solution of (6.158)–(6.160). The proof resembles that of Step 6 of Theorem 6.3.1. Since $z_p^{\nu_k} = 0$ for $p \in I_1$, and $\langle a^p, x^* \rangle = b_p$ for $p \in I_2$, we have that

$$\begin{aligned} \langle z^{\nu_k}, Ax^{\nu_k} - b \rangle &= \langle A^{\mathrm{T}} z^{\nu_k}, x^{\nu_k} - x^* \rangle \\ &= -\langle \nabla f(x^{\nu_k}), x^{\nu_k} - x^* \rangle \\ &= -D_f(x^*, x^{\nu_k}) + f(x^*) - f(x^{\nu_k}). \end{aligned} \tag{6.196}$$

By Definition 6.8.2 the right-hand side of the last equation tends to zero, as $k \to \infty$, and thus, by continuity of $f$ and by Proposition 6.8.4 we have that, for any $x \in Q \cap \widetilde{S}$,

$$f(x) \geq \lim_{k\to\infty} L(x^{\nu_k}, z^{\nu_k}) = \lim_{k\to\infty} (f(x^{\nu_k}) + \langle z^{\nu_k}, Ax^{\nu_k} - b\rangle) = f(x^*), \tag{6.197}$$

showing that $x^*$ is indeed a minimizer of problem (6.158)–(6.160).

Since by strict convexity of $f$ problem (6.158)–(6.160) has at most one solution, it follows that $\{x^\nu\}$ has a unique accumulation point $x^*$, i.e., a limit point $x^* \in \widetilde{S}$. ∎

Algorithm 6.8.1 can be modified to become more suitable for application to the case of linear equations instead of linear inequality constraints, i.e., to the problem

$$\text{Minimize} \quad f(x) \tag{6.198}$$

$$\text{s.t.} \quad Ax = b, \tag{6.199}$$

$$x \in \widetilde{S}. \tag{6.200}$$

The modified algorithm eliminates dual iterates $z^\nu$ and has, for all $\nu \geq 0$, $c_\nu = \beta_\nu$. It generalizes Bregman's original algorithm for equality constraints by allowing underrelaxation parameters and by encompassing extended functions $f \in \mathcal{EB}(S)$. Theorem 6.8.1 also holds in this case. Note that the sequence $\{x^\nu\}$ produced by the modified algorithm will, in general, be different from the sequence obtained by applying Algorithm 6.8.1 after conversion of each equation into two inequalities, unless unity relaxation ($\alpha_\nu = 1$, for all $\nu \geq 0$) is used.

**Algorithm 6.8.2 Underrelaxed Bregman's Algorithm for Extended Bregman Functions and Linear Equality Constraints**

**Step 0:** (Initialization.) $x^0 \in Z_0$ is arbitrary, where $Z_0$ is defined as in (6.16), and $z^0$ is such that $\nabla f(x^0) = -A^T z^0$.

**Step 1:** (Iterative step.) Given $\nu$ and $x^\nu$, take an $\alpha_\nu$ such that (6.162) holds, pick a control index and solve for the next iterate $x^{\nu+1}$ the system:

$$\nabla f(x^{\nu+1}) = \nabla f(x^\nu) + \beta_\nu a^{i(\nu)}, \tag{6.201}$$

$$\langle x^{\nu+1}, a^{i(\nu)}\rangle = \alpha^\nu b_{i(\nu)} + (1 - \alpha_\nu)\langle a^{i(\nu)}, x^\nu\rangle. \tag{6.202}$$

**Control:** The sequence $\{i(\nu)\}_{\nu=0}^\infty$ is almost cyclic on $I$.

## 6.9 The Hybrid Algorithm: A Computational Simplification

From a computational point of view, the difficult part in Algorithms 6.3.1, 6.4.1, 6.6.1, 6.8.1, and 6.8.2 lies in the projection operation at each iterative step. The system of equations (6.155)–(6.156) consists of $n+1$ equations, $n$ of which are nonlinear in the $n+1$ variables $(x, s)$. If this system has to be solved numerically in each iteration, then the computational burden might reduce the efficiency of the algorithm. Numerical errors in the calculation of the projections may also cause the practical algorithm to deviate from the conceptual one.

As we have seen in Section 6.5 the special case when $f(x) = \frac{1}{2} \parallel x \parallel^2$, which leads to orthogonal projections onto hyperplanes, always leads to a closed-form formula that replaces the system (6.155)–(6.156). The main question that we address here is whether it is possible in the general case to have such algorithms with a closed-form formula in the inner loop. To motivate the subsequent development of hybrid algorithms (the reason for this name will be explained below) let us look at the case when $f(x) = -\operatorname{ent} x$, the negative of Shannon's entropy (see Example 2.1.2).

The iterative step of Algorithm 6.6.1 calls for the solution of (6.112)–(6.113) and finding $x^{\nu+1}$ as the unrelaxed entropy projection of $x^\nu$ onto the $i(\nu)$th hyperplane $H_{i(\nu)}$, i.e.,

$$x_j^{\nu+1} = x_j^\nu \exp(c_\nu a_j^{i(\nu)}), \quad j = 1, 2, \ldots, n. \tag{6.203}$$

For the same $i(\nu)$ and the same $x^\nu$, the iterative step of MART (Algorithm 6.6.2) without relaxation, i.e., with $\lambda_\nu = 1$, has the form

$$x_j^{\nu+1} = x_j^\nu \left( \frac{b_{i(\nu)}}{\langle a^{i(\nu)}, x^\nu \rangle} \right)^{a_j^{i(\nu)}}, \quad j = 1, 2, \ldots, n. \tag{6.204}$$

Is there any connection between the iterative steps (6.203) and (6.204)? The answer is more difficult to detect than to prove once detected. In order to perform (6.203) $c_\nu$ has to be calculated by solving (6.203) together with the additional equation

$$\langle x^{\nu+1}, a^{i(\nu)} \rangle = b_{i(\nu)}. \tag{6.205}$$

Substitution leads to the following single nonlinear equation for the single unknown $c_\nu$,

$$\sum_{j=1}^{n} a_j^{i(\nu)} x_j^\nu \exp(c_\nu a_j^{i(\nu)}) - b_{i(\nu)} = 0. \tag{6.206}$$

Letting

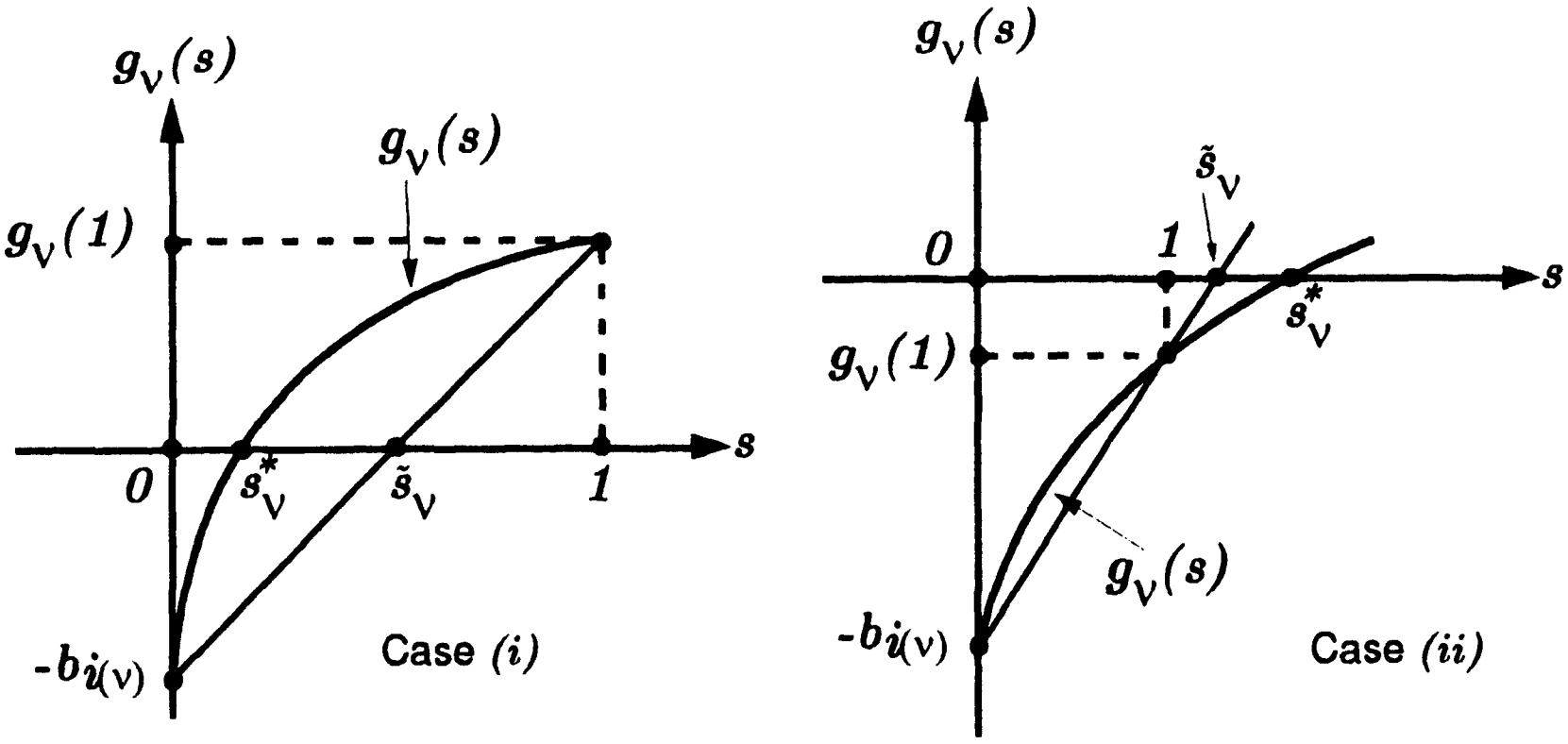

**Figure 6.6** The exact root $s^*_\nu$ and the secant root $\tilde{s}_\nu$ in case (*i*) when $g_\nu(1) > 0$, or in case (*ii*) when $g_\nu(1) < 0$.

$$\exp c_\nu \doteq s_\nu, \tag{6.207}$$

and using $s$ for $s_\nu$ we define the function

$$g_\nu(s) \doteq \sum_{j=1}^{n} a_j^{i(\nu)} x_j^\nu s^{a_j^{i(\nu)}} - b_{i(\nu)}, \tag{6.208}$$

and, of course, a positive root of this function is all that is needed to determine $c_\nu$. Now, $g_\nu(0) = -b_{i(\nu)} < 0$ (we assume that Assumptions 6.6.1–6.6.3 hold) and $\lim_{s\to\infty} g_\nu(s) = +\infty$. The derivates are $g'_\nu(s) > 0$ and $g''_\nu(s) \le 0$ for all $s \ge 0$, because of Assumption 6.6.2 and since $x^\nu > 0$ for all $\nu \ge 0$. Thus, there exists a unique $s^* > 0$ for which $g_\nu(s^*) = 0$, and $g_\nu$ is monotonically increasing and concave for $s \ge 0$ . Algorithm 6.6.1 requires the calculation of this root at each iterative step.

To see how MART fits into the picture we consider the line through the points $(0, -b_i)$ and $(1, g_\nu(1))$ in the plane of the graph of $g_\nu(s)$. This secant line to the graph of $g_\nu(s)$ intersects the $s$-axis at the point $\tilde{s}_\nu$ which, as a simple calculation confirms, is given by

$$\tilde{s}_\nu = \frac{b_{i(\nu)}}{\langle a^{i(\nu)}, x^\nu \rangle}, \tag{6.209}$$

in either of the two possible cases, $g_\nu(1) > 0$ or $g_\nu(1) < 0$.

Figure 6.6 describes the situation for both cases. Using $\tilde{s}_\nu$ for $s_\nu$ in (6.207) and substituting the resulting $\tilde{c}_\nu$ (i.e., $\tilde{c}_\nu = \log \tilde{s}_\nu$) into (6.203) for $c_\nu$ yields exactly the iterative step (6.204). Using $s^*_\nu$ in a similar way yields, from (6.203), an iterative step of an entropic Bregman projection.

This analysis shows that MART, which has a convergence proof of its own, is actually a kind of approximation to Bregman's algorithm (Algo-rithm 6.6.1). These considerations motivate us to replace $c_\nu$ in Bregman's algorithm for entropy maximization over linear inequalities by a MART-step. The algorithm we obtain is a computationally simplified algorithm, which we call a *hybrid algorithm* because it retains the overall structure of Bregman's algorithm (Algorithm 6.6.1) but replaces the computational burden involved in an inner loop calculation of $c_\nu$ by a closed-form formula. The use of closed-form formulae also facilitates the parallel implementation of such algorithms in situations when multiple $c_\nu$'s can be computed simultaneously for different constraints. Multiple processors can then evaluate in parallel the closed-form expressions.

We now present this idea in a more generalized way so that it will be applicable to cases other than Shannon's entropy. Equation (6.155) can sometimes be solved explicitly for $x$ in terms of $y, a$, and $s$. Since this is the case for all examples in which we are interested, we restrict our discussion to this situation. Therefore, we make the following assumption (which is indirectly an assumption on the Bregman function $f$).

**Assumption 6.9.1** *There exists a function $h : \mathbb{R}^{2n+1} \rightarrow \mathbb{R}^n$, called a solution function, such that $x = h(y, a, s)$ solves equation (6.155).*

Substitution of $x$ from (6.155) into (6.156) leads to a single nonlinear equation in $s$. Defining $g : \mathbb{R} \rightarrow \mathbb{R}$ by $g(s) \doteq \langle a, h(y, a, s)\rangle$ reduces the task of solving (6.155)–(6.156) to finding a root $s^*$ of the equation

$$g(s) = b. \tag{6.210}$$

Such a root provides $s^* = \pi_H(y)$ where $H$ is the hyperplane (6.153). As noted above, there may be several roots $s^*$ but only one will be such that

$$h(y, a, s^*) \in \widetilde{S}, \tag{6.211}$$

so that $x^* = h(y, a, s^*)$ will obey (6.155)–(6.157). Note also that

$$g(0) = \langle a, y\rangle. \tag{6.212}$$

Solving the single nonlinear equation $g(s) = b$ at each iterative step of Algorithm 6.8.1 instead of the system (6.155)–(6.157) might be difficult. However, the flexibility obtained through the introduction of the underrelaxation parameters $\{\alpha_\nu\}$ into Bregman's method (Algorithm 6.8.1) sometimes cancels out the need to solve iteratively (6.210). This may be achieved by choosing a particular underrelaxation strategy that enables a closed-form formula for the calculation of $\beta_\nu$ in each iterative step. If in addition to the previous assumptions Assumption 6.9.1 also holds, then the calculation of $\beta_\nu$ in equation (6.170) amounts to finding a root of

$$g_\nu(s) = \alpha_\nu b_{i(\nu)} + (1 - \alpha_\nu)\langle a^{i(\nu)}, x^\nu\rangle, \tag{6.213}$$

where

$$g_\nu(s) \doteq \langle a^{i(\nu)}, h(x^\nu, a^{i(\nu)}, s)\rangle. \tag{6.214}$$

The restriction (6.162) on $\alpha_\nu$ means that we are essentially seeking an $s$ such that $g_\nu(s)$ is in the interval between the values $\langle a^{i(\nu)}, x^\nu\rangle$ and $b_{i(\nu)}$. If $g_\nu(s)$ is continuous and monotone near zero, as it is in the examples below, then in view of (6.212) we seek some $s$ between zero and the solution $s^*$ of $g_\nu(s) = b_{i(\nu)}$. Formally, we need an $s$ such that

$$\epsilon \le \alpha_\nu = \frac{g_\nu(s) - \langle a^{i(\nu)}, x^\nu\rangle}{b_{i(\nu)} - \langle a^{i(\nu)}, x^\nu\rangle} \le 1, \tag{6.215}$$

for some arbitrary, fixed for all $\nu \ge 0$, $\epsilon > 0$, when

$$\langle a^{i(\nu)}, x^\nu\rangle \ne b_{i(\nu)}. \tag{6.216}$$

If equality holds in (6.216) then $\beta_\nu = 0$ and $x^{\nu+1} = x^\nu$. Therefore, if we are able to present an explicit formula for the calculation of $s$, i.e.,

$$s = \varphi(a, y, b), \tag{6.217}$$

such that for

$$s_\nu \doteq \varphi(a^{i(\nu)}, x^\nu, b_{i(\nu)}), \tag{6.218}$$

condition (6.215) is satisfied, then we have a closed-form formula to replace the calculations involved in projections and a new, computationally simpler algorithm emerges as follows.

### Algorithm 6.9.1 The Hybrid Algorithm

**Step 0:** (Initialization.) $x^0 \in Z_0$ is arbitrary, where $Z_0$ is as defined in (6.16), and $z^0$ is such that $\nabla f(x^0) = -A^T z^0$.

**Step 1:** (Iterative step.) Given $\nu$, $x^\nu$, and $z^\nu$, pick a control index $i(\nu)$ and take an $\alpha_\nu$ such that (6.215) holds with $s = s_\nu$ given explicitly by (6.218). Calculate $x^{\nu+1}$ and $z^{\nu+1}$ from

$$x^{\nu+1} = h(x^\nu, a^{i(\nu)}, d_\nu), \tag{6.219}$$

and

$$z^{\nu+1} = z^\nu - d_\nu e^{i(\nu)}, \tag{6.220}$$

where

$$d_\nu \doteq \min(z^\nu_{i(\nu)}, \mu_\nu), \tag{6.221}$$

and

$$\mu_\nu \doteq \varphi(a^{i(\nu)}, x^\nu, b_{i(\nu)}). \tag{6.222}$$

**Control:** The sequence $\{i(\nu)\}_{\nu=0}^\infty$ is almost cyclic on $I$.

The functions $h$ and $\varphi$ appearing above are those introduced in Assumption 6.9.1 and in (6.217), respectively. The term *hybrid* for Algorithm 6.9.1 is used because the algorithm has the same overall structure as Algorithm 6.8.1 but the projection coefficient $\beta_\nu$ is replaced by another coefficient $\mu_\nu$, which has a closed-form formula and whose nature is reminiscent of that of the MART algorithm.

The hybrid algorithm becomes a particular case of Algorithm 6.8.1, and Theorem 6.8.1 can be used to guarantee its convergence if the underrelaxation parameters are specifically chosen as

$$\alpha_\nu \doteq \frac{g_\nu(\varphi(a^{i(\nu)}, x^\nu, b_{i(\nu)})) - \langle a^{i(\nu)}, x^\nu \rangle}{b_{i(\nu)} - \langle a^{i(\nu)}, x^\nu \rangle}, \tag{6.223}$$

and condition (6.215) holds.

When applying Algorithm 6.9.1, some additional conditions on the matrix $A$ and the vector $b$ of (6.159) may be required to insure our ability to find a formula (6.217) such that (6.215) would hold. This is demonstrated in the following subsections where specific realizations of the hybrid algorithm are presented.

The extension to functions $f \in \mathcal{EB}(S)$ and the construction of the hybrid algorithm are motivated by entropy maximization problems. Such problems, which seek to maximize an entropy function over linear constraints (equality, inequality, or interval constraints), arise in various fields of applications. These include transportation planning, statistics, linear numerical analysis, chemistry, geometric programming, image reconstruction from projections, image restoration, pattern recognition, spectral analysis, and others. Numerous measures of entropy have been proposed in the different fields of applications. Following, we use as examples the previously discussed Shannon's entropy, $\operatorname{ent} x$, and Burg's entropy, and introduce another type of entropy called Rényi's entropy.

### 6.9.1 Hybrid algorithms for Shannon's entropy

Maximization of the entropy function $\operatorname{ent} x$ (defined in Example 2.1.2) gives rise to problem (6.158)–(6.160) or (6.198)–(6.200) with objective function $f(x) = \sum_{j=1}^n x_j \log x_j$. From Lemma 2.1.3 this is a Bregman function, thus $f \in \mathcal{EB}(S)$ with the zone $S = \operatorname{int} \mathbb{R}_+^n$, $\widetilde{S} = \overline{S}$ and $T = \emptyset$ , i.e., the singular set is empty. Theorem 6.8.1 applies and Algorithm 6.8.1 may be used with the iterative step for the primal iterates having the form

$$x_j^{\nu+1} = x_j^\nu \exp(c_\nu a_j^{i(\nu)}), \quad j = 1, 2, \ldots, n. \tag{6.224}$$

So, for this function, $h_j(y, a, s) \doteq y_j \exp(s a_j)$, $j = 1, 2, \ldots, n$, is the $j$th component of a solution function in the sense of Assumption 6.9.1 and $g(s) = \sum_{j=1}^n a_j y_j \exp(s a_j)$. If the additional conditions (compare with As-sumptions 6.6.2 and 6.6.3)

$$| a_j^i | \leq 1, \quad \text{and} \quad a_j^i b_i \geq 0, \quad \text{for } i = 1, 2, \ldots, m, \;\; j = 1, 2, \ldots, n, \qquad (6.225)$$

are imposed on the entries of the matrix $A$ and the vector $b$, then a closed-form formula of the form (6.217) is given by

$$\varphi(a, y, b) = \begin{cases} \lambda \log \left( \dfrac{b}{\langle a, y \rangle} \right), \text{ if } b > 0, \\ \\ \lambda \log \left( \dfrac{\langle a, y \rangle}{b} \right), \text{ if } b < 0, \end{cases} \qquad (6.226)$$

with $0 < \lambda < 1$. The case $b = 0$ can be disregarded because the nonnegativity of the primal iterates implies that if the right-hand side is zero then the variables corresponding to nonzero coefficients in such a row must be zero. Therefore, such rows can be deleted from the constraints beforehand.

In fact, the relaxation parameter $\lambda$ may vary from iteration to iteration as long as

$$0 < \eta \leq \lambda_\nu \leq 1, \qquad (6.227)$$

for an arbitrary but fixed $\eta$. The resulting hybrid algorithm follows.

**Algorithm 6.9.2 A Hybrid Algorithm for Shannon's Entropy**

**Step 0:** (Initialization.) $z^0 \in \mathbb{R}_+^m$ is arbitrary and $x^0 \in \mathbb{R}_+^n$ is such that

$$x_j^0 = \exp(-(A^{\mathrm{T}} z^0)_j - 1), \quad j = 1, 2, \ldots, n. \qquad (6.228)$$

**Step 1:** (Iterative step.) Given $\nu$, $x^\nu$, and $z^\nu$, pick a control index $i(\nu)$ and choose a relaxation parameter $\lambda_\nu$ that satisfies (6.227). Calculate $x^{\nu+1}$ and $z^{\nu+1}$ by

$$x_j^{\nu+1} = x_j^\nu \exp(\ell_\nu a_j^{i(\nu)}), \quad j = 1, 2, \ldots, n, \qquad (6.229)$$

$$z^{\nu+1} = z^\nu - \ell_\nu e^{i(\nu)}, \qquad (6.230)$$

where $e^{i(\nu)}$ is the $i(\nu)$th standard basis vector and

$$\ell_\nu \doteq \min(M_\nu, z_{i(\nu)}^\nu), \qquad (6.231)$$

$$M_\nu \doteq \begin{cases} \lambda_\nu \log \left( \dfrac{b_{i(\nu)}}{\langle a^{i(\nu)}, x^\nu \rangle} \right), \text{ if } b_{i(\nu)} > 0, \\ \\ \lambda_\nu \log \left( \dfrac{\langle a^{i(\nu)}, x^\nu \rangle}{b_{i(\nu)}} \right), \text{ if } b_{i(\nu)} < 0. \end{cases} \qquad (6.232)$$

**Control:** The sequence $\{i(\nu)\}_{\nu=0}^\infty$ is almost cyclic on $I$.

Adapting Algorithm 6.9.2 for linear equality constraints and using Assumptions 6.6.1–6.6.3 yield precisely MART as given in Algorithm 6.6.2. This MART was actually the forerunner of hybrid algorithms for entropy maximization that inspired the development of Algorithm 6.9.1.

### 6.9.2 Algorithms for the Burg entropy function

We revisit Burg's entropy, defined in (6.149). The function $h(x) \doteq -B(x)$ is not a Bregman function in the sense of Definition 2.1.1 because of the singularities occurring on the boundary of its zone int $\mathbb{R}^n_+$. It is however an extended Bregman function, i.e., $h \in \mathcal{EB}(S)$, and therefore Algorithms 6.8.1 and 6.9.1 can be applied. We describe now how this is done and present the resulting algorithms.

**Proposition 6.9.1** *$h \in \mathcal{EB}(S)$ with $S = \text{int } \mathbb{R}^n_+$, $\overline{S} = \mathbb{R}^n_+$, $\widetilde{S} = S$, and $T = \mathbb{R}^n_+ \backslash \text{int } \mathbb{R}^n_+$.*

**Proof** Conditions $(i)$, $(ii)$, and $(iii)$ of Definition 6.8.2 hold, thus we have

$$D_h(x, y) = \sum_{j=1}^{n} \log\left(\frac{y_j}{x_j}\right) + \sum_{j=1}^{n}\left(\frac{x_j}{y_j}\right) - n, \tag{6.233}$$

and therefore,

$$\lim_{\|x\| \to \infty} D_h(x, y) = +\infty, \qquad \lim_{\|y\| \to \infty} D_h(x, y) = +\infty, \tag{6.234}$$

for all $y$ and all $x$, respectively, in int $\mathbb{R}^n_+$. Also, for every $j$, $j = 1, 2, \ldots, n$,

$$\lim_{x_j \to 0^+} D_h(x, y) = +\infty, \quad \lim_{y_j \to 0^+} D_h(x, y) = +\infty, \tag{6.235}$$

for all $y$ and all $x$, respectively, in int $\mathbb{R}^n_+$. It follows that $L_1^h(y, \alpha)$ and $L_2^h(x, \alpha)$ satisfy $(iv)$ of Definition 6.8.2. Conditions $(v)$ and $(vi)$ follow from equation (6.233). ■

Now let us consider problem (6.158)–(6.160) with $f(x) = h(x)$. The feasible set of this problem is $F \doteq Q \cap \text{int } \mathbb{R}^n_+$, where $Q$ is given by (6.161). Denoting $N(A) \doteq \{x \in \mathbb{R}^n \mid Ax \leq 0\}$, we give first a necessary and sufficient condition for the solvability of (6.158)–(6.160) with $f(x) = h(x)$.

**Proposition 6.9.2** *If the set $F \neq \emptyset$ and the matrix $A \neq 0$, then a necessary and sufficient condition for the solution of (6.158)–(6.160) with $f(x) = h(x)$ to exist is*

$$N(A) \cap \mathbb{R}^n_+ = \{0\}. \tag{6.236}$$

**Proof** To prove sufficiency assume that $h(x)$ is unbounded from below on the set $F$, i.e., that $\inf \{h(x) \mid x \in F\} = -\infty$ and take a sequence $\{x^\nu\}_{\nu=0}^{\infty}$ such that $x^\nu \in F$ and $h(x^\nu) = -\nu$, for all $\nu \geq 0$. Such a sequence exists since $\{h(x) \mid x \in F\}$ is an interval, i.e., a connected set on the real line. This sequence must be unbounded because

$$\exp(-h(x)) = \prod_{j=1}^{n} x_j \leq \left( \max_{1 \leq j \leq n} x_j \right)^n, \tag{6.237}$$

and therefore, $\max_{1 \leq j \leq n} x_j^\nu \geq \exp(\nu/n)$, for all $\nu \geq 0$. Take a subsequence of $\{x^\nu / \parallel x^\nu \parallel\}$ which converges to a $v \in \mathbb{R}^n$. This limit must have the properties $v \geq 0$ and $\parallel v \parallel = 1$. However, since $Ax^\nu \leq b$, $Ax^\nu / \parallel x^\nu \parallel \leq b / \parallel x^\nu \parallel$, and the unboundedness of $\{x^\nu\}_{\nu=0}^{\infty}$ implies that $\parallel x^\nu \parallel \to \infty$, thus $Av \leq 0$. This is in contradiction to (6.236) because, since $\parallel v \parallel = 1$, it must be that $v \neq 0$.

To prove necessity suppose that $x^*$ solves (6.158)–(6.160) and that there exists some $v \geq 0$, $v \neq 0$, with $Av \leq 0$. Then, $x^* + tv \in F$ for all $t \geq 0$, and $\lim_{t \to \infty} h(x^* + tv) = -\infty$, which contradicts the assumption that $x^*$ is a solution. ■

We therefore assume henceforth that $a^i \neq 0$ for all $i = 1, 2, \ldots, m$, that $F \neq \emptyset$ (feasibility assumption), and that (6.236) holds. Next we show the validity of Assumptions 6.8.1 and 6.9.1 in this case.

**Proposition 6.9.3** *Assumptions 6.8.1 and 6.9.1 hold.*

**Proof** Regarding Assumption 6.8.1 for $(i)$ see Proposition 6.9.1; $(iii)$ holds by choice; and $(iv)$ is the feasibility assumption $F \neq \emptyset$ made in order to guarantee solvability according to Proposition 6.9.2. As observed immediately after Definition 6.8.4, when $\widetilde{S} = S$ (as is the case here; see Proposition 6.9.1), then $f$ is strongly zone consistent with respect to any hyperplane $H$ for which $H \cap S \neq \emptyset$. This settles $(ii)$ of Assumption 6.8.1. The solution function required in Assumption 6.9.1 is given by

$$h_j(y, a, s) = \frac{y_j}{1 - s y_j a_j}, \tag{6.238}$$

and thus all assumptions hold. ■

Given a hyperplane $H$ as in (6.153), and a $y \in \operatorname{int} \mathbb{R}^n_+$, the Burg's entropy projection of $y$ onto $H$ is the unique solution $x^*$ of (6.155)–(6.157). The associated projection coefficient $s^* = \pi_H(y)$ is calculated in the following way. Define

$$g(s) \doteq \sum_{j=1}^{n} \frac{a_j y_j}{1 - s y_j a_j}, \tag{6.239}$$

and observe that $g(s) = b$ has many solutions. In order to ensure that $P_H(y) \in S$, i.e., that $h(y, a, s) > 0$, $s^*$ has to be in the interval

$$r < s^* < t, \tag{6.240}$$

where

$$r \doteq \max\{1/a_j y_j \mid 1 \le j \le n, \quad a_j < 0\}, \tag{6.241}$$

and if $a_j \geq 0$ for all $j$, set $r = -\infty$ and

$$t \doteq \min\{1/a_j y_j \mid 1 \le j \le n, \quad a_j > 0\}, \tag{6.242}$$

and if $a_j \leq 0$ for all $j$, set $t = +\infty$.

**Proposition 6.9.4** *If $H$ is a hyperplane as in (6.153) with $a \neq 0$, and such that $H^+ \doteq H \cap \operatorname{int} \mathbb{R}^n_+ \neq \emptyset$, and if $y > 0$, then there exists a unique pair $(x^*, s^*)$ such that $g(s^*) = b$, equation (6.240) holds, and $x^* > 0$ is the Burg's entropy projection of $y$ onto $H$.*

**Proof** This follows directly from Lemma 6.8.1 and Proposition 6.9.1. But to better understand this issue we again analyze the specific solution at hand. The explicit form of $(x^*, s^*)$ can be obtained from (6.155)–(6.157) as the solution of the system

$$x_j^* = \frac{y_j}{1 - s^* y_j a_j}, \quad j = 1, 2, \ldots, n, \tag{6.243}$$

$$\langle x^*, a \rangle = b, \tag{6.244}$$

where $s^*$ obeys the additional restriction (6.240). These equations lead to the single equation

$$\sum_{j=1}^{n} a_j \frac{y_j}{1 - s^* y_j a_j} = b, \tag{6.245}$$

and therefore we consider the behavior of the function

$$g(s) \doteq \sum_{j=1}^{n} \frac{\alpha_j}{(1 - s\alpha_j)} = \sum_{j=1}^{n} \left( \frac{1}{\alpha_j} - s \right)^{-1}, \tag{6.246}$$

of the single variable $s$, $r < s < t$, where $\alpha_j \doteq a_j y_j$, for all $j = 1, 2, \ldots, n$. We observe that we may assume $\alpha_j \neq 0$ since $\alpha_j = 0$ contributes nothing to the first sum in (6.246). The derivative

$$g'(s) = \sum_{j=1}^{n} \left( \frac{1}{\alpha_j} - s \right)^{-2} \tag{6.247}$$

is always nonnegative, thus $g(s)$ is monotonically increasing. At the points $s_j = 1/\alpha_j$, $j = 1, 2, \ldots, n$, $g(s)$ is discontinuous in such a manner that

$$g(s) = b \tag{6.248}$$

has exactly one solution $s_j^*$ between each two consecutive $s_j$'s. To verify the existence of a solution to (6.248) we distinguish between three cases.

$-\infty < r < t < \infty$. In this case $\lim_{s \to r^+} g(s) = -\infty$ and $\lim_{s \to t^-} g(s) = +\infty$. Hence, by continuity of $g$ over the open interval $(r, t)$ there is a solution to (6.248).

$-\infty = r < t < \infty$. In this case $a_j \geq 0$, for all $j = 1, 2, \ldots, n$, $\lim_{s \to r^+} g(s) = 0$, and $\lim_{s \to t^-} g(s) = \infty$. By assumption, $H^+ \neq \emptyset$, so that $b$ in the definition of $H$ must be positive, hence (6.248) has a (finite) solution.

$-\infty < r < t = \infty$. This case is handled in the same way as the second case but with $b < 0$.

To conclude our analysis we verify that $x^*$ belongs to int $\mathbb{R}_+^n$. Indeed, if $a_j = 0$ then $x_j^* = y_j > 0$. If $a_j > 0$ then $1 - s^* y_j a_j > 0$, which shows that $x_j^* > 0$. The same argument works if $a_j < 0$. ■

The following algorithm is the underrelaxed Bregman algorithm for problem (6.158)–(6.160) with $f(x) = h(x)$.

**Algorithm 6.9.3 Underrelaxed Bregman Algorithm for Burg's Entropy Over Linear Inequalities**

**Step 0:** (Initialization.) $z^0 > 0$ and $x^0$ are such that $x_j^0 = 1/(A^T z^0)_j$, for all $j = 1, 2, \ldots, n$.

**Step 1:** (Iterative step.) Given $\nu$, $x^\nu$, and $z^\nu$, take an $\alpha_\nu$ such that (6.162) holds, pick a control index $i(\nu)$, and calculate $x^{\nu+1}$ and $z^{\nu+1}$ from:

$$x_j^{\nu+1} = \frac{x_j^\nu}{1 - c_\nu x_j^\nu a_j^{i(\nu)}}, \quad j = 1, 2, \ldots, n, \tag{6.249}$$

$$z^{\nu+1} = z^\nu - c_\nu e^{i(\nu)}, \tag{6.250}$$

where $c_\nu$ is given by (6.169) and $\beta_\nu$ is the unique solution $s$ of

$$\sum_{j=1}^n \frac{a_j^{i(\nu)} x_j^\nu}{1 - s x_j^\nu a_j^{i(\nu)}} = \alpha_\nu b_{i(\nu)} + (1 - \alpha_\nu) \langle a^{i(\nu)}, x^\nu \rangle, \tag{6.251}$$

for which $r_\nu < \beta_\nu < t_\nu$, where

$$r_\nu \doteq \max\{1/a_j^{i(\nu)} x_j^\nu \mid 1 \leq j \leq n, \quad a_j^{i(\nu)} < 0\}, \tag{6.252}$$

and if $a_j^{i(\nu)} \geq 0$ for all $j$, then $r_\nu = -\infty$, and

$$t_\nu \doteq \min\{1/a_j^{i(\nu)} x_j^\nu \mid 1 \leq j \leq n, \quad a_j^{i(\nu)} > 0\}, \tag{6.253}$$

and if $a_j^{i(\nu)} \leq 0$, for all $j$, then $t_\nu = +\infty$.

**Control:** The sequence $\{i(\nu)\}_{\nu=0}^{\infty}$ is almost cyclic on $I$.

We conclude our treatment of problems with Burg's entropy by discussing the hybrid algorithm for such entropy maximization. If we impose the additional condition

$$a_j^i b_i \geq 0, \quad i = 1, 2, \ldots, m, \quad j = 1, 2, \ldots, n, \tag{6.254}$$

then an appropriate $\varphi$ function for the hybrid algorithm turns out to be

$$\varphi(a, y, b) \doteq \begin{cases} \lambda(1 - \dfrac{\langle a, y \rangle}{b})t, & \text{if } b > 0, \\ \lambda(1 - \dfrac{\langle a, y \rangle}{b})r, & \text{if } b < 0, \end{cases} \tag{6.255}$$

where $t$ and $r$ are given by (6.241) and (6.242) and $0 < \lambda < 1$.

**Algorithm 6.9.4 Hybrid Algorithm for Burg's Entropy Maximization Over Inequalities**

**Step 0:** (Initialization.) $z^0 > 0$ and $x^0$ are such that

$$x_j^0 = 1/(A^{\mathrm{T}} z^0)_j, \quad j = 1, 2, \ldots, n. \tag{6.256}$$

**Step 1:** (Iterative step.) Given $\nu$, $x^\nu$, and $z^\nu$, take an $\alpha_\nu$ such that (6.162) holds, pick a control index $i(\nu)$, and calculate $x^{\nu+1}$ and $z^{\nu+1}$ from:

$$x_j^{\nu+1} = \frac{x_j^\nu}{1 - c_\nu x_j^\nu a_j^{i(\nu)}}, \quad j = 1, 2, \ldots, n, \tag{6.257}$$

and

$$z^{\nu+1} = z^\nu - c_\nu e^{i(\nu)}, \tag{6.258}$$

where $c_\nu$ is given by (6.169) and $\beta_\nu$ is given by the closed-form formula

$$\beta_\nu \doteq \lambda_\nu(1 - \frac{\langle a^{i(\nu)}, x^\nu \rangle}{b_{i(\nu)}})\theta_\nu, \tag{6.259}$$

where

$$\theta_\nu \doteq \begin{cases} t_\nu, & \text{if } b_{i(\nu)} > 0, \\ r_\nu, & \text{if } b_{i(\nu)} < 0. \end{cases} \tag{6.260}$$

**Relaxation parameters:** The sequence $\{\lambda_\nu\}$ satisfies (6.227).

**Control:** The sequence $\{i(\nu)\}_{\nu=0}^{\infty}$ is almost cyclic on $I$.

In order to establish convergence of any sequence $\{x^\nu\}$ generated by Algorithm 6.9.4 to a point $x^*$ that is a solution of problem (6.158)–(6.160) with $f(x) = h(x)$, it is enough to verify that the function $\varphi(a, y, b)$ defined by (6.255) satisfies condition (6.215) with $s = s_\nu$ as given by (6.218).

**Proposition 6.9.5** *Under condition (6.254), there exists an $\epsilon > 0$ such that*

$$\epsilon \leq \alpha_\nu \leq 1, \text{ for all } \nu \geq 0, \tag{6.261}$$

*where $\alpha_\nu$ are given by (6.223) with $\varphi$ as in (6.255).*

**Proof** We consider in detail the case $b_{i(\nu)} > 0$. From (6.214), (6.223), (6.238), and (6.255) we obtain, following some arithmetic manipulations,

$$\alpha_\nu = \sum_{j=1}^{n} \frac{a_j^{i(\nu)} x_j^\nu}{\langle a^{i(\nu)}, x^\nu \rangle + \sigma_{j\nu}}, \tag{6.262}$$

where, for all $j = 1, 2, \ldots, n$ and all $\nu \geq 0$,

$$\sigma_{j\nu} \doteq -b_{i(\nu)}(1 - \frac{1}{a_j^{i(\nu)} x_j^\nu t_\nu \lambda_\nu}). \tag{6.263}$$

From $b_{i(\nu)} > 0$ and (6.254) we have that

$$a_j^{i(\nu)} \geq 0, \quad j = 1, 2, \ldots, n. \tag{6.264}$$

From (6.227) and (6.253) we get

$$\frac{1}{a_j^{i(\nu)} x_j^\nu t_\nu \lambda_\nu} \geq 1, \tag{6.265}$$

thus $\sigma_{j\nu} \geq 0$ for all $j = 1, 2, \ldots, n$ and all $\nu \geq 0$, which proves, in view of (6.262), that $\alpha_\nu \leq 1$ for all $\nu \geq 0$. Also, since we assumed that $a^i \neq 0$ for all $i = 1, 2, \ldots, m$, we get $0 < \alpha_\nu$ for all $\nu \geq 0$.

As noted before, Proposition 6.8.5 holds even if $0 < \alpha_\nu \leq 1$. Therefore, $\{x^\nu\}$ is bounded and bounded away from the singularity set $T$, meaning here that there exist $\underline{q}$, $\overline{q} \in \mathbb{R}$ such that

$$0 < \underline{q} \leq x_j^\nu \leq \overline{q}, \text{ for all } j = 1, 2, \ldots, n, \text{ and all } \nu \geq 0. \tag{6.266}$$

Define,

$$\underline{p} \doteq \min\{\, |a_j^i| \mid 1 \leq i \leq m,\ 1 \leq j \leq n,\ a_j^i \neq 0\}, \tag{6.267}$$

$$\overline{p} \doteq \max\{\, |a_j^i| \mid 1 \leq i \leq m,\ 1 \leq j \leq n\}, \tag{6.268}$$

and

$$\delta \doteq \max\{\mid b_i \mid \mid 1 \leq i \leq m\}. \tag{6.269}$$

Considering only summands of (6.262) for which $a_j^{i(\nu)} \neq 0$, and using (6.227) where $\eta$ is defined, and (6.253) we conclude that

$$\frac{1}{\mid a_j^{i(\nu)} x_j^\nu t_\nu \lambda_\nu \mid} \leq \frac{\overline{p}\,\overline{q}}{\eta \underline{p}\,\underline{q}}, \tag{6.270}$$

because

$$\frac{1}{t_\nu} = \max\{a_j^{i(\nu)} x_j^\nu \mid 1 \leq j \leq n,\ \ a_j^{i(\nu)} > 0\}. \tag{6.271}$$

From (6.263), (6.269), and (6.270) we get

$$\sigma_{j\nu} \leq \delta(1 + \frac{\overline{p}\,\overline{q}}{\eta \underline{p}\,\underline{q}}) \doteq \omega, \tag{6.272}$$

and also that $\underline{p}\,\underline{q} \leq \left|\langle a^{i(\nu)}, x^\nu \rangle\right| \leq n \overline{p}\,\overline{q}$, and therefore, from (6.262) and in view of (6.264), we obtain finally that, for all $\nu \geq 0$,

$$\alpha_\nu \geq \frac{\underline{p}\,\underline{q}}{n \overline{p}\,\overline{q} + \omega} \doteq \epsilon. \tag{6.273}$$

This choice of $\epsilon$ completes the proof. The case $b_{i(\nu)} < 0$ can be handled in a similar way by using $r_\nu$ instead of $t_\nu$. ■

### 6.9.3 Rényi's entropy function

Rényi's entropy function of order $\alpha$ is defined by

$$H_\alpha(x) \doteq (1-\alpha)^{-1} \log \sum_{j=1}^{n} x_j^\alpha, \tag{6.274}$$

for $0 < \alpha < 1$. We consider the optimization problems (6.158)–(6.160) or (6.198)–(6.200) with $f(x) = g(x) \doteq -H_\alpha(x)$, for a fixed $\alpha$. This function is not a Bregman function because of a singularity that appears at the origin, which is a boundary point of its zone. Rényi's entropy is thus another example of an extended Bregman function.

**Proposition 6.9.6** *We have $g \in \mathcal{EB}(S)$ with zone $S = \mathbb{R}_+^n$, $\overline{S} = \mathbb{R}_+^n$, $\widetilde{S} = \{x \in \mathbb{R}_+^n \mid x \neq 0\}$ and $T = \{0\}$.*

**Proof** Using an auxiliary function $a : \widetilde{S} \to \mathbb{R}$ defined by

$$a(x) \doteq (1-\alpha) \sum_{j=1}^{n} x_j^\alpha, \tag{6.275}$$

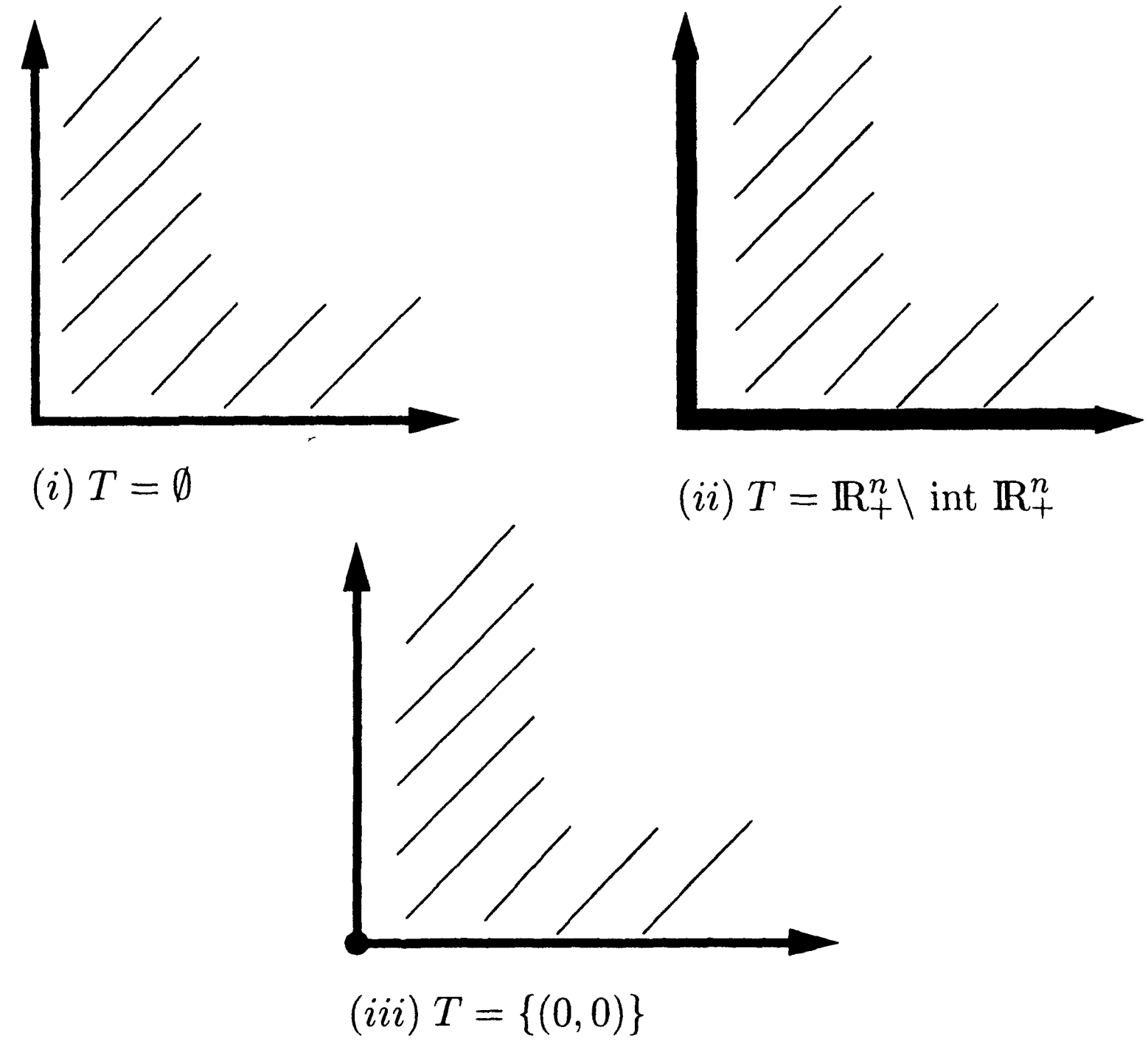

**Figure 6.7** The singularity set $T$ in three different cases: $(i)$ Shannon's entropy, $(ii)$ Burg's entropy, $(iii)$ Rényi's entropy.

the generalized distance takes the form

$$D_g(x,y) = \frac{\alpha \sum_{j=1}^{n} y_j^{\alpha-1} x_j}{a(y)} - (1-\alpha)^{-1} \log\left(\frac{a(x)}{a(y)}\right) + \frac{\alpha}{1-\alpha}. \tag{6.276}$$

The conditions of Definition 6.8.2 can be verified in a similar way to the proof of Proposition 6.9.1. ■

It is worthwhile to note that in the special case of Rényi's entropy the singularity at $x = 0$ may be circumvented by replacing in (6.158) or in (6.198) the objective function by

$$\psi(x) \doteq \sum_{j=1}^{n} x_j^{\alpha}. \tag{6.277}$$

This function is a Bregman function with zone $S = \mathbb{R}_+^n$; equivalently, $\psi(x) \in \mathcal{EB}(S)$, with $S = \widetilde{S} = \overline{S} = \mathbb{R}_+^n$, and $T = \emptyset$, thus Algorithms 6.8.1 and 6.9.1 are applicable. This alternative is preferable for any practical pur-pose of solving a linearly constrained entropy optimization problem with Rényi's entropy. The explicit formulae for implementing Algorithms 6.8.1 and 6.9.1 with Rényi's entropy or with the function $\psi(x)$, can be written down and all necessary verifications can be made. Calculations are somewhat cumbersome and the resulting formulae are quite complicated, so they are not included here.

Another example of an extended Bregman function with a nonempty singularity set is

$$q(x) \doteq \sum_{j=1}^{n}(x_j^{\alpha} + x_j^{-\beta}), \tag{6.278}$$

for $\alpha > 1$ and $\beta > 0$. In this case, $S = \text{int } \mathbb{R}^n_+$, $\overline{S} = \mathbb{R}^n_+$, $\widetilde{S} = S$ and $T = \mathbb{R}^n_+ \backslash \text{int } \mathbb{R}^n_+$.

In Figure 6.7 we describe graphically three singularity sets. In $(i)$ the set $T$ is empty and the zone $S = \text{int } \mathbb{R}^n_+$. This situation applies to Shannon's entropy as treated in Section 6.9.1. In $(ii)$ $T = \mathbb{R}^n_+ \backslash \text{int } \mathbb{R}^n_+$ and the zone is $S = \text{int } \mathbb{R}^n_+$, and this applies to Burg's entropy discussed in Proposition 6.9.1. The singularity set $T = \{(0,0)\}$ in $(iii)$ with zone $S = \mathbb{R}^n_+$ and $\widetilde{S} = \{x \in \mathbb{R}^n_+ \mid x \neq 0\}$ is related to Rényi's entropy function (see Proposition 6.9.6).

## 6.10 Notes and References

**6.1** The term *row-action* was introduced by Censor (1981). The regularization approach has its roots in Tikhonov's work; see, e.g., Tikhonov and Arsenin (1977), Bertero, De Mol, and Pike (1985), and references therein. An extensive and up-to-date study on least squares methods is Björck (1996). Application of row-action methods to (6.7)–(6.8) was suggested by Herman, Hurwitz, and Lent (1980). Kovarik (1977) studied other choices for $f$ and $g$ in (6.7)–(6.8) and Frieden (1972, 1975) used the entropy regularization approach where both $f$ and $g$ are of the form discussed in Example 2.1.2. The distance reduction approach is related to minimum distance estimation of statistical analysis; see, e.g., Titterington, Smith, and Makov (1985). The Kullback-Liebler distance reduction problem and the behavior of the EM algorithm when applied to it are discussed by Eggermont (1990), Iusem (1991b) and Snyder, Schultz, and O'Sullivan (1992).

**6.2–6.4** The collection of row-action methods presented here is not exhaustive. Various other well-known methods, such as Richardson's method (see, e.g., Young (1971)), also have row-action implementations. The Karush-Kuhn-Tucker conditions and the basics of duality theory can be found in a text on mathematical programming, e.g., Avriel (1976) or Luenberger (1984). Much of the material presented in this and subsequent sections has its roots in Bregman's pioneering work, Bregman (1967b). Censor and Lent (1981) seem to have done the first follow-up work on Bregman's study. They formally defined Bregman functions, studied Bregman projections, and proposed and established convergence of Algorithm 6.4.1 for interval constraints. De Pierro and Iusem (1986) introduced underrelaxation parameters into Algorithm 6.3.1 and recently Iusem and Zenios (1995) did the same for Algorithm 6.4.1. Iusem (1991a) continues the study of these algorithms by analyzing their dual convergence and the rate of primal convergence. Closely related to the material presented here, and in several important respects complementing it, are the works of Tseng (1990, 1991), Tseng and Bertsekas (1987, 1991), and Luo and Tseng (1992a, 1992b). See also Bertsekas and Tsitsiklis (1989, Chapter 5). For the special case $f(y) = \frac{1}{2} \parallel y \parallel^2$, Algorithm 6.4.1 coincides with Algorithm Scheme I of Herman and Lent (1978). The latter, which actually inspired the development of the algorithm for interval constraints (Algorithm 6.4.1) by Censor and Lent (1981), gives rise to the algorithms ART2 and ART4 (discussed in Section 6.5). Lamond and Stewart (1981) showed that most of the independently discovered matrix balancing methods, used in transportation planning and in other fields (see Chapter 9), are in fact special cases of Algorithms 6.3.1 and 6.4.1. However, they did not show how the latter algorithms are related to MART (Algorithm 6.6.2) and left this as an open problem.

**6.5** Kaczmarz (1937) published his algorithm for a square nonsingular constraint matrix and without relaxation parameters. The method was extensively studied since then and we only sample the vast literature on it. Herman, Lent, and Lutz (1978) introduced relaxation parameters and applied the method to the fully discretized model of the problem of image reconstruction from projections; see Chapter 10. Tanabe (1971, 1974) studied the behavior of Kaczmarz's method and other linear stationary iterative processes in the singular and general ($m \neq n$) cases; see also Dax (1990) and Tewarson (1972). Consult, e.g., Schott (1990), Maess (1988) and references therein for further related works by researchers of the German school, including Berg, Maess, Peters, Schott, and others. Censor, Eggermont, and Gordon (1983) analyzed the behavior of the method when strong underrelaxation is used and the system of equations is inconsistent. This was studied also by Elsner, Koltracht, and Lancaster (1991). The block-Kaczmarz algorithm was first published by Elfving (1980a) and by Eggermont, Herman and Lent (1981). Related to this is also the work of Aharoni, Duchet and Wajnryb (1984) showing that any sequence of points obtained by successive projections onto any finite family of hyperplanes in $\mathbb{R}^n$ is bounded.

Algorithm 6.5.2 appeared in Hildreth (1957) and in D'Esopo (1959); see also Luenberger (1969, Section 10.10). Its primal-dual nature and the almost cyclic control were studied in Lent and Censor (1980). A simultaneous version of Hildreth's algorithm was proposed and analyzed by Iusem and De Pierro (1987). Sugimoto, Fukushima, and Ibaraki (1995) investigated theoretically and experimentally a simultaneous ART4 algorithm. Dykstra's (1983) algorithm, further studied by Boyle and Dykstra (1986), Han (1988), Iusem and De Pierro (1991), Bauschke and Borwein (1994), and Deutsch and Hundal (1994) can be considered a generalization of Hildreth's method for convex constraints sets. See also the related work of He and Stoer (1992). Fukushima (1990) studies a conjugate gradient algorithm approach for norm minimization over linear inequalities.

**6.6** Entropy functions find various uses in many independent branches of science; see, e.g., Aczél and Daróczy (1975), Kapur and Kesavan (1992) and references mentioned therein, or Kapur (1983). The MART algorithm (Algorithm 6.6.2) was invented and reinvented in several fields. When the constraints in (6.108) are the so called transportation constraints, i.e., $A = \Phi$, where $\Phi$ is as in (9.33), or even for any 0–1 matrix $A$, i.e., all entries of $A$ are either zero or one, then Algorithm 6.6.2 coincides with Bregman's algorithm (Algorithm 6.6.1). This case of MART has been called the RAS algorithm in the matrix balancing literature (see Chapter 9). Bregman (1967a) studied this case and attributed the method to G.V. Sheleikhovskii.

Another early reference is Kruithof (1937); see also Krupp (1979). Lent (1977) proved convergence of MART for a general (not necessarily 0–1) matrix $A$ and introduced the underrelaxation parameters. Elfving (1980b) studied MART and the fully simultaneous version of block-MART and showed how both are obtained by applying the nonlinear SOR method; see, e.g., Ortega and Rheinboldt (1970), to the Karush-Kuhn-Tucker conditions of the problem (6.106)–(6.107) with $Q = Q_1$ as in (6.108). He has also rate of convergence results, applications, and numerical results. Further analysis of MART can be found in Lent and Censor (1991). An approximate version of MART was proposed in Censor, Lakshminarayanan, and Lent (1979).

**6.7** The block-MART algorithm (Algorithm 6.7.1) was studied by Censor and Segman (1987). It is similar to the method of *generalized iterative scaling* of Darroch and Ratcliff (1982). A detailed comparison between these two algorithms appears in Censor and Segman (1987). Byrne (1993) studied a fully simultaneous MART algorithm that he calls SMART and which closely resembles the algorithm obtained from Algorithm 6.7.1 when $M = 1$, $I_1 = I$, i.e., all constraint equations are lumped into a single block. For a proof of the arithmetic-geometric mean inequality (6.130) see, e.g., Peressini, Sullivan, and Uhl (1988).

**6.8–6.9** The extension of the family of Bregman functions to include functions with a nonempty singularity set is motivated by the special case of Burg's entropy. See, e.g., Johnson and Shore (1984) for a comparative discussion of Burg's entropy and Shannon's entropy. Censor and Lent (1987) studied the ad-hoc application of Bregman's algorithm (Algorithm 6.3.1) to linear inequalities constrained Burg's entropy maximization. See also the remarks in Section 2.6 on the recent extensions of the family of Bregman functions. The underrelaxed algorithm for Burg's entropy over inequality constraints (Algorithm 6.9.3) as well as the hybrid algorithm (Algorithm 6.9.4) come from Censor, De Pierro, and Iusem (1991). The general treatment of extended Bregman functions and the hybrid algorithm (Algorithm 6.9.1) presented here were developed by Censor, De Pierro and Iusem (1986). The relationship between MART (Algorithm 6.6.2) and Bregman's algorithm, which Lamond and Stewart (1981) left open, was studied in Censor et al. (1990). For a discussion of Rényi's and other entropy functions see, e.g., Aczél and Daróczy (1975) and Kapur and Kesavan (1992). Auslender (1992) identifies two fundamental properties of the dual function that allow a more general outlook on some of the primal-dual algorithms discussed here.

## CHAPTER 7

# Model Decomposition Algorithms

---

> Out of intense complexities intense simplicities emerge.
>
> Winston Churchill

In the previous chapters we studied algorithms that facilitate parallel computations due to the structure of their operations. Sometimes, however, it is the structure of the model of the problem at hand (that is suitable for some decomposition) that leads to efficient parallel computations. Large-scale problems often display some characteristic sparsity pattern which is amenable to decomposition. In time-dependent optimization problems, for example, one has to optimize a given system at different points in time. Successive time periods are linked through the flow of inventory. It might be possible to partially optimize the operations of the system for each time period, while maintaining a level of inventory consistent with the optimal operating schedule of successive time periods. In large spatial systems (e.g., transportation or telecommunication problems) one has to optimize distinct geographical regions. Adjacent regions are linked through trading and the flow of traffic. Totally decentralized optimization is not possible, but it might be possible to partially optimize each region separately, while restricting the trading between adjacent regions to be consistent with each region's optimal state.

Optimization algorithms have been devised over the last fifty years specifically to deal with such problems; the two most noteworthy examples are the Dantzig-Wolfe decomposition and Benders decomposition. An important book by Lasdon (1970) discusses decomposition algorithms for large-scale optimization problems. A common feature of these algorithms is that they solve a sequence of (smaller) *subproblems* to optimize the distinct components (time periods, regions, etc.), while a coordinating *master* problem synthesizes these solutions into an estimate of the overall optimum. The current solution estimate of the master program then defines a new subproblem and the process repeats iteratively. The subproblems and the master problem are much smaller than the original program. Hence, even if the decomposition algorithm requires several iterations between the master and the subproblems to reach a solution, it is usually faster than algorithms that attack the original undecomposed problem. We call such algorithms here *model decomposition algorithms*, as they do not solve the original model of the problem directly, but instead solve a modified decomposed variant of the model.

In the early days of parallel optimization, it was anticipated that parallelism would substantially speed up model decomposition algorithms. This has not been the case, however. The speedups observed with decomposition algorithms such as those of Dantzig-Wolfe or Benders have been modest. Why is this so? The most successful attempts to parallelize these algorithms (see references in Section 7.3) solve the subproblems in parallel, but solve the master program on a single processor. There have also been cases where the solution of the master program was parallelized, but the efficiency was low. Thus the master program becomes a serial bottleneck: as the decomposition algorithm iterates the master programs increase in size and the serial bottleneck becomes more restrictive, as Amdahl's law (Definition 1.4.6) dictates.

With a view toward parallelism other model decomposition algorithms have recently been designed and implemented, which either use a very simple coordination phase that does not create any serial bottleneck, or have a coordination phase that is itself suitable for parallel computations. In this chapter we will discuss one such algorithm. Section 7.1 contains preliminary discussion on model decompositions and discusses parallel decompositions based on linearization or diagonal-quadratic approximations. Section 7.2 discusses the Linear-Quadratic Penalty (LQP) algorithm for large-scale structured problems. Notes and references are given in Section 7.3.

## 7.1 General Framework of Model Decompositions

Consider the minimization of a convex, continuously differentiable, block-separable function $F : \mathbb{R}^{nK} \to \mathbb{R}$, written as $F(\boldsymbol{x}) \doteq \sum_{k=1}^{K} f_k(x^k)$ where $f_k : \mathbb{R}^n \to \mathbb{R}$ for all $k = 1, 2, \ldots, K$. The vector $\boldsymbol{x} \in \mathbb{R}^{nK}$ is the concatenation of $K$ subvectors $\boldsymbol{x} = \left((x^1)^T, (x^2)^T, \ldots, (x^K)^T\right)^T$ where $x^k \in \mathbb{R}^n$ for all $k = 1, 2, \ldots, K$. (Boldface letters denote vectors in the product space $\mathbb{R}^{nK}$.) Consider now the following constrained optimization problem:

**Problem** $[\mathcal{P}]$:

$$\text{Minimize} \quad F(\boldsymbol{x}) \tag{7.1}$$

$$\text{s.t.} \quad x^k \in X_k \subseteq \mathbb{R}^n, \qquad \text{for all } k = 1, 2, \ldots, K, \tag{7.2}$$

$$\boldsymbol{x} \in \Omega \subseteq \mathbb{R}^{nK}. \tag{7.3}$$

The sets $X_k, k = 1, 2, \ldots, K$, and $\Omega$ are assumed to be closed and convex. Figure 7.1 illustrates the structure of this problem in two dimensions.

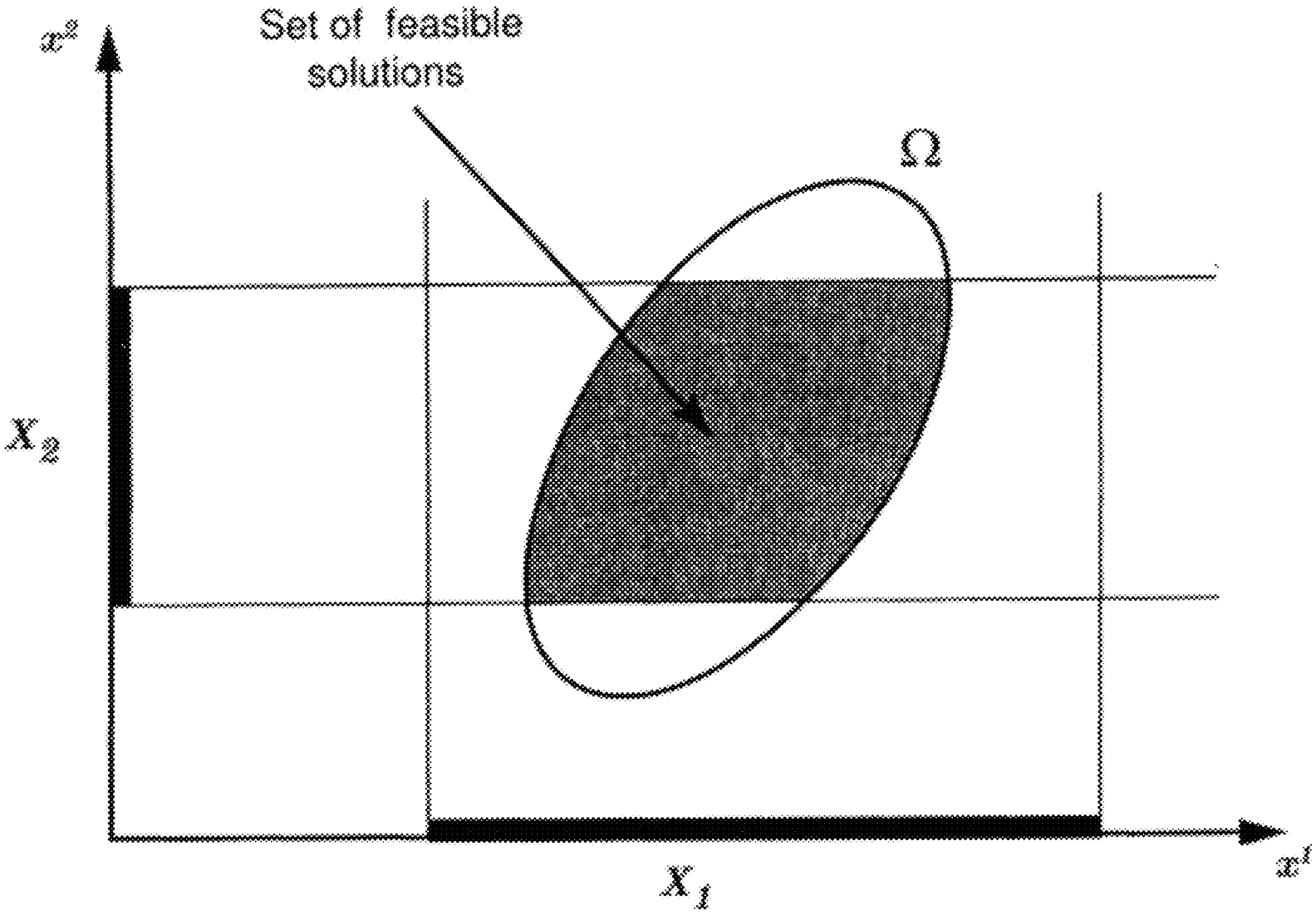

**Figure 7.1** Constraint sets and set of feasible solutions of problem $[\mathcal{P}]$ in $\mathbb{R}^2$ with $k = 2$.

If the product set $X_1 \times X_2 \times \cdots \times X_K \subseteq \Omega$, then problem $[\mathcal{P}]$ can be solved by simply ignoring the constraints $x \in \Omega$ and solving $K$ independent subproblems in each of the $x^k$ vector variables. In this respect, the constraints $x \in \Omega$ are *complicating* (or *coupling*) constraints. When the complicating constraints cannot be ignored, a model decomposition applies a *modifier* to problem $[\mathcal{P}]$ to obtain a problem $[\mathcal{P}']$ in which the complicating constraints are not explicitly present. It then employs a suitable algorithm to solve $[\mathcal{P}']$. If the solution to $[\mathcal{P}']$ is *sufficiently close* to a solution of the original problem $[\mathcal{P}]$ then the process terminates. Otherwise, the current solution is used to construct a new modified problem and the process repeats. Figure 1.2 illustrates the model decomposition algorithmic framework (see Chapter 1).

It is the judicious combination of a modifier and a suitable algorithm for the solution of the modified problem $[\mathcal{P}']$ that leads to a decomposition of the original problem $[\mathcal{P}]$ suitable for parallel computations. We discuss in this section modifiers and algorithms suitable for solving the modified problems.

### 7.1.1 Problem modifiers

We present now two modifiers for problem $[\mathcal{P}]$, drawing on the general theory developed in Chapter 4.

### Modifier I: Penalty or Barrier Functions

The first modifier eliminates the complicating constraints $\boldsymbol{x} \in \Omega$ by using a *penalty* or a *barrier* function. Using a penalty function $p : \mathbb{R}^{nK} \to \mathbb{R}$ with respect to the set $\Omega$ defining the complicating contraints (see Definition 4.1.1), the modified problem can be written as:

**Problem** $[\mathcal{P}']$:

$$\text{Minimize} \quad F(\boldsymbol{x}) + cp(\boldsymbol{x}) \tag{7.4}$$

$$\text{s.t.} \quad x^k \in X_k, \qquad \text{for all } k = 1, 2, \ldots, K. \tag{7.5}$$

We know (see Section 7.3) that it is possible to construct penalty functions that are *exact*, i.e., there exists some constant $\bar{c} > 0$, such that for $c > \bar{c}$ any solution of $[\mathcal{P}']$ is also a solution to $[\mathcal{P}]$. Hence, a solution to $[\mathcal{P}]$ can be obtained by solving $[\mathcal{P}']$. Note that $[\mathcal{P}']$ has a simpler constraint set than $[\mathcal{P}]$ because the complicating constraints $\boldsymbol{x} \in \Omega$ have been removed. However, problem $[\mathcal{P}']$ still cannot be solved by solving $K$ independent subproblems, since the function $p$ is not necessarily block-separable (see Definition 1.3.1 and Section 4.3 for definitions of separability). The next section explores algorithms that induce separability of this function.

Consider now situations when the set $\Omega$ has a nonempty interior. Such sets arise when inequality constraints are used to define them, e.g.,

$$\Omega = \{\boldsymbol{x} \mid g_l(\boldsymbol{x}) \leq 0, \text{ for all } l = 1, 2, \ldots, L\}, \tag{7.6}$$

where $g_l : \mathbb{R}^{nK} \to \mathbb{R}$. In this case we can use a *barrier function* (see Definition 4.2.1) to establish a barrier on the boundary of $\Omega$ so that the iterates of an algorithm that starts with an interior point remain in the interior of the set, therefore satisfying the constraints $\boldsymbol{x} \in \Omega$. For example, a barrier function for the set $\Omega$ defined by (7.6) can be constructed with the aid of Burg's entropy (6.149) (see Example 4.2.2) as:

$$q(\boldsymbol{x}) = -\sum_{l=1}^{L} \log g_l(\boldsymbol{x}),$$

where $g_l$ for $l = 1, 2, \ldots, L$ are the functions used in the definition of $\Omega$. With the use of such a barrier function the modified problem is written as:

**Problem** $[\mathcal{P}']$:

$$\text{Minimize} \quad F(\boldsymbol{x}) + cq(\boldsymbol{x}) \tag{7.7}$$

$$\text{s.t.} \quad x^k \in X_k, \qquad \text{for all } k = 1, 2, \ldots, K. \tag{7.8}$$

A solution to the problem $[\mathcal{P}]$ can be approximated by solving the modified problem with the barrier function for a sufficiently small value of the param-eter $c$. It is also possible to solve the barrier-modified problem repeatedly for a sequence of barrier parameters $\{c_\nu\}$, such that $c_\nu > c_{\nu+1} > 0$. If $\{\boldsymbol{x}^\nu\}$ denotes the sequence of solutions of these barrier problems, then it is known that $\{\boldsymbol{x}^\nu\}$ converges to a solution of $[\mathcal{P}]$ as $c_\nu \to 0$ (see Theorem 4.2.1).

Like the penalty-modified problem $[\mathcal{P}']$, the barrier-modified problem has a simpler constraint structure than problem $[\mathcal{P}]$. However, it still cannot be decomposed into independent components since the barrier function is not necessarily block-separable. The algorithms of the next section can be used to induce separability of penalty and barrier functions.

### Modifier II: Variable Splitting and Augmented Lagrangian

The second modifier first *replicates* (or *splits*) the components $x^k$ of the vector $\boldsymbol{x}$ into two copies, one of which is constrained to belong to the set $X_k$ and the other constrained to satisfy the complicating constraints. Let $\boldsymbol{z} \in \mathbb{R}^{nK}$ denote the replication of $\boldsymbol{x}$, where $\boldsymbol{z} = \left((z^1)^{\mathrm{T}}, (z^2)^{\mathrm{T}}, \ldots, (z^K)^{\mathrm{T}}\right)^{\mathrm{T}}$ and the vector $z^k \in \mathbb{R}^n$ for all $k = 1, 2, \ldots, K$. Consider now the equivalent split-variable formulation of $[\mathcal{P}]$:

**Problem** [Split-$\mathcal{P}$]:

$$\text{Minimize} \quad \sum_{k=1}^{K} f_k(x^k) \tag{7.9}$$

$$\text{s.t.} \quad x^k \in X_k, \qquad \text{for all } k = 1, 2, \ldots, K, \tag{7.10}$$

$$\boldsymbol{z} \in \Omega, \tag{7.11}$$

$$z^k = x^k, \qquad \text{for all } k = 1, 2, \ldots, K. \tag{7.12}$$

The constraints $z^k = x^k$ link the variables that appear in the constraint sets $X_k$ with the variables that appear in the set $\Omega$. An augmented Lagrangian formulation (see Section 4.4) is now used to eliminate these complicating constraints. We let $\boldsymbol{\pi} \doteq \left((\pi^1)^{\mathrm{T}}, (\pi^2)^{\mathrm{T}}, \ldots, (\pi^K)^{\mathrm{T}}\right)^{\mathrm{T}}$ where $\pi^k \in \mathbb{R}^n$ denotes the Lagrange multiplier vector for the complicating constraints $z^k = x^k$, and let $c > 0$ be a constant. Then a partial augmented Lagrangian for (7.9)–(7.12) can be written as:

$$\mathcal{L}_c(\boldsymbol{x}, \boldsymbol{z}, \boldsymbol{\pi}) = \sum_{k=1}^{K} f_k(x^k) + \sum_{k=1}^{K} \langle \pi^k, z^k - x^k \rangle + \frac{c}{2} \sum_{k=1}^{K} \| z^k - x^k \|^2 . \tag{7.13}$$

A solution to problem [Split-$\mathcal{P}$] can be obtained by solving the dual problem:

**Problem** $[\mathcal{P}'']$:

$$\underset{\boldsymbol{\pi} \in \mathbb{R}^{nK}}{\text{Maximize}} \ \varphi_c(\boldsymbol{\pi}), \tag{7.14}$$

where $\varphi_c(\boldsymbol{\pi}) \doteq \min_{x^k \in X_k, \, \boldsymbol{z} \in \Omega} \mathcal{L}_c(\boldsymbol{x}, \boldsymbol{z}, \boldsymbol{\pi})$. This is the modified problem whose solution yields a solution of $[\mathcal{P}]$.

An algorithm for solving convex optimization problems using augmented Lagrangians is the *method of multipliers*, which is an instance of the augmented Lagrangian algorithmic scheme (Algorithm 4.4.1). It proceeds by minimizing the augmented Lagrangian for a fixed value of the Lagrange multiplier vector, followed by a simple update of this vector. Using the method of multipliers to solve the dual problem $[\mathcal{P}'']$ we obtain the following algorithmic scheme:

**Algorithm 7.1.1 Method of Multipliers for Solving the Modified Problem $[\mathcal{P}'']$.**

**Step 0:** (Initialization.) Set $\nu = 0$. Let $\boldsymbol{\pi}^0$ be an arbitrary Lagrange multiplier vector.

**Step 1:** (Minimizing the augmented Lagrangian.)

$$\left(\boldsymbol{x}^{\nu+1}, \boldsymbol{z}^{\nu+1}\right) = \underset{x^k \in X_k,\ \boldsymbol{z} \in \Omega}{\operatorname{argmin}} \ \mathcal{L}_c\left(\boldsymbol{x}, \boldsymbol{z}, \boldsymbol{\pi}^{\nu}\right) \tag{7.15}$$

**Step 2:** (Updating the Lagrange multiplier vector.) For $k = 1, 2, \ldots, K$, update:

$$(\pi^k)^{\nu+1} = (\pi^k)^{\nu} + c\left((z^k)^{\nu+1} - (x^k)^{\nu+1}\right). \tag{7.16}$$

**Step 3:** Replace $\nu \leftarrow \nu + 1$ and return to Step 1.

The minimization problem in Step 1 has a block-decomposable constraint set. The problem, however, still cannot be decomposed into $K$ independent subproblems since the augmented Lagrangian is not block-separable due to the cross-products $\langle z^k, x^k \rangle$ in the quadratic term of (7.13). The next section explores algorithms that induce separability of this term. Step 2 consists of simple vector operations that can be executed very efficiently on parallel architectures.

### 7.1.2 Solution algorithms

Both modified problems $[\mathcal{P}']$ and $[\mathcal{P}'']$ have a block-decomposable constraint set, but the objective function is not block-separable. These problems can be written in the general form:

$$\text{Minimize} \quad \Phi(\boldsymbol{x}) \doteq \Phi(x^1, x^2, \ldots, x^K) \tag{7.17}$$

$$\text{s.t.} \quad \boldsymbol{x} \in X \doteq X_1 \times X_2 \times \cdots \times X_K. \tag{7.18}$$

We consider in this section two solution algorithms that—when applied to problems of this form—give rise to block-separable functions and thus decompose the problem into $K$ subproblems, which can then be solved in parallel. The first algorithm uses linear approximations to the nonlinear function $\Phi$. These linear approximations are block-separable. The second algorithm uses a diagonal approximation for cases where $\Phi$ is a nonseparable quadratic function, such as the last summand appearing in the augmented Lagrangian (7.13).

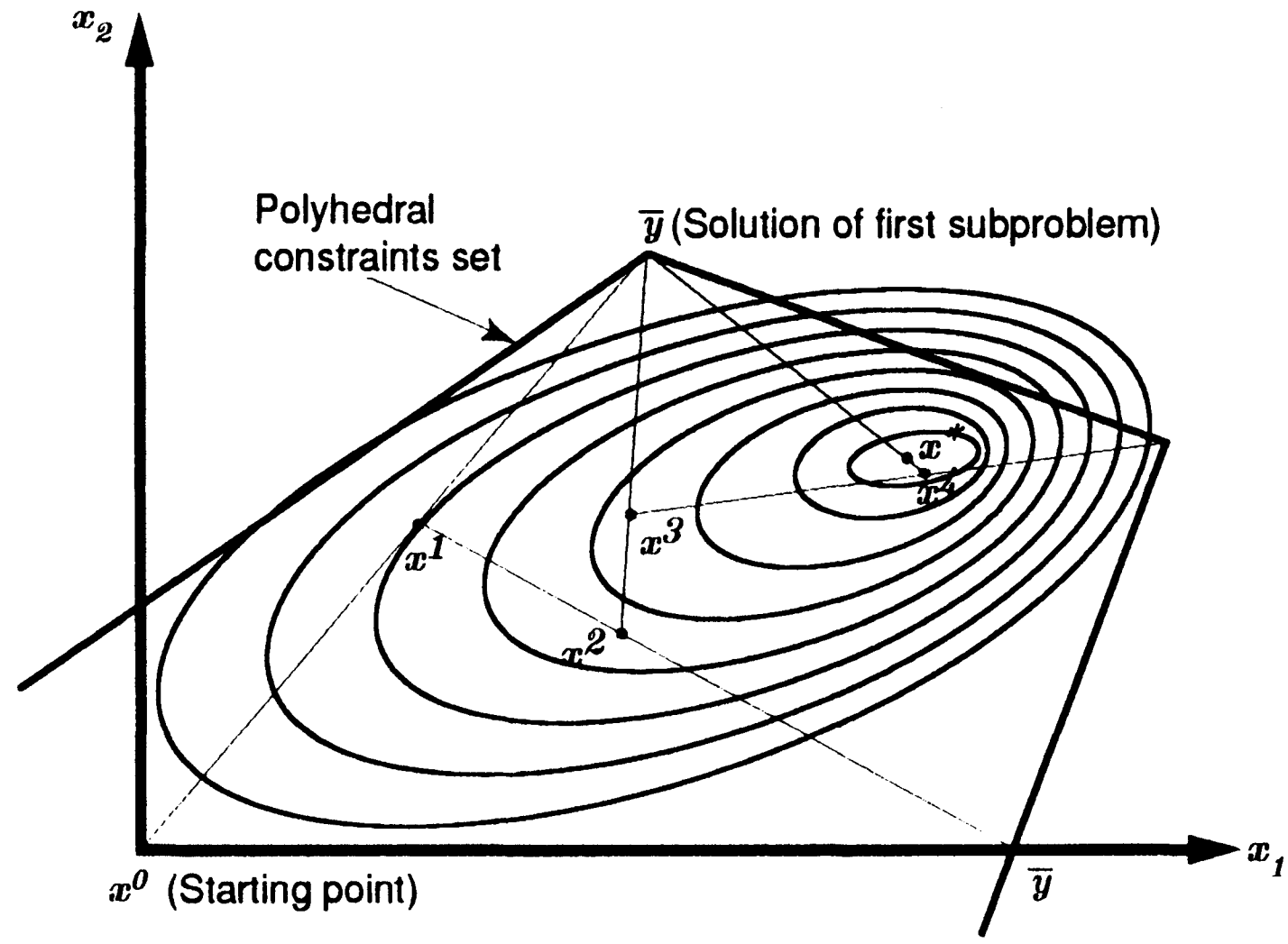

**Figure 7.2** Illustration of the Frank-Wolfe linearization algorithm for a function $\Phi(x_1, x_2)$. (The curves denote level sets.)

### Solution Algorithms Based on Linearization

One of the earliest algorithms suggested for the solution of nonlinear optimization problems using linear approximations is the Frank-Wolfe algorithm. It uses Taylor's expansion formula to obtain a first-order approximation $\hat{\Phi}$ to $\Phi$ around the current iterate $\boldsymbol{x}^\nu$, i.e.,

$$\hat{\Phi}(\boldsymbol{x}) = \Phi(\boldsymbol{x}^\nu) + \langle \nabla\Phi(\boldsymbol{x}^\nu), \boldsymbol{x} - \boldsymbol{x}^\nu \rangle,$$

ignoring second- and higher order terms. The Frank-Wolfe algorithm now minimizes this linear function, subject to the original constraints. The solution of this linear program, $\overline{\boldsymbol{y}}$, is a vertex of the constraint set, which determines a direction of descent for the original nonlinear function, given by $\boldsymbol{p} = \overline{\boldsymbol{y}} - \boldsymbol{x}^\nu$. The algorithm then performs a one-dimensional search along this direction to determine the step length where the nonlinear function attains its minimum. Figure 7.2 illustrates the algorithm.

Applied to problem (7.17)–(7.18) the Frank-Wolfe algorithm is formally stated below:

**Algorithm 7.1.2 The Frank-Wolfe Algorithm**

**Step 0:** (Initialization.) Set $\nu = 0$. Let $\boldsymbol{x}^0 \in X$ be an arbitrary vector.

**Step 1:** (Solving linearized subproblem.) Evaluate the first-order approximation of $\Phi(\boldsymbol{x})$ at $\boldsymbol{x}^\nu$, i.e., $\hat{\Phi}(\boldsymbol{x}) = \Phi(\boldsymbol{x}^\nu) + \langle \nabla\Phi(\boldsymbol{x}^\nu), \boldsymbol{x} - \boldsymbol{x}^\nu \rangle$, and compute a direction of descent of $\Phi$ by solving the linear program:

$$\underset{\boldsymbol{y} \in \mathbb{R}^{nK}}{\text{Minimize}} \quad \langle \nabla\Phi(\boldsymbol{x}^\nu), \boldsymbol{y} - \boldsymbol{x}^\nu \rangle \tag{7.19}$$

$$\text{s.t.} \qquad \boldsymbol{y} \in X. \tag{7.20}$$

Let $\overline{\boldsymbol{y}}$ denote the optimal vertex of the linear program. Then $\boldsymbol{p} = \overline{\boldsymbol{y}} - \boldsymbol{x}^\nu$ is a direction of descent for $\Phi$.

**Step 2:** (Linesearch.) Compute a step length $\alpha^*$ along the direction $\boldsymbol{p}$ that minimizes the nonlinear function $\Phi(\boldsymbol{x})$ by solving the one-dimensional nonlinear program:

$$\alpha^* = \underset{0 \le \alpha \le 1}{\text{argmin}}\, \Phi(\boldsymbol{x}^\nu + \alpha \boldsymbol{p}). \tag{7.21}$$

**Step 3:** (Updating the iterate.) Let $\boldsymbol{x}^{\nu+1} = \boldsymbol{x}^\nu + \alpha^* \boldsymbol{p}$ and return to Step 1.

The subproblem in Step 1 is a linear programming problem over a Cartesian product of linear constraint sets. It can be solved by solving $K$ independent linear programs in the variables $x^k$, $k = 1, 2, \ldots, K$, because although $\Phi(\boldsymbol{x})$ need not to be a block-separable function, the objective function in (7.19) is. Denoting the subvector of $\nabla\Phi(\boldsymbol{x})$, which is calculated with respect to $x^k$, by $\nabla_{x^k}\Phi(\boldsymbol{x})$ we obtain

$$\langle \nabla\Phi(\boldsymbol{x}^\nu), \boldsymbol{y} - \boldsymbol{x}^\nu \rangle = \sum_{k=1}^{K} \langle \nabla_{x^k}\Phi((x^k)^\nu), y^k - (x^k)^\nu \rangle. \tag{7.22}$$

The independent linear subproblems are then:

$$\underset{y^k \in \mathbb{R}^n}{\text{Minimize}} \quad \langle \nabla_{x^k}\Phi((x^k)^\nu), y^k - (x^k)^\nu \rangle \tag{7.23}$$

$$\text{s.t.} \qquad y^k \in X_k. \tag{7.24}$$

The problem in Step 2 is a nonlinear program in a single bounded variable. For large-scale problems the calculations of Step 2 are insignificant in comparison to the calculations of Step 1, and the algorithm parallelizes efficiently.

It is well-known that the Frank-Wolfe algorithm can *zigzag* toward the solution. If the optimal solution lies on a facet of the constraint set, the algorithm cannot reach it by a linesearch between one vertex (i.e., $\overline{\boldsymbol{y}}$) and an interior point (i.e., $\boldsymbol{x}^\nu$). A mild version of this effect is illustrated in Figure 7.2 where the optimal solution is close to a facet. The *simplicial decomposition* algorithm (described next) avoids the zigzagging effect of the Frank-Wolfe algorithm. It uses information contained in multiple vertices generated during successive iterations of the algorithm and solves a larger nonlinear master program rather than the simple linesearch. We describe next the simplicial decomposition algorithm.

Let $\mathcal{Y} = \{\overline{\boldsymbol{y}}^1, \overline{\boldsymbol{y}}^2, \ldots, \overline{\boldsymbol{y}}^\nu\}$ denote the set of vertices of the feasible region $X$ generated during the first $\nu$ iterations of the Frank-Wolfe algorithm. The convex hull of $\mathcal{Y}$ (i.e., the set defined by all convex combinations of the vertices) is

$$\text{conv}(\mathcal{Y}) = \{\boldsymbol{y} = \sum_{l=1}^{\nu} w_l \overline{\boldsymbol{y}}^l \mid \overline{\boldsymbol{y}}^l \in \mathcal{Y}, w_l \geq 0,\ l = 1, 2, \ldots, \nu,\ \sum_{l=1}^{\nu} w_l = 1\}, \tag{7.25}$$

which is a subset of the feasible set $X$. Simplicial decomposition generates vertices as in the Frank-Wolfe algorithm by solving linear programming subproblems. It then optimizes the original objective function by solving a *master* program over the set conv $(\mathcal{Y})$. The dimension of the master program is usually much smaller than the dimension of the original program. The master program also has a simple constraint structure—consisting of a single equality constraint and bounds on the variables— which can be exploited to convert the problem into a locally unconstrained program. Now the problem can be solved using standard unconstrained optimization techniques, and it can also be solved inexactly. Vertices that do not contribute to the representation of the optimal solution of the master program can be removed, thereby reducing the size of the master program.

Applied to problem (7.17)–(7.18), the simplicial decomposition algorithm is described next.

#### Algorithm 7.1.3 The Simplicial Decomposition Algorithm

**Step 0:** (Initialization.) Set $\nu = 0$. Let $\boldsymbol{x}^0 \in X$ be an arbitrary vector. Let $\mathcal{Y} = \emptyset$ be the set of vertices and its cardinality $v = 0$.

**Step 1:** (Solving linearized subproblem.) Evaluate the first-order approximation $\hat{\Phi}(\boldsymbol{x})$ at $\boldsymbol{x}^\nu$, i.e., $\hat{\Phi}(\boldsymbol{x}) = \Phi(\boldsymbol{x}^\nu) + \langle \nabla\Phi(\boldsymbol{x}^\nu), \boldsymbol{x} - \boldsymbol{x}^\nu \rangle$, and compute a direction of descent of $\Phi$ by solving the linear programming problem:

$$\underset{\boldsymbol{y} \in \mathbb{R}^{nK}}{\text{Minimize}} \quad \langle \nabla\Phi(\boldsymbol{x}^\nu), \boldsymbol{y} - \boldsymbol{x}^\nu \rangle \tag{7.26}$$

$$\text{s.t.} \quad \boldsymbol{y} \in X. \tag{7.27}$$

Let $\overline{\boldsymbol{y}}^\nu$ denote the optimal solution of this linear program. Update the set of vertices $\mathcal{Y} \leftarrow \mathcal{Y} \cup \{\overline{\boldsymbol{y}}^\nu\}$, and its cardinality $v \leftarrow v + 1$.

**Step 2:** (Solving the nonlinear master program.) Optimize the nonlinear function $\Phi$ over the convex hull of $\mathcal{Y}$ (equation (7.25)). That is, compute

$$w^* = \underset{w \in W_v}{\operatorname{argmin}} \, \Phi\left(\sum_{l=1}^{v} w_l \overline{\boldsymbol{y}}^l\right), \tag{7.28}$$

where $\overline{\boldsymbol{y}}^l \in \mathcal{Y}$, for all $l = 1, 2, \ldots, v$, and

$$W_v \doteq \{w = (w_l) \in \mathbb{R}^v \mid \sum_{l=1}^{v} w_l = 1, \; w_l \geq 0, \; \text{for all } l = 1, 2, \ldots, v\}.$$

**Step 3:** (Updating the iterate.) Let $\boldsymbol{x}^{\nu+1} = \sum_{l=1}^{v} w_l^* \overline{\boldsymbol{y}}^l$. Update the set of vertices $\mathcal{Y}$ by deleting from it any vertices with zero weight in the representation of $\boldsymbol{x}^{\nu+1}$, i.e., set $\mathcal{Y} \leftarrow \mathcal{Y} \backslash \{\overline{\boldsymbol{y}}^l \mid w_l^* = 0, 1 \leq l \leq v\}$ and let $v = \operatorname{card}(\mathcal{Y})$. Set $\nu \leftarrow \nu + 1$ and return to Step 1.

At Step 1 the algorithm solves a linear program with a Cartesian product of linear constraint sets. This subproblem can be decomposed, and its components solved independently and in parallel:

**Step 1 (alternate):** (Solving the decomposed linearized subproblems.) Let $\nabla_{x^k} \Phi((x^k)^\nu)$ denote the subvector of the gradient vector $\nabla \Phi(\boldsymbol{x})$ corresponding to the $k$th block of $\boldsymbol{x}$ evaluated at the current iterate $\boldsymbol{x}^\nu$. For each $k = 1, 2, \ldots, K$, solve

$$\underset{y^k \in \mathbb{R}^n}{\text{Minimize}} \quad \langle \nabla_{x^k} \Phi((x^k)^\nu), y^k - (x^k)^\nu \rangle \tag{7.29}$$

$$\text{s.t.} \qquad y^k \in X_k. \tag{7.30}$$

Let $(\overline{y}^k)^\nu$ denote the optimal solution of this linear program and form $\overline{\boldsymbol{y}}^\nu$ as the concatenation of the subvectors $\{(\overline{y}^l)^\nu \mid l = 1, 2, \ldots, K\}$.

The nonlinear master program in Step 2 is much smaller in size than the original problem. Typically, the number of vertices upon termination of the algorithm does not exceed one hundred. Furthermore, it has a simple structure, i.e., a simplex equality constraint, namely $\sum_{l=1}^{v} w_l = 1$, and bounds on the variables $0 \leq w_l \leq 1$. This structure can be exploited by designing an unconstrained optimization procedure for its solution as follows. Using the simplex equality constraint we can substitute $w_v$ with $w_v = 1 - \sum_{l=1}^{v-1} w_l$. Then we can write the master program (7.28) as:

$$w^* = \underset{\{0 \leq w_l \leq 1 \mid l = 1, 2, \ldots, v-1\}}{\operatorname{argmin}} \Phi(\overline{\boldsymbol{y}}^v + \sum_{l=1}^{v-1} w_l (\overline{\boldsymbol{y}}^l - \overline{\boldsymbol{y}}^v)). \tag{7.31}$$

Recall that at the current iteration we have $v-1$ *active* vertices (i.e., $w_l > 0$ for $l = 1, \ldots, v-1$) and the last vertex $\overline{\boldsymbol{y}}^v$ lies along a direction of descent. Hence the nonlinear master program is locally unconstrained in the neighborhood of the current iterate $\boldsymbol{x}^\nu$. Any unconstrained optimization algorithm can be used to compute a descent direction, followed by a simple test to determine the maximum allowable step that will keep the $w$'s within the bounds. Note that the evaluation of the objective function involves operations on dense vectors, i.e, the vectors $\overline{\boldsymbol{y}}^l$, $l = 1, 2, \ldots, v$. Such operations parallelize naturally and the solution of the master program is also amenable to parallel computations.

### Solution Algorithms Based on Diagonalization

Consider now a special structure of the objective function (7.17) arising from the modified problem $[\mathcal{P}'']$. In particular, we assume that $\Phi(\cdot)$ can be written as

$$\Phi(\boldsymbol{x}, \boldsymbol{z}) \doteq \sum_{k=1}^{K} \| x^k - z^k \|^2 . \tag{7.32}$$

This is the structure of the quadratic terms of the augmented Lagrangian (7.13). (We consider only the quadratic terms of the augmented Lagrangian, since these are the nonseparable terms that prevent us from decomposing the minimization of Step 1 of Algorithm 7.1.1 into $K$ independent subproblems.) We will approximate this nonseparable function using a separable quadratic function. The term *diagonal quadratic approximation* is also used, which indicates that the Hessian matrix of $\Phi(\boldsymbol{x}, \boldsymbol{z})$ is approximated by a diagonal matrix.

The terms in (7.32) can be expanded as

$$\| x^k - z^k \|^2 = \| x^k \|^2 + \| z^k \|^2 - 2\langle x^k, z^k \rangle,$$

and we only discuss the cross-product terms $\langle x^k, z^k \rangle$ for $k = 1, 2, \ldots, K$. Using Taylor's expansion formula we obtain a first-order approximation of the cross-product terms around the current iterate $(\boldsymbol{x}^\nu, \boldsymbol{z}^\nu)$ as:

$$\langle x^k, z^k \rangle \approx \langle x^k, (z^k)^\nu \rangle - \langle (x^k)^\nu, (z^k)^\nu \rangle + \langle (x^k)^\nu, z^k \rangle,$$

for $k = 1, 2, \ldots, K$. With this approximation of its cross-product terms, the function (7.32) is approximated by the expression:

$$\sum_{k=1}^{K} \| x^k - (z^k)^\nu \|^2 - \sum_{k=1}^{K} \| (x^k)^\nu - (z^k)^\nu \|^2 + \sum_{k=1}^{K} \| (x^k)^\nu - z^k \|^2 .$$

A solution algorithm, based on diagonalization solves problem $[\mathcal{P}'']$ using the method of multipliers (Algorithm 7.1.1) but instead of minimizing the augmented Lagrangian in Step 1 it minimizes the diagonal quadratic approximation. This approximation is block-separable into the variable blocks $x^k, k = 1, 2, \ldots, K$, and the minimization is decomposed into $K$ independent problems. These problems can be solved in parallel.

## 7.2 The Linear-Quadratic Penalty (LQP) Algorithm

Many different model decomposition algorithms can be designed using problem modifiers based on penalty or barrier methods and then followed by the use of linearization. We discuss here one such algorithm based on a linear-quadratic penalty (LQP) function. The algorithm has been shown to be efficient for large-scale, structured optimization problems. It has also been implemented on different parallel architectures.

We consider the problem:

$$\text{Minimize} \quad f_0(\boldsymbol{x}) \tag{7.33}$$

$$\text{s.t.} \quad f_k(x^k) \leq 0, \quad \text{for all } k = 1, 2, \ldots, K, \tag{7.34}$$

$$g_l(\boldsymbol{x}) \leq 0, \quad \text{for all } l = 1, 2, \ldots, L. \tag{7.35}$$

The functions $f_k : \mathbb{R}^n \rightarrow \mathbb{R}$, $k = 0, 1, \ldots, K$ and $g_l : \mathbb{R}^{nK} \rightarrow \mathbb{R}$, for $l = 1, 2, \ldots, L$ are convex and continuously differentiable. The constraints (7.34) decompose into blocks, one for each subvector $x^k$, while constraints (7.35) are complicating. Let $X = X_1 \times X_2 \times \cdots \times X_K$ be a product of the sets $X_k = \{x^k \in \mathbb{R}^n \mid f_k(x^k) \leq 0\}$ for all $k = 1, 2, \ldots, K$ and assume that $X$ is a compact set. We further make the following assumptions:

**Assumption 7.2.1** *Problem (7.33)-(7.35) has a nonempty and compact optimal solutions set.*

**Assumption 7.2.2** *Problem (7.33)-(7.35) has at least one feasible solution that satisfies all the constraints with strict inequality.*

Under these assumptions a Kuhn-Tucker vector (i.e., a Lagrange multiplier vector as defined in Rockafellar (1970, p. 274)) exists for problem (7.33)–(7.35).

Consider now the $\ell_1$-norm penalty function $p : \mathbb{R} \rightarrow \mathbb{R}$ given by:

$$p(t) \doteq \mu \max(0, t), \tag{7.36}$$

where $t$ is a scalar variable, and $\mu$ is a positive constant (see Figure 7.3). We want to obtain a solution to (7.33)–(7.35) by solving the following exact penalty problem:

$$\min_{\boldsymbol{x} \in X} \left( f_0(\boldsymbol{x}) + \sum_{l=1}^{L} p(g_l(\boldsymbol{x})) \right). \tag{7.37}$$

It is known (see, e.g., Bertsekas (1982, Chapter 4)) that under the assumptions stated above there exists a penalty parameter $\mu$ for which the optimal solutions to (7.37) and (7.33)–(7.35) coincide. In particular, if the penalty parameter is larger than a *threshold* value $\mu^*$ given by the largest component of a Lagrange multiplier vector of (7.33)–(7.35), then a solution to (7.37) is also a solution to (7.33)–(7.35). However the $\ell_1$-norm exact penalty function is nondifferentiable, and this precludes the application of gradient-based descent methods for the solution of the exact penalty problem. In order to gain access to gradient-based minimization techniques for solving (7.37) we consider an $\epsilon$-smoothing of the function $p$ around $t = 0$. In particular we introduce the following *$\epsilon$-smoothed function* $p_\epsilon$:

$$p_\epsilon(t) \doteq \begin{cases} 0, & \text{if } t < 0, \\ \mu \dfrac{t^2}{2\epsilon}, & \text{if } 0 \le t \le \epsilon, \\ \mu(t - \dfrac{\epsilon}{2}), & \text{if } t > \epsilon, \end{cases} \tag{7.38}$$

where $\epsilon$ is a positive scalar; see Figure 7.3. It is easy to see that

$$\lim_{\epsilon \longrightarrow 0^+} p_\epsilon(t) = p(t), \tag{7.39}$$

and, furthermore, this convergence is *uniform*, that is, given an arbitrary $\eta > 0$ we can find a $\delta > 0$ such that if $0 < \epsilon < \delta$ then $\mid p_\epsilon(t) - p(t) \mid < \eta$, for all $t$ (i.e., $\delta$ does not depend on $t$).

With the introduction of $\epsilon$-smoothing we obtain a continuously differentiable penalty function that can be optimized using algorithms such as the simplicial decomposition Algorithm 7.1.3. Furthermore, the $\epsilon$-smoothing device provides a natural mechanism for handling "soft" constraints that need not be satisfied exactly.

The use of $\epsilon$-smoothing might introduce some approximation error because an optimum point of the original problem is not necessarily an optimum of the smoothed penalty problem, even for penalty parameter values larger than the threshold value. However, an a priori upper bound to this error can be computed as a function of the penalty parameters. It is also possible to compute a solution that is feasible to within any given $\epsilon > 0$ for given values of the penalty parameter. Such a solution is termed *$\epsilon$-feasible.*

**Definition 7.2.1** **($\epsilon$-feasibility)** *Given some $\epsilon > 0$, a vector $\tilde{x} \in X$ is $\epsilon$-feasible for problem (7.33)–(7.35) if $g_l(\tilde{x}) \le \epsilon$, for all $l = 1, 2, \ldots, L$.*

We first describe the linear-quadratic penalty (LQP) algorithm that uses the smoothed penalty function $p_\epsilon(t)$ instead of the exact penalty $p(t)$ in (7.37). We then proceed with the analysis of the properties of $\epsilon$-smoothing, and derive bounds on the difference between the solution of the smoothed

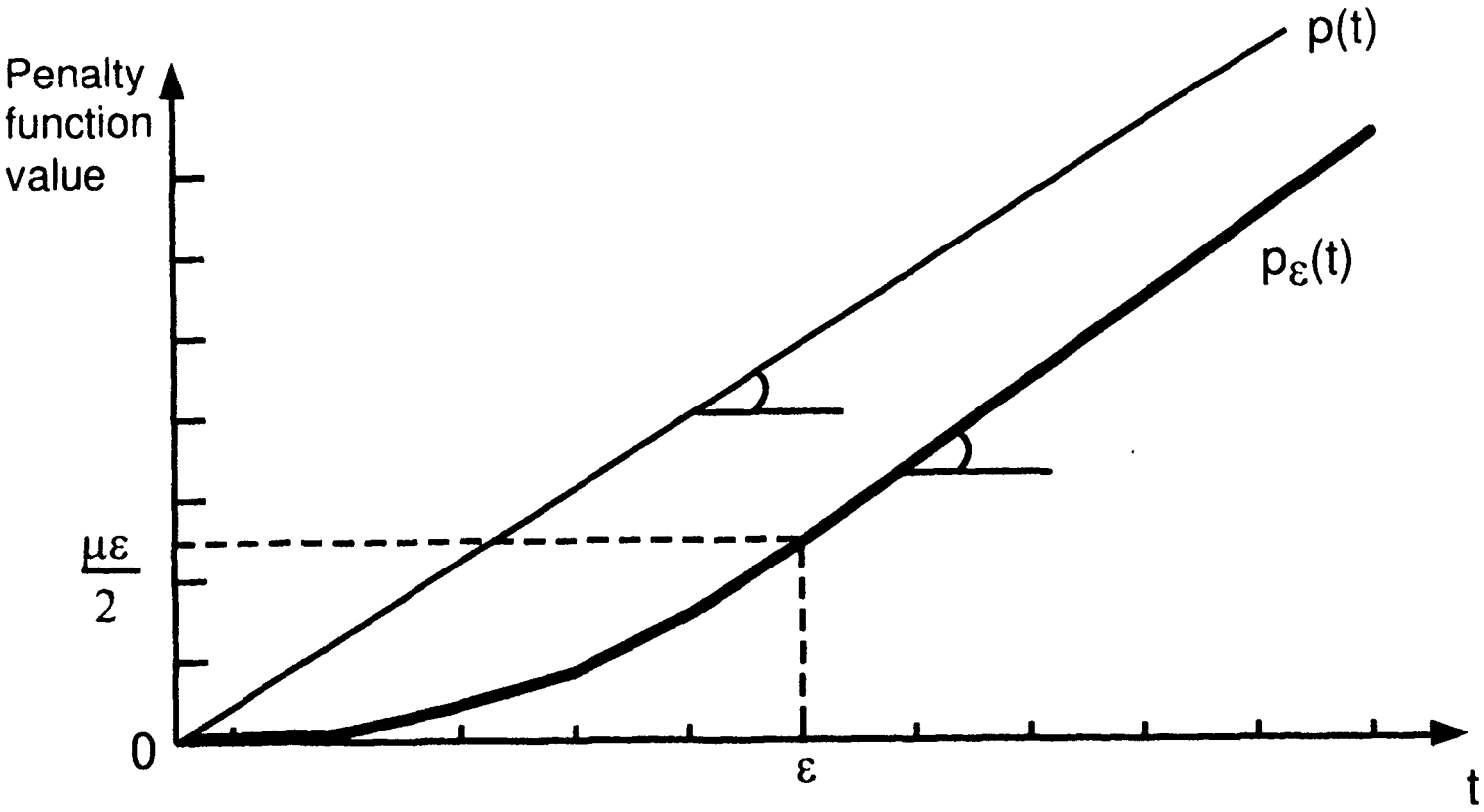

**Figure 7.3** The $\ell_1$-norm penalty function and the linear-quadratic smooth penalty function defined by (7.38).

penalty problem and the original problem (7.33)–(7.35).

Define first the objective function for the exact penalty problem

$$F(\boldsymbol{x}, \mu) \doteq f_0(\boldsymbol{x}) + \sum_{l=1}^{L} p(g_l(\boldsymbol{x})), \tag{7.40}$$

and then express the objective function for the $\epsilon$-smoothed penalty problem

$$\tilde{F}(\boldsymbol{x}, \mu, \epsilon) \doteq f_0(\boldsymbol{x}) + \sum_{l=1}^{L} p_\epsilon(g_l(\boldsymbol{x})). \tag{7.41}$$

The LQP algorithm can now be described in detail:

**Algorithm 7.2.1 The Linear-Quadratic Penalty (LQP) Algorithm**

**Step 0:** (Initialization.) Set $\nu = 0$, set the initial parameter values $\mu_0 > 0$, $\epsilon_0 > 0$, and define some $\epsilon_{min} > 0$.

**Step 1:** Solve the problem

$$\min_{\boldsymbol{x} \in X} \tilde{F}(\boldsymbol{x}, \mu_\nu, \epsilon_\nu).$$

and let $\tilde{\boldsymbol{x}}^\nu$ denote its optimal solution.

**Step 2:** If $\tilde{\boldsymbol{x}}^\nu$ is $\epsilon_\nu$-feasible and $\epsilon_\nu \leq \epsilon_{min}$, then stop. Otherwise, update the penalty parameters $\mu_\nu$ and $\epsilon_\nu$ according to the rules given below, set $\nu \leftarrow \nu + 1$ and go to Step 1.

The parameter $\epsilon_{min} > 0$ is a user-determined final feasibility tolerance. At Step 1 the algorithm solves the modified penalty problem, whereby the complicating constraints have been placed in the objective function, and at Step 2 the modifier is updated. Two situations may arise:

- The point $\tilde{x}^\nu$ is $\epsilon_\nu$-feasible but $\epsilon_\nu > \epsilon_{min}$. In this case the penalty parameter remains unchanged, i.e., $\mu_{\nu+1} = \mu_\nu$ and $\epsilon_\nu$ is decreased by replacing it with $\epsilon_{\nu+1} = \eta_1 \max_{1\le l\le L} \tilde{u}_l$ where $\tilde{u}_l \doteq g_l(\tilde{x}^\nu)$ for all $l = 1, 2, \ldots, L$, and $0 < \eta_1 < 1$ is a user-specified fixed parameter.
- The point $\tilde{x}^\nu$ is not $\epsilon_\nu$-feasible. This indicates that the value of the penalty parameter $\mu_\nu$ is not large enough and must be increased by replacing it with $\mu_{\nu+1} = \eta_2\mu_\nu$ where $\eta_2 > 1$ is a user-specified fixed parameter. The feasibility tolerance remains unchanged, i.e., $\epsilon_{\nu+1} = \epsilon_\nu$.

### 7.2.1 Analysis of the $\epsilon$-smoothed linear-quadratic penalty function

We now analyze the approximation error introduced with the use of $\epsilon$-smoothing. The following result gives an a priori upper bound on the difference between the exact and the smoothed penalty functions.

**Proposition 7.2.1** *Let the functions $F$ and $\tilde{F}$ be defined by (7.40) and (7.41), respectively. Then*

$$0 \le F(x, \mu) - \tilde{F}(x, \mu, \epsilon) \le L\mu\frac{\epsilon}{2}, \tag{7.42}$$

*for any $x \in \mathbb{R}^{nK}$, $\mu > 0$ and $\epsilon \ge 0$.*

**Proof** From the definitions of $p$ and $p_\epsilon$ we have

$$0 \le p(g_l(x)) - p_\epsilon(g_l(x)) \le \mu\frac{\epsilon}{2} \text{ for all } l = 1, 2, \ldots, L, \text{ and } x \in \mathbb{R}^{nK}. \tag{7.43}$$

Adding these inequalities up, for $l = 1, 2, \ldots, L$, we obtain

$$0 \le \sum_{l=1}^{L} p(g_l(x)) - \sum_{l=1}^{L} p_\epsilon(g_l(x)) \le L\mu\frac{\epsilon}{2}, \text{ for all } x \in \mathbb{R}^{nK}, \tag{7.44}$$

and the result follows from the definitions of $F$ and $\tilde{F}$. ∎

The LQP algorithm solves in Step 1 the smoothed penalty problem

$$\min_{x\in X} \tilde{F}(x, \mu, \epsilon), \tag{7.45}$$

instead of solving the exact penalty problem (7.37). The following results give a priori upper bounds on the error incurred by solving the smooth penalty problem (7.45) in lieu of the nondifferentiable penalty problem (7.37).

**Proposition 7.2.2** *Let $\boldsymbol{x}^* \in X$ be an optimal solution of (7.33)–(7.35) and $\tilde{\boldsymbol{x}} \in X$ be an optimal solution of (7.45) for some $\mu$ and $\epsilon$. Then*

$$0 \leq F(\boldsymbol{x}^*, \mu) - \tilde{F}(\tilde{\boldsymbol{x}}, \mu, \epsilon) \leq L\mu\frac{\epsilon}{2}. \tag{7.46}$$

**Proof** From Proposition 7.2.1 we have

$$F(\boldsymbol{x}, \mu) \leq \tilde{F}(\boldsymbol{x}, \mu, \epsilon) + L\mu\frac{\epsilon}{2}. \tag{7.47}$$

Taking the infimum, we obtain

$$\inf_{\boldsymbol{x}\in X} F(\boldsymbol{x}, \mu) \leq \inf_{\boldsymbol{x}\in X} \tilde{F}(\boldsymbol{x}, \mu, \epsilon) + L\mu\frac{\epsilon}{2}, \tag{7.48}$$

which proves the right-hand side inequality. The left-hand side inequality can be similarly proved. ■

The next proposition tells us that the difference between the optimal values of the exact penalty problem and the smoothed penalty problem can be controlled through the parameter $\epsilon$ provided that the solution of the penalty problem is $\epsilon$-feasible.

**Proposition 7.2.3** *Let $\boldsymbol{x}^* \in X$ be an optimal solution of (7.33)–(7.35) and $\tilde{\boldsymbol{x}} \in X$ be an optimal solution of (7.45) for some $\mu$ and $\epsilon$. Furthermore let $\tilde{\boldsymbol{x}}$ be $\epsilon$-feasible. Then*

$$0 \leq f_0(\boldsymbol{x}^*) - f_0(\tilde{\boldsymbol{x}}) \leq L\mu\epsilon. \tag{7.49}$$

**Proof** Since $\tilde{\boldsymbol{x}}$ is $\epsilon$-feasible, we have by the definition of $p_\epsilon$ that

$$\sum_{l=1}^{L} p_\epsilon(g_l(\tilde{\boldsymbol{x}})) \leq L\mu\frac{\epsilon}{2}. \tag{7.50}$$

Also, $\boldsymbol{x}^*$ is a solution to (7.33)–(7.35), which implies that

$$\sum_{l=1}^{L} p(g_l(\boldsymbol{x}^*)) = 0. \tag{7.51}$$

From Proposition 7.2.2 we have

$$0 \leq (f_0(\boldsymbol{x}^*) + \sum_{l=1}^{L} p(g_l(\boldsymbol{x}^*))) - (f_0(\tilde{\boldsymbol{x}}) + \sum_{l=1}^{L} p_\epsilon(g_l(\tilde{\boldsymbol{x}}))) \leq L\mu\frac{\epsilon}{2}. \tag{7.52}$$

Substituting (7.50) and (7.51) into (7.52) and rearranging terms, the result is established. ■

The next question is how to specify conditions on the penalty parameter $\mu$ to ensure that $\epsilon$-feasibility is achieved. For the exact penalty function case it is known (see, for example, Bertsekas (1975, Proposition 1)) that for $\mu$ larger than a threshold $\mu^*$ given by $\mu^* = \max_{l=1,2,\ldots,L} z_l^*$ where $z^* = (z_l^*)$ is a Lagrange multiplier vector for the constraints (7.35), an optimal solution $\boldsymbol{x}^*$ to (7.33)–(7.35) is also an optimal solution for (7.37). This result provides the motivation for exact penalty methods. The question we address next is how this result is affected by the $\epsilon$-smoothing of the exact penalty function. The answer is that although an optimal solution to (7.33)–(7.35) does not necessarily coincide with an optimal solution to (7.45), its suboptimality with respect to the latter can be bounded as a function of the penalty and smoothing parameters.

To simplify notation we introduce, for each constraint function $g_l(\boldsymbol{x})$, the function $g_l^+(\boldsymbol{x}) \doteq \max(0, g_l(\boldsymbol{x}))$, thus $p(g_l(\boldsymbol{x})) = \mu g_l^+(\boldsymbol{x})$. Obviously if $g_l(\boldsymbol{x}) > 0$ then $g_l^+(\boldsymbol{x}) = g_l(\boldsymbol{x})$.

**Proposition 7.2.4** *Let $\boldsymbol{x}^*$ be an optimal solution for (7.33)–(7.35) and let $\begin{pmatrix} y^* \\ z^* \end{pmatrix} \in \mathbb{R}^{K+L}$ be a Lagrange multiplier vector. Then, for some $\epsilon \geq 0$,*

$$\tilde{F}(\boldsymbol{x}^*, \mu, \epsilon) \leq \tilde{F}(\boldsymbol{x}, \mu, \epsilon) + L\mu\frac{\epsilon}{2}, \text{ for all } \boldsymbol{x} \in X, \tag{7.53}$$

*provided that $\mu \geq z_l^*$, for all $l = 1, 2, \ldots, L$.*

**Proof** Since $\boldsymbol{x}^*$ is an optimal solution for (7.33)–(7.35) and $(y^*, z^*)$ is a Lagrange multiplier vector whose existence is guaranteed by Assumption 7.2.2, we have (see, e.g., Rockafellar (1970, Theorem 28.3)):

$$\nabla f_0(\boldsymbol{x}^*) = -\sum_{k=1}^{K} y_k^* \nabla f_k((x^*)^k) - \sum_{l=1}^{L} z_l^* \nabla g_l(\boldsymbol{x}^*), \tag{7.54}$$

$$y_k^* f_k((x^*)^k) = 0 \text{ for all } k = 1, 2, \ldots, K, \ \ z_l^* g_l(\boldsymbol{x}^*) = 0 \text{ for all } l = 1, \ldots, L, \tag{7.55}$$

$$y_k^* \geq 0 \text{ for all } k = 1, 2, \ldots, K, \ \ z_l^* \geq 0 \text{ for all } l = 1, 2, \ldots, L, \tag{7.56}$$

$$f_k((x^*)^k) \leq 0 \text{ for all } k = 1, 2, \ldots, K, \ \ g_l(\boldsymbol{x}^*) \leq 0 \text{ for all } l = 1, 2, \ldots, L. \tag{7.57}$$

It is known that convexity and differentiability of any function $h$ guarantee that $h(x) \geq h(x^*) + \langle \nabla h(x^*), x - x^* \rangle$; see, e.g., Luenberger (1984). Therefore, using (7.54)–(7.55) in the definition of $F$, and by convexity and differentiability of the functions $f_k$, $k = 0, 1, 2, \ldots, K$ and $g_l$, $l = 1, 2, \ldots, L$, we obtain

$$F(\boldsymbol{x}, \mu) \geq f_0(\boldsymbol{x}^*) + \langle \nabla f_0(\boldsymbol{x}^*), \boldsymbol{x} - \boldsymbol{x}^* \rangle + \mu \sum_{l=1}^{L} g_l^+(\boldsymbol{x})$$

$$
\begin{aligned}
&= f_0(\boldsymbol{x}^*) - \langle \sum_{k=1}^{K} y_k^* \nabla f_k((x^*)^k), \boldsymbol{x} - \boldsymbol{x}^* \rangle \\
&\quad - \langle \sum_{l=1}^{L} z_l^* \nabla g_l(\boldsymbol{x}^*), \boldsymbol{x} - \boldsymbol{x}^* \rangle + \mu \sum_{l=1}^{L} g_l^+(\boldsymbol{x}) \\
&\geq f_0(\boldsymbol{x}^*) - \sum_{k=1}^{K} y_k^* (f_k(x^k) - f_k((x^*)^k)) \\
&\quad - \sum_{l=1}^{L} z_l^* (g_l(\boldsymbol{x}) - g_l(\boldsymbol{x}^*)) + \mu \sum_{l=1}^{L} g_l^+(\boldsymbol{x}) \\
&= f_0(\boldsymbol{x}^*) - \sum_{k=1}^{K} y_k^* f_k(x^k) - \sum_{l=1}^{L} z_l^* g_l(x^l) + \mu \sum_{l=1}^{L} g_l^+(\boldsymbol{x}).
\end{aligned}
$$

The first equality follows from (7.54); the second inequality follows from the convexity and differentiability of the functions $f_k$ for all $k = 1, 2, \ldots, K$, and $g_l$ for all $l = 1, 2, \ldots, L$; the second equality follows from (7.55). Since $g_l(\boldsymbol{x}) \leq g_l^+(\boldsymbol{x})$, we have:

$$
F(\boldsymbol{x}, \mu) \geq f_0(\boldsymbol{x}^*) - \sum_{k=1}^{K} y_k^* f_k(x^k) + \sum_{l=1}^{L} (\mu - z_l^*) g_l^+(\boldsymbol{x}).
$$

Now, since $\mu \geq z_l^*$ for all $l = 1, \ldots, L$, $y_k^* \geq 0$, and $f_k(x^k) \leq 0$ for all $k = 1, 2, \ldots, K$ because $\boldsymbol{x} \in X$, we get

$$
F(\boldsymbol{x}, \mu) \geq f_0(\boldsymbol{x}^*). \tag{7.58}
$$

But from Proposition 7.2.1 we have

$$
F(\boldsymbol{x}, \mu) - \tilde{F}(\boldsymbol{x}, \mu, \epsilon) \leq L\mu \frac{\epsilon}{2}. \tag{7.59}
$$

Rewriting (7.58) as

$$
f_0(\boldsymbol{x}^*) - F(\boldsymbol{x}, \mu) \leq 0, \tag{7.60}
$$

observing that $f_0(\boldsymbol{x}^*) = \tilde{F}(\boldsymbol{x}^*, \mu, \epsilon)$ since $\boldsymbol{x}^*$ is feasible, and adding up inequalities (7.59) and (7.60), we obtain

$$
\tilde{F}(\boldsymbol{x}^*, \mu, \epsilon) - \tilde{F}(\boldsymbol{x}, \mu, \epsilon) \leq L\mu \frac{\epsilon}{2}, \tag{7.61}
$$

which establishes the result. ■

We see that even if the optimal solution $x^*$ of (7.33)–(7.35) is not an optimal solution of (7.45), the objective function value computed at $x^*$ is greater than the true optimal value by at most $L\mu\epsilon/2$. As $\epsilon$ vanishes we recover the classical assertion that an optimal solution to the original problem coincides with an optimal solution to the penalty problem for some value of the penalty parameter $\mu$.

In the next proposition we obtain, as an immediate corollary from Proposition 7.2.4, a bound on the difference between the optimal value of the original problem (7.33)–(7.35) and the optimal value of the smoothed penalty problem.

**Proposition 7.2.5** *Let $\bar{\boldsymbol{x}}$ be a minimum point of $\min_{\boldsymbol{x}\in X} \tilde{F}(\boldsymbol{x},\mu,\epsilon)$, let $\boldsymbol{x}^*$ be an optimal point for the original problem (7.33)–(7.35), and let $z^* \in \mathbb{R}^L$ be a Lagrange multiplier vector associated with the inequalities $g_l(\boldsymbol{x}) \leq 0$ for all $l = 1, 2, \ldots, L$. Then*

$$0 \leq f_0(\boldsymbol{x}^*) - \tilde{F}(\bar{\boldsymbol{x}},\mu,\epsilon) \leq L\mu\frac{\epsilon}{2}, \tag{7.62}$$

*provided that $\mu > z_l^*$ for $l = 1, 2, \ldots, L$.*

We can now prove the following result:

**Proposition 7.2.6** *Let $\{\epsilon_\nu\}_{\nu=0}^{\infty}$ be a sequence of positive numbers with $\lim_{\nu\to\infty} \epsilon_\nu = 0$, and assume that for some $\mu > 0$, $\boldsymbol{x}^\nu$ is a solution to*

$$\min_{\boldsymbol{x}\in X} \tilde{F}(\boldsymbol{x},\mu,\epsilon_\nu). \tag{7.63}$$

*Also let $\bar{\boldsymbol{x}}$ be an accumulation point of the sequence $\{\boldsymbol{x}^\nu\}$, then*

$$F(\bar{\boldsymbol{x}},\mu) = F(\boldsymbol{x}(\mu),\mu), \tag{7.64}$$

*where $\boldsymbol{x}(\mu)$ is an optimal solution to $\min_{\boldsymbol{x}\in X} F(\boldsymbol{x},\mu)$.*

**Proof** The result is obtained from the continuity of $F$, the uniform convergence of (7.39), and inequalities (7.46) as follows. From (7.46) we have, for all $\nu \geq 0$, $F(\boldsymbol{x}(\mu),\mu) \geq \tilde{F}(\boldsymbol{x}^\nu,\mu,\epsilon_\nu)$. Using (7.39), remembering that the convergence is uniform, and by the definitions of $\tilde{F}$ and $F$, we obtain

$$\lim_{s\to\infty} \tilde{F}(\boldsymbol{x}^{\nu_s},\mu,\epsilon_{\nu_s}) = F(\bar{\boldsymbol{x}},\mu), \tag{7.65}$$

where $\lim_{s\to\infty} \boldsymbol{x}^{\nu_s} = \bar{\boldsymbol{x}}$ is the subsequence converging to $\bar{\boldsymbol{x}}$. Now, (7.46) yields $F(\boldsymbol{x}(\mu),\mu) \geq \tilde{F}(\boldsymbol{x}^{\nu_s},\mu,\epsilon_{\nu_s})$ which, together with (7.65) gives

$$F(\bar{\boldsymbol{x}},\mu) \geq F(\boldsymbol{x}(\mu),\mu).$$

The opposite inequality follows from the definition of $\boldsymbol{x}(\mu)$. ■

Therefore, as the smoothing parameter $\epsilon$ tends to 0, and if the penalty parameter $\mu$ is specified larger than the threshold value $\mu^*$, then the so-lutions to (7.45) and (7.33)–(7.35) would coincide. Hence, we obtain a solution to the original problem (7.33)–(7.35) by solving smooth penalty problems for a decreasing sequence of smoothing parameters $\epsilon_\nu$ and increasing penalty parameters $\mu_\nu$, since the threshold value $\mu^*$ is not known in practice. Also, driving $\epsilon$ to zero is not desirable since we would recover the nondifferentiable $\ell_1$ penalty function. Therefore, we cannot expect to obtain a solution to the original problem by solving a linear-quadratic penalty problem, but we obtain an $\epsilon$-feasible solution instead. In the next section we show that $\epsilon$-feasibility is attained if the penalized problem (7.45) is solved using a penalty parameter equal to a constant multiple of the threshold for the nondifferentiable $\ell_1$ penalty function.

### 7.2.2 $\epsilon$-exactness properties of the LQP function

We study now the approximate *exactness* properties that the smooth LQP penalty function inherits from its nondifferentiable counterpart. We show, in particular, that $\epsilon$-feasibility is achieved for some (finite) value of the penalty parameter.

First we look at the optimality conditions of problems (7.33)–(7.35) and (7.45), respectively.

#### Optimality Conditions for (7.33)–(7.35)

As stated earlier in the proof of Proposition 7.2.4, for $\boldsymbol{x}^*$ to be an optimal solution of (7.33)–(7.35) and for $(y^*, z^*)$ to be a Lagrange multiplier vector, it is necessary and sufficient that $(\boldsymbol{x}^*, y^*, z^*)$ satisfy (7.54)–(7.55). The existence of $(y^*, z^*)$ is guaranteed under Assumption 7.2.2.

#### Optimality Conditions for (7.45)

We introduce $\mathcal{A} \doteq \{l \mid g_l(\boldsymbol{x}) \leq \epsilon,\ 1 \leq l \leq L\}$ to denote the set of indices corresponding to constraints that are satisfied or violated to within $\epsilon$ at $\boldsymbol{x}$, and $\mathcal{V} \doteq \{l \mid g_l(\boldsymbol{x}) > \epsilon,\ 1 \leq l \leq L\}$, for the set of indices corresponding to constraints violated beyond $\epsilon$ at $\boldsymbol{x}$. Again using Rockafellar (1970, Theorem 28.3) we obtain: for $\boldsymbol{x}^*$ to be an optimal solution of (7.45) and for $y^*$ to be a Lagrange multiplier vector, it is necessary and sufficient that $(\boldsymbol{x}^*, y^*)$ satisfy the conditions

$$\nabla f_0(\boldsymbol{x}^*) = -\sum_{l \in \mathcal{A}} \frac{\mu g_l^+(\boldsymbol{x}^*)}{\epsilon} \nabla g_l(\boldsymbol{x}^*) - \sum_{l \in \mathcal{V}} \mu \nabla g_l(\boldsymbol{x}^*) - \sum_{k=1}^{K} y_k^* \nabla f_k((x^*)^k), \tag{7.66}$$

$$y_k^* f_k((x^*)^k) = 0 \text{ for all } k = 1, 2, \ldots, K, \tag{7.67}$$

$$y_k^* \geq 0 \text{ for all } k = 1, 2, \ldots, K, \tag{7.68}$$

$$f_k((x^*)^k) \leq 0 \text{ for all } k = 1, 2, \ldots, K. \tag{7.69}$$

The existence of $y^*$ is guaranteed under Assumption 7.2.2.

These optimality conditions for problems (7.33)–(7.35) and (7.45) have some features in common. Let us take a pair of vectors $(\tilde{\boldsymbol{x}}, \tilde{y})$ satisfying conditions (7.66)–(7.68) for some $\mu$ and $\epsilon$. Consider also an estimate of the Lagrange multiplier $z$, denoted by $\tilde{z}$, obtained from

$$\tilde{z}_l = \frac{dp_\epsilon}{dt}(g_l(\tilde{\boldsymbol{x}})) \text{ for all } l = 1, 2, \ldots, L. \tag{7.70}$$

Then the triplet $(\tilde{\boldsymbol{x}}, \tilde{y}, \tilde{z})$ would satisfy the optimality conditions for (7.33)–(7.35), with the exception of the complementary slackness conditions given by $z_l^* g_l(\boldsymbol{x}^*) = 0$ for all $l = 1, 2, \ldots, L$, and the feasibility conditions given by $g_l(\boldsymbol{x}^*) \leq 0$ for all $l = 1, 2, \ldots, L$. If $\mu$ and $\epsilon$ are chosen such that $\tilde{\boldsymbol{x}}$ is $\epsilon$-feasible, then the error in complementary slackness is bounded by the quantity $\mu\epsilon$. To see this, observe that all constraints are satisfied to within $\epsilon$ and that estimates for the Lagrange multiplier vector are given by

$$\tilde{z}_l = \frac{\mu \max(0, \tilde{u}_l)}{\epsilon} \quad \text{for all } l = 1, 2, \ldots, L,$$

where $\tilde{u}_l \doteq g_l(\tilde{\boldsymbol{x}})$. Therefore, we have

$$\tilde{z}_l \tilde{u}_l \leq \max(0, \mu \frac{\tilde{u}_l^2}{\epsilon}) \text{ for all } l = 1, 2, \ldots, L.$$

However, since $\tilde{u}_l \leq \epsilon$ for $l = 1, 2, \ldots, L$, it follows that

$$\tilde{z}_l \tilde{u}_l \leq \mu\epsilon \;\text{ for all } l = 1, 2, \ldots, L, \tag{7.71}$$

and the assertion is verified.

These observations provide the justification for computing an $\epsilon$-feasible solution to the original problem by solving the smooth penalty problem (7.45).

In the remainder of this section we study the conditions under which a solution that is $\epsilon$-feasible for the original problem can be obtained by solving the linear-quadratic penalty problem. This result characterizes the threshold value of the penalty parameter $\mu$ and provides a sufficient condition for $\epsilon$-feasibility.

**Proposition 7.2.7** *Let $\boldsymbol{x}^*$ be an optimal solution to (7.33)–(7.35), and let $\begin{pmatrix} y^* \\ z^* \end{pmatrix} \in \mathbb{R}^{K+L}$ be a Lagrange multiplier vector. Let $\tilde{\boldsymbol{x}}$ be an optimal solution to the smooth penalty problem (7.45) for some $\mu$ and $\epsilon$. Then $\tilde{\boldsymbol{x}}$ is $\epsilon$-feasible if $\mu$ is chosen such that $\mu \geq \mu^*/\kappa$, where $\mu^*$ is a threshold value of the penalty parameter with $\mu^* > z_l^*$ for all $l = 1, 2, \ldots, L$, and $\kappa$ is a constant given by*

$$\kappa \doteq \frac{1 - \sqrt{L}}{1 - L}. \tag{7.72}$$

**Proof** Recall that $\mathcal{A}$ is the set of indices corresponding to constraints that are satisfied or violated to within $\epsilon$ at $\tilde{\boldsymbol{x}}$, and $\mathcal{V}$ is the set of indices corresponding to constraints violated beyond $\epsilon$ at $\tilde{\boldsymbol{x}}$. We will assume that card($\mathcal{V}$) (the cardinality of the set $\mathcal{V}$) is at least one (otherwise the proposition holds trivially) and argue by negation.

We write the objective function of the smooth penalty problem (7.45) at the assumed minimum $\tilde{\boldsymbol{x}}$ as

$$f_0(\tilde{\boldsymbol{x}}) + \mu \sum_{l \in \mathcal{V}} (g_l(\tilde{\boldsymbol{x}}) - \frac{\epsilon}{2}) + \mu \sum_{l \in \mathcal{A}} \frac{g_l^+(\tilde{\boldsymbol{x}})^2}{2\epsilon}. \tag{7.73}$$

This can be so written because when $l \in \mathcal{V}$ then $g_l(\tilde{\boldsymbol{x}}) > \epsilon > 0$, thus, $g_l^+(\tilde{\boldsymbol{x}}) = g_l(\tilde{\boldsymbol{x}})$; however if $l \in \mathcal{A}$ then $g_l(\tilde{\boldsymbol{x}}) \leq \epsilon$, except that when $g_l(\tilde{\boldsymbol{x}}) \leq 0$ then $p_\epsilon(g_l(\tilde{\boldsymbol{x}})) = 0$ by definition.

We consider first the linear term (i.e., the second summand) in (7.73), and rewrite it as

$$\mu \sum_{l \in \mathcal{V}} (g_l^+(\tilde{\boldsymbol{x}}) - \frac{\epsilon}{2}) = \mu^* \sum_{l \in \mathcal{V}} g_l^+(\tilde{\boldsymbol{x}}) + (\mu \sum_{l \in \mathcal{V}} (g_l^+(\tilde{\boldsymbol{x}}) - \frac{\epsilon}{2}) - \mu^* \sum_{l \in \mathcal{V}} g_l^+(\tilde{\boldsymbol{x}})).$$

Considering the term in parentheses for $\mu = \mu^*/\kappa$ we get

$$\mu \sum_{l \in \mathcal{V}} (g_l^+(\tilde{\boldsymbol{x}}) - \frac{\epsilon}{2}) - \mu^* \sum_{l \in \mathcal{V}} g_l^+(\tilde{\boldsymbol{x}}) = \mu^* \sum_{l \in \mathcal{V}} (\frac{(g_l^+(\tilde{\boldsymbol{x}}) - \frac{\epsilon}{2})}{\kappa} - g_l^+(\tilde{\boldsymbol{x}})). \tag{7.74}$$

Now defining $t \doteq g_l^+(\tilde{\boldsymbol{x}})$, the term under the summation on the right-hand side can be written as:

$$\frac{t - \frac{\epsilon}{2}}{\kappa} - t = t(\frac{1}{\kappa} - 1) - \frac{\epsilon}{2\kappa}.$$

Since $t = g_l(\tilde{\boldsymbol{x}}) > \epsilon$, for all $l \in \mathcal{V}$, and since from the definition of $\kappa$ in (7.72) we have $0 < \kappa < \frac{1}{2}$ for all $L > 1$, we obtain

$$\frac{t - \frac{\epsilon}{2}}{\kappa} - t > \frac{\epsilon}{2\kappa} - \epsilon > 0.$$

Therefore, from (7.74) we have

$$\mu \sum_{l \in \mathcal{V}} (g_l^+(\tilde{\boldsymbol{x}}) - \frac{\epsilon}{2}) > \mu^* (\frac{\epsilon}{2\kappa} - \epsilon) \operatorname{card}(\mathcal{V}) + \mu^* \sum_{l \in \mathcal{V}} g_l^+(\tilde{\boldsymbol{x}}). \tag{7.75}$$

Now we consider the quadratic term in (7.73) which can be rewritten as

$$\mu \sum_{l \in \mathcal{A}} \frac{g_l^+(\tilde{\boldsymbol{x}})^2}{2\epsilon} = \mu^* \sum_{l \in \mathcal{A}} g_l^+(\tilde{\boldsymbol{x}}) + (\mu \sum_{l \in \mathcal{A}} \frac{g_l^+(\tilde{\boldsymbol{x}})^2}{2\epsilon} - \mu^* \sum_{l \in \mathcal{A}} g_l^+(\tilde{\boldsymbol{x}})).$$

Again considering the term in parentheses, for $\mu = \mu^*/\kappa$ we get

$$\mu \sum_{l \in \mathcal{A}} \frac{g_l^+(\tilde{\boldsymbol{x}})^2}{2\epsilon} - \mu^* \sum_{l \in \mathcal{A}} g_l^+(\tilde{\boldsymbol{x}}) = \mu^* \sum_{l \in \mathcal{A}} (\frac{g_l^+(\tilde{\boldsymbol{x}})^2}{2\epsilon\kappa} - g_l^+(\tilde{\boldsymbol{x}})). \tag{7.76}$$

Using $t \doteq g_l^+(\tilde{\boldsymbol{x}})$ we see that the term in the summation of the right-hand side has the form $(t^2/2\epsilon\kappa) - t$, which is a convex function in $t$. Its minimum is attained at $t^* = \epsilon\kappa$, at which point the function takes the value $-\epsilon\kappa/2$. Since this lower bound is negative by the positivity of $\epsilon$ and $\kappa$, the following inequality follows from (7.76),

$$\mu \sum_{l \in \mathcal{A}} \frac{g_l(\tilde{\boldsymbol{x}})^2}{2\epsilon} \geq \mu^* \sum_{l \in \mathcal{A}} g_l^+(\tilde{\boldsymbol{x}}) - \mu^* \frac{\epsilon\kappa}{2} \operatorname{card}(\mathcal{A}). \tag{7.77}$$

Now combining inequalities (7.75) and (7.77) we get the following result for the objective function value of the smooth penalty problem evaluated at $\tilde{\boldsymbol{x}}$, whose derivation is explained below:

$$\begin{aligned}
& f_0(\tilde{\boldsymbol{x}}) + \mu \sum_{l \in \mathcal{V}} (g_l(\tilde{\boldsymbol{x}}) - \frac{\epsilon}{2}) + \mu \sum_{l \in \mathcal{A}} \frac{g_l^+(\tilde{\boldsymbol{x}})^2}{2\epsilon} \\
& > f_0(\tilde{\boldsymbol{x}}) + \mu^* \sum_{l=1}^{L} g_l^+(\tilde{\boldsymbol{x}}) + \mu^* (\frac{\epsilon}{2\kappa} - \epsilon) - (L-1)\mu^* (\frac{\epsilon\kappa}{2}) \\
& = f_0(\tilde{\boldsymbol{x}}) + \mu^* \sum_{l=1}^{L} g_l^+(\tilde{\boldsymbol{x}}) \\
& > f_0(\tilde{\boldsymbol{x}}) + \sum_{k=1}^{K} y_k^* f_k(\tilde{x}^k) + \sum_{l=1}^{L} z_l^* g_l(\tilde{\boldsymbol{x}}) \\
& \geq f_0(\boldsymbol{x}^*) + \sum_{k=1}^{K} y_k^* f_k((x^*)^k) + \sum_{l=1}^{L} z_l^* g_l(\boldsymbol{x}^*) \\
& = f_0(\boldsymbol{x}^*) = f_0(\boldsymbol{x}^*) + \sum_{l=1}^{L} p_\epsilon(g_l(\boldsymbol{x}^*)).
\end{aligned}$$

The first inequality follows by adding (7.75) and (7.77), based on the assumption that $\operatorname{card}(\mathcal{V}) \geq 1$, and using the fact that $\operatorname{card}(\mathcal{A}) \leq L - 1$. The next equality follows from $(\epsilon/2\kappa) - \epsilon - (L-1)\epsilon\kappa/2 = 0$, by definition of $\kappa$. The subsequent inequality follows from the definition of $\mu^*$ and from (7.67). The third inequality along with the subsequent equality follow from the definition of a Lagrange multiplier vector. The last equality follows from the feasibility of $x^*$. But this result is a contradiction since $\tilde{\boldsymbol{x}}$ was assumed to be an optimal solution to (7.45). Hence the assumption

card($\mathcal{V}$) $\geq 1$ is violated, i.e., the set $\mathcal{V}$ must be empty. ■

Therefore, for a given problem the threshold value of $\mu$ required to attain $\epsilon$-feasibility is a constant independent of $\epsilon$. The existence of a threshold value of the penalty parameter $\mu$ indicates that an $\epsilon$-feasible solution can be computed in a finite number of minimizations. An important consequence of Proposition 7.2.7 is that we can now characterize the conditions under which the upper bound on the difference between the optimal objective values of problems (7.33)–(7.35) and (7.45), as stated in Proposition 7.2.3, is achieved. This result is now summarized in the following proposition.

**Proposition 7.2.8** *Let $x^*$ be an optimal solution for (7.33)–(7.35) and let $\begin{pmatrix} y^* \\ z^* \end{pmatrix} \in \mathbb{R}^{K+L}$ be the Lagrange multiplier vector. Let $\mu^*$ be such that $\mu^* > z_l^*$ for all $l = 1, 2, \ldots, L$. Furthermore let $\tilde{x}$ be an optimal solution to (7.45), with $\mu \doteq \mu^*/\kappa$ where $\kappa \doteq (1 - \sqrt{L})/(1 - L)$. Then, $0 \leq f_0(x^*) - f_0(\tilde{x}) \leq L\mu\epsilon$.*

**Proof** Since $\kappa < 1$ it follows that $\mu > \mu^*$. Hence, solving (7.37) using $\mu = \mu^*/\kappa$ produces an optimal solution to the original problem. The rest of the proof follows from Proposition 7.2.3. ■

Propositions 7.2.7 and 7.2.8 motivate the procedure for controlling the accuracy of the solution in the LQP algorithm (Algorithm 7.2.1). The procedure starts with appropriately chosen $\mu_\nu$ and $\epsilon_\nu$. The value of the penalty parameter $\mu_\nu$ is increased according to some criteria, while the smoothing parameter $\epsilon_\nu$ can be decreased after each penalty minimization in Step 1 is completed. More precisely, if the solution of the penalty problem is $\epsilon$-feasible, this is an indication that the penalty parameter $\mu_\nu$ is large enough. Therefore, if the smoothing parameter $\epsilon$ is below the final accuracy tolerance, $\epsilon_{min}$, the algorithm terminates. In this case Proposition 7.2.8 provides an upper bound on the difference between the optimal value of the original problem and the optimal value of the smoothed penalty problem. If the solution of the penalty problem is not $\epsilon$-feasible, then the penalty parameter $\mu_\nu$ is not large enough. Another round of penalty minimization is carried out with a larger value of the penalty parameter. The smoothing parameter $\epsilon_\nu$ can be left unchanged in this case.

## 7.3 Notes and References

The development of algorithms for the decomposition of structured optimization models into smaller and simpler subproblems dates back to the early days of linear programming. See Dantzig (1963) for an early account. A classification of decomposition algorithms is given by Geoffrion (1970), and a textbook treatment of algorithms for large-scale systems developed through the 1960s was provided by Lasdon (1970). The 1980s have witnessed several efforts in parallelizing the most popular decomposition algorithms, including the Dantzig-Wolfe algorithm in Dantzig and Wolfe (1960), and Benders decomposition in Benders (1962); see also Geoffrion (1972). For an early discussion see the introduction in Meyer and Zenios (1988). For parallelizations of the Dantzig-Wolfe algorithm see Ho, Lee and Sundarraj (1988), and for parallelizations of Benders decomposition see Ariyawansa (1991), Ariyawansa and Hudson (1991), Dantzig, Ho, and Infanger (1991), Nielsen and Zenios (1994), and Qi and Zenios (1994).

Significant recent developments that fit into the general framework of model decomposition methods are the parallel *constraint distribution* and parallel *variable distribution* methods of Ferris and Mangasarian (1991, 1994). Both methods have been proven effective for the solution of large-scale linear programs and were also efficient when implemented on parallel machines. See also Mangasarian (1995) for the related gradient distribution method for unconstrained optimization and Ferris (1994) for extensions to convex quadratic programs.

**7.1** For an alternative classification of decomposition algorithms to the one taken in this section see Geoffrion (1970).

**7.1.1** Penalty methods for the solution of constrained optimization problems are attributed to Courant (1962), and the barrier method was first suggested by Caroll (1961). A textbook treatment of penalty methods can be found in Bertsekas (1982). Barrier, or interior penalty methods were developed and popularized in Fiacco and McCormick (1968) and McCormick (1983). For an introduction to both methods see Luenberger (1984). A discussion of exact penalty methods can be found in Han and Mangasarian (1979). The use of split-variable formulations is fairly standard in large-scale optimization, e.g., Bertsekas and Tsitsiklis (1989, p.231). For applications see Rockafellar and Wets (1991), Mulvey and Vladimirou (1989, 1991), and Nielsen and Zenios (1993a, 1996a). The augmented Lagrangian was introduced in Arrow, Hurwicz, and Uzawa (1958) as a means for the convexification of nonconvex problems, and was further extended and analyzed by Rockafellar (1974, 1976a). An extensive treatment of this topic can be found in Bertsekas (1982). The method of multipliers was introduced, independently, by Hestenes (1969) and Powell (1969). Consult Section 4.5 for more notes and references on penalty and barrier methods.

**7.1.2** The Frank-Wolfe algorithm was developed by Frank and Wolfe (1956), and the simplicial decomposition algorithm was first suggested as a generalization to Frank-Wolfe by Holloway (1974). The representation of the master program as a nonlinear program with a simplex constraint is due to Von Hohenbalken (1977); see also Von Hohenbalken (1975) and Grinold (1982). A memory-efficient variant of the algorithm that maintains only a limited number of vertices is due to Hearn, Lawphongpanich, and Ventura (1984). Mulvey, Zenios, and Ahlfeld (1990) developed numerical procedures for solving the master program and specialized the algorithm for large-scale network problems. Nonlinear programs with a single equality constraint, and bounded variables (such as the master program in Algorithm 7.1.3), appear in many applications. Algorithms for their solution have been suggested by Helgason, Kennington, and Lass (1980) and Tseng (1990); and parallel algorithms were developed, implemented, and compared by Nielsen and Zenios (1992c).

The use of linearization algorithms in conjunction with barrier or penalty functions for the decomposition of structured programs for parallel computing was suggested by Schultz and Meyer (1991) and Pinar and Zenios (1992). Both studies report encouraging computational results on parallel machines for large-scale problems.

Diagonal quadratic approximations to nonseparable quadratic functions, and in particular for the augmented Lagrangian, were suggested by Stephanopoulos and Westerberg (1975) and further extended by Tatjewski (1989). The use of split-variable formulations in conjunction with augmented Lagrangians and diagonal quadratic approximations for the solution of large-scale structured programs was suggested by Mulvey and Ruszczyński (1992), and was subsequently applied to the solution of stochastic programming problems by Mulvey and Ruszczyński (1994) and Berger, Mulvey, and Ruszczyński (1994). The same references also report encouraging computational results with the use of the algorithm for solving large-scale problems on a distributed network of workstations.

**7.2** Smoothing approximations to exact penalty functions were introduced by Bertsekas (1973), who also suggested their use for the solution of minimax optimization problems. Zang (1980) also discusses smoothing techniques for the solution of minimax problems. Madsen and Nielsen (1993) introduced smoothing for the solution of $\ell_1$ estimation problems. We point out that the linear-quadratic penalty function is used in statistics to develop robust estimation procedures, and it is known as the Huber estimator; see Huber (1981). Linear-quadratic programming problems such as those arising from the smoothing approximations of this section are analyzed in Rockafellar (1987, 1990). The LQP algorithm for large-scale problems was developed by Zenios, Pinar, and Dembo (1994). Its properties were analyzed in Pinar and Zenios (1994b), and the algorithm was used for the parallel decomposition of multicommodity flow problems in Pinar and Zenios (1992, 1994a).

**7.2.2** The analysis of the properties of the LQP function is due to Pinar and Zenios (1994b). It is based on earlier works by Bertsekas (1975),

Charalambous (1978), and Charalambous and Conn (1978) who analyzed exact penalty functions, and on the work of Truemper (1975) who analyzed the quadratic penalty function.

CHAPTER 8

# Decompositions in Interior Point Algorithms

In 1984 N. Karmarkar of AT&T Bell Laboratories introduced an interior point algorithm for linear programming. The algorithm has a polynomial time complexity bound and has been established, following extensive computational experimentation, as a viable competitor to the classic simplex algorithm for the solution of large-scale linear programs. A flurry of research activities followed Karmarkar's work and several interior point algorithms were developed for linear programming, convex quadratic programming, convex programming, linear complementarity problems, and nonlinear complementarity problems.

Different mathematical tools have been employed to develop and analyze these algorithms. Depending on the theory underlying the mathematical analysis the algorithms are classified as potential reduction algorithms, path following algorithms, barrier function algorithms, affine scaling algorithms, and projective scaling algorithms.

It has been gradually realized that some of the algorithms can be obtained as special cases of a unified interior point method. It is beyond the scope of this chapter to cover all these developments, and relevant references are given in Section 8.4.

One common feature of the interior point algorithms is that they compute a search direction by solving a system of equations that involves the matrix $AA^{\mathrm{T}}$, where $A$ is the constraint matrix. For problems with structured constraints matrices it is possible to solve the system of equations using parallel matrix factorization procedures. Hence, interior point algorithms are well suited for parallel computations. In this respect they draw heavily from developments in parallel numerical linear algebra.

In this chapter we discuss one particular interior point algorithm—the *primal-dual path following algorithm*—for both linear programs and for quadratic programming problems. This is the algorithm implemented in some of the most efficient and widely used software systems such as OB1 (Optimization with Barrier 1) and its quadratic programming extension, OBN, the linear programming code ALPO, and LOQO. We show that, for structured problems, it is possible to develop parallel matrix factorization procedures for solving the required systems of equations. This results in improvement in the performance of an already efficient algorithm via the use of parallelism. Sections 8.1 and 8.2 develop the algorithms for linear and quadratic programs, respectively. The parallel matrix factorization techniques are the topic of Section 8.3. Notes and references in Section 8.4 conclude this chapter.

## 8.1 The Primal-Dual Path Following Algorithm for Linear Programming

Consider the *primal* linear programming problem:

$$\text{Minimize} \quad \langle c, x\rangle \tag{8.1}$$

$$\text{s.t.} \quad Ax = b, \tag{8.2}$$

$$x \geq 0. \tag{8.3}$$

The $m \times n$ constraint matrix $A$ that is assumed to have full row rank, the cost vector $c \in \mathbb{R}^n$, and the vector of resource coefficients $b \in \mathbb{R}^m$ are given. Associated with this problem is its *dual* linear programming problem, defined as follows:

$$\text{Maximize} \quad \langle b, y\rangle \tag{8.4}$$

$$\text{s.t.} \quad A^T y + z = c, \tag{8.5}$$

$$z \geq 0. \tag{8.6}$$

The dual formulation uses the vector of *slack variables* $z \in \mathbb{R}^n$ in order to express its constraints as equalities.

The primal-dual path following algorithm applies first a logarithmic barrier function (see Example 4.2.2) to formulate the following pair of problems: a logarithmic barrier problem for the primal program, given by,

$$\text{Minimize} \quad \langle c, x\rangle - \mu \sum_{j=1}^{n} \log x_j \tag{8.7}$$

$$\text{s.t.} \quad Ax = b, \tag{8.8}$$

where log is the natural logarithmic function and a barrier problem for the dual program

$$\text{Maximize} \quad \langle b, y\rangle + \mu \sum_{j=1}^{n} \log z_j \tag{8.9}$$

$$\text{s.t.} \quad A^T y + z = c. \tag{8.10}$$

The algorithm then solves these two barrier problems for different values of the parameter $\mu$, as $\mu \to 0$, and obtains sequences of feasible points that converge to the solutions of the original primal and dual problems, respectively (see Section 4.2).

Let us now fix a value of $\mu > 0$, which is the same for the primal and dual barrier problems, at a point that satisfies $x > 0$ and $z > 0$. We can then obtain the first-order conditions for simultaneous optimality of both the primal and the dual barrier problems. In particular, the solution to the primal barrier problem is the unique point where the gradient of the associated Lagrangian,

$$L_P(x,y) \doteq \langle c,x\rangle - \mu\sum_{j=1}^{n}\log x_j + \langle y, b-Ax\rangle,$$

vanishes. This is the *critical point* of the Lagrangian, and it is given by the solution of the following system:

$$c - \mu X^{-1}\mathbf{1} - A^{\mathrm{T}}y = 0, \tag{8.11}$$

$$Ax = b. \tag{8.12}$$

Here $X^{-1}$ denotes the inverse of the $n \times n$ diagonal matrix $X$ defined as $X \doteq \text{diag}\,(x_1, x_2, \ldots, x_n)$, and $\mathbf{1} \in \mathbb{R}^n$ is the all ones vector. The solution to the dual barrier problem is the unique point where the gradient of the Lagrangian

$$L_D(y,z,x) \doteq \langle b,y\rangle + \mu\sum_{j=1}^{n}\log\, z_j + \langle x, c - A^{\mathrm{T}}y - z\rangle$$

vanishes. This point is obtained by solving:

$$Ax = b, \tag{8.13}$$

$$\mu Z^{-1}\mathbf{1} - x = 0, \tag{8.14}$$

$$A^{\mathrm{T}}y + z = c, \tag{8.15}$$

where $Z^{-1}$ is the inverse of the diagonal matrix $Z \doteq \text{diag}\,(z_1, z_2, \ldots, z_n)$.

Combining the first-order conditions for the primal and dual problems, and using the substitution $\mu Z^{-1}\mathbf{1} = X\mathbf{1}$ from (8.14) we obtain the following system of equations that characterizes, simultaneously, the optimum of both the primal and the dual barrier problems:

$$Ax = b, \tag{8.16}$$

$$A^{\mathrm{T}}y + z = c, \tag{8.17}$$

$$XZ\mathbf{1} = \mu\mathbf{1}. \tag{8.18}$$

Equation (8.16) is part of the primal feasibility requirement and (8.17) is part of the dual feasibility requirement. Equation (8.18), relating the dual slack variables $z$ with the primal vector $x$, is called *$\mu$-complementarity.* The duality theory of mathematical programming guarantees that if either the primal or the dual polytope, defined by the feasibility sets of (8.2)–(8.3) and (8.5)–(8.6) respectively, is bounded and has a nonempty interior, then there exists a unique solution to (8.16)–(8.18) for which $x \geq 0$, $z \geq 0$, for all $\mu > 0$. See, e.g., Goldfarb and Todd (1989, p. 156), Fiacco and McCormick (1968, Chapter 3), or Kojima et al. (1991).

For a given $\mu > 0$ the solution to system (8.16)–(8.18) is denoted by $(x_\mu, y_\mu, z_\mu)$ and the one-parameter family of solutions $\{(x_\mu, y_\mu, z_\mu)\}$ is called the *central path.* As $\mu$ tends to zero this path converges to $(x^*, y^*, z^*)$, where $x^*$ is optimal for the original primal problem (8.1)–(8.3), and $(y^*, z^*)$ is optimal for the dual problem (8.4)–(8.6). The objective values $\langle c, x_\mu \rangle$ and $\langle b, y_\mu \rangle$ converge, as $\mu$ tends to zero, to the common objective value of the primal and dual problems.

A primal-dual path following algorithm iteratively tracks the central path. It starts from some $x > 0, z > 0$, and a positive $\mu$, and takes a single step of Newton's method on the system (8.16)–(8.18) to find a point close to the central path. It then reduces $\mu$ and takes another Newton step starting from the previous iterate, and so forth.

To solve the system (8.16)–(8.18) we start from a point $(x, y, z)$ with $x > 0,\ z > 0$, and move to a new point $(x+\Delta x, y+\Delta y, z+\Delta z)$. Substituting the new point in (8.16)–(8.18) we obtain a system of equations that is complicated by the fact that (8.18) is nonlinear. The nonlinear equation $(X + \Delta X)(Z + \Delta Z)\mathbf{1} = \mu\mathbf{1}$ is then linearized by ignoring the cross-product term $\Delta X \Delta Z$. The resulting system of linear equations is written as

$$A\Delta x = \rho, \tag{8.19}$$

$$A^{\mathrm{T}}\Delta y + \Delta z = \sigma, \tag{8.20}$$

$$Z\Delta x + X\Delta z = \phi. \tag{8.21}$$

where $\rho \doteq b - Ax$, $\sigma \doteq c - A^{\mathrm{T}}y - z$, and $\phi \doteq \mu\mathbf{1} - XZ\mathbf{1}$.

To solve this system we use first (8.21) to obtain

$$\Delta z = X^{-1}(\phi - Z\Delta x), \tag{8.22}$$

and substitute it into (8.20) to obtain

$$-ZX^{-1}\Delta x + A^{\mathrm{T}}\Delta y = \sigma - X^{-1}\phi. \tag{8.23}$$

Extracting $\Delta x$ from (8.23) we obtain

$$\Delta x = -XZ^{-1}(\sigma - X^{-1}\phi - A^{\mathrm{T}}\Delta y), \tag{8.24}$$

which, upon substitution into (8.19), yields

$$(AXZ^{-1}A^{\mathrm{T}})\Delta y = \left(\rho + AXZ^{-1}(\sigma - X^{-1}\phi)\right). \tag{8.25}$$

Let $\Theta \doteq XZ^{-1}$ and $\psi \doteq \rho + A\Theta(\sigma - X^{-1}\phi)$. The *dual step direction* $\Delta y$ is then obtained by solving the system

$$(A\Theta A^{\mathrm{T}})\Delta y = \psi. \tag{8.26}$$

This is a symmetric positive definite system. (Recall that $A$ has full row rank, and $\Theta$ is positive definite since $X$ and $Z$ are positive definite.) It can be solved using any direct or iterative method for solving systems of equations. The solution of this system, using parallel matrix factorization techniques, is the topic of Section 8.3. Having obtained the dual step direction $\Delta y$ we can substitute it in (8.24) to obtain the *primal step direction* $\Delta x$, and finally substitute the primal step direction in (8.22) to obtain the step direction for the slack variables.

The primal-dual path following algorithm is thus summarized as follows.

#### Algorithm 8.1.1 The Primal-Dual Path Following Algorithm for Linear Programs

**Step 0:** (Initialization.) Start with a triplet $(x^0, y^0, z^0)$ satisfying $x^0 > 0$, $z^0 > 0$, and any $\mu_0 > 0$. Let $\nu = 0$.

**Step 1:** (Iterative Step.) Given $(x^\nu, y^\nu, z^\nu)$ and a barrier parameter $\mu_\nu$, let

$$\begin{aligned}
X_\nu &\doteq \mathrm{diag}\,(x_1^\nu, x_2^\nu, \ldots, x_n^\nu),\\
Z_\nu &\doteq \mathrm{diag}\,(z_1^\nu, z_2^\nu, \ldots, z_n^\nu),\\
\Theta_\nu &\doteq X_\nu Z_\nu^{-1},\\
\rho^\nu &\doteq b - Ax^\nu,\\
\sigma^\nu &\doteq c - A^{\mathrm{T}}y^\nu - z^\nu,\\
\phi^\nu &\doteq \mu_\nu \mathbf{1} - X_\nu Z_\nu \mathbf{1},\\
\psi^\nu &\doteq \rho^\nu + A\Theta_\nu(\sigma^\nu - X_\nu^{-1}\phi^\nu),
\end{aligned}$$

and solve for the dual step direction $\Delta y^\nu$:

$$(A\Theta_\nu A^{\mathrm{T}})\Delta y^\nu = \psi^\nu. \tag{8.27}$$

Compute the primal step direction $\Delta x^\nu$:

$$\Delta x^\nu = -\Theta_\nu(\sigma^\nu - X_\nu^{-1}\phi^\nu - A^{\mathrm{T}}\Delta y^\nu), \tag{8.28}$$

and the *slack variable direction* $\Delta z^\nu$ from

$$\Delta z^\nu = X_\nu^{-1}(\phi^\nu - Z_\nu \Delta x^\nu). \tag{8.29}$$

Update:

$$x^{\nu+1} = x^\nu + \alpha_{P,\nu}\Delta x^\nu, \tag{8.30}$$
$$y^{\nu+1} = y^\nu + \alpha_{D,\nu}\Delta y^\nu, \tag{8.31}$$
$$z^{\nu+1} = z^\nu + \alpha_{D,\nu}\Delta z^\nu, \tag{8.32}$$

where $\alpha_{P,\nu}$, $\alpha_{D,\nu}$ are primal and dual step lengths, respectively, restricted to the interval $(0,1]$. Reduce $\mu_\nu$ to $\mu_{\nu+1}$, update $\nu \leftarrow \nu + 1$, and repeat the iterative step.

We discuss now the choices of the step lengths $\alpha_{P,\nu}$ and $\alpha_{D,\nu}$ and of the barrier parameter $\mu_\nu$. Other features of the algorithm that need to be properly addressed for an efficient implementation are the treatment of upper bounds, the treatment of free variables, and the computation of an initial feasible solution. These technical issues are not discussed here, as they are not essential to understanding the algorithm and do not alter its parallel implementation. References to the relevant literature are given in Section 8.4.

### 8.1.1 Choosing the step lengths

The primal and dual step lengths are chosen in such a way that both the primal and the slack variables remain positive. To determine the maximum step length that will keep the iterates positive define first

$$\alpha_{P,j} \doteq \max\{\alpha > 0 \mid x_j^\nu + \alpha \Delta x_j^\nu \geq 0\}, \quad \text{for all } j = 1, 2, \ldots, n,$$
$$\hat{\alpha}_{P,\nu} \doteq \min_{1\leq j\leq n} \alpha_{P,j},$$

and

$$\alpha_{D,j} \doteq \max\{\alpha > 0 \mid z_j^\nu + \alpha \Delta z_j^\nu \geq 0\}, \quad \text{for all } j = 1, 2, \ldots, n,$$
$$\hat{\alpha}_{D,\nu} \doteq \min_{1\leq j\leq n} \alpha_{D,j}.$$

Then the step lengths are defined by $\alpha_{P,\nu} \doteq \delta\hat{\alpha}_{P,\nu}$ and $\alpha_{D,\nu} \doteq \delta\hat{\alpha}_{D,\nu}$, where $0 < \delta < 1$. A typical value for $\delta$, used in practice, is 0.9995.

### 8.1.2 Choosing the barrier parameter

The barrier parameter $\mu_\nu$ determines the severity of the barrier term. For large $\mu_\nu$ the search direction points away from the boundary of the feasible region, and therefore large step lengths are allowed. As a solution is approached the value of the barrier parameter is reduced so that the iterates can approach the boundary, since it is known that the solutions to a linear program are on the boundary.

An acceptable (theoretically) value for the initial barrier parameter $\mu_0$ is $2^{O(L)}$, where $O(L)$ is the "order of greatness" of $L$ which, in turn, is a measure of the size of the problem at hand. $K = O(L)$ means that there exists a constant $M$ such that $| K | \leq ML$, for all problems above a certain size $L$. The penalty parameter $\mu_\nu$ is reduced at every step by the recursion

$$\mu_{\nu+1} = (1 - \frac{0.1}{\sqrt{n}})\mu_\nu.$$

This choice of the barrier parameter assures polynomial convergence of the algorithm, but the algorithm is very slow in practice. An updating formula for setting the barrier parameter, which works well in practice, is given by

$$\mu_\nu = \frac{\langle c, x^\nu \rangle - \langle b, y^\nu \rangle}{n^2}.$$

Thus $\mu_\nu$ is large when far from the optimum (as measured by a large duality gap $\langle c, x^\nu \rangle - \langle b, y^\nu \rangle$), and is reduced rapidly as the optimum is approached.

## 8.2 The Primal-Dual Path Following Algorithm for Quadratic Programming

A primal-dual path following algorithm for the quadratic programming problem is derived using the ideas discussed in the previous section for the linear program. Namely, a barrier function is employed to eliminate the inequality constraints and then Newton's method is used to solve the optimality conditions of the barrier problem. For separable quadratic functions the algorithms for the linear and the quadratic problems are virtually identical.

Consider the *primal* quadratic programming problem

$$\text{Minimize} \quad \langle c, x \rangle + \frac{1}{2}\langle x, Qx \rangle \tag{8.33}$$

$$\text{s.t.} \quad Ax = b, \tag{8.34}$$

$$x \geq 0. \tag{8.35}$$

The $m \times n$ constraint matrix $A$ has full row rank, $c \in \mathbb{R}^n$ is the cost vector, and $b \in \mathbb{R}^m$ is the vector of resource coefficients. The $n \times n$ matrix $Q$ is positive definite.

Associated with this problem is a *dual* quadratic program, defined by

$$\text{Maximize} \quad \langle b, y \rangle - \frac{1}{2}\langle x, Qx \rangle \tag{8.36}$$

$$\text{s.t.} \quad A^T y + z - Qx = c, \tag{8.37}$$

$$z \geq 0. \tag{8.38}$$

The dual formulation uses slack variables $z \in \mathbb{R}^n$ in order to express constraints as equalities.

The logarithmic barrier problem for the primal program is given by

$$\text{Minimize} \quad \langle c, x\rangle + \frac{1}{2}\langle x, Qx\rangle - \mu \sum_{j=1}^{n} \log x_j \tag{8.39}$$

$$\text{s.t.} \quad Ax = b. \tag{8.40}$$

Similarly, the barrier problem for the dual program is given by

$$\text{Maximize} \quad \langle b, y\rangle - \frac{1}{2}\langle x, Qx\rangle + \mu \sum_{j=1}^{n} \log z_j \tag{8.41}$$

$$\text{s.t.} \quad A^{\mathrm{T}} y + z - Qx = c. \tag{8.42}$$

We can solve the barrier problems (8.39)–(8.40) or (8.41)–(8.42) for various values of the parameter $\mu$, as $\mu \to 0$, and obtain sequences of feasible points that converge to the solutions of the original primal and dual quadratic programs, respectively.

Let us fix a value of $\mu > 0$ that is the same for the primal and dual barrier problems at a point that satisfies $x > 0$ and $z > 0$. We can then obtain the first-order conditions for simultaneous optimality of the primal and the dual barrier problems. In particular, the solution to the primal barrier problem is the unique critical point of the associated Lagrangian

$$L_P(x, y) \doteq \langle c, x\rangle + \frac{1}{2}\langle x, Qx\rangle - \mu \sum_{j=1}^{n} \log x_j + \langle y, b - Ax\rangle.$$

The critical point is obtained by setting the gradient of the associated Lagrangian equal to zero, i.e., it is the solution of the following system:

$$c + Qx - \mu X^{-1}\mathbf{1} - A^{\mathrm{T}} y = 0, \tag{8.43}$$

$$Ax = b. \tag{8.44}$$

Similarily, the solution to the dual barrier problem is the unique critical point of the Lagrangian

$$L_D(y, z, x) \doteq \langle b, y\rangle - \frac{1}{2}\langle x, Qx\rangle + \mu \sum_{j=1}^{n} \log z_j + \langle x, c + Qx - A^{\mathrm{T}} y - z\rangle.$$

The first-order conditions are:

$$Ax = b, \tag{8.45}$$

$$\mu Z^{-1}\mathbf{1} - x = 0, \tag{8.46}$$
$$-Qx + A^{\mathrm{T}}y + z = c. \tag{8.47}$$

Doing some simple algebraic manipulations on the first-order conditions for the primal and dual barrier problems we obtain the following system of equations that characterizes, simultaneously, the optimum of both problems (8.39)–(8.40) and (8.41)–(8.42):

$$Ax = b, \tag{8.48}$$
$$-Qx + A^{\mathrm{T}}y + z = c, \tag{8.49}$$
$$XZ\mathbf{1} = \mu\mathbf{1}. \tag{8.50}$$

Similar to the linear programming case, it is known that if either the primal or the dual polytope that describes the feasible region is bounded and has a nonempty interior, then there exists a unique solution to (8.48)–(8.50) in the domain given by $x \geq 0$, $z \geq 0$ for all $\mu > 0$.

As in the linear programming case, the primal-dual path following algorithm for quadratic programs tracks the central path, defined by solving the system (8.48)–(8.50), for decreasing values of $\mu$. To solve this system we start from $(x, y, z)$ with $x > 0$, $z > 0$, and move to a new point $(x + \Delta x, y + \Delta y, z + \Delta z)$. Substituting the new point in (8.48)–(8.50) we obtain a system of equations that is complicated by the nonlinear equation (8.50). This equation, i.e., $(X + \Delta X)(Z + \Delta Z)\mathbf{1} = \mu\mathbf{1}$, is again linearized by ignoring the cross-product term $\Delta X \Delta Z$. The resulting system of linear equations is written as

$$A\Delta x = \rho, \tag{8.51}$$
$$-Q\Delta x + A^{\mathrm{T}}\Delta y + \Delta z = \sigma, \tag{8.52}$$
$$Z\Delta x + X\Delta z = \phi, \tag{8.53}$$

where $\rho \doteq b - Ax$, $\sigma \doteq c + Qx - A^{\mathrm{T}}y - z$, and $\phi \doteq \mu\mathbf{1} - XZ\mathbf{1}$. To solve this system we use first (8.53) to obtain

$$\Delta z = X^{-1}(\phi - Z\Delta x), \tag{8.54}$$

and then we eliminate $\Delta z$ in (8.52) to obtain

$$-(Q + ZX^{-1})\Delta x + A^{\mathrm{T}}\Delta y = \sigma - X^{-1}\phi. \tag{8.55}$$

Solving this system for $\Delta x$ we get

$$\Delta x = -\Theta(\sigma - X^{-1}\phi - A^{\mathrm{T}}\Delta y), \tag{8.56}$$

where $\Theta \doteq (Q + ZX^{-1})^{-1}$. Using now (8.51) to eliminate $\Delta x$ we obtain the system of equations for the dual step direction $\Delta y$ :

$$(A\Theta A^{\mathrm{T}})\Delta y = \psi, \tag{8.57}$$

where $\psi \doteq \rho + A\Theta(\sigma - X^{-1}\phi)$.

The only difference between the dual step direction calculation for the linear program (equation (8.26)) and for the quadratic program (equation (8.57)) is the definition of the matrix $\Theta$: for linear programs it is $\Theta = XZ^{-1}$ while for quadratic programs $\Theta = (Q + ZX^{-1})^{-1}$. When $Q$ is diagonal the computation of $\Theta$ is easily obtained as the inverse of a positive diagonal matrix. The complexity of the computation of $\Delta y$ is then the same for linear and quadratic programs. When $Q$ is not diagonal then $\Theta$ is more difficult to compute and may even be dense, in which case the matrix $A\Theta A^{\mathrm{T}}$ will also be dense.

The primal-dual path following algorithm is thus summarized as follows.

**Algorithm 8.2.1 The Primal-Dual Path Following Algorithm for Quadratic Programs**

**Step 0:** (Initialization.) Start with a triplet $(x^0, y^0, z^0)$ satisfying $x^0 > 0$, $z^0 > 0$, and any $\mu_0 > 0$. Let $\nu = 0$.

**Step 1:** (Iterative Step.) Given $(x^\nu, y^\nu, z^\nu)$ and a barrier parameter $\mu_\nu$, let

$$\begin{aligned}
X_\nu &\doteq \mathrm{diag}\,(x_1^\nu, x_2^\nu, \ldots, x_n^\nu),\\
Z_\nu &\doteq \mathrm{diag}\,(z_1^\nu, z_2^\nu, \ldots, z_n^\nu),\\
\Theta_\nu &\doteq (Q + Z_\nu X_\nu^{-1})^{-1},\\
\rho^\nu &\doteq b - Ax^\nu,\\
\sigma^\nu &\doteq c + Qx^\nu - A^{\mathrm{T}}y^\nu - z^\nu,\\
\phi^\nu &\doteq \mu_\nu \mathbf{1} - X_\nu Z_\nu \mathbf{1},\\
\psi^\nu &\doteq \rho^\nu + A\Theta_\nu(\sigma^\nu - X_\nu^{-1}\phi^\nu),
\end{aligned}$$

and solve for the dual step direction $\Delta y^\nu$:

$$(A\Theta A^{\mathrm{T}})\Delta y^\nu = \psi^\nu. \tag{8.58}$$

Compute the primal step direction $\Delta x^\nu$:

$$\Delta x^\nu = -\,\Theta_\nu(\sigma^\nu - X_\nu^{-1}\phi^\nu - A^{\mathrm{T}}\Delta y^\nu), \tag{8.59}$$

and the slack step direction $\Delta z^\nu$ from

$$\Delta z^\nu = X_\nu^{-1}(\phi^\nu - Z_\nu \Delta x^\nu). \tag{8.60}$$

Update:

$$x^{\nu+1} = x^\nu + \alpha_{Q,\nu}\Delta x^\nu, \tag{8.61}$$

$$y^{\nu+1} = y^\nu + \alpha_{Q,\nu}\Delta y^\nu, \tag{8.62}$$

$$z^{\nu+1} = z^\nu + \alpha_{Q,\nu}\Delta z^\nu, \tag{8.63}$$

where $\alpha_{Q,\nu}$ is a step length restricted to the interval $(0,1]$. Reduce $\mu_\nu$ to $\mu_{\nu+1}$, update $\nu \leftarrow \nu+1$, and repeat the iterative step.

There are two differences between Algorithm 8.1.1 for linear programming and the quadratic programming Algorithm 8.2.1: one is in the definition of the matrix $\Theta$; the other is in the definition of the step length. In particular, for linear programming it is possible to estimate different step lengths ($\alpha_{P,\nu}$, and $\alpha_{D,\nu}$) that will preserve positivity of the primal and dual variables, respectively, because the primal problem is a problem only in $x$, while the dual problem is a problem only in $y, z$. For the quadratic programming problem, however, the primal variables appear in the constraints of the dual problem. Hence the step length that will preserve positivity of both primal and dual variables is computed as $\alpha_{Q,\nu} = \min(\alpha_{P,\nu}, \alpha_{D,\nu})$, where $\alpha_{P,\nu}$ and $\alpha_{D,\nu}$ are determined as in Section 8.1.1.

In the next section we discuss the parallel solution of systems of the form given by (8.27) and (8.58) when the constraint matrix $A$ is structured and when the objective function is either linear or separable quadratic (i.e., when $\Theta$ is diagonal). The solutions of these systems are the computationally intensive parts of Algorithms 8.1.1 and 8.2.1. They usually account for more than 90 percent of the execution time of the algorithm. The calculations of $\Delta x^\nu$ and $\Delta z^\nu$ involve simple matrix-vector products and vector additions, and can be implemented efficiently on parallel machines; they are also well suited for vector processing.

## 8.3 Parallel Matrix Factorization Procedures for the Interior Point Algorithm

We consider now the solution of the linear systems (8.27) or (8.58) for structured problems. (To simplify the presentation we discuss the solution of these systems during a generic iteration of the interior point algorithm and, therefore, we drop the iteration index $\nu$.) In particular, we assume that the constraint matrix $A$ has the following block-angular structure:

$$A = \begin{pmatrix} A_0 & & & \\ T_1 & W_1 & & \\ \vdots & & \ddots & \\ T_N & & & W_N \end{pmatrix}. \tag{8.64}$$

The matrix $A_0$ is of dimension $m_0 \times n_0$. The matrices $W_l$ are of dimensions $m_l \times n_l$, for $l = 1, 2, \ldots, N$. The matrices $T_l$ are of dimensions $m_l \times n_0$. It is assumed that $A_0$ and $W_l$ have full row rank, with $m_l \leq n_l$ for all $l = 1, 2, \ldots, N$, and that $n_0 \leq \sum_{l=1}^N n_l$.

The block-angular structure arises in the deterministic equivalent formulation of two-stage stochastic programs with recourse (see Section 13.3.3). This is also the constraint matrix of the dual formulation of the multicommodity network flow problem of Section 12.2.2. Stochastic programming problems and multicommodity network flow problems encountered in practical applications are of extremely large size. Hence, matrices of this form occasionally have millions of columns and hundreds of thousands of rows. Furthermore, while these matrices are sparse—as illustrated above—the product matrix $A\Theta A^{\mathrm{T}}$ could be completely dense due to the presence of the matrices $T_l$, $l = 1, 2, \ldots, N$.

This section develops a matrix factorization procedure for solving the linear systems of equations for the calculation of $\Delta y$, that exploits the special structure of the constraint matrix $A$. The vectors $\Delta y$ and $\psi$ in equations (8.27) or (8.58) are written as concatenations of subvectors, i.e., as $((\Delta y^0)^{\mathrm{T}}, (\Delta y^1)^{\mathrm{T}}, \ldots, (\Delta y^N)^{\mathrm{T}})^{\mathrm{T}}$ and $((\psi^0)^{\mathrm{T}}, (\psi^1)^{\mathrm{T}}, \ldots, (\psi^N)^{\mathrm{T}})^{\mathrm{T}}$, respectively, with $\Delta y^l$, $\psi^l \in \mathbb{R}^{m_l}$ for $l = 0, 1, \ldots, N$. $\Theta$ is assumed to be diagonal, as is the case for linear programs and for separable quadratic programs.

### 8.3.1 The matrix factorization procedure for the dual step direction calculation

The procedure for solving (8.27) or (8.58) is based on the use of a generalized version of the Sherman-Morrison-Woodbury formula. This formula is stated in the following lemma, where $I$ denotes the identity matrix.

**Lemma 8.3.1** *For matrices $A \in \mathbb{R}^{n\times n}$, $U \in \mathbb{R}^{n\times K}$, and $V \in \mathbb{R}^{n\times K}$ such that both $A$ and $(I + V^{\mathrm{T}}A^{-1}U)$ are invertible,*

$$(A + UV^{\mathrm{T}})^{-1} = A^{-1} - A^{-1}U(I + V^{\mathrm{T}}A^{-1}U)^{-1}V^{\mathrm{T}}A^{-1}.$$

**Proof** We verify the formula by multiplying its two sides by $A + UV^{\mathrm{T}}$:

$$\begin{aligned}
&(A + UV^{\mathrm{T}})^{-1}(A + UV^{\mathrm{T}}) \\
&= I + A^{-1}UV^{\mathrm{T}} - A^{-1}U(I + V^{\mathrm{T}}A^{-1}U)^{-1}V^{\mathrm{T}} \\
&\quad -A^{-1}U(I + V^{\mathrm{T}}A^{-1}U)^{-1}V^{\mathrm{T}}A^{-1}UV^{\mathrm{T}} \\
&= I + A^{-1}UV^{\mathrm{T}} - A^{-1}U(I + V^{\mathrm{T}}A^{-1}U)^{-1}(I + V^{\mathrm{T}}A^{-1}U)V^{\mathrm{T}} \\
&= I + A^{-1}UV^{\mathrm{T}} - A^{-1}UV^{\mathrm{T}} = I.
\end{aligned}$$

Having obtained the identity matrix completes the proof. ∎

The following theorem provides the foundation for the matrix factorization routine that solves the systems (8.27) and (8.58) for $\Delta y$.

**Theorem 8.3.1** *Let $M \doteq A\Theta A^{\mathrm{T}}$, where $\Theta$ is a diagonal matrix, and let $S \doteq diag\,(S_0, S_1, S_2, \ldots, S_N)$, $S_l \doteq W_l\Theta_l W_l^{\mathrm{T}} \in \mathbb{R}^{m_l\times m_l}, l = 1, 2, \ldots, N$, where the matrix $S_0 \doteq I$ is an $m_0 \times m_0$ identity matrix, and for $l =$*

*$0, 1, 2, \ldots, N$, $\Theta_l \in \mathbb{R}^{n_l \times n_l}$ is the diagonal submatrix of $\Theta$ corresponding to the $l$th block. Also, denoting $(\Theta_0^{-1})^2$ by $\Theta_0^{-2}$, let*

$$G_1 \doteq \Theta_0^{-2} + A_0^{\mathrm{T}} A_0 + \sum_{l=1}^{N} T_l^{\mathrm{T}} S_l^{-1} T_l, \tag{8.65}$$

$$G \doteq \begin{pmatrix} G_1 & A_0^{\mathrm{T}} \\ -A_0 & 0 \end{pmatrix}, U \doteq \begin{pmatrix} A_0 & I \\ T_1 & 0 \\ \vdots & \vdots \\ T_N & 0 \end{pmatrix}, V \doteq \begin{pmatrix} A_0 & -I \\ T_1 & 0 \\ \vdots & \vdots \\ T_N & 0 \end{pmatrix}.$$

*If $A_0$ and $W_l$, $l = 1, 2, \ldots, N$, have full row rank then the matrices $M$ and $G_2 \doteq -A_0 G_1^{-1} A_0^{\mathrm{T}}$ are invertible, and*

$$M^{-1} = S^{-1} - S^{-1} U G^{-1} V^{\mathrm{T}} S^{-1}. \tag{8.66}$$

**Proof** The matrix $\Theta$ is invertible for both the linear and the quadratic programs. By assumption $W_l$ has full row rank and hence $S_l$ is invertible for $l = 1, 2, \ldots, N$. It follows that $S$ is also invertible. Let

$$\Phi \doteq \begin{pmatrix} \Theta_0 & 0 \\ 0 & I \end{pmatrix},$$

$\overline{U} \doteq U\Phi$ and $\overline{V} \doteq V\Phi$. Then $M$ can be written as $M = S + \overline{U}\,\overline{V}^{\mathrm{T}}$.

In order to employ Lemma 8.3.1 to calculate the inverse of $M$, the matrix $I + \overline{V}^{\mathrm{T}} S^{-1} \overline{U}$ must be invertible. We prove that this is indeed the case. $I + \overline{V}^{\mathrm{T}} S^{-1} \overline{U} = \Phi(\Phi^{-2} + V^{\mathrm{T}} S^{-1} U)\Phi$ and hence it is invertible if and only if $\Phi$ and $(\Phi^{-2} + V^{\mathrm{T}} S^{-1} U)$ are invertible. $\Phi$ is invertible, so we check if $(\Phi^{-2} + V^{\mathrm{T}} S^{-1} U)$ is invertible. We can write

$$\begin{aligned} G &= \begin{pmatrix} \Theta_0^{-2} + A_0^{\mathrm{T}} A_0 + \sum_{l=1}^{N} T_l^{\mathrm{T}} S_l^{-1} T_l & A_0^{\mathrm{T}} \\ -A_0 & I - I \end{pmatrix} \\ &= \begin{pmatrix} \Theta_0 & 0 \\ 0 & I \end{pmatrix}^{-2} + \begin{pmatrix} A_0^{\mathrm{T}} A_0 + \sum_{l=1}^{N} T_l^{\mathrm{T}} S_l^{-1} T_l & A_0^{\mathrm{T}} \\ -A_0 & -I \end{pmatrix} \\ &= \Phi^{-2} + V^{\mathrm{T}} S^{-1} U. \end{aligned}$$

Therefore, $I + \overline{V}^{\mathrm{T}} S^{-1} \overline{U}$ is invertible if and only if $G$ is invertible.

By construction $\Theta_0^{-2}$ and $A_0^{\mathrm{T}} A_0$ are positive definite and symmetric. $T_l^{\mathrm{T}} S_l^{-1} T_l$ is positive definite and symmetric for all $l = 1, 2, \ldots, N$, since $S_l$ are positive definite. Hence $G_1$ is positive definite and symmetric, being the sum of positive definite symmetric matrices, and it is invertible. Its inverse matrix $G_1^{-1}$ is symmetric, has rank $n_0$, and can be written as $G_1^{-1} = G_1^{-\frac{1}{2}} G_1^{-\frac{1}{2}}$ where $G_1^{-\frac{1}{2}}$ is also symmetric. By assumption $A_0$ has full row rank, hence $A_0 G_1^{-\frac{1}{2}}$ has full row rank and $G_2 = -A_0 G_1^{-1} A_0^{\mathrm{T}}$ is invertible. The rank of $G$ is $m_0 + n_0$ so that $G$ is invertible. Hence $I + \overline{V}^{\mathrm{T}} S^{-1} \overline{U}$ is invertible and Lemma 8.3.1 can be applied to invert $M$:

$$M^{-1} = (S + \overline{U}\,\overline{V}^{\mathrm{T}})^{-1} = S^{-1} - S^{-1} U G^{-1} V^{\mathrm{T}} S^{-1}.$$

This completes the proof. ■

It is easy to verify, using equation (8.66), that the solution of the linear system $(A\Theta A^{\mathrm{T}})\Delta y = \psi$ is given by $\Delta y = p - r$ where $p$ solves the system $Sp = \psi$, and $r$ is obtained from the system

$$Gq = V^{\mathrm{T}} p \text{ and } Sr = Uq. \tag{8.67}$$

The vector $p$ can be computed componentwise by solving $S_l p^l = \psi^l$, for $l = 1, 2, \ldots, N$. In order to solve for $q$ we exploit the block structure of $G$ and write:

$$Gq = \begin{pmatrix} G_1 & A_0^{\mathrm{T}} \\ -A_0 & 0 \end{pmatrix} \begin{pmatrix} q^1 \\ q^2 \end{pmatrix} = \begin{pmatrix} \hat{p}^1 \\ \hat{p}^2 \end{pmatrix} \text{ where } \begin{pmatrix} \hat{p}^1 \\ \hat{p}^2 \end{pmatrix} \doteq V^{\mathrm{T}} p. \tag{8.68}$$

Hence, we get

$$q^2 = -G_2^{-1}(\hat{p}^2 + A_0 G_1^{-1} \hat{p}^1), \tag{8.69}$$

$$q^1 = G_1^{-1}(\hat{p}^1 - A_0^{\mathrm{T}} q^2). \tag{8.70}$$

Once $q$ is known, $r$ can be computed componentwise from (8.67). The procedure for calculating $\Delta y$ can be summarized as follows. (We adopt here the notation $A_{i\cdot}$ for the $i$th row of $A$, and $A_{\cdot j}$ for the $j$th column of $A$.)

#### Procedure 8.3.1 Matrix factorization for dual step direction calculation

**Step 1:** *Solve* $Sp = \psi$.

**Step 2:** *(Solve* $Gq = V^{\mathrm{T}} p$.*)*

- **a.** *For all* $l = 1, 2, \ldots, N$, *solve* $S_l (u^l)^i = (T_l)_{\cdot i}$, *for* $(u^l)^i$, $i = 1, 2, \ldots, n_0$ *thus computing the columns of the matrix* $S_l^{-1} T_l$.
- **b.** *For all* $l = 1, 2, \ldots, N$ *multiply* $T_l^{\mathrm{T}} (u^l)^i$, *for* $i = 1, 2, \ldots, n_0$ *to form* $T_l^{\mathrm{T}} S_l^{-1} T_l$. *Form* $G_1$ *(eq. (8.65)). Compute* $\hat{p}^1, \hat{p}^2$ *(eq. (8.68)).*

**c.** *Solve $G_1 u = \hat{p}^1$ for $u$ and set $v = \hat{p}^2 + A_0 u$ (eq. (8.69)).*

**d.** *Form $G_2$ by solving $G_1 w^i = (A_0^{\mathrm{T}})_{\cdot i}$ for $w^i$, for $i = 1, 2, \ldots, m_0$, and setting $G_2 = -A_0[w^1\, w^2 \cdots w^{m_0}]$.*

**e.** *Solve $G_2 q^2 = -v$ for $q^2$, and solve $G_1 q^1 = \hat{p}^1 - A_0^{\mathrm{T}} q^2$ for $q^1$, (eq. (8.70)).*

**Step 3:** *(Solve $Sr = Uq$.) Set $r^0 = A_0 q^1 + q^2$ and for all indices $l = 1, 2, \ldots, N$ solve $S_l r^l = T_l q^1$ for $r^l$.*

**Step 4:** *(Form $\Delta y$.) For $l = 1, 2, \ldots, N$, set $\Delta y^l = p^l - r^l$.*

## Decompositions for Parallel Computing

Procedure 8.3.1 is well suited for parallel implementation because it relies largely on matrix (sub)block computations that can be performed independent of one another. The computation begins with the submatrices $T_l$, $W_l$ and $\Theta_l$ and the vector segment $\psi^l$ located at the $l$th processor. Processor $l$ can compute $S_l$ and proceed independently with all computations involving only local data.

Interprocessor data communication is necessary at only three instances in the parallel algorithm. After forming the terms $T_l^{\mathrm{T}} S_l^{-1} T_l$ in Steps 2a–2b of Procedure 8.3.1 in parallel, the processors communicate to form the matrix $G_1$ and the vectors $\hat{p}^1$ and $\hat{p}^2$ in Step 2b. The results can be accumulated at a single processor, designated as the *master.* The computations involving the dense matrix $G_1$ can be done serially on the master processor. Steps 2d–2e, which involve the dense matrix $G_2$, are also done serially on the master processor.

The master processor must then broadcast the computed vector $q$ to all other nodes. The remaining Steps 3 and 4 require only the distributed data $S_l$, $T_l$, and $p^l$ on the $l$th processor and may be carried out with full parallelism. A final communication step accumulates the partial vectors $\Delta y^l$ at the master processor. This vector can then be made available to all processors for use in the calculation of the directions $\Delta x$ and $\Delta z$ of an interior point algorithm.

An alternative parallel implementation of Procedure 8.3.1 would distribute the matrices $G_1$ and $G_2$ to all processors and let the processors proceed locally (and redundantly) with all calculations involving these matrices. The approach based on a master processor described above uses an *all-to-one* communication step to accumulate the dense matrices at the master, followed by a *one-to-all* communication step to distribute the results to the processors. The alternative approach suggested here combines these two communication steps into a single *all-to-all* communication step that distributes the dense matrices to all processors. Both of these alternatives are very efficient on present day distributed memory machines with high-bandwidth communication networks (see the results in Section 15.5.1).

Yet another alternative approach is to distribute the dense matrices across processors and use parallel dense linear algebra techniques for all calculations that involve these matrices. This approach is more suitable to shared-memory, tightly coupled multiprocessors, or when the dense matrices $G_1$ and $G_2$ are large. The parallel procedure that uses the master processor is summarized as follows.

**Procedure 8.3.2 Parallel matrix factorization for dual step direction calculation**

*Begin with the following data distribution: Processor $l$ holds $S_l$, $T_l$, and $\psi^l$, $l = 1, 2, \ldots, N$. A designated master processor also holds $A_0$, $S_0$, $\Theta_0$, and $\psi^0$.*

**Step 1:** *(Solve $Sp = \psi$ in parallel.) Solve $S_0 p^0 = \psi^0$ on the master processor. Solve in parallel, on processors $l = 1, 2, \ldots, N$, the system $S_l p^l = \psi^l$ for $p^l$.*

**Step 2:** *(Solve $Gq = V^{\mathrm{T}} p$ in parallel.)*

- **a.** *Solve in parallel, on processors $l = 1, 2, \ldots, N$, the linear systems $S_l (u^l)^i = (T_l)_{\cdot i}$ for $(u^l)^i$, $i = 1, 2, \ldots, n_0$.*
- **b.** *Multiply $T_l^{\mathrm{T}} (u^l)^i$ for $i = 1, 2, \ldots, n_0$ in parallel, on processors $l = 1, 2, \ldots, N$. Communicate to form $G_1$ and $\hat{p}^1$ on the master processor and compute $\hat{p}^2$ serially on the master processor.*
- **c.** *Solve $G_1 u = \hat{p}^1$ for $u$ serially on the master processor, and set $v = \hat{p}^2 + A_0 u$.*
- **d.** *Form $G_2$ by solving serially on the master processor the linear systems $(G_1) w^i = (A_0^{\mathrm{T}})_{\cdot i}$ for $w^i$, for $i = 1, 2, \ldots, m_0$, and by setting $G_2 = -A_0 [w^1\, w^2 \cdots w^{m_0}]$.*
- **e.** *Solve serially, on the master processor, the system $G_2 q^2 = -v$ for $q^2$ and the system $G_1 q^1 = \hat{p}^1 - A_0^{\mathrm{T}} q^2$ for $q^1$. Communicate to distribute $q^1$ to all processors.*

**Step 3:** *(Solve $Sr = Uq$ in parallel.) Set $r^0 = A_0 q^1 + q^2$ on the master processor. Solve the systems $S_l r^l = T_l q^1$ for $r^l$ in parallel, on all processors $l = 1, 2, \ldots, N$.*

**Step 4:** *(Form $\Delta y$ in parallel.) Set $\Delta y^0 = p^0 - r^0$ on the master processor. Set $\Delta y^l = p^l - r^l$ on processors $l = 1, 2, \ldots, N$. Communicate to gather the vector $\Delta y$ on the master processor.*

Sections 14.6 and 15.5 give details on the implementation of this procedure on different parallel architectures and report numerical results.

## 8.4 Notes and References

The seminal contribution in interior point methods was made by Karmarkar (1984), where the projective scaling algorithm was proposed and its polynomial complexity was established. His paper and claims about superior computational performance of this algorithm over the classic simplex algorithm for large-scale problems, fueled the extensive research that followed. It is now recognized, though, that the closely related affine scaling algorithm was proposed earlier by Dikin (1967, 1974). For a general introduction see Goldfarb and Todd (1989), Gonzaga (1992), Kojima et al. (1991), and Marsten et al. (1990).

The linear programming software system OB1 (Optimization with Barrier 1) is described in Lustig, Marsten, and Shanno (1991). OBN is the quadratic programming optimization version of OB1. ALPO (Another Linear Programming Optimizer) is described in Vanderbei (1993), and LOQO (Linear or Quadratic Optimizer) is described in Vanderbei (1992) and Vanderbei and Carpenter (1993).

**8.1** The primal-dual path following algorithm for linear programs was developed, independently, by Monteiro and Adler (1989a) and Kojima, Mizuno, and Yoshise (1989). Earlier work on the primal-dual central path was developed by Meggido (1989). Numerical issues for the efficient implementation of the algorithm were resolved by McShane, Monma, and Shanno (1989), Lustig, Marsten, and Shanno (1991) and Vanderbei (1993). See also Choi, Monma, and Shanno (1990) and Mehrotra (1992).

**8.2** The primal-dual path following algorithm for quadratic programming problems is described in Monteiro and Adler (1989b, 1990), Vanderbei and Carpenter (1993), and Carpenter et al. (1993).

**8.3** For a general discussion of direct matrix inversion methods see Householder (1975). Parallel implementation of interior point algorithms for dense problems is discussed in Qi and Zenios (1994) and Eckstein et al. (1992). The exploitation of matrix structure in computing the dual step direction for interior point algorithms has been addressed in Loute and Vial (1992), Choi and Goldfarb (1993), Czyzyk, Fourer, and Mehrotra (1995), and De Silva and Abramson (1996). The use of the Sherman-Morrison-Woodbury updating formula for computing dual steps for stochastic programs was suggested by Birge and Qi (1988), and was further extended by Birge and Holmes (1992) who also performed extensive numerical experiments with this method. The extension to multistage stochastic programming problems is discussed in Holmes (1993). Its use for solving multicommodity network flow problems has been suggested by Choi and Goldfarb (1990, 1993) and Lustig and Li (1992). The parallel matrix factorization procedure and its implementation on hypercubes and other parallel machines, was developed by Jessup, Yang, and Zenios (1994a). Yang and Zenios (1996) discuss the application of the parallel matrix factorization procedure within an interior point algorithm for the solution of large-scale stochastic programs. Ruszczyński (1993) provides a survey on the use of interior point methods for solving stochastic programming problems, and Vladimirou and Zenios (1996) provide a survey of parallel methods for stochastic programming, including a discussion on the parallelization of interior point methods.

# Part III

# APPLICATIONS

> The sciences do not try to explain, they hardly even try to interpret, they mainly make models. By a model is meant a mathematical construct which, with the addition of certain verbal interpretations, describes observed phenomena. The justification of such a mathematical construct is solely and precisely that it is expected to work.
>
> John Von Neumann

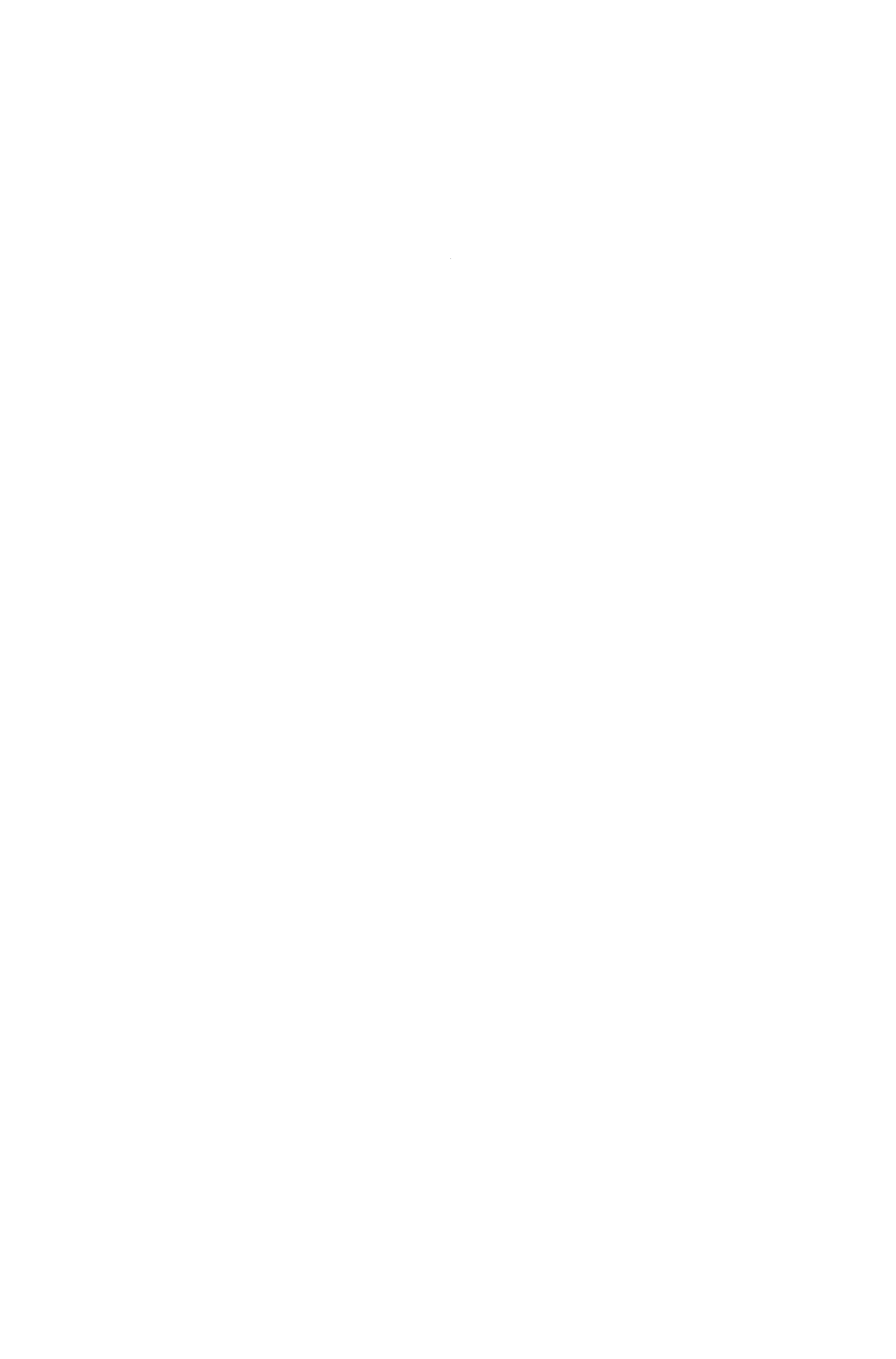

CHAPTER 9

# Matrix Estimation Problems

The problem of data inconsistency and measurement errors in the economic sciences has traditionally been a stumbling block in the development of formal economic models. Oscar Morgenstern devoted his 1963 book to the problem of inconsistency in economic measurements. The Cambridge Growth Project under the direction of Sir Richard Stone addressed (among other issues) this problem with particular focus on the estimation of tables of economic data that satisfy internal consistency conditions. One important outcome of the project was the book by Bacharach (1970) on the estimation of input-output tables.

The problem of adjusting the entries of a large matrix to satisfy consistency requirements occurs frequently in economics as well as in urban planning, statistics, demography, and stochastic modeling. It is typically posed as follows:

> Given a rectangular matrix $A$, determine a matrix $X$ that is *close* to $A$ and satisfies a given set of linear restrictions on its entries.

A well-studied instance of this problem—occurring in transportation planning and input-output analysis—requires that $A$ be adjusted so that the row and column totals equal fixed positive values. A related problem occurring in developmental economics requires that the row and column totals (of a square matrix) be equal, but not necessarily to prespecified values.

These problems are known as *balancing* the matrix $A$. The terms *matrix estimation* or *matrix adjustment* have also been used to describe the problem of calculating a balanced matrix. Because we discuss several balancing problems with different consistency requirements the definition of a balanced matrix is problem dependent. That is, a matrix is defined to be *balanced* if it satisfies the given set of linear restrictions of the problem.

This chapter addresses problems arising in the estimation of matrices. Section 9.1 discusses applications of this problem in diverse real-world settings. Formal mathematical models are given in Section 9.2, and Section 9.3 develops iterative row-action algorithms for solving these models. Notes and references in Section 9.4 conclude this chapter.

## 9.1 Applications of Matrix Balancing

In this section we discuss several applications wherein problems of estimation of balanced matrices arise. These applications are drawn from economics, transportation and telecommunications, statistics, demography and stochastic modeling.

### 9.1.1 Economics: social accounting matrices (SAMs)

A social accounting matrix (SAM) is a square matrix $A$ whose entries represent the flow of funds among the national income accounts of a country's economy at a fixed point in time. Each index of a row or of a column in $A$ represents an account or agent in the economy. Entry $a_{ij}$ is positive if agent $j$ receives funds from agent $i$. Thus a SAM is a snapshot of the critical variables in a general equilibrium model describing the circular flow of financial transactions in an economy. For balancing problems arising from estimating SAMs, the balance conditions are the a priori accounting identities that each agent's total expenditures must equal his total receipts. That is, for each index $i$, the sum of the entries in row $i$ of the matrix $A$ must equal the sum of the entries in column $i$.

The agents of an economy include institutions, factors of production, households, and the rest of the world (to account for transactions with the economies of other countries). Figure 9.1 shows the major relationships among accounts in a simplified SAM. Briefly, the production activities generate added value that flows to the factors of production—land, labor, and capital. Factor income is the primary source of income for institutions—households, government, and firms—that purchase goods and services supplied by productive activities, thereby completing the cycle. This simple model produces a SAM with three agents and three nonzero entries. Of course, to be useful for equilibrium modeling, this highly aggregated model must be disaggregated into subaccounts for each sector of the economy. Estimation of disaggregated SAMs with more than 250 agents has been reported in the literature.

SAMs are used as the database for complex economy-wide general equilibrium models. The entries of the SAM provide a convenient database used to estimate the parameters of an equilibrium model. Researchers at the World Bank have developed specialized modeling systems based on SAM databases; similar models have been developed by the United Nations Statistical Office, the Cambridge Growth Project, and the statistical bureaus of developing and industrialized countries.

Inconsistent data is an inherent problem when statistical methods are used to estimate underlying economic models. In particular, the direct estimate of a SAM is never balanced. The following quote of Sir Richard Stone summarizes the sources of inconsistency in SAM modeling:

> ...it is impossible to establish by direct estimation a system of national accounts free of statistical discrepancies, residual

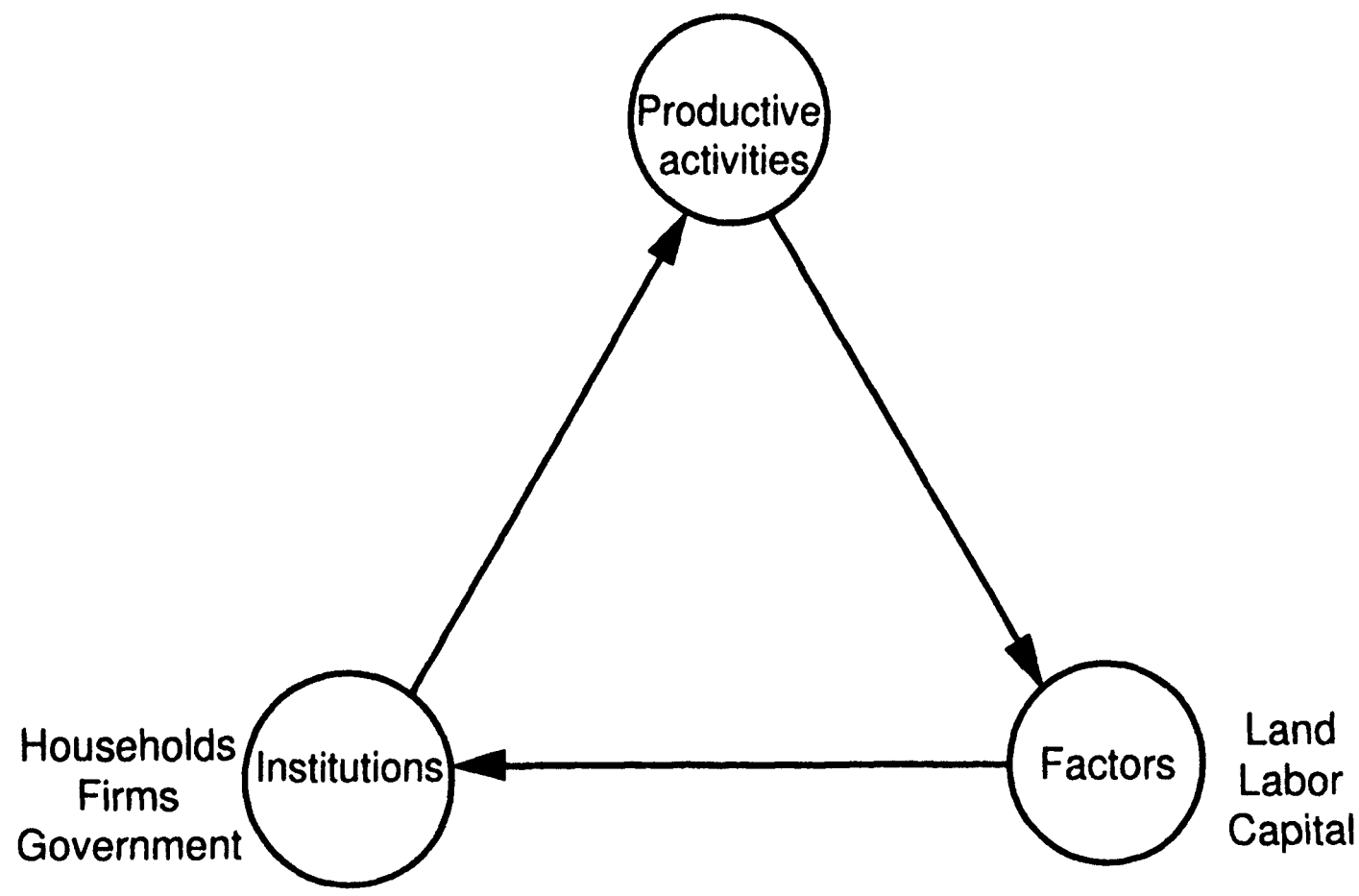

**Figure 9.1** A simplified Social Accounting Matrix (SAM).

> errors, unidentified items, balancing entries and the like since the information available is in some degree incomplete, inconsistent and unreliable. Accordingly, the task of measurement is not finished when the initial estimates have been made and remains incomplete until final estimates have been obtained which satisfy the constraints that hold between their true values. (Van der Ploeg 1982, p.186)

Therefore, the raw estimates of a SAM must be adjusted so that the consistency requirements are satisfied. The problem of adjusting the initial matrix so that the row and column sums are equal is an example of matrix balancing. This problem motivates much of the work on matrix balancing for economic modeling.

A matrix balancing problem also arises when *partial survey* methods are used to estimate a SAM. Frequently estimates of the total expenditures and receipts for each agent in an economy are available but current data for the individual transactions between the agents are not. If a complete and balanced SAM is available from an earlier period it must be updated to reflect the recent index totals. The problem is then to adjust the entries of the old matrix $A$ so that the row and column totals equal the given fixed amounts. A similar balancing problem occurs when the entries of an input-output matrix must be updated to be consistent with exogenous estimates of the total levels of primary inputs and final demands.

### 9.1.2 Transportation: estimating Origin-Destination Matrices

Urban transportation analysis often models traffic flow over physical transportation systems involving, for example, urban highways, subways, buses, or airplanes. The transportation system is represented as a *network* with nodes corresponding to intersections and arcs corresponding to travel links. (Refer to Section 12.1 for the definition of a network and further discussion of transportation models.) Each arc has an associated *travel cost function* representing the disutility of travel time as a function of the total flow along that arc. A user entering the transportation network must choose a path from the origin node to the destination node. The matrix representing the rate of traffic flow between all pairs of nodes is called the *origin-destination matrix*. Each entry of the matrix indicates the fraction of the total network traffic that travels between an origin-destination pair.

We refer here to data estimation problems arising in the transportation of physical quantities over a network, such as those found in a highway or railroad system. However similar problems arise in the transport of messages over a telecommunication network or the transport of signals over computer networks.

The traffic assignment problem requires determination of the flow rate along each arc and the travel time between each origin-destination pair given the network representing the system, the travel cost functions, and the origin-destination matrix. For example, under the assumption that users choose the path with the lowest total travel time, the system is in equilibrium if no motorist can lower his or her travel time by altering his route. The resulting steady-state distribution of traffic flows is called a *user equilibrium*. We assume that the model represents the network at, for example, a peak travel time and that the number of trips is fixed. That is, the total travel demand is represented by a fixed $n \times m$ origin-destination matrix $A$, where $a_{ij}$ is the rate of flow from origin zone $i$ to destination zone $j$.

The origin-destination matrix $A$ is a required input for a wide variety of transportation planning models. Matrix balancing problems arise when origin-destination matrices are estimated from observed traffic flows on selected sets of arcs. For example, data might be collected measuring the total flow out of each origin node and into each destination node. When the estimated origin-destination matrix is used as an input into the underlying traffic assignment model, it should ensure that the traffic flows produced as outputs will equal the observed traffic flows. Since there may be many such matrices, the estimated matrix is chosen by minimizing some properly defined distance between $A$ and a fixed reference origin-destination matrix $A_0$, which could be an estimate based on a previous study.

In practice the total flow out of each origin and into each destination is frequently used for observed traffic flows. If there is a single route connecting each origin and destination, then the resulting problem of estimating an origin-destination matrix can be formulated as a matrix balancing problem.

### 9.1.3 Statistics: estimating contingency tables

Contingency tables are used in the classification of items by several criteria, each of which is partitioned into a finite number of categories. For example, suppose that $N$ individuals in a population are classified according to the criteria marital status and age, which are divided into $m$ and $n$ categories, respectively. The $mn$ distinct classifications are called *cells* and the resulting $m \times n$ contingency table contains the number of individuals in each cell.

For many applications it is important to estimate the underlying cell probabilities, i.e., the fraction of the population $N$ that is classified in each cell. Determining the cell probabilities directly by population sampling is rarely possible because in general such a procedure is prohibitively expensive. Thus cell probabilities must be derived from partial observations. A balancing problem arises when prior values for the cell probabilities are revised so that they become consistent with known information. For example, suppose that the marginal distribution of each criterion within the population is known—i.e., the total number of people in each one of the $n$ age categories is known, together with the total number of people in each one of the $m$ categories indicating marital status—and that prior cell values are available. The prior values could be obtained from an out-of-date general census, or from another population with similar ethnic and socioeconomic characteristics, or from a contingency table derived from sampling. The prior cell values must then be adjusted so that they are consistent with the known marginal probabilities.

Similar models are extended to four-way tables—tables whereby each cell is characterized by four characteristics (e.g., marital status, age, income, and education). In the estimation of four-way tables different marginal totals may be available (e.g., one-, two- and three-way marginal totals can be given).

### 9.1.4 Demography: modeling interregional migration

Estimating interregional migration flows based on partial and outdated information is a recurrent problem in demography that arises when figures are available for net in- and out-migration from every region, and an estimate is needed of the interregional migration patterns. In the United States, for example, flow matrices with detailed migration characteristics become available once every ten years from the general census of the population. In the interim, however, net migration estimates for every region are available as byproducts of annual population estimates. (Net migration is the difference between total population change and changes from births and deaths). An updated migration matrix that reconciles the out-of-date migration patterns with the more recent net figures is therefore required.

Demographers have postulated several models for estimating interregional migration including gravity models, Markov or fixed transition prob-ability models, and doubly constrained *minimum-information* models.

The doubly constrained minimum-information model leads to a matrix balancing problem, defined as follows. Let $M$ denote the total population in all regions, and for each $i = 1, 2, \ldots, n$ let $O_i$ and $I_i$ be, respectively, the net out-migration and in-migration for region $i$. Let $x_{ij}$, $i, j = 1, 2, \ldots, n$, $i \neq j$ be the model-estimated probabilities that any individual in the system is a migrant from region $i$ to region $j$, and let $a_{ij}$ be some prior estimate of the probabilities. The minimum-information model can be written as:

$$\underset{x}{\text{Minimize}} \quad \sum_{i=1}^{n} \sum_{\substack{j=1 \\ j \neq i}}^{n} x_{ij} \left( \log \left( \frac{x_{ij}}{a_{ij}} \right) - 1 \right) \tag{9.1}$$

$$\text{s.t.} \quad \sum_{\substack{j=1 \\ j \neq i}}^{n} x_{ij} = O_i / M, \qquad \text{for } i = 1, 2, \ldots, n, \tag{9.2}$$

$$\sum_{\substack{i=1 \\ i \neq j}}^{n} x_{ij} = I_j / M, \qquad \text{for } j = 1, 2, \ldots, n. \tag{9.3}$$

An additional constraint is often imposed to fix the total distance $D$ traveled by all migrants. If $d_{ij}$ is the distance between regions $i$ and $j$, this additional constraint is $\sum_{i=1}^{n} \sum_{\substack{j=1 \\ j \neq i}}^{n} d_{ij} x_{ij} = D$.

### 9.1.5 Stochastic modeling: estimating transition probabilities

A problem that appears in the estimation of transition probabilities from macrodata is the following:

> Given the proportion of observations in $n$ alternative states during $T$ time periods, estimate the probability of transition between any two states at the next time period.

Problems of this sort typically appear in marketing research. Observations are available for the proportion of customers using brand $i$ where $i = 1, 2, \ldots, n$ during the time periods $t = 1, 2, \ldots, T$ when a market survey was conducted. Assuming that transition probabilities are constant over the time interval of the survey, we want to estimate the underlying transition probability matrix. Let $x_i^t$ be the probability of a customer using the $i$th brand during period $t$, and let $a_{ij}$ be the probability of transition from $i$ to $j$. Furthermore, let $x^t$ be the vector $x^t = (x_i^t)_{i=1}^n$ and let $A = (a_{ij})$ be an $n \times n$ matrix. The probability vector at time $t + 1$ can be obtained from that at time $t$ by applying the transpose of the transition probability matrix $A$ as follows:

$$x^{t+1} = A^{\mathrm{T}} x^t. \tag{9.4}$$

In practice, one is given the observations $x^t$, $t = 1, 2, \ldots, N$, and would like to compute the transition probability matrix $A$ that satisfies the relation (9.4). There are $n(n-1)$ probabilities $a_{ij}$ to be estimated and a solution exists if the number of observations is sufficiently small. In general $N > n$ and no solution exists. For any estimated matrix $A$ we may expect at most $N$ nonzero discrepancy vectors $x^{t+1} - A^{\mathrm{T}}x^t$. A quadratic programming model for estimating the transition probability matrix $A$ minimizes a quadratic form of the discrepancy vectors:

$$\sum_{t=1}^{T} \langle x^{t+1} - A^{\mathrm{T}}x^t,\ Q(x^{t+1} - A^{\mathrm{T}}x^t)\rangle, \tag{9.5}$$

where $Q$ is an appropriate symmetric positive semidefinite matrix. In addition we require that the sum of the entries in every row of the estimated matrix equal one (i.e., there is always a transition to some state) and that all entries are nonnegative. The matrix estimation problem that arises in this case can be modeled as the problem of minimizing (9.5) subject to the constraints

$$\sum_{j=1}^{n} a_{ij} = 1, \quad \text{for } i = 1, 2, \ldots, n, \quad \text{and} \tag{9.6}$$

$$a_{ij} \geq 0, \quad \text{for all } i, j. \tag{9.7}$$

The resulting problem is a simple version of a matrix balancing problem with constraints only on the row sums. The estimated matrix $A$ is not related to any prior estimate of the transition probabilities, instead, the model minimizes a quadratic term of discrepancy vectors that are related to the computed matrix $A$.

A typical real-world example of estimating transition probabilities is the following. Consider the market shares of new and used cars among car buyers over a period of five years. We want to determine the probability of a used-car buyer switching over to a new car and vice versa. One is therefore estimating the entries in a matrix of transition probabilities subject to the constraints that all entries are between 0 and 1 and that the sum of the probabilities over all states is equal to 1. Estimation of transition probabilities may also be required for a set of data on market shares for various other commodities, such as different cigarette or soft drink brands.

## 9.2 Mathematical Models for Matrix Balancing

Many different mathematical formulations for balancing problems have appeared in the literature. While some are equivalent, others are genuinely different due to the underlying application. In this section we review the diverse approaches to modeling matrix balancing problems and give some of the most popular mathematical formulations.

### 9.2.1 Matrix estimation formulations

We consider first a general formulation of the problem. Let $X = (x_{ij})$ be an $m \times n$ matrix and denote

$$\mathcal{A} \doteq \{(i,j) \mid i = 1, 2, \ldots, m, \text{ and } j = 1, 2, \ldots, n\}.$$

An *integral of the matrix* $X$ is a sum $\sum_{(i,j)\in T} x_{ij}$ where $T$ is a given subset of $\mathcal{A}$. The matrix balancing problem is to estimate a matrix $X$ such that some of its integrals will obey certain conditions and the matrix itself will be related in a specific way to another given matrix $A = (a_{ij})$.

We refer to such problems as *matrix estimation* or *matrix construction* problems. Matrix balancing is a special case of matrix estimation. Other related problems that can be viewed as special cases of matrix estimation appear in the literature under the terms matrix scaling (similarity scaling, equivalence scaling, truncated scaling), constrained matrix problems, generalized scaling, and fair-share matrix allocation.

Various conditions on the integrals of the matrix $X$ have been considered in the examples mentioned earlier in this chapter. These include (but are not restricted to): fixing the row and column sums to preassigned values; constraining the individual entries of $X$ to lie within specified lower and upper bounds; confining the row and column sums to lie within specified lower and upper bounds and fixing the sum of all the entries of $X$ to a preassigned value; forcing row sums to equal column sums.

The precise relationship between the constructed matrix $X$ and the given matrix $A$ is what differentiates the various instances of matrix estimation problems. These relationships are of three general types.

*The form relationship:* Here one imposes a certain form relation that dictates the desired form of $X$. One such relation is to demand the existence of vectors $\lambda \in \mathbb{R}^m$ and $\mu \in \mathbb{R}^n$ and a real $\delta \in \mathbb{R}$ such that for a particular $m \times n$ weight matrix $W$ the form relation

$$X = A + \text{diag}(\lambda)\, W + W \,\text{diag}(\mu) + \delta \,\text{diag}(\lambda)\, W \,\text{diag}(\mu) \qquad (9.8)$$

would hold. With the choice $W = A$, $\delta = 1$ and the substitution $r = \lambda + \mathbf{1}$ and $s = \mu + \mathbf{1}$ ($\mathbf{1}$ is the vector all of whose components equal 1), relation (9.8) becomes

$$X = \text{diag}(r)\, A \,\text{diag}(s). \qquad (9.9)$$

The problem of estimating a matrix $X$ which is related to a given matrix $A$ via this particular form relationship is known as the *equivalence scaling* problem.

Yet another form relation between $X$ and $A$ is to demand the existence of a positive vector $z \in \mathbb{R}^k$ (i.e., $z_j > 0,\ j = 1, 2, \ldots, k$) such that the lexicographic vectorial reorderings (see Section 9.3) $x \in \mathbb{R}^{mn}$ and $a \in \mathbb{R}^{mn}$ of $X$ and $A$, respectively, are related by

$$x_q = a_q \left( \prod_{p=1}^{k} z_p^{\phi_{pq}} \right), \quad q = 1, \ldots, mn, \tag{9.10}$$

where $\Phi = (\phi_{pq})$ is a $k \times mn$ coefficients matrix of a system of linear equality or linear inequality integral constraints on $X$.

*The axiomatic approach:* This approach specifies a list of axioms that the relationship between $X$ and $A$ should satisfy, and then establishes a specific form relation that is proven to fulfill the given axioms.

*Distance optimization:* In this approach the constructed matrix $X$ should be as close as possible to the original matrix $A$, subject to constraints on the integrals of the matrix. The notion *close* is defined by some distance function $f(X; A)$ that measures the *distance* between $X$ and $A$. The choice $f(X; A) \doteq \| X - A \|_F^2$, where $\| \cdot \|_F$ denotes the Frobenius norm, i.e., $\| A \|_F^2 \doteq \sum_{i=1}^m \sum_{j=1}^n | a_{ij} |^2$, leads to a linearly constrained quadratic optimization problem. Another commonly used objective is the negative entropy function (see also Example 2.1.2)

$$\sum_{(i,j) \in \mathcal{A}} x_{ij} \left( \log \left( \frac{x_{ij}}{a_{ij}} \right) - 1 \right). \tag{9.11}$$

It is important to point out that connections exist between some of the seemingly different approaches to matrix estimation, particularly between several form relations and the entropy optimization problem (see discussion in Section 9.4).

Most of the matrix balancing applications encountered in practice can be formulated as one of the following two matrix estimation problems, with simple matrix integrals. The first problem is:

**Problem 9.2.1** *Given an $m \times n$ nonnegative matrix $A = (a_{ij})$ and positive vectors $u \in \mathbb{R}^m$ and $v \in \mathbb{R}^n$, determine a "nearby" nonnegative matrix $X = (x_{ij})$ (of the same dimensions) such that*

$$\sum_{j=1}^{n} x_{ij} = u_i, \quad \text{for } i = 1, 2, \ldots, m, \tag{9.12}$$

$$\sum_{i=1}^{m} x_{ij} = v_j, \quad \text{for } j = 1, 2, \ldots, n, \tag{9.13}$$

*and $x_{ij} > 0$ only if $a_{ij} > 0$.*

The estimation of a social accounting matrix to agree with prespecified row and column sums is an example of this problem. The estimation of origin-destination matrices, contingency tables, and migration patterns (discussed in the previous section) are also examples of this problem.

The second problem is defined as follows:

**Problem 9.2.2** *Given an $n \times n$ nonnegative matrix $A = (a_{ij})$, determine a "nearby" nonnegative matrix $X = (x_{ij})$ (of the same dimensions) such that*

$$\sum_{j=1}^{n} x_{ij} = \sum_{j=1}^{n} x_{ji}, \quad \text{for } i = 1, \ldots, n, \tag{9.14}$$

*and $x_{ij} > 0$ only if $a_{ij} > 0$.*

The estimation of a social accounting matrix so that row sums equal column sums is an example of this problem.

In principle, there are infinitely many matrices satisfying the consistency restrictions (9.12)–(9.13) or (9.14). For both problems we need to define the notion of a nearby matrix, and different models and algorithmic approaches follow naturally from different definitions. In the next section we develop specific models based on the distance optimization approach. First we describe the network structure of the two basic matrix balancing Problems, 9.2.1 and 9.2.2 (see also Section 12.1) as this structure is instructive for understanding the algorithms.

### 9.2.2 Network structure of matrix balancing problems

The connection of matrix balancing problems to network models is established by associating with the matrix $A$ a directed graph $\mathcal{G} = (\mathcal{N}, \mathcal{A})$ with node set $\mathcal{N}$ and arc set $\mathcal{A}$. Each row $i$ and each column $j$ of the matrix is associated with a node; a directed arc $(i, j)$ exists between the nodes $i$ and $j$ associated with the $i$th row and the $j$th column respectively, when the matrix entry $a_{ij}$ is nonzero. The notation $\delta_i^+ \doteq \{j \mid (i, j) \in \mathcal{A}\}$ denotes the set of nodes having an incident arc with origin node $i$, and $\delta_j^- \doteq \{i \mid (i, j) \in \mathcal{A}\}$ denotes the set of nodes having an incident arc with destination node $j$.

We define first the graph for Problem 9.2.1.

**Definition 9.2.1 (Bipartite graph of matrix A)** *Given an $m \times n$ nonnegative matrix $A = (a_{ij})$, the graph $\mathcal{G} = (\mathcal{N}, \mathcal{A})$ is called the bipartite graph of $A$ if the set of arcs is $\mathcal{A} \doteq \{(i, j) \mid a_{ij} > 0\}$ and the set of nodes $\mathcal{N}$ can be partitioned into two subsets $\mathcal{N}_1$ and $\mathcal{N}_2$, of cardinality $m$ and $n$, respectively, so that for every arc $(i, j) \in \mathcal{A}$, $i \in \mathcal{N}_1$ and $j \in \mathcal{N}_2$.*

A graph in which the nodes are divided into two subsets with all arcs leading from nodes of one subset to those of the other is commonly referred to as *bipartite* (see an example in Figure 9.2). On bipartite graphs it is possible to define an optimization problem known as the *transporta-*

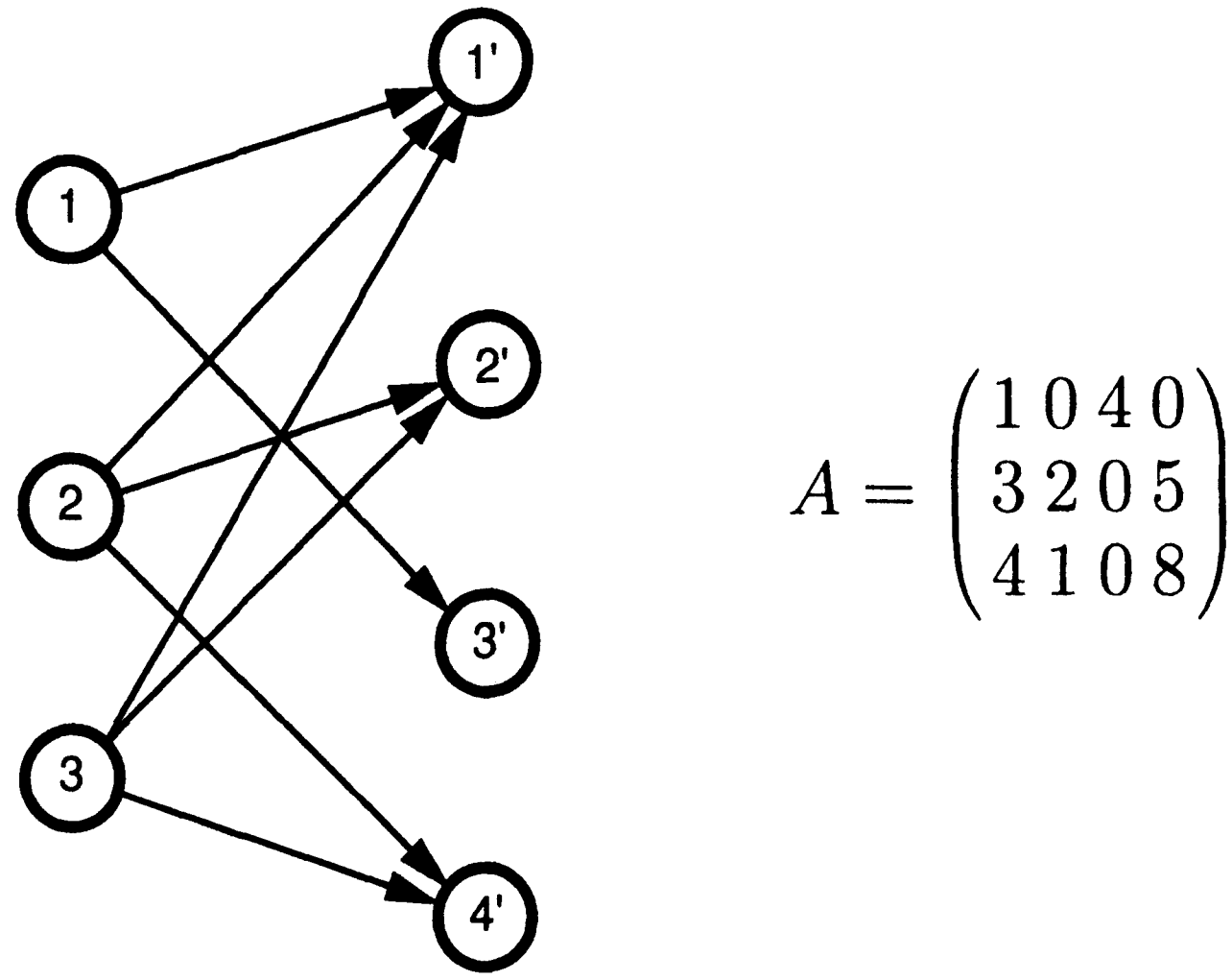

$$A = \begin{pmatrix} 1\,0\,4\,0 \\ 3\,2\,0\,5 \\ 4\,1\,0\,8 \end{pmatrix}$$

**Figure 9.2** The bipartite graph of a matrix $A$.

*tion problem* (see also Section 12.2.1). With each arc $(i,j)$ we associate a variable $x_{ij}$ that denotes flow from node $i \in \mathcal{N}_1$ to node $j \in \mathcal{N}_2$, and denote by $f_{ij}(x_{ij})$ the cost of sending $x_{ij}$ units along arc $(i,j)$. We use $i = 1,2,\ldots,m$, to denote the nodes in $\mathcal{N}_1$, and $j = 1,2,\ldots,n$, to denote the nodes in $\mathcal{N}_2$. We define now an optimization problem that minimizes the total cost $\sum_{(i,j)\in\mathcal{A}} f_{ij}(x_{ij})$, subject to conditions on the conservation of total flow at each node:

$$\sum_{j\in\delta_i^+} x_{ij} = u_i, \quad i = 1,2,\ldots,m, \tag{9.15}$$

$$\sum_{i\in\delta_j^-} x_{ij} = v_j, \quad j = 1,2,\ldots,n, \tag{9.16}$$

$$x_{ij} \geq 0, \text{ for all } (i,j) \in \mathcal{A}. \tag{9.17}$$

These are precisely the consistency conditions required for Problem 9.2.1. The positive elements of matrix $A$ are viewed as flows on the arcs of $\mathcal{G}$, that is, $a_{ij} > 0$ denotes flow on arc $(i,j)$. Thus the relationship between Problem 9.2.1 and transportation optimization problems is established.

We turn now to Problem 9.2.2. The natural graph for this problem is a *transshipment graph* with $n$ nodes and an arc for every nonzero entry of $A$ as defined next (see Figure 9.3).

**Definition 9.2.2 (Transshipment graph of matrix A)** *Given a square $n\times n$ nonnegative matrix $A$, the graph $\mathcal{G} = (\mathcal{N},\mathcal{A})$ is called the transshipment graph of $A$ if $\mathcal{N} \doteq \{1,2,\ldots,n\}$ and $\mathcal{A} \doteq \{(i,j) \mid a_{ij} > 0\}$.*

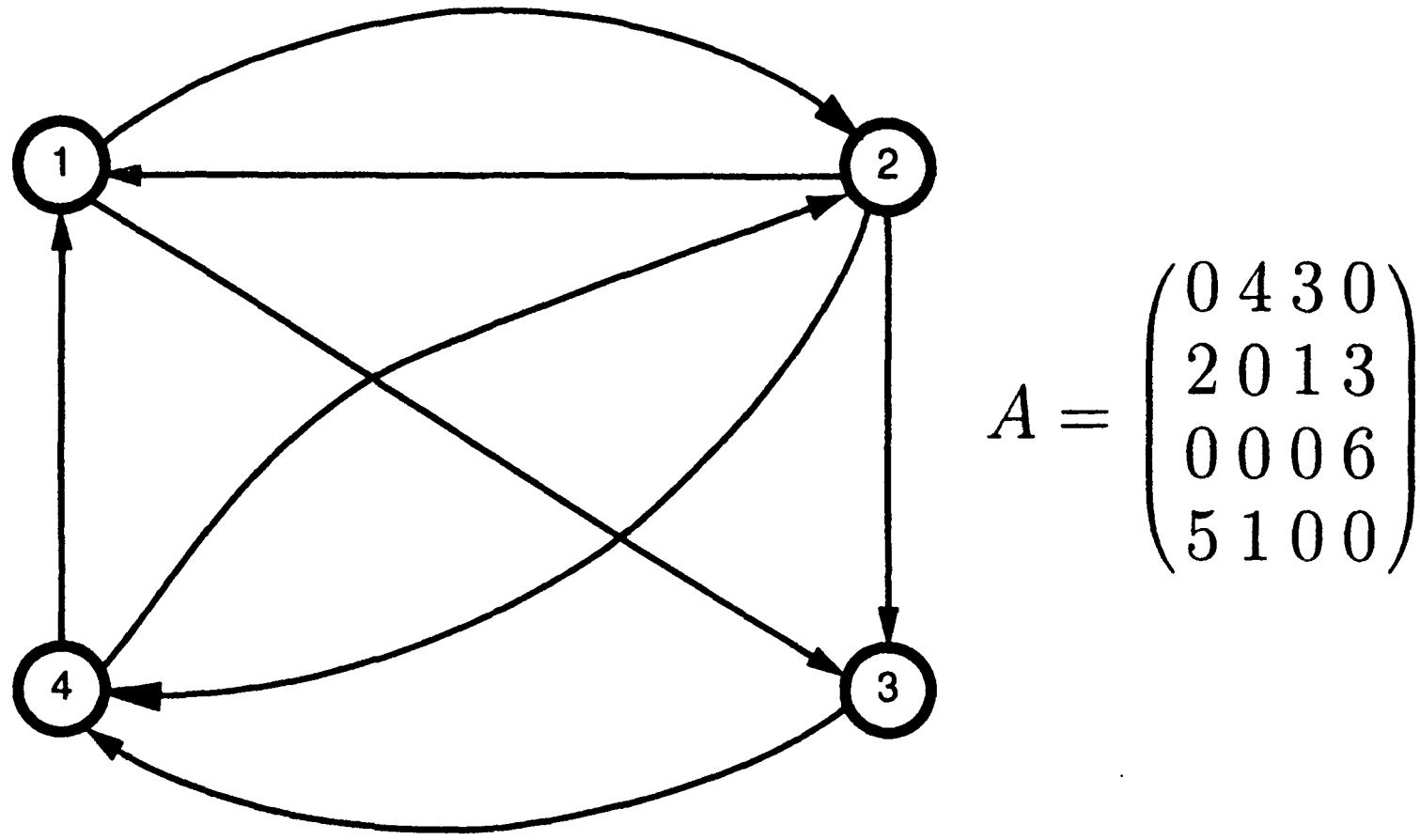

**Figure 9.3** The transshipment graph of a matrix $A$.

Similar to the transportation optimization problem defined over a bipartite graph, we can define a transshipment optimization problem over a transshipment graph. In this graph model, the entries of $A$ are flows on the arcs of $\mathcal{G}$. The conservation of flow conditions of the transshipment problem are precisely the consistency requirements for Problem 9.2.2, which can be written as $\sum_{j\in\delta_i^+} x_{ij} = \sum_{j\in\delta_i^-} x_{ji}, \quad i = 1, 2, \ldots, n.$

### 9.2.3 Entropy optimization models for matrix balancing

We consider now entropy optimization formulations of Problems 9.2.1 and 9.2.2. That is, we consider the estimated and balanced matrix $X$ to be nearby the matrix $A$ if the entropic distance between them is minimized. Entropy optimization formulations have been justified in the literature as resulting in balanced matrices that are the least biased, or maximally uncommitted with respect to missing information from the original matrix $A$. Some axiomatic approaches also lead naturally to entropy optimization formulations. Further justification for the use of entropy maximization models can be deduced from Section 10.3 if we note the similarity between matrix balancing and fully discretized models for image reconstruction from two orthogonal projections. As a result matrix balancing problems formulated as entropy optimization models over network flow constraints are some of the most widely used in practice. The entropy optimization models for Problems 9.2.1 and 9.2.2 are defined as follows. Note that in all models the optimization problem does not define variables for $(i, j) \notin \mathcal{A}$ (i.e., when $a_{ij} = 0$) and the corresponding entries of the matrix $X$ are set to zero.

**Problem 9.2.3** *Given an $m \times n$ nonnegative matrix $A = (a_{ij})$ and positive vectors $u = (u_i) \in \mathbb{R}^m$ and $v = (v_j) \in \mathbb{R}^n$, determine the matrix $X = (x_{ij})$* *that solves the optimization problem:*

$$\text{Minimize} \quad \sum_{(i,j)\in\mathcal{A}} x_{ij}\left(\log\left(\frac{x_{ij}}{a_{ij}}\right)-1\right) \tag{9.18}$$

$$\text{s.t.} \quad \sum_{j\in\delta_i^+} x_{ij} = u_i, \qquad \text{for } i=1,2,\ldots,m, \tag{9.19}$$

$$\sum_{i\in\delta_j^-} x_{ij} = v_j, \qquad \text{for } j=1,2,\ldots,n, \tag{9.20}$$

$$x_{ij} \geq 0, \qquad \text{for all } (i,j)\in\mathcal{A}. \tag{9.21}$$

**Problem 9.2.4** *Given an $n \times n$ nonnegative matrix $A = (a_{ij})$, determine the matrix $X = (x_{ij})$ that solves the optimization problem:*

$$\text{Minimize} \quad \sum_{(i,j)\in\mathcal{A}} x_{ij}\left(\log\left(\frac{x_{ij}}{a_{ij}}\right)-1\right) \tag{9.22}$$

$$\text{s.t.} \quad \sum_{j\in\delta_i^+} x_{ij} = \sum_{j\in\delta_i^-} x_{ji}, \qquad \text{for } i=1,2,\ldots,n, \tag{9.23}$$

$$x_{ij} \geq 0, \qquad \text{for all } (i,j)\in\mathcal{A}. \tag{9.24}$$

In practical applications when the prescribed row and column sums are unreliable or when they cannot be specified more precisely than being confined to certain intervals, we cannot use these formulations. Therefore we now consider extensions of these problems to situations in which the given row and column sums are confined to an interval or in which the strict agreement between row and column sums is relaxed such that their differences fall within a specified interval. These are known as *interval-constrained matrix balancing* problems.

**Problem 9.2.5** *Given an $m\times n$ nonnegative matrix $A = (a_{ij})$, nonnegative vectors $\underline{u} = (\underline{u}_i)$, $\overline{u} = (\overline{u}_i) \in \mathbb{R}^m$ such that $0 \leq \underline{u} \leq \overline{u}$, and nonnegative vectors $\underline{v} = (\underline{v}_j)$, $\overline{v} = (\overline{v}_j) \in \mathbb{R}^n$ such that $0 \leq \underline{v} \leq \overline{v}$, determine the matrix $X = (x_{ij})$ that solves the optimization problem:*

$$\text{Minimize} \quad \sum_{(i,j)\in\mathcal{A}} x_{ij}\left(\log\left(\frac{x_{ij}}{a_{ij}}\right)-1\right) \tag{9.25}$$

$$\text{s.t.} \quad \underline{u}_i \leq \sum_{j\in\delta_i^+} x_{ij} \leq \overline{u}_i, \qquad \text{for } i=1,2,\ldots,m, \tag{9.26}$$

$$\underline{v}_j \leq \sum_{i\in\delta_j^-} x_{ij} \leq \overline{v}_j, \qquad \text{for } j=1,2,\ldots,n, \tag{9.27}$$

$$x_{ij} \geq 0, \qquad \text{for all } (i,j)\in\mathcal{A}. \tag{9.28}$$

**Problem 9.2.6** *Given an $n \times n$ nonnegative matrix $A = (a_{ij})$ and a positive tolerance vector $\epsilon = (\epsilon_i) \in \mathbb{R}^n$, $\epsilon > 0$, determine the matrix $X = (x_{ij})$ that solves the optimization problem:*

$$\text{Minimize} \quad \sum_{(i,j)\in\mathcal{A}} x_{ij}\left(\ln\left(\frac{x_{ij}}{a_{ij}}\right) - 1\right) \tag{9.29}$$

$$\text{s.t.} \quad -\epsilon_i \le \sum_{j\in\delta_i^+} x_{ij} - \sum_{j\in\delta_i^-} x_{ji} \le \epsilon_i, \quad \text{for } i = 1,2,\ldots,n, \tag{9.30}$$

$$x_{ij} \ge 0, \quad \text{for all } (i,j) \in \mathcal{A}. \tag{9.31}$$

## 9.3 Iterative Algorithms for Matrix Balancing

We apply now the row-action iterative algorithms of Chapter 6 to the entropy optimization formulations of the matrix balancing problems. The entropy function is a Bregman function (see Lemma 2.1.3), and it is optimized over equality constraints (for Problems 9.2.3 and 9.2.4) and interval constraints (for Problems 9.2.5 and 9.2.6). Hence, we can develop specializations of the general row-action Algorithm 6.4.1. We develop here the algorithms for the interval-constrained Problems 9.2.5 and 9.2.6. Then we obtain, as special cases, the algorithms for the equality-constrained Problems 9.2.3 and 9.2.4.

### 9.3.1 The range-RAS algorithm (RRAS)

The range-RAS (RRAS) algorithm is designed to solve Problem 9.2.5 with interval constraints on the row and column sums of $X$. We call the algorithm range-RAS to indicate its relation with the well-known RAS algorithm (see Sections 9.3.2 and 9.4). It is an adaptation of the iterative row-action Algorithm 6.4.1 for interval convex programming. While the general algorithm is a dual ascent method that does not perform exact minimization in the dual problem, in the particular case of RRAS the dual ascent is indeed exact, due to the fact that the constraint matrix (denoted by $\Phi$ below) is a $0-1$ matrix, i.e., all its entries are either zero or one.

It is useful to reformulate the constraints of Problem 9.2.5 by lexicographically rearranging the $m \times n$ matrix $X$ into an $mn$-dimensional vector $x = (x_s), s = 1,2,\ldots,mn$, where $s \doteq (i-1)m + j$. The matrix $A$ is similarly arranged into a vector $(a_s)$, $s = 1,2,\ldots,mn$. The constraints then take the form

$$\begin{pmatrix} \underline{u} \\ \underline{v} \end{pmatrix} \le \Phi x \le \begin{pmatrix} \overline{u} \\ \overline{v} \end{pmatrix}. \tag{9.32}$$

Using the index sets $\delta_i^+$ and $\delta_j^-$ defined in the previous section the matrix $\Phi$ is defined as

$$\Phi \doteq \begin{pmatrix} R \\ C \end{pmatrix}, \tag{9.33}$$

where the $m \times mn$ matrix $R$ is given by

$$R \doteq \begin{pmatrix} \alpha_{11} \cdots \alpha_{1n} & | & & | & & | & \\ & | & \alpha_{21} \cdots \alpha_{2n} & | & & | & \\ & | & & | & \ddots & | & \\ & | & & | & & | & \alpha_{m1} \cdots \alpha_{mn} \end{pmatrix} \tag{9.34}$$

with entries

$$\alpha_{ij} = \begin{cases} 1, & \text{if } j \in \delta_i^+, \\ 0, & \text{otherwise.} \end{cases} \tag{9.35}$$

The remaining entries of $R$ are all zero. $C$ is the $n \times mn$ matrix:

$$C \doteq \begin{pmatrix} \beta_{11} & & & | & \beta_{21} & & & | & & | & \beta_{m1} & & \\ & \ddots & & | & & \ddots & & | & \cdots & | & & \ddots & \\ & & \beta_{1n} & | & & & \beta_{2n} & | & & | & & & \beta_{mn} \end{pmatrix} \tag{9.36}$$

with entries $\beta_{ij}$ given by

$$\beta_{ij} = \begin{cases} 1, & \text{if } i \in \delta_j^-, \\ 0, & \text{otherwise,} \end{cases} \tag{9.37}$$

and the remaining entries of $C$ all zero.

In Algorithm 6.4.1 we identify $\Phi$ as above as the constraint matrix, let $\gamma \doteq \begin{pmatrix} \underline{u} \\ \underline{v} \end{pmatrix}$ and $\delta \doteq \begin{pmatrix} \overline{u} \\ \overline{v} \end{pmatrix}$, and take the objective function $f$ to be

$$f(x) = \sum_{s=1}^{mn} x_s \left( \log \left( \frac{x_s}{a_s} \right) - 1 \right). \tag{9.38}$$

Then the $\nu$th iterative step derived from Algorithm 6.4.1 takes the form

$$x_s^{\nu+1} = x_s^\nu \exp\left(c_\nu \phi_s^{l(\nu)}\right), \quad s = 1, 2, \ldots, mn, \tag{9.39}$$

$$z_t^{\nu+1} = \begin{cases} z_t^\nu, & \text{if } t \neq l(\nu), \\ z_t^\nu - c_\nu, & \text{if } t = l(\nu). \end{cases} \tag{9.40}$$

Here $\{l(\nu)\}$ is the control sequence over the set $\{1, 2, \ldots, m+n\}$ according to which the rows of $\Phi$ are chosen, and $\phi^l = (\phi_s^l) = (\phi_{ls})$, $s = 1, 2, \ldots, mn$, is the $l$th column of $\Phi^T$, the transpose of $\Phi$. We assume the cyclic control sequence $l(\nu) = \nu \bmod (m+n) + 1$ and henceforth abbreviate $l \equiv l(\nu)$. $\{z^\nu\}$ is the sequence of dual vectors.

The parameter $c_\nu$ is computed, for every $\nu \geq 0$, by

$$c_\nu = \text{mid}\,(z_l^\nu, \Gamma_\nu, \Delta_\nu). \tag{9.41}$$

The parameters $\Gamma_\nu$ and $\Delta_\nu$ are determined by solving, for $l = l(\nu)$, the systems:

$$y_s = x_s^\nu \exp(\Gamma_\nu \phi_{ls}), \; s = 1, 2, \ldots, mn, \tag{9.42}$$

$$\sum_{s=1}^{mn} y_s \phi_{ls} = \gamma_l, \tag{9.43}$$

and

$$w_s = x_s^\nu \exp(\Delta_\nu \phi_{ls}), \; s = 1, 2, \ldots, mn, \tag{9.44}$$

$$\sum_{s=1}^{mn} w_s \phi_{ls} = \delta_l. \tag{9.45}$$

The solution for $\Gamma_\nu$ is obtained as follows. Consider the first $m$ row indices of $\Phi$, $l = 1, 2, \ldots, m$ that correspond to row indices $i = 1, 2, \ldots, m$ of the matrix $R$. Substitute (9.42) into (9.43) to get $\sum_{s=1}^{mn} x_s^\nu \exp(\Gamma_\nu \phi_{is}) \phi_{is} = \gamma_i$. Use now the fact that for the indices $l \doteq i = 1, 2, \ldots, m$ we have $\gamma_l = \underline{u}_i$, and $\phi_{is}$ takes the value 1 if the arc in $\mathcal{A}$ with lexicographic order $s$ has origin node $i$, or 0 otherwise as specified by (9.35). Then, reverting to the graph notation of Problem 9.2.5, we can write the equation for $\Gamma_\nu$ as $(\exp \Gamma_\nu) \sum_{(i,j)\in\mathcal{A}} x_{ij}^\nu = \underline{u}_i$, or, in the equivalent form $(\exp \Gamma_\nu) \sum_{j\in\delta_i^+} x_{ij}^\nu = \underline{u}_i$, from which we obtain $\Gamma_\nu$. The expression for $\Delta_\nu$, for $l = 1, 2, \ldots, m$, can be derived in a similar way from (9.44)–(9.45).

Consider now the derivations for $\Gamma_\nu$ and $\Delta_\nu$ for row indices of $\Phi$ given by $l = m+1, m+2, \ldots, m+n$. These indices correspond to row indices $j = 1, 2, \ldots, n$ of the matrix $C$. Substituting again (9.42) into (9.43) we get $\sum_{s=1}^{mn} x_s^\nu \exp(\Gamma_\nu \phi_{is}) \phi_{is} = \gamma_i$. Use now the fact that for all the indices $l = m+1, m+2, \ldots, m+n$ we have $\gamma_l = \underline{v}_j$ for an appropriate index $j$, and $\phi_{ls}$ takes the value 1 if the arc in $\mathcal{A}$ with lexicographic order $s$ has destination node $j$, or 0 otherwise as specified by (9.37). Then we can write the equation for $\Gamma_\nu$ as $(\exp \Gamma_\nu) \sum_{(i,j)\in\mathcal{A}} x_{ij}^\nu = \underline{v}_j$, or in the equivalent form $(\exp \Gamma_\nu) \sum_{i\in\delta_j^-} x_{ij}^\nu = \underline{v}_j$ from which we obtain $\Gamma_\nu$. The expressions for $\Delta_\nu$ can be derived in a similar way.

To avoid the exponentiation terms (appearing in equation (9.39)) in the execution of the algorithm we work with the exponents of the parameters $\Gamma_\nu$ and $\Delta_\nu$. The parameter $c_\nu$ (eq. (9.41)) can be obtained using the fact that for any real triplet $x, y, z$,

$$\log(\text{mid}\ (\exp x, \exp y, \exp z)) = \text{mid}\ (x, y, z).$$

We also use $\underline{\rho}_i^\nu, \overline{\rho}_i^\nu$ to denote $\exp \Gamma_\nu$ and $\exp \Delta_\nu$ respectively, when row-index $l$ of $\Phi$ corresponds to the $i$th row of the matrix $R$, and $\underline{\sigma}_j^\nu, \overline{\sigma}_j^\nu$ to denote $\exp \Gamma_\nu$ and $\exp \Delta_\nu$ respectively, when row-index $l$ of $\Phi$ corresponds to the $j$th row of the matrix $C$. Then the row-action algorithm for the interval-constrained Problem 9.2.5 can be summarized as follows:

### Algorithm 9.3.1 The range-RAS (RRAS) Algorithm

**Step 0:** (Initialization.) Set $\nu = 0$ and $x_{ij}^0 = a_{ij}$ for all $i = 1, 2, \ldots, m$, and $j = 1, 2, \ldots, n$. Set $\rho_i^0 = 1$ for $i = 1, 2, \ldots, m$, and $\sigma_j^0 = 1$ for all $j = 1, 2, \ldots, n$.

**Step 1:** (Row scaling.) For $i = 1, 2, \ldots, m$ define

$$\underline{\rho}_i^\nu \doteq \frac{\underline{u}_i}{\sum_{j \in \delta_i^+} x_{ij}^\nu}, \tag{9.46}$$

$$\overline{\rho}_i^\nu \doteq \frac{\overline{u}_i}{\sum_{j \in \delta_i^+} x_{ij}^\nu}, \tag{9.47}$$

compute

$$\Delta\rho_i^\nu \doteq \text{mid}\,(\rho_i^\nu, \underline{\rho}_i^\nu, \overline{\rho}_i^\nu), \tag{9.48}$$

and update:

$$x_{ij}^\nu \leftarrow x_{ij}^\nu \Delta\rho_i^\nu, \quad i = 1, 2, \ldots, m, \; j = 1, 2, \ldots, n, \tag{9.49}$$

$$\rho_i^{\nu+1} = \frac{\rho_i^\nu}{\Delta\rho_i^\nu}, \quad i = 1, 2, \ldots, m. \tag{9.50}$$

**Step 2:** (Column scaling.) For $j = 1, 2, \ldots, n$ define

$$\underline{\sigma}_j^\nu \doteq \frac{\underline{v}_j}{\sum_{i \in \delta_j^-} x_{ij}^\nu}, \tag{9.51}$$

$$\overline{\sigma}_j^\nu \doteq \frac{\overline{v}_j}{\sum_{i \in \delta_j^-} x_{ij}^\nu}, \tag{9.52}$$

compute

$$\Delta\sigma_j^\nu \doteq \text{mid}\left(\sigma_j^\nu, \underline{\sigma}_j^\nu, \overline{\sigma}_j^\nu\right), \tag{9.53}$$

and update:

$$x_{ij}^{\nu+1} = x_{ij}^\nu \Delta\sigma_j^\nu, \quad i = 1, 2, \ldots, m, \; j = 1, 2, \ldots, n, \tag{9.54}$$

$$\sigma_j^{\nu+1} = \frac{\sigma_j^\nu}{\Delta\sigma_j^\nu}, \quad j = 1, 2, \ldots, n. \tag{9.55}$$

**Step 3:** Replace $\nu \leftarrow \nu + 1$ and return to Step 1.

### Decompositions for Parallel Computing

The range-RAS algorithm decomposes naturally for parallel computations. The row scaling operations at Step 1 can proceed concurrently for all rows $i = 1, 2, \ldots, m$. The calculation of the scaling parameters $\underline{\rho}_i^\nu, \overline{\rho}_i^\nu$ involves only the flows of arcs incident to the $i$th node. Since no two nodes corresponding to the matrix rows have arcs in common it is possible to compute the scaling parameters $\underline{\rho}_i^\nu$ and $\overline{\rho}_i^\nu$ for all values of $i$ independent of each other. The algorithm could utilize up to $m$ processors in calculating these parameters. Similarly the flow variables $x_{ij}^\nu$ can be updated concurrently for all arcs incident to node $i$ and for all indices $i$, using as many processors as there are arcs. The scaling parameters $\rho_i^{\nu+1}$ can also be updated concurrently for all $i$. The column scaling operations parallelize in a similar fashion. The parallel implementation of RRAS is similar to that of the transportation algorithms in Chapter 12. Chapter 14 gives the data structures for the parallel implementation of these algorithms.

### 9.3.2 The RAS scaling algorithm

If the input to the RRAS algorithm is $\underline{u} = \overline{u}$ and $\underline{v} = \overline{v}$ then for all $\nu$ $\underline{\rho}_i^\nu = \overline{\rho}_i^\nu$ and $\underline{\sigma}_j^\nu = \overline{\sigma}_j^\nu$ will always be equal to $\Delta\rho_i^\nu$ and $\Delta\sigma_j^\nu$, respectively. The scaling parameters $\rho_i^\nu$ and $\sigma_j^\nu$ are therefore not used explicitly and the algorithm coincides with the classic RAS algorithm for Problem 9.2.3 with fixed row and column sums, and has the following form:

### Algorithm 9.3.2 The RAS Algorithm

**Step 0:** (Initialization.) Set $\nu = 0$ and $x_{ij}^0 = a_{ij}$ for all $i = 1, 2, \ldots, m$, and $j = 1, 2, \ldots, n$.

**Step 1:** (Row scaling.) For $i = 1, 2, \ldots, m$, define

$$\rho_i^\nu \doteq \frac{u_i}{\sum_{j \in \delta_i^+} x_{ij}^\nu}, \tag{9.56}$$

and update:

$$x_{ij}^\nu \leftarrow x_{ij}^\nu \rho_i^\nu, \quad i = 1, 2, \ldots, m, \; j = 1, 2, \ldots, n. \tag{9.57}$$

**Step 2:** (Column scaling.) For $j = 1, 2, \ldots, n$, define

$$\sigma_j^\nu \doteq \frac{v_j}{\sum_{i \in \delta_j^-} x_{ij}^\nu}, \tag{9.58}$$

and update:

$$x_{ij}^{\nu+1} = x_{ij}^{\nu}\sigma_j^{\nu}, \quad i = 1,2,\ldots,m, \; j = 1,2,\ldots,n. \tag{9.59}$$

**Step 3:** Replace $\nu \leftarrow \nu + 1$ and return to Step 1.

Within the bipartite graph model of Problem 9.2.3 the RAS algorithm iterates through the nodes on each side of the bipartite graph $\mathcal{G}$, and multiplies the flows on all arcs incident to a node by the positive constant necessary to force the sum of the flows to equal the desired total. The algorithm can be implemented in parallel in the same way as RRAS.

### 9.3.3 The range-DSS algorithm (RDSS)

The range-DSS (RDSS) algorithm solves Problem 9.2.6. We call the algorithm range-DSS to denote its relation to the diagonal similarity scaling (DSS) algorithm of Section 9.3.4. Here we derive range-DSS from the general Algorithm 6.4.1.

The interval constraints of Problem 9.2.6 may be written as

$$-\epsilon \leq (R - C)x \leq \epsilon, \tag{9.60}$$

where $R - C$ is the difference matrix of $R$ and $C$ given by (9.34) and (9.36), where both $R$ and $C$ are, for Problem 9.2.6, of dimension $n \times n^2$. That is, $R - C$ is as follows:

$$\left(\begin{array}{cc|ccc|c|ccc}
0 & \alpha_{12}\cdots\alpha_{1n} & -\beta_{21} & & & & -\beta_{n1} & & \\
-\beta_{12} & & \alpha_{21} & 0\cdots & \alpha_{2n} & & & -\beta_{n2} & \\
\ddots & & & \ddots & & \cdots & & & \ddots \\
 & -\beta_{1n} & & & -\beta_{2n} & & \alpha_{n1} & \alpha_{n2}\cdots & 0
\end{array}\right). \tag{9.61}$$

Let us denote

$$R - C = \Psi = (\psi_{is}), \tag{9.62}$$

and let $\psi^i = (\psi_s^i) = (\psi_{is})$ for $s = 1,2,\ldots,n^2$ be the $i$th column of $\Psi^T$, the transpose of $\Psi$. When applied to the entropy objective function of Problem 9.2.6 with the interval constraints above, the iterative step derived from Algorithm 6.4.1 takes the form

$$x_s^{\nu+1} = x_s^{\nu}\exp\left(c_\nu \psi_s^{i(\nu)}\right), \; s = 1,2,\ldots,n^2, \tag{9.63}$$

$$z_t^{\nu+1} = \begin{cases} z_t^{\nu}, & \text{if } t \neq i(\nu), \\ z_t^{\nu} - c_\nu, & \text{if } t = i(\nu). \end{cases} \tag{9.64}$$

Here $\{i(\nu)\}$ is the control sequence over the set $\{1,2,\ldots,n\}$ according to which the rows of $\Psi$ are chosen. We assume the cyclic control given by $i(\nu) = \nu(\bmod n) + 1$ and henceforth abbreviate it as $i \equiv i(\nu)$. $\{z^\nu\}$ is again a sequence of dual vectors. In this iteration the parameter $c_\nu$ is

$$c_\nu = \text{mid } (z_i^\nu, \underline{d}^\nu, \overline{d}^\nu), \tag{9.65}$$

where $\underline{d}^\nu$ and $\overline{d}^\nu$ are determined by solving, for $i = i(\nu)$, the systems:

$$y_s = x_s^\nu \exp\left(\underline{d}_\nu \psi_{is}\right), \ s = 1, 2, \ldots, n^2, \tag{9.66}$$

$$\sum_{s=1}^{n^2} y_s \psi_{is} = -\epsilon_i, \tag{9.67}$$

and

$$w_s = x_s^\nu \exp\left(\overline{d}_\nu \psi_{is}\right), \ s = 1, 2, \ldots, n^2, \tag{9.68}$$

$$\sum_{s=1}^{n^2} w_s \psi_{is} = \epsilon_i. \tag{9.69}$$

In order to solve these systems consider the graph representation of the problem, i.e., the transshipment graph of Definition 9.2.2. From (9.66) and (9.67), and denoting $\exp \underline{d}_\nu \doteq \underline{\alpha}_\nu$, we get $\sum_{s=1}^{n^2} \psi_{is} x_s^\nu (\underline{\alpha}_\nu)^{\psi_{is}} = -\epsilon_i$. We use now the fact that $\psi_{is}$ takes the value $+1$ if the arc in $\mathcal{A}$ with lexicographic order $s$ has origin node $i$, the value $-1$ if the $s$th arc has destination node $i$, or 0 otherwise, as specified by (9.61), and the expressions (9.35) and (9.37). Then the equation for $\underline{\alpha}_\nu$ can be written as:

$$p_i^\nu \underline{\alpha}_\nu - q_i^\nu \frac{1}{\underline{\alpha}_\nu} = -\epsilon_i, \tag{9.70}$$

where

$$p_i^\nu \doteq \sum_{j \in \delta_i^+} x_{ij}^\nu, \quad q_i^\nu \doteq \sum_{j \in \delta_i^-} x_{ji}^\nu. \tag{9.71}$$

Similarly, for $\overline{\alpha}_\nu \doteq \exp \overline{d}_\nu$ we obtain from (9.68) and (9.69),

$$p_i^\nu \overline{\alpha}_\nu - q_i^\nu \frac{1}{\overline{\alpha}_\nu} = \epsilon_i. \tag{9.72}$$

Taking the nonnegative solutions of the quadratic equations (9.70) and (9.72) we obtain the values of $\underline{d}_\nu$ and $\overline{d}_\nu$ respectively, which go into the *mid* operator in (9.65) to obtain $c_\nu$. Using $c_\nu$ in (9.63)–(9.64) completes the iterative step. The RDSS algorithm is thus formulated as follows.

**Algorithm 9.3.3 The range-DSS (RDSS) Algorithm**

**Step 0:** (Initialization.) Set $\nu = 0$ and $x_{ij}^0 = a_{ij}$ for all $i, j = 1, 2, \ldots, n$. Set $\rho_i^0 = 1$ for $i = 1, 2, \ldots, n$.

**Step 1:** (Computation of scaling parameters.) Choose a control index from a cyclic control sequence $\{i(\nu)\}$ with $i \doteq i(\nu) = \nu(\mathrm{mod}\, n) + 1$, and calculate the sums

$$p_i^\nu \doteq \sum_{j \in \delta_i^+} x_{ij}^\nu \,, \quad q_i^\nu \doteq \sum_{j \in \delta_i^-} x_{ji}^\nu. \tag{9.73}$$

Compute $\underline{\alpha}_\nu$ and $\overline{\alpha}_\nu$ as the nonnegative roots of (9.70) and (9.72), respectively.

**Step 2:** (Update.) Calculate

$$\Delta\rho_i^\nu = \text{mid } (\rho_i^\nu, \underline{\alpha}_\nu, \overline{\alpha}_\nu)\,, \tag{9.74}$$

and update:

$$x_{ij}^{\nu+1} = x_{ij}^\nu \Delta\rho_i^\nu, \text{ for all } j \in \delta_i^+, \tag{9.75}$$

$$x_{ij}^{\nu+1} = \frac{x_{ij}^\nu}{\Delta\rho_i^\nu}, \text{ for all } j \in \delta_i^-, \tag{9.76}$$

$$x_{ij}^{\nu+1} = x_{ij}^\nu, \text{ otherwise.} \tag{9.77}$$

$$\rho_i^{\nu+1} = \frac{\rho_i^\nu}{\Delta\rho_i^\nu}. \tag{9.78}$$

**Step 3:** Replace $\nu \leftarrow \nu + 1$ and return to Step 1.

## Decompositions for Parallel Computing

At the $\nu$th iteration the algorithm computes at Step 1 the scaling parameters $\underline{\alpha}_\nu$, $\overline{\alpha}^\nu$ for a given row/column of the matrix, whose index is $i$. It then updates at Step 2 the flow of all arcs of the underlying graph, incident to the node associated with the $i$th row/column. The calculation of the scaling parameters is essentially a serial operation since a single pair of parameters is calculated for a given index $i$. It is possible however, to exploit parallelism in the calculation of the sums in equation (9.73) using as many processors as the number of arcs incident to node $i$ to compute the sums. This step requires communication among processors. Once the scaling parameters are computed the calculations of Step 2 can be executed in parallel, utilizing as many processors as the number of arcs incident to the $i$th node.

It is also possible to operate on multiple nodes/rows concurrently. This can be achieved by selecting those pairs of node/row indices that correspond to matrix entries with zero value. For example, if the nodes corresponding to rows/columns $i$ and $j$ do not have an arc in common, then the scaling parameters for $i$ and $j$ can be computed concurrently and the flows of all arcs incident to both $i$ and $j$ can be updated in parallel. Identifying sets of nodes that do not have arcs in common is equivalent to solving the classic problem of *coloring* the graph. Steps 1 and 2 are then executed concurrently on all nodes with the same color. Note that a control sequence that selects nodes based on their color is allowed, since the convergence of the algorithm does not depend on the cyclic control sequence specified for simplicity in Step 1. Instead, any other almost cyclic sequence can be specified that chooses nodes with the same color.

### 9.3.4 The diagonal similarity scaling (DSS) algorithm

In the special case when the RDSS algorithm is applied to Problem 9.2.6 with $\epsilon = 0$ we obtain the *diagonal similarity scaling* (DSS) algorithm that solves Problem 9.2.4. In this case $\underline{\alpha}_\nu = \overline{\alpha}_\nu = \alpha_\nu$ and therefore $\underline{d}_\nu = \overline{d}_\nu = d_\nu = c_\nu$. Since $\alpha_\nu = \exp c_\nu$ the iteration (9.63) becomes

$$x_s^{\nu+1} = x_s^\nu (\alpha_\nu)^{\psi_{is}}.$$

Since then also $\alpha_\nu = \sqrt{q_i^\nu / p_i^\nu}$, (from (9.70) or (9.72)), the DSS algorithm can be written as follows.

**Algorithm 9.3.4 The DSS Algorithm**

**Step 0:** (Initialization.) Set $\nu = 0$ and $x_{ij}^0 = a_{ij}$ for all $i, j = 1, 2, \dots, n$.

**Step 1:** (Computation of scaling parameters.) Choose a control index from a cyclic control sequence $\{i(\nu)\}$ with $i \doteq i(\nu) = \nu (\bmod n) + 1$ and calculate the sums

$$p_i^\nu \doteq \sum_{j \in \delta_i^+} x_{ij}^\nu \,, \quad q_i^\nu \doteq \sum_{j \in \delta_i^-} x_{ji}^\nu. \tag{9.79}$$

Compute $\alpha_\nu = \sqrt{\dfrac{q_i^\nu}{p_i^\nu}}$.

**Step 2:** (Update.) Update:

$$x_{ij}^{\nu+1} = x_{ij}^\nu \alpha_\nu, \text{ for all } j \in \delta_i^+, \tag{9.80}$$

$$x_{ij}^{\nu+1} = \frac{x_{ij}^\nu}{\alpha_\nu}, \text{ for all } j \in \delta_i^-, \tag{9.81}$$

$$x_{ij}^{\nu+1} = x_{ij}^\nu, \text{ otherwise}. \tag{9.82}$$

**Step 3:** Replace $\nu \leftarrow \nu + 1$ and return to Step 1.

With the transshipment graph model (Definition 9.2.2) of the problem, the DSS algorithm scans the nodes of $\mathcal{G}$ and adjusts the flows on arcs directed into and out of that node so that conservation of flow at that node is satisfied. The algorithm can be implemented in parallel in a way similar to RDSS.

## 9.4 Notes and References

For a general introduction to models and algorithms for matrix balancing refer to the survey by Schneider and Zenios (1990). A comprehensive treatment of the problem of matrix estimation in the context of economic modeling is given by Bacharach (1970). Morgenstern (1963) gives an interesting account of the sources of errors in economic observations and the problems thus created for formal economic modeling. This problem is prevalent in economics: Morgenstern's book was first published in condensed form in 1953; an expanded version was published in 1963; and it was reprinted "on demand" as recently as 1985.

**9.1** The volume by Pyatt and Round (1985) contains a collection of papers that give an introduction to Social Accounting Matrices (SAM). The use of SAMs as the database for estimating parameters in economic equilibrium models is discussed, for example, in Bacharach (1970), Baker, Van der Ploeg, and Weale (1984), Dervis, De Melo, and Robinson (1982), Van Tongeren (1986), and Stone (1985). Computerized modeling systems for the estimation of Social Accounting Matrices were developed for The World Bank by Kendrick and Drud (1985) and Zenios, Drud, and Mulvey (1989). The estimation of input-output tables is treated in Miller and Blair (1985). See also Morrison and Thumann (1980), Harrigan and Buchanan (1984) or Jensen and McGaurr (1977).

The estimation of origin-destination tables for transportation applications is discussed in Sheffi (1985) and in papers by Abdfulaal and Le Blanc (1979), Le Blanc and Farhangian (1982), Nguyen (1974), Carey, Hendrickson, and Siddharathan (1981), Jefferson and Scott (1979), McNeil (1983), and Van Zuylen and Willumsen (1980). The model arising when the aggregate flow between an origin-destination pair is given by some demand function is addressed by Nguyen (1974) and Ben-Akiva (1987).

Deming and Stephan (1940) and Stephan (1942) describe applications of matrix balancing procedures for deriving probability estimates from the United States 1940 census data. Friedlander (1961) discusses similar applications by the British government. Ireland and Kullback (1968) discuss algorithms for estimating the underlying cell probabilities of two-way and four-way contingency tables. Darroch and Ratcliff (1982) discuss further applications in statistics.

Interregional migration models are developed by Plane (1982) and Eriksson (1980). Plane (1982) develops the minimum-information model and shows that it contains as special cases the so-called *gravity* and the *fixed-rate* models. Applications of matrix estimation in stochastic modeling using market-share data are discussed in Theil and Rey (1966) and Theil (1967).

**9.2** The literature on models for matrix estimation is vast and our notes here are necessarily incomplete. Furthermore, relationships exist between the various models of matrix estimation. We do not point out all such relationships. The problem of matrix balancing, estimation, or adjustment is discussed in Bachem and Korte (1979, 1981a, 1981b), Erlander, Jörnsten, and Lundgren (1985), and Schneider and Zenios (1990).

The most general form relation is (9.8), proposed by Bachem and Korte (1981a). For a recent discussion and further pointers to the literature on the *equivalence scaling problem* and its relation to the form relation (9.9) see Schneider (1989, 1990). Multiplicative form relations such as (9.10) are discussed by Darroch and Ratcliff (1982), Rothblum (1989, 1992), and Rothblum and Zenios (1992). For special cases of the coefficient matrix $\Phi$, this form relationship also leads to the equivalence scaling problem; see Rothblum (1989). The axiomatic approach has been pursued by Balinski and Demange (1989a, 1989b).

Distance optimization formulations using either the Frobenius distance norm or the entropic norm or combinations of both are discussed in Cottle, Duvall, and Zikan (1986), Zenios, Drud, and Mulvey (1989), Bachem and Korte (1981a, 1981b), and Schneider and Zenios (1990). See also Rothblum and Schneider (1989) on the use of optimization for matrix scaling problems. Interval-constrained formulations were proposed independently and in slightly different forms by Balinksi and Demange (1989a, 1989b) and Censor and Zenios (1991). Another instance of matrix estimation problems is that of *max-scaling*, discussed in Schneider and Schneider (1988) and Rothblum, Schneider, and Schneider (1994). Max-scaling problems are analogous to the problems discussed in this section, but the precise relationships are not explored here.

**9.3** The iterative algorithms for the interval-constrained entropy optimization formulations were developed by Censor and Zenios (1991) who also developed a massively parallel implementation of range-RAS.

The RAS algorithm is an old method that has been independently discovered and analyzed by researchers working in different areas: Kruithof (1937) proposed the method for predicting traffic flows in telecommunication networks and he called it the *method of twin factors*; Deming and Stephan (1940) used the algorithm that they called *method of iterative proportions* to find an approximate solution to the least-squares formulation of problems arising in the estimation of contingency tables; see also Stephan (1942). Other early results are contained in Bacharach (1970), Bregman (1967a), and Sinkhorn (1964). Bregman attributes the method to the Russian architect Sheleikhovskii in the 1930s. The MART algorithm (Algorithm 6.6.2) when applied to $0-1$ matrices coincides with RAS.

The parallel implementation of RAS on different architectures is discussed in Zenios (1990) and Zenios and Iu (1990). For a discussion of the graph coloring problem and a heuristic algorithm for its solution see, e.g., Christofides (1971, 1975).

The DSS algorithm (for problems without bounds on the variables) was first suggested by Schneider and Schneider (1988), and was further extended to the bounded case by Schneider (1989, 1990). Parallel implementations of range-DSS or DSS have not appeared in the literature. Parallelism can be introduced for these problems using the coloring of the underlying graph, as suggested in the context of nonlinear network optimization in Zenios and Mulvey (1988a).

CHAPTER 10

# Image Reconstruction from Projections

---

For many significant problems in diverse fields of applications in science and technology an object—described by a vector or a function $x$—is related to some data $y$ through a relation $\mathcal{O}x = y$. The recovery of $x$ requires essentially the inversion, in some well-specified sense, of the operator $\mathcal{O}$. Such problems are known as *inversion problems.* The inversion problems discussed in this chapter are those of *image reconstruction from projections.* We must caution the reader that *projections* in this context have a different meaning from the *projections* onto convex sets defined and used in Chapters 2 and 5. This double usage has historical roots and cannot be changed but the precise meaning can always be understood from the context.

Image reconstruction problems differ widely depending on the area of application in which they arise. Even within one field—such as medicine or industry—many diverse reconstruction problems occur because of the physically different methods of data collection. In spite of such differences there is a common mathematical nature to reconstruction problems: there is an unknown (two- or three-dimensional) distribution of some physical parameter. This parameter could be, for example, the linear attenuation coefficients of x-rays in human tissue, or the attenuation coefficients of the material in a nuclear reactor, or the density of electrons in the sun's corona. A finite number of line integrals of this parameter can be estimated from physical measurements, and an estimate of the distribution of the original parameter is desired. For example, in the x-ray transmission computerized tomography (CT) case the total attenuation of the x-ray beam between a source and a detector is approximately the integral of the linear attenuation coefficient along the line between the source and the detector. A word of caution regarding notation in this chapter: $x$ and $y$ here sometimes denote Cartesian coordinates that are real variables of a function $f$ and sometimes, even in the same section, they are used also to denote vectors or functions that describe an object $x$ and data $y$. There should be no confusion of meaning since it is always understood from the context.

Let us denote by $f(x, y)$ the *picture* that has to be reconstructed, i.e., the two dimensional function that represents the spatial distribution of the physical parameter of the application at hand. Although $f$ is unknown a priori, in most applications it is known that the distribution is spatially bounded, so that $f$ vanishes outside a finite region $\Omega$ of the plane. It is sometimes more convenient to express $f$ in terms of polar coordinates $(r, \phi)$ rather than Cartesian coordinates $(x, y)$ (see Figure 10.1). In this case we write

$$f(x, y) = f(r \cos \phi, r \sin \phi). \tag{10.1}$$

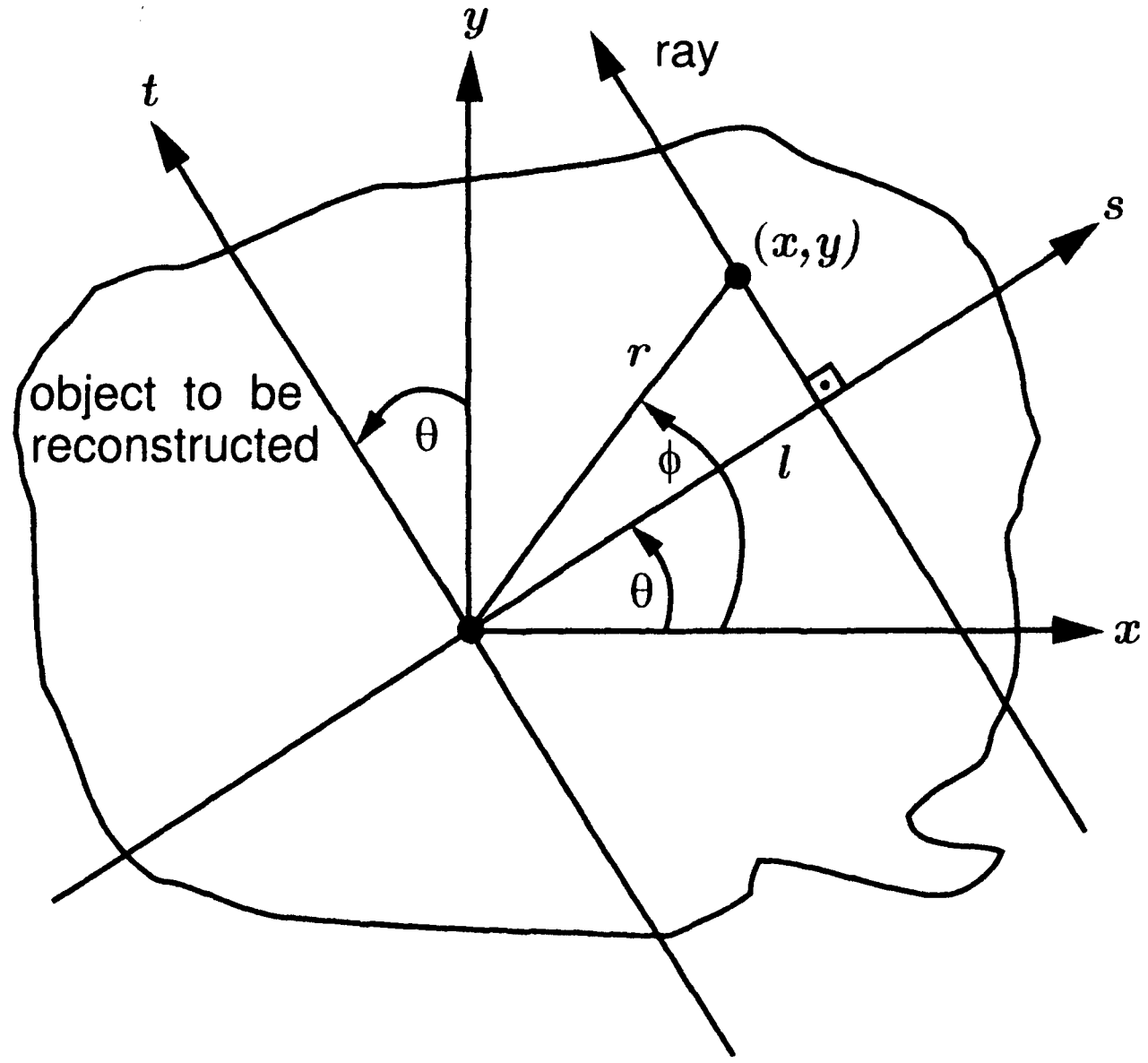

**Figure 10.1** The position of a point in the plane is specified by $(x, y)$ or by $(r \cos \phi, r \sin \phi)$. A line (ray) is specified by its distance $l$ from the origin and its orientation $\theta$.

A line (ray) in the plane is specified by two parameters: its (signed) distance $l$ from the origin and its angle $\theta$ with respect to the $y$-axis. We denote by $p(l, \theta)$ that function of two variables whose value at any $(l, \theta)$ is defined as the line integral of $f$ along the line $L$ specified by $l$ and $\theta$,

$$p(l, \theta) = \int_L f(x, y) dz, \tag{10.2}$$

where the limits of integration depend, in general, on $l, \theta$, and $\Omega$, and $dz$ means integration along the line $L$. Note that for any $l$ and $\theta$, $(l, \theta)$ and $(-l, \theta+\pi)$ represent the same line in the plane, so that $p(l, \theta) = p(-l, \theta+\pi)$.

It is important to note that the arguments of the function $p(l, \theta)$ are *not* the familiar common polar coordinates. To see this, consider two different lines through the origin, with angles $\theta_1$ and $\theta_2$. Generally the integrals of $f$ are different along these two lines so that $p(0, \theta_1)$ is not the same as $p(0, \theta_2)$, which does not make sense if the arguments are interpreted as polar coordinates.

The problem of reconstruction from projections is to find $f(x, y)$ given $p(l, \theta)$. As a mathematical problem, where $f$ and $p$ are functions (unknown and known, respectively), finding $f$ requires solving the integral equation (10.2). In fact the solution (a form of which is given below by equation (10.5)) was published by Radon in 1917 as an inversion formula expressing $f$ in terms of $p$.

The mathematical problem posed above, and solved by Radon, represents an idealized abstraction of the problem as it occurs in practical applications. In practice we are given discrete projection data that are estimates of $p$ for a finite number of rays, and we want to find a two-dimensional image function that is a reconstructed estimate of the unknown object. In addition to the discreteness and limited precision of measured data, we note that a variety of physical problems can lead to significant nonlinearity in the relationship between data and the original object. Because of its mathematical tractability, however, the linear model expressed by (10.2) is the preferred starting point for the derivation of algorithms for image reconstruction.

This chapter only samples the vast literature on the topic focusing on directions that lead to the employment of iterative algorithms. Work in such directions can evidently benefit from parallel algorithmic developments of the kind studied in this book. In Section 10.1 we develop the fully discretized model for transmission CT and describe approaches toward its solution. A fully discretized model for positron emission tomography (PET) is presented in Section 10.2. In Section 10.3 we show how entropy maximization can be justified in image reconstruction and in Section 10.4 the algebraic reconstruction techniques (ART), some of which were studied mathematically in Chapters 5 and 6, are revisited. In Section 10.5 we examine an algorithmic framework for iterative data refinement (IDR) that helps to reduce the effects of discrepancy between an idealized model and a real-world problem of image reconstruction. Iterative algorithms, in spite of their versatility, are not a panacea for all inversion problems. Some guidelines for deciding when their use is appropriate are presented in Section 10.6. Further information and references are discussed in Section 10.7.

## 10.1 Transform Methods and the Fully Discretized Model

Two fundamentally different approaches to the image reconstruction problem exist. In one approach the problem is formulated for continuous functions $f$ and $p$ and an inversion formula is derived in this *continuous model.* Since this approach involves the mathematical study of transforms and their inverses, it is often called the *transform methods approach.* The alternative approach is to discretize the functions involved in the problem at the outset. The object and measurements then become vectors in finite-dimensional Euclidean space and methods from different branches of mathematics, mainly linear algebra and optimization theory, are used to study the model and construct algorithms for solving the problem. This approach is called *finite series-expansion* or the *fully discretized model.* After a brief discussion of Radon's inversion formula we develop the fully discretized model.

In the two-dimensional case we have the following result, which is Radon's inversion formula. Let $f$ be a continuous bounded function of two polar variables such that $f(r, \phi) = 0$, for all $r \geq E$, for some $E > 0$. For any pair of real numbers $(l, \theta)$ define the *Radon transform* $\mathcal{R}f$ of $f$ by

$$[\mathcal{R}f](l, \theta) \doteq \int_{-\infty}^{\infty} f((l^2 + z^2)^{1/2},\ \theta + \arctan(z/l))dz, \text{ if } l \neq 0, \tag{10.3}$$

and

$$[\mathcal{R}f](0, \theta) \doteq \int_{-\infty}^{\infty} f(z,\ \theta + \pi/2)dz. \tag{10.4}$$

Under some physically reasonable technical conditions specified in Radon's theorem we have for all $(r, \phi)$ the following inversion formula, which gives an explicit expression for the inverse Radon transform $\mathcal{R}^{-1}$,

$$[\mathcal{R}^{-1}\mathcal{R}f](r, \phi) = f(r, \phi) = \frac{1}{2\pi^2} \int_0^{\pi} \int_{-E}^{E} [\mathcal{D}_1\mathcal{R}f](l, \theta) \frac{1}{r\cos(\theta - \phi) - l} dl\, d\theta, \tag{10.5}$$

where $\mathcal{D}_1\mathcal{R}f$ denotes the partial derivative of $\mathcal{R}f$ with respect to its first variable, i.e., for any function $p(l, \theta)$ of two variables

$$[\mathcal{D}_1 p](l, \theta) = \lim_{\Delta l \to 0} \frac{p(l + \Delta l, \theta) - p(l, \theta)}{\Delta l}, \tag{10.6}$$

assuming that the limit exists.

The fully discretized model in the series-expansion approach to the image reconstruction problem of x-ray transmission is formulated in the following way. A Cartesian grid of square picture elements, called *pixels*, is introduced into the region of interest so that it covers the whole picture that has to be reconstructed. The pixels are numbered in some agreed manner, say from 1 (top left corner pixel) to $n$ (bottom right corner pixel) (see Figure 10.2).

The x-ray attenuation function is assumed to take a constant value $x_j$ throughout the $j$th pixel, for $j = 1, 2, \ldots, n$. Sources and detectors that transmit and receive some energy are assumed to be points and the rays between them are assumed to be lines. Further, assume that the length of intersection of the $i$th ray with the $j$th pixel, denoted by $a_{ij}$ for all $i = 1, 2, \ldots, m,\ \ j = 1, 2, \ldots, n$, represents the weight of the contribution of the $j$th pixel to the total attenuation of energy along the $i$th ray.

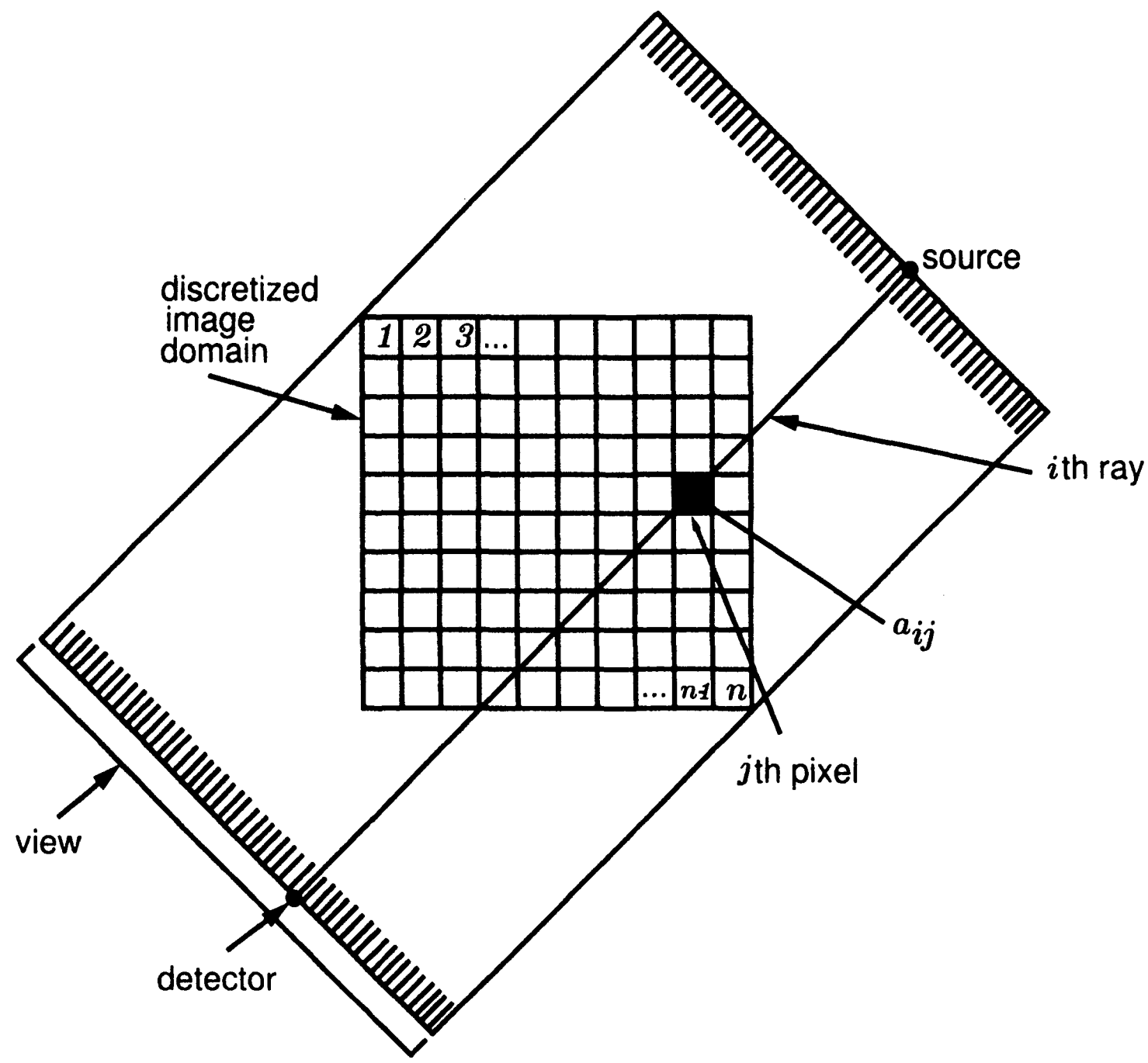

**Figure 10.2** The fully discretized model for transmission tomography image reconstruction.

The physical measurement of the total attenuation along the $i$th ray, denoted by $y_i$, represents the line integral of the unknown attenuation function along the path of the ray. Therefore, in this fully discretized model, the line integral turns out to be a finite sum and the entire model is described by a system of linear equations

$$\sum_{j=1}^{n} x_j a_{ij} = y_i, \qquad i = 1, 2, \ldots, m. \tag{10.7}$$

In matrix notation we write (10.7) as

$$y = Ax, \tag{10.8}$$

where $y = (y_i) \in \mathbb{R}^m$ is the *measurement vector*, $x = (x_j) \in \mathbb{R}^n$ is the *image vector*, and the $m \times n$ matrix $A = (a_{ij})$ is the *projection matrix*.

We may describe the above model in a different way. Let $\{b_j(r, \phi)\}_{j=1}^{n}$

be a set of *basis functions* in polar coordinates in the plane given by

$$b_j(r,\phi) \doteq \begin{cases} 1, & \text{if } (r,\phi) \text{ belongs to the } j\text{th pixel}, \\ 0, & \text{otherwise}, \end{cases} \tag{10.9}$$

and call

$$\hat{f}(r,\phi) \doteq \sum_{j=1}^{n} x_j b_j(r,\phi) \tag{10.10}$$

the *digitization* of $f(r,\phi)$ with respect to the basis functions $\{b_j\}_{j=1}^n$, where $x_j$ are the coefficients of the expansion. Let $\{\mathcal{R}_i\}_{i=1}^m$ be a set of linear and continuous functionals which assign to any picture $f(r,\phi)$ real numbers $\{\mathcal{R}_i f\}_{i=1}^m$. In our case, $\mathcal{R}_i f$ is the line integral of $f(r,\phi)$ along the $i$th ray. Now, $y_i$ is only an approximation of $\mathcal{R}_i f$ because of the inaccuracy in the physical measurements. $\mathcal{R}_i f$ is close to $\mathcal{R}_i \hat{f}$ because of the continuity of $\mathcal{R}_i$ and because we want to use $\hat{f}$ as our approximation for $f$. Using the linearity of $\mathcal{R}_i$, we may then write

$$y_i \simeq \mathcal{R}_i f \simeq \mathcal{R}_i \hat{f} = \sum_{j=1}^{n} x_j \mathcal{R}_i b_j(r,\phi) = \sum_{j=1}^{n} x_j a_{ij}, \tag{10.11}$$

where $a_{ij} \doteq \mathcal{R}_i b_j(r,\phi)$, thus arriving again at (10.7).

The flexibility of the series-expansion approach can be appreciated from this description. First, the functionals $\mathcal{R}_i$ need not be the line integrals we used; they can be any other functionals that may arise in different reconstruction problems if their basic properties make them amenable to a similar modeling process. As an example, we just mention the case of emission computerized tomography (ECT) where the functionals $\mathcal{R}_i$ are the *attenuated line integrals* and the series-expansion approach applied to calculating simultaneously activity and attenuation coefficients leads to a system of nonlinear equations instead of the linear system (10.8) (see Section 10.7 for references). Also, the family of basis functions may be chosen differently than in (10.9).

When considering the fully discretized model in image reconstruction we have to bear in mind some special features of the problem. The system (10.7) is extremely large—with $n$ (number of pixels) and $m$ (number of rays) of the order of magnitude of $10^5$ each—in order to produce images with good resolution. The matrix $A$ of the system is very sparse with less than 1 percent of its entries nonzero, because only few pixels have a nonempty intersection with each particular ray. The system is sometimes underdetermined due to lack of information, but often it is greatly overdetermined in which case it is most probably inconsistent (i.e., for this data vector $y$, there does not exist an $x$ that satisfies $Ax = y$). Moreover we might have reason to believe that the exact algebraic solution of the system, even if it did exist and could be computed accurately, is no more desirable in terms of the reconstruction problem than some other, differently defined "solution." Such a belief may stem from evidence of measurement inaccuracy, or noise corruption of data, or the fact that the original problem has undergone discretization.

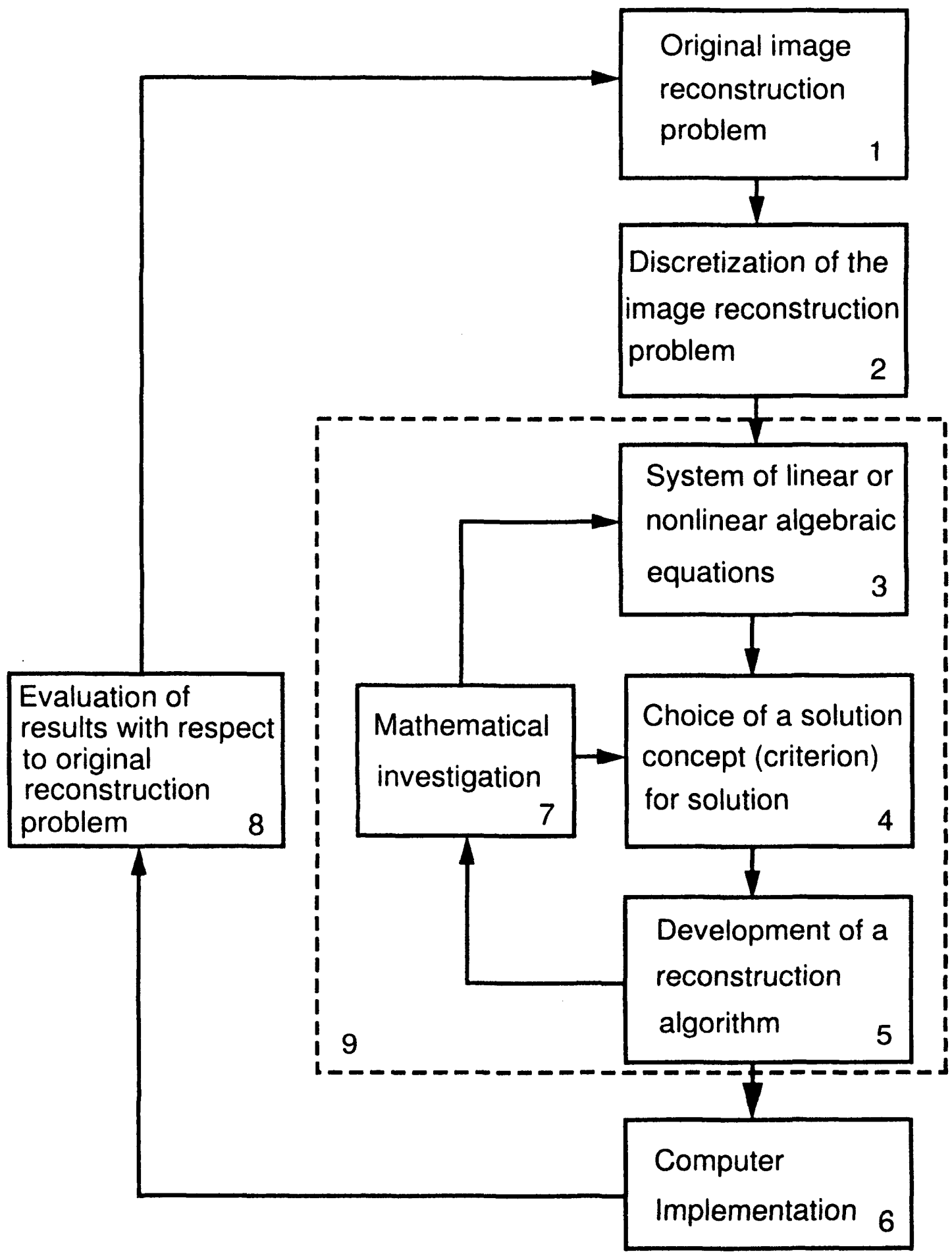

**Figure 10.3** Methodology of the series-expansion approach.

The methodology of the series-expansion approach can be described in general terms by the block-diagram in Figure 10.3. The original image reconstruction problem must first be formulated (box 1). At this stage, the physics of the particular problem (e.g., x-ray transmission, single-photon or positron-annihilation emission, ultrasound, etc.) is taken into account and some idealizing assumptions are usually made. The scanning geometry is also recorded here.

The next step in the series-expansion approach is the discretization (box 2) of the region of interest, and it is here that we depart from the approach of transform methods. In transform methods the problem is handled in its continuous formulation up to the implementation phase (where the necessary approximations are made to the derived inversion formula), whereas here the discretization leads to a system of algebraic equations—linear for x-ray transmission image reconstruction and, possibly, nonlinear for other reconstruction problems (box 3).

In view of the system's special features, such as dimensions, sparsity, inconsistency, ill-conditioning, etc., one has to choose a *solution concept* according to which a solution for the system will be sought (box 4). Most often an optimization criterion is set up, with the system of equations or some system of inequalities derived from it as the set of constraints over which the optimization is performed. This phase is important since the choice here determines the set of mathematical solutions that will be considered as solutions of the underlying image reconstruction problem. Generally, a variety of solution concepts is possible but the final choice rests with the user and should usually be made with reference to the specific reconstruction problem at hand. The principal approaches to choosing a solution concept are the following:

*The feasibility approach*: Here one seeks a solution that lies within a specific vicinity of all hyperplanes defined by the equations of (10.8). This situation is described in Section 6.1 and demonstrated here in Figure 10.4.

*The optimization approach*: An exogenous (i.e., independent of the measurement process) objective function is specified, according to which a particular element will be singled out from the *feasible region* that is a set $Q$ in $\mathbb{R}^n$. This feasible region may be composed of inequalities derived from (10.8) and additional inequalities providing a priori information about the desirable solution. Examples of optimization problems in image reconstruction are discussed later (see Sections 10.2 and 10.3).

*The distance minimization approach*: In this strategy the generalized distance $D_f(y, Ax)$, with respect to some Bregman function $f \in \mathcal{B}(S)$ (see Chapter 2), between the data $y$ and the range $\mathcal{R}(A) = \{Ax \mid x \in \mathbb{R}^n\}$ of the matrix $A$ is minimized. Thus a solution is required for the unconstrained (apart from nonnegativity constraints) optimization problem:

$$\text{Minimize} \quad D_f(y, Ax) \tag{10.12}$$

$$\text{s.t.} \quad x \geq 0. \tag{10.13}$$

For $f(x) = \frac{1}{2} \parallel x \parallel^2$, the distance in (10.12) is $\frac{1}{2} \parallel y - Ax \parallel^2$, yielding the

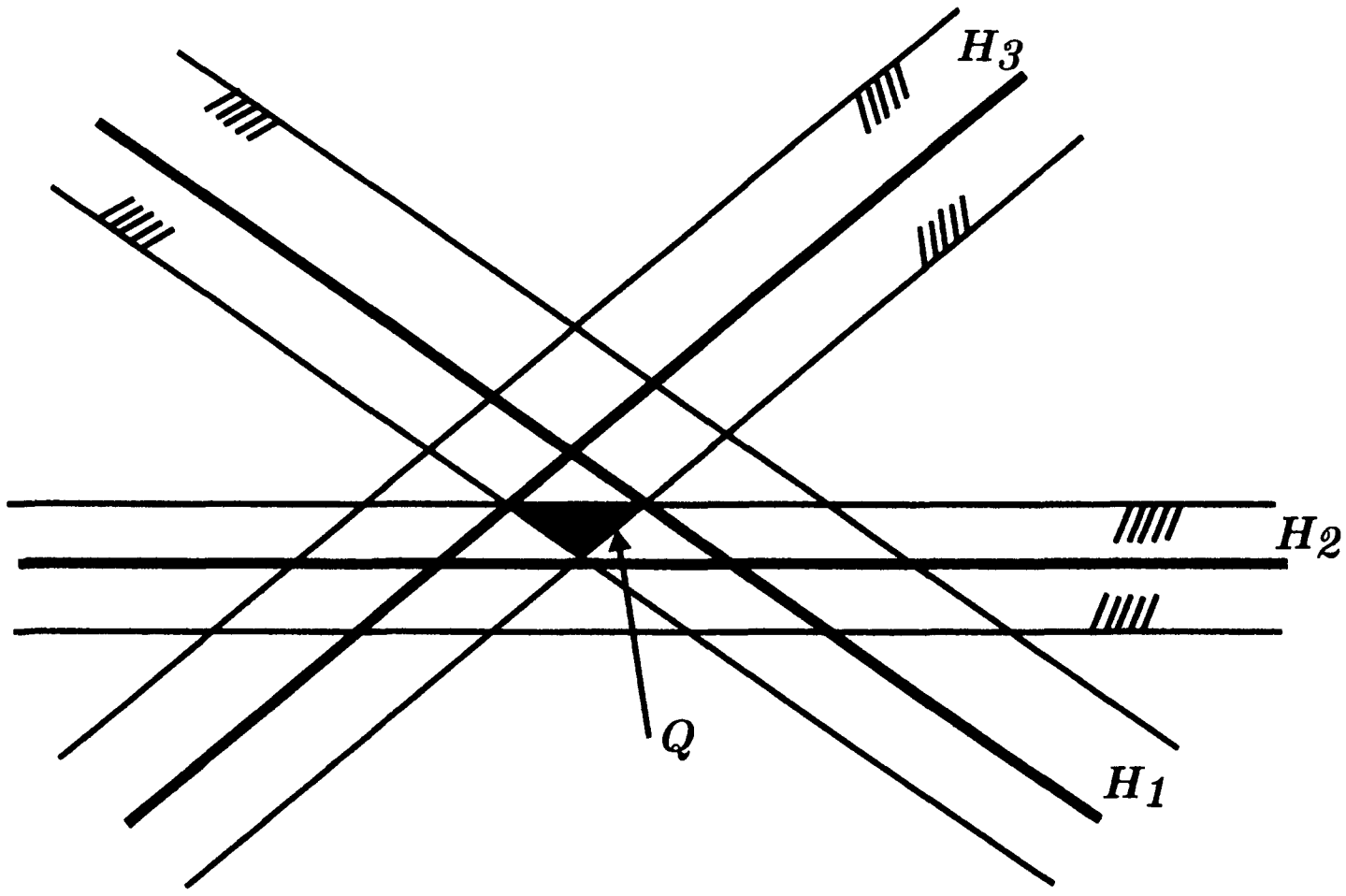

**Figure 10.4** The feasibility approach: the linear system of equations represented by the three hyperplanes $H_1$, $H_2$ and $H_3$ is inconsistent but any point in the intersection $Q$ of the hyperslabs is considered a solution to the image reconstruction problem.

classic least-squares minimization problem. For $f(x) = \sum_{j=1}^{n} x_j \log x_j$, the generalized distance in (10.12) is, up to an additive constant, the Kullback-Leibler information divergence (see Examples 2.1.1 and 2.1.2, respectively).

*The regularization approach*: This approach combines the two previous approaches. The original problem:

$$\text{Solve} \quad Ax = y \tag{10.14}$$

$$\text{s.t.} \quad x \in Q, \tag{10.15}$$

where the linear system is derived from (10.7) and the constraints set $Q \subseteq \mathbb{R}^n$ describes some additionally available information, is treated by interpreting it as

$$\text{Minimize} \quad g(x) + \alpha^2 D_f(y, Ax) \tag{10.16}$$

$$\text{s.t.} \quad x \in Q. \tag{10.17}$$

Here $g(x)$ is a user-chosen objective function, $D_f(y, Ax)$ is as in the previous approach for some user-chosen Bregman function $f \in \mathcal{B}(S)$, $\alpha^2$ is a user-determined regularization parameter, and the unique solution of (10.16)–(10.17) is referred to as the *regularized solution* of (10.14)–(10.15).

The next step is the development of a reconstruction algorithm (box

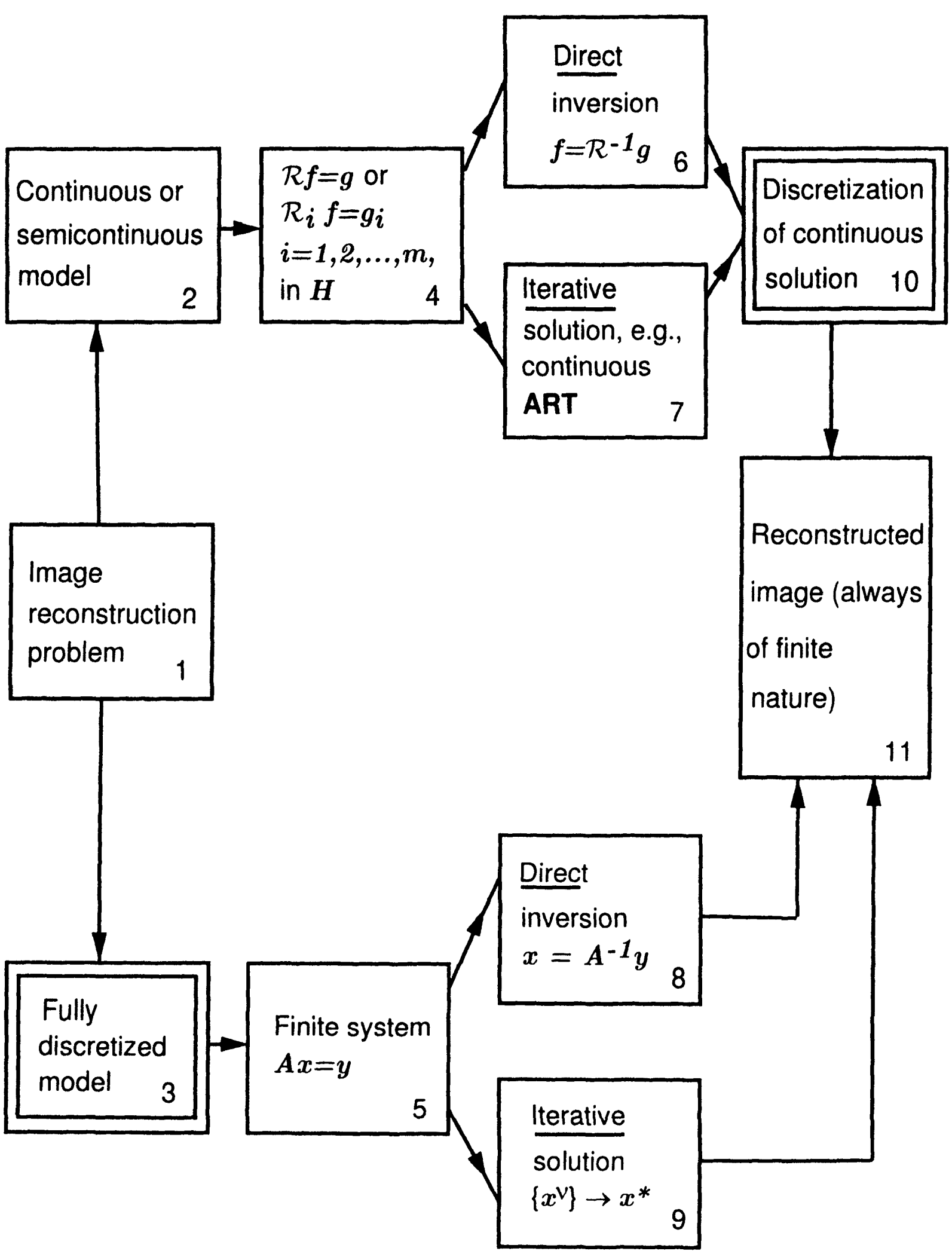

**Figure 10.5** Transform methods versus finite series-expansion methods. The double-framed boxes, 3 and 10, indicate the stage where discretization is introduced in the two approaches. The path connecting boxes 1, 3, 5, 8 or 9, and 11 represents the full discretization approach.

5 in Figure 10.3). Such an algorithm should be capable of handling the problem—which is usually of huge dimensions and very sparse—and it should also be efficiently implementable on a computer in a way that takes advantage of the computer architecture, especially if multiple processors are available.

The results of the computer implementation are then evaluated not only in view of the mathematical problem (box 7) but mainly with respect to the original reconstruction problem (box 8). Figure 10.3 illustrates graphically the steps in the approach. The broken line (box 9) in the figure represents the mathematical efforts involved; the tools for which are found in linear algebra, optimization theory, and numerical analysis.

To get a general idea of the two fundamentally different approaches to the image reconstruction problem we now refer to Figure 10.5. The approach generally termed *transform methods* handles the continuous model $\mathcal{R}f = g$, or a semicontinuous model

$$\mathcal{R}_i f = g_i, \quad i = 1, 2, \ldots, m, \tag{10.18}$$

where $\mathcal{R}_i f = [\mathcal{R}f](l_i, \theta_i)$, i.e., it is the value of $\mathcal{R}f$ at the point $(l_i, \theta_i)$. The model is treated in an appropriate Hilbert space and either direct inversion or an iterative solution method is used (boxes 6 or 7 in Figure 10.5). All variants of the well-known *convolution-backprojection* (CBP) method (also called *filtered backprojection*) would naturally come under the heading of direct inversion. Studies related to the iterative solution approach to the semicontinuous problem (box 7 in Figure 10.5) are also available. Numerical implementation then demands that solution formulae be discretized, often requiring that interpolation be used to fill in values that are needed for the numerical computation.

In contrast to this, the so-called *finite series-expansion* approach to image reconstruction from projections starts off with a fully discretized model (box 3 in Figure 10.5), which leads to a problem of solving, in some sense, a linear or nonlinear system of equations or inequalities in the $n$-dimensional Euclidean space $\mathbb{R}^n$. After a solution concept is chosen an iterative algorithm is usually used, and when iterations are stopped the current iterate $x^\nu \in \mathbb{R}^n$ is taken as the reconstructed image. This approach is represented by the path connecting boxes 1, 3, 5, 9 and 11 in Figure 10.5.

The basic difference between these two approaches lies in whether discretization takes place at the very beginning or at the very end (represented in Figure 10.5 by a double-framed box, i.e., boxes 3 and 10).

## 10.2 A Fully Discretized Model for Positron Emission Tomography

In this section we examine another image reconstruction problem in which fully discretized modeling plays an important role. The problem we ex-amine concerns positron emission tomography (PET), a process whereby a positron-emitting substance is introduced into an organism. Our problem is to determine the concentration of this substance, after some fixed time period, at each point within a certain region of the body. The region of interest is discretized into a finite number of volume elements (*voxels*) throughout which the concentration is assumed constant. Each voxel emits positrons according to a Poisson process at a rate proportionate to the concentration of the substance in the voxel. These rates are the unknowns. An emitted positron combines with an electron in a mutual annihilation event that produces two photons that move in opposite directions along the same straight line (where each such line has a random orientation in space). A ring of detectors is placed around the organism, and it is possible to count the number of events detected in the volume between each pair of detectors (called a tube). Unless the body is completely surrounded by detectors (in which case the ring is a sphere, which is an impractical arrangement, at least in medical tomography) some photons will go undetected (see Figure 10.6).

Let $N_j$, $j = 1, 2, \dots, n$ be the numbers of positrons emitted during a time interval $[0, \tau]$ in voxel $j$. The $N_j$ are independent Poisson variables with means $\mu_j$. Let $p_{ij}$ be the probability of the event that "an emission in voxel $j$ is detected in tube $i$". The $p_{ij}$ depend on the geometry of the system and can be determined from it because the direction of the moving photons has a uniform radial distribution in $\mathbb{R}^3$, that is, the event considered is equivalent to an event in which a line through voxel $j$ intersects the two detectors that constitute tube $i$. Let $N_{ij}$ be the number of events detected in tube $i$, due to photons emitted in voxel $j$, during the time interval $[0, \tau]$. It follows that the $N_{ij}$ are independent Poisson variables with means $\mu_j p_{ij}$. Let $\overline{N}_j$ be the number of emissions in voxel $j$ detected in some tube during the time interval $[0, \tau]$. Then $\overline{N}_j$ are independent Poisson variables with means

$$x_j = \mu_j \sum_{i=1}^{m} p_{ij}. \tag{10.19}$$

Finally, let $y_i$ be the number of photons detected in tube $i$ during the time interval $[0, \tau]$. We have $y_i = \sum_{j=1}^{n} N_{ij}$, so that the $y_i$ are independent Poisson variables with means

$$\xi_i = \sum_{j=1}^{n} \mu_j p_{ij}. \tag{10.20}$$

The $y_i$ and the $p_{ij}$ are the known data, and the $x_j$ from which the $\mu_j$ can be easily obtained using (10.19) are the unknowns. Consider the likelihood function for $\xi$:

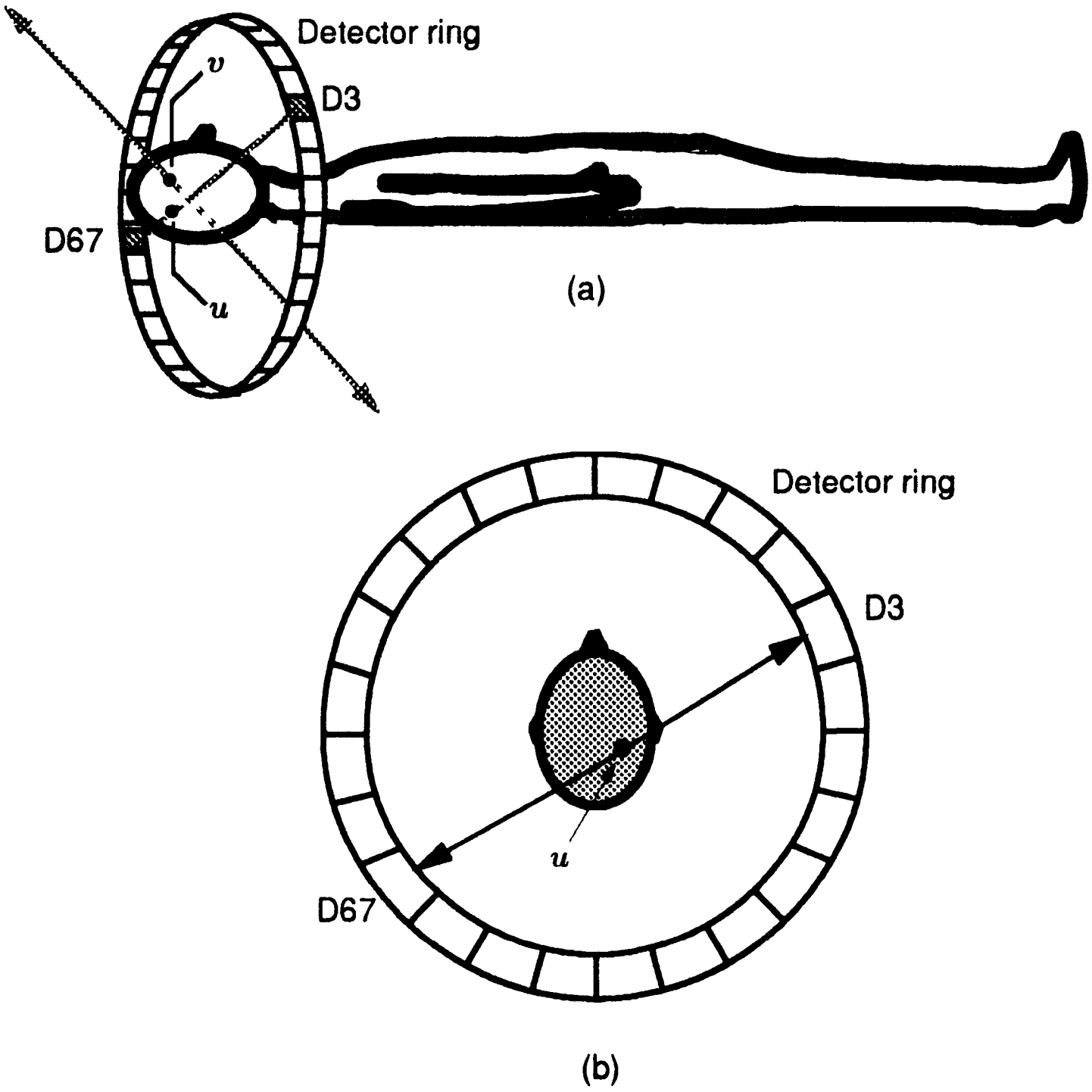

**Figure 10.6** Illustration of PET: (a) Two annihilations: one, at $u$, is detected in the tube (D3, D67) and the other, at $v$, passes undetected because the photon path does not intersect the detector ring. (b) View of the detector ring plane of (a). (Drawn after a figure in Vardi, Shepp, and Kaufman (1985).)

$$L(\xi) \doteq \mathrm{P}(y \mid \xi) = \prod_{i=1}^{m} e^{-\xi_i} \frac{\xi_i^{y_i}}{y_i!}, \tag{10.21}$$

where $y$ and $\xi$ are the vectors with components $y_i$ and $\xi_i$ respectively, the factorial is denoted by ! , and $P(y \mid \xi)$ is the conditional probability of $y$ given $\xi$. Then, the log-likelihood function is

$$\log\ L(\xi) = -\sum_{i=1}^{m} (\xi_i - y_i \log \xi_i + \log (y_i!)). \tag{10.22}$$

Let $\mu$ and $p^i$ be the vectors with components $\mu_j$ and $p_{ij}$, respectively. From (10.20) we have $\xi_i = \langle p^i, \mu \rangle$ so that (10.22) becomes

$$\log\ L(\xi) = -\sum_{i=1}^{m}(\langle p^i, \mu\rangle - y_i\ \log\ \langle p^i, \mu\rangle +\ \log\ (y_i!)). \tag{10.23}$$

The maximum-likelihood estimator of $\xi$ is given by the solution of the problem: max $L(\xi)$, such that $\xi \geq 0$. This maximization problem is equivalent to minimization of $-\log L(\xi)$ such that $\xi \geq 0$, which is equivalent to

$$\text{Minimize} \quad \sum_{i=1}^{m}(\langle p^i, \mu\rangle - y_i \log \langle p^i, \mu\rangle) \tag{10.24}$$

$$\text{s.t.} \quad \mu \geq 0. \tag{10.25}$$

Define

$$a_{ij} \doteq \frac{p_{ij}}{\sum_{i=1}^{m} p_{ij}}, \tag{10.26}$$

so that $\sum_{i=1}^{m} a_{ij} = 1$. Let $x$ and $a^i$ be the vectors with components $x_j$ and $a_{ij}$, respectively. It follows from (10.19) and (10.26) that $\langle p^i, \mu\rangle = \langle a^i, x\rangle$. So problem (10.24)–(10.25) is equivalent to

$$\text{Minimize} \quad \sum_{i=1}^{m}(\langle a^i, x\rangle - y_i \log \langle a^i, x\rangle) \tag{10.27}$$

$$\text{s.t.} \quad x \geq 0, \tag{10.28}$$

where $0 \leq a_{ij} \leq 1$ for all $i, j$, and $\sum_{i=1}^{m} a_{ij} = 1$, and $a^i \neq 0$ for all $i$. Since the function $\phi(t) \doteq t - y_i \log t$ attains its minimum at $t = y_i$, it follows that, if there exists an $x \geq 0$ such that $Ax = y$ (where $A$ is the matrix with rows $a^i$) then such an $x$ is a solution of (10.27)–(10.28). So (10.27)–(10.28) can be seen as a generalization of $Ax = y$, $x \geq 0$. If this system is feasible its solutions are the same as those of (10.27)–(10.28). If it is infeasible solutions of (10.27)–(10.28) can be considered to be approximate solutions of the linear systems in the sense that they minimize $D_f(y, Ax)$ where $D_f$ is the Kullback-Leibler information divergence (see Example 2.1.2).

This model is an approximate representation of the physics of positron emission tomography (PET). In reality, positrons travel some distance before hitting an electron; the angle between emitted photons is not exactly $180^o$; and some photons are absorbed by the tissue surrounding the voxel where they were emitted and therefore do not reach the detection ring. To incorporate these and other factors the model needs to be extended.

### 10.2.1 The expectation-maximization algorithm

The *expectation-maximization* (EM) algorithm is an iterative technique for computing maximum-likelihood estimators. We discuss here its implementation for the model discussed in Section 10.2, which deals with a finite number of independent Poisson variables.

Suppose that we have an approximation $x^\nu$ for $x$, that is, $x_j^\nu$ is an approximation of the expected value $E[\overline{N}_j]$. Since $\overline{N}_j = \sum_{i=1}^m N_{ij}$ and since $y_i = \sum_{j=1}^n N_{ij}$, we consider the conditional expectation value

$$E[\overline{N}_j \mid y] = \sum_{i=1}^{m} E[N_{ij} \mid y] = \sum_{i=1}^{m} E[N_{ij} \mid y_i], \tag{10.29}$$

using the independence of the $y_i$ in the last equality. From (10.29)

$$E[\overline{N}_j \mid y] = \sum_{i=1}^{m} E\left[N_{ij} \mid \sum_{j=1}^{n} N_{ij} = y_i\right]. \tag{10.30}$$

We take this expression to be the next iterate $x_j^{\nu+1}$. However, if (as assumed above) we let $x_j^\nu = E[\overline{N}_j]$ then $E[N_{ij}]$ would be equal to $a_{ij}x_j^\nu$. Given independent Poisson variables $X_1, \dots, X_n$ with means $\zeta_1, \dots, \zeta_n$, the random variable $X_j$ for which $\sum_{l=1}^n X_l = \omega$ is a binomial random variable, i.e., its mathematical expectation is given by

$$E\left[X_j \mid \sum_{l=1}^{n} X_l = \omega\right] = \frac{\omega \zeta_j}{\sum_{l=1}^{n} \zeta_l}. \tag{10.31}$$

In our case, taking $x_j^\nu$ as an approximation for $E[\overline{N}_j]$ and using (10.31) in (10.30) gives

$$x_j^{\nu+1} = E[\overline{N}_j \mid y] \simeq \sum_{i=1}^{m} \frac{y_i x_j^\nu a_{ij}}{\langle a^i, x^\nu \rangle} = x_j^\nu \sum_{i=1}^{m} \frac{y_i a_{ij}}{\langle a^i, x^\nu \rangle}. \tag{10.32}$$

This version of the EM algorithm can be applied in the case of PET. A relaxed version of the iterative step (10.32) for the solution of the same problem is given by

$$x_j^{\nu+1} = x_j^\nu \left( \sum_{i=1}^{m} \frac{y_i a_{ij}}{\langle a^i, x^\nu \rangle} \right)^{\alpha}, \tag{10.33}$$

where $\alpha$ is a relaxation parameter. Algorithm (10.32) converges to a solution of (10.27)–(10.28) from any positive starting vector $x^0 \in \mathbb{R}^n_{++}$. The relaxed EM algorithm defined by (10.33) converges to a solution of (10.27)–(10.28) for $\alpha \in (0, 2)$. If $\alpha \in (0, 1]$ convergence holds from any $x^0 \in \mathbb{R}^n_{++}$; if $\alpha \in (1, 2)$ then the algorithm converges only locally.

## 10.3 A Justification for Entropy Maximization in Image Reconstruction

Going back to the fully discretized model of transmission computerized tomography, suppose that we are interested in the solution of the system

$$Ax = y, \quad x \geq 0, \tag{10.34}$$

with $A \in \mathbb{R}^{m \times n}, y \in \mathbb{R}^m$, and assume that $y$ is a vector of measurements subject to experimental errors. System (10.34) might have no solutions or infinitely many solutions, and we want to find the most "likely" one, in some sense to be made precise later on. We start by replacing $Ax = y$ with

$$Ax = y + r, \tag{10.35}$$

where $r$ is a random error vector in $\mathbb{R}^m$. We consider also $x$ and $y$ as random vector variables, and thus the most likely solution is the one that maximizes $P(x \mid y)$, i.e., the conditional probability of $x$ given the measurements $y$ where $x$ and $y$ are related through (10.35). A good point of departure is Bayes' formula

$$P(x \text{ and } y) = P(x \mid y)P(y) = P(y \mid x)P(x), \tag{10.36}$$

where $P(x)$ and $P(y)$ are the a priori probabilities in the usual Bayesian notation. From (10.36) we get

$$P(x \mid y) = \frac{P(y \mid x)P(x)}{P(y)}, \tag{10.37}$$

and our desired maximization of $P(x \mid y)$ is equivalent by (10.37) to

$$\max_x P(y \mid x)P(x). \tag{10.38}$$

We must therefore specify $P(y \mid x)$ and $P(x)$. Let us look first at $P(x)$. We view the pixels $j = 1, 2, \ldots, n$ in the fully discretized model (see Figure 10.2) as containers among which a total of $N$ particles are to be distributed, with $N_j$ particles falling in pixel $j$, so that $N = \sum_{j=1}^n N_j$. We will also assume that for all $j$ the $j$th pixel is divided into $m_j$ subpixels, so that the total number of subpixels is given by $M = \sum_{j=1}^n m_j$.

Consider now two possible configurations: fine, corresponding to the distribution of particles in the subpixels and coarse, corresponding to the distributions of particles in the pixels. The total numbers of possible fine and coarse configurations are clearly $N^M$ and $N^n$, respectively. We make the following assumptions on the probabilities of these configurations.

**Assumption 10.3.1**

(*i*) *All fine configurations are equally likely.*

(*ii*) $x_j$ *is proportional to the fraction of all particles that fall in pixel* $j$*, i.e.,* $x_j = \alpha N_j/N$*, where* $\alpha$ *is a scale factor.*

It is clear that the number of coarse configurations such that exactly $N_j$ particles fall in pixel $j$, i.e., those which produce a coarse distribution $\overline{N} = (N_1, \ldots, N_n)$, is given by the multinomial coefficient

$$\binom{N}{N_1, \ldots, N_n} = \frac{N!}{\prod_{j=1}^{n} N_j!}. \tag{10.39}$$

Also, the number of fine configurations in pixel $j$ is $m_j^{N_j}$. Therefore, the number of fine configurations compatible with a given coarse distribution $\overline{N}$ is the number $V(\overline{N})$ given by

$$V(\overline{N}) = \frac{N!}{\prod_{j=1}^{n} N_j!} \prod_{j=1}^{n} m_j^{N_j} = N! \prod_{j=1}^{n} \frac{m_j^{N_j}}{N_j!}. \tag{10.40}$$

By Assumption 10.3.1($i$), the probability of coarse distribution $\overline{N}$ is given by

$$P(\overline{N}) = \frac{V(\overline{N})}{M^N}. \tag{10.41}$$

Next we define, for all $j = 1, 2, \ldots, n$, $z_j \doteq N_j/N$, $q_j \doteq m_j/M$, and $u_j \doteq \alpha q_j$. Using $q_j$ and $z_j$ in (10.40), taking logarithms in (10.41), and applying Stirling's approximation formula (i.e., $N! \approx (N/e)^N$) we get

$$\log\left(P(\overline{N})\right) = -N \sum_{j=1}^{n} z_j \log\left(\frac{z_j}{q_j}\right). \tag{10.42}$$

Note that, from the definitions of $M$ and $N$ and from Assumption 10.3.1($ii$), we have

$$\sum_{j=1}^{n} x_j = \sum_{j=1}^{n} u_j = \alpha. \tag{10.43}$$

Since $x_j = \alpha z_j$ (10.42) can be rewritten as

$$\log P(x) = -\frac{N}{\alpha} \sum_{j=1}^{n} x_j \log\left(\frac{x_j}{u_j}\right). \tag{10.44}$$

Next we must specify $u_j$, which is proportional to the number of subpixels in pixel $j$. If we have an estimate $\overline{x} = (\overline{x}_j) \in \mathbb{R}^n$ of the solution, then it is reasonable to take $u_j \doteq \overline{x}_j$. In the absence of such information, the best option is to take all the $u_j$'s equal, i.e., $u_j = \alpha/n$ so that all pixels contain the same number of subpixels. Note that the maximum of the expression in (10.44) is attained at $x = u$, but $P(x)$ is the prior probability of $x$ ignoring the information given by $y$. In order to specify $P(y \mid x)$ we need to make the following assumption on the error vector $r = (r_i) \in \mathbb{R}^m$.

**Assumption 10.3.2** *The $r_i$'s are independent and normally distributed with mean zero and variance $\sigma$, i.e., for every $i = 1, 2, \ldots, m$,*

$$P(r_i) \doteq \beta \exp\left(\frac{-r_i^2}{2\sigma^2}\right), \tag{10.45}$$

*where $\beta$ is some scaling factor.*

By (10.35) and (10.45) we have

$$P(y \mid x) = \beta \exp\left(-\frac{\| Ax - y \|^2}{2\sigma^2}\right), \tag{10.46}$$

and combining (10.44) and (10.46) gives

$$\log\left(P(y \mid x)P(x)\right) = -\frac{N}{\alpha}\sum_{j=1}^{n} x_j \log\left(\frac{x_j}{u_j}\right) - \frac{1}{2\sigma^2} \| Ax - y \|^2 - \log\beta, \tag{10.47}$$

so that the solution vector $x^*$, which maximizes $P(x \mid y)$, is the vector that minimizes the function

$$f(x) \doteq \sum_{j=1}^{n} x_j \log\left(\frac{x_j}{u_j}\right) + \mu \| Ax - y \|^2, \tag{10.48}$$

where $\mu \doteq \alpha/(2N\sigma^2)$. The objective function $f$ in (10.48) is strictly convex and coercive (i.e., $\lim_{\|x\| \to \infty} f(x)/ \| x \| = +\infty$) so that (10.48) has a unique minimizer. Reintroducing $r$, the minimization of (10.48) can be written as

$$\underset{x,r}{\text{Minimize}} \quad \sum_{j=1}^{n} x_j \log\left(\frac{x_j}{u_j}\right) + \mu \| r \|^2 \tag{10.49}$$

$$\text{s.t.} \quad Ax = y + r, \tag{10.50}$$

$$x \geq 0. \tag{10.51}$$

Consider now the Kullback-Leibler information divergence, defined by (2.9). In view of (10.43), it can be written as $D(x, u) = \sum_{j=1}^{n} x_j \log(x_j/u_j)$, so that (10.49) is equivalent to

$$\text{Minimize } D(x,u) + \mu \parallel r \parallel^2 . \tag{10.52}$$

As mentioned before, in the absence of a good estimate for the solution of the problem we should take $u_j = \alpha/N$, and thus in view of (10.43) we get $D(x,u) = \sum_{j=1}^{n} x_j \log x_j - \alpha \log(\alpha/N)$, so that (10.52) is equivalent to

$$\text{Minimize } \sum_{j=1}^{n} x_j \log x_j + \mu \parallel r \parallel^2 . \tag{10.53}$$

If the error is known to be small (i.e., if $\sigma$ is small) then the second term in the objective of (10.49) is negligible and (10.49)–(10.51) becomes

$$\text{Minimize } \sum_{j=1}^{n} x_j \log\left(\frac{x_j}{u_j}\right) \tag{10.54}$$

$$\text{s.t. } Ax = y, \tag{10.55}$$

$$x \geq 0, \tag{10.56}$$

where $u$ is our initial estimate. In the absence of such an estimate we may assume that $u_j = 1$ for all $j = 1, 2, \ldots, n$, and the problem reduces to

$$\text{Minimize } \sum_{j=1}^{n} x_j \log x_j \tag{10.57}$$

$$\text{s.t. } Ax = y, \tag{10.58}$$

$$x \geq 0. \tag{10.59}$$

We remark that MART (Algorithm 6.6.2) is specifically devised to solve such linearly constrained entropy optimization problems. If we want to keep the error term as in (10.49)–(10.51), we can still use Bregman's method, (Algorithms 6.3.1 or 6.8.1) because the objective function in (10.49) is a Bregman function with zone $\mathbb{R}^n_{++} \times \mathbb{R}^m$; but no MART-type hybrid algorithm is known for solving the resulting mixed entropic-quadratic optimization problem.

## 10.4 Algebraic Reconstruction Technique (ART) for Systems of Equations

In this section we review the *variable-block Algebraic Reconstruction Technique* (ART) algorithmic scheme for solving systems of linear equations that arise from the fully discretized model of transmission tomography. Mathematically speaking, this scheme is a special instance of the more general *block-iterative projections* (BIP) method (see Section 5.6). However, the variable-block ART scheme deserves to be examined separately because of its importance to image reconstruction from projections, and it includes the classical row-action ART, the fixed-block ART (i.e., block-Kaczmarz algorithm), as well as the Cimmino version of SIRT (Simultaneous Iterative Reconstruction Technique), as special cases. It serves as a unifying framework for the description of all these popular iterative reconstruction techniques and it enhances them in that it allows the processing of blocks (i.e., groups of equations) that need not be fixed in advance, but may rather change dynamically throughout iterations. The number of blocks, their sizes, and their structure may all vary provided that the weights attached to the equations do not vanish.

How iterative reconstruction algorithms behave when the underlying system of equations is inconsistent is an interesting question because noise and other inaccuracies in data would usually make any consistency assumption unrealistic. Results related to the behavior of ART, fixed-block ART, and SIRT when applied to an inconsistent system are available in the literature, but the question is still open, in general, for the variable-block ART. Another important issue is the implementation and practical performance of block-iterative reconstruction algorithms with fixed or variable blocks. While research on the algorithms in terms of quality of reconstructed images is important there is also the issue of their potential use on parallel machines, which as yet has not been fully explored.

Consider a system of linear equations obtained from the fully discretized model for transmission tomography image reconstruction, which has the form (10.8). Unless otherwise stated we assume that the system is consistent, i.e., given $A$ and $y$, the set $\{x \in \mathbb{R}^n \mid Ax = y\}$ is nonempty.

Let $w = (w(i)) \in \mathbb{R}^m$ be a *weight vector* with $w(i) \geq 0$ for all $i$, and $\sum_{i=1}^m w(i) = 1$. Recall (Section 5.6) that a sequence $\{w^\nu\}_{\nu=0}^\infty$ of weight vectors is called *fair* if for every $i$ there exist infinitely many values of $\nu$ for which $w^\nu(i) > 0$. This guarantees that no equation is given a zero weight indefinitely. To prevent equations from fading away by positive but steadily diminishing weights we must require also that for every $i$ the stronger condition

$$\sum_{\nu=0}^{\infty} w^\nu(i) = +\infty, \tag{10.60}$$

holds (see Theorem 5.6.1).

The sequence of weight vectors $\{w^\nu\}$ actually determines which block is used at each iteration by attaching positive weights to equations that belong to the current block and zero weights to the rest.

### Algorithm 10.4.1 Variable-block ART—The General Scheme

**Step 0:** (Initialization.) $x^0 \in \mathbb{R}^n$ is arbitrary.

**Step 1:** (Iterative step.)

$$x^{\nu+1} = x^{\nu} + \lambda_{\nu} \sum_{i=1}^{m} w^{\nu}(i) \left( \frac{y_i - \langle a^i, x^{\nu} \rangle}{\| a^i \|^2} \right) a^i. \tag{10.61}$$

**Relaxation parameters:** $\{\lambda_{\nu}\}_{\nu=0}^{\infty}$ is a sequence of relaxation parameters chosen by the user such that $\epsilon_1 \leq \lambda_{\nu} \leq 2 - \epsilon_2$ for all $\nu \geq 0$, with some arbitrary but fixed $\epsilon_1, \epsilon_2 > 0$.

As a direct corollary of Theorem 5.6.1 we obtain:

**Theorem 10.4.1** *If the system of equations in (10.8) is consistent and $\{w^{\nu}\}$ are weight vectors with property (10.60), then any sequence $\{x^{\nu}\}$ generated by Algorithm 10.4.1 converges to a solution of the system (10.8).*

The variable-block ART algorithm allows processing of the information contained in groups of equations, called blocks. The number of blocks, their sizes (i.e., the number of equations in each block), and their specific structure (i.e., which equations are assigned to each block) may all vary from one iterative step to the next. The following special cases of Algorithm 10.4.1 have been studied separately in the past.

*Row-action ART* This classic sequential ART is obtained by choosing the weight vectors as $w^{\nu} = e^{i(\nu)}$ where $e^l \in \mathbb{R}^n$ is the $l$th standard basis vector (having one in its $l$th coordinate and zeros elsewhere). Each block contains a single equation and the index sequence $\{i(\nu)\}$ with $i \leq i(\nu) \leq m$ for all $\nu$ is the *control sequence* of the algorithm that determines the index of the single equation upon which the algorithm operates at the $\nu$th iterative step. The iterative steps take the form

$$x_j^{\nu+1} = x_j^{\nu} + \lambda_{\nu} \frac{b_{i(\nu)} - \langle a^{i(\nu)}, x^{\nu} \rangle}{\| a^{i(\nu)} \|^2} a_j^{i(\nu)}, \text{ for all } j = 1, 2, \ldots, n. \tag{10.62}$$

Depending on the initialization vector $x^0$ this is either Algorithm 5.4.3 or Algorithm 6.5.1.

*Cimmino-type SIRT* This is a fully simultaneous reconstruction algorithm in which all equations are lumped, with a fixed system of weights, into a single block and are acted upon simultaneously in every step. The iteration formula is precisely (10.61) with $w^{\nu} = w$ for all $\nu \geq 0$, for some fixed $w$ and $w(i) \neq 0$ for all $i = 1, 2, \ldots, m$.

*Fixed-block ART* Here the index set $I = \{1, 2, \ldots, m\}$ is partitioned as $I = I_1 \cup I_2 \cup \cdots \cup I_M$ into $M$ blocks. $\{t(\nu)\}$ is a control sequence over the set $\{1, 2, \ldots, M\}$ of block indices and the weight vectors are of the form $w^{\nu} = \sum_{i \in I_{t(\nu)}} w^{\nu}(i) e^i$. This guarantees that equations outside the $t(\nu)$th block are not operated upon in the $\nu$th iterative step. The block-Kaczmarz procedure is thus obtained.

In addition to these well-known special cases, the variable-block ART enables the implementation of dynamic ART algorithms in which the block-formation strategy might vary during the iterative process itself. Such variations in the block-formation strategy may be used to accelerate the initial convergence toward an acceptable reconstructed image. It is possible under the regime of the variable-block ART to perform *multilevel image reconstruction.* By multilevel approaches we mean methods whereby the system of equations to be solved and/or the way the system is organized into blocks can be changed repeatedly during the iterative process.

Because of their computational appeal (simplicity, row-action nature, etc.) and their efficacious performance in some specific situations, ART methods can also be applied to inconsistent systems of equations, i.e., systems for which there is no common point to all hyperplanes represented by the equations of the system. Their surprisingly good practical performance is confirmed by convergence results on the behavior of the algorithms when applied to such systems. For example, the row-action Kaczmarz algorithm is known to be *cyclically convergent*, which means that the subsequences of iterates lying on each of the hyperplanes of the system converge separately although the whole sequence does not. This also holds for the fixed-block Kaczmarz method. The fully simultaneous Cimmino method also converges locally to a weighted least-squares solution if the system is inconsistent. As yet no study exists that in any way synthesizes these results into a definitive statement on the behavior of the variable-block ART method for inconsistent systems.

## 10.5 Iterative Data Refinement in Image Reconstruction

Measurement is fundamental to many scientific disciplines; however the measuring device often provides only approximate data. The measured data will be called *actual data*; and data which are to some extent idealized, i.e., assumed to obey certain idealizing assumptions, will be called *ideal data.* The discrepancy between the actual data and the ideal data can sometimes be estimated from the actual data (based on knowledge of the measuring process), leading to a better approximation of the ideal data. This new approximation can then be used to estimate the new discrepancy, and the process can be repeated. Typically, our measurement tools are insufficient to obtain the ideal data exactly, but the original discrepancy may be significantly reduced with a few iterative steps. We term this process *iterative data refinement* (IDR) and examine it in this section. We discuss first the fundamentals, then we study applications to medical imaging (see Section 10.5.2).

### 10.5.1 The fundamentals of iterative data refinement

In Figure 10.7 we introduce the IDR approach and notation. The object that is measured is denoted by $x$. We have a *recovery operator* $S$, which from the data $\hat{y}$ that would be produced by an idealized measuring device produces an $\hat{x}$ that is an acceptable approximation of $x$. Our actual measuring device produces data $\tilde{y}$, and the overall aim is to refine this $\tilde{y}$ to bring it closer to the ideal $\hat{y}$ before applying the recovery operator $S$. The operator $U$ (located in the box outlined with broken lines in Figure 10.7) is intended to provide an approximation to $\hat{y} - \mu\tilde{y}$, when applied to an approximation of $\hat{y}$ (here $\mu$ is a constant, which is sometimes referred to as the *relaxation parameter*). Let $y^0$ be an initial approximation to $\hat{y}$. The hope is that $y^1 = \mu\tilde{y} + Uy^0$ is a better approximation to $\hat{y}$ than $y^0$, and that $y^2 = \mu\tilde{y} + Uy^1$ is a better approximation to $\hat{y}$ than $y^1$, and so on. One way of obtaining such an operator $U$ is detailed inside the broken-line box in Figure 10.7. There the recovery operator $S$ is applied to our current approximation $y^\nu$ of $\hat{y}$ to produce an approximation $x^\nu$ of $\hat{x}$. Then the idealized and actual measuring devices are simulated by applying operators $A$ and $B$ respectively, to $x^\nu$. We use $Ax^\nu - \mu Bx^\nu$ as an estimate of $\hat{y} - \mu\tilde{y}$. Thus in this case $U \doteq (A - \mu B)S$.

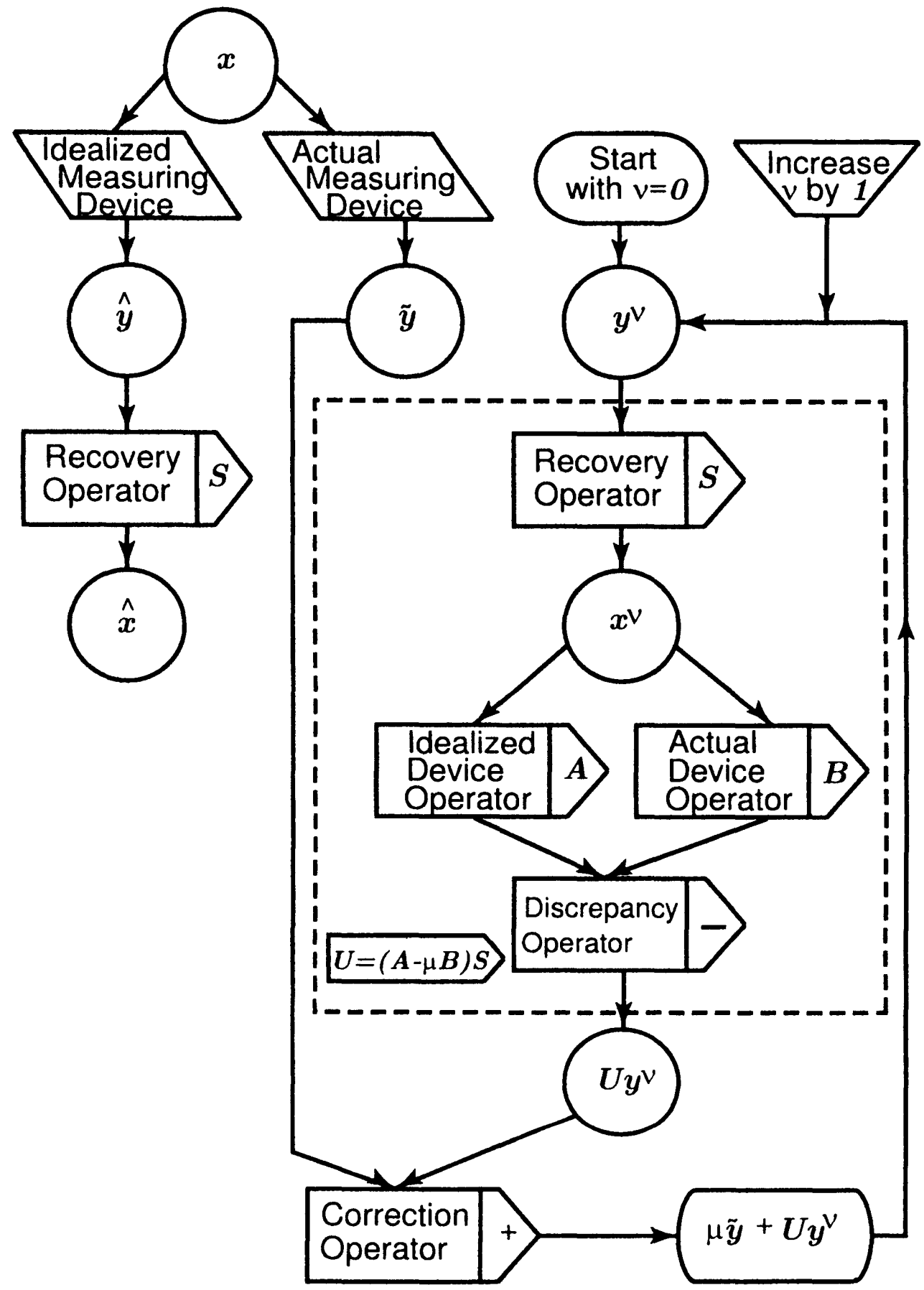

**Figure 10.7** Outline of iterative data refinement (IDR).

More formally, we think of the normed space $Y$ as the space of all possible measurements (both actual and ideal), and of $\| y' - y \|$ as the distance between measurements $y'$ and $y$ in $Y$. The iterative data refinement (IDR) algorithm associates with every 4-tuple $(\tilde{y}, y^0, U, \mu)$ a sequence $\{y^\nu\}_{\nu=0}^\infty$ generated by

$$y^{\nu+1} = \mu\tilde{y} + Uy^\nu, \tag{10.63}$$

where $\tilde{y}$ and $y^0$ are elements of $Y$, $U$ is an operator from $Y$ into $Y$, and $\mu$ is a scalar, A natural (although not necessary) choice for the initial iterate $y^0$ would be the actual data $\tilde{y}$.

The purpose of this process is the following. The actual measurement $\tilde{y}$ is an approximation to some ideal measurement $\hat{y}$. The operator $U$ is supposed to provide an approximation of $\hat{y} - \mu\tilde{y}$. We expect this approximation to be good, provided the argument of $U$ is near $\hat{y}$.

In practical situations IDR does not necessarily produce a convergent sequence. It is nevertheless enlightening to recall first its behavior in the much-studied convergent case.

**Proposition 10.5.1** *If the sequence $\{y^\nu\}$ produced by IDR for the 4-tuple $(\tilde{y}, y^0, U, \mu)$ converges to a $y^*$ such that $U$ is continuous at $y^*$ and if*

$$Uy^* = \hat{y} - \mu\tilde{y}, \tag{10.64}$$

*then $y^* = \hat{y}$.*

**Proof** By letting $\nu$ on both sides of equation (10.63) go to infinity we get

$$y^* = \mu\tilde{y} + Uy^*. \tag{10.65}$$

This and equation (10.64) imply that $y^* = \hat{y}$. ∎

This proposition says that if the sequence converges at all and it behaves as expected near the limit point, then this limit is indeed the desired point. An essential problem associated with IDR is to find an appropriate operator $U$ that ensures convergence. In many practical problems such an operator does not exist. One way of ensuring convergence is to assume that $U$ is linear and has bounded norm. This leads to the following well-known special case.

**Proposition 10.5.2** *If $U$ is a linear operator with bounded norm $\| U \| < 1$, then IDR converges to the point*

$$y^* = \mu(I - U)^{-1}\tilde{y}, \tag{10.66}$$

*regardless of $y^0$, where $I$ is the identity operator and $(I-U)^{-1}$ denotes the inverse operator of $(I-U)$. If in addition equation (10.64) is satisfied, then $y^* = \hat{y}$.*

**Proof** If $U$ is a linear operator, then it is easy to prove by induction that for all $\nu \geq 0$,

$$y^\nu = \mu \sum_{i=0}^{\nu-1} U^i \tilde{y} + U^\nu y^0, \tag{10.67}$$

where $U^i$ is the $i$th power of $U$, i.e., $U^i = UU^{i-1}$ for all $i \geq 0$, and $U^0 = I$, the identity. If $\| U \| < 1$ then (see, e.g., Kolmogorov and Fomin (1957, Theorem 8, p.102)

$$(I-U)^{-1} = \sum_{i=0}^{\infty} U^i. \tag{10.68}$$

Also, the second term on the right-hand side of equation (10.67) goes to zero as $\nu$ goes to infinity. This, combined with equation (10.68), proves equation (10.66). ■

In many practical situations $U$ may be nonlinear and satisfaction of equation (10.64) cannot be guaranteed; thus we cannot expect convergence to $\hat{y}$. Instead, we aim at a much weaker goal of merely reducing the distance to $\hat{y}$.

An operator $U$: $Y \to Y$ is called a *contraction* on $Y$ if there exists a number $0 \leq \beta < 1$ such that

$$\| Uy' - Uy \| \leq \beta \| y' - y \|, \tag{10.69}$$

for all $y'$ and $y$ in $Y$. To describe precisely the behavior of IDR in practical situations, we introduce the notion of *targeted contraction.*

**Definition 10.5.1 (Targeted Contraction)** *Given any pair of points $u^* \in Y$ and $y^* \in Y$, an operator $U$ from $Y$ into $Y$ is called a targeted contraction with respect to the pair $(u^*, y^*)$ on a set $\Omega \subseteq Y$ if there exists a number $\beta$, $0 \leq \beta < 1$, such that*

$$\| Uy - u^* \| \leq \beta \| y - y^* \|, \tag{10.70}$$

*for all $y \in \Omega$.*

This definition means that $\Omega$ consists of all points $y \in Y$ with the property that, as they get closer to $y^*$, their images under $U$ (i.e., the points $Uy$) tend to $u^*$. Using $y^* = \hat{y}$ and $u^* = \hat{y} - \mu\tilde{y}$ in the inequality (10.70) we may define for every real number $\beta$, $0 \leq \beta < 1$, the set

$$N_\beta \doteq \{y \mid \| \hat{y} - \mu\tilde{y} - Uy \| \leq \beta \| \hat{y} - y \|\}. \tag{10.71}$$

Then if a sequence $\{y^\nu\}$ is produced by IDR for the 4-tuple $(\tilde{y}, y^0, U, \mu)$, and if $y^\nu \in N_\beta$ for some $0 \leq \beta < 1$, then

$$\| \hat{y} - y^{\nu+1} \| \leq \beta \| \hat{y} - y^\nu \| . \tag{10.72}$$

This means that as long as the points $y^0, y^1, y^2, \ldots$ belong to $N_\beta$, for some $0 \leq \beta < 1$, each element of the sequence is nearer to the ideal measurement $\hat{y}$ than the previous element. This simple observation is the backbone of practical IDR. The choice of $U$ that will make $y^0, y^1, y^2, \ldots$ behave as described above is not a trivial matter, especially in the case of ill-posed problems. The question of how much constraining or a priori knowledge should be incorporated into the process to achieve a reasonable output is very much application dependent.

In typical applications of IDR, the measurements $y$ serve as an input to a (often linear) procedure $S$, which is supposed to produce a result $x$. To formalize this, let the normed space $X$ denote a space of all possible results, and let $S$ denote an operator, usually linear, from $Y$ into $X$.

We think of $\hat{x} = S\hat{y}$ as the ideal result and assume that, provided $y$ is near $\hat{y}$, $Sy$ is near $\hat{x}$. Suppose now that we have a not necessarily linear operator $T$ from $X$ into $Y$, which when applied to an approximation of $\hat{x}$ is supposed to provide an approximation of $\hat{y} - \mu\tilde{y}$. Then we see that the composition $TS$ can serve as the operator $U$. Furthermore, the IDR procedure leads to the following procedure for *iterative result refinement* (IRR) described in next two propositions, which are reformulations of the corresponding statements for IDR.

**Proposition 10.5.3** *Let $\{y^\nu\}$ be the sequence produced by the IDR algorithm for the 4-tuple $(\tilde{y}, y^0, TS, \mu)$ where $S$ is a linear operator from $Y$ into $X$, and let for all $\nu \geq 0$,*

$$x^\nu = Sy^\nu. \tag{10.73}$$

*Also, let $\tilde{x} = S\tilde{y}$. Then, for all $\nu \geq 0$ we have*

$$x^{\nu+1} = \mu\tilde{x} + STx^\nu. \tag{10.74}$$

If we consider the elements $\tilde{x} = S\tilde{y}$ and $\hat{x} = S\hat{y}$ of $X$ and the scalar $\mu$ fixed, we may define for any $\alpha$, $0 \leq \alpha < 1$ the set

$$M_\alpha \doteq \{x \mid \| \hat{x} - \mu\tilde{x} - STx \| \leq \alpha \| \hat{x} - x \|\}. \tag{10.75}$$

This makes the operator $TS$ a targeted contraction with respect to the pair $(\hat{x} - \mu\tilde{x}, \hat{x})$ on the set $M_\alpha$, and we immediately conclude the following.

**Proposition 10.5.4** *Let $\{y^\nu\}$ be the sequence produced by the IDR algorithm for the 4-tuple $(\tilde{y}, y^0, TS, \mu)$ where $S$ is a linear operator from $Y$ into*

*X; let $\{x^\nu\}$ be defined by equation (10.73); and let $0 \le \alpha < 1$. If $x^\nu \in M_\alpha$ then*

$$\| \hat{x} - x^{\nu+1} \| \le \alpha \| \hat{x} - x^\nu \| . \tag{10.76}$$

Iterative data refinement and iterative result refinement are similar processes. We will see later why in practice one may be preferable to the other. Here we note only that the two sequences $\{x^\nu\}$ and $\{y^\nu\}$ might have different convergence properties. If $S$ is a linear operator with a nontrivial null space, i.e., $\mathcal{N}(S) \doteq \{y \in Y \mid Sy = 0\} \neq \{0\}$, then it is possible that $\{y^\nu\}$ does not converge while $\{x^\nu\}$ does.

We now look further into the possible structure of $T$. Assume we have a model of the ideal measurement process (expressed by an operator $A$ from $X$ into $Y$) and of the actual measurement process (expressed by an operator $B$ from $X$ into $Y$). Our intent is that if $x$ is near $\hat{x}$, then $Ax$ is near $\hat{y}$ and $Bx$ is near $\tilde{y}$. We may then define

$$T \doteq A - \mu B. \tag{10.77}$$

Since $S$ is used to recover the ideal result from the ideal measurement and $A$ is supposed to represent the ideal measurement process, $S$ and $A$ are "nearly inverses" of each other. More formally we obtain the following proposition.

**Proposition 10.5.5** *Let $\{y^\nu\}$ be the sequence produced by the IDR algorithm for the 4-tuple $(\tilde{y}, y^0, (A - \mu B)S, \mu)$, where $S$ is an operator from $Y$ into $X$ and $A$ and $B$ are operators from $X$ into $Y$ (all of them can be nonlinear). If for some $\nu \ge 0$*

$$y^\nu = ASy^\nu, \tag{10.78}$$

*then*

$$y^{\nu+1} = y^\nu + \mu(\tilde{y} - BSy^\nu). \tag{10.79}$$

In equation (10.79) $y^\nu$ is updated to $y^{\nu+1}$ by adding a correction term $\tilde{y} - BSy^\nu$ multiplied by the relaxation parameter $\mu$. It is often the case that the practical behavior of such an algorithm is improved by the use of variable relaxation parameters, i.e., by using a sequence of real numbers $\{\mu_\nu\}$ and defining the sequence $\{y^\nu\}$ by

$$y^{\nu+1} = y^\nu + \mu_\nu(\tilde{y} - BSy^\nu). \tag{10.80}$$

We now state the corresponding IRR algorithm under a suitable assumption on $A$ and $S$.

**Proposition 10.5.6** *Let $\{y^\nu\}$ be the sequence produced by the IDR algorithm for the 4-tuple $(\tilde{y}, y^0, (A - \mu B)S, \mu)$ where $S$ is an operator from $Y$* *into $X$ and $A$ and $B$ are operators from $X$ into $Y$. Let $\{x^\nu\}$ be defined by equation (10.73) and let $\tilde{x} = S\tilde{y}$. If for some $\nu \geq 0$,*

$$x^\nu = SAx^\nu, \tag{10.81}$$

*then*

$$x^{\nu+1} = x^\nu + \mu(\tilde{x} - SBx^\nu). \tag{10.82}$$

One can again use variable relaxation parameters $\{\lambda_\nu\}$ and replace the iterative step in equation (10.82) with

$$x^{\nu+1} = x^\nu + \lambda_\nu(\tilde{x} - SBx^\nu). \tag{10.83}$$

More specifically we consider the case in which $S$ is the inverse of the operator $A$, and therefore $SAx = x$ for all $x$ in $X$, and equation (10.81) is satisfied for all $x^\nu$. Writing $S$ as $A^{-1}$ equation (10.82) yields

$$x^{\nu+1} = x^\nu + \mu(\tilde{x} - A^{-1}Bx^\nu). \tag{10.84}$$

The process expressed by equation (10.84) can also be derived using a standard method of numerical analysis. Given a function $F : \mathbb{R}^n \to \mathbb{R}^n$ and some nonsingular matrix $A$, the *parallel-chord method* (see, e.g., Ortega and Rheinboldt (1970)) aims at finding a root of $F(x)$ by the iterative process

$$x^{\nu+1} = x^\nu - A^{-1}F(x^\nu). \tag{10.85}$$

The following proposition relates IDR to the parallel-chord method.

**Proposition 10.5.7** *Let $\{y^\nu\}$ be the sequence produced by the IDR algorithm for the 4-tuple $(\tilde{y}, y^0, (A-\mu B)A^{-1}, \mu)$ where $A$ is an invertible linear operator from $X$ into $Y$, and let $x^\nu = A^{-1}y^\nu$. Then the sequence $\{x^\nu\}$ satisfies the parallel-chord iteration formula (10.85) with*

$$F(x) \doteq \mu(Bx - \tilde{y}). \tag{10.86}$$

It is interesting to reexamine Proposition 10.5.2 in light of Proposition 10.5.7. In Proposition 10.5.7 $U$ is $(A-\mu B)A^{-1}$, which is linear if $B$ is linear. Hence the convergence of IDR to $(BA^{-1})^{-1}\tilde{y}$ is guaranteed by Proposition 10.5.2 provided that $\| I - \mu BA^{-1} \| < 1$. Furthermore if $A\hat{x} = \hat{y}$ and $B\hat{x} = \tilde{y}$ (i.e., the idealized device operator and the actual device operator work as intended on the idealized result $\hat{x}$), then $\tilde{y} = (BA^{-1})\hat{y}$ and IDR converges to the idealized measurement $\hat{y}$.

Several standard techniques of numerical analysis are special cases of the IDR. There is also a class of algorithms for constrained *iterative signal restoration* that are special cases of IDR. In the restoration problem the space $X$ is a subspace of $Y$ and the idealized measurement process is merely identification (e.g., $\hat{y}$ is $x$ in Figure 10.7); and the actual measurement process involves some *distortion*, producing a distorted version $\tilde{y}$ of $\hat{y} = x$. The restoration of $\hat{y}$ from $\tilde{y}$ via an iterative procedure is facilitated by prior knowledge of some of the properties of the desired solution. This can be expressed by defining a constraint operator $C$ from $Y$ into $X$ such that $x = Cx$ for all $x$ in $X$. Note that $C$ can take the role of the recovery operator in Figure 10.7 since $\hat{y} = x$ is in $X$ and $\hat{x} = C\hat{y} = x$, which corresponds to a perfect recovery. In this application, the idealized device operator of Figure 10.7 should be the identity $I$. If we use $D$ to denote the operator that represents the distortion (the actual device operator of Figure 10.7) and let $y^0 = \tilde{y}$ then equation (10.63) translates into

$$y^{\nu+1} = \mu y^0 + (I - \mu D)Cy^\nu, \tag{10.87}$$

which is exactly the general procedure for iterative signal restoration.

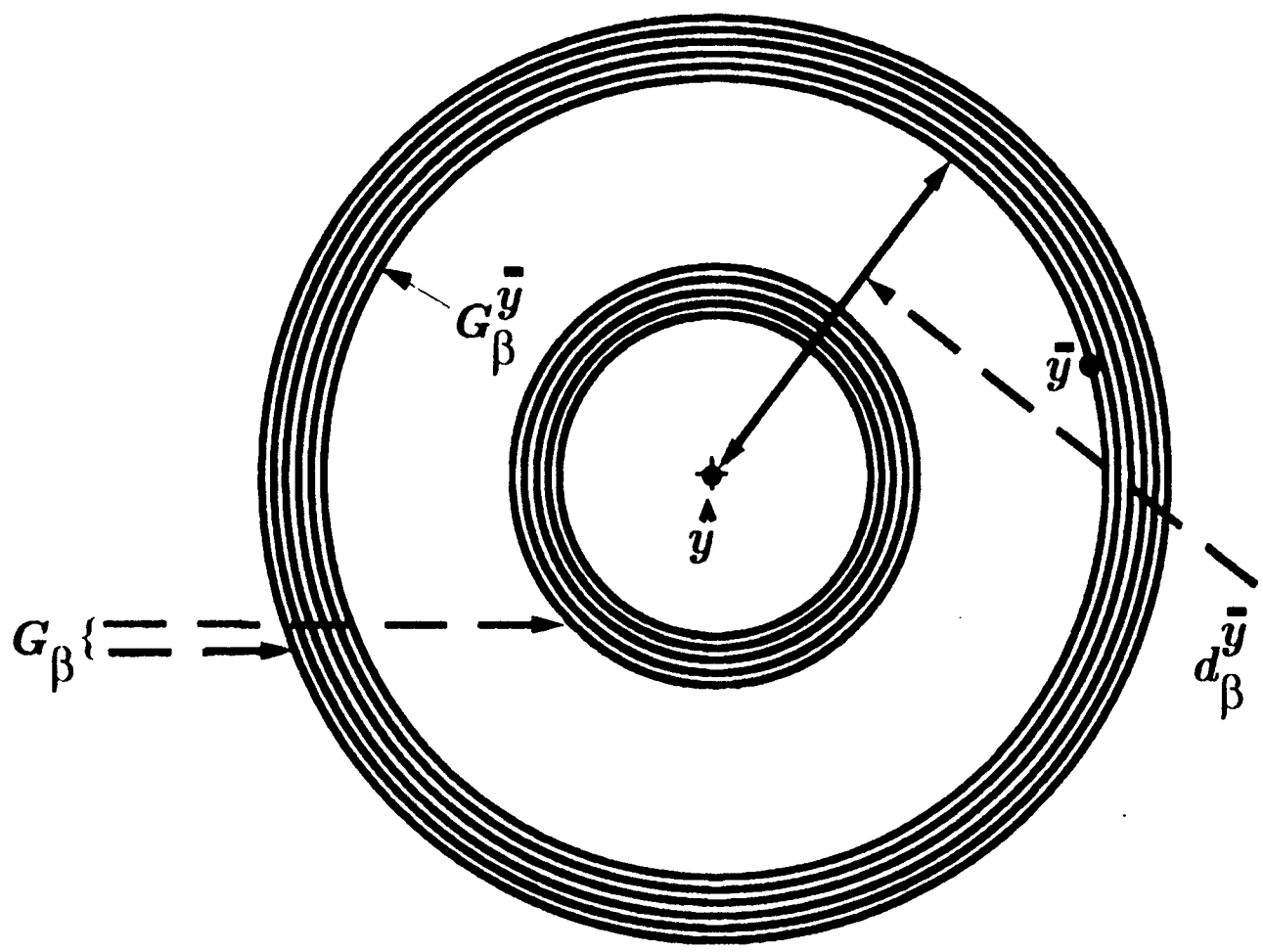

**Figure 10.8** The set $G_\beta$ is a union of disjoint spherical shells whose common center is $\hat{y}$.

We now return to the general form of IDR as given by (10.63). By strengthening the assumptions we derive—in Proposition 10.5.9—a sharpening of the result, which contains (10.72). Although the practical usefulness of the following material depends on the ability to verify the stronger assumptions stated below, the results further clarify the behavior of IDR. For any real $0 \leq \beta < 1$ define the set

$$G_\beta \doteq \{y \mid \sup_{\|z-\hat{y}\|=\|y-\hat{y}\|} \| \hat{y} - \mu\tilde{y} - Uz \| \leq \beta \| y - \hat{y} \|\}. \tag{10.88}$$

This means that we look at all points $z$ on a sphere through $y$ centered at $\hat{y}$, and that this $y$ belongs to the set $G_\beta$ if and only if the above supremum taken over all points $z$ on the sphere is bounded from above by $\beta$ times the radius of the sphere. The set $G_\beta$ is always a (possibly infinite) union of spheres centered at $\hat{y}$, since for any $y$ in $G_\beta$ all points of equal distance from $\hat{y}$ are also in $G_\beta$ (see Figure 10.8). It is easy to see that $G_\beta \subseteq N_\beta$, with $N_\beta$ as in equation (10.71).

Moreover, for each $\beta$, $0 < \beta < 1$, $G_\beta$ is the union of a (possibly infinite) family of disjoint sets, each one of which is a connected component (i.e., a shell) of $G_\beta$. We use $G_\beta^{\overline{y}}$ to denote the component of $G_\beta$ that contains the point $\overline{y}$. For any $\overline{y}$ in $G_\beta$,

$$d_\beta^{\overline{y}} \doteq \inf_{y \in G_\beta} \| y - \hat{y} \|, \tag{10.89}$$

is the distance from $\hat{y}$ to the point in $G_\beta^{\overline{y}}$ that is nearest to $\hat{y}$. The following result says that the IDR algorithm (i.e., equation (10.63)) starting with an element $y^0 \in G_\beta$ (with some $0 < \beta < 1$) will either converge to the ideal $\hat{y}$, or will, at least, get nearer to $\hat{y}$ than any element of $G_\beta^{y^0}$.

**Proposition 10.5.8** *Let $\{y^\nu\}$ be the sequence produced by the IDR algorithm for the 4-tuple $(\tilde{y}, y^0, U, \mu)$, and assume that $y^0 \in G_\beta$, for some $0 < \beta < 1$. Then either (i) or (ii) holds:*

*(i) For all $\nu \geq 0$, $y^\nu \in G_\beta$ and*

$$\lim_{\nu \to \infty} y^\nu = \hat{y}. \tag{10.90}$$

*(ii) There exists an index $\nu_0 \geq 0$ such that $y^0, y^1, y^2, \ldots, y^{\nu_0}$ are all in $G_\beta^{y^0}$ and*

$$\| y^{\nu_0+1} - \hat{y} \| < d_\beta^{y^0}. \tag{10.91}$$

**Proof** Observe that it is not true that $(i)$ implies $(ii)$, because $(i)$ may hold with $d_\beta^{y^0} = 0$ and (10.91) cannot hold. To prove the proposition we show that if $(ii)$ fails to hold then $(i)$ must. If $(ii)$ does not hold then $y^\nu \in G_\beta^{y^0}$ for all $\nu \geq 0$. Since $G_\beta^{y^0} \subseteq G_\beta \subseteq N_\beta$, (10.72) holds implying that $\| y^\nu - \hat{y} \| \leq (\beta)^\nu \| y^0 - \hat{y} \|$, which implies $(i)$. ■

The assumption $y^0 \in G_\beta$ in this proposition means that

$$\{y \mid \| y - \hat{y} \| = \| y^0 - \hat{y} \|\} \subseteq N_\beta, \tag{10.92}$$

i.e., that $N_\beta$ contains all points that are at the same distance from $\hat{y}$ as $y^0$ is. This is much stronger than the assumption that led to (10.72).

Finally, we give a sufficient condition for the convergence of $\{y^\nu\}$ to $\hat{y}$.

**Proposition 10.5.9** *Let $\{y^\nu\}$ be the sequence produced by the IDR algorithm for the 4-tuple $(\tilde{y}, y^0, U, \mu)$. If there exists a number $\beta$, $0 \le \beta < 1$ such that $y^0 \in G_\beta$ and*

$$\{y \mid 0 < \| y - \hat{y} \| < \| y^0 - \hat{y} \|\} \subseteq G_\beta, \tag{10.93}$$

*then $\{y^\nu\}$ converges to $\hat{y}$.*

**Proof** This is an immediate consequence of Proposition 10.5.8. ■

In practical applications the conditions of Proposition 10.5.9 are unlikely to be satisfied because $U$ provides only an estimate of $\hat{y} - \mu\tilde{y}$. Thus the most one can hope for is applicability of the second part of Proposition 10.5.8, which guarantees *initial improvement* provided $y^0 \in G_\beta$, but will not guarantee convergence to the ideal $\hat{y}$. We will illustrate this below on a small example from *beam hardening correction,* after we have shown how IDR incorporates procedures proposed in different areas of medical imaging.

Before beginning discussion of these applications, we summarize the notation that we will use. $Y$ is a space of measurements, whose elements include $\hat{y}$ (the ideal measurement), $\tilde{y}$ (the actual measurement), and $y^\nu$ for all $\nu \ge 0$ (the $\nu$th iterative refinement). $U$ is the operator from $Y$ into $Y$ used in the iterative data refinement procedure. $S$ is a usually linear operator from $Y$ into $X$, the space of results. $X$ includes elements $\hat{x} = S\hat{y}$ (the ideal result), $\tilde{x} = S\tilde{y}$, and $x^\nu = Sy^\nu$ for $\nu \ge 0$. For any operator $T$ from $X$ into $Y$, $U = TS$ can be used to define an IDR procedure. One choice for $T$ is $A - \mu B$, where $A$ and $B$ represent the ideal and actual measurement processes, respectively.

### 10.5.2 Applications in medical imaging

In this section we discuss techniques used in three different areas of medical imaging with image reconstruction from projections, and show that each is a special case of IDR. In the notation of the last section, $X$ is a Banach space of functions of two variables defined over a fixed compact subset $Q$ of the plane. We refer to elements $x$ of $X$ as *images.* The set $Y$ is $\mathbb{R}^m$, the $m$-dimensional Euclidean space. The ideal measurement process $A$ produces an $m$-dimensional vector $Ax$ whose $i$th component

$$(Ax)_i = \int_{L_i} x dz \tag{10.94}$$

is a line integral of $x$ along a line $L_i$ (this is the same as equation (10.2) but in different notation).

There are many processes (both linear and nonlinear) for estimating $x$ from $Ax$. Provided that there are many uniformly spaced lines $L_i$, these processes produce efficacious results. We shall not discuss here the details of the processes, but simply use the symbol $S$ to denote the operator from $Y$ into $X$ that is used to estimate $x$ from $Ax$. Unless otherwise stated, we assume $S$ to be linear.

Since $Y$ is finite dimensional and $S$ is linear, $SY$ is a finite-dimensional subspace of $X$. We distinguish between two cases depending on whether the elements of $SY$ are treated as continuous or discretized. In the continuous case we use a closed-form formula, which for any $y$ in $Y$ and for any point of $Q$ gives the value of the function $Sy$ at that point in terms of the $m$ components $\{y_i \mid 1 \leq i \leq m\}$ of $y$ and the two coordinates of the point. If we let $c_i$ denote the element of $X$ obtained by applying $S$ to the standard basis vector whose $i$th component is 1 and all other components are 0, then by the linearity of $S$,

$$Sy = \sum_{i=1}^{m} y_i c_i. \tag{10.95}$$

Thus $\{c_i \mid 1 \leq i \leq m\}$ is a *natural basis* for $SY$. Note that because $A^{-1}$ in general does not exist, different choices for $S$ may lead to very different $c_i$'s.

There are situations in which the natural basis representation is not appropriate, for example when we wish to display images. Displayed images are uniformly gray within *pixels*, which cover the support $Q$ of functions in $X$. Suppose that there are $n$ such pixels in total; let $b_j$ denote the function whose value is 1 in the $j$th pixel and 0 everywhere else, and assume that these $b_j$'s are elements of $X$. Certain reconstruction methods produce pictures that are in the subset of $X$ generated by this *pixel basis*, i.e.,

$$Sy = \sum_{j=1}^{n} x_j b_j. \tag{10.96}$$

Note that in such a case $S$ has a representation as an $m \times n$ matrix that produces for any $m$-dimensional measurement vector $y$ an $n$-dimensional image vector whose $j$th component is $x_j$; this is the fully discretized approach discussed in Section 10.1 (compare with equation (10.10)).

While a pixel basis is necessary for display and has been popular in image reconstruction in general, its use with IDR is a potential source of an additional error that does not arise when the natural basis is used. Forcing of the result of an image reconstruction method based on a closed-form formula to be representable in the pixel basis introduces an additional approximation.

The three applications that we briefly discuss next within the framework of IDR are: beam hardening correction in x-ray computerized tomography; attenuation correction in single photon emission computerized tomography; and image reconstruction with incomplete data.

**Beam Hardening Correction in x-ray Computerized Tomography**

The linear x-ray attenuation at a point inside a cross section of the human body depends on the position of the point and on the energy $E$ of the x-ray. We use $x_E$ to denote the element of $X$ that represents the linear attenuation distribution at energy level $E$. If we had a monoenergetic x-ray source at energy $\overline{E}$ we could estimate $(Ax_{\overline{E}})_i$ by placing the x-ray source and detector on the two sides of the body along the line $L_i$, and calculating $-\log(I_i'/I_i)$ where $I_i$ and $I_i'$ are the measured intensities of the x-ray beam entering and exiting (respectively) the body along the line $L_i$. These are the ideal measurements

$$\hat{y} = Ax_{\overline{E}}, \tag{10.97}$$

from which $x_{\overline{E}}$ can be reconstructed using $S$.

The x-rays used in practice for diagnostic CT imaging are polyenergetic. Since the linear attenuation coefficient of various biological materials is highly dependent on the photon energy, the spectrum of the x-rays is modified as they travel through the body cross section. Consequently the attenuation due to the tissue at a point will be different for rays from different directions. Thus if we evaluate $-\log(I_i'/I_i)$ for a polyenergetic x-ray beam along $L_i$, the result will not be $\hat{y}_i$ (or the line integral of any two-dimensional distribution) but rather an estimate for the more complicated integral

$$p_i = -\log \int_0^\infty \tau(E) \exp\left(-\int_{L_i} x_E dl\right) dE, \tag{10.98}$$

where $\tau(E)$ represents the fractional number of photons at energy $E$ in the (measured) spectrum of the x-ray source. Currently we do not have any rigorous mathematical approaches for the recovery of $x_{\overline{E}}$ from these components $p_i, i = 1, 2, \ldots, m$, of the vector $p$.

Since $p$ does not uniquely determine $\hat{y}$, the practical question is: "Given $p$, can we approximate $\hat{y}$ well enough so that it leads to diagnostically useful CT reconstructions?" The answer is affirmative. One class of methods for obtaining such an approximation are called post-reconstruction beam hardening correction methods. Such methods proceed in two stages.

*Initial stage:* A real-valued function $q$ is specified and $q(p_i)$ is used as the estimate $y_i^0 = \tilde{y}_i$, of $\hat{y}_i$, for $1 \le i \le m$. A reconstruction is performed which is referred to as the *initial reconstruction.*

*Iterative stage:* A sequence of corrections is applied, where each correction is calculated on the basis of earlier reconstructions. The process stops with the *final reconstruction* when further iterations do not change the reconstruction significantly.

The exact nature of the function $q$ in the initial stage is irrelevant to our discussion here. There have been a number of different proposals for performing the iterative stage; following we describe one method by way of an appropriate example.

| | $L_3$ ↓ | $L_4$ ↓ |
|---|---|---|
| $L_1$ → | $a$ 1 | $b$ 2 |
| $L_2$ → | $c$ 3 | $c$ 4 |

**Figure 10.9** Object used in illustrative example of the IDR procedure.

The method that we describe makes use of the pixel basis (see equation (10.96)) and the iterative step expressed in (10.74) with $\mu = 1$. The operator $T$ is of the form $A - B$, where $A$ is the operator of (10.94) (the line integrals of the ideal measurement process) and $B$ is an approximation of the actual measurement process (equation (10.98)) as modified by $q$. One difficulty in defining $B$ is due to the fact that $x^\nu$ is an approximation to $x_{\overline{E}}$, but equation (10.98) requires $x_E$ for all values of $E$.

To overcome this difficulty we can define, for every energy level $E$ a function $g_E$ such that if $x_j$ is the value of $x_{\overline{E}}$ in the $j$th pixel, then $g_E(x_j)$ approximates the value of $x_E$ in the $j$th pixel. Denoting by $a_{ij}$ the length of intersection of line $L_i$ with the $j$th pixel, we define the operator $B$, for a given function $q$ by

$$(B(\sum_{j=1}^{n} x_j b_j))_i \doteq q(-\log \int_0^\infty \tau(E) \exp\left(-\sum_{j=1}^{n} a_{ij} g_E(x_j)\right) dE). \quad (10.99)$$

Similarly,

$$(A(\sum_{j=1}^{n} x_j b_j))_i \doteq \sum_{j=1}^{n} a_{ij} x_j. \quad (10.100)$$

We set $T = A - B$ and use equation (10.74). Note that since $B$ is not linear $T$ and $U$ are not linear, and the results that are based on assuming linearity (such as Proposition 10.5.2) are not applicable.

Let us look at a detailed example of the behavior of IDR on the beam hardening correction problem. Our aim is merely to illustrate the IDR procedure and we use an extremely simple example. Use of IDR for beam hardening correction on realistic examples has been reported in the literature. We consider a very simple object comprised of three different tissue types $a, b$, and $c$ occupying pixels as indicated in Figure 10.9. Pixels are assumed to be of unit size. We use $\mu_E^t$ to denote the linear attenuation of tissue $t$ at energy $E$. Table 10.1 gives the assumed values of $\mu_E^t$ for tissues $a, b$, and $c$, and energies 30, 65, and 100. (The example is for illustration only and we do not specify units.) If we choose $\overline{E} = 65$, then $x_{\overline{E}}$ is described by the vector whose transpose is (6.0, 6.0, 2.0, 2.0) (see equation (10.96)). The projection operator $A$ can be represented by the matrix

$$A = \begin{pmatrix} 1 & 1 & 0 & 0 \\ 0 & 0 & 1 & 1 \\ 1 & 0 & 1 & 0 \\ 0 & 1 & 0 & 1 \end{pmatrix} \tag{10.101}$$

and so $\hat{y} = (12.0, 4.0, 8.0, 8.0)^{\mathrm{T}}$. An appropriate reconstruction operator $S$ is given by the *generalized inverse* $A^\dagger$ of the matrix $A$ in (10.101), namely

$$A^\dagger = \frac{1}{8} \begin{pmatrix} 3 & -1 & 3 & -1 \\ 3 & -1 & -1 & 3 \\ -1 & 3 & 3 & -1 \\ -1 & 3 & -1 & 3 \end{pmatrix}. \tag{10.102}$$

**Table 10.1** Values of the linear attenuation coefficients $\mu_E^t$.

| $E$ | $t = a$ | $t = b$ | $t = c$ |
|---|---|---|---|
| 30 | 7.2 | 7.5 | 2.4 |
| 65 | 6.0 | 6.0 | 2.0 |
| 100 | 3.6 | 5.4 | 1.2 |

Note that multiplying $\hat{y}$ with this matrix does give us the original picture $x_{\overline{E}}$. Let us assume that in the polyenergetic x-rays one-third of the photons are at energy 100 and two-thirds of the photons are at energy 30. Using the identity function for $q$, equation (10.98) yields

$$y^0 = (10.0929, 3.3327, 5.8833, 6.683)^{\mathrm{T}}. \tag{10.103}$$

This is a rather rough approximation of $\hat{y}$. In fact, $\| \hat{y} - y^0 \| = 2.950$.

To apply IDR we need to define $g_E$ for $E = 30$ and $E = 100$. We have chosen

$$g_{30}(x_j) \doteq \frac{x_j - 2.0}{4.0} 4.95 + 2.4, \tag{10.104}$$

$$g_{100}(x_j) \doteq \frac{x_j - 2.0}{4.0} 3.3 + 1.2. \tag{10.105}$$

Observe that these choices of $g_{30}$ and $g_{100}$ assign the correct value to tissue $c$, and the average of the correct values for tissues $a$ and $b$. In view of Table 10.1, we see that we cannot possibly assign correct values for both the tissues $a$ and $b$, since if $x_j = 6.0$ we cannot tell whether the tissue in the $j$th pixel is of type $a$ or $b$. This implies that the actual device operator defined by equation (10.99) will not produce the same result as the actual measuring device, even at $x = \hat{x} = x_{\overline{E}}$, i.e., that $Bx = B\hat{x} \neq \tilde{y}$. In view of this, we cannot expect convergence of IDR to $\hat{y}$, but only an initial improvement (see the remarks after Proposition 10.5.9).

**Table 10.2** Errors in the IDR sequence $\{y^\nu\}$ generated for the illustrative example.

| $\nu$ | $\| \hat{y} - y^\nu \|$ |
|---|---|
| 0 | 2.950 |
| 1 | 1.488 |
| 2 | 1.461 |
| 3 | 1.468 |

Using the operator $T$ as defined by equations (10.77), (10.99), and (10.100) with $\mu = 1$ and the operator $S$ given by (10.102), the beginning of the sequence $\{y^\nu\}$ produced by IDR for the 4-tuple $(y^0, y^0, TS, 1)$ is

$$y^1 = (11.6778,\ 3.9329,\ 6.7904,\ 8.801)^{\mathrm{T}},$$
$$y^2 = (11.9437,\ 3.9920,\ 6.9286,\ 8.9944)^{\mathrm{T}},$$
$$y^3 = (11.9898,\ 3.9988,\ 6.9495,\ 9.0261)^{\mathrm{T}}.$$

The corresponding values of the errors are given in Table 10.2. Note that the initially rapid improvement from $\nu = 0$ to $\nu = 1$ is substantially reduced from $\nu = 1$ to $\nu = 2$, and the error actually gets slightly worse from $\nu = 2$ to $\nu = 3$. In view of the result stated in (10.72), this implies that $y^2$ is not an element of $N_\beta$ for any $\beta$ less than 1. While our example is unrealistically simple, this behavior matches the behavior of beam hardening correction methods reported in the literature.

### Attenuation Correction in Single Photon Emission Computed Tomography (SPECT)

In SPECT the function $x$ that we wish to reconstruct is the distribution of activity of some radiating material in a cross section of the human body. Using an appropriate detector placed outside the body we could measure line integrals of $x$, if it were not for the fact that the body attenuates the emitted radiation between the source point inside the body and the detector.

Assuming that the attenuation function $z$ is known, let $x(i, s)$ and $z(i, s)$ denote the activity and attenuation respectively, along the line $L_i$ at a distance $s$ from the detector. Then the actual measurement can be represented by the operator $B$ from $X$ into $Y$, given by

$$(Bx)_i = \int_0^\infty x(i,s) \exp\left(-\int_0^s z(i,s')\, ds'\right) ds. \tag{10.106}$$

One method proposed for finding the ideal measurements (the line integrals of $x$) is given by (10.80). The values $\mu_1 = 2/3$ and $\mu_\nu = 1/3$ for $\nu > 1$ have been suggested in the literature.

An interesting point is the possibility of considering the operator $BS$ as a single operator, thus avoiding the need for explicitly producing a pixel-based version of each $Sy^\nu$. Note that since $BS$ is a linear operator from $\mathbb{R}^m$ into $\mathbb{R}^m$ it can be described by an $m \times m$ matrix whose $(i, i')$th element in our natural pixel notation is

$$(BS)_{i,i'} = (Bc_{i'})_i, \tag{10.107}$$

(see equations (10.95) and (10.106)). The matrix $BS$ may have special properties that allow its elements to be calculated efficiently. It has been reported in the literature that using this method on realistic SPECT examples produced the fifth iterate $x^5$ as an excellent approximation to $\hat{x}$.

**Image Reconstruction with Incomplete Data**

Let us return to transmission CT. If there are many uniformly spaced lines $L_i, 1 \le i \le m$ along which integrals are taken, then there are efficient reconstruction operators $S$ available in the literature for recovering $x$ from $Ax$. In certain applications the collected data set is incomplete in the following sense. Let $I = \{1, 2, \ldots, m\}$ and $I = K \cup M$, where the sets of indices $K$ and $M$ have no elements in common. The data that are collected provide us with a good approximation of $(Ax)_i$ for $i$ in $K$, but do not provide us with $(Ax)_i$ for $i$ in $M$. Various methods have been proposed, which are essentially IDR approaches, for estimating the whole vector $Ax$ from its known components. Translated into our notation the approach is the following. Let

$$\tilde{y}_i \doteq \begin{cases} \text{measured value of } (Ax)_i, & \text{if } i \in K, \\ 0, & \text{if } i \in M. \end{cases} \tag{10.108}$$

Accordingly (see the discussion before equation (10.77)), we define the operator $B$ from $X$ into $Y$ by

$$(Bx)_i \doteq \begin{cases} (Ax)_i, & \text{if } i \in K, \\ 0, & \text{if } i \in M. \end{cases} \tag{10.109}$$

If we now define $T$ by equation (10.77) with $\mu = 1$ we get, for any $x$ in $X$,

$$(Tx)_i \doteq \begin{cases} 0, & \text{if } i \in K, \\ (Ax)_i, & \text{if } i \in M. \end{cases} \tag{10.110}$$

The sequence $\{y^\nu\}$ produced by IDR for the 4-tuple $(\tilde{y}, \tilde{y}, TS, 1)$ has the property

$$y_i^{\nu+1} = \begin{cases} \tilde{y}_i, & \text{if } i \in K, \\ (TSy^\nu)_i, & \text{if } i \in M. \end{cases} \tag{10.111}$$

It has been reported in the literature that in this application area the incorporation of prior knowledge (such as $x$ being a positive-valued function) can improve the practical usefulness of the procedure. An operator $S$ that incorporates such knowledge is often nonlinear. However, in realistic simulations of various applications the $Sy^\nu$ provide very good approximation of $\hat{x}$ in about 5 to 10 iterations.

## 10.6 On the Selective Use of Iterative Algorithms for Inversion Problems

In this section we discuss the selective use of iterative algorithms for inversion problems. We distinguish between situations where the *discretization first* approach is appropriate, and those where *transform methods*, which rely on analytic inversion and discretize only the final inversion formulae, are apt to perform better. The discussion is motivated and supported by examples from image reconstruction from projections and from radiation therapy treatment planning (see Chapter 11).

The forerunner of iterative reconstruction algorithms was ART (Algebraic Reconstruction Technique), which was published in 1970 by Gordon, Bender, and Herman (1970); independently employed in a computerized tomography (CT) head scanner in 1972 by Hounsfield (1972); and subsequently recognized as Kaczmarz's 1937 method for solving systems of linear equations (Algorithm 5.4.3). The controversy between iterative algorithms and transform methods seems to have been resolved in the practice of x-ray transmission CT for medical purposes, where versions of the transform method called *Filtered Backprojection* (FBP) have been universally adopted in most (although not all) cases. Iterative reconstruction algorithms continue however to be recognized as useful tools in various specific situations in x-ray CT in medical imaging and in nonmedical applications of image reconstruction. While they have been variously criticized (e.g., Shepp and Kruskal (1978, Section 3), advocated (e.g., Herman and Lent (1976a, Section 6)), and compared with transform methods (e.g., Censor (1983, section VIII), or Censor and Herman (1987, Section 4)), iterative reconstruction algorithms kept appearing as tools of choice in many application areas of image reconstruction, until and including recent papers regarding the EM algorithm for *positron emission tomography* (PET), (see Sections 10.2 and 10.2.1).

Iterative inversion or analytic inversion: when is one likely to perform better? The diagram in Figure 10.10 illustrates why iterative reconstruction algorithms are not always suitable in *x-ray CT medical imaging* and why many other inversion problems in image reconstruction, radiation therapy treatment planning and other fields call for their use. All these situ-ations can benefit immensely from the use of optimization or feasibility algorithms that exploit parallelism, since traditionally these applications have been computationally very intensive.

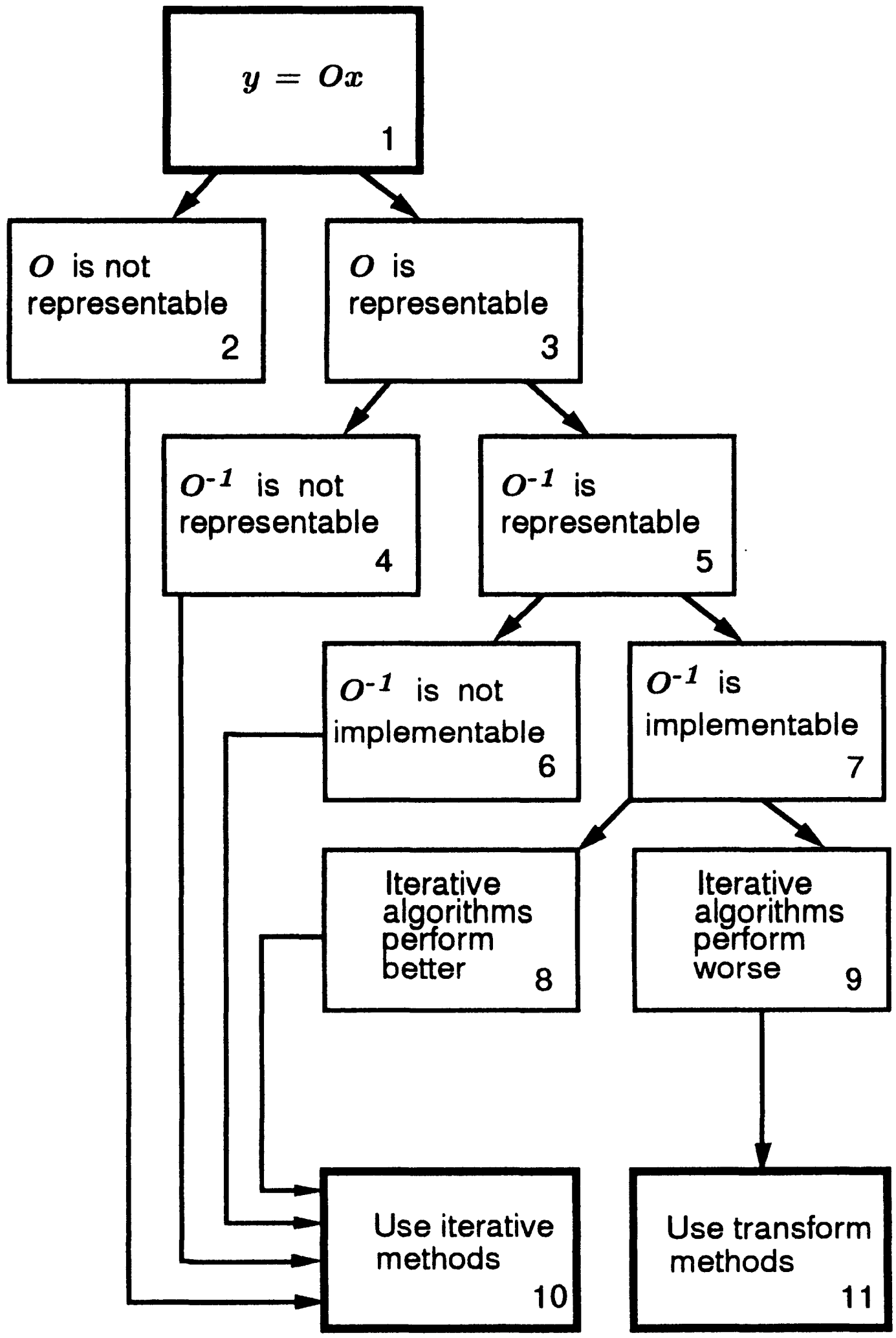

**Figure 10.10** Selective use of iterative algorithms in inversion problems.

In the diagram, $y = \mathcal{O}x$ is a functional equation that relates an unknown object $x$ to some measurements or other data $y$ through an operator $\mathcal{O}$. The object $x$ and the data $y$ belong to appropriate function spaces, and the fundamental inversion problem is to recover $x$ or an acceptable approximation thereof. The diagram identifies in general terms situations in which transform methods are not likely to perform well. In such cases full discretization of the model and use of iterative algorithms for the resulting optimization or feasibility problems are potentially advantageous and sometimes the only way to treat the problem. The right-hand side of the diagram (all even-numbered boxes) describes a logical path necessary for successful application of a transform method. The default options—boxes 2, 4, 6, 8—all lead to the consideration of iterative algorithms. We now explain the phrases in the diagram.

*Representable* means that an operator (either $\mathcal{O}$ or its inverse $\mathcal{O}^{-1}$) can be represented in a mathematically closed-form analytic formula. *Implementable* means that all conditions which are required either mathematically or practically for efficient implementation of the computations involved in $\mathcal{O}^{-1}$ are met. *Perform better* and *perform worse* are meant in a general sense wherein the quality of the results in terms of the specific real-world application as well as the efficiency and speed of computation are evaluated.

An example that fits the description of box 2, where $\mathcal{O}$ is not representable, is the inverse problem of radiation therapy treatment planning discussed in the next chapter. Even if $\mathcal{O}$ is representable (box 3), there is no guarantee that its inverse $\mathcal{O}^{-1}$ is also representable. In x-ray CT, $\mathcal{O}$ is the *Radon transform* (equations (10.3) and (10.4)) and an analytic derivation of its inverse is given by (10.5), thus it fits into box 5. An example that fits the description of box 4, where the analytic inversion has not yet been solved although the forward operator $\mathcal{O}$ is representable, is the attenuated Radon transform of the image reconstruction problem of *single photon emission computerized tomography* (SPECT) (see, e.g., Budinger, Gullberg, and Huesman (1979) or Knoll (1983)). The analytic inverse of this transform is known only when the attenuation coefficient has constant value throughout the object.

To efficiently implement a reconstruction method such as filtered backprojection, which is based on Radon's inversion formula, considerable data have to be available. Because interpolation in the data space is a required intermediate operation, insufficient or sparse data cause deterioration of the overall reconstruction. This is an example of the apparent dichotomy between boxes 6 and 7 in Figure 10.10. Thus, for example, *limited range* (Tuy (1981)) or *limited views* (see, e.g., Natterer (1986, Chapter VI)) image reconstruction or reconstruction from a very small number of views would often fall into box 6. When both a transform method and an iterative algorithm are implementable one may perform better (boxes 8 and 9) depending on characteristics of the true object $x$—such as whether its true image has high or low contrast. If the quality of $x$ generated by the iterative algorithm is poor, i.e., it is not useful for the particular application, then the passage to box 11 is prescribed. If however, the recovery by an iterative algorithm is better or no worse than recovery by a transform method but demands more computation time then the iterative algorithm approach can still be pursued. This situation would especially benefit from the use of parallel computations.

## 10.7 Notes and References

Many texts on image reconstruction from projections are available, varying according to emphasis, i.e., physical, radiological, or mathematical. A representative sample includes: Herman (1980), Barrett and Swindell (1981), Kak and Slaney (1988), Hurt (1989), Bates and McDonnell (1986), Deans (1983), Natterer (1986), Hendee (1983), and Louis (1989). Many book collections of articles have appeared over the years from the early days of the field's development, see Ter-Pogossian et al. (1977), Herman (1979b), and Shepp (1983) until the more recent volumes by Stark (1987), Newhouse (1988), Herman, Louis, and Natterer (1990), Censor, Elfving and Herman (1990), Udupa and Herman (1991), Nolet (1987), and many more—too numerous to be listed here. Additional notes and references relevant to some of the algorithms discussed here appear in Sections 5.11 and 6.10.

**10.1** The material in this section is based on papers by Gordon and Herman (1974), Herman and Lent (1976a), Lewitt (1983), Herman (1980), Censor (1983, 1988), and Censor and Herman (1987). Radon's inversion formula was published in Radon (1917). Transform methods are discussed by Lewitt (1983) as well as in the books mentioned above. Other review articles on transform methods, with various emphases, are Smith, Solmon, and Wagner (1977), Lindgren and Rattey (1981), Kak and Roberts (1986), Smith and Keinert (1985) and Louis (1992). Lewitt (1992) developed a theory that uses spherically symmetric functions, which are the so-called generalized Kaiser-Bessel window functions (Lewitt (1990)), often referred to as *blobs*. It has recently been demonstrated by Matey et al. (1994), for the case of fully three-dimensional PET reconstruction (see Section 10.2), that such basis functions indeed lead to statistically significant improvements in the task-oriented performance of reconstruction methods based on full discretization of the model.

**10.2** The maximum likelihood approach, see, e.g., Titterington, Smith and Makov (1985), for the fully discretized model in PET, was proposed by Rockmore and Macovski (1976). Shepp and Vardi (1982) and Vardi, Shepp, and Kaufman (1985) made this approach popular in image reconstruction by suggesting the EM algorithm (equation (10.32)) for the solution of the problem. The EM algorithm has been used in statistics for quite some time; see Dempster, Laird and Rubin (1977). The literature on the EM algorithm is vast, reporting both experimental results and mathematical studies. Our presentation is taken from Iusem (1991b). A general review of the class of multiplicative iterative methods, to which the EM algorithm belongs, is given by De Pierro (1991). Figure 10.6 is drawn after a figure in Vardi, Shepp, and Kaufman (1985) where more details on this topic and further references may be found.

**10.3** For a discussion of maximum entropy methods see Jaynes (1982). For a recent text on entropy optimization see Kapur and Kesavan (1992). The justification for entropy maximization in image reconstruction given here is taken from Iusem (1994) and is based on Elfving (1989) and Iusem (1988). Consult also Skilling and Gull (1984) for a thorough investigation of entropy expressions and many relevant references. The Bayes formula, and other notions from probability theory, that underlie the development in this section can be found in books on probability theory; see, e.g., Uspensky (1937), or Mullins and Rosen (1971). The Bayesian approach was studied as a methodology in image reconstruction; see Herman et al. (1979) and Hanson (1987). Bregman's method for the mixed entropic-quadratic problem (10.49)–(10.51) was analyzed by Elfving (1989). This method was applied to matrix balancing from noisy or incomplete data, in which case the constraint matrix $A$ is a $0-1$ matrix, by Zenios and Zenios (1992). Regarding hybrid algorithms, see also Elfving (1989) where a single step of Newton's method (rather than a secant step as in MART) was used. A recent review on iterative reconstruction algorithms based on cross-entropy minimization is Byrne (1996a), and a different kind of entropy maximization reconstruction algorithm, called MENT, is proposed in Minerbo (1979), and further studied by Dusaussoy and Abdou (1991). Herman (1982) did an experimental case study on the maximum entropy criterion in image reconstruction.

**10.4** The ART algorithm was invented for image reconstruction from projections by Gordon, Bender, and Herman (1970) and was, at about that time, independently utilized by Hounsfield (1972) in the first commercially available CT scanner produced by EMI, Ltd; see also Hounsfield (1973). The building of this machine, which marked the beginning of a revolution in diagnostic radiology, earned Hounsfield in 1979 a Nobel prize (shared with Cormack). ART was further studied within the field of image reconstruction from projections by Herman, Lent, and Rowland (1973), Herman, Lent, and Lutz (1978), and many others. The SIRT algorithm for image reconstruction is due to Gilbert (1972), and was further studied by Lakshminarayanan and Lent (1979). Van der Sluis and Van der Vorst (1987, 1990) studied SIRT and compared it with conjugate gradient type methods; see also Van der Vorst (1994). The block-ART, i.e., block-Kaczmarz algorithm with fixed blocks of equations was proposed by Eggermont, Herman, and Lent (1981) and applied to a real problem of tomographic image re-construction by Altschuler et al. (1980). Björck and Elfving (1979) and Elfving (1980a) also studied projection methods and block-iterative methods with regard to the image reconstruction problem. Earlier theoretical work was done by Peters (1976). Cimmino's algorithm was proposed in Cimmino (1938); see also Gastinel (1970). Many of these algorithms were brought into a unified framework of quadratic optimization and tested experimentally by Herman and Lent (1976b) and Artzy, Elfving, and Herman (1979); see also the survey by Herman and Lent (1976a). Dax's review (1990) of linear stationary iterative processes for solving singular unstructured systems of linear equations puts the methods of Kaczmarz and Cimmino in perspective. He discusses their relation to the SOR and the Jacobi methods. Herman and Levkowitz (1987) investigated the effect of different fixed-block sizes on the initial performance of the block-Kaczmarz algorithm. In spite of the limited scope of their reported experiments, Herman and Levkowitz noticed significant dependence on the block size. They found that implementations with block sizes that are intermediate to the extremes of row-action (one equation per block) and fully simultaneous (all equations in one block)—show better initial performance according to several calculated measures of merit. Further results concerning block implementations of an iterative reconstruction algorithm were published by Hudson and Larkin (1994), who used a block version of the EM algorithm (equation (10.32)). Byrne (1996b) has recently investigated block-iterative versions of MART and of the EM algorithm in a unified framework. However, a judicious choice of the order of the equations, and the choice of optimal relaxation parameters, even in row-action implementation of ART can lead to very good results. This has been recently demonstrated by a special data-access ordering developed by Herman and Meyer (1993). A recent tutorial by Herman (1996) reviews recent developments which improve the efficacy and reduce computational costs of reconstructions produced by ART algorithms.

ART has been used in many reconstruction problems outside medical imaging, see, e.g., Nolet (1987) and a recent application to turbulent field diagnostics by Stricker, Censor and Zakharin (1996). An analysis and generalizations of ART and Cimmino's method to operator equations in function space were given by Nashed (1981).

**10.5** The IDR method presented in this section was proposed by Censor, Elfving and Herman (1985). Constrained iterative signal restoration methods mentioned here which are special cases of IDR are surveyed by Shafer, Mersereau and Richards (1981). Material on the beam hardening problem in x-ray computerized tomography can be found, e.g., in Herman (1980) and in Stonestrom, Alvarez and Macovski (1981). The method described here is from Herman (1979a), more can be found in Herman and Trivedi (1983). Walters et al. (1981) did attenuation correction in SPECT which we recognize as being actually a process of IDR. Image reconstruction from incomplete data has received much attention and is the topic of many articles. The following authors have proposed what are essentially IDR approaches: Tam and Perez-Mendez (1983), Tuy (1983) and Nassi et al. (1982).

Proposition 10.5.9 is in essence Theorem 10.1.2 of Ortega and Rheinboldt (1970). IDR is closely related to the *Defect Correction Principle* of Stetter (1978); see the review paper by Böhmer, Hemker and Stetter (1984) and the collection of papers in Böhmer and Stetter (1984). Ro, Herman and Joseph (1989) and Ro et al.(1992) have applied IDR to image reconstruction in magnetic resonance imaging (MRI). Herman and Ro (1990) showed that well-known image recovery methods such as the Gerchberg-Saxton algorithm and the error reduction and hybrid input-output methods of Fienup (see Herman and Ro (1990) for references) are also encompassed by IDR.

For the notion of *generalized inverse* of a matrix (such as (10.102)) see, e.g., Ben-Israel and Greville (1974). The notion of a *targeted contraction* (Definition 10.5.1) generalizes the well-known notion of a *contraction*, see, e.g., Luenberger (1969). It also generalizes the notion of *pseudocontraction* used by Bertsekas and Tsitsiklis (1989). See also Istrăţescu (1981, Definition 6.6.4).

CHAPTER 11

# The Inverse Problem in Radiation Therapy Treatment Planning

In this chapter we introduce the problem of therapy treatment planning in radiation oncology. There are some good reasons to construct a fully discretized model of this problem, to which we can apply the various row-action algorithms and block-iterative methods discussed in earlier chapters. We will explain our reasons and also present some computational results.

Radiation therapy maintains a prominent place in the treatment of neoplastic lesions (tumors). Over the past ninety years, medical experience and laboratory studies have established the principles whereby radiation can be used as an effective modality for oncological treatment.

Radiation *dose* is a measure of the actual energy absorbed in a biological (or other) absorber from a given x-ray beam. It is measured in units of energy per mass. The dose unit is called *gray* and is abbreviated by Gy, and 1 Gy $\doteq$ 1 Joule/Kg (where Joule is the standard unit of energy). It is very difficult to define quantitatively the actual biological damage produced by the dose but since the dose is closely related to the amount of ionization produced we can roughly estimate the degree of biological damage. Although all viable tissues are affected to some degree by ionizing radiation, a dose below a certain empirical limit, or tolerance value, does not cause irreparable damage to normal tissue. However, a dose of radiation that exceeds the tolerance can be expected to cause permanent tissue damage, and the tolerance value varies with tissue type. For any given organ, the risk of injury depends on the total and incremental dose, and on the fraction of the organ irradiated. Similarly, empirical values have been determined for doses lethal to different tumors. Generally the range of dose that will destroy a tumor without an unacceptable risk of injury to normal tissue is very narrow.

The radiation dose delivered to any point within a region is determined by two components, the *primary* and the *scatter* radiation. The primary radiation causes that portion of the dose which results from interactions lying directly along the path of the undeviated radiation beam. All other contributions to the dose at a specific point come from radiation deflected, i.e., scattered, to that point from within the radiation field.

After the lesion is diagnosed, a decision must be made on treatment modality, e.g., surgery, chemotherapy, radiation therapy, or some combination of these. If a tumor is of the size and type likely to respond to radiation therapy, the radiation therapist (physician) delineates the volume to be treated and prescribes the dose he believes most effective to the lesion. Together with the radiation therapist, the radiation dosimetrist (technician) must ensure that this physician-specified dose is correctly delivered to the site of the lesion, without endangering the function of critical organs or normal tissue.

Radiation can be delivered in two ways: by *brachytherapy*, i.e., the direct implantation of radioactive sources into the lesion, and by *teletherapy*, i.e., the use of beams of penetrating radiation directed at the lesion from an external source. When radioactive sources are implanted into tumors the dose is determined by the strength of the sources, the geometric arrangement of these sources within the tumor region, and the duration of insertion. Brachytherapy calculations are relatively straightforward and show a predictable dose gradient from the high-dose regions around the sources to the small-dose regions at increasing distances from the sources.

When radiation is delivered by beams from an external source (e.g., cobalt-60 or a linear accelerator), the situation becomes quite complex. Beams from the external radiation source are shaped and directed at a specific target region, namely, the volume of the neoplastic lesion and the immediately surrounding volume that is likely to contain microscopic diseases. The beams must be aimed and contoured so that when their dose contributions are added, the total dose is lethal to the neoplastic cells in the target volume but not to the healthy tissue. The clinical simplicity, reproducibility, and ability to plan and deliver a uniform dose throughout the tumor volume makes teletherapy the more widely used modality of radiation therapy. In this chapter we are concerned only with teletherapy treatment planning.

The teletherapy radiation source is mounted in a gantry (or arm) that rotates in a circular arc in some plane at a specified distance from the axis of rotation. A patient is positioned on the therapy table so that the center of rotation is placed in the target volume. The target can be irradiated by positioning the gantry at a number of fixed angles around the arc of rotation or by rotating the gantry continuously along segments of the arc.

For teletherapy planning, the contours of the patient boundary, the target, and the critical organs are delineated in one or more contiguous transverse sections of the patient. Only where the anatomy changes rapidly from section to section (as between the head and neck or the neck and shoulders) are several contiguous sections required for accurate calculation of delivered dose, called *dosimetry*. We restrict ourselves to two-dimensional treatment planning, meaning that one patient section is coincident with the plane of gantry motion and all radiation rays are confined to that plane. Extensions to three-dimensional radiation therapy treatment are possible.

The radiation therapist chooses a treatment dose with a high probability of eradicating the tumor without permanent normal tissue injury. The dosimetrist then arranges the treatment fields to cover the target volume while trying to minimize the amount of normal tissue within the irradiated field. Radiation beam modifiers, called wedges or compensators, are added to increase the dose homogeneity across the target volume. The dosimetrist tries to achieve: a uniform dose distribution across the tumor; a minimum dose to the critical organs and other nondiseased tissues; and a minimum of complexity in the technical setup. Most often the experience and judgment of the dosimetrist quickly lead to an efficacious setup.

There are a number of cases, however, that present a challenge, for example when tumors are juxtaposed against critical organs such as the spinal cord, optic nerve, or kidney. The dosimetrist must then consider many possibilities, including fixed fields, arc rotation, beam modifiers, or a combination of these. Achieving a satisfactory treatment plan can therefore become time consuming and tedious. How can we determine if a feasible treatment plan that satisfies all of the radiation therapist's specifications even exists? Should we accept a plan after a preset number of tries if there is still the possibility of a significantly better plan? How do we know when we are close to an optimal plan in these complex cases?

In this chapter we describe a mathematical modeling approach to radiation therapy treatment planning (RTTP) that aims toward automating the selection of an acceptable setup once the physician's specifications are given. The mathematical formulation of the RTTP problem leads to a mathematical problem of operator inversion, which can be approached by either continuous analytical methods or by fully discretizing the model at the outset and resorting to linear algebraic or optimization techniques in a finite-dimensional space. This dichotomy between solution approaches is reminiscent of the situation in image reconstruction from projections (see Chapter 10). The reasons that we favor the full discretization approach will become clear as the discussion unfolds.

## 11.1 Problem Definition and the Continuous Model

Let the plane of the gantry motion be parametrized by polar coordinates $r$ and $\theta$ where $r$ is the distance of a point from the center of rotation of the gantry and $\theta$ is an angle measured clockwise from the positive $y$-axis (see Figure 11.1).

Let $D(r,\theta)$ be a real-valued nonnegative function of the polar coordinates whose value describes the dose absorbed at a point in the patient section coincident with the plane of gantry motion. $D(r,\theta)$ is referred to as the *dose distribution*. We define a *ray* to be a directed line along which radiation travels away from the source, where *source* is a short term for *teletherapy source position*. A ray is specified by the gantry angle $u$ of its source on the gantry circle and by the length $w$ determined by $w = d \tan A$, where $d$ is the radius of the gantry circle and $A$ is the angle of the ray measured clockwise from the central ray. Each ray is outgoing from a source. Thus rays $(u, 0)$ and $(u + \pi, 0)$ represent coincident lines but distinct rays since they have opposite directions. The domain of $u$ is $0 \leq u < 2\pi$, and the domain of $w$ is $-W \leq w \leq W$ where $2W$ is the *beam width.*

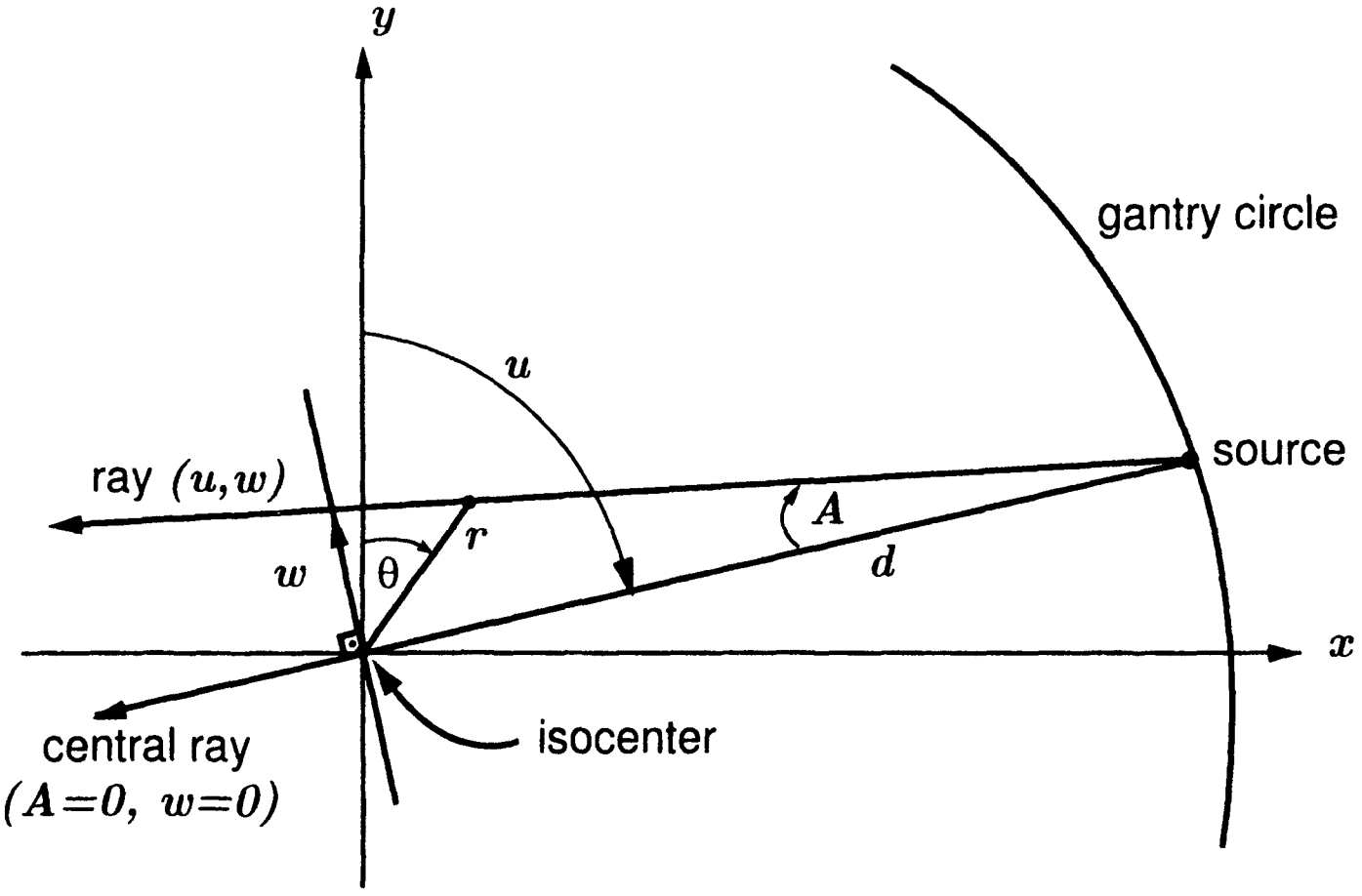

**Figure 11.1** Patient space in the plane of the gantry circle. $(r, \theta)$ denotes a point in patient space. Rays are defined by the gantry angle $u$ of the source and the distance $w$ in the "beam window plane" (shown perpendicular to the central ray at the isocenter). Here ray $(u, w)$ is shown to pass through the point at $(r, \theta)$.

The real-valued nonnegative function $\rho(u, w)$ for $0 \leq u < 2\pi$ and $-W \leq w \leq W$ represents the radiation intensity along the ray $(u, w)$ due to a point source on the gantry circle. We refer to $\rho(u, w)$ as the *radiation intensity distribution.* The continuous version of the *forward* and *inverse* problems of radiation therapy treatment planning can now be formulated.

### 11.1.1 The continuous forward problem

The continuous forward problem of RTTP is the following. Assume that the cross section $\Omega$ of the patient and its radiation absorption characteristics are known. Given a radiation intensity distribution $\rho(u, w)$ for $0 \leq u < 2\pi$ and $-W \leq w \leq W$, find the dose distribution $D(r, \theta)$ for all $(r, \theta) \in \Omega$ from the formula

$$D(r, \theta) = \Delta[\rho(u, w)](r, \theta), \tag{11.1}$$

where $\Delta$ is the *dose operator.* This operator relates the dose distribution to the radiation intensity distribution.

In other words, the forward problem amounts to the calculation of the total dose absorbed at each point of a patient section when all parameters of each radiation beam are specified and the description of the patient section is known. The difficulties associated with the forward problem (i.e., *dose calculation*, or *dosimetry*) stem from the fact that there exists no closed-form analytic representation of the dose operator $\Delta$ that will enable us to use equation (11.1) for the calculation of $D(r, \theta)$. Although the interaction between radiation and tissue is measured and understood at the atomic level, the situation is so complex that, to solve the forward problem in practice, a good state-of-the-art computer program, which represents a *computational approximation* of the operator $\Delta$ and which enables reasonably good dose calculations, must be used.

Let us elaborate on what we intend by stating "there exists no closed-form analytic representation of the dose operator $\Delta$." We actually mean the following: If certain simplifying assumptions are made about the physics of the model as well as the particulars of the desired dose distribution, then it is sometimes possible to express the dose operator in a closed-form analytic formula. This has been done by various authors and some references are given in Section 11.7.

In current practice of RTTP, when dose calculations are performed to verify the dose distribution that will result from a proposed treatment plan, the goal is to obtain results that are as accurate as possible. To achieve this, various empirical data, which are often condensed in look-up tables, are incorporated into the forward calculation. Thus, the true forward calculation, or true dose operator, is not represented by a closed-form analytic relation between the radiation intensity distribution $\rho(u, w)$ and the dose distribution $D(r, \theta)$, but by a software package that calculates $D(r, \theta)$ from $\rho(u, w)$. These arguments are well expressed in the words of Dr. M. Goitein, from Massachusetts General Hospital in Boston. Describing the 3D-CATP (Three-Dimensional Computer-Aided Treatment Planning) program he says:

> Dose calculations are the most critical and most CPU-intensive calculations the system produces. The system must interrogate every pixel in the path of the beam. In order to derive the dose at any point, the system has to take into account the beam's energy profile, and the density of the tissue (gray scale value of the pixels) preceding that point. The system can then look up the dose in a depth-dose curve table. The system also takes into account a scattering factor, through another table look-up. None of the individual calculations are particularly complex, but there are a tremendous number of calculations required. Once the dose distribution calculations have been completed for individual treatments, they can be combined into an overall dose distribution.

Thus, what we really mean by saying that there is no closed-form ana-lytic expression for $\Delta$ is that we choose to adhere to the software representation rather than compromise on the accuracy of the forward calculations by allowing simplifying assumptions that might lead to a closed-form analytic mathematical formula.

### 11.1.2 The continuous inverse problem

The *inverse problem* of radiation therapy is the treatment planning problem:

> Given a description of the patient section, the dose prescribed for the target, and the maximum permissible doses to the target, critical organs, and other tissues, calculate the external configuration and relative intensities of radiation sources (i.e., the radiation field) that will deliver the specified radiation doses (or some acceptable approximation thereof).

Assuming that the cross section $\Omega$ of the patient and its radiation absorption characteristics are known, and given a prescribed dose distribution $D(r,\theta)$, the problem is to find a radiation intensity distribution $\rho(u,w)$ such that equation (11.1) holds, or $\rho(u,w) = \Delta^{-1}[D(r,\theta)]$ where $\Delta^{-1}$ is the inverse operator of $\Delta$. This is the inversion problem that we want to solve, in a computationally tractable way, although no closed-form analytic mathematical representation is available for the dose operation $\Delta$. The dose at $(r,\theta)$ is the sum of the dose contributions from the sources at all the different gantry angles. Thus

$$D(r,\theta) = \sum_{i=1}^{S} y_i D_i(r,\theta), \tag{11.2}$$

where, for each $i = 1,2,\ldots,S$, the value $D_i(r,\theta)$ is the dose deposited at point $(r,\theta)$ by a beam of unit intensity from the $i$th source, and $y_i$ is the total intensity of the $i$th source. It is also possible to interpret $D_i(r,\theta)$ as the dose per unit time deposited at $(r,\theta)$ by the $i$th beam, and $y_i$ as the time the $i$th beam is kept on.

The dose can be further partitioned into two components that are due to primary and scattered radiation. Thus

$$D_i(r,\theta) = D_i^{(pr)}(r,\theta) + D_i^{(sc)}(r,\theta). \tag{11.3}$$

The primary dose $D_i^{(pr)}$ is physically due to the first interaction of beam photons with the tissue medium, and is delivered along the fan of rays of the $i$th beam profile. The value of the primary dose at a point in the patient depends on the distance from the source to the patient surface, the depth of the point, the electron density distribution, the angle between the central ray (i.e., the ray through the center of the gantry circle) and the ray to the point, and the value of the beam profile for the latter ray. The scatter dose $D_i^{(sc)}$ is due to radiation scattered to a point after beam photons first interact with other points in the medium. That is, the scatter dose is due to secondary interactions. The scatter dose at any position depends on the depth of this location in the tissue and the area of irradiated tissue. Both primary and scatter dose also depend on the beam energy, the incident photon spectrum, and the beam modifiers.

For the purposes of our study, it will be assumed here that the dose $D_i(r,\theta)$, and its components $D_i^{(pr)}(r,\theta)$ and $D_i^{(sc)}(r,\theta)$ for each source $i$, can be calculated accurately once the beam parameters and patient section information are specified. That is, we assume that we can solve the forward problem and calculate $D(r,\theta)$ accurately from (11.2) and (11.3). This assumption is confirmed by innumerable direct measurements in water and tissue-equivalent phantoms.

Whereas a dose distribution that solves the forward problem is always obtained for a specified radiation intensity field, the inverse problem may have no solution at all, since some prescribed dose distributions may be unobtainable from any radiation field. Therefore, we do not aim at a solution of the continuous inversion problem, but look instead at a *feasibility formulation* that relaxes the equality in (11.1) in the following manner.

Let $\overline{D} = \overline{D}(r,\theta)$ and $\underline{D} = \underline{D}(r,\theta)$ be two dose distribution functions whose values represent upper and lower bounds, respectively, on the permitted and required dose inside the patient's crosssection. A radiation therapist defines $\overline{D}$ and $\underline{D}$ for each given case and will accept as a solution to the RTTP problem any radiation intensity distribution $\rho(u,w)$ that satisfies

$$\underline{D}(r,\theta) \leq \Delta[\rho(u,w)](r,\theta) \leq \overline{D}(r,\theta), \quad \text{for all } (r,\theta) \in \Omega. \tag{11.4}$$

In target regions (tumors) the lower bound $\underline{D}$ is usually the important factor because the dose there should exceed some given value. In critical organs and other healthy tissues $\underline{D}(r,\theta) = 0$, so that $\overline{D}(r,\theta)$ is the dose that cannot be exceeded. Any solution $\rho(u,w)$ that fulfills (11.4), for given $\overline{D}$ and $\underline{D}$, is a *feasible solution* to the RTTP problem.

## 11.2 Discretization of the Feasibility Problem

In the approach presented here, we adhere to the computerized calculation of the dose operator $\Delta$. Full discretization of the problem at the outset is used to circumvent the difficulties associated with the inversion of $\Delta$. We also neglect the effect of scatter. The patient's cross section $\Omega$ is discretized into a grid of points represented by $\{(r_j,\theta_j) \mid j = 1,2,\ldots,J\}$. Define $\Delta_j\rho$ by

$$\Delta_j\rho \doteq [\Delta\rho](r_j,\theta_j) \tag{11.5}$$

and call $\Delta_j$ a *dose functional*, for every $j = 1, 2, \ldots, J$. Acting on a radiation intensity distribution $\rho(u, w)$, the functional $\Delta_j$ provides $\Delta_j \rho$, which is the dose absorbed at the $j$th grid point of the patient's cross section $\Omega$ due to the radiation intensity field $\rho$. The dose distributions $\overline{D}$ and $\underline{D}$ are specified at the grid points by giving, for all $j = 1, 2, \ldots, J$,

$$\overline{D}(r_j, \theta_j) = \overline{D}_j, \quad \underline{D}(r_j, \theta_j) = \underline{D}_j, \tag{11.6}$$

thus converting (11.4) into a finite system of *interval inequalities*

$$\underline{D}_j \leq \Delta_j \rho \leq \overline{D}_j, \quad j = 1, 2, \ldots, J. \tag{11.7}$$

Denoting hereafter by $\overline{D}$ ($\underline{D}$) the $J$-dimensional column vector whose $j$th element is $\overline{D}_j$ ($\underline{D}_j$), the inverse problem of RTTP is restated as follows:

> Given vectors $\overline{D}$ and $\underline{D}$ of permitted and required doses at $J$ grid points in the patient's cross section $\Omega$, find a radiation intensity distribution $\rho = \rho(u, w)$ such that (11.7) holds.

In continuing the discretization process of the problem it is assumed that a set of $I$ *basis radiation intensity fields* is fixed and that their non-negative linear combinations can give adequate approximations to any radiation intensity field we wish to specify. This is done by discretizing the region $0 \leq u < 2\pi$, $-W \leq w \leq W$ in the $(u, w)$-plane into a grid of points given by $\{(u_i, w_i) \mid i = 1, 2, \ldots, I\}$. A radiation intensity distribution

$$\sigma_i(u, w) \doteq \begin{cases} 1, & \text{if } (u, w) = (u_i, w_i), \\ 0, & \text{otherwise}, \end{cases} \tag{11.8}$$

is a *unit intensity ray* and serves as a member of the set of basis intensity fields, $i = 1, 2, \ldots, I$. A desired radiation intensity distribution $\rho$ that solves (11.7) is approximated by

$$\hat{\rho}(u, w) = \sum_{i=1}^{I} x_i \sigma_i(u, w), \tag{11.9}$$

where $x_i$ is the intensity of the $i$th ray, and it is required that $x_i \geq 0$, for all $i = 1, 2, \ldots, I$. Once the grid points are fixed, any radiation intensity distribution $\hat{\rho}$ that can be presented as a nonnegative linear combination of the rays is uniquely determined by the coefficients $x_i$, $1 \leq i \leq I$. The vector $x = (x_i)$, in the $I$-dimensional Euclidean space $\mathbb{R}^I$ is referred to as the *radiation vector* or *basic solution.*

Further, assume that the dose functionals $\Delta_j$ are linear and continuous. This assumption cannot be mathematically verified due to the absence of an analytic representation of $\Delta$ or $\Delta_j$, but it is a reasonable assumption based on the empirical knowledge of $\Delta_j$. Using linearity and continuity of all $\Delta_j$'s, we can write $\Delta_j\rho \simeq \Delta_j\hat{\rho} = \sum_{i=1}^{I} x_i\Delta_j\sigma_i$. For $j = 1, 2, \ldots, J$, and $i = 1, 2, \ldots, I$, denote by

$$d_{ij} \doteq \Delta_j\sigma_i \tag{11.10}$$

the dose deposited at the $j$th point $(r_j, \theta_j)$ in the patient's cross section $\Omega$ due to a unit intensity ray $\sigma_i(u, w)$. The *fully discretized feasibility inverse problem* of RTTP then becomes the linear interval feasibility problem of finding a vector $x \in \mathbb{R}^I$ such that

$$\begin{aligned} \underline{D}_j \leq \sum_{i=1}^{I} x_i d_{ij} \leq \overline{D}_j, \quad j = 1, 2, \ldots, J, \\ x_i \geq 0, \quad i = 1, 2, \ldots, I. \end{aligned} \tag{11.11}$$

Let the set of pixels in the discretized patient cross section be denoted by $N = \{1, 2, \ldots, J\}$. Organs within the patient section are then defined as subsets of $N$. The subsets $B_k \subset N$, where $k = 1, 2, \ldots, K$ denote $K$ critical organs to be spared from excessive radiation. Let the values $b_k$ denote the corresponding upper bounds on the dose permitted in each critical organ. The subsets $T_q \subset N$ where $q = 1, 2, \ldots, Q$ denote $Q$ target regions. Let the values $t_q$ denote the corresponding prescribed lower bounds for the absorbed dose in each. All the $B_k$ and $T_q$ are pairwise disjoint. The set of pixels inside the patient section that are not in any $B_k$ or $T_q$ are called the *complement*, denoted as the subset $C \subset N$, and $c$ is the upper bound for the total permitted dose there. It is assumed that the definition of all subsets $B_k, T_q$, and $C$ and the prescription of all $b_k, t_q$, and $c$ are given by the radiotherapist as input data for the discretized treatment planning problem.

Problem (11.11) then becomes the following system of linear inequalities, which we call the *basic model*:

$$\sum_{i=1}^{I} d_{ij}x_i \leq b_k, \quad \text{for all } j \in B_k,\ k = 1, 2, \ldots, K, \tag{11.12}$$

$$t_q \leq \sum_{i=1}^{I} d_{ij}x_i, \quad \text{for all } j \in T_q,\ q = 1, 2, \ldots, Q, \tag{11.13}$$

$$\sum_{i=1}^{I} d_{ij}x_i \leq c, \quad \text{for all } j \in C, \tag{11.14}$$

$$x_i \geq 0, \quad \text{for all } i = 1, 2, \ldots, I. \tag{11.15}$$

With $b_k, t_q$, and $c$ given and the $d_{ij}$'s calculated from (11.10), the mathematical question represented by the basic model (11.12)–(11.15) is to find a nonnegative solution vector $x = (x_i)$ for a system of linear inequalities. This is a *linear feasibility problem.* Any vector $x$ whose components $x_i$ satisfy the basic model is acceptable in terms of the dose specifications set forth. Consult Figure 11.2 in which the individual ray intensities $x_i$ are double-indexed to identify the number of the source from which they are derived in the discretization process. With the discretization of source positions on the gantry arc as described above, we may consider a *discretized ray space*, namely, a rectangular grid of $I$ points (in the nonnegative orthant of the ray space), with $S$ equally spaced locations in the interval $0 \leq u < 2\pi$ along the $u$-axis and $M$ equally spaced locations placed in the interval $-W \leq w \leq W$ along the $w$-axis. We enclose the points in the ray space in pixels and assign intensity $x_i$ to the $i$th pixel, where $i = i(s, m)$ according to the lexicographic ordering formula $i = m + (s - 1)M$ (see Figure 11.3).

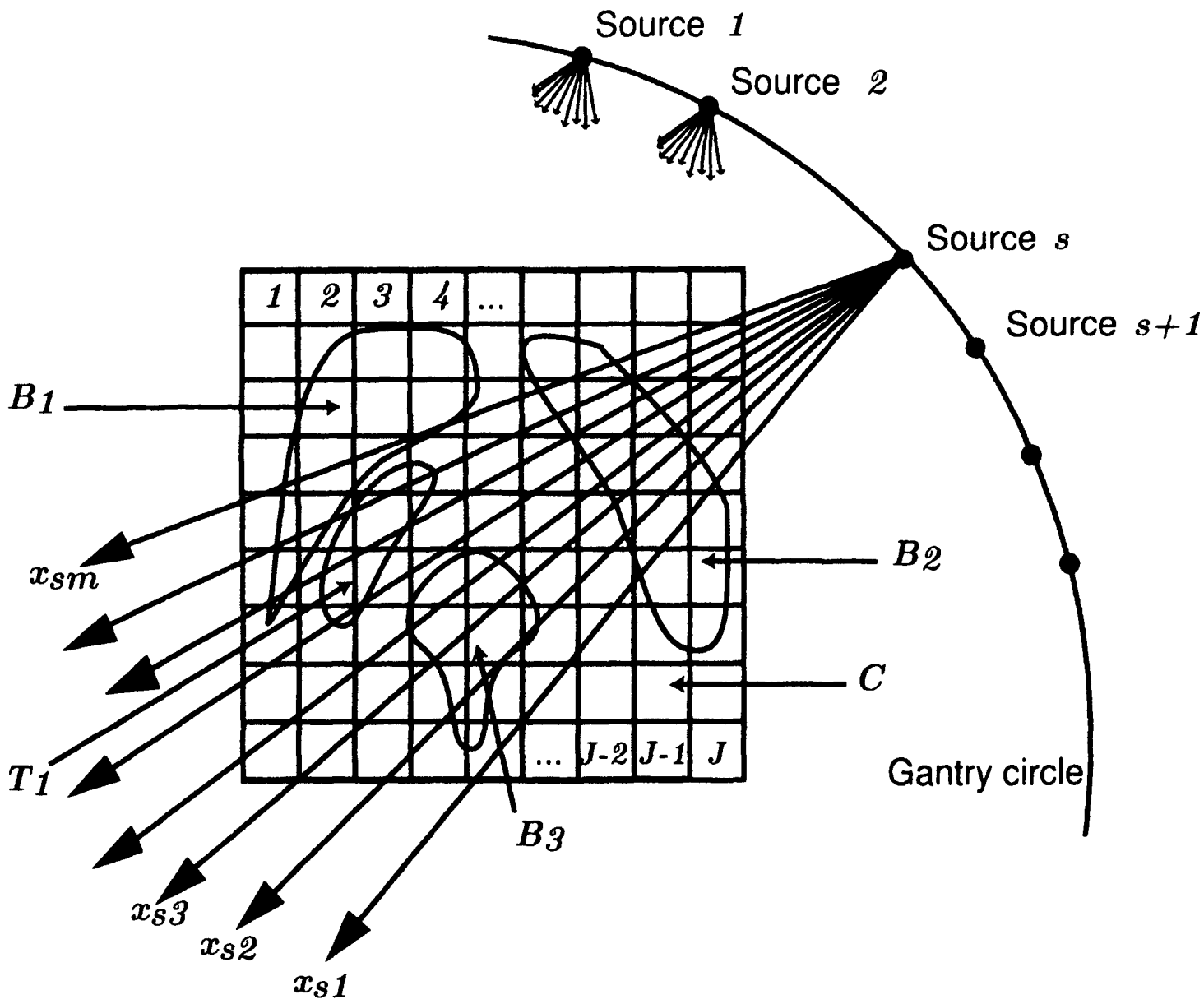

**Figure 11.2** Discretization of patient section into pixels, and discretization of the radiation field into sources and rays.

For this problem to be well defined the quantities $d_{ij}$ must be precalculated, i.e., by apportioning the calculated dose per unit beam intensity between the rays representing that beam in the discretization. A forward problem solver in the form of a computer program is assumed available to calculate each specific dose $d_s(j)$, which is the dose per unit time, absorbed at pixel $j$ of the patient section when the radiation field is due to a radi-

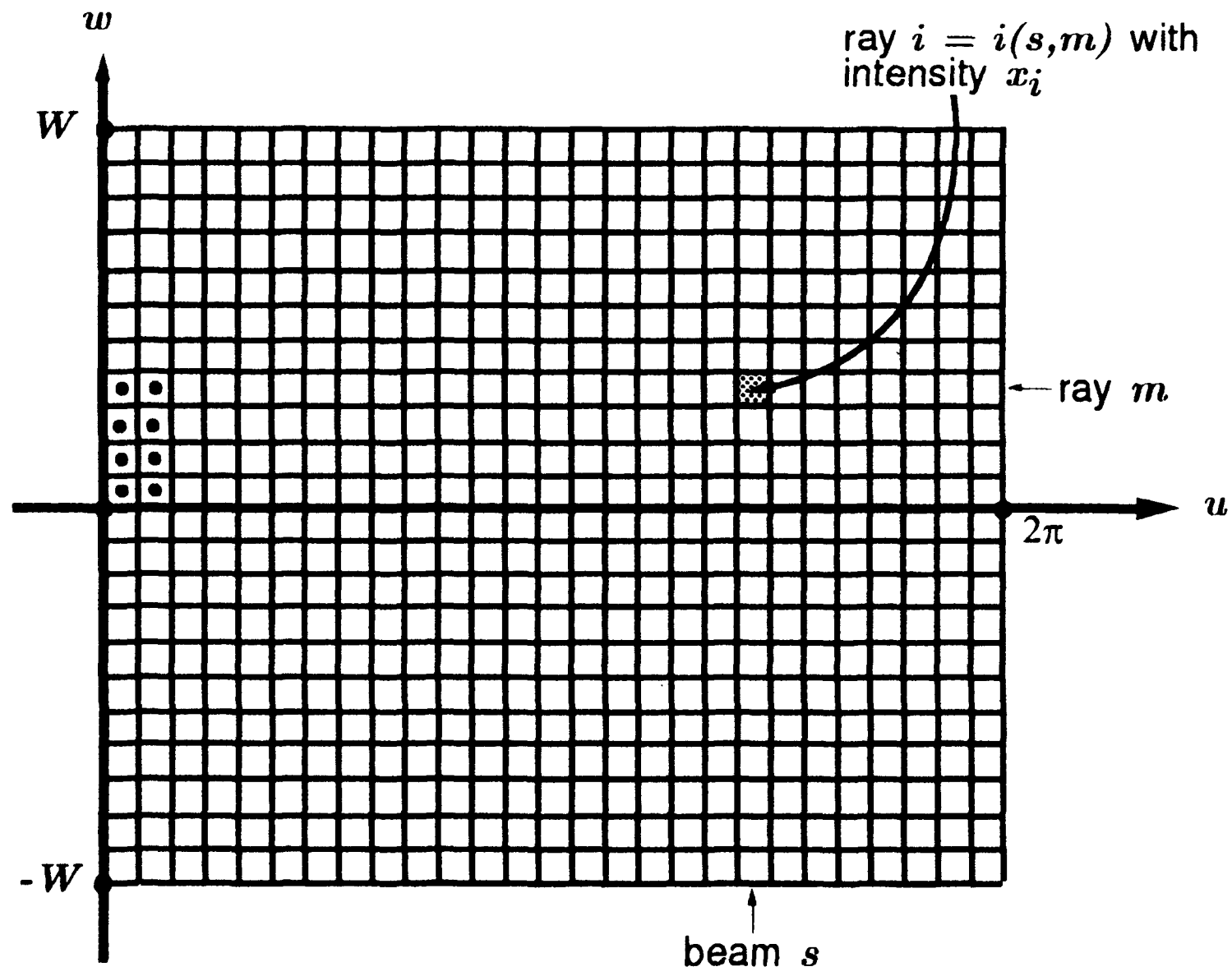

**Figure 11.3** Discretized ray space wherein each column represents a single beam.

ation beam from source $s$. For our discretized radiation field, the gantry angle $u$ assumes only the values

$$u_s = 2\pi\frac{s-1}{S}, \quad s = 1, 2, \ldots, S. \tag{11.16}$$

We must now attribute the specific dose at pixel $j$ to a contribution from each of the rays in our discretized ray space. We thus need to calculate the numbers $d_{ij}$, $i = 1, 2, \ldots, I$, $j = 1, 2, \ldots, J$, which represent the specific dose absorbed in the $j$th pixel due to radiation from the $i$th ray alone. We have already decided to consider only the primary component of the radiation dose; therefore we use a dose apportionment scheme which is now defined and explained (see Figure 11.4).

Assume that a source of radiation is placed at gantry angle $u_s$, and that the value $d_s(j)$ has been calculated for each pixel $j$ by means of a forward problem solver. This value is apportioned among a chosen set of rays, and at present, only rays that pass through or immediately straddle the pixel are used for apportionment. A possible apportionment scheme results from the formula

$$d_{ij} = \frac{1/\sigma_{ij}}{\sum_p 1/\sigma_{pj}} d_s(j), \tag{11.17}$$

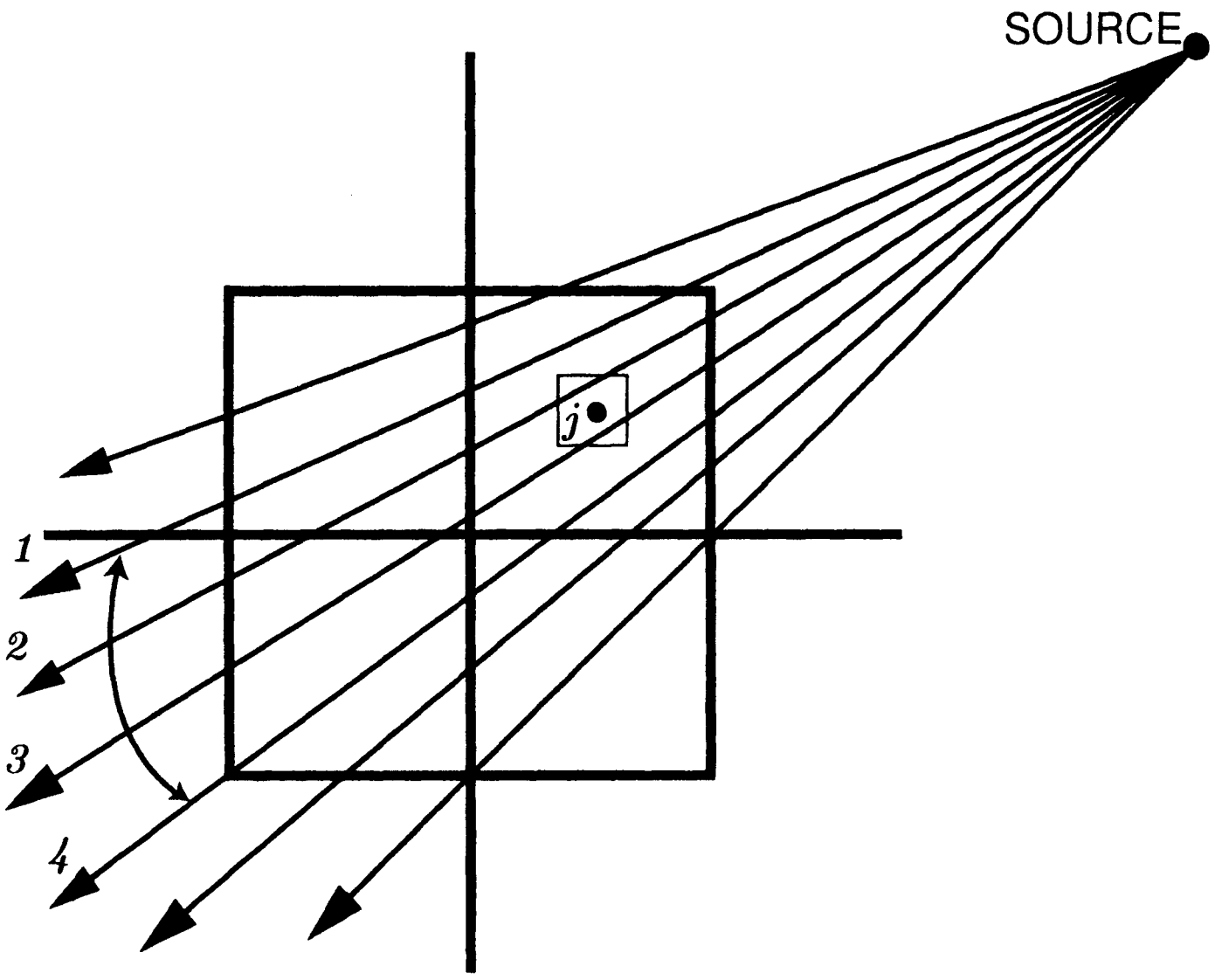

**Figure 11.4** Specific dose absorbed in pixel $j$ due to a source is apportioned among rays labeled 1 to 4 (in this case) to obtain the quantities $d_{ij}$.

where $\sigma_{ij}$ is the distance from the center of the $j$th pixel to the $i$th ray. The summation over $p$ goes up to the number $M_s(j)$ of rays from source $s$ at gantry angle $u_s$ that pass through or immediately straddle the $j$th pixel. If the pixel center lies on a ray whose index is $t$ so that $\sigma_{tj} = 0$, then we take $d_{tj} = d_s(j)$ and set $d_{ij} = 0$ for rays $i \neq t$. In computational trials, however, the apportionment scheme of (11.17) was found to be too sensitive to small changes in the planar discretization and ray-sampling schema. A more satisfactory apportioning scheme was found to be

$$d_{ij} = \frac{d_s(j)}{M_s(j)}. \tag{11.18}$$

This simple scheme is the one used to obtain the numerical results demonstrated in Section 11.5.

## 11.3 Computational Inversion of the Data

A procedure that solves the fully discretized inversion problem formulated above can be constructed in the following manner:

**Step 1.** Read as input the parameters necessary for the forward problem solver (dose calculation program). These are the beam width, beam energy, calculation tables, etc. Read as input the patient boundary contour and the contours of the internal structures (targets and critical organs). Read as input the physician's prescription for the dose distribution.

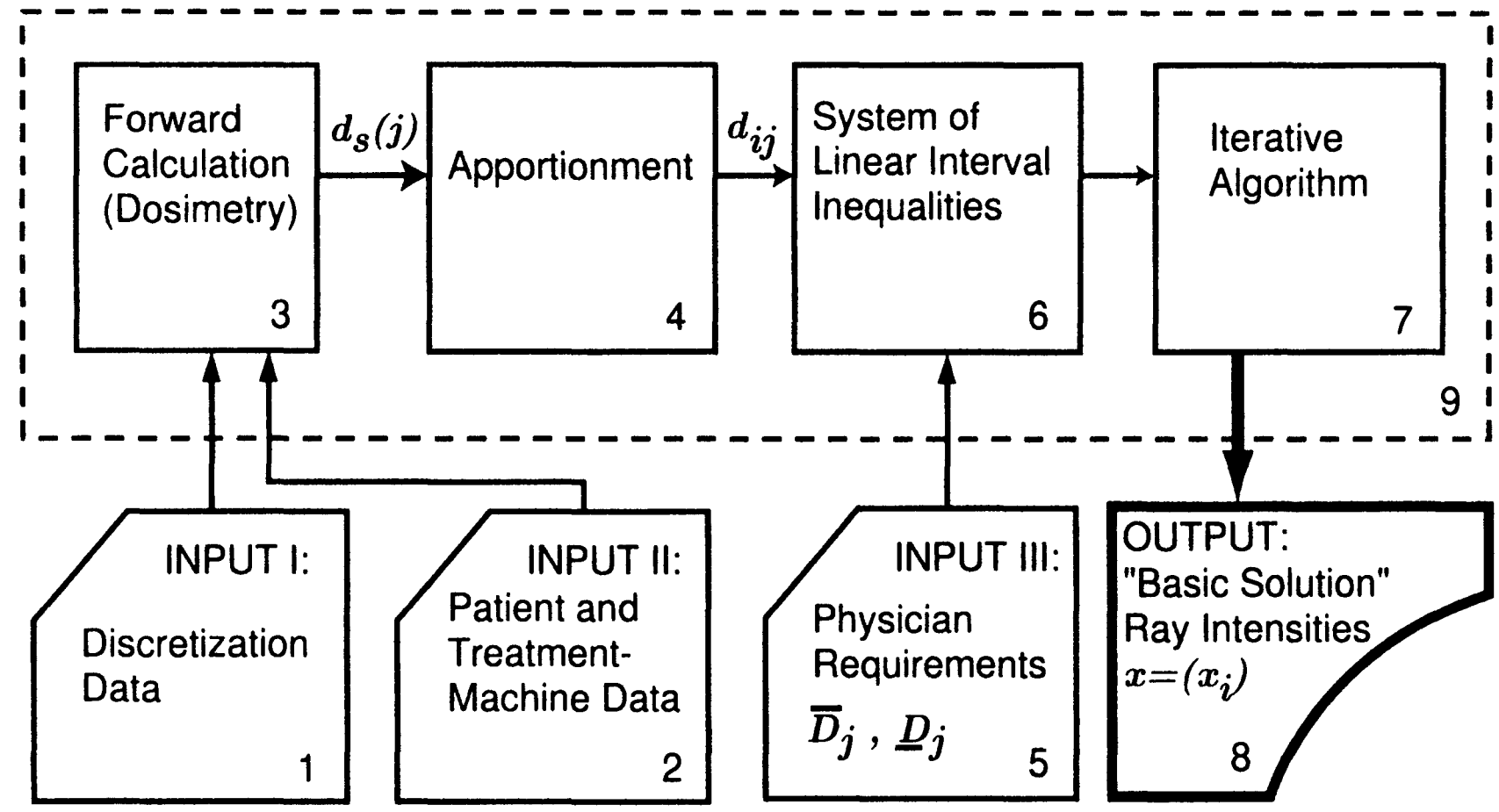

**Figure 11.5** Derivation of the *basic solution* of ray intensities from physician requirements, patient and machine data, and discretization.

**Step 2.** Choose a discretization grid for the cross section, and determine whether each pixel lies within a tumor region, within a critical organ, or elsewhere. Assign to each pixel the appropriate values of the physician-chosen bounds on radiation dose.

**Step 3.** Choose the number $S$ of possible gantry (source) positions. Calculate and save the specific dose distribution $(d_s(j))_{j=1}^{J}$ for each gantry position of the beam. No blocks, wedges, or filters are used to modify the beam profile. Since scatter radiation is omitted here, this forward dose calculation is quite fast.

**Step 4.** Choose the number of rays per beam for discretization.

**Step 5.** For each pixel, the specific dose due to the source at gantry angle $u_s$, $s = 1, 2, \ldots, S$ is apportioned among those rays passing through or straddling the pixel. This dose apportionment is applied to the dose distributions calculated in Step 3 for all the sources. In this way the numbers $d_{ij}$ are calculated.

**Step 6.** Use an iterative algorithm to solve the resulting linear feasibility problem. Initially all rays are assigned zero intensity. During the iterative process, the algorithm gradually changes these intensities as it tries to satisfy the prescribed physician constraints. Rays that do not cross any tumor region are permanently assigned zero intensity.

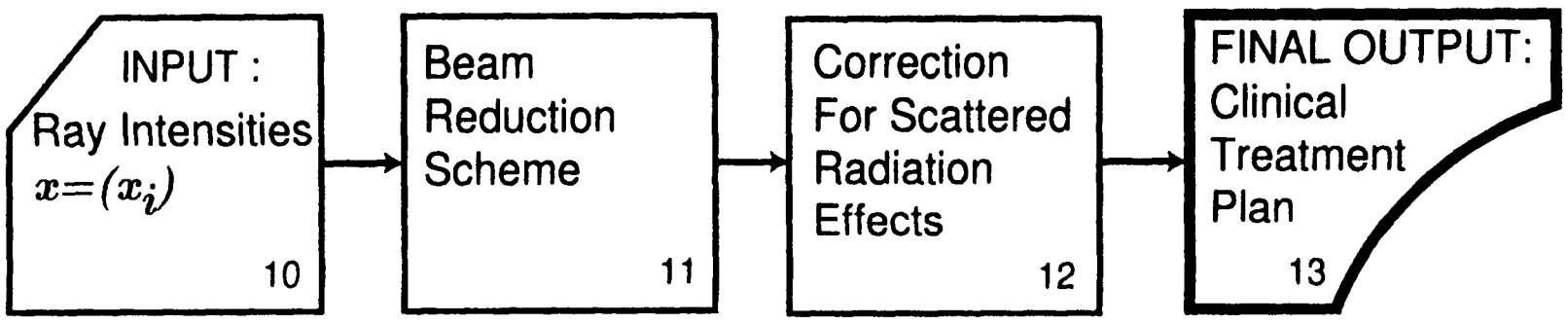

**Figure 11.6** Derivation of the *clinical treatment plan* from a *basic solution* of ray intensities.

Figures 11.5 and 11.6 describe the overall process. All boxes in these diagrams are numbered for easy reference. Box 3 represents and assumes the availability of a state-of-the-art computer program for forward calculation. Given discretization data for pixels, beams, and rays (box 1), and data on the patient's cross section and the treatment machine parameters (box 2), it calculates the numbers $d_s(j)$, $s = 1, 2, \ldots, S$, $j = 1, 2, \ldots, J$. Each $d_s(j)$ is the dose absorbed at the $j$th location of the cross section (i.e., at $(r_j, \theta_j)$) due to a unit intensity of radiation from source $s$ on the gantry circle. The dose apportionment scheme of box 4 distributes these values among individual rays. The resulting $d_{ij}$'s are the coefficients of the system (11.11), which also needs the physician requirements from box 5. An iterative algorithm (box 7) then produces an approximate basic solution of ray intensities $x = (x_i)_{i=1}^{I}$.

If a treatment machine that could deliver pencil-thin single rays of controlled intensity existed, then the basic solution could have been implemented clinically. Since this is not yet the case, we use this solution as input to the process described in Figure 11.6, which first employs a beam reduction scheme that extracts from the basic solution a clinically acceptable treatment plan. After reducing the number of beams, one would need to correct the plan to incorporate the effect of scattered radiation because initially the calculations in the system of Figure 11.5 were for primary dose only.

## 11.4 Consequences and Limitations

If a mathematically feasible solution to the radiation therapy treatment planning problem exists, i.e., if the system (11.12)–(11.15) is feasible, then the final set of ray intensities, obtained after stopping the algorithm, provides an approximation to a mathematically feasible radiation field. On the other hand, if a feasible solution to the treatment planning problem does not exist, i.e., if the system (11.12)–(11.15) is not consistent, then the final set of ray intensities could provide a possible compromise assuming that the nonnegativity constraints (11.15) are satisfied. With the final set of ray intensities denoted by the vector $x^* = (x_i^*)$, the question is how to implement this derived radiation field as a treatment plan.

A straightforward implementation of radiation field $x^*$ calls for a machine that will be able to emit controlled amounts of radiation in pencil-thin rays and repeat this in a multitude of directions. The number of directions, determined by the fineness of the grid for radiation field discretization, may range into the thousands. Available radiation treatment machines, however, use finite-width beam sources. Although a beam profile may be modified to some extent by using certain accessories such as blocks, wedges, and compensators, we cannot expect present day machines to implement accurately our finely discretized radiation fields.

While the technology may sufficiently improve to meet such a challenge, we are studying the mathematical problem of how to extract from the ray intensities solution vector $x^*$ a set of physical parameters such as: number of sources (beams), gantry angles, specification of beam modifiers, and intensities of sources (beams) that will constitute a clinically acceptable and implementable treatment plan. Figure 11.3 shows that the discretized radiation field $x^*$ can be considered a distribution of values $x_i^*$ (or $x_{sm}^*$) over the pixels $(s, m)$ that represent the rays in this figure. To be clinically implementable, a derived radiation field should have nonzero values in only a few (say up to 10) columns in Figure 11.3 (each representing a single source or gantry position), and only a linear or piecewise linear change of ray intensities over a succession of rays within a column (to allow wedge or block beam modifications). To meet these requirements, a method must be devised to extract a blockwise varying radiation field from the finely discretized radiation field.

At present we examine the distribution of ray intensities $x_{sm}^*$ (or $x_i^*$) and choose those source positions $u_s$ that provide the greatest ray intensity contributions. We then repeat the feasibility calculation described above for a limited number (not necessarily uniformly spaced) of source positions around the gantry arc.

## 11.5 Experimental Results

The system (11.11) constitutes the mathematical model (box 6 in Figure 11.5) that we adopt for computing a solution of the inverse problem in RTTP. The physician's requirements are given by the lower and upper bounds, $\underline{D}_j$ and $\overline{D}_j$, $j = 1, 2, \ldots, J$. The coefficients $d_{ij}$ are calculated with the aid of a state-of-the-art forward problem solver, according to the principles stated above. We apply to this problem the fully simultaneous Cimmino algorithm. This is the block-AMS algorithm (Algorithm 5.6.2) in which there is only one block ($M = 1$) that contains all inequality constraints. The Cimmino algorithm (for box 7 in Figure 11.5) was tested for the two-dimensional case of a single transverse section of a patient with beams confined to the plane of the section. No algorithmic changes are required to extend the applicability to three-dimensional problems in radiation therapy, although numbering of volume elements (voxels) and rays (the input to the implementation), and the graphic presentation of the calculated dose in a three-dimensional volume (the output of the implementation) would differ.

To use the algorithm the clinician must specify: the contours for the patient, target, and critical organs; the dose prescription for the target, together with a minimum and a maximum dose over it; the maximum dose allowed for specified critical organs; and the maximum dose allowed for the remaining noncritical tissue. In addition, the clinician may assign relative weights (of importance) to each region within the cross section, such that the algorithm generates a solution that least violates the constraints of the heavily weighted regions and distributes proportionately more error to the least weighted regions.

Figure 11.7 shows a computed tomography (CT) section of a patient's skull. In Figure 11.8 the picture is discretized into squares with dimensions $0.85 \times 0.85$ cm (not drawn to scale). The region marked by x is the target containing the tumor. The dose to be delivered to the target is between 50 and 53 Gy. The region marked with 2 corresponds to the brain stem, which is a critical organ. The dose there is not to exceed 40 Gy. Regions marked with 3 and 4 are the left and right eyes respectively, and are not to receive any radiation. Finally, the dose delivered to the remaining unlabelled tissue should not exceed 45 Gy (see Figure 11.9).

The constraint (or regional) weights for the Cimmino algorithm are chosen to be fixed, i.e., in Algorithm 5.6.2 the weights $w^{\nu}(i) = w(i)$ are independent of the iteration index $\nu$ and have the same value for all pixels $i$ that belong to a specific organ region in the patient's cross section. Note that in the fully discretized model of RTTP each inequality represents a pixel in the cross section, whereas in image reconstruction discretized models (Chapter 10) each constraint (equality or inequality) comes from an individual ray traversing the object. The weights $w(i)$ chosen here are: ten for the target region and for the eyes, three for the brain stem, and one for the unlabelled tissue. We start the discretization of the model with 24 uniformly spaced beam positions, calculate the intensities of 50 rays for each beam (a total of $50 \times 24$ rays), and then plot the resulting dose distribution. The 50 rays of each beam are set to span the target, which means that rays which do not intersect the target are set, and kept throughout, at zero intensity.

In clinical situations, complex multibeam treatment plans are generally infeasible and the number of beam positions must be limited to fewer than five or six. To this end, we describe our beam reduction scheme (box 11 in Figure 11.6), which leads to our clinical treatment plan (box 13 in Figure 11.6). Correction for scattered radiation effects (box 12 in Figure 11.6) has not been incorporated in the experiments described here. Ideally, one would hope that a reduced number of beams, suitable for a practical clinical treatment plan, could be read directly from the basic solution vector of ray intensities $x^* = (x_i^*)$. Since we are unaware of a method that will directly reduce beams, we repeatedly eliminate certain beams from the original 24 beams, used for the initial discretization of the full radiation field, and then repeat the computational process (described in Figure 11.5) until a satisfactory plan is achieved. Thus, we use the Cimmino algorithm to select systematically and iteratively a small number of beam positions that can deliver a clinically acceptable dose distribution to the patient. This is done in the following manner. Once we have at hand a basic solution $x^* = (x_i^*)$ with the $24 \times 50$ ray intensities $x_i^*$, we calculate for each beam the sum of the individual ray intensities (this number is proportionate to the dose contribution of the whole beam) and the sum of the absolute values of the differences of adjacent ray intensities (this number is indicative of the beam smoothness across its profile). Because the ray intensity sums of the 24 beams vary considerably, we separate the beams into two groups. The high ray-intensity group has 11 beams; thus the other 13 beams can be eliminated. Applying the Cimmino algorithm a second time to the 11 remaining beams, on the basis of ray intensity sums alone we can determine six beams with higher such sums. When the Cimmino algorithm is applied a third time—to the remaining six beams—the resulting ray intensity sums are all very similar and moreover, the beams are adjacent along an arc from 45° to 120° (see Table 11.1). We then use the criterion of smoothness across the target, and eliminate two beams with significantly large adjacent difference sums. A fourth application of Cimmino's algorithm to the four remaining beams shows that one beam (the one at 45°) can be eliminated on the basis of ray sum and smoothness. The fifth application of Cimmino's algorithm to three beams shows virtually the same ray sums and adjacent difference sums for each beam; and a sixth application of the algorithm gives the ray intensities for a two-beam treatment plan.

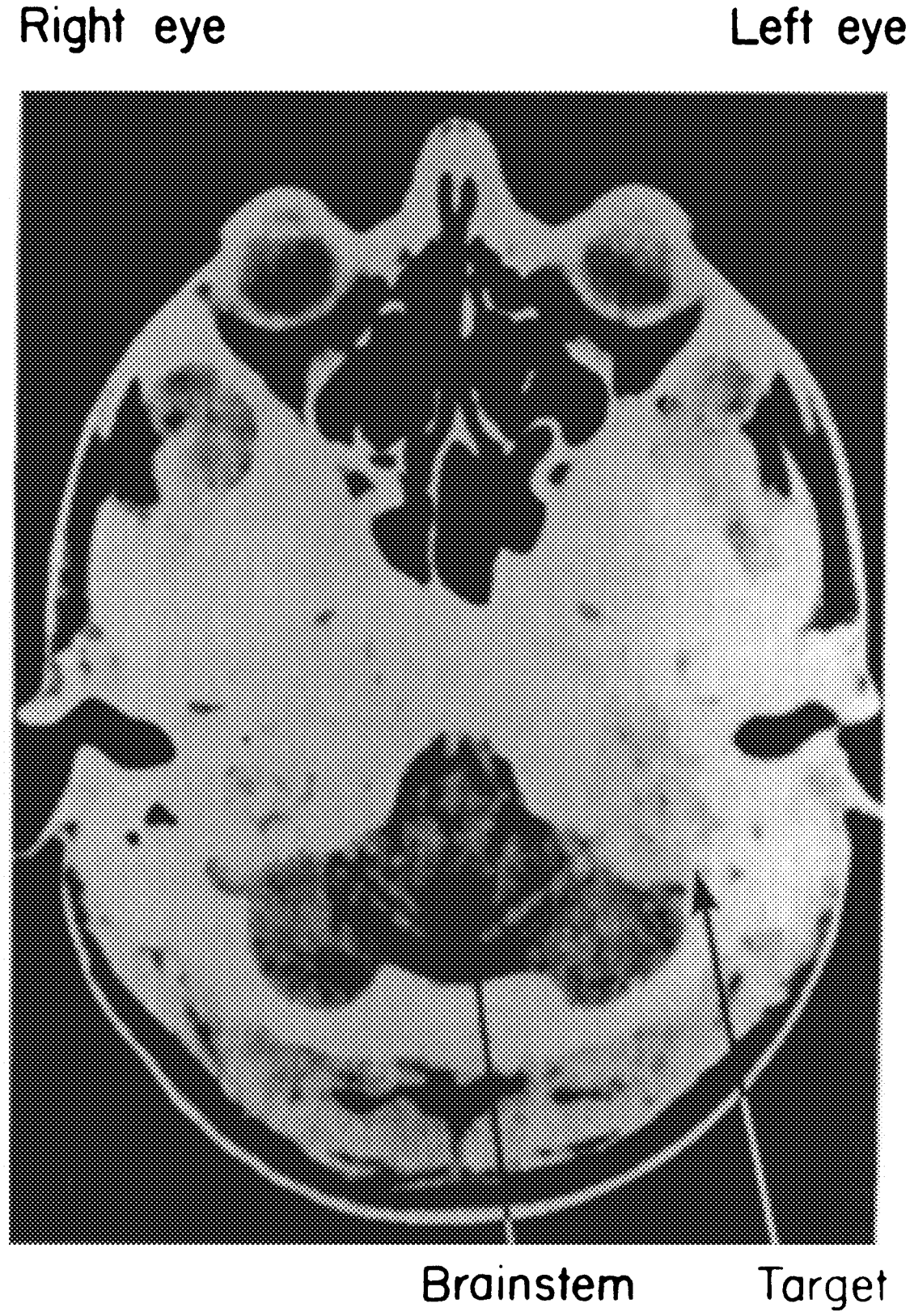

**Figure 11.7** A CT section of the patient to be treated. Our experimental results are demonstrated on this case. (Reproduced from Censor, Altschuler, and Powlis 1988).

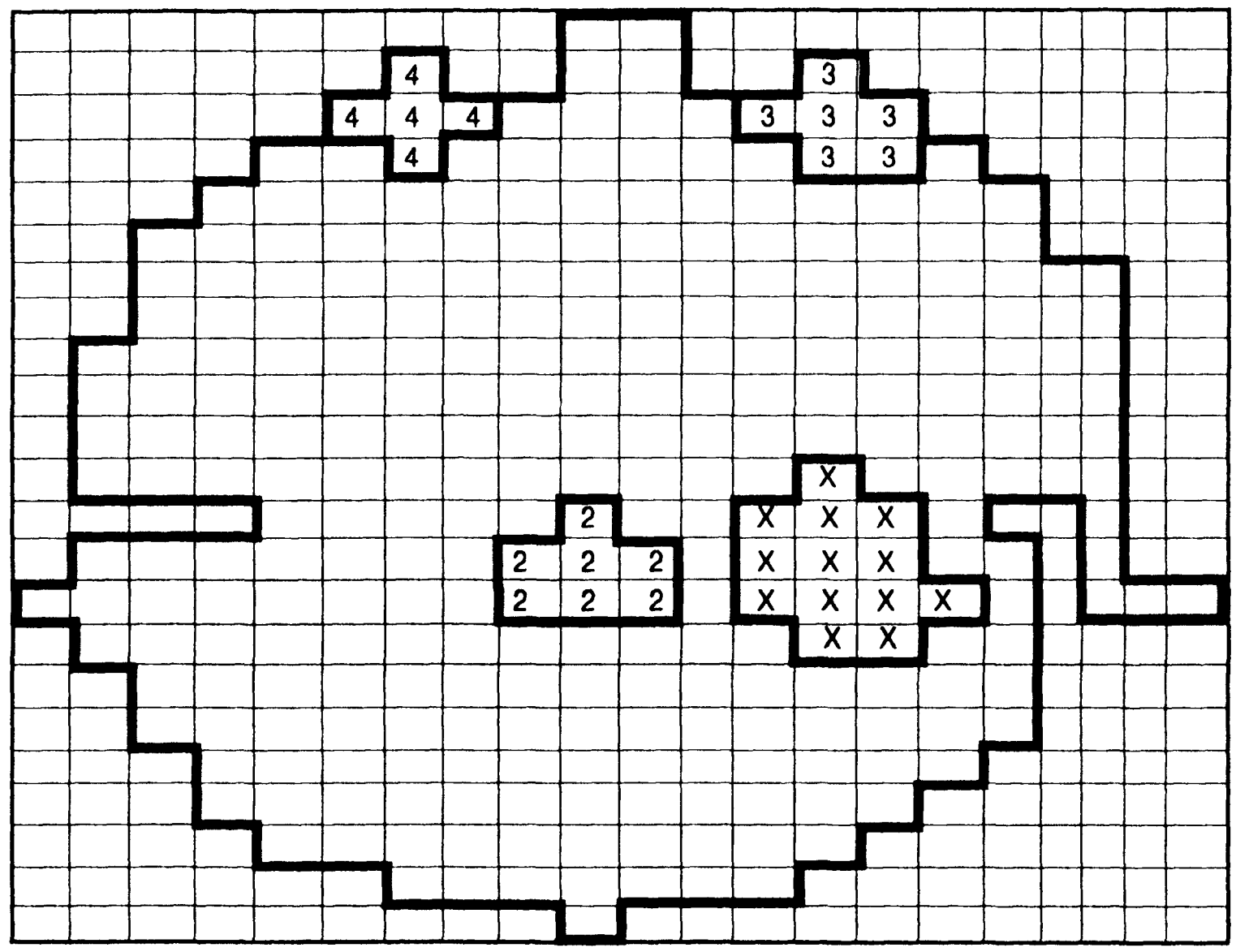

**Figure 11.8** The discretized section derived from Figure 11.7 (with the pixels compressed, for convenience of display, in the vertical direction).

In this way we have used two empirical criteria to reduce the number of beams in the final treatment plan. The more important criterion is a threshold on the sum of ray intensities of the beams and it can be used if the ray intensity sums are sufficiently different. The secondary criterion is the choice of beams with least sum of absolute values of differences of intensities of adjacent rays. This criterion helps to differentiate between beams with similar sums of ray intensities.

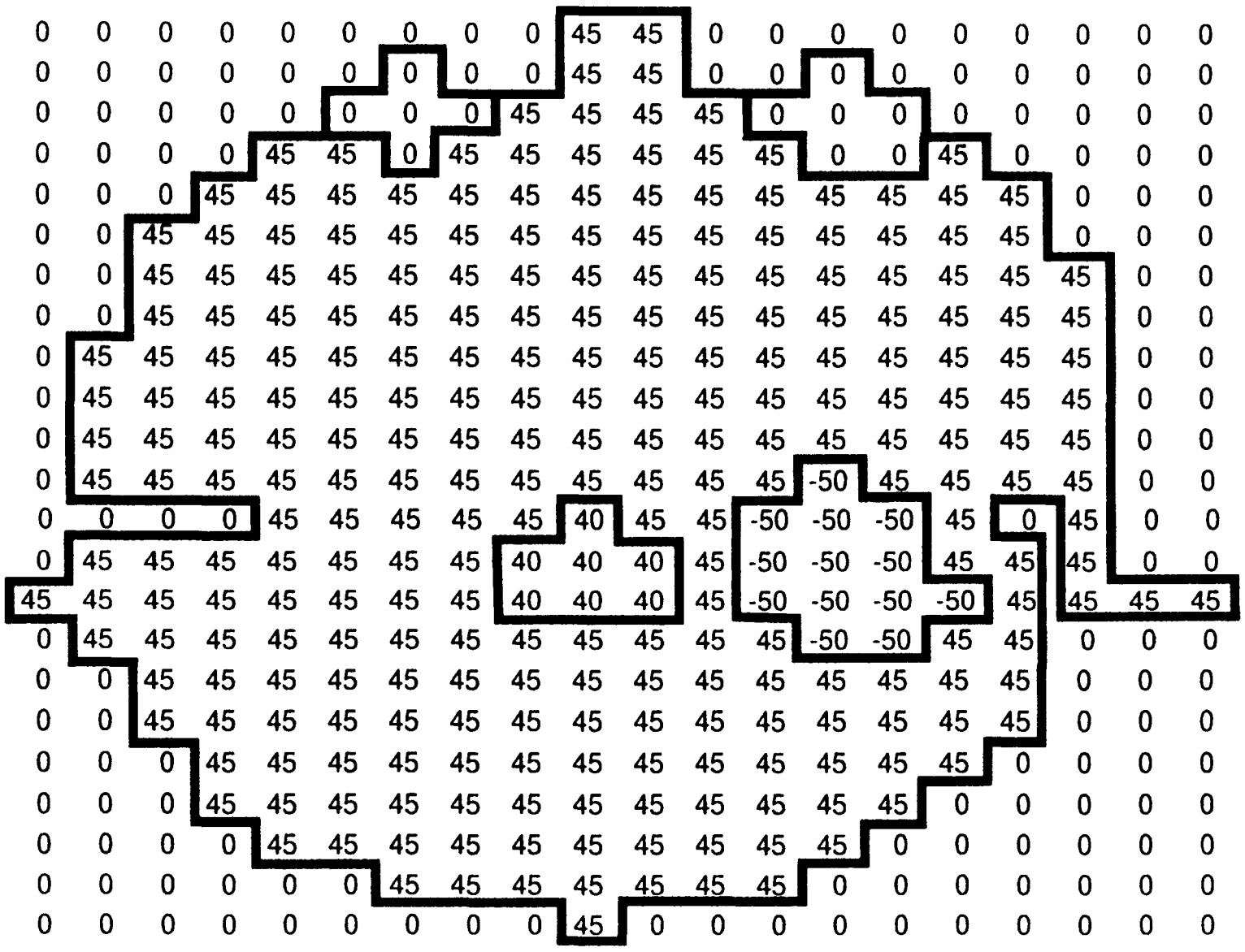

| 0 | 0 | 0 | 0 | 0 | 0 | 0 | 0 | 0 | 45 | 45 | 0 | 0 | 0 | 0 | 0 | 0 | 0 | 0 | 0 |
|---|---|---|---|---|---|---|---|---|---|---|---|---|---|---|---|---|---|---|---|
| 0 | 0 | 0 | 0 | 0 | 0 | 0 | 0 | 0 | 45 | 45 | 0 | 0 | 0 | 0 | 0 | 0 | 0 | 0 | 0 |
| 0 | 0 | 0 | 0 | 0 | 0 | 0 | 0 | 45 | 45 | 45 | 45 | 0 | 0 | 0 | 0 | 0 | 0 | 0 | 0 |
| 0 | 0 | 0 | 0 | 45 | 45 | 0 | 45 | 45 | 45 | 45 | 45 | 45 | 0 | 0 | 45 | 0 | 0 | 0 | 0 |
| 0 | 0 | 0 | 45 | 45 | 45 | 45 | 45 | 45 | 45 | 45 | 45 | 45 | 45 | 45 | 45 | 45 | 0 | 0 | 0 |
| 0 | 0 | 45 | 45 | 45 | 45 | 45 | 45 | 45 | 45 | 45 | 45 | 45 | 45 | 45 | 45 | 45 | 0 | 0 | 0 |
| 0 | 0 | 45 | 45 | 45 | 45 | 45 | 45 | 45 | 45 | 45 | 45 | 45 | 45 | 45 | 45 | 45 | 45 | 0 | 0 |
| 0 | 0 | 45 | 45 | 45 | 45 | 45 | 45 | 45 | 45 | 45 | 45 | 45 | 45 | 45 | 45 | 45 | 45 | 0 | 0 |
| 0 | 45 | 45 | 45 | 45 | 45 | 45 | 45 | 45 | 45 | 45 | 45 | 45 | 45 | 45 | 45 | 45 | 45 | 0 | 0 |
| 0 | 45 | 45 | 45 | 45 | 45 | 45 | 45 | 45 | 45 | 45 | 45 | 45 | 45 | 45 | 45 | 45 | 45 | 0 | 0 |
| 0 | 45 | 45 | 45 | 45 | 45 | 45 | 45 | 45 | 45 | 45 | 45 | 45 | 45 | 45 | 45 | 45 | 45 | 0 | 0 |
| 0 | 45 | 45 | 45 | 45 | 45 | 45 | 45 | 45 | 45 | 45 | 45 | 45 | -50 | 45 | 45 | 45 | 45 | 0 | 0 |
| 0 | 0 | 0 | 0 | 45 | 45 | 45 | 45 | 45 | 40 | 45 | 45 | -50 | -50 | -50 | 45 | 0 | 45 | 0 | 0 |
| 0 | 45 | 45 | 45 | 45 | 45 | 45 | 45 | 40 | 40 | 40 | 45 | -50 | -50 | -50 | 45 | 45 | 45 | 0 | 0 |
| 45 | 45 | 45 | 45 | 45 | 45 | 45 | 45 | 40 | 40 | 40 | 45 | -50 | -50 | -50 | -50 | 45 | 45 | 45 | 45 |
| 0 | 45 | 45 | 45 | 45 | 45 | 45 | 45 | 45 | 45 | 45 | 45 | 45 | -50 | -50 | 45 | 45 | 0 | 0 | 0 |
| 0 | 0 | 45 | 45 | 45 | 45 | 45 | 45 | 45 | 45 | 45 | 45 | 45 | 45 | 45 | 45 | 45 | 0 | 0 | 0 |
| 0 | 0 | 45 | 45 | 45 | 45 | 45 | 45 | 45 | 45 | 45 | 45 | 45 | 45 | 45 | 45 | 45 | 0 | 0 | 0 |
| 0 | 0 | 0 | 45 | 45 | 45 | 45 | 45 | 45 | 45 | 45 | 45 | 45 | 45 | 45 | 45 | 0 | 0 | 0 | 0 |
| 0 | 0 | 0 | 45 | 45 | 45 | 45 | 45 | 45 | 45 | 45 | 45 | 45 | 45 | 45 | 0 | 0 | 0 | 0 | 0 |
| 0 | 0 | 0 | 0 | 45 | 45 | 45 | 45 | 45 | 45 | 45 | 45 | 45 | 45 | 0 | 0 | 0 | 0 | 0 | 0 |
| 0 | 0 | 0 | 0 | 0 | 0 | 45 | 45 | 45 | 45 | 45 | 45 | 45 | 0 | 0 | 0 | 0 | 0 | 0 | 0 |
| 0 | 0 | 0 | 0 | 0 | 0 | 0 | 0 | 0 | 45 | 0 | 0 | 0 | 0 | 0 | 0 | 0 | 0 | 0 | 0 |

**Figure 11.9** Physician-prescribed treatment requirements. (We use a minus sign to indicate the target pixels containing the tumor.)

When employing the full set (24 in our case) of beams there are no violations of the dose prescription, indicating that a feasible solution exists. As more beams are eliminated with our empirical criteria, errors begin to occur, indicating that the dose prescription cannot be satisfied with few beams (see Table 11.2). Figure 11.10 shows the dose distribution and the dose prescription violations obtained with the six-beam configuration of our empirical method; and Figure 11.11 shows the same information for the two-beam configuration. A trained dosimetrist, who did not use our model, chose a two-beam configuration of 60° and 120°—which are exactly the same angles derived by our empirical method for the two-beam configuration.

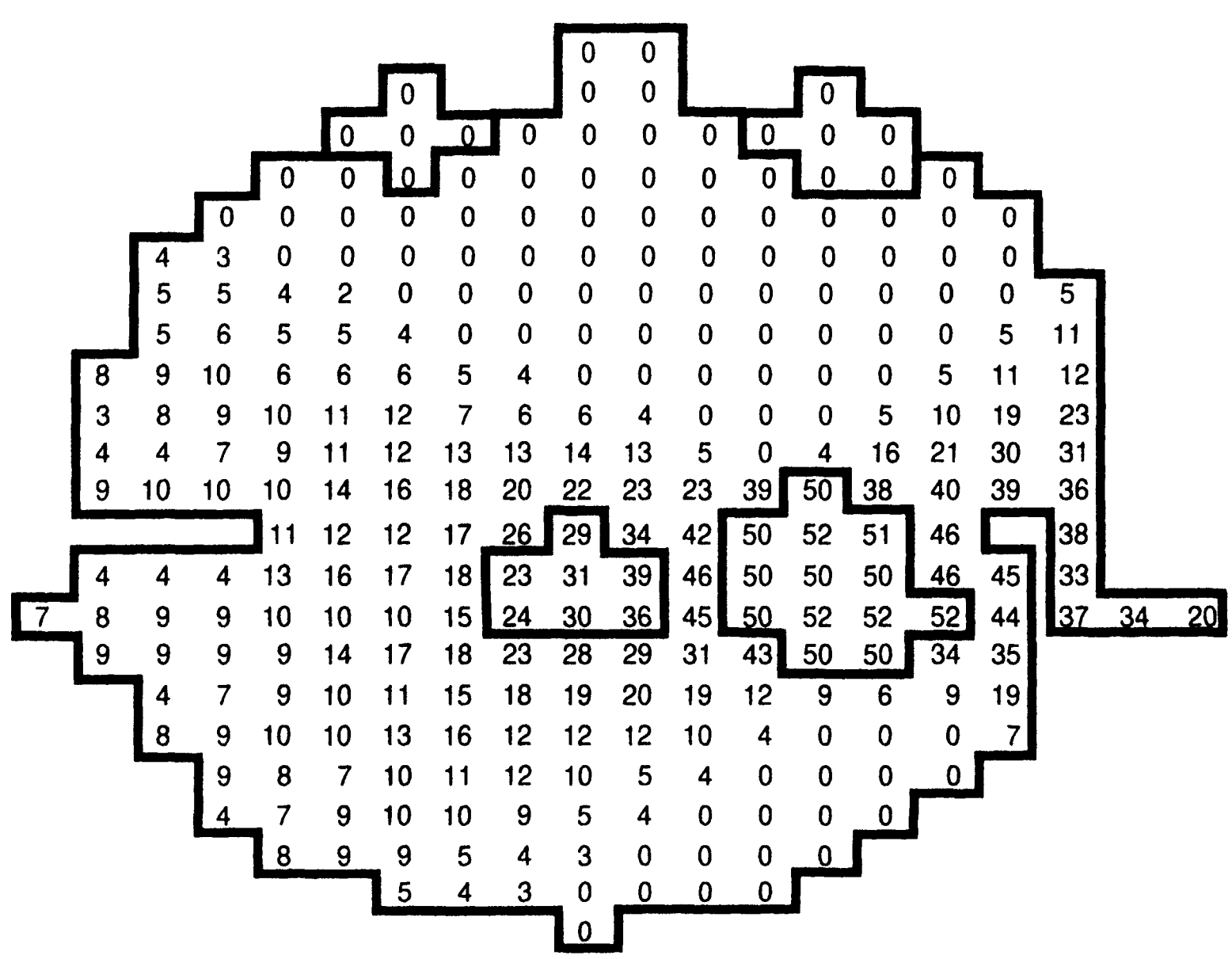

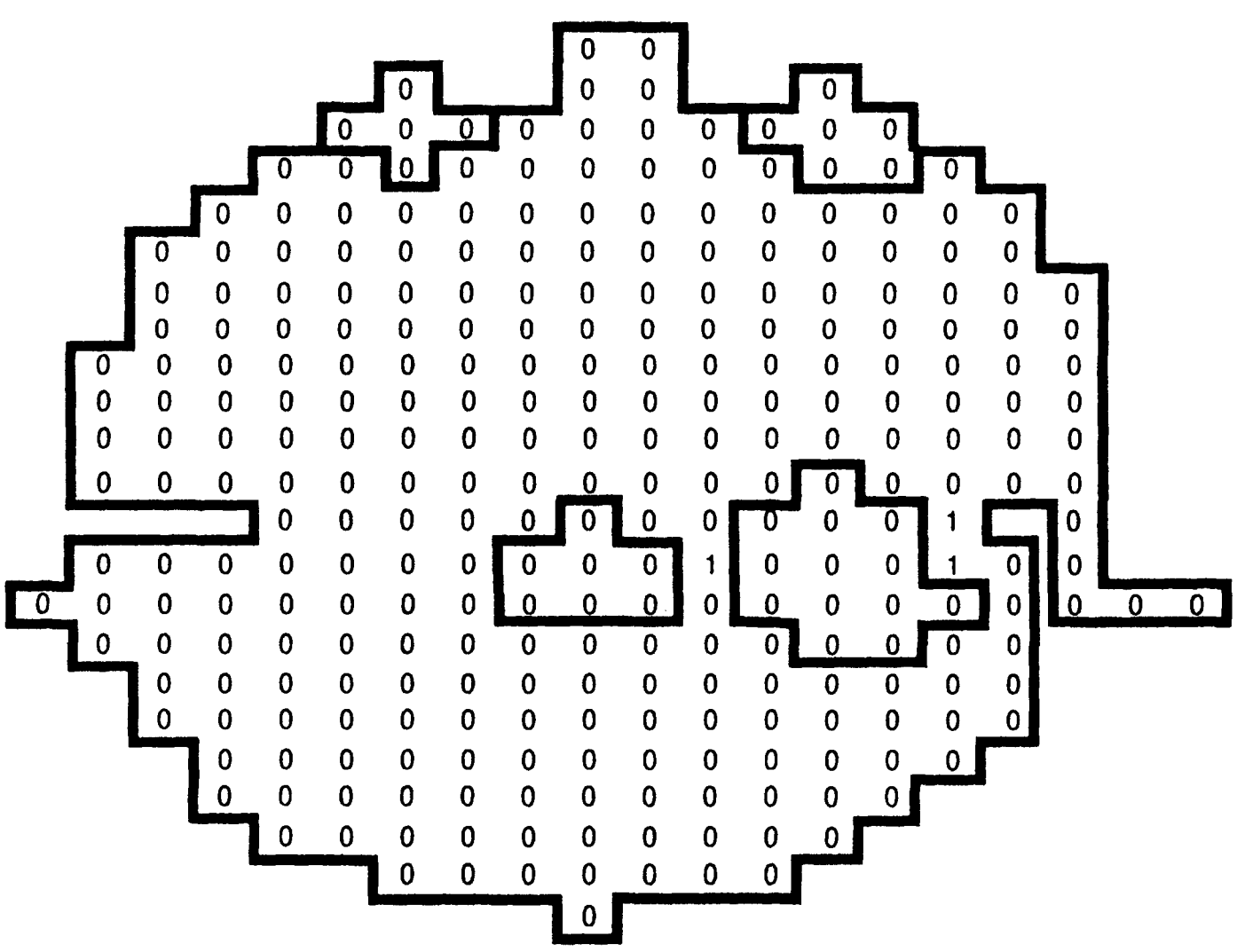

**Figure 11.10** Dose distribution (top figure) and dose prescription violations (bottom figure) for the six-beam configuration.

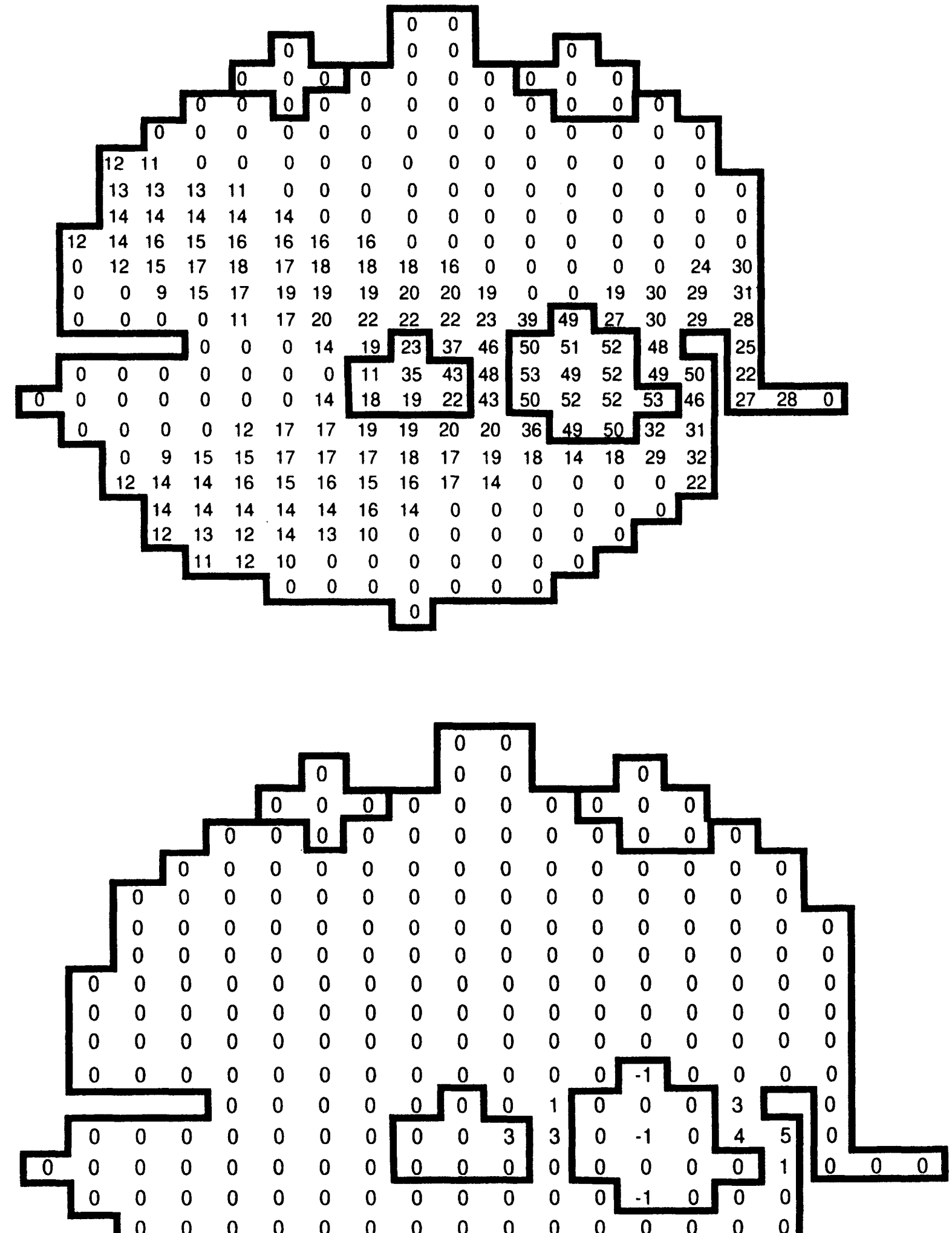

**Figure 11.11** Dose distribution (top figure) and dose prescription violations (bottom figure) for the two-beam configuration.

**Table 11.1** Repeated application of Cimmino's algorithm to choose beams. (The two leftmost columns in the table display the beam numbers and angles respectively.)

| | | Sum of ray intensities | | | | | | Absolute sum of adjacent ray intensity differences | | | | | |
|---|---|---|---|---|---|---|---|---|---|---|---|---|---|
| | | Number of beams | | | | | | Number of beams | | | | | |
| | | 24 | 11 | 6 | 4 | 3 | 2 | 24 | 11 | 6 | 4 | 3 | 2 |
| 1 | 0 | 10 | - | - | - | - | - | 4 | - | - | - | - | - |
| 2 | 15 | 126 | - | - | - | - | - | 28 | - | - | - | - | - |
| 3 | 30 | 165 | 133 | - | - | - | - | 22 | 24 | - | - | - | - |
| 4 | 45 | 179 | 147 | 143 | 137 | - | - | 34 | 28 | 28 | 34 | - | - |
| 5 | 60 | 184 | 151 | 147 | 146 | 145 | 150 | 28 | 26 | 24 | 26 | 32 | 40 |
| 6 | 75 | 184 | 155 | 151 | - | - | - | 38 | 34 | 36 | - | - | - |
| 7 | 90 | 184 | 150 | 152 | - | - | - | 46 | 36 | 38 | - | - | - |
| 8 | 105 | 187 | 157 | 150 | 152 | 147 | - | 34 | 30 | 28 | 32 | 34 | - |
| 9 | 120 | 186 | 154 | 148 | 143 | 149 | 149 | 30 | 26 | 28 | 24 | 30 | 34 |
| 10 | 135 | 153 | 112 | - | - | - | - | 36 | 28 | - | - | - | - |
| 11 | 150 | 65 | - | - | - | - | - | 32 | - | - | - | - | - |
| 12 | 165 | 145 | 104 | - | - | - | - | 38 | 32 | - | - | - | - |
| 13 | 180 | 48 | - | - | - | - | - | 28 | - | - | - | - | - |
| 14 | 195 | 133 | - | - | - | - | - | 24 | - | - | - | - | - |
| 15 | 210 | 148 | 124 | - | - | - | - | 24 | 26 | - | - | - | - |
| 16 | 225 | 147 | 125 | - | - | - | - | 28 | 24 | - | - | - | - |
| 17 | 240 | 139 | - | - | - | - | - | 22 | - | - | - | - | - |
| 18 | 255 | 138 | - | - | - | - | - | 28 | - | - | - | - | - |
| 19 | 270 | 134 | - | - | - | - | - | 32 | - | - | - | - | - |
| 20 | 285 | 134 | - | - | - | - | - | 22 | - | - | - | - | - |
| 21 | 300 | 129 | - | - | - | - | - | 20 | - | - | - | - | - |
| 22 | 315 | 126 | - | - | - | - | - | 22 | - | - | - | - | - |
| 23 | 330 | 21 | - | - | - | - | - | 12 | - | - | - | - | - |
| 24 | 345 | 109 | - | - | - | - | - | 28 | - | - | - | - | - |

**Table 11.2** Increasing error with decreasing number of beams—for the procedure described in text to eliminate beams from a treatment plan.

| Number of beams | Number of pixels with errors | Absolute sum of errors |
|---|---|---|
| 24 | 2 | 2 |
| 6 | 3 | 3 |
| 4 | 3 | 5 |
| 3 | 4 | 5 |
| 2 | 10 | 22 |

## 11.6 Combination of Plans in Radiotherapy

We conclude our chapter on radiation therapy treatment planning with a discussion of how to systematically combine several given treatment plans, which we term *plan combination* in RTTP. Plan combination gives rise to a problem that is neither of the forward nor of the inverse type discussed earlier, but one that addresses a situation when for a specific clinical case, a set of several treatment plans is proposed, for which both the forward and the inverse problems have already been solved with whatever available methods. This means that for each plan both a description of the dose distribution and a detailed scheme of machine setup parameters for dose delivery are available. Additionally, it is assumed that although all proposed plans approximate the desired and prescribed dose distribution, they all violate the prescribed dose in at least one significant region of the volume to be treated.

Such situations arise in clinically complicated cases, so that even an experienced dosimetrist will have difficulty in finding an acceptable treatment plan. Alternatively, it is conceivable that a set of such inaccurate treatment plans will be generated intentionally in some time- or effort-saving crude planning method. Both possibilities are especially plausible for truly three-dimensional planning, and present a dilemma to the radiotherapist. He must either choose a single plan from the set of plans, or attempt to create new plans, probably based on relaxed requirements. We will now discuss treatment plans as vectors in the Euclidean space, and define their *equivalence*, *acceptability* and *realizability*. A simple linear algebraic model for combining them is utilized in order to derive, from the given set of approximate plans, a *combined treatment plan* that will be both acceptable and technically realizable. In the event that such a combined plan does not exist, the alternatives for relaxing the treatment requirements will be systematically considered.

### 11.6.1 Basic definitions and mathematical modeling

Let $\Omega_j$, $j = 1, 2, \ldots, n$, be $n$ mutually disjoint subsets of a cross section $\Omega$ such that $\Omega = \cup_{j=1}^{n} \Omega_j$. These could be exactly identified with structures such as tumor regions and critical organs, but could also be (as will be understood later) any other regions in $\Omega$ for which the dose distributions in any of the proposed plans do not conform with the requirements.

**Definition 11.6.1 (Treatment Plan)** *A vector $x = (x_j)_{j=1}^{n} \in \mathbb{R}^n$ is called a treatment plan vector relative to a given cross section $\Omega$, and its partition $\{\Omega_j\}$, $j = 1, 2, \ldots, n$, if $x_j$ is nonnegative and it represents the total dose absorbed everywhere inside the region $\Omega_j$ for all $j = 1, 2, \ldots, n$.*

**Definition 11.6.2 (Radio-equivalence)** *Two radiation therapy treatment plans that give rise to two treatment plan vectors $x^1, x^2$ are called radio-equivalent if $x^1 = x^2$ (component wise) regardless of the actual treatment* *machine setups and the types of machines that yield them.*

Let $l = (l_j)$ and $u = (u_j)$ be two given treatment plan vectors in which $l_j$ and $u_j$ are respectively, the lower and upper bounds on the dose at region $\Omega_j$. These vectors are assumed to be prescribed by the radiotherapist.

**Definition 11.6.3 (Acceptable Treatment Plan)** *Given $\Omega$, subsets $\Omega_j$, $j = 1, 2, \ldots, n$, $l$ and $u$, a treatment plan vector $x$ is called acceptable if regardless of how $x$ is delivered, $l_j \leq x_j \leq u_j$, $j = 1, 2, \ldots, n$.*

**Definition 11.6.4 (Realizable Treatment Plan)** *A treatment plan vector $x$ is called realizable if there exists (in-house) a clinical machine setup that can deliver it.*

Realizability is a subjective and user-dependent concept that may vary according to equipment availability and treatment planning capabilities. The zero vector $0 \in \mathbb{R}^n$ is always realizable, but not acceptable unless $l = 0$. It is henceforth assumed that by utilization of whatever planning procedure (or procedures) and by using whatever equipment, a finite set of treatment plan vectors $\{x^1, x^2, \ldots, x^m\}$ has been generated for a given clinical case such that every $x^i = (x_1^i, x_2^i, \ldots, x_n^i)$ is realizable in the given clinical environment, but none of the $x^i$ is acceptable. This means that for each $i$, $1 \leq i \leq m$ there exists at least one $j$, $1 \leq j \leq n$, such that either $x_j^i < l_j$ or $x_j^i > u_j$.

A nonnegative combination of treatment plan vectors is a treatment plan vector $y = \sum_{i=1}^m a_i x^i$ such that the linear combination coefficients $a_1, a_2, \ldots, a_m$, are nonnegative real numbers. From the nature of RTTP, a nonnegative combination of realizable vectors is also realizable.

The treatment plan vectors $l$ and $u$, prescribed for a given clinical case, determine a box $B$ (henceforth called *treatment box*) , in $\mathbb{R}^n$ defined as

$$B = \{x \in \mathbb{R}^n \mid l_j \leq x_j \leq u_j, \; j = 1, 2, \ldots, n\}. \tag{11.19}$$

Obviously, every treatment plan vector that belongs to $B$ is acceptable.

The set $\{x^1, x^2, \ldots, x^m\}$ of realizable but unacceptable treatment plan vectors, together with the always realizable zero vector generate a polyhedral convex cone $C$ with apex at the origin, given by the convex hull of all nonnegative linear combinations of the vectors $\{x^i\}_{i=1}^m$, i.e.,

$$C = \mathrm{conv}\{y \in \mathbb{R}^n \mid y = \sum_{i=1}^m a_i x^i, \; a_i \geq 0, \; i = 1, 2, \ldots, m\}. \tag{11.20}$$

Following these definitions, the linear algebraic model for plan combination in RTTP can be formulated as follows:

> Given $\Omega$, $\Omega_j$, $j = 1, 2, \ldots, n$, $l$ and $u$, and given a set of realizable treatment plan vectors denoted by $\{x^1, x^2, \ldots, x^m\}$, find a treatment plan vector $y$ such that $y \in B \cap C$.

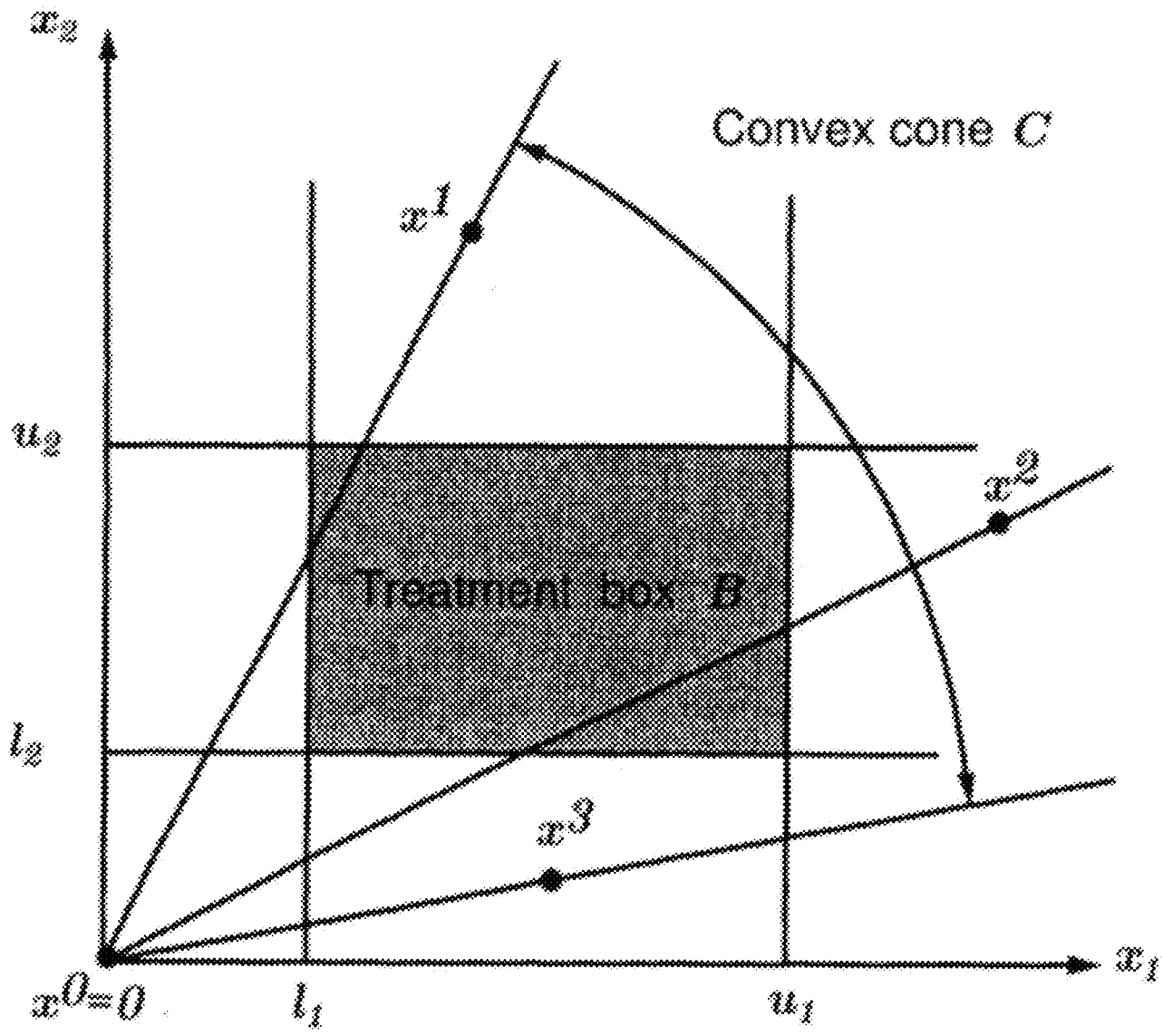

**Figure 11.12** A plan combination problem with $n = 2$ regions and a set $\{x^1, x^2, x^3\}$ of realizable but unacceptable treatment plans. In this drawing, $B \cap C \neq \emptyset$ and the problem is feasible.

Since the plan combination problem is a linear feasibility problem (see Chapter 5), we call any treatment plan vector $y$ that solves it a *feasible solution.* Figure 11.12 depicts a typical feasible plan combination problem with $m = 3$ and $n = 2$. The overall methodology considered in this section is summarized in the flow chart presented in Figure 11.13 whose meaning becomes clear after reading the next subsections.

### 11.6.2 The feasible case

The formulation of the plan combination problem makes a thorough analysis possible and suggests several options for radiation therapy treatment. The problem poses first the question of feasibility, i.e., does at least one acceptable linear combination of realizable plans exist? This question can be answered using linear programming methods.

First, we give a dual formulation of the plan combination problem. Let us organize the given set $\{x^1, x^2, \ldots, x^m\}$ of realizable but probably unacceptable treatment plan vectors in an $n \times m$ matrix $X$ whose $i$th column is $x^i$. The transpose of the $j$th row of $X$ is denoted by $r^j$, i.e., $r^j = (x_j^1, x_j^2, \ldots, x_j^m)^{\mathrm{T}}$. The dual formulation of the plan combination problem of finding a vector $y \in B \cap C$ is the following linear feasibility problem:

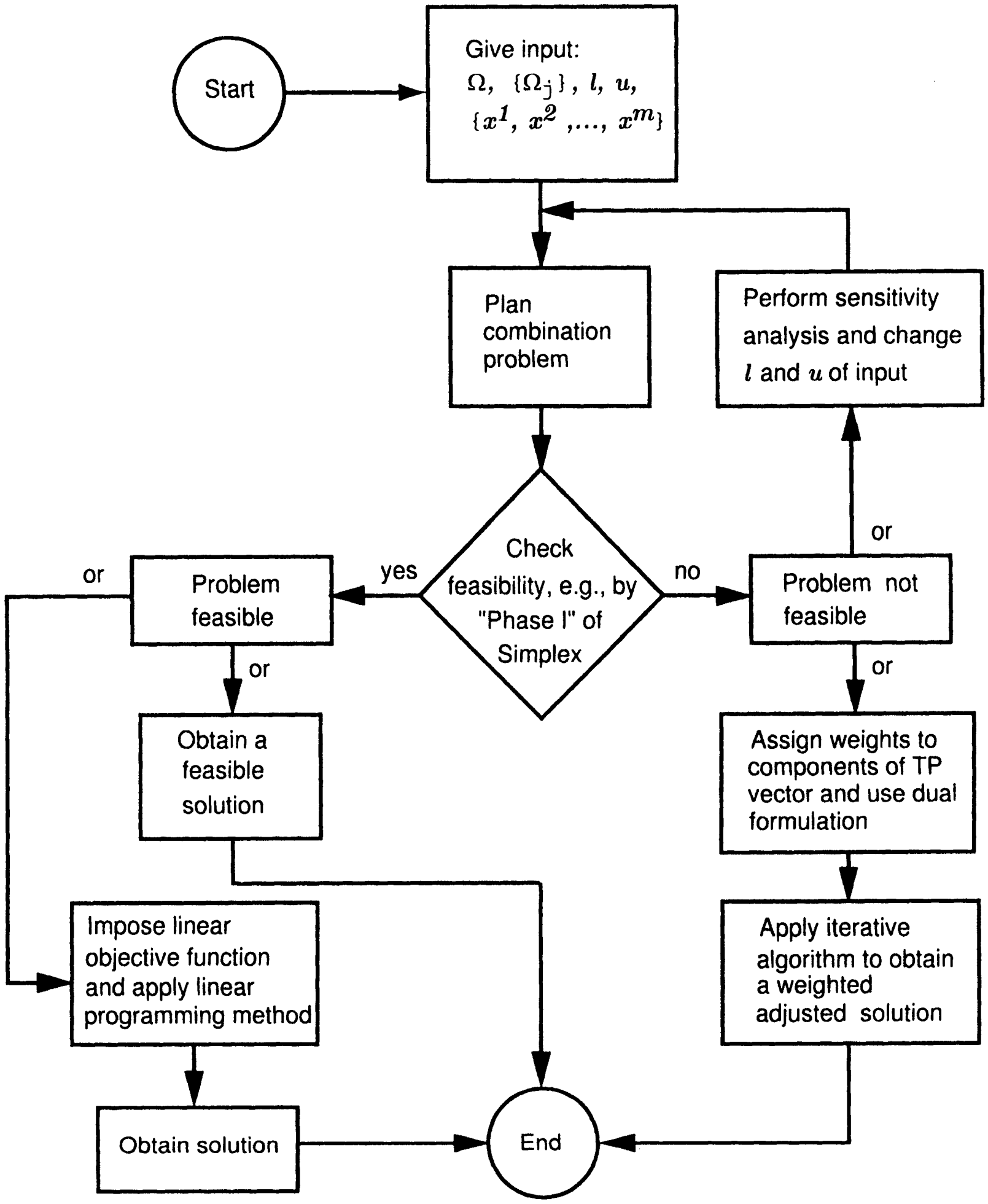

**Figure 11.13** Methodological flow chart of plan combination.

Given $l, u$ and $r^j$, $j = 1, 2, \ldots, n$, find a vector $a = (a_i) \in \mathbb{R}^m$ such that

$$\begin{aligned} &l_j \leq \langle a, r^j \rangle \leq u_j, \; j = 1, 2, \ldots, n, \\ &a_i \geq 0, \qquad\qquad i = 1, 2, \ldots, m. \end{aligned} \tag{11.21}$$

The duality between the plan combination problem and (11.21) is obvious, and it is also clear that the plan combination problem is feasible if and only if (11.21) is.

The simplex method of linear programming is a practical tool for solving (11.21). Phase 1 of the simplex method is concerned in particular with finding a feasible solution to a linear programming problem and it answers the question whether (11.21) is feasible or not. For the remainder of this section we assume feasibility, i.e., that $B \cap C \neq \emptyset$.

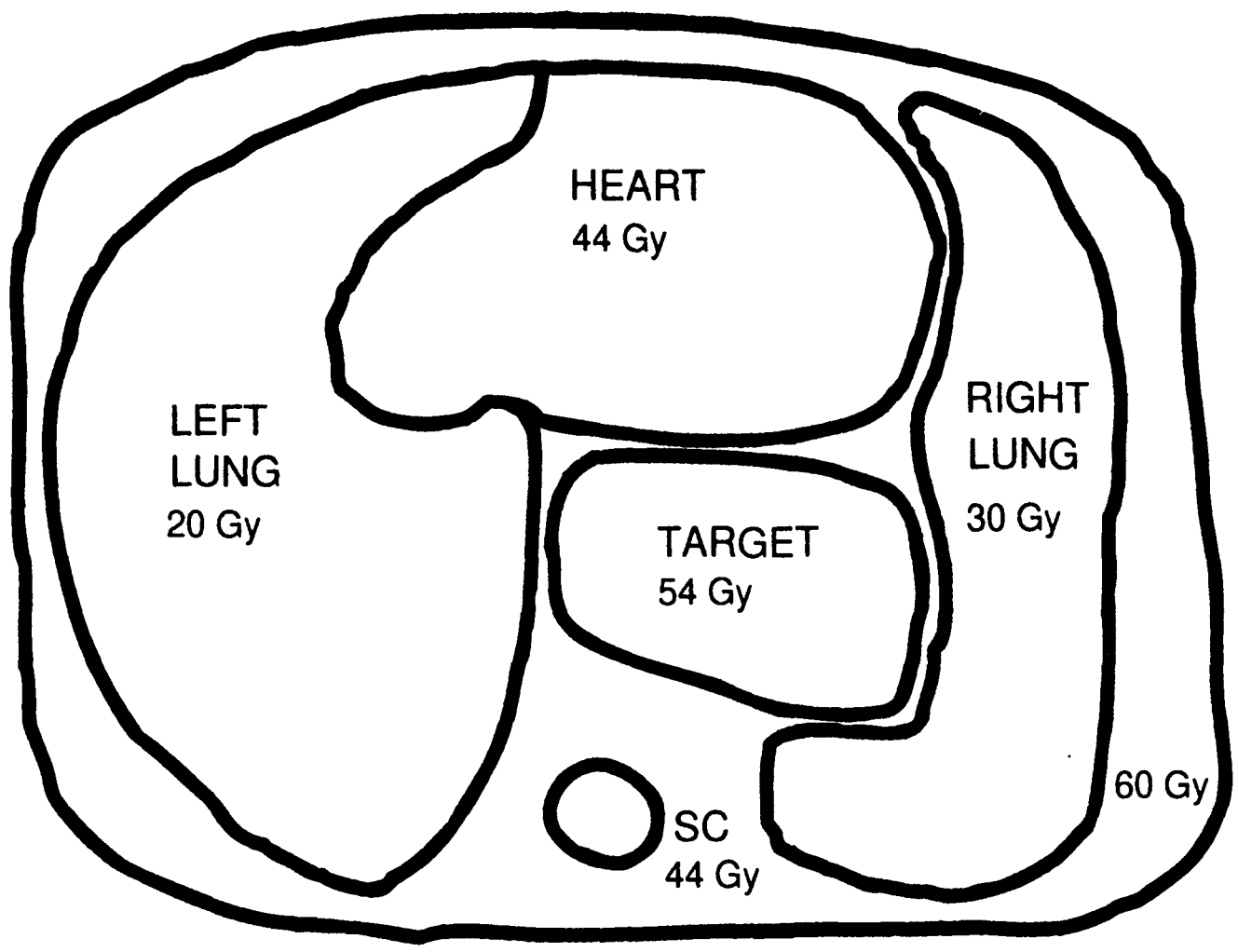

**Figure 11.14** A patient cross-section contours and dose constraints.

Therefore, in this case more information can be extracted, and probably several feasible treatment plan vectors can be suggested. The physician may be interested in a treatment plan that is not only feasible, but is also optimal in some sense. For example, in addition to feasibility we may wish to minimize the total radiation dose $f(a)$ given to the patient, where $f(a) \doteq \sum_{j=1}^{n} \sum_{i=1}^{m} a_i x_j^i$; or alternatively, to maximize the total dose to the target $g(a)$, where $g(a) \doteq \sum_{j \in S} \sum_{i=1}^{m} a_i x_j^i$, and $S$ is a subset of $\{1, 2, \ldots, n\}$ that consists of all the indices corresponding to the target. Another criterion for selecting one realizable and acceptable treatment plan vector over another may be that a combined treatment plan that uses fewer beams, i.e., $a_i = 0$, for some or several indices $i \in \{1, 2, \ldots, m\}$ is preferable, perhaps for clinical reasons.

These are all linear programming problems and the simplex method or an interior point algorithm (see Chapter 8) can be applied. Since in our case the polyhedral set described by (11.21) is bounded, it is impossible to have an unbounded solution. The following example illustrates this approach.

**Example 11.6.1 (Plan combination)** *Suppose that the target for the radiation treatment is in the chest, surrounded by the two lungs, the heart and the spinal cord, as illustrated in Figure 11.14. Let the vectors of lower and upper bounds be* $l = (0, 0, 0, 0, 0, 54)^{\mathrm{T}}$ *and* $u = (20, 30, 44, 44, 60, 100)^{\mathrm{T}}$, *respectively. The entries in each vector are doses measured in Gy units and correspond to the organs as follows: the first entry corresponds to the right lung, the second entry to the left lung, the third to the heart, the fourth to the spinal cord, the fifth to the complementary tissue and the sixth entry corresponds to the target. Let $x^1, x^2, x^3$ be three given vectors of realizable but unacceptable treatment plans: $x^1 = (15, 25, 40, 43, 45, 50)^{\mathrm{T}}$, $x^2 = (30, 40, 50, 40, 60, 60)^{\mathrm{T}}$, $x^3 = (40, 40, 60, 55, 60, 55)^{\mathrm{T}}$. Let $a_1, a_2$, and $a_3$ be the corresponding coefficients that we seek for the linear combination. Then the linear feasibility problem for this example is to find nonnegative numbers $a_1, a_2$, and $a_3$ such that*

$$\begin{aligned}
0 &\le 15a_1 + 30a_2 + 40a_3 \le 20, \\
0 &\le 25a_1 + 40a_2 + 40a_3 \le 30, \\
0 &\le 40a_1 + 50a_2 + 60a_3 \le 44, \\
0 &\le 43a_1 + 40a_2 + 55a_3 \le 44, \\
0 &\le 45a_1 + 60a_2 + 60a_3 \le 60, \\
54 &\le 50a_1 + 60a_2 + 55a_3 \le 100.
\end{aligned} \tag{11.22}$$

*We impose an objective function $f(a)$ to minimize the total radiation given to the organs surrounding the target tissue. By way of example we take $f(a) = 168a_1 + 220a_2 + 255a_3$, which has to be minimized subject to $a \ge 0$ and to all inequalities in (11.22). A solution with a linear programming algorithm yields the following coefficients (all numbers are rounded): $a_1 = 0.83$, $a_2 = 0.21$, $a_3 = 0.0$. The treatment vector corresponding to the amount of radiation each organ receives, if this combined plan is used, is $y = (18.72, 29.10, 43.62, 44, 49.86, 54)^{\mathrm{T}}$. The objective function value is 185.31 Gy. Moreover, elementary sensitivity analysis shows that if we decrease the lower bounds of the radiation given to the target by one unit (i.e., from 54 Gy to 53 Gy), then the total radiation to all but the target organ will decrease by 4.72 Gy.*

*A different objective function, which might at times be more desirable, is one that maximizes the total radiation $g(a) = 50a_1 + 60a_2 + 55a_3$ given to the target tissue, subject to the constraints. The solution obtained for this problem is: $a_1 = 0.8$, $a_2 = 0.24$, $a_3 = 0$, with total radiation to the target 54.4 Gy, and with a combined treatment plan vector $y = (19.2, 29.6, 44, 44, 50.4, 54.5)^{\mathrm{T}}$.*

### 11.6.3 The infeasible case

When no feasible solution exists (see, for example, Figure 11.15), it is not necessarily a hopeless situation, and there are several possible alternatives. The analysis of Phase 1 of the simplex method provides some insight to the infeasibility issue, as it yields important information about the model, such as which constraints are violated if the current (infeasible) solution is used, and which constraints are tight for this solution, namely constraints that are met with equality. Using this information, the physician may decide to reformulate the dose requirements $l$ and $u$ only in so far as necessary to arrive at a feasible combined treatment plan. Such a possibility is illustrated in the next example.

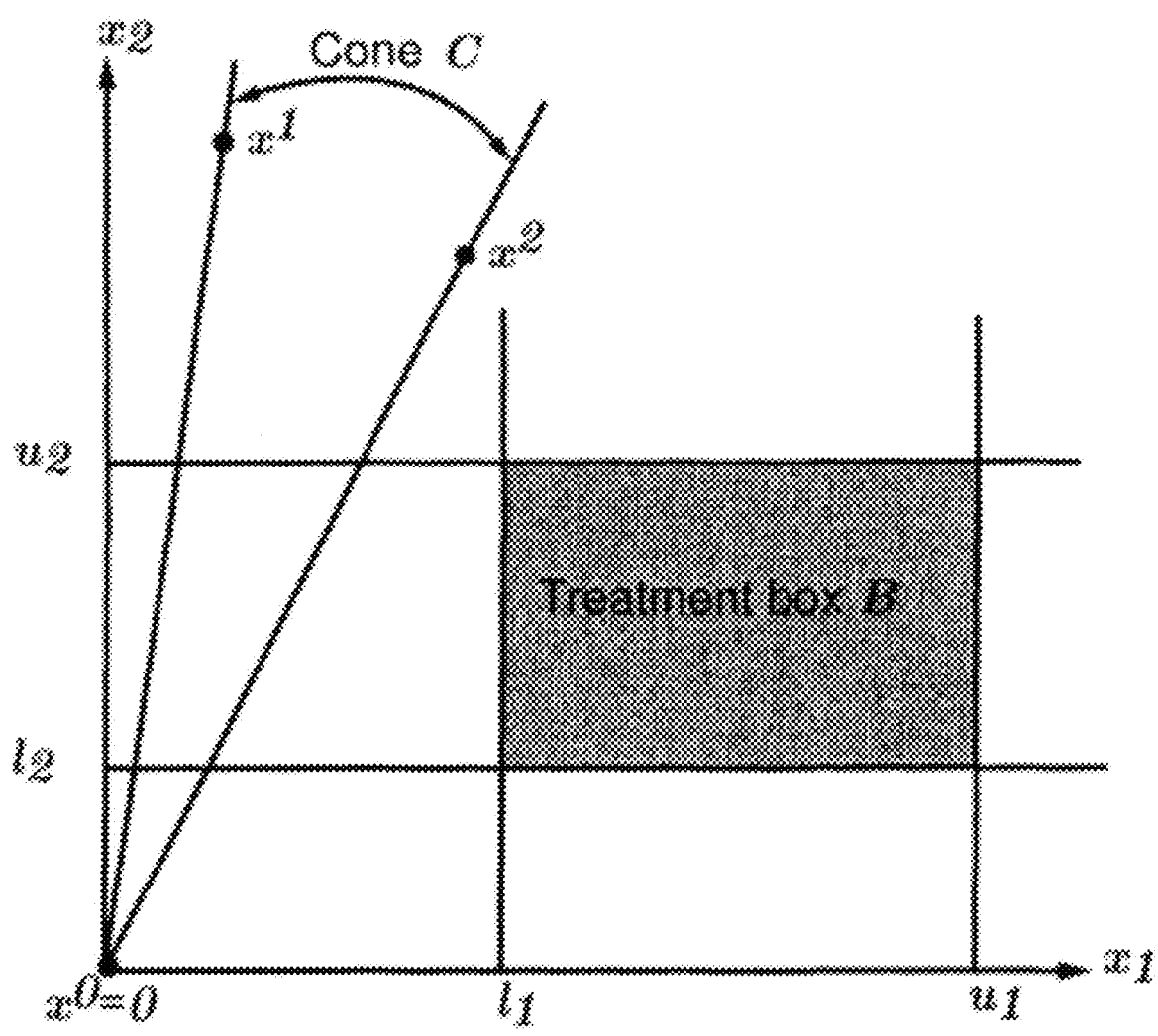

**Figure 11.15** $B \cap C = \emptyset$, an infeasible plan combination problem.

**Example 11.6.2 (An infeasible case)** *Suppose that for the same patient as in Example 11.6.1 we have the following realizable but unacceptable plans* $x^1 = (15, 25, 40, 45, 45, 55)^T$, $x^2 = (30, 40, 50, 50, 60, 50)^T$, *and* $x^3 = (40, 40, 60, 55, 60, 40)^T$. *Application of Phase 1 of the simplex method yields an infeasible solution with* $a_1 = 0.9777$, $a_2 = 0.0$, $a_3 = 0.0$, *giving rise to the treatment vector:* $y = (14.67, 24.44, 39.11, 44, 44, 53.78)^T$. *The violated constraint is the sixth one, corresponding to the target region. It is violated (underdosed) by 0.2222 Gy. In such a situation, the physician may decide that a relatively small decrease in the lower bound of the radiation requirement for the target can—and should—be done. Thus to obtain a feasible solution, the lower bound should be reduced from 54 Gy to, say, 53.7 Gy. An alternative change might be an increase in the upper bound of the constraint that contributes to the infeasibility, namely the fourth constraint. A simple computation indicates that the problem becomes feasible if the upper bound (44 Gy) is replaced by 44.2 Gy, which results in the combined treatment vector* $y = (14.73, 24.55, 39.27, 44.18, 54)^T$.

There is an alternative approach for the infeasible case. Rather than relaxing the vectors $l$ and/or $u$ as dictated by the sensitivity analysis of the simplex method, the physician is called upon to provide additional input. He is required to point out the relative importance (priorities) of delivering the prescribed doses to the various regions. An iterative algorithm is then applied to the dual formulation of the plan combination problem (11.21) and the result is called an *adjusted solution* to the original infeasible plan combination problem.

Various iterative algorithms are applicable to (11.21), many of which can be derived as special cases of the block-iterative projections (BIP) method (Algorithm 5.6.1). Some, such as ART3 (Algorithm 5.10.2) or ARM (Algorithm 5.10.1) are specifically tailored to handle an interval linear feasibility problem. However since many of these algorithms rely on feasibility of (11.21) for their convergence, we construct here an ad hoc fully simultaneous version of ART3 and demonstrate its capability to (asymptotically) generate adjusted solutions to the original infeasible plan combination problem.

We use the notation of Section 5.10. For each $j$, an interval linear inequality of (11.21) represents in $\mathbb{R}^m$ a hyperslab, whose bounding hyperplanes are $L_j \doteq \{a \in \mathbb{R}^m \mid \langle a, r^j \rangle = l_j\}$ and $U_j \doteq \{a \in \mathbb{R}^m \mid \langle a, r^j \rangle = u_j\}$. The *median hyperplane* is, for each $j$, $H_j \doteq \{a \in \mathbb{R}^m \mid \langle a, r^j \rangle = \frac{1}{2}(l_j + u_j)\}$, and the *half-width* of the $j$th hyperslab is $w_j = (u_j - l_j)/2 \parallel r^j \parallel$. The signed distance of a point $z \in \mathbb{R}^m$ from $H_j$ is given by

$$\hat{d}(z, H_j) = \frac{\langle z, r^j \rangle - \frac{1}{2}(l_j + u_j)}{\parallel r^j \parallel}. \tag{11.23}$$

Let $\{a^\nu\}_{\nu=0}^{\infty}$ be a sequence of iterates. We denote the distances of an element of the sequence from the median hyperplane by

$$\hat{d}_j \doteq d(a^\nu, H_j). \tag{11.24}$$

The sign +1 or −1 of $\hat{d}_j$ is denoted by $t_j \doteq \text{sign}\,(\hat{d}_j)$. With these definitions a simultaneous ART3 algorithm that aims at solving the feasibility problem (11.21) can be constructed. Assume that a set of real numbers $\pi_j$ for all $j = 1, 2, \ldots, n$ is given, such that $0 < \pi_j < 1$ for all $j$, and $\sum_{j=1}^{n} \pi_j = 1$. By assigning such a set of priorities $\pi_j$ the therapist indicates the relative importance attributed to each face of the treatment box $B$, defined by (11.19), thus to each interval of the system (11.21). Analogous to a typical step of ART3 we define, for a given iterate $a^\nu \in \mathbb{R}^m$ and with respect to each hyperslab $Q_j = \{a \in \mathbb{R}^m \mid l_j \leq \langle a, r^j \rangle \leq u_j\}$, the step sizes

$$s_{\nu,j} = \begin{cases} 0, & \text{if } \mid \hat{d}_j \mid \leq w_j, \\ \alpha_\nu(\hat{d}_j - t_j w_j), & \text{if } w_j < \mid \hat{d}_j \mid < 2w_j, \\ \hat{d}_j, & \text{if } 2w_j \leq \mid \hat{d}_j \mid, \end{cases} \tag{11.25}$$

where $\hat{d}_j$ and $t_j$ are determined from the formulae given above. Factors that determinine the value of the constants $\alpha_\nu$ are explained later.

**Algorithm 11.6.1 Simultaneous ART3**

**Step 0:** (Initialization.) $a^0 \in \mathbb{R}^m$ is arbitrary.

**Step 1:** (Iterative step.) Given the current iterate $a^\nu$, calculate for each index $j$, $1 \le j \le n$,

$$a^{\nu+1,j} = a^\nu - s_{\nu,j} \frac{r^j}{\| r^j \|}, \tag{11.26}$$

then calculate

$$a^{\nu+1} = \sum_{j=1}^{n} \pi_j a^{\nu+1,j}. \tag{11.27}$$

It is straightforward to verify that each intermediate iterate $a^{\nu+1,j}$ is a relaxed orthogonal projection of $a^\nu$ onto $Q_j$. Thus if $\lambda_\nu$ denotes relaxation parameters, then equation (11.26) can be rewritten as

$$a^{\nu+1,j} = a^\nu + \lambda_\nu \left( P_{Q_j}(a^\nu) - a^\nu \right), \tag{11.28}$$

where,

$$\lambda_\nu \doteq \begin{cases} 1, & \text{if } \mid \hat{d}_j \mid \le w_j, \\ \alpha_\nu, & \text{if } w_j < \mid \hat{d}_j \mid < 2w_j, \\ 1, & \text{if } 2w_j \le \mid \hat{d}_j \mid . \end{cases} \tag{11.29}$$

For $\alpha_\nu = \alpha = 2$ the iterative step (11.26) coincides with that of the original ART3 algorithm. However, to ensure applicability of the theory to the present simultaneous ART3 algorithm, the relaxation parameters $\lambda_\nu$ in (11.28) must fulfill the condition that for all $\nu \ge 0$,

$$\epsilon_1 \le \lambda_\nu \le 2 - \epsilon_2, \tag{11.30}$$

with any fixed $\epsilon_1, \epsilon_2 > 0$.

If this condition is adhered to, then the simultaneous ART3 algorithm becomes a special case of the nonlinear fully simultaneous Cimmino algorithm, i.e., the BIP algorithm (Algorithm 5.6.1) with all indices of all constraints lumped into a single block. Then the theory guarantees convergence of the simultaneous ART3 algorithm. The algorithm converges globally to a solution of (11.21) if (11.21) is feasible. If (11.21) is infeasible, then the iterates $\{a^\nu\}$ generated by the simultaneous ART3 algorithm can be shown to converge locally to a point that minimizes the convex combination with $\{\pi_j\}$ as coefficients of the squares of the distances to the hyperslabs $Q_j$, $j = 1, 2, \ldots, n$.

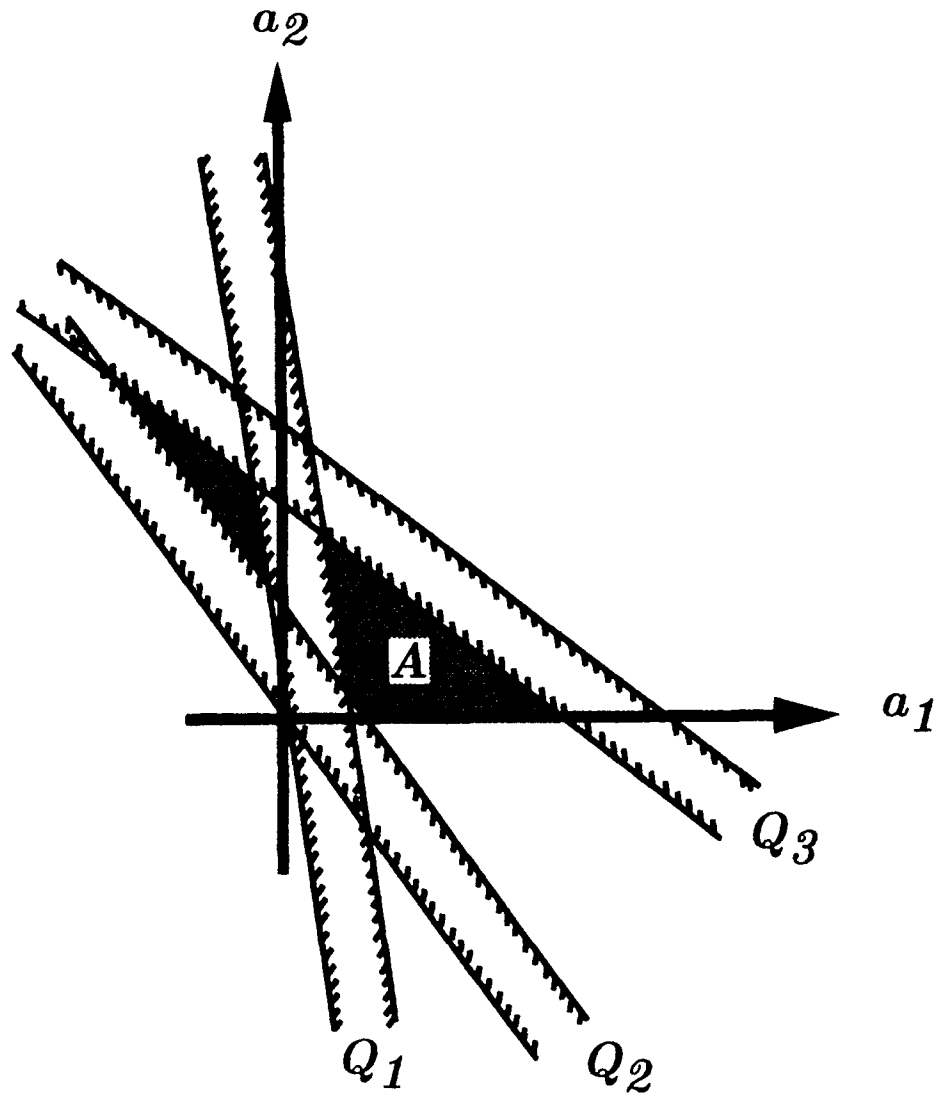

**Figure 11.16** Dual space description of an infeasible plan combination problem.

In this sense we call such a limit point an *adjusted solution* to the infeasible original plan combination problem. It violates some or all constraints which define the treatment box $B$, but at the same time it is a weighted least squares solution to the linear feasibility problem – the weights being the prescribed priorities $\{\pi_j\}_{j=1}^n$. Nonnegativity of the iterates $\{a^\nu\}$ generated by the simultaneous ART3 algorithm is not guaranteed even if the initial point $a^0 \in \mathbb{R}^m$ is chosen so that $a^0 \geq 0$. To overcome this difficulty we perform orthogonal projections onto the nonnegative orthant between every two consecutive iterations of the simultaneous ART3 algorithm. This is done by replacing $a^{\nu+1}$ in (11.27) by $\hat{a}^{\nu+1}$ and then adding to the iterative step the operation

$$a_i^{\nu+1} = \begin{cases} \hat{a}_i^{\nu+1}, & \text{if } \hat{a}_i^{\nu+1} \geq 0, \\ 0, & \text{if } \hat{a}_i^{\nu+1} < 0. \end{cases} \tag{11.31}$$

Algorithm 11.6.1 with (11.30)–(11.31) converges to a nonnegative weighted least squares solution of the system (11.21), i.e., to a point in the region marked $A$ in Figure 11.16.

**Example 11.6.3** *Using the data in Example 11.6.2, we run a computer program that executes the simultaneous ART3 algorithm. We stopped after $\mu$ iterations, where $\mu$ is the smallest integer such that the inequality* $\| a^{\mu-1} - a^{\mu} \| / \| a^{\mu} \| < 10^{-5}$ *holds. If equal weights are given to all organs, then the solution is* $a = (0.98, 0.0013, 0.001)^{\mathrm{T}}$, *and the resulting plan is* $y = (14.75, 24.54, 39.25, 44.13, 44.15, 53.9)^{\mathrm{T}}$. *We see that the problematic constraints, namely the fourth (which is associated with the spinal cord) and the sixth (the target) are violated. Nevertheless, y is a least-squares point as described above. If different weights are given, the algorithm converges to a weighted least-squares point. For example, if the weights given in this example are* (0.05, 0.05, 0.05, 0.05, 0.05, 0.75) *then the corresponding treatment plan is* $y = (14.8, 24.62, 39.35, 44.23, 44.27, 53.99)^{\mathrm{T}}$, *resulting from the solution* $a = (0.98, 0.002, 0.0018)^{\mathrm{T}}$.

It should be noted that the two alternatives suggested for the infeasible case, may be mathematically and practically equivalent for some cases. The same feasible solution obtained by moving an upper or lower bound, as described in the first option for the infeasible case, can be obtained using the simultaneous ART3 with appropriate weights.

## 11.7 Notes and References

For a treatise on the physics of radiology see, e.g., Johns and Cunningham (1974) or De Vita, Hellman, and Rosenberg (1989). Some condensed details about radiation quantities related to this chapter can be found in Barrett and Swindell (1981, Appendix D). See also the computationally oriented book of Wood (1981), or a conference proceedings such as Bruinvis et al. (1987). An up-to-date review of the field (including 130 references) is found in Brahme (1995b). This work contains also a concise description of the biological objective function optimization approach, recently pursued with promising results by Brahme and his coworkers. Much can be learned about recent progress from a special issue of the *International Journal of Imaging Systems and Technology*, Brahme (1995a).

**11.1** We presented the continuous model and formulated the continuous forward and inverse problems in RTTP in general terms only. Work in this direction (as opposed to the fully discretized approach that we have adopted here) originates in the paper by Brahme, Roos and Lax (1982), followed by the papers of Cormack (1987), Cormack and Cormack (1987), and Cormack and Quinto (1989, 1990). The forward problem, represented by (11.1), i.e., the process commonly called dose calculations, has received great attention for many years now and its complexity brought about many methods and formulae. Several relevant references are: Geijn (1972), Cunningham, Shrivastava, and Wilkinson (1972), Sontag and Cunningham (1978), Lane, Bloch and Davis (1974), Cunningham (1972), and the review by Goitein ([illegible]). Dose calculations continue to form a substantial research issue, see, e.g., Altschuler, Sontag and Bloch (1987), and Näfstadius, Brahme and Nordell (1984). The quotation of Dr. M. Goitein is taken from a report by Hindus (1988); see also Barrett and Swindell (1981, Appendix D.5.3). The feasibility formulation (11.4) was proposed by Altschuler and Censor (1984).

**11.2** Mathematical methods (mostly optimization theory techniques) are used in radiation therapy in ways that differ conceptually from the approach presented here. These other approaches can be generally identified as follows.

*Comparison among rival plans.* Here several treatment plans are compared, based on their score with respect to some predetermined quality index. The treatment plans are all fixed, and the selection of a plan depends largely on the definition of the quality index. Various quality-index functions were proposed and advocated on different grounds. These include, by way of example, the complication probability factor (CPF) of Dritschilo et al. (1978) and Wolbarst et al. (1980), the normal tissue dose (NTD), the optimum target dose (OTD), and the critical organ dose (COD), in Kartha et al. (1982); see also Wolbarst, Chin, and Svensson (1982).

*Optimization with respect to a few setup parameters.* Here the compared treatment plans are not fixed, but most setup parameters are kept unchanged while a few (one, two, or more) are allowed to take several possible values which will result in different treatment plans. A criterion function is set up, and a mathematical programming (i.e., optimization) problem is solved to yield an optimal choice of the unfixed variables. Usually linear programming is used. Although such an approach is logically also a comparison among rival plans, it is apparent that because of the greater freedom in the selection process (i.e., the process of mathematical optimization) many more options are compared. Into this category fall the work of Gallagher (1967), Bahr et al. (1968), Redpath, Vickery and Wright (1976), McDonald and Rubin (1977), and Starkschall (1984). The implementation of some of these methods, such as in McDonald and Rubin (1977) and in Starkschall (1984) constrains the dose distribution at only a small number of points (typically less than 20) chosen by the dosimetrist.

Additional work with optimization techniques includes Legrass et al. (1986), Fymat et al. (1988), the method of Webb (1989) with simulated annealing, Swan (1981), Burkard et al. (1994), and many others. For general reviews of the status of involvement of mathematical modeling see, e.g., Aird (1989) and Goitein (1990).

Another attempt at approaching the problem mathematically is given in the papers of Ebert (1977a, 1977b). His models are different from those studied here. The S-model assumes a source of radiation that moves around the patient section in orbits, while his F-model allows some parameters to vary as others remain fixed and therefore falls into the second class of methods described above.

Among more recent advances we refer the reader to the works of Brahme, Lind, and Källman (1990), Lind (1990), and Lind and Källman (1990). Their approach is to apply various iterative techniques to a discretized version of a Fredholm integral equation of the first kind. See also Söderström (1995). Raphael (1991, 1992a, 1992b) in his dissertation and subsequent publications formulated the RTTP problem as constrained optimization in appropriate $\mathcal{L}^2$ function space and applied to it iterative algorithms. His work and the work of Brahme and co-workers mentioned above, along with the full discretization approach presented here, show in more than one way the usefulness of iterative techniques for the solution of the inverse problem in RTTP.

The approach that we presented in this chapter differs fundamentally from all these approaches. We do not compare fixed rival plans. We do not select a plan by allowing one or few parameters to change sequentially. We rather address the inverse problem within a general framework where both the patient section and the sources are fully discretized and the full discretization is the tool for handling the modeling difficulties inherent in the continuous model.

**11.3–11.5** The presentation in these sections, as well as in the previous one, is based on works of Altschuler, Censor, and Powlis (1988, 1992), Censor, Altschuler, and Powlis (1988), Censor, Powlis, and Altschuler (1987), and Powlis et al. (1989). A review appears also in Altschuler, Censor, and Powlis (1992). A description of a new generation of radiation therapy machines is given by Brahme (1987). Experiments with a similar approach are reported in Altschuler, Powlis, and Censor (1985) and Powlis, Altschuler, and Censor (1985) where the method is applied to uniform beams which are not further discretized into rays.

**11.6** The material of this section comes from Censor et al. (1988) and Censor and Schwarz (1989). It should be noted that the plan combination approach is closely related to the full discretization philosophy given earlier. The inverse problem of RTTP, in its fully discretized formulation is actually a high-resolution plan combination problem.

CHAPTER 12

# Multicommodity Network Flow Problems

> However, the production and distribution of commodities is entirely unorganized so that everybody must live in fear of being eliminated from the economic cycle, in this way suffering from the want of everything.
>
> I trust that posterity will read these statements with a feeling of proud and justified superiority.
>
> Albert Einstein, *Message in the Time-Capsule*, 1939

The movement of some commodity between two points on a network can be represented using *network flow models.* Such models are prevalent in the real world; and the transportation of products over a railroad system is the prototypical example of a network flow problem. Other examples include the transportation of messages over a telecommunications network, the transmission of electricity over a power distribution grid, the scheduling of water releases through a hydroelectric power generation network of dams and pipelines, the transfer of funds between bank accounts, the exchange of funds among different currencies, and so on. Occasionally, multiple but distinct commodities share the resources of a common underlying network and the resulting models become more complicated. A common example is the movement of traffic along a highway system, where vehicles with different destinations and of different sizes share the same roads.

The study of network flow problems predates the developments in mathematical programming. Its origins can be traced to the work of Gustav Kirchoff who studied the equilibrium flow of current in electrical networks, and to D. König's research on linear graphs. The study of the effective distribution of products over a transportation system was initiated, independently, by the Russian mathematician L.V. Kantorovich and the Dutch economist T.C. Koopmans in the 1940s. They studied the problem of allocating tasks to machines, and the problem of optimum utilization of transportation networks for the distribution of cargo. The classic formulation of the *transportation problem* is often referred to as the Hitchcock-Koopmans transportation problem, after the work of F.L. Hitchcock who first formulated it.

The field underwent major developments in the 1940s when G.B. Dantzig introduced the simplex algorithm. The fact that network problems are prevalent in real-world applications, together with the observation that they are usually extremely large, prompted the development of specialized algorithms for their efficient solution. The study of networks has occasionally motivated further developments in linear programming. For example, the use of price-directed decomposition methods for linear programming has its origins in network optimization algorithms.

Over the last half century, network optimization has led the developments in large-scale optimization. Some of the first genuinely parallel algorithms, such as the auction algorithm for the assignment problem formulated by D.P. Bertsekas in the 1970s, were developed specifically for network optimization. S.A. Zenios and R. Lasken in the 1980s were the first to successfully implement optimization algorithms for network problems on massively parallel computers, the Connection Machines CM–1 and CM–2.

This chapter addresses a class of network flow problems that are well-suited to parallel computations and solution by the algorithms introduced in Chapters 6 and 7. We focus especially on network optimization problems with quadratic or entropic objective functions, so that the algorithms of Chapter 6 are directly applicable. Formulations with linear cost functions can be solved by embedding the algorithms of this chapter within the PMD algorithm framework (see Chapter 3). This combination of algorithms is further investigated empirically in Section 15.6.

Notation is established in Section 12.1. Section 12.2 formulates the classic transportation problem and the multicommodity network flow problem, and some sample applications are discussed in Section 12.3. Section 12.4 develops iterative algorithms for quadratic transportation problems and multicommodity network flow problems. Section 12.5 develops a model decomposition algorithm for the multicommodity network flow problem. Notes and references in Section 12.6 conclude this chapter.

## 12.1 Preliminaries

A *directed network*, or *directed graph* $\mathcal{G} = (\mathcal{N}, \mathcal{A})$ consists of a collection $\mathcal{N}$ of *nodes* together with a set $\mathcal{A}$ of ordered pairs $(i, j)$ of nodes. The elements of $\mathcal{A}$ are referred to as *arcs*, and for each arc $(i, j)$ the node $i$ is termed the *origin*, and node $j$ is termed the *destination*. A graph is *bipartite* if the set of nodes $\mathcal{N}$ can be partitioned into two subsets with all arcs in $\mathcal{A}$ leading from nodes of one subset to those of the other. The first subset is composed of the *origin* nodes, and the second subset contains the *destination* nodes.

With every directed network $(\mathcal{N}, \mathcal{A})$ we associate a *node-arc incidence matrix* as follows: we create a matrix with number of rows equal to the number of nodes, and number of columns equal to the number of arcs. Each column of the matrix has two nonzero entries. One entry is $+1$ and the other is a negative real number denoted by $-m_{ij}$. For a given column, corresponding to arc $(i, j)$, the $i$th row has the $+1$ entry, and the $j$th row has the entry $-m_{ij}$. Conversely, a network can be obtained from each matrix that has two nonzero entries per column. The network is obtained by associating each row of the matrix with a node, and by associating an arc between two nodes when the corresponding rows have nonzero entries. Each arc has its origin at the node corresponding to the row with the $+1$ entry, and its destination at the node corresponding to the row with the other nonzero entry. When the entries $m_{ij}$ are all equal to 1 the network problem is called *pure*. A network that is not pure is called *generalized.* All information related to the structure of a network is captured in the node-arc incidence matrix. Figure 12.1 illustrates a simple network.

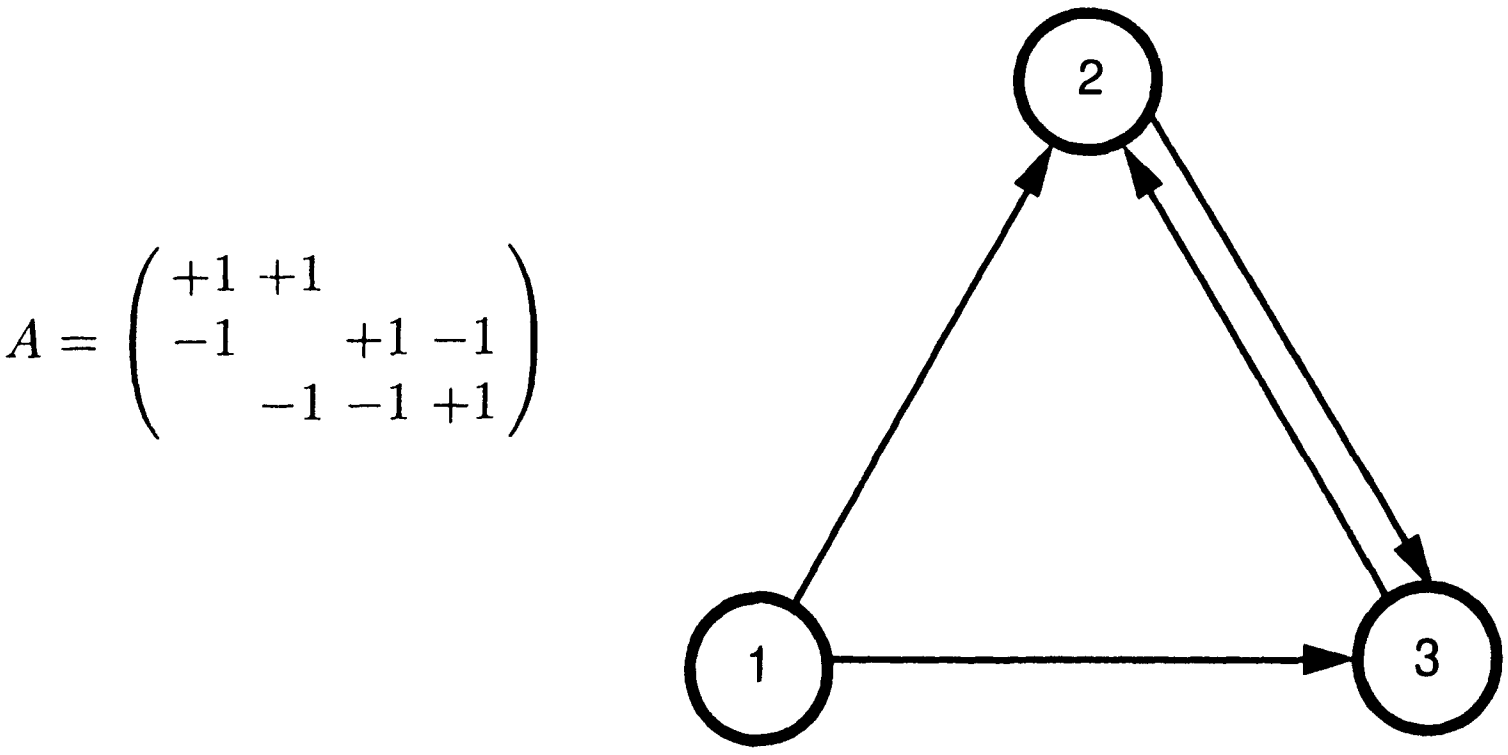

**Figure 12.1** An example of a simple node-arc incidence matrix, and the associated network structure.

With every arc $(i,j)$ we associate a variable $x_{ij}$ to denote the *flow* of some commodity from node $i$ to node $j$. When each arc has an *upper bound* $u_{ij}$ on the flow on arc $(i,j)$, then the network is *capacitated.* As in Section 9.2.1 the notation $\delta_i^+ \doteq \{j \mid (i,j) \in \mathcal{A}\}$ denotes the set of nodes having an arc with origin node $i$, and $\delta_j^- \doteq \{i \mid (i,j) \in \mathcal{A}\}$ denotes the set of nodes with destination node $j$.

## 12.2 Problem Formulations

The following problem will serve as the general formulation for multicommodity network flow problems:

$$\text{Minimize} \quad F(x) \tag{12.1}$$

$$\text{s.t.} \quad Ax = b, \tag{12.2}$$

$$Gx \leq U, \tag{12.3}$$

$$0 \leq x \leq u, \tag{12.4}$$

where:

$F : \mathbb{R}^n \to \mathbb{R}$ is the objective function, assumed to be convex and at least once continuously differentiable.

$x \in \mathbb{R}^n$ is the vector of decision variables representing flows on the network.

$A$ is an $m \times n$ constraints matrix with network structure. It could be the node-arc incidence matrix of a network flow problem, or it could be a *block-diagonal* matrix where each block is a node-arc incidence matrix. Such block-diagonal matrices appear in multicommodity network flow (Section 12.2.2), stochastic network flow (Section 13.7), and time-staged problems. We assume that $A \doteq \text{diag}\,(A_1, A_2, \ldots, A_K)$, and that each $A_k$, $k = 1, 2, \ldots, K$ is a node–arc incidence matrix.

$G$ is the $m_c \times n$ matrix of complicating (i.e., nonnetwork) constraints.

$u \in \mathbb{R}^n$ is the vector of upper bounds on the variables.

$b \in \mathbb{R}^m, U \in \mathbb{R}^{m_c}$ are the right-hand side vectors of the constraints.

The constraints $Ax = b$ express *conservation of flow* conditions. That is, the total flow into a given node $i$ minus the total flow out of the node must be equal to the exogenous coefficient $b_i$. The following standing assumption is made throughout this section.

**Assumption 12.2.1** *The set of feasible solutions to problem (12.1)–(12.4) is nonempty.*

Specializations of problem (12.1)–(12.4) yield several well-known formulations for network flow problems. Two such formulations, i.e., the transportation problem and the multicommodity network flow problem, are given next. The split-variable formulation of the deterministic equivalent stochastic network problems (discussed in Section 13.7) can also be cast in the form of (12.1)–(12.4) .

### 12.2.1 Transportation problems

Consider the special case of problem (12.1)–(12.4) without the complicating constraints $Gx \leq U$ and with only a single block $K = 1$ in the matrix $A$, i.e., $A \doteq A_1$. Assume, further, that the network corresponding to this node-arc incidence matrix is bipartite, with a subset $\{i \mid i = 1, 2, \ldots, m_O\} \subset \mathcal{N}$ of origin nodes, and a subset $\{j' \mid j = 1, 2, \ldots, m_D\} \subset \mathcal{N}$ of destination nodes. Problem (12.1)–(12.4) becomes the classic transportation problem. The cardinality of the set $\mathcal{N}$ is $m = m_O + m_D$. The node-arc incidence matrix of the transportation problem can be partitioned as

$$A = \begin{pmatrix} S \\ D \end{pmatrix}, \tag{12.5}$$

where $S$ is the $m_O \times n$ matrix

$$S \doteq \begin{pmatrix} \alpha_{11} \dots \alpha_{1m_D} & & & \\ & \alpha_{21} \dots \alpha_{2m_D} & & \\ & & \ddots & \\ & & & \alpha_{m_O 1} \dots \alpha_{m_O m_D} \end{pmatrix} \tag{12.6}$$

with entries $\alpha_{ij}$ given by

$$\alpha_{ij} \doteq \begin{cases} 1, & \text{if } j \in \delta_i^+, \\ 0, & \text{otherwise.} \end{cases} \tag{12.7}$$

The remaining entries of $S$ are all zero. $D$ is the $m_D \times n$ matrix

$$D \doteq \begin{pmatrix} \beta_{11} & & \beta_{21} & & & \beta_{m_O 1} & \\ & \ddots & & \ddots & \dots & & \ddots \\ & \beta_{1m_D} & & \beta_{2m_D} & & & \beta_{m_O m_D} \end{pmatrix} \tag{12.8}$$

with entries $\beta_{ij}$ given by

$$\beta_{ij} \doteq \begin{cases} -1, & \text{if } i \in \delta_j^-, \\ 0, & \text{otherwise.} \end{cases} \tag{12.9}$$

The remaining entries of $D$ are all zero. Since the $\beta_{ij}$'s are either zero or $-1$ we have a pure network problem. The right-hand side vector $b \in \mathbb{R}^m$ can be partitioned similarly as

$$b = \begin{pmatrix} s \\ d \end{pmatrix},$$

where $s \in \mathbb{R}^{m_O}$ is termed the *supply* vector and $d \in \mathbb{R}^{m_D}$ is the *demand* vector. The decision variables vector $x$ is a lexicographic ordering of the flow variables $x_{ij}$ for all arcs in $\mathcal{A}$. Problem (12.1), (12.2), and (12.4) can now be rewritten in the form:

$$\text{Minimize} \quad F(x) \tag{12.10}$$

$$\text{s.t.} \quad Sx = s, \tag{12.11}$$

$$Dx = d, \tag{12.12}$$

$$0 \leq x \leq u. \tag{12.13}$$

This is the *transportation problem.* The objective function is usually separable, that is, it can be expressed as the sum of functions of a single variable in the form $F(x) \doteq \sum_{(i,j)\in\mathcal{A}} f_{ij}(x_{ij})$, where $f_{ij} : \mathbb{R} \to \mathbb{R}$ for all $(i,j) \in \mathcal{A}$ are also convex and at least once continuously differentiable. The separable transportation problem is expressed in a componentwise notation as follows:

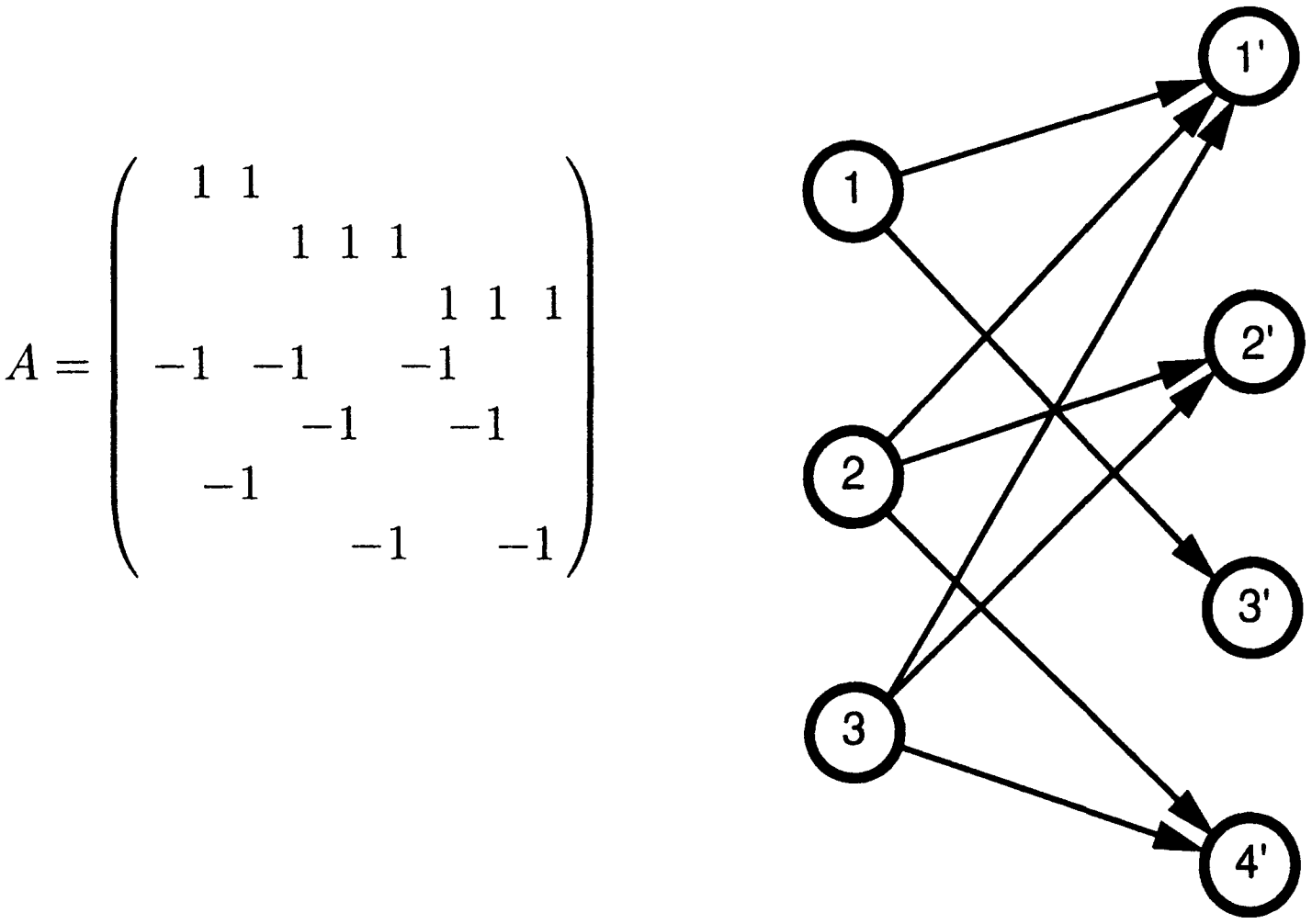

**Figure 12.2** An example of a node-arc incidence matrix for a simple transportation network and the associated bipartite graph.

$$\text{Minimize} \quad \sum_{(i,j)\in\mathcal{A}} f_{ij}(x_{ij}) \tag{12.14}$$

$$\text{s.t.} \quad \sum_{j\in\delta_i^+} x_{ij} = s_i, \qquad i = 1,2,\ldots,m_O, \tag{12.15}$$

$$\sum_{i\in\delta_j^-} x_{ij} = d_j, \qquad j = 1,2,\ldots,m_D, \tag{12.16}$$

$$0 \le x_{ij} \le u_{ij}, \qquad \text{for all } (i,j) \in \mathcal{A}. \tag{12.17}$$

Figure 12.2 illustrates the graph of a transportation problem, and its node-arc incidence matrix.

### 12.2.2 Multicommodity network flow problems

In many applications of network models several distinct commodities share the resources of a common network. The flow of each commodity on the network is governed by its own conservation of flow equations, but the commodities also share some common resources. Typically, each arc is utilized by multiple commodities that share its limited capacity. The interaction among commodities must be modeled explicitly.

The *multicommodity network flow* problem can be cast in the form of problem (12.1)–(12.4). Assume that there are $K > 1$ distinct commodities, and let

$$A \doteq \text{diag}\,(A_1, A_2, \ldots, A_K),$$

where $A_k$, $k = 1, 2, \ldots, K$ denotes the node-arc incidence matrix for the $k$th commodity. All matrices $A_k$ are assumed to be of equal dimensions $m \times n$, and hence the matrix $A$ is of dimension $Km \times Kn$. The flow vector of the $k$th commodity is denoted by $x^k \in \mathbb{R}^n$, and the upper bound on this vector is denoted by $u^k \in \mathbb{R}^n$. The resource vector for the $k$th commodity is denoted by $b^k \in \mathbb{R}^m$. We also use boldface $\boldsymbol{x} \in \mathbb{R}^{nK}$ to denote the concatenated vector of the flow vectors of all commodities, i.e., $\boldsymbol{x} = ((x^1)^T, \ldots, (x^K)^T)^T$. Let the components of $U \in \mathbb{R}^n$ denote the total capacities of the network arcs that are shared by the multiple commodities, and let $F_k(x^k)$ denote the cost function for the $k$th commodity. The multicommodity network flow problem is then written as follows:

$$\text{Minimize} \quad \sum_{k=1}^{K} F_k(x^k) \tag{12.18}$$

$$\text{s.t.} \quad A_k x^k = b^k, \qquad k = 1, 2, \ldots, K, \tag{12.19}$$

$$0 \le x^k \le u^k, \qquad k = 1, 2, \ldots, K, \tag{12.20}$$

$$0 \le G\boldsymbol{x} \le U. \tag{12.21}$$

The matrix $G$ is of dimension $n \times Kn$ and it operates on the $Kn$-vector $\boldsymbol{x}$. Hence, it links the flow conservation constraints that are otherwise independent for each commodity. $G$ has the *generalized upper bounding* (GUB) structure $G \doteq (I_n \mid \cdots \mid I_n)$, where $I_n$ denotes the $n \times n$ identity matrix. $U \in \mathbb{R}^n$ denotes the upper bounds on the joint capacity constraints.

### The Multicommodity Transportation Problem

The multicommodity network flow problem (12.18)–(12.21) can be specialized further to the case wherein each matrix $A_k$ is of the form given in (12.5), corresponding to the node-arc incidence matrix of a transportation problem. Each matrix is written as

$$A_k = \begin{pmatrix} S_k \\ D_k \end{pmatrix}, \tag{12.22}$$

where each $S_k$ and $D_k$ for all $k = 1, 2, \ldots, K$ are of the form (12.6) and (12.8), respectively. The resource vector $b^k$ is also partitioned as

$$b^k = \begin{pmatrix} s^k \\ d^k \end{pmatrix},$$

where $s^k \in \mathbb{R}^{m_O}$ and $d^k \in \mathbb{R}^{m_D}$ are the respective supply and demand vectors of the $k$th commodity. If we define matrices $\tilde{S} \doteq \text{diag}\,(S_1, S_2, \ldots, S_K)$, and $\tilde{D} \doteq \text{diag}\,(D_1, D_2, \ldots, D_K)$, and vectors $\tilde{\boldsymbol{s}} = \left((s^1)^{\mathrm{T}}, \ldots, (s^K)^{\mathrm{T}}\right)^{\mathrm{T}}$, $\tilde{\boldsymbol{d}} = \left((d^1)^{\mathrm{T}}, \ldots, (d^K)^{\mathrm{T}}\right)^{\mathrm{T}}$, and $\tilde{\boldsymbol{u}} = \left((u^1)^{\mathrm{T}}, (u^2)^{\mathrm{T}}, \ldots, (u^K)^{\mathrm{T}}\right)^{\mathrm{T}}$, then we can express the *multicommodity transportation problem* in matrix form as

$$
\begin{aligned}
&\text{Minimize} && F(\boldsymbol{x}) \\
&\text{s.t.} && \tilde{S}\boldsymbol{x} = \tilde{\boldsymbol{s}}, && (12.23)\\
& && \tilde{D}\boldsymbol{x} = \tilde{\boldsymbol{d}}, && (12.24)\\
& && 0 \le \boldsymbol{x} \le \tilde{\boldsymbol{u}}, && (12.25)\\
& && G\boldsymbol{x} \le U. && (12.26)
\end{aligned}
$$

In order to express the problem componentwise some additional notation is needed. For each commodity $k = 1, 2, \ldots, K$, let

$x^k = (x^k_{ij}) \in \mathbb{R}^n$, $(i,j) \in \mathcal{A}$ be the vector of flows,
$u^k = (u^k_{ij}) \in \mathbb{R}^n$, $(i,j) \in \mathcal{A}$ be the vector of upper bounds on the flows,
$s^k = (s^k_i) \in \mathbb{R}^{m_O}$, $i = 1, 2, \ldots, m_O$ be the vector of supplies, and
$d^k = (d^k_j) \in \mathbb{R}^{m_D}$, $j = 1, 2, \ldots, m_D$ be the vector of demands.

We also let $U = (U_{ij}) \in \mathbb{R}^n$, $(i,j) \in \mathcal{A}$ be the vector of mutual arc capacities. The objective function is usually separable and of the form

$$
F(\boldsymbol{x}) \doteq \sum_{k=1}^{K} F_k(x^k) \doteq \sum_{k=1}^{K} \sum_{(i,j)\in\mathcal{A}} f^k_{ij}(x^k_{ij}).
$$

Thus the multicommodity transportation problem can be rewritten componentwise as follows:

$$
\begin{aligned}
&\text{Minimize} && \sum_{k=1}^{K} \sum_{(i,j)\in\mathcal{A}} f^k_{ij}(x^k_{ij}) && (12.27)\\
&\text{s.t.} && \sum_{j\in\delta^+_i} x^k_{ij} = s^k_i, \quad i = 1, 2, \ldots, m_O,\ k = 1, 2, \ldots, K, && (12.28)\\
& && \sum_{i\in\delta^-_j} x^k_{ij} = d^k_j, \quad j = 1, 2, \ldots, m_D,\ k = 1, 2, \ldots, K, && (12.29)\\
& && 0 \le x^k_{ij} \le u^k_{ij}, \qquad \text{for all } (i,j) \in \mathcal{A}, && (12.30)\\
& && \sum_{k=1}^{K} x^k_{ij} \le U_{ij}, \qquad \text{for all } (i,j) \in \mathcal{A}. && (12.31)
\end{aligned}
$$

The constraints (12.28) represent conservation of flow at all origin nodes and for all commodities. Similarly, the constraints (12.29) represent conservation of flow at the destination nodes. The constraints (12.31) impose an upper bound on the total flow on each arc, and they are the *generalized upper bounding* (GUB) constraints. Figure 12.3 illustrates the structure of the constraints matrix of the multicommodity network flow problem for the constraints (12.28)–(12.29) and (12.31)

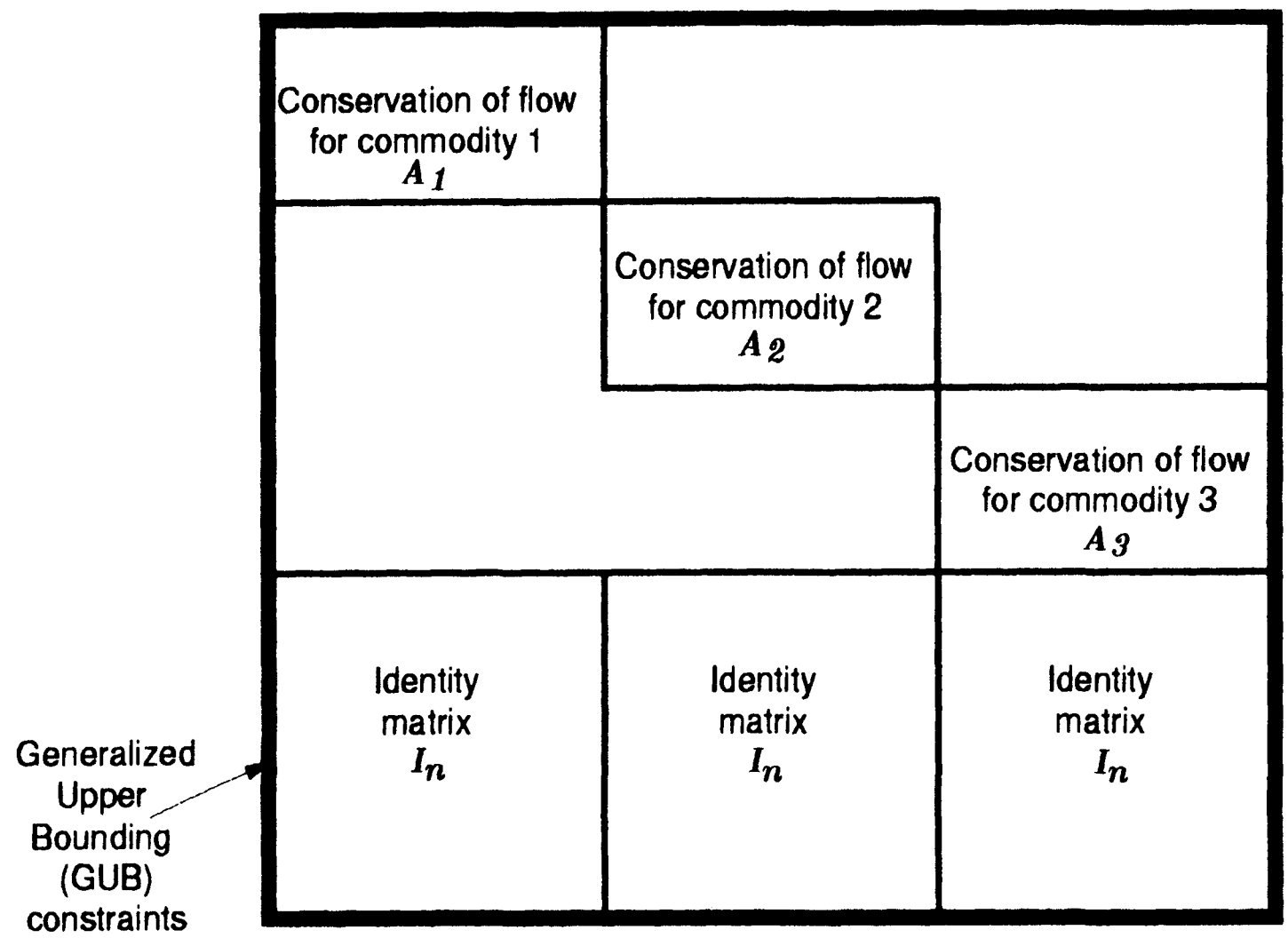

**Figure 12.3** Structure of the constraint matrix for multicommodity transportation problems.

## 12.3 Sample Applications

We now discuss three applications of network flow models. The domain of applications of network models is extensive, and this section is merely illustrative. We consider three diverse problem areas: finance, air-traffic scheduling, and transportation. The matrix balancing problem discussed in Chapter 9 is also amenable to formulation as a nonlinear transportation problem (see Section 9.2). Section 12.6 gives additional references that point readers to the vast literature on network applications.

### 12.3.1 Example one: Covering positions in stock options

Regulations for stock option investments require that a *margin* be provided before a position in a stock option is taken. This margin is a form of insurance. The total margin requirements can be difficult to calculate when an investor takes a complex position in stocks and options. In particular, different options on the same security may be paired together in order to reduce the margin requirements, since some options act as insurance (hedge) against each other.

The calculation of minimum margin can be formulated as a transportation problem. Assume that there are $Z_j$ long stocks of the $j$th stock, and $Y_i$ short stocks of the $i$th stock. (A position in stocks is *long* if there is an agreement to buy stocks at some future point in time at a predetermined price, while a position is *short* if there is an agreement to sell stocks at some future time at a predetermined price.) Let $m_{ij}$ denote the margin required for a position involving one share or option of security $i$ paired with one option or share of security $j$. The minimum margin can be determined by deciding which securities should be paired together. Denote by $x_{ij}$ the number of contracts paired between securities $i$ and $j$. The following transportation model determines the minimum margin:

$$\begin{aligned}
\text{Minimize} \quad & \sum_i \sum_j m_{ij} x_{ij} & \\
\text{s.t.} \quad & \sum_i x_{ij} = Z_j, & \text{for all } j, \\
& \sum_j x_{ij} = Y_i, & \text{for all } i, \\
& x_{ij} \geq 0, & \text{for all } i \text{ and } j.
\end{aligned}$$

### 12.3.2 Example two: Air-traffic control

The allocation of aircraft to high-altitude jet routes is an important aspect of the air-traffic control system. Aircraft are assigned jet routes in a sequence that allows them to reach their destination in the shortest period of time. However, delays may be necessary to avoid congestion on the preferred routes, especially during peak travel hours.

Airports and control towers form the nodes of the air-traffic control network. The flight paths, between these control points, are the arcs of the network. Conservation of flow constraints dictate that the required number of flights is allowed to depart from each origin airport, and that the scheduled number of flights is allowed to land at each destination airport. In addition, all aircraft that move into the jurisdiction of a control tower are eventually passed on to the jurisdiction of some other tower or to their destination airport.

The congestion of a jet route is a complicated nonlinear function of variables that depend on all aircraft flying on the route. Upper bounds are imposed on the total number of flights to maintain acceptable safety limits. The network optimization model has dual, and conflicting, objectives: allocate aircraft to the most cost-effective routes, and allocate aircraft to different routes so that congestion is reduced.

For military applications the air-traffic control problem usually encompasses both the problem of routing the aircraft as well as the problem of scheduling cargo on the aircraft. Cargo of different types (e.g., personnel, vehicles, supplies) is loaded on planes and then routed to their required destination. The network model becomes a multicommodity flow problem, since the different types of cargo share the limited capacity of each aircraft.

### 12.3.3 Example three: Routing of traffic

In many applications of multicommodity flow models, commodities are distinguished by their different origins and destinations. While the same entity flows on each arc of the network, each unit of flow is characterized by its origin/destination pair. Hence, flow cannot be simply aggregated at each node in order to satisfy conservation of flow, but instead, conservation has to be satisfied by all flows with the same origin/destination pair. This situation arises in routing applications when different nodes on a network are sending information or products to each other.

In communication networks (such as telephone or computer networks) messages are transmitted between different origins and destinations sharing the same arcs. In these applications arcs are the communication links, e.g., wires. Origin and destination nodes represent different calling stations, storage devices, computer terminals, and the like. Node conservation conditions at each node of the network ensure that no messages are lost and that all messages are routed to their destination. The capacities of the transmission lines define the joint capacity constraints of the multicommodity flow problem. The capacity of a transmission line could be finite. For example, if a communication link is established using a switching device, then there is a limit on the total number of messages that can be routed through the switch per unit time interval. In other applications the capacity is virtually unlimited, but very long delays can occur when the communication link is congested. Delays are modeled using nonseparable functions. The delay (therefore, cost) experienced by the $k$th commodity, if it were the only commodity using the network, is given by a separable nonlinear function. If $x^k$ denotes the units of the $k$th commodity using the transmission link, then the delay is the function of a single variable $f_k(x^k)$. The delay due to the presence of $K$ distinct commodities is a function of the total flow on the transmission link, and it takes the nonseparable form $F(\sum_{k=1}^{K} x^k)$. If the cost function is ignored the problem is decomposable by commodity, since the flow vector for each commodity satisfies different constraints. A different problem could be solved for each commodity. However, the cost function for a given commodity depends on the value of the flow vector of the remaining commodities, and the problems for each commodity are linked through the objective function. Such applications are usually called *traffic assignment problems.*

Another routing application arises in railroad transportation, when goods are transported through a network of tracks, yards, and junction points. Cars are loaded at each yard, and they can carry goods that have different destinations. These goods share the total capacity of the railroad system, determined by the number of cars that can be loaded onto a train. Conservation of flow at yards and junctions ensures that all goods are delivered to their destinations. The capacity of the train defines the joint capacity of the multicommodity flow problem.

Similar situations arise in the distribution of products from plants to retail stores using a fleet of trucks. Warehouses are used for temporary storage of the products to balance seasonal variations in demand. Multiple products may share the same storage space and distribution vehicles. Nodes of the underlying network model represent plants, warehouses, and retail stores; arcs represent transportation links such as truck routes. Conservation of flow at the nodes ensures that products reach the retail stores without any losses at the storage sites. The capacity of the warehouses defines the joint capacity constraints.

## 12.4 Iterative Algorithms for Multicommodity Network Flow Problems

In this section we develop algorithms for the multicommodity transportation problem. They are derived from the general row-action Algorithm 6.4.1 but they are formulated specifically to accommodate the unique structure of this problem. An algorithm for the single commodity quadratic transportation problem is developed first, and is then extended to solve generalized network problems with entropic objective functions. This algorithm can of course be simplified to solve pure networks with entropic objectives simply by setting all arc multipliers equal to 1. Finally the algorithms are extended to solve multicommodity network flow problems. We point out that problems with linear objective functions can be solved by embedding the algorithms of this section within the PMD framework (see Chapter 3, Section 15.6).

### 12.4.1 Row-action algorithm for quadratic transportation problems

Consider problem (12.14)–(12.17) with the quadratic objective function:

$$F(x) \doteq \sum_{(i,j)\in\mathcal{A}} \left(\frac{1}{2} w_{ij} x_{ij}^2 + c_{ij} x_{ij}\right), \tag{12.32}$$

where $\{w_{ij}\}$ and $\{c_{ij}\}$ are given positive real numbers.

The row-action iterative Algorithm 6.4.1 is applied to this problem, iterating on the constraints matrix one row at a time. We derive first the iterative step of the algorithm when the chosen row corresponds to the flow conservation constraints at the origin nodes, i.e., equation (12.15). The projection $y'$, of the current iterate $y$, onto the $i$th hyperplane $H(s^i, s_i) \doteq \{x \mid \sum_{j\in\delta_i^+} x_{ij} = s_i\}$ is the solution of the system

$$\nabla F(y') = \nabla F(y) + \beta_i s^i, \tag{12.33}$$

$$y' \in H(s^i, s_i), \tag{12.34}$$

where $\beta_i$ is the parameter associated with the generalized projection. (Recall that $s^i$ denotes the transpose of the $i$th column of the matrix $S$ in (12.6), while $s_i$ denotes the $i$th component of the supply vector.) If $y \in H(s^i, s_i)$, then $\beta_i = 0$ and $y' = y$. Solving the system given above means that if the current iterate $y$ does not satisfy flow conservation, i.e., $y \notin H(s^i, s_i)$, we compute a $y'$ that is related to $y$ via (12.33) and that satisfies flow conservation. Evaluating (12.33) when $F$ is the quadratic objective function (12.32), and using the structure of the $i$th row of the constraints matrix $s^i$ that corresponds to conservation of flow constraints at the origin node $i$, we get

$$y'_{ij} = y_{ij} + \frac{\beta_i}{w_{ij}}, \quad \text{for all } j \in \delta_i^+. \tag{12.35}$$

Equation (12.34) means that the vector $y' = (y'_{ij})$ must satisfy the flow conservation constraint (12.15) of the $i$th origin node. Substituting therefore the expressions for $y'_{ij}$ from (12.35) into the flow conservation constraint (12.15) we obtain

$$\sum_{j \in \delta_i^+} (y_{ij} + \frac{\beta_i}{w_{ij}}) = s_i. \tag{12.36}$$

Hence,

$$\beta_i = \frac{s_i - \sum_{j \in \delta_i^+} y_{ij}}{\sum_{j \in \delta_i^+} 1/w_{ij}}. \tag{12.37}$$

Using this value of $\beta_i$ in (12.35) yields the formula for the projected flow vector $y'$ on arcs incident to node $i$. The dual variable for the $i$th origin node (denoted by $\pi_i^O$) is updated by subtracting $\beta_i$ from its current value, $\pi_i^O \leftarrow \pi_i^O - \beta_i$ (refer to the iterative step of Algorithm 6.4.1). Similar algebraic manipulations lead to the updating formula for flows incident to the destination nodes.

Now we develop the projections onto the positive orthant $\mathbb{R}_+^n$ and on the upper bounds (12.17). Denote by $y_{ij}$ the current value of the flow variable on arc $(i, j)$. If $y_{ij} < 0$ then the iterative step takes the form

$$\nabla F(y') = \nabla F(y) + \beta e^t, \tag{12.38}$$

$$y' \in H(e^t, 0), \tag{12.39}$$

where $e^t$ is the standard basis vector having 1 in the $t$th coordinate and zeros elsewhere. Here $t$ corresponds to the lexicographic index of arc $(i, j)$.

Evaluating the first expression, for the quadratic objective function, and substituting it into the second expression we get

$$0 = y_{ij} + \frac{\beta}{w_{ij}}, \tag{12.40}$$

and the parameter associated with the generalized projection is

$$\beta = -w_{ij} y_{ij}. \tag{12.41}$$

Similarly, if $y_{ij} > u_{ij}$ we compute the generalized projection parameter

$$\beta = (u_{ij} - y_{ij}) w_{ij}, \tag{12.42}$$

and the primal and dual updating steps (see (6.61)) take the form

$$\begin{aligned} \Delta &\doteq \text{mid}\ (r_{ij},\ \ -w_{ij}y_{ij},\ \ (u_{ij} - y_{ij})w_{ij}\ ), \\ y'_{ij} &= y_{ij} + \frac{\Delta}{w_{ij}}, \\ r_{ij} &\leftarrow r_{ij} - \Delta, \end{aligned}$$

where $r_{ij}$ are the components of the current iterate of the dual variables associated with the constraints (12.17). It is easy to verify that $y_{ij} < 0$ implies that $y'_{ij} = 0$, and that $y_{ij} > u_{ij}$ implies $y'_{ij} = u_{ij}$.

The complete algorithm is summarized below. Relaxation parameters are allowed. To guarantee the asymptotic convergence of the algorithm all relaxation parameters $\lambda$ must be confined to an interval $0 < \epsilon \leq \lambda \leq 2 - \epsilon$, for some arbitrary fixed positive $\epsilon$.

**Algorithm 12.4.1 Row-Action Algorithm for Quadratic Transportation Problems**

**Step 0:** (Initialization.) Set $\nu = 0$. Take an arbitrary flow vector $x^0$, and initial dual prices vectors for origin nodes $(\pi^O)^0 = ((\pi_i^O)^0)$, for destination nodes $(\pi^D)^0 = ((\pi_j^D)^0)$, and for the simple upper bounds $r^0 = (r_{ij}^0)$ that satisfy the initialization step of Algorithm 6.4.1 for the current problem with the quadratic objective function (12.32). The matrix $\Phi$ in (6.60) represents the constraints (12.15)–(12.17). That is, for all $(i,j) \in \mathcal{A}$,

$$x_{ij}^0 = -\frac{1}{w_{ij}} \left( c_{ij} + (\pi_i^O)^0 + (\pi_j^D)^0 + r_{ij}^0 \right). \tag{12.43}$$

**Step 1:** (Iterative step for origin nodes.) For $i = 1, 2, \ldots, m_O$, pick relaxation parameters $\lambda_i^\nu$, calculate scaling parameters $\rho_i^\nu$, and update the primal and dual vectors as follows:

$$\rho_i^\nu = \frac{\lambda_i^\nu}{\sum_{j \in \delta_i^+} 1/w_{ij}} (s_i - \sum_{j \in \delta_i^+} x_{ij}^\nu), \tag{12.44}$$

$$x_{ij}^\nu \leftarrow x_{ij}^\nu + \frac{\rho_i^\nu}{w_{ij}}, \quad \text{for all } j \in \delta_i^+, \tag{12.45}$$

$$(\pi^O)_i^{\nu+1} = (\pi^O)_i^\nu - \rho_i^\nu. \tag{12.46}$$

**Step 2:** (Iterative step for destination nodes.) For $j = 1, 2, \ldots, m_D$, pick relaxation parameters $\lambda_j^\nu$, calculate scaling parameters $\sigma_j^\nu$, and update the primal and dual vectors as follows:

$$\sigma_j^\nu = \frac{\lambda_j^\nu}{\sum_{i \in \delta_j^-} 1/w_{ij}} (d_j - \sum_{i \in \delta_j^-} x_{ij}^\nu), \tag{12.47}$$

$$x_{ij}^\nu \leftarrow x_{ij}^\nu + \frac{\sigma_j^\nu}{w_{ij}}, \quad \text{for all } i \in \delta_j^-, \tag{12.48}$$

$$(\pi^D)_j^{\nu+1} = (\pi^D)_j^\nu - \sigma_j^\nu. \tag{12.49}$$

**Step 3:** (Iterative step for bounds on flows.) For all $(i,j) \in \mathcal{A}$, pick a relaxation parameter $\lambda_{ij}$, and calculate:

$$\Delta_{ij}^\nu \doteq \text{mid}\,(r_{ij}^\nu,\ \lambda_{ij} w_{ij}(u_{ij} - x_{ij}^\nu),\ -\lambda_{ij} w_{ij} x_{ij}^\nu), \tag{12.50}$$

$$x_{ij}^{\nu+1} = x_{ij}^\nu + \frac{\Delta_{ij}^\nu}{w_{ij}}, \tag{12.51}$$

$$r_{ij}^{\nu+1} = r_{ij}^\nu - \Delta_{ij}^\nu. \tag{12.52}$$

**Step 4:** Let $\nu \leftarrow \nu + 1$ and return to Step 1.

### Decompositions for Parallel Computing

The algorithm decomposes naturally for parallel computations. The execution of Step 1 can proceed concurrently for all origin nodes $i = 1, 2, \ldots, m_O$. The calculation of the parameter $\rho_i^\nu$ involves only the current flow on arcs incident to node $i$, and since no two origin nodes have arcs in common it is possible to compute the parameters $\rho_i^\nu$, for all values of $i$, independent of each other. The algorithm could utilize up to $m_O$ processors in calculating these parameters for all indices $i$ simultaneously. Additional processors can be utilized to compute the sum of the flows on all arcs incident to each node. Similarly, the flow variables can be updated concurrently for all arcs with destination nodes $j \in \delta_i^+$, for all indices $i$, using up to $n$ processors. The dual prices can also be updated concurrently for all $i$ using $m_O$ processors. Step 2 of the algorithm parallelizes in a way similar to Step 1. Step 3 can be executed concurrently for all arcs. The calculation of the parameters $\Delta_{ij}$ and the updating of $x_{ij}^{\nu+1}$ and $r_{ij}^{\nu+1}$ need only the current values of $x_{ij}^{\nu}$ and $r_{ij}^{\nu}$ respectively, for the given arc. Hence the flows and the dual variables for all arcs can be updated concurrently using up to $n$ processors. Chapter 14 gives data-structures for the parallel implementation of the algorithm.

### 12.4.2 Extensions to generalized networks

The algorithm for pure quadratic transportation problems developed above, uses simple linear expressions for computing the scaling parameters and updating the primal variable. The general row-action Algorithm 6.4.1 requires the solution of a nonlinear system of equations in order to find the generalized projection parameter at the iterative step (6.61) (refer, for example, to the solution of (12.33)–(12.34)). For a quadratic objective function the system becomes linear and it is possible to obtain a closed-form solution for the parameter, denoted by $\beta$, as in the previous section. Closed-form solutions can also be obtained with other nonlinear functions (whose gradients are also nonlinear) when the constraints are those of a pure network, because the constraint coefficients of such a network are either 0 or 1. For generalized network problems, however, it is not always possible to solve analytically for the generalized projection parameter and each step of the algorithm would require an iterative procedure to solve the nonlinear system of equations. If the objective function $F$ is twice continuously differentiable these equations are smooth and Newton's method can be used.

In some cases it is possible to obtain in closed-form, an *approximate* solution to the nonlinear system of equations and still preserve the overall asymptotic convergence of the row-action algorithm (see Section 6.9). Obtaining such closed-form solutions has important implications for the efficient implementation of an algorithm. First, the computationally expensive part of invoking a nonlinear equation solver is avoided. This may result in a slight increase in the total number of iterations of the row-action algorithm before a target accuracy of the solution can be achieved, however numerical experiences indicate that the increased number of iterations does not adversely affect the total solution time since each iteration is executed very efficiently. Second, it is possible to implement the algorithm on massively parallel machines, even if they are SIMD computers, since the analytic expression for the generalized projection parameter can be evaluated simultaneously for multiple rows of the constraint set (see also Chapter 13, p. 393). We illustrate here the derivation of the row-action algorithm for generalized networks, using a closed-form approximate calculation for the generalized projection parameter.

Consider the generalized transportation problem, with the entropy objective function

$$\text{Minimize} \quad F(x) \doteq \sum_{(i,j)\in\mathcal{A}} x_{ij}\left(\log\left(\frac{x_{ij}}{a_{ij}}\right) - 1\right) \tag{12.53}$$

$$\text{s.t.} \quad \sum_{j\in\delta_i^+} x_{ij} \leq s_i, \qquad i = 1,2,\ldots,m_O, \tag{12.54}$$

$$\sum_{i\in\delta_j^-} m_{ij}x_{ij} = d_j, \qquad j = 1,2,\ldots,m_D, \tag{12.55}$$

$$0 \leq x_{ij} \leq u_{ij}, \qquad \text{for all } (i,j) \in \mathcal{A}, \tag{12.56}$$

where $\{a_{ij}\}$ are given positive real numbers. This problem is written in matrix form as

$$\text{Minimize} \quad F(x) \tag{12.57}$$

$$\text{s.t.} \quad Sx \leq s, \tag{12.58}$$

$$Dx = d, \tag{12.59}$$

$$0 \leq x \leq u. \tag{12.60}$$

$S$ is as defined in (12.6). $D$ is structurally identical to the matrix defined in (12.8), but with entries given by

$$\beta_{ij} \doteq \begin{cases} m_{ij}, & \text{if } i \in \delta_j^-, \\ 0, & \text{otherwise.} \end{cases} \tag{12.61}$$

Let $\Phi \doteq \begin{pmatrix} S \\ D \\ I \end{pmatrix}$ where $I$ is the $n \times n$ identity matrix. Let also

$$\gamma \doteq \begin{pmatrix} -\infty \\ d \\ 0 \end{pmatrix} \quad \text{and} \quad \delta \doteq \begin{pmatrix} s \\ d \\ u \end{pmatrix}$$

Problem (12.57)–(12.60) can be written in the form

$$\text{Minimize} \quad F(x)$$

$$\text{s.t.} \quad \gamma \leq \Phi x \leq \delta. \tag{12.62}$$

Denote by $z \in \mathbb{R}^{m_O+m_D+n}$ the vector of dual variables

$$z = \left((\pi^O)^T, (\pi^D)^T, (r)^T\right)^T$$

and apply the general row-action Algorithm 6.4.1 to this formulation. The iterative step (6.61), applied to the $l$th row $1 \leq l \leq m_O$ of the constraints matrix takes the form

$$\begin{aligned} y'_t &= y_t \exp(\beta \phi^l_t), \quad t = 1, 2, \ldots, n, \\ z &\leftarrow z - \beta e^l. \end{aligned} \tag{12.63}$$

For $l = 1, 2, \ldots, m_O$ we deal with rows of (12.58), and the generalized projection parameter $\beta$ is obtained from the iterative step of Algorithm 6.4.1 as

$$\beta = \min(\pi^O_l, \Delta), \tag{12.64}$$

where $\Delta$ is obtained by solving the system

$$y'_t = y_t \exp(\Delta \phi^l_t), \quad t = 1, 2, \ldots, n, \tag{12.65}$$

$$\langle y', \phi^l \rangle = s_l. \tag{12.66}$$

Substituting $y'$ from the first equation into the second, and using the facts that $\phi^l$ is the $l$th column of $S^T$ when $l = i = 1, 2, \ldots, m_O$ with values $\alpha_{ij} \in \{0, 1\}$ as defined in (12.7), and that $(y_t)$ is the lexicographically ordered vector $(y_{ij})$, for $j \in \delta^+_i$, we get $\exp \Delta = s_i / \sum_{j \in \delta^+_i} y_{ij}$. This value of $\Delta$ goes into (12.64) to get $\beta$, which is then used in (12.63) to complete the iterative step.

For $l = m_O + 1, m_O + 2, m_O + 3, \ldots, m_O + m_D$ we deal with rows of (12.59), and the generalized projection parameter is obtained from the iterative step of Algorithm 6.4.1 by solving the system

$$y'_t = y_t \exp(\beta \phi^l_t), \quad t = 1, 2, \ldots, n, \tag{12.67}$$

$$\langle y', \phi^l \rangle = d_l. \tag{12.68}$$

Substituting $y'$ from the first equation into the second and using the facts that $\phi^l$ is the $i \doteq (l - m_O)$th column of $D^T$ with values $\beta_{ij}$ as given in (12.61), and that $(y_t)$ is the lexicographically ordered vector $(y_{ij})$ for $i \in \delta^-_j$, we get the following nonlinear equation in $\beta$

$$\sum_{i \in \delta^-_j} m_{ij} x_{ij} \exp(\beta m_{ij}) = d_j. \tag{12.69}$$

Denote now $\omega \doteq \exp \beta$ and rewrite (12.69) as a nonlinear equation in $\omega$

$$\psi(\omega) \doteq \sum_{i \in \delta^-_j} m_{ij} x_{ij} \omega^{m_{ij}} - d_j = 0. \tag{12.70}$$

Solving this equation for $\omega$ requires an iterative procedure (e.g., Newton's algorithm). However, it is possible to simplify the algorithm by computing a MART-like step (see Section 6.9), which eliminates the need to solve this equation exactly. An approximate solution of (12.70) can be obtained in closed-form, by taking one secant step. Consider the line through $(0, -d_j)$ and $(1, \psi(1))$, as illustrated in Figure 12.4. At the $\nu$th iteration this line intersects the $\omega$-axis at

$$\hat{\omega}_\nu = \frac{d_j}{\sum_{i \in \delta_j^-} m_{ij} x_{ij}^\nu}. \tag{12.71}$$

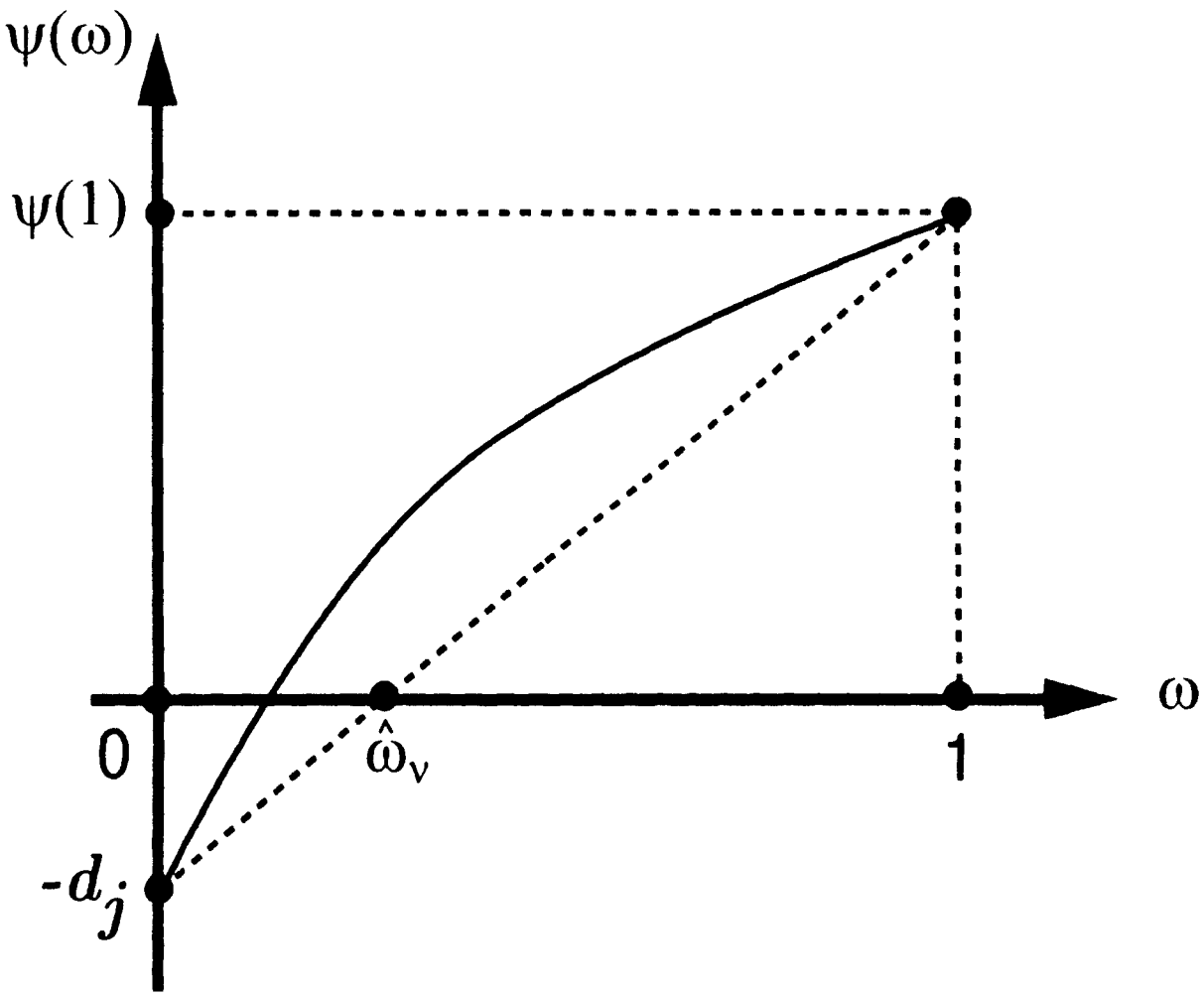

**Figure 12.4** Approximate solution $\hat{\omega}_\nu$ for the projection parameter of the entropy optimization algorithm for generalized networks.

We use this value of $\hat{\omega}_\nu$ as an approximation to the exact solution of (12.70) for the $\nu$th iteration. This approximation preserves asymptotic convergence of the algorithm (see Section 12.6). The value $\hat{\omega}_\nu$ is used to set $\beta_\nu = \log \hat{\omega}_\nu$ and, in turn, to complete the iterative step in (12.63).

In order to avoid the exponential terms the algorithm works with the exponents of the dual prices, which are denoted by $\overline{\pi}_i^O, \overline{\pi}_j^D$ and $\overline{r}_{ij}$, and makes use of the equality $\log(\min(\exp \alpha, \exp \beta)) = \min(\alpha, \beta)$. The algorithm also permits a user-chosen relaxation parameter $\lambda$, and to preserve asymptotic convergence it is required that $0 < \epsilon \leq \lambda \leq 1$. Before stating the algorithm we explain how the formula in the initialization step is derived. According to (6.60), using the entropy function (12.53), the matrix $\Phi$ of (12.62) with the appropriate $S$ and $D$, and the dual vector $z$, we write

$$\log\left(\frac{x_{ij}^0}{a_{ij}}\right) = \left(-\,(S^{\mathrm{T}} \mid D^{\mathrm{T}} \mid I)\begin{pmatrix} (\pi^O)^0 \\ (\pi^D)^0 \\ r^0 \end{pmatrix}\right)_{ij}.$$

With $\alpha_{ij}$ as in (12.7) and $\beta_{ij}$ as in (12.61) we obtain, after exponentiation,

$$x_{ij}^0 = \frac{a_{ij}((\overline{\pi}^D)_j^0)^{m_{ij}}}{(\overline{\pi}^O)_i^0 \overline{r}_{ij}^0}.$$

The entropy optimization algorithm for generalized transportation problems is stated as follows.

**Algorithm 12.4.2 Row-action Algorithm for Entropy Optimization Generalized Transportation Problems**

**Step 0:** (Initialization.) Set $\nu = 0$. Get $x^0 > 0, (\overline{\pi}^O)^0, (\overline{\pi}^D)^0$, and $\overline{r}^0$ that satisfy the initialization step of Algorithm 6.4.1 for the present problem with the entropy function (12.53), that is,

$$x_{ij}^0 = \frac{a_{ij}((\overline{\pi}^D)_j^0)^{m_{ij}}}{(\overline{\pi}^O)_i^0 \overline{r}_{ij}^0}. \tag{12.72}$$

**Step 1:** (Iterative step for origin nodes.) For $i = 1, 2, \ldots, m_O$, pick relaxation parameters $\lambda_i^\nu$, and calculate:

$$\rho_i^\nu = \left( \frac{s_i}{\sum_{j \in \delta_i^+} x_{ij}^\nu} \right)^{\lambda_i^\nu}, \tag{12.73}$$

$$\Delta\rho_i^\nu \doteq \min\left((\overline{\pi}^O)_i^\nu, \rho_i^\nu\right), \tag{12.74}$$

$$x_{ij}^\nu \leftarrow x_{ij}^\nu \Delta\rho_i^\nu, \quad \text{for all } j \in \delta_i^+, \tag{12.75}$$

$$(\overline{\pi}^O)_i^{\nu+1} = \frac{(\overline{\pi}^O)_i^\nu}{\Delta\rho_i^\nu}. \tag{12.76}$$

**Step 2:** (Iterative step for destination nodes.) For $j = 1, 2, \ldots, m_D$, pick relaxation parameters $\lambda_j^\nu$, and calculate:

$$\sigma_j^\nu = \left( \frac{d_j}{\sum_{i \in \delta_j^-} m_{ij} x_{ij}^\nu} \right)^{\lambda_j^\nu}, \tag{12.77}$$

$$x_{ij}^\nu \leftarrow x_{ij}^\nu (\sigma_j^\nu)^{m_{ij}}, \quad \text{for all } i \in \delta_j^-, \tag{12.78}$$

$$(\overline{\pi}^D)_j^{\nu+1} = \frac{(\overline{\pi}^D)_j^\nu}{\sigma_j^\nu}. \tag{12.79}$$

**Step 3:** (Iterative step for bounds on flows.) For all $(i,j) \in \mathcal{A}$, pick relaxation parameters $\lambda_{ij}$, and calculate:

$$\Delta_{ij}^{\nu} \doteq \min\left(\overline{r}_{ij}^{\nu}, \left(\frac{u_{ij}}{x_{ij}^{\nu}}\right)^{\lambda_{ij}}\right), \tag{12.80}$$

$$x_{ij}^{\nu+1} = x_{ij}^{\nu}\Delta_{ij}^{\nu}, \tag{12.81}$$

$$\overline{r}_{ij}^{\nu+1} = \frac{\overline{r}_{ij}^{\nu}}{\Delta_{ij}^{\nu}}. \tag{12.82}$$

**Step 4:** Replace $\nu \leftarrow \nu + 1$ and return to Step 1.

Note that in Step 3 it is not necessary to perform a projection on the nonnegative orthant since $x^0 > 0$ and the multiplicative adjustment of $x$ in Steps 1 and 2 preserve the positivity of the primal variables.

An algorithm for pure entropy optimization problems can be derived from Algorithm 12.4.2, simply by setting the multiplier $m_{ij}$ equal to one, in equations (12.72), (12.77), and (12.78). When dealing with pure networks we can also require that the inequality constraint (12.58) be an equality. The resulting algorithm is, with the exception of Step 3 for the variable bounds, identical to the RAS algorithm (Algorithm 9.3.2) for matrix balancing.

### Decompositions for Parallel Computing

The calculations involved in this algorithm are identical in structure to the calculations of Algorithm 12.4.1, and the ideas described above for parallelizing the algorithm for the pure network problems apply to generalized network problems as well. Chapter 14 gives data-structures for the parallel implementation of the algorithm.

## 12.4.3 Row-action algorithm for quadratic multicommodity transportation problems

We turn now to the multicommodity transportation problem (12.27)–(12.31), with the quadratic objective function

$$F(x) \doteq \sum_{k=1}^{K} \sum_{(i,j)\in\mathcal{A}} \left(\frac{1}{2} w_{ij}^k (x_{ij}^k)^2 + c_{ij}^k x_{ij}^k\right), \tag{12.83}$$

where $\{w_{ij}^k\}$ and $\{c_{ij}^k\}$ are given positive real numbers. We apply the row-action Algorithm 6.4.1 to this problem. For any fixed commodity index $k$ the constraints (12.28)–(12.30) are identical to the constraints of the transportation problem (12.15)–(12.17). Hence, for these constraints the iterative steps of the algorithm are identical to the steps derived earlier for the quadratic transportation problem, and they are repeated for each commodity $k$.

To formulate the algorithm for the specific multicommodity problem we need to derive the iterative step for rows corresponding to the GUB constraints (12.31). We introduce the vector $\psi \in \mathbb{R}^n$ to denote dual prices for the GUB constraints. The $l$th component of this vector is the dual price of the GUB constraint of arc $(i,j)$, where the mapping of pairs $(i,j)$ onto a single index $l$ is done lexicographically according to $l \doteq (i-1)m+j$. This component is denoted by $\psi_{ij}$. Assume that the current and next iterates are $y = (y_{ij}^k)$ and $y' = (y_{ij}^k)'$ respectively, for $(i,j) \in \mathcal{A}$, $k = 1,2,\ldots,K$. If $y$ is such that $\sum_{k=1}^K y_{ij}^k < 0$ then the point $y$ is projected onto the hyperplane $H(g^l, 0) \doteq \{x \mid \sum_{k=1}^K x_{ij}^k = 0\}$ where $l$ denotes the lexicographic index of the GUB constraint (12.31) corresponding to the $(i,j)$th arc and $g^l$ is the $l$th row of matrix $G \doteq (I_n \mid I_n \mid \cdots \mid I_n)$ (see (12.21)). The iterative step of the algorithm takes the form

$$\nabla F(y') = \nabla F(y) + \Gamma e^l, \tag{12.84}$$

$$y' \in H(g^l, 0). \tag{12.85}$$

Evaluating (12.84) when $F$ is the quadratic objective function (12.83), and using $\Gamma$ to denote the projection parameter when projecting on the lower bounds, we obtain the primal updating step as: $(y_{ij}^k)' = y_{ij}^k + \Gamma/w_{ij}^k$.

Since $y'$ must be such that $\sum_{k=1}^{K} (y_{ij}^k)' = 0$ we obtain

$$\Gamma = -\frac{1}{\sum_{k=1}^{K} 1/w_{ij}^k} \sum_{k=1}^{K} y_{ij}^k. \tag{12.86}$$

Similarly, if $y$ is such that $\sum_{k=1}^K y_{ij}^k > U_{ij}$, then $y$ is projected onto the hyperplane denoted by $H(g^l, U_{ij}) \doteq \{x \mid \sum_{k=1}^K x_{ij}^k = U_{ij}\}$ where $l$ and $g^l$ are the same as above. The primal updating step is, therefore, $(y_{ij}^k)' = y_{ij}^k + \Delta/w_{ij}^k$, where $\Delta$ denotes the projection parameter. Since $y'$ must satisfy the constraint $\sum_{k=1}^K x_{ij}^k = U_{ij}$ we obtain

$$\Delta = \frac{1}{\sum_{k=1}^{K} 1/w_{ij}^k} \left( U_{ij} - \sum_{k=1}^{K} y_{ij}^k \right). \tag{12.87}$$

The projection parametern $\beta$, which is used when projecting onto the interval GUB constraints, is now calculated by

$$\beta = \text{mid}\ (\psi_{ij}, \Gamma, \Delta),$$

where $\psi_{ij}$ is the dual price for the GUB constraint corresponding to the $(i,j)$th arc. The iterative step for an interval of the GUB constraints takes the form:

$$(y_{ij}^k)' = y_{ij}^k + \frac{\beta}{w_{ij}^k}, \tag{12.88}$$

$$\psi_{ij} \leftarrow \psi_{ij} - \beta. \tag{12.89}$$

The complete algorithm is summarized below. The dual prices of origin and destination nodes are denoted by $\pi^{O,k}$ and $\pi^{D,k}$ respectively, for each commodity $k$. The reduced cost for the bound on the flow of the $k$th commodity on arc $(i,j)$ is denoted by $r_{ij}^k$.

### Algorithm 12.4.3 Row-action Algorithm for Quadratic Multicommodity Transportation Problems

**Step 0:** (Initialization.) $\nu = 0$. Set the initial dual prices $(\pi^{O,k})^0 = 0$, $(\pi^{D,k})^0 = 0$, $(r^k)^0 = 0$, for $k = 1, 2, \ldots, K$, and $\psi^0 = 0$, and get an initial primal vector that satisfies the initialization step of Algorithm 6.4.1 for the present problem with the quadratic function (12.83), i.e.,

$$(x_{ij}^k)^0 = -\frac{c_{ij}^k}{w_{ij}^k},\ k = 1, 2, \ldots, K,\ (i,j) \in \mathcal{A}. \tag{12.90}$$

**Step 1:** (Solve the single commodity problems.) For $k = 1, 2, \ldots, K$ :

**Step 1.1:** (Iterative step for origin nodes.) For $i = 1, 2, \ldots, m_O$, calculate:

$$\beta_\nu = \frac{1}{\sum_{j \in \delta_i^+} 1/w_{ij}^k}\left(s_i^k - \sum_{j \in \delta_i^+} (x_{ij}^k)^\nu\right), \tag{12.91}$$

$$(x_{ij}^k)^\nu \leftarrow (x_{ij}^k)^\nu + \frac{\beta_\nu}{w_{ij}^k},\ \text{for all } j \in \delta_i^+, \tag{12.92}$$

$$(\pi^{O,k})_i^{\nu+1} = (\pi^{O,k})_i^\nu - \beta_\nu. \tag{12.93}$$

**Step 1.2:** (Iterative step for destination nodes.) For $j = 1, 2, \ldots, m_D$, calculate:

$$\beta_\nu = \frac{1}{\sum_{i \in \delta_j^-} 1/w_{ij}^k}\left(d_j^k - \sum_{i \in \delta_j^-} (x_{ij}^k)^\nu\right), \tag{12.94}$$

$$(x_{ij}^k)^\nu \leftarrow (x_{ij}^k)^\nu + \frac{\beta_\nu}{w_{ij}^k},\ \text{for all } i \in \delta_j^-, \tag{12.95}$$

$$(\pi^{D,k})_j^{\nu+1} = (\pi^{D,k})_j^\nu - \beta_\nu. \tag{12.96}$$

**Step 1.3:** (Iterative step for bounds on flows.) For all $(i,j) \in \mathcal{A}$, calculate:

$$\beta_\nu = \text{mid}\left((r_{ij}^k)^\nu,\ w_{ij}^k(u_{ij}^k - (x_{ij}^k)^\nu),\ -w_{ij}^k(x_{ij}^k)^\nu\right), \tag{12.97}$$

$$(x_{ij}^k)^\nu \leftarrow (x_{ij}^k)^\nu + \frac{\beta_\nu}{w_{ij}^k}, \tag{12.98}$$

$$(r_{ij}^k)^{\nu+1} = (r_{ij}^k)^\nu - \beta_\nu. \tag{12.99}$$

**Step 2:** (Iterative step for generalized upper bound GUB constraints.) For all $(i,j) \in \mathcal{A}$, calculate:

$$\Gamma_\nu = -\frac{1}{\sum_{k=1}^K 1/w_{ij}^k} \sum_{k=1}^K (x_{ij}^k)^\nu, \tag{12.100}$$

$$\Delta_\nu = \frac{1}{\sum_{k=1}^K 1/w_{ij}^k} (U_{ij} - \sum_{k=1}^K (x_{ij}^k)^\nu), \tag{12.101}$$

$$\beta_\nu = \text{mid}\ (\psi_{ij}^\nu, \Gamma_\nu, \Delta_\nu), \tag{12.102}$$

$$(x_{ij}^k)^{\nu+1} = (x_{ij}^k)^\nu + \frac{\beta_\nu}{w_{ij}^k},\ k = 1, 2, \ldots, K, \tag{12.103}$$

$$\psi_{ij}^{\nu+1} = \psi_{ij}^\nu - \beta_\nu. \tag{12.104}$$

**Step 3:** Let $\nu \leftarrow \nu + 1$ and return to Step 1.

### Decompositions for Parallel Computing

The calculations in Steps 1.1 – 1.3 of the algorithm are structurally identical to the calculations of Algorithms 12.4.1 and 12.4.2 for transportation problems. Hence, the ideas discussed in Section 12.4.1 for parallel computation are applicable to the multicommodity network flow problem as well. Furthermore, the iteration in Step 1 over the commodity index $k$ can proceed concurrently for multiple commodities. Hence, Step 1.1 can be executed concurrently for all origin nodes and for all commodities utilizing $Km_O$ processors in calculating the projection parameters and updating dual variables, and utilizing $Kn$ processors for the flow updates. Step 1.2 parallelizes in a similar fashion. Step 1.3 can be executed concurrently for all arcs and for all commodities, utilizing $Kn$ processors.

Step 2 can be executed concurrently for all arcs. For each arc the algorithm computes the total flow of all commodities, and computes the parameters $\Gamma, \Delta$, and $\beta$ utilizing $n$ processors. The flows can be subsequently updated concurrently for all arcs and all commodities, using $Kn$ processors. The reduced cost vector $\psi$ can be updated concurrently for all arcs, again utilizing $n$ processors.

## 12.5 A Model Decomposition Algorithm for Multicommodity Network Flow Problems

We develop now a model decomposition algorithm for the multicommodity network flow problem, which is based on the linear-quadratic penalty (LQP) function (see Chapter 7). It decomposes the model described in (12.18)–(12.21) by removing the inequalities $Gx \leq U$ from the feasible set and appending them to the objective via a penalty function (see Chapter 4).

To describe the algorithm we distinguish between the set of solutions that satisfy only the flow conservation constraints,

$$X \doteq \{x \mid Ax = b,\ 0 \leq x \leq u\},$$

and the set of solutions that satisfy the complicating constraints,

$$\Omega \doteq \{x \mid\ Gx \leq U\}.$$

The constraint set $X$, for problem (12.18)–(12.21) is a Cartesian product. Using the notation introduced in Section 12.2.2 we write $X = \prod_{k=1}^{K} X_k$, where $X_k \doteq \{x^k \mid A_k x^k = b^k,\ 0 \leq x^k \leq u^k\}$.

### 12.5.1 The linear-quadratic penalty (LQP) algorithm

The complicating GUB constraints $Gx \leq U$ are removed from the constraint set, and a penalty function is added to the objective to penalize violations of these constraints. The penalty function is of the linear-quadratic form

$$\phi_\epsilon(t) \doteq \begin{cases} 0, & \text{if } t < 0, \\ \dfrac{t^2}{2\epsilon}, & \text{if } 0 \leq t \leq \epsilon, \\ (t - \dfrac{\epsilon}{2}), & \text{if } t > \epsilon, \end{cases} \tag{12.105}$$

where $\epsilon > 0$ is a user-specified parameter. In lieu of solving (12.18)–(12.21) consider instead the following smoothed penalty version

$$\text{Minimize} \quad \sum_{k=1}^{K} F_k(x^k) + \mu \sum_{l=1}^{n} \phi_\epsilon(t_l) \tag{12.106}$$

$$\text{s.t.} \quad x \in X. \tag{12.107}$$

where $t_l \doteq (Gx - U)_l$ for $l = 1, 2, \ldots, n$, and $\mu > 0$ is a *penalty* parameter. The LQP algorithm solves a sequence of problems of this form, using increasing values of the penalty parameter $\mu$ and decreasing values of the approximation error $\epsilon$ (see Section 7.2). Convergence of the algorithm follows from Theorem 4.1.1, since condition (*ii*) is satisfied. The algorithm is summarized as follows.

**Algorithm 12.5.1 The LQP Algorithm for Multicommodity Network Flow Problems**

**Step 0:** (Initialization.) Set $\nu = 0$. Choose parameters $\mu_0 > 0$, $\epsilon_0 > 0$.

**Step 1:** Using the violations $t_l = (Gx - U)_l$ for all $l = 1, 2, \ldots, n$, for evaluating the penalty function, solve the problem:

$$\text{Minimize} \quad \sum_{k=1}^{K} F_k(x^k) + \mu_\nu \sum_{l=1}^{n} \phi_{\epsilon_\nu}(t_l) \tag{12.108}$$

$$\text{s.t.} \quad \boldsymbol{x} \in X. \tag{12.109}$$

Let $(\boldsymbol{x}^*)^\nu$ denote the optimal solution.

**Step 2:** If $(\boldsymbol{x}^*)^\nu$ satisfies some termination criteria then stop. Otherwise update the parameters $\mu_\nu$ and $\epsilon_\nu$, set $\nu \leftarrow \nu + 1$, and go back to Step 1.

(Details on this general algorithm, such as procedures for updating the penalty parameters, are given in Section 7.2.)

A starting point for initializing any iterative algorithm for solving the minimization problem in Step 1 can be obtained by solving the following network problems for all $k = 1, 2, \ldots, K$,

$$\text{Minimize} \quad F_k(x^k) \tag{12.110}$$

$$\text{s.t.} \quad A_k x^k = b^k, \tag{12.111}$$

$$0 \le x^k \le u^k. \tag{12.112}$$

If $(x^{k*})$, $k = 1, 2, \ldots, K$ denote the optimal solutions, then the concatenated $\boldsymbol{x}^0 = ((x^{1*})^{\mathrm{T}}, \ldots, (x^{K*})^{\mathrm{T}})^{\mathrm{T}}$ is a suitable starting point. Of course, if this point satisfies the complicating GUB constraints, i.e., $\boldsymbol{x}^0 \in \Omega$, then the original problem has been solved and there is no need to proceed further with the LQP algorithm.

Most of the computational effort of the algorithm is invested in solving the smoothed penalty problem in Step 1. The objective function is continuously differentiable, and the feasible set is the Cartesian product set $X = \prod_{k=1}^{K} X_k$. Each set $X_k$ corresponds to network flow constraints. The smoothed penalty problem can be solved using specialized implementations of nonlinear programming algorithms for network problems, such as the truncated Newton method. However, the penalty function $\phi_{\epsilon_\nu}(t_l)$ is not block-separable in the $x^k$ (recall that $t_l = (Gx - U)_l$). The fact that $X$ is a Cartesian product motivates the use of an algorithm for executing Step 1 based on a linearization of the penalty function. In this way the objective function becomes block-separable in the $x^k$'s, and the nonlinear program (12.108)–(12.109) decomposes into $K$ independent linear network problems. The simplicial decomposition algorithm (Algorithm 7.1.3) can be used to linearize the smoothed penalty function, and hence decompose the problem by commodity.

### Decompositions for Parallel Computing

The use of a linearization algorithm (such as simplicial decomposition) to execute Step 1 of the LQP algorithm results in a decomposition into $K$ linear network subproblems, which can be solved in parallel. Simplicial decomposition, however, has a coordinating master program that uses the solutions of the linear subproblems to obtain an estimate of the minimum of the nonlinear penalty function. As discussed earlier (Sections 7.1 and 7.2), this master program is also amenable to parallel computations, and takes only a small fraction of the execution time of simplicial decomposition. Finally, updating the penalty parameters in Step 2 of the LQP algorithm involves simple scalar calculations and is computationally insignificant. The combination of the LQP algorithm with simplicial decomposition for solving multicommodity network flow problems is a procedure that parallelizes very well, either when the number of blocks $K$ is large, or when each block is large. Numerical experiences with this algorithm are reported in Section 15.4.

## 12.6 Notes and References

The seminal papers on network flows are due to Kantorovich (1939, 1942), Hitchcock (1941), and Koopmans (1949) wherein the linear transportation problem was first formulated for production scheduling and transportation applications. Dantzig (1951) developed specializations of the simplex algorithm for the transportation problem. See also Dantzig (1963, Chapter 14). Ford and Fulkerson (1962) developed primal-dual algorithms that exploit combinatorial properties of the problem. Their book fostered the rapid growth of the study of network problems. More recent books that treat theory and algorithms in network optimization are Christofides (1975), Kennington and Helgason (1980), Jensen and Barnes (1980), Tarjan (1983), Rockafellar (1984), Bertsekas (1991), Murty (1992), and Ahuja, Magnanti, and Orlin (1993). The last two books also treat applications extensively.

The visual representations of the graphs of network problems is an effective way to communicate mathematical programs to users. Glover, Klingman, and Phillips (1990) develop network-related formulations, which they call *netforms.* A textbook by the same authors (1992) describes exclusively model formulations and applications. Sheffi (1985) examines network flow problems in transportation and discusses models and algorithms for their solution.

For a survey of the state-of-the-art in network optimization see Ahuja, Magnanti, and Orlin (1989), and for a survey of research in parallel computing for network optimization see Bertsekas et al. (1995). Survey articles for the multicommodity network flow problem are found in Assad (1978) and Kennington (1978). A survey of algorithms and models for network optimization problems with nonlinear objective functions is given in Dembo, Mulvey and Zenios (1989). Pinar and Zenios (1993) survey and compare parallel algorithms for multicommodity flow problems.

The first parallel algorithm for a network optimization problem was developed by Bertsekas (1979), and was further refined and tested by Bertsekas (1985, 1988), Phillips and Zenios (1989), Wein and Zenios (1991), Kempka, Kennington, and Zaki (1991), and Zaki (1995). It is known as the *auction algorithm*—due to its dual price adjustment mechanism, which is similar to that of an auction—and is applicable to assignment problems.

The first implementation of an optimization algorithm on a massively parallel computer was that of the relaxation algorithm for strictly convex network problems by Zenios and Lasken (1988a) on the Connection Machine CM–1. They developed the data-structures for representing network problems on SIMD architectures. Further experiments using the Connection Machine CM–2 were reported in Zenios and Lasken (1988b).

**12.1–12.2** The notation used in this chapter is fairly standard in the literature, and is based on notation introduced in Ford and Fulkerson (1962) and developed further in Kennington and Helgason (1980).

**12.3** The application of linear programming to determine the minimum margin requirements for covering positions in options is suggested in Rudd and Schroeder (1982); see also Cox and Rubinstein (1985, pp. 107–109). The network formulation is adopted from Dahl, Meeraus and Zenios (1993). The use of network models for air-traffic control has been discussed by several authors; see Ferguson and Dantzig (1957), Mulvey and Zenios (1987a), Booth and Harvey (1983), Zenios (1991a), Bielli et al. (1969), Odoni (1987), and Richetta (1991). Bertocchi and Zenios (1993) discuss the use of parallel computers for problems in air-traffic control. See Del Balzo (1989) for a general discussion on the air-traffic system. Routing applications and traffic assignment problems are discussed in Sheffi (1985).

**12.4.1** The row-action algorithms for nonlinear transportation problems were developed by Zenios and Censor (1991a). They also discuss the solution of pure and generalized transportation problems with quadratic and entropic objective functions, and develop implementations of the algorithms on massively parallel SIMD architectures. The introduction of underrelaxation parameters in Algorithm 6.4.1, which allows us to use them here, was made by Iusem and Zenios (1995). The use of a secant approximation for the estimation of the generalized projection parameter was studied by Censor et al. (1990). As discussed earlier (Section 6.9) they established that the asymptotic convergence of the row-action algorithm is preserved even with the use of an approximation for the projection parameter. Related primal-dual algorithms for the transportation problem were developed by Ohuchi and Kaji (1984), Zenios and Mulvey (1986b), and Bertsekas, Hossein, and Tseng (1987). A similar algorithm for the quadratic transportation problem was developed by Bachem and Korte (1978). These algorithms are known as *relaxation algorithms*, and they use a dual formulation of the network optimization problem that differs slightly from the one adopted by the row-action algorithms. In particular, the row-action algorithms dualize both the flow conservation constraints and the variable bounds, while the relaxation methods dualize only the flow conservation constraints and explicitly bound the primal variables in the calculation of the Lagrangian function. A discussion of the differences between row-action methods and relaxation algorithms is given by Bertsekas et al. (1995). A primal Newton type algorithm for nonlinear transportation problems, which is quite distinct from the iterative algorithms discussed here, is due to Klincewicz (1989). A survey of algorithms for nonlinear network optimization problems is given in Dembo, Mulvey, and Zenios (1989).

**12.4.3** The row-action algorithm for the multicommodity network flow problem was developed by Zenios (1991b). The same reference also discusses the solution of pure and generalized multicommodity network flow problems with quadratic and entropic objective functions, and develops data-structures for the implementation of the algorithms on massively parallel SIMD architectures. A related algorithm was developed by Nagamochi, Fukushima, and Ibaraki (1990), who adopted the dual formulation of the relaxation algorithms (see the notes and references discussion for subsection 12.4.1 above). The row-action algorithm for the multicommodity network flow problem was further studied by Censor, Chajakis, and Zenios (1995), who developed alternative parallel implementations on shared-memory multiprocessors.

**12.5** The linear-quadratic penalty (LQP) algorithm for problems with network structures was proposed by Zenios, Pinar, and Dembo (1994), and was applied to solve multicommodity network flow problems by Pinar and Zenios (1992). Data-structures for a massively parallel implementation for the solution of large-scale problems were developed by Pinar and Zenios (1994a). A comparative study of algorithms, including the algorithm of this section, is found in Pinar and Zenios (1993). For references on simplicial decomposition see Section 7.3.

Truncated Newton algorithms for solving large-scale nonlinear programs with network constraints, such as the one appearing in Step 1 of the LQP algorithm, have been developed by Escudero (1986), Dembo (1987), Ahlfeld et al. (1987), and Zenios and Pinar (1992).

Other authors have also suggested the decomposition of multicommodity network flow problems using a combination of penalty methods with linearization techniques. Schultz and Meyer (1991) developed an interior point method using a shifted barrier function for structured optimization problems. They applied their algorithm to the solution of large-scale multicommodity network flow problems using parallel computing, with encouraging computational results. Mulvey and Ruszczyński (1992), based on earlier work by Stephanopoulos and Westerberg (1975) and Ruszczyński (1989), used an augmented Lagrangian formulation followed by a diagonalization of the augmented Lagrangian using a quadratic approximation. They applied this algorithm to the solution of structured nonlinear programs arising in stochastic programming. Encouraging computational results for the solution of large-scale problems in a distributed computing environment of high-performance workstations were obtained by Berger, Mulvey, and Ruszczyński (1994).

The idea of linearizing or diagonalizing a nonseparable function is extensively used in solving *traffic assignment* problems, whereby multiple distinct commodities are linked via a nonseparable delay function and rather than through generalized upper bounding constraints. The book by Sheffi (1985) discusses several such algorithms, including the well-known Frank-Wolfe algorithm. The linearization of traffic assignment problems to induce separability and use parallel computations was first explored by Chen and Meyer (1988). It is also discussed in Bertsekas and Tsitsiklis (1989, pp. 414–422).

CHAPTER 13

# Planning Under Uncertainty

> The winners of tomorrow will deal *proactively* with chaos, will look at the chaos per se as the source for market advantage, not as a problem to be got around. Chaos and uncertainty are (will be) market opportunities for the wise.
>
> Tom Peters, *Thriving on Chaos*, 1987

The analyst who attempts to build a mathematical model for a real-world system is often faced with the problem of uncertain, noisy, incomplete or erroneous data. This is true for several application domains. In business applications noisy data are prevalent. Returns of financial instruments, demand for a firm's products, the cost of fuel, and consumption of power and other resources are examples of model data that are known with some probabilistic distribution at best. In social sciences data are often incomplete—for example, partial census surveys are carried out periodically in lieu of a complete census of the population. In the physical sciences and engineering data are usually subject to measurement errors, as in models of image restoration from remote sensing experiments.

For some applications not much is lost by assuming that the value of the uncertain data is known and then developing a deterministic mathematical programming model. Worst case or mean values can be used in this respect because they provide reasonable approximations when either the level of uncertainty is low, or when the uncertain parameters have a minor impact on the system we want to model. For many applications, however, uncertainty plays a key role in the performance of the real-world system: worst case analysis often leads to conservative and potentially expensive solutions, and solving the *mean value problem*, i.e., a problem where all random variables are replaced by their mean values can even lead to nonsensical solutions since the mean of a random variable might not be a value that can be realized in practice.

A general approach to dealing with uncertainty is to assign to the unknown parameters a probability distribution, which should then be incorporated into an appropriate mathematical programming model. This chapter addresses the problem of planning under uncertainty and develops mathematical programming formulations using *stochastic linear programming* *models* and *robust optimization models.* Section 13.2 discusses a classic example to highlight the issues involved. Sections 13.3 and 13.4 give mathematical programming models. Section 13.5 discusses diverse real-world applications and 13.6 focuses on an application from financial planning. Section 13.7 introduces a broad class of stochastic programming applications, namely that of stochastic network problems. Section 13.8 gives a row-action iterative algorithm for solving this special class of problems. Notes and references in Section 13.9 conclude this chapter.

## 13.1 Preliminaries

We introduce first some basic definitions on probability spaces that are needed throughout this section. Additional background material on probability theory can be found, e.g., in Billingsley (1995) and, with emphasis on stochastic programming, in Kall (1976) and Wets (1989). In this chapter boldface Greek characters are used to denote random vectors which belong to some probability space as defined below.

Let $\Omega$ be an arbitrary space or set of points. A *$\sigma$-field* for $\Omega$ is a family $\Sigma$ of subsets of $\Omega$ such that $\Omega$ itself, the complement with respect to $\Omega$ of any set in $\Sigma$, and any union of countably many sets in $\Sigma$ are all in $\Sigma$. The members of $\Sigma$ are called *measurable* sets, or *events* in the language of probability theory. The set $\Omega$ with the $\sigma$-field $\Sigma$ is called a *measurable space* and is denoted by $(\Omega, \Sigma)$.

Let $\Omega$ be a (linear) vector space and $\Sigma$ a $\sigma$-field. A *probability measure* $P$ on $(\Omega, \Sigma)$ is a real-valued function defined over the family $\Sigma$, which satisfies the following conditions: (*i*) $0 \leq P(A) \leq 1$ for $A \in \Sigma$; (*ii*) $P(\emptyset) = 0$ and $P(\Omega) = 1$; and (*iii*) if $\{A_k\}_{k=1}^{\infty}$ is a sequence of disjoint sets $A_k \in \Sigma$ and if $\cup_{k=1}^{\infty} A_k \in \Sigma$ then $P(\cup_{k=1}^{\infty} A_k) = \sum_{k=1}^{\infty} P(A_k)$. The triplet $(\Omega, \Sigma, P)$ is called a *probability space.* The *support* of $(\Omega, \Sigma, P)$ is the smallest subset of $\Omega$ with probability 1. If the support is a countable set then the probability measure is said to be discrete. The term *scenario* is used for the elements of $\Omega$ of a probability space with a discrete distribution.

A proposition is said to hold *almost surely* (abbreviated a.s.) or *$P$-almost surely* if it holds on a subset $A \subseteq \Omega$ with $P(A) = 1$. The *expected value* of a random variable $\mathcal{Q}$ on $(\Omega, \Sigma, P)$ is the Stieltjes integral of $\mathcal{Q}$ with respect to the measure $P$:

$$E[\mathcal{Q}] \doteq \int \mathcal{Q} dP = \int_{\Omega} \mathcal{Q}(\omega) dP(\omega).$$

The expectation of a constant function is also constant and it is easy to see that $E[\mathcal{Q}_1 + \mathcal{Q}_2] = E[\mathcal{Q}_1] + E[\mathcal{Q}_2]$. The $k$th *moment* of $\mathcal{Q}$ is the expected value of $\mathcal{Q}^k$, i.e., $E[\mathcal{Q}^k] \doteq \int_{\Omega} \mathcal{Q}^k(\omega) dP(\omega)$. The *variance* of the random variable $\mathcal{Q}$ is defined as $\mathrm{Var}[\mathcal{Q}] \doteq E[\mathcal{Q}^2] - (E[\mathcal{Q}])^2$.

Finally, we give a formal but restricted (to our needs) definition of a conditional expectation. Let $(\Omega, \Sigma, P)$ be a probability space and suppose that $\mathcal{A}_1, \mathcal{A}_2, \ldots, \mathcal{A}_K$ is a finite partition of the set $\Omega$. From this partition we form a $\sigma$-field $\mathcal{A}$ which is a subfield of $\Sigma$. Then the *conditional expectation* of the random variable $\mathcal{Q}(\boldsymbol{\omega})$ on $(\Omega, \Sigma, P)$ given $\mathcal{A}$ at $\boldsymbol{\omega}$ is denoted by $E[\mathcal{Q} \mid \mathcal{A}]$ and defined as

$$E[\mathcal{Q} \mid \mathcal{A}] \doteq \frac{1}{P(\mathcal{A}_i)} \int_{\mathcal{A}_i} \mathcal{Q}(\boldsymbol{\omega}) dP(\boldsymbol{\omega})$$

for $\boldsymbol{\omega} \in A_i$, assuming that $P(\mathcal{A}_i) > 0$.

## 13.2 The Newsboy Problem

To develop an understanding of stochastic programming problems we consider first the following simple problem of planning under uncertainty. On a street corner a young entrepreneur is selling newspapers that he buys from a local distributor each morning. He sells these papers for a profit $p^+$ per unit, and any papers that remain at the end of the day are sold as scrap paper, in which case a net loss $p^-$ is realized per unit. The demand for newspapers is a random variable $\boldsymbol{\omega}$ which belongs to a probability space with support denoted by $\Omega \doteq \{\boldsymbol{\omega} \in \mathbb{R} \mid 0 \leq \boldsymbol{\omega} < \infty\}$ and probability distribution function $P(\boldsymbol{\omega})$. The problem is to choose the *optimal* number of papers $x$ that should be bought from the local distributor.

An approach to modeling this situation is to consider a *policy* $x$ as optimal if it maximizes the *expected* profit. Profit is a function of the policy and the demand random variable $\boldsymbol{\omega}$. Let $F(x, \boldsymbol{\omega})$ be the profit function:

$$F(x, \boldsymbol{\omega}) \doteq \begin{cases} p^+ x & \text{if } x \leq \boldsymbol{\omega}, \\ p^+ \boldsymbol{\omega} - p^-(x - \boldsymbol{\omega}) & \text{if } x > \boldsymbol{\omega}. \end{cases}$$

The expected value of the profit function is the Stieltjes integral with respect to the distribution function:

$$\begin{aligned} E\left[F(x, \boldsymbol{\omega})\right] &= \int_{\Omega} F(x, \boldsymbol{\omega}) dP(\boldsymbol{\omega}) \\ &= \int_0^x \left(p^+ \boldsymbol{\omega} - p^-(x - \boldsymbol{\omega})\right) dP(\boldsymbol{\omega}) + \int_x^\infty p^+ x dP(\boldsymbol{\omega}), \end{aligned}$$

and the mathematical model for the newsboy problem is the following optimization problem with respect to $x$,

$$\text{Maximize} \quad E\left[F(x, \boldsymbol{\omega})\right] \tag{13.1}$$

$$\text{s.t.} \quad x \geq 0. \tag{13.2}$$

This is a simple example of a problem of planning under uncertainty. It is represented by an *adaptive* model, since decisions adapt as more information becomes available, i.e., as newspapers are sold to customers that arrive during the day. The model has *fixed recourse*, meaning that the reaction to the observed demand is fixed. That is, the number of newspapers sold for a profit is uniquely determined by the number of customers. The same is true for the surplus created at the end of the day, which is sold for scrap at a loss. Other forms of recourse action might have been possible, such as purchasing additional papers at a higher cost later during the day, or returning newspapers before the end of the day, at a value higher than that of scrap paper. This simple, fixed recourse model does not allow for such considerations, and it also assumes that all risk preferences are captured by the expected value of the profit. Higher moments of the distribution of the profit function $F(x, \boldsymbol{\omega})$ are ignored. The next section presents mathematical models for planning under uncertainty in more complicated settings.

## 13.3 Stochastic Programming Problems

The following problem is the general formulation of *stochastic programming*:

$$\begin{aligned} \text{Minimize} \quad & E\left[f_0(x, \boldsymbol{\omega})\right] \\ \text{s.t.} \quad & E\left[f_i(x, \boldsymbol{\omega})\right] = 0, \qquad i = 1, 2, \ldots, m, \\ & x \in X \subset \mathbb{R}^n. \end{aligned} \tag{13.3}$$

The following notation is used: $\boldsymbol{\omega}$ is a random vector with support $\Omega \subset \mathbb{R}^N$ and $P \doteq P(\boldsymbol{\omega})$ is a probability distribution function on $\mathbb{R}^N$. Also $f_0 : \mathbb{R}^n \times \Omega \rightarrow \mathbb{R} \cup \{+\infty\}$, $f_i : \mathbb{R}^n \times \Omega \rightarrow \mathbb{R}$, $i = 1, 2, \ldots, m$, and $X$ is a closed set. Inequality constraints can be incorporated into this formulation with the use of slack variables.

The expectation functions

$$E[f_i(x, \boldsymbol{\omega})] \doteq \int_\Omega f_i(x, \boldsymbol{\omega}) dP(\boldsymbol{\omega}), \tag{13.4}$$

are assumed finite for all $i = 0, 1, \ldots, m$ unless the set $\{\boldsymbol{\omega} \mid f_0(x, \boldsymbol{\omega}) = +\infty\}$ has a nonzero probability, in which case $E[f_0(x, \boldsymbol{\omega})] = +\infty$. The feasibility set

$$X \cap \{x \mid E[f_i(x, \boldsymbol{\omega})] = 0,\ i = 1, 2, \ldots, m\} \cap \{x \mid E[f_0(x, \boldsymbol{\omega})] < +\infty\}$$

is assumed to be nonempty.

The model (13.3) is a nonlinear programming problem, whose constraints and objective functions are represented by integrals. Much of the theory of stochastic programming is concerned with identifying the properties of these integral functions and devising suitable approximation schemes for their evaluation. Optimality conditions are derived from those for nonlinear programming, with the aid of subdifferential calculus for the expectation functions. However, the computation of solutions for these nonlinear programs poses serious challenges, since evaluation of the integrals can be an extremely difficult task, especially when the expectation functionals are multidimensional. There are even cases when the integrands are neither differentiable, nor convex nor even continuous. A broad class of stochastic programming models, however, can be formulated as large-scale linear or nonlinear programs with a specially structured constraints matrix. Most of the work on parallel computing for stochastic programming focuses on the development of decomposition algorithms that exploit this special structure. In the next subsections we look at further refinements of the general stochastic programming formulation.

### 13.3.1 Anticipative models

Consider now the following situation. A decision $x$ must be made in an uncertain world where the uncertainty is described by the random vector $\boldsymbol{\omega}$. The decision does not in any way depend on future observations, but prudent planning has to anticipate possible future realizations of the random vector.

In anticipative models feasibility is expressed in terms of *probabilistic* (or *chance*) *constraints*. For example, a *reliability* level $\alpha$ where $0 < \alpha \leq 1$ is specified and constraints are expressed in the form

$$P\{\boldsymbol{\omega} \mid g_i(x, \boldsymbol{\omega}) = 0, \;\; i = 1, 2, \ldots, m\} \geq \alpha,$$

where $g_i : \mathbb{R}^n \times \Omega \rightarrow \mathbb{R}, \;\; i = 1, 2, \ldots, m$. This constraint can be cast in the form of the general model (13.3) by defining $f_i$ as follows

$$f_i(x, \boldsymbol{\omega}) \doteq \begin{cases} \alpha - 1 & \text{if } g_i(x, \boldsymbol{\omega}) = 0, \\ \alpha & \text{otherwise.} \end{cases}$$

The objective function may also be of a reliability type, such as $P\{\boldsymbol{\omega} \mid g_0(x, \boldsymbol{\omega}) \leq \gamma\}$, where $\gamma$ is a constant.

In summary, an anticipative model selects a policy that leads to some desirable characteristics of the constraint and objective functionals under the realizations of the random vector. In the example above it is desirable that the probability of a constraint violation is less than the prespecified threshold value $\alpha$. The precise value of $\alpha$ depends on the application at hand, the cost of constraint violation, and other similar considerations.

### 13.3.2 Adaptive models

In an adaptive model, observations related to uncertainty become available before a decision $x$ is made, such that optimization takes place in a learning environment. It is understood that observations provide only partial information about the random variables because otherwise the model would simply wait to observe the values of the random variables, and then make a decision $x$ by solving a deterministic mathematical program. In contrast to this approach we have the other extreme situation where all observations are made after the decision $x$ has been made, and the model becomes anticipative.

Let $\mathcal{A}$ be the collection of all the relevant information that could become available by making an observation. This $\mathcal{A}$ is a subfield of the $\sigma$-field of all possible events, generated from the support set $\Omega$ of the random vector $\boldsymbol{\omega}$. The decisions $x$ depend on the events that could be observed, and $x$ is termed *$\mathcal{A}$-adapted* or *$\mathcal{A}$-measurable.* Using the conditional expectation with respect to $\mathcal{A}$, $E[\,\cdot \mid \mathcal{A}]$, the adaptive stochastic program can be written as:

$$
\begin{aligned}
\text{Minimize} \quad & E[f_0(x(\boldsymbol{\omega}), \boldsymbol{\omega}) \mid \mathcal{A}] \\
\text{s.t.} \quad & E[f_i(x(\boldsymbol{\omega}), \boldsymbol{\omega}) \mid \mathcal{A}] = 0, \qquad i = 1, 2, \ldots, m, \\
& x(\boldsymbol{\omega}) \in X, \qquad \text{almost surely.}
\end{aligned} \tag{13.5}
$$

where the mapping $x : \Omega \to X$ is such that $x(\boldsymbol{\omega})$ is $\mathcal{A}$-measurable. This problem can be addressed by solving for every $\boldsymbol{\omega}$ the following deterministic programs

$$\text{Minimize} \quad E[f_0(x, \cdot) \mid \mathcal{A}](\boldsymbol{\omega}) \tag{13.6}$$

$$\text{s.t.} \quad E[f_i(x, \cdot) \mid \mathcal{A}](\boldsymbol{\omega}) = 0, \qquad i = 1, 2, \ldots, m, \tag{13.7}$$

$$x \in X. \tag{13.8}$$

Each such problem for a given $\boldsymbol{\omega}$ is of the canonical form (13.3).

The two extreme cases (i.e., complete information with $\mathcal{A} = \Sigma$, or no information at all) deserve special mention. The case of no information reduces the model to the form of the anticipative model; when there is complete information the model (13.5) is known as the *distribution model.* The goal in this later case is to characterize the distribution of the optimal objective function value. The precise values of the objective function and the optimal policy $x$ are determined after realizations of the random vector $\boldsymbol{\omega}$ are observed. The most interesting situations arise when partial information becomes available after some decisions have been made, and models to address such situations are discussed next.

### 13.3.3 Recourse models

The *recourse problem* combines the anticipative and adaptive models in a common mathematical framework. The problem seeks a policy that not only anticipates future observations but also takes into account that obser-vations are made about uncertainty, and thus can adapt by taking *recourse* decisions. For example, a portfolio manager specifies the composition of a portfolio considering both future movements of stock prices (anticipation) and that the portfolio will be rebalanced as prices change (adaptation).

The two-stage version of this model has been studied extensively. It is amenable to formulations as a large-scale deterministic nonlinear program with a special structure of the constraints matrix. These formulations yield naturally to solution via decomposition algorithms and parallel computations. To formulate the two-stage stochastic program with recourse we need two vectors for decision variables to distinguish between the anticipative policy and the adaptive policy. The following notation is used.

$x \in \mathbb{R}^{n_0}$ denotes the vector of *first-stage* decisions. These decisions are made before the random variables are observed and are anticipative.

$y \in \mathbb{R}^{n_1}$ denotes the vector of *second-stage* decisions. These decisions are made after the random variables have been observed and are adaptive. They are constrained by decisions made at the first-stage, and depend on the realization of the random vector.

We formulate the *second-stage problem* in the following manner. Once a first-stage decision $x$ has been made, some realization of the random vector can be observed. Let $q(y, \boldsymbol{\omega})$ denote the second-stage cost function, and let $\{T(\boldsymbol{\omega}), W(\boldsymbol{\omega}), h(\boldsymbol{\omega}) \mid \boldsymbol{\omega} \in \Omega\}$ be the model parameters. Those parameters are functions of the random vector $\boldsymbol{\omega}$ and are, therefore, random parameters. $T$ is the *technology matrix* of dimension $m_1 \times n_0$. It contains the technology coefficients that convert the first-stage decision $x$ into resources for the second-stage problem. $W$ is the *recourse matrix* of dimension $m_1 \times n_1$. $h$ is the second-stage resource vector of dimension $m_1$.

The second-stage problem seeks a policy $y$ that optimizes the cost of the second-stage decision for a given value of the first-stage decision $x$. We denote the optimal value of the second-stage problem by $\mathcal{Q}(x, \boldsymbol{\omega})$. This value depends on the random parameters and on the value of the first-stage variables $x$. $\mathcal{Q}(x, \boldsymbol{\omega})$ is the optimal value, for any given $\Omega$, of the following nonlinear program

$$\begin{aligned} \text{Minimize} \quad & q(y, \boldsymbol{\omega}) \\ \text{s.t.} \quad & W(\boldsymbol{\omega})y = h(\boldsymbol{\omega}) - T(\boldsymbol{\omega})x, \\ & y \in \mathbb{R}^{n_1}_{+}. \end{aligned} \tag{13.9}$$

If this second-stage problem is infeasible then we set $\mathcal{Q}(x, \boldsymbol{\omega}) \doteq +\infty$. The model (13.9) is an *adaptation* model in which $y$ is the *recourse* decision and $\mathcal{Q}(x, \boldsymbol{\omega})$ is the *recourse cost function.*

The two-stage stochastic program with recourse is an optimization problem in the first-stage variables $x$, which optimizes the sum of the cost of the first-stage decisions, $f(x)$, and the expected cost of the second-stage decisions. It is written as follows.

$$\begin{array}{lll} \text{Minimize} & f(x) + E[\mathcal{Q}(x,\boldsymbol{\omega})] & \\ \text{s.t.} & \langle a^i, x\rangle = b_i, \qquad i = 1,2,\ldots,m_0, & \\ & x \in \mathbb{R}^{n_0}_+, & \end{array} \tag{13.10}$$

where $a^i$ denotes the transpose of the $i$th row of the $m_0 \times n_0$ matrix $A$, and $b_i$ is the $i$th component of the $m_0$-vector $b$. $\langle a^i, x\rangle = b_i$, $i = 1,2,\ldots,m_0$ are linear restrictions on the first-stage variables. This model can be cast in the general formulation (13.3) simply by denoting $f_0(x,\boldsymbol{\omega}) \doteq f(x) + \mathcal{Q}(x,\boldsymbol{\omega})$, and $f_i(x,\boldsymbol{\omega}) \doteq \langle a^i, x\rangle - b_i$.

A formulation that combines (13.9) and (13.10) is the following:

$$\begin{array}{ll} \text{Minimize} & \Big(f(x) + E[\text{Min}_{y\in\mathbb{R}^{n_1}_+} \{q(y,\boldsymbol{\omega}) \mid T(\boldsymbol{\omega})x + W(\boldsymbol{\omega})y = h(\boldsymbol{\omega})\}]\Big) \\ \text{s.t.} & Ax = b, \\ & x \in \mathbb{R}^{n_0}_+, \end{array} \tag{13.11}$$

where *Min* denotes the minimal function value.

Let $K_1 \doteq \{x \in \mathbb{R}^{n_0}_+ \mid Ax = b\}$, denote the feasible set for the first-stage problem. Let also $K_2 \doteq \{x \in \mathbb{R}^{n_0} \mid E[\mathcal{Q}(x,\boldsymbol{\omega})] < +\infty\}$ denote the set of *induced constraints*. This is the set of first-stage decisions $x$ for which the second-stage problem is feasible. Problem (13.10) is said to have *complete recourse* if $K_2 = \mathbb{R}^{n_0}$, that is, if the second-stage problem is feasible for any value of $x$. The problem has *relatively complete recourse* if $K_1 \subseteq K_2$, that is, if the second-stage problem is feasible for any value of the first-stage variables that satisfies the first-stage constraints. *Simple recourse* refers to the case when the resource matrix $W(\boldsymbol{\omega}) = I$ and the recourse constraints take the simple form $Iy^+ - Iy^- = h(\boldsymbol{\omega}) - T(\boldsymbol{\omega})x$, where $I$ is the identity matrix, and the recourse vector $y$ is written as $y \doteq y^+ - y^-$ with $y^+ \geq 0$, $y^- \geq 0$.

### Deterministic Equivalent Formulation

We consider now the case where the random vector $\boldsymbol{\omega}$ has a discrete and finite distribution, with support $\Omega = \{\omega^1, \omega^2, \ldots, \omega^N\}$. In this case the set $\Omega$ is called a *scenario set.* Denote by $p_s$ the probability of realization of the $s$th scenario $\omega^s$. That is, for every $s = 1,2,\ldots,N$,

$$\begin{aligned} p_s &\doteq \text{Prob}\ (\boldsymbol{\omega} = \omega^s) \\ &= \text{Prob}\left\{(q(y,\boldsymbol{\omega}), W(\boldsymbol{\omega}), h(\boldsymbol{\omega}), T(\boldsymbol{\omega})) \;=\; (q(y,\omega^s), W(\omega^s), h(\omega^s), T(\omega^s))\right\}. \end{aligned}$$

It is assumed that $p_s > 0$ for all $\omega^s \in \Omega$, and that $\sum_{s=1}^N p_s = 1$.

The expected value of the second-stage optimization problem can be expressed as

$$E[Q(x,\omega)] = \sum_{s=1}^{N} p_s Q(x,\omega^s). \tag{13.12}$$

For each realization of the random vector $\omega^s \in \Omega$ a different second-stage decision is made, which is denoted by $y^s$. The resulting second-stage problems can then be written as

$$\begin{aligned} \text{Minimize} \quad & q(y^s,\omega^s) \\ \text{s.t.} \quad & W(\omega^s)y^s = h(\omega^s) - T(\omega^s)x, \\ & y^s \in \mathbb{R}^{n_1}_{+}. \end{aligned} \tag{13.13}$$

Combining now (13.12) and (13.13) we reformulate the stochastic nonlinear program (13.11) as the following *large-scale deterministic equivalent* nonlinear program

$$\text{Minimize} \quad f(x) + \sum_{s=1}^{N} p_s q(y^s,\omega^s) \tag{13.14}$$

$$\text{s.t.} \quad Ax = b, \tag{13.15}$$

$$T(\omega^s)x + W(\omega^s)y^s = h(\omega^s) \quad \text{for all } \omega^s \in \Omega, \tag{13.16}$$

$$x \in \mathbb{R}^{n_0}_{+}, \tag{13.17}$$

$$y^s \in \mathbb{R}^{n_1}_{+}. \tag{13.18}$$

The constraints (13.15)–(13.18) for this deterministic equivalent program can be combined into a matrix equation with block-angular structure

$$\begin{pmatrix} A & & & & \\ T(\omega^1) & W(\omega^1) & & & \\ T(\omega^2) & & W(\omega^2) & & \\ \vdots & & & \ddots & \\ T(\omega^N) & & & & W(\omega^N) \end{pmatrix} \begin{pmatrix} x \\ y^1 \\ y^2 \\ \vdots \\ y^N \end{pmatrix} = \begin{pmatrix} b \\ h(\omega^1) \\ h(\omega^2) \\ \vdots \\ h(\omega^N) \end{pmatrix}. \tag{13.19}$$

### Split-Variable Formulation

The system of linear equations in (13.19) can be rewritten in a form that is, for some algorithms, more amenable to decomposition and parallel computations. In particular, in the absence of the $x$ variables the system (13.19) becomes block-diagonal. The *split-variable* formulation replicates the first-stage variable vector $x$ into a set of vectors $x^s \in \mathbb{R}^{n_0}$ for each $\omega^s \in \Omega$. Once a different first-stage decision is allowed for each scenario, the stochastic program decomposes into $S$ independent problems. Of course, the first-stage variables must be *nonanticipative*, that is they cannot depend on scenarios that have not yet been observed when the first-stage decisions are made. This requirement is enforced by adding the restrictions that $x^1 = x^2 = \ldots = x^S$. The split-variable formulation is then equivalent to the original stochastic problem (13.14)–(13.18), for which equation (13.19) can be written in the equivalent form

$$\begin{pmatrix} A & & & & & \\ T(\omega^1) & W(\omega^1) & & & & \\ & & A & & & \\ & & T(\omega^2) & W(\omega^2) & & \\ & & & & \ddots & \\ & & & & A & \\ & & & & T(\omega^N) & W(\omega^N) \\ I & & -I & & & \\ \vdots & & & & \ddots & \\ I & & & & -I & \end{pmatrix} \begin{pmatrix} x^1 \\ y^1 \\ x^2 \\ y^2 \\ \vdots \\ x^N \\ y^N \end{pmatrix} = \begin{pmatrix} b \\ h(\omega^1) \\ b \\ h(\omega^2) \\ \vdots \\ b \\ h(\omega^N) \\ 0 \\ \vdots \\ 0 \end{pmatrix} \tag{13.20}$$

## Multistage Recourse Problems

The recourse problem is not restricted to the two-stage formulation. It is possible that observations are made at $K$ different stages and are captured in the information sets $\{\mathcal{A}_t\}_{t=1}^K$ with $\mathcal{A}_1 \subset \mathcal{A}_2 \cdots \subset \mathcal{A}_K$. These sets are subfields of the underlying $\sigma$-field $\Sigma$ of all possible observations. A multistage stochastic program with recourse will have a recourse problem at stage $\tau$ conditioned on the information provided by $\mathcal{A}_\tau$, which includes all information provided by the information sets $\mathcal{A}_t$, for $t = 1, 2, \ldots, \tau$. The program also anticipates the information in $\mathcal{A}_t$, for $t = \tau + 1, \ldots, K$.

Let the random vector $\boldsymbol{\omega}$ have support $\Omega = \Omega_1 \times \Omega_2 \times \cdots \times \Omega_K$, which is the product set of all individual support sets $\Omega_t$, $t = 1, 2, \ldots, K$. $\boldsymbol{\omega}$ is written componentwise as $\boldsymbol{\omega} = (\boldsymbol{\omega}^1, \ldots, \boldsymbol{\omega}^K)$. Denote the first-stage variable vector by $y^0$. For each stage $t = 1, 2, \ldots, K$ define the recourse variable vector $y^t \in \mathbb{R}^{n_t}$, the random cost function $q_t(y^t, \boldsymbol{\omega}^t)$, and the random parameters $\{W_t(\boldsymbol{\omega}^t), h_t(\boldsymbol{\omega}^t), T_t(\boldsymbol{\omega}^t) \mid \boldsymbol{\omega}^t \in \Omega_t\}$.

The multistage program, which extends the two-stage model (13.11), is formulated as the following nested optimization problem

$$\begin{aligned} \underset{y_0 \in \mathbb{R}_+^{n_0}}{\text{Minimize}} \quad & f(y^0) + E\left[ \underset{y^1 \in \mathbb{R}_+^{n_1}}{\text{Min}} q_1(y^1, \boldsymbol{\omega}^1) + \cdots E\left[ \underset{y^K \in \mathbb{R}_+^{n_K}}{\text{Min}} q_K(y^K, \boldsymbol{\omega}^K) \right] \cdots \right] \\ \text{s.t.} \quad & T_1(\boldsymbol{\omega}^1) y^0 + W_1(\boldsymbol{\omega}^1) y^1 = h_1(\boldsymbol{\omega}^1), \\ & \qquad \vdots \\ & T_K(\boldsymbol{\omega}^K) y^{K-1} + W_K(\boldsymbol{\omega}^K) y^K = h_K(\boldsymbol{\omega}^K). \end{aligned} \tag{13.21}$$

For the case of discrete and finitely distributed probability distributions it is again possible to formulate the multistage model into a deterministic equivalent large-scale nonlinear program. Section 13.9 provides references to the literature on multistage programs.

## 13.4 Robust Optimization Problems

This section considers an alternative approach to handling uncertain or noisy data. This approach is applicable to optimization models that have two distinct components: a *structural* component that is fixed and free of any noise in its input data; and a *control* component that is subject to noisy input data. In some cases the robust optimization model is identical to a two-stage stochastic program with recourse. But it also allows additional flexibility in dealing with noise. In order to define the model we introduce two sets of variables

$x \in \mathbb{R}^{n_0}$ denotes the vector of decision variables that depend only on the noise-free structural constraints. These are the *design* variables whose values are independent of realizations of the noisy parameters.

$y \in \mathbb{R}^{n_1}$ denotes the vector of *control* variables that can be adjusted once the uncertain parameters are observed. Their optimal values depend both on the realization of uncertain parameters, and on the optimal values of the design variables.

This terminology is borrowed from the flexibility analysis of production and distribution systems. The design variables determine the structure of the system and the size of production modules; the control variables are used to adjust the mode and level of production in response to disruptions in the system, changes in demand or production yield, and so on.

The optimization model we are interested in has the following structure

$$\text{Minimize} \quad \langle c, x\rangle + \langle d, y\rangle \tag{13.22}$$

$$\text{s.t.} \quad Ax = b, \tag{13.23}$$

$$Bx + Cy = e, \tag{13.24}$$

$$x \in \mathbb{R}^{n_0}_{+}, \tag{13.25}$$

$$y \in \mathbb{R}^{n_1}_{+}, \tag{13.26}$$

where $b, c, d, e$ are given vectors and $A, B, C$ are given matrices. Equation (13.23) denotes the *structural constraints* that are free of noise. Equation (13.24) denotes the *control constraints.* The coefficients of these constraints, i.e., the elements of $B, C,$ and $e$ are subject to noise. The cost vector $d$ is also subject to noise, while $A, b,$ and $c$ are not.

To define the robust optimization problem we introduce an index set $\Omega \doteq \{1, 2, \ldots, S\}$. With each index $s \in \Omega$ we associate the scenario set $\{d(s), B(s), C(s), e(s)\}$ of realizations of the control coefficients. Reference to an index $s$ imply reference to the scenario set associated with this index. The probability of the $s$th scenario is $p_s$, and $\sum_{s\in\Omega} p_s = 1$. Now the following question is posed: What are the desirable characteristics of a solution to problem (13.22)–(13.26) when the coefficients of the constraints (13.24) take values from some given set of scenarios? The solution is considered robust with respect to optimality if it remains *close* to optimal for any realization of the scenario index $s \in \Omega$. The problem is then termed *solution robust.* The solution is robust with respect to feasibility if it remains *almost* feasible for any realization of $s$. The problem is then termed *model robust.* The concepts of close and almost are precisely defined later through the choice of appropriate norms.

It is unlikely that a solution to the mathematical program will remain both feasible and optimal for all realizations of $s$. If the system being modeled has substantial built-in redundancies, then it might be possible to find solutions that remain both feasible and optimal. Otherwise a model is needed that permits a trade-off between solution and model robustness. The model developed next formalizes a way to measure this trade-off.

First let us introduce a set $\{y^1, y^2, \ldots, y^N\}$ of *control variables* for each scenario $s \in \Omega$, and another set $\{z^1, z^2, \ldots, z^N\}$ of *feasibility error vectors* that measure the infeasibility of the control constraints under each scenario.

The real-valued objective function $\xi(x, y) = \langle c, x\rangle + \langle d, y\rangle$ is a random variable taking the value $\xi_s(x, y^s) \doteq \langle c, x\rangle + \langle d(s), y^s\rangle$ with probability $p_s$. Hence, there is no longer a simple single choice for an aggregate objective function. The expected value

$$\sigma(\cdot) \doteq \sum_{s\in\Omega} p_s \xi_s(\cdot) \tag{13.27}$$

is precisely the objective function used in the stochastic programming formulations studied in the previous section. Another choice is to employ worst case analysis and minimize the maximal value. The objective function is then defined by

$$\sigma(\cdot) \doteq \max_{s\in\Omega} \xi_s(\cdot). \tag{13.28}$$

The robust optimization formulation also allows the introduction of higher moments of the distribution of $\xi(\cdot)$ in the optimization model. Indeed, the introduction of higher moments is one of the features of robust optimization that distinguishes it from the stochastic programming model of the previous sections. For example, we could use a nonlinear utility function that embodies a trade-off between mean value and variability in this mean value. If $\mathcal{U}(\xi_s)$ denotes the utility of $\xi_s$, then the function

$$\sigma(\cdot) \doteq \sum_{s\in\Omega} p_s \mathcal{U}(\xi_s(\cdot)) \tag{13.29}$$

captures the risk preference of the user. A popular choice of utility functions, for portfolio management applications, is the logarithmic function $\mathcal{U}(\xi_s) = \log \xi_s$. The general robust optimization model includes a term $\sigma(x, y^1, y^2, \ldots, y^N)$ in the objective function to denote the dependence of the function value on the scenario index $s$. This term controls solution robustness, and can take different forms depending on the application. The examples mentioned above are some popular choices.

The robust optimization model introduces a second term in the objective function to control model robustness. This term is a feasibility penalty function, denoted by $\rho(z^1, z^2, \ldots, z^N)$, and it is used to penalize violations of the control constraints under some of the scenarios. The introduction of this penalty function also distinguishes the robust optimization model from the stochastic programming approach for dealing with noisy data. In particular, the model recognizes that it may not always be possible to arrive at a feasible solution to a problem under all scenarios. Infeasibilities will inevitably arise, and they will be dealt with outside the optimization model. The robust optimization model generates solutions that present the modeler with the fewest infeasibilities to be dealt with outside the model.

The specific choice of penalty function is problem dependent, and it also has implications for the choice of a solution algorithm. Two suitable penalty functions are the following:

$\rho(z^1, z^2, \ldots, z^S) \doteq \sum_{s\in\Omega} p_s \|z^s\|^2$. This quadratic penalty function (i.e., a weighted $\ell_2$-norm) is applicable to equality control constraints where both positive and negative violations of the constraints are equally undesirable. The resulting quadratic programming problem is twice continuously differentiable, and can be solved using standard quadratic programming algorithms, although it is typically large scale.

$\rho(z^1, z^2, \ldots, z^S) \doteq \sum_{s\in\Omega} p_s \max(0, (\max_j z_j^s))$. This penalty function is applicable to inequality control constraints when only positive violations are of interest (negative values of some $z_j$ indicate slack in the inequality constraints). With this choice of penalty function, however, the resulting mathematical program is nondifferentiable. It is possible to use an $\epsilon$-smoothing of the exact penalty function, and employ the Linear Quadratic Penalty (LQP) algorithm (Chapter 7). The result is a differentiable problem that is easier to solve and produces a solution that lies within $\epsilon$ from the solution of the nondifferentiable problem.

The robust optimization model takes a multicriteria objective form. A goal programming weight parameter $\lambda$ is used to derive a spectrum of answers that trade off solution for model robustness. The general formulation of the robust optimization model is stated as follows:

$$
\begin{aligned}
\text{Minimize} \quad & \sigma(x, y^1, y^2, \ldots, y^N) + \lambda \rho(z^1, z^2, \ldots, z^N) \\
\text{s.t.} \quad & Ax = b, \\
& B(s)x + C(s)y^s + z^s = e(s) \qquad \text{for all } s \in \Omega, \\
& x \in \mathbb{R}_+^{n_0}, \\
& y^s \in \mathbb{R}_+^{n_1}, \\
& z^s \in \mathbb{R}_+^{m}.
\end{aligned}
$$

## 13.5 Applications

In this section we discuss real-world applications where uncertainty is prevalent and wherein it is handled using the models introduced above. We give first (subsection 13.5.1) an illustration of the robust optimization framework, using the classic diet problem as an example. The other subsections discuss models for planning under uncertainty in production and inventory management (subsection 13.5.2) and models for matrix balancing (subsection 13.5.3), with a new section devoted to models for financial planning under uncertainty (Section 13.6).

### 13.5.1 Robust optimization for the diet problem

The well-known diet problem is used here as an example to illustrate the feature of *model robustness.* This feature is particularly interesting in the context of optimization formulations, since feasibility has traditionally been overemphasized.

The problem is to find a minimum-cost diet that will satisfy certain nutritional requirements. The origins of this problem date back to the 1940s and to the works of G.J. Stigler and G.B. Dantzig where it was soon recognized as a problem of robust optimization, since the nutritional content of some food products may not be certain. Dantzig was still intrigued by this ambiguity when he wrote, several decades later:

> When is an apple an apple and what do you mean by its cost and nutritional content? For example, when you say *apple* do you mean a Jonathan, or McIntosh, or Northern Spy, or Ontario, or Winesap, or Winter Banana? You see, it can make a difference, for the amount of ascorbic acid (vitamin C) can vary from 2.0 to 20.8 units per 100 grams depending upon the type of apple. (Dantzig, 1990.)

The standard linear programming formulation assumes an average nutritional content for each food product and produces a diet. However, as consumers buy food products of varying nutritional content they will soon build a deficit or surplus of some nutrients. This situation may be irrelevant for a healthy individual over long periods of time, or it may require remedial action in the form of vitamin supplements.

We develop here a robust optimization formulation for this problem. Let $x_i$ denote the total cost of food type $i$ in the diet; let $a_{ij}$ denote the contents of food type $i$ in nutrient $j$ per unit cost; and let $b_j$ be the required daily allowance of nutrient $j$. Let $c$ denote one specific nutrient, e.g., vitamin C, from the set of nutrients and let $A$ denote, for example, apples, from the set of foods. A point estimate for the content of apples in vitamin C is $a_{Ac} = \alpha$ per unit cost. For the sake of the example assume that the coefficient can take any value $\{a^s_{Ac}\}$ for $s$ in a scenario set $\Omega$. The robust optimization formulation of the diet problem is then

$$\text{Minimize} \quad \sum_i x_i + \lambda \frac{1}{S} \sum_{s \in \Omega} (b_c - \sum_{i \neq A} a_{ic} x_i - a^s_{Ac} x_A)^2 \tag{13.30}$$

$$\text{s.t.} \quad \sum_i a_{ij} x_i = b_j \qquad \text{for all } j \neq c, \tag{13.31}$$

$$x_i \geq 0 \qquad \text{for all } i. \tag{13.32}$$

The weight $\lambda$ is used to trade off feasibility robustness with cost. For $\lambda = 0$, and also allowing the index $j$ in (13.31) to take the value $j = c$, with $a_{Ac} = \alpha$, we obtain the classic linear programming formulation.

The diets obtained with larger values of $\lambda$ have vitamin C content that varies very little with respect to the quality of apples. Figure 13.1 illustrates the error in vitamin C intake due to alternative robust optimization solutions under different scenarios of vitamin C content. This figure illustrates the efficacy of the robust optimization model in hedging against alternative realizations of the data. For example, if an error of $\pm 2\epsilon$ units in total vitamin C intake is acceptable, no remedial action will be needed for the robust optimization dieter—no matter what quality of apples are included in the diet. On the other hand, the linear programming dieter will need remedial treatment for several of the available apple qualities (note, from Figure 13.1, that the error in vitamin C intake exceeds the allowable margin of $\pm 2\epsilon$ under seven scenarios).

The robust optimization diet is somewhat more expensive than the diet produced by the linear program. Figure 13.2 shows the increase in the cost of the diet as it becomes more robust with respect to nutritional content. This simple example clarifies the meaning of a robust solution and shows that robust solutions are indeed possible, but at some cost.

### 13.5.2 Robust optimization for planning capacity expansion

Manufacturing and service firms have to plan for capacity expansion in order to meet increasing demand for their products and services over time. Demand, however, is usually highly uncertain. Product demand depends on general economic conditions, competition, technological changes, and the general business cycle. Demand may also exhibit seasonal variations, which are particularly difficult to address in the context of service oper-ations since an inventory of services cannot be created during periods of low demand. Public utility companies—power and water distribution—face both seasonal and daily variations. The hotel and travel businesses are also highly seasonal. Stochasticities are not restricted to demand alone; equipment failure, delivery of material by suppliers, and routine maintenance operations also show varying degrees of uncertainty.

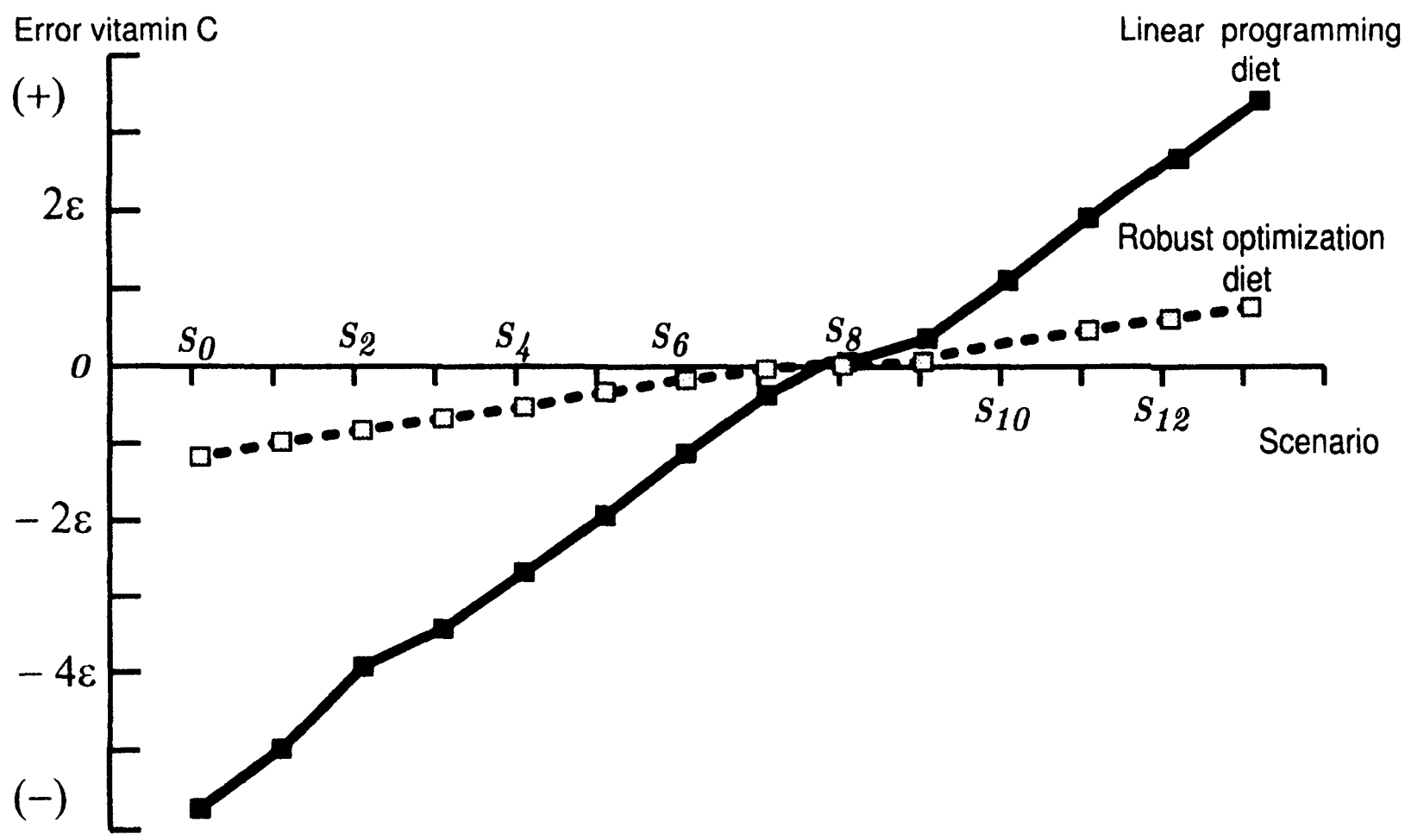

**Figure 13.1** Error (negative for deficit, positive for surplus) of the dieter's intake of vitamin C. The leftmost scenario corresponds to the lowest vitamin C content, while the rightmost scenario corresponds to the highest vitamin C content. The diet obtained with the linear programming formulation (i.e., $\lambda = 0$) is very sensitive to the vitamin content of the food products, whereas the diet obtained with the robust optimization model (i.e., $\lambda = 5$) is much less sensitive.

Planning for capacity expansion must account for the stochastic aspects of the system. A conservative plan may be developed by estimating capacity based on worst case analysis, wherein capacity must be sufficient to meet the demands of the last customer during the peak season. This strategy is quite expensive, and the prudent manager must carefully weigh the marginal profit from the last sale against the cost of maintaining underutilized manufacturing or service facilities. The United States automobile industry arrived at this conclusion in the early eighties and slashed capacity to increase utilization of their facilities, even at the cost of losing some potential customers. The decision made front-page news in the *The Wall Street Journal* on October 7, 1986. The Ford's chief financial officer summarized the essence of this strategy:

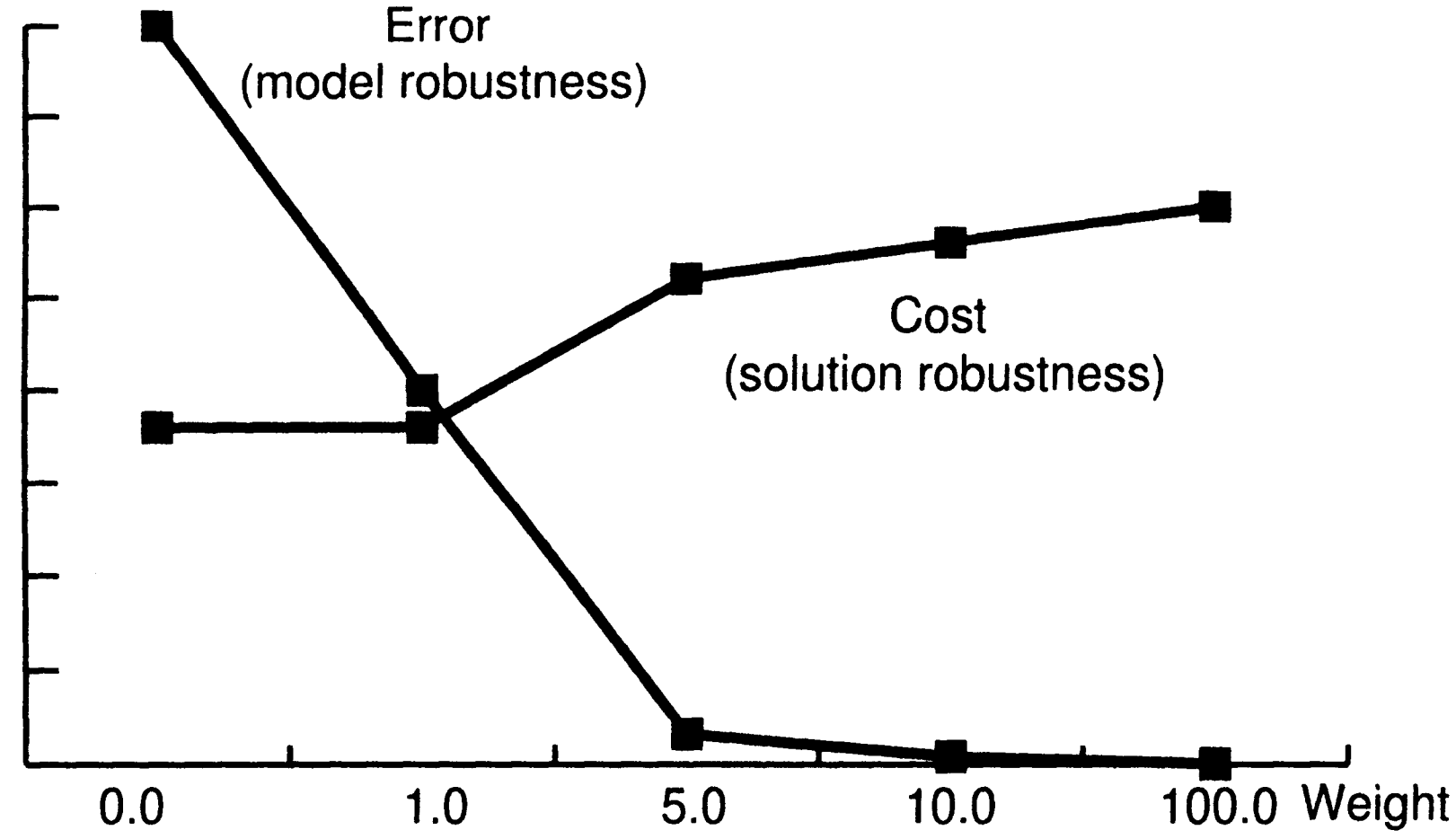

**Figure 13.2** Trade-off between cost (i.e., solution robustness) and expected error in the vitamin C contents (i.e., model robustness) for diets obtained using increasing values of the weight parameter $\lambda$ in the robust optimization model.

> "We arrived at a conscious willingness to give up the last vehicle we needed in peak years."

In this section we formulate a stochastic programming model that plans capacity expansion for a manufacturing firm.

### The Multiperiod Stochastic Program

We consider, for simplicity, the example of a firm producing a single product at multiple manufacturing sites. The model deals with uncertain product demand, and can be easily extended to handle multiple products, which is a more realistic situation. We begin with a stochastic programming formulation such that solution robustness is dealt with as an extension of the stochastic model.

The model makes capacity expansion decisions at $K$ plant sites during a planning horizon of $T$ time periods. There are $H$ different decisions that can be made with respect to each plant site, denoted by $h = 0, 1, 2, \ldots, H$ where $h = 1$ signifies the current state of the plant, $h = 0$ indicates that the plant is shut down and $h = 2, 3, \ldots, H$ signal various retooling options that are possible at each plant. Stochasticity in demand is dealt with by postulating a set of scenarios $\Omega = \{1, 2, \ldots, N\}$.

The model makes first-stage decisions on capacity planning, i.e., which plants to shut down or to retool, and which plants to maintain at their current status. The recourse decisions are the production levels for each time period and under the different scenarios.

### Notation

We define first the parameters of the model: use $k = 1, 2, \ldots, K$ to denote plant sites, $h = 0, 1, 2, \ldots, H$ to denote plant configurations, $s \in \Omega$ to denote scenarios, and $t = 1, 2, \ldots, T$ to denote time periods.

$p_{st}$ : the probability that scenario $s$ occurs at time period $t$.

$d_{st}$ : the demand for the product under scenario $s$, during time period $t$.

$\alpha$ : the fraction of unmet demand that is translated to profit by being diverted to other products of the same firm.

The capacity parameters are:

$U_{kht}$ : the capacity available at site $k$ under configuration $h$ during the $t$th time period.

$a_{kh}$ : the production coefficients, indicating the capacity required to produce one unit of the product at site $k$ under configuration $h$.

$L_{kht}$ : the capacity lost during the retooling process if site $k$ is retooled into configuration $h$ during the $t$th time period.

The cost coefficients, given next, are discounted to the initial time period, using an appropriate discount rate:

$F_{kht}$ : the fixed cost for changing the configuration of site $k$ from $h = 1$ to $h = 0, 2, \ldots, H$ during the $t$th time period.

$r_{kht}$ : the marginal contribution of producing and selling one unit of the product at site $k$, using configuration $h$ at the $t$th time period.

$r$ : the marginal contribution realized when there is unmet demand and a fraction $\alpha$ of it is diverted to other products.

Now we define the decision variables. There are two sets of continuous variables that denote the schedule of production and the level of unmet demand. Two sets of integer variables are used to denote retooling decisions and plant configurations.

$x_{khst}$ : the number of units produced at site $k$, operating using configuration $h$ under scenario $s$ during the $t$th time period.

$z_{st}$ : the number of units of unmet demand under scenario $s$ during the $t$th time period.

$$y_{kht} \doteq \begin{cases} 1 & \text{if site } k \text{ is in configuration } h \text{ at time period } t, \\ 0 & \text{otherwise.} \end{cases}$$

$$w_{kht} \doteq \begin{cases} 1 & \text{if site } k \text{ is retooled to configuration } h \text{ at time period } t, \\ 0 & \text{otherwise.} \end{cases}$$

The continuous variables are constrained to be nonegative. The integer variables are binary, that is $y_{kht}, w_{kht}$ are either 0 or 1.

### Model Formulation

We now define the model by specifying precisely the objective function and the constraints.

*Objective function:* The objective function that must be maximized has two terms that account for direct profits from sales of the product and from indirect profits from diverted demand, and a third term that accounts for the retooling cost. It takes the form:

$$\sum_{s=1}^{N}\sum_{t=1}^{T} p_{st}\left(\sum_{k=1}^{K}\sum_{h=0}^{H} r_{kht}x_{khst}\right) + \alpha r \sum_{s=1}^{N}\sum_{t=1}^{T} p_{st}z_{st} - \sum_{t=1}^{T}\sum_{k=1}^{K}\sum_{h=0}^{H} F_{kht}w_{kht}. \tag{13.33}$$

*Demand constraints:* For each time period and for each scenario the total production from all sites under all configurations plus the unmet demand diverted to other products is equal to the total realized demand.

$$\sum_{k=1}^{K}\sum_{h=0}^{H} x_{khst} + z_{st} = d_{st} \text{ for all } s = 1, 2, \ldots, N,\ t = 1, 2, \ldots, T. \tag{13.34}$$

*Capacity constraints:* The total production capacity utilized at each site cannot exceed the capacity available at the site under the given configuration, and taking into account losses of production capacity due to retooling operations. This is described for all $k = 1, 2, \ldots, K$, all $h = 0, 2, \ldots, H$, all $s = 1, 2, \ldots, N$, and all $t = 1, 2, \ldots, T$ by

$$a_{kh}x_{khst} \leq U_{kht}y_{kht} - L_{kht}. \tag{13.35}$$

*Retooling constraints:* We consider now the logical conditions between the retooling decisions and the plant configurations. A plant cannot be in a given configuration $h \neq 1$ unless it has first been retooled from its original configuration $h = 1$. For the first time period this condition is imposed by the constraint $y_{kh1} \leq w_{kh1}$. For subsequent time periods we impose the constraints $y_{kht} - y_{kh(t-1)} \leq w_{kh1}$ for all $k = 1, 2, \ldots, K$, all $h = 0, \ldots, H$, and all $t = 1, 2, \ldots, T$.

Operational considerations usually dictate that a plant cannot be retooled more than once during the planning horizon. Hence, we add the constraint:

$$\sum_{h=0}^{H}\sum_{t=1}^{T} w_{kht} \leq 1 \text{ for all } k = 1, 2, \ldots, K. \tag{13.36}$$

Finally, we require that each plant operate under some configuration or be shut down. These considerations are imposed by the constraints:

$$\sum_{h=0}^{H} y_{kht} = 1 \text{ for all } k = 1, 2, \ldots, K,\ t = 1, 2, \ldots, T, \tag{13.37}$$

$$w_{k0t} \leq y_{k0t} \text{ for all } k = 1, 2, \ldots, K,\ t = 1, 2, \ldots, T. \tag{13.38}$$

### Robustness Considerations

In practical capacity expansion applications, the expected cost or profit of the decision is not the only objective. Because large amounts of capital and other resources are involved, and the careers of many employees are at stake, some form of risk measure should be incorporated into the model with the objective of reducing it. Models of robust optimization—and especially those with solution robustness—provide the framework for dealing with risk.

In order to produce capacity expansion plans that are solution robust we define first a scenario-dependent measure of marginal profit from production $P_{st}$ at each time period $t$. We also define the total profit $P_s$ under each scenario $s$ as the total marginal profit from production during the planning horizon, less the fixed cost of the capacity expansion plan. We then impose additional constraints on the optimization model (13.34)–(13.38) so that some acceptable level of profit is realized under all scenarios.

The marginal profit realized under scenario $s$, from a given capacity plan and production schedule, is given by:

$$P_{st} \doteq \sum_{k=1}^{K} \sum_{h=0}^{H} r_{kht} x_{khst} + \alpha r \sum_{t=1}^{T} z_{st}. \tag{13.39}$$

The total profit for a given scenario $s$, accounting for the fixed cost of a capacity expansion plan, is given by:

$$P_s \doteq \sum_{t=1}^{T} P_{st} - \sum_{t=1}^{T} \sum_{k=1}^{K} \sum_{h=0}^{H} F_{kht} w_{kht}. \tag{13.40}$$

One way to introduce solution robustness now is to obtain solutions that have maximum expected profit for a given level of variance of profit. This approach, which is common in the finance literature, is known as the *Markowitz criterion.* A range of such solutions can be obtained by trading off expected profit for variance, that is, reducing the level of acceptable variance which also reduces the expected profit. This trade-off can be achieved by setting up the objective of the optimization model as follows:

$$\text{Maximize } \left(E[P_s] - \lambda \text{Var}[P_s]\right), \tag{13.41}$$

where $E[\cdot]$ and $\text{Var}[\cdot]$ denote the expected value and the variance of the random variable respectively, and $\lambda$ is a user specified parameter. Large values of $\lambda$ reduce variance at the expense of reduced profits, while smaller values allow the variance to increase, producing higher returns at a penalty of increased risks.

Calculation of the variance requires evaluation of a quadratic function and the resulting optimization program becomes a nonlinear program with continuous and integer variables. Such problems are very difficult to solve with currently available software systems, and they are likely to remain so in the future. Moreover, the use of a variance term as a measure of risk is inappropriate when the distribution of profits is not symmetric. The decision-maker wishes to reduce only downside risk, but by reducing variance the model reduces both upside potential and downside risk. An alternative formulation for robust capacity expansion planning imposes constraints that limit downside risk alone and gives rise to a linear optimization program. A target level of profit $\rho$ is set and linear constraints are added to the optimization program of the previous section such that the profit is greater than $\rho$ for all scenarios, i.e., $P_s \geq \rho$ for all $s = 1, 2, \ldots, N$. If $\rho$ is assigned a large negative value then these constraints are not binding and the model obtains a capacity expansion plan that maximizes expected profit. As $\rho$ is increased the model is forced to seek capacity plans that are guaranteed to have a profit of at least $\rho$ under all scenarios. Such plans, however, have reduced expected profit.

### 13.5.3 Robust optimization for matrix balancing

The problem of matrix balancing was examined earlier (Chapter 9), and entropy optimization models were developed wherein the problem data were assumed consistent, so that the sets of feasible solutions to Problems 9.2.3 and 9.2.4 were nonempty. For inconsistent data an interval-constrained formulation was proposed that allowed the solution of the matrix balancing problem in a way that satisfied the constraints within an error of $\pm\epsilon$. Here we develop a robust optimization model that provides an alternative way to overcome difficulties due to data inconsistency.

Consider the following equality-constrained optimization problem

$$\text{Minimize} \quad \sum_{i=1}^{m} \sum_{j=1}^{n} x_{ij} \log \left( \frac{x_{ij}}{a_{ij}} \right) \tag{13.42}$$

$$\text{s.t.} \quad \sum_{j=1}^{n} x_{ij} = u_i \qquad \text{for } i = 1, 2, \ldots, m, \tag{13.43}$$

$$\sum_{i=1}^{m} x_{ij} = v_j \qquad \text{for } j = 1, 2, \ldots, n, \tag{13.44}$$

$$x \geq 0. \tag{13.45}$$

When the observation vectors $u$ and $v$ are noisy it is possible that this optimization problem has no solution. Clearly, if $\sum_{i=1}^{m} u_i \neq \sum_{j=1}^{n} v_j$ then the optimization problem has no feasible solution. Several approaches can be pursued in order to overcome this difficulty. Tradition suggests that the vectors $u$ and $v$ be first scaled using the transformation $u_i \leftarrow u_i \left( (\sum_j v_j)/(\sum_i u_i) \right)$ for all $i$ so that feasibility is restored. Alternatively, the *interval constrained* formulation suggested in Problem 9.2.5 may be used, becasue with sufficiently large values of the interval parameter $\epsilon$ the program becomes feasible.

The robust optimization formulation of the matrix estimation problem explicitly accounts for potential infeasibilities in the linear constraints. It then introduces a penalty term in the objective function that minimizes a norm of the infeasibilities. Let $y \in \mathbb{R}^m$ and $z \in \mathbb{R}^n$ denote the infeasibility vectors for the constraints (13.43) and (13.44) respectively. The robust optimization model is written as

$$\text{Minimize} \quad \sum_{i=1}^{m}\sum_{j=1}^{n} x_{ij}\log\left(\frac{x_{ij}}{a_{ij}}\right) + \frac{\lambda}{2}\left(\sum_{i=1}^{m} y_i^2 + \sum_{j=1}^{n} z_j^2\right) \tag{13.46}$$

$$\text{s.t.} \quad \sum_{j=1}^{n} x_{ij} - y_i = u_i \quad \text{for } i = 1, 2, \ldots, m, \tag{13.47}$$

$$\sum_{i=1}^{m} x_{ij} - z_j = v_j \quad \text{for } j = 1, 2, \ldots, n, \tag{13.48}$$

$$x \geq 0. \tag{13.49}$$

This formulation is a direct application of the robust optimization model with quadratic penalty for infeasibilities. It is possible to arrive at a similar mathematical formulation, beginning with statistical arguments on the desirable properties of the balanced matrix. In particular, this formulation results in a Bayesian estimate of the matrix, that is, it maximizes (the logarithm of) the probability of the matrix $X = (x_{ij})$, conditional on the noisy observations $\{u_i, v_j\}$. The entropy term in the objective function estimates the matrix that is the least biased (or maximally uncommitted) with respect to missing information, conditioned on the observations $\{a_{ij}\}$. The quadratic terms are the logarithms of the probability distribution function of the error (i.e., noise) term, conditioned on the matrix $X = (x_{ij})$, assuming that the errors are normally distributed with mean zero and standard deviations that are identical for all observations.

### Row-action Algorithm for Robust Optimization of Matrix Balancing

Both the quadratic and entropy terms in the objective function of the robust optimization models are Bregman functions (see Section 2.1) and so is their sum. It is easy to verify that both functions have the strong zone consistency property with respect to the hyperplanes specified by the equality constraints (13.47)–(13.48). Hence, we can apply a row-action algorithm (Algorithm 6.4.1) to solve this model.

Consider first the application of the iterative step of Algorithm 6.4.1 to the $i$th row of the constraints (13.47). At the $\nu$th iteration it takes the form

$$x_{ij}^{\nu+1} = x_{ij}^{\nu} \exp \beta \quad \text{for all } j = 1, 2, \ldots, n, \tag{13.50}$$

$$y_i^{\nu+1} = y_i^{\nu} - \frac{\beta}{\lambda}. \tag{13.51}$$

It is also important to observe that since the algorithm is applied to equality constraints there is no need to explicitly update the dual variables and so they are omitted from Algorithm 13.5.1 below. The projection parameter $\beta$ is calculated such that $x_{ij}^{\nu+1}, y_i^{\nu+1}$ satisfy the $i$th constraint. Hence, substituting (13.50)–(13.51) to the $i$th equation of (13.47) we obtain

$$(\sum_{j=1}^{n} x_{ij}^{\nu} \exp \beta) - y_i^{\nu} + \frac{\beta}{\lambda} = u_i. \tag{13.52}$$

Let $\Psi(\beta)$ be the nonlinear function

$$\Psi(\beta) \doteq (\sum_{j=1}^{n} x_{ij}^{\nu} \exp \beta) - y_i^{\nu} + \frac{\beta}{\lambda} - u_i. \tag{13.53}$$

We seek a $\beta^*$ such that $\Psi(\beta^*) = 0$. We may use any nonlinear equation solver (e.g., Newton's method) to solve $\Psi(\beta) = 0$. However, a sufficient approximation to $\beta^*$ can be obtained by taking a single Newton's step, starting from $\beta = 0$, similar to the use of secant approximations discussed earlier (Sections 6.9 and 12.4.2). The asymptotic convergence of the row-action algorithm is preserved if this approximation is used in the iterative step instead of the exact value $\beta^*$.

It is worth mentioning that this approximation can be calculated in closed form, whereas the exact calculation of $\beta^*$ requires the use of an iterative procedure. Obtaining closed-form solutions to the nonlinear system of equations has important implications for the efficient implementation of the algorithm (see Section 12.4.2). The approximate solution of the nonlinear equation is obtained as

$$\tilde{\beta} = -\frac{\Psi(0)}{\Psi'(0)}, \tag{13.54}$$

where $\Psi'$ denotes first derivative with respect to $\beta$. Straightforward calculations yield

$$\tilde{\beta} = -\frac{(\sum_{j=1}^{n} x_{ij}^{\nu}) - y_i^{\nu} - u_i}{\frac{1}{\lambda} + \sum_{j=1}^{n} x_{ij}^{\nu}}. \tag{13.55}$$

This value of $\tilde{\beta}$ is used in (13.50)–(13.51) to complete the iterative step over the constraint equations (13.47). Following a similar argument we obtain the iterative step of the row-action algorithm over the constraint equations (13.48). The complete algorithm is summarized below.

**Algorithm 13.5.1 Row-Action Algorithm for Robust Matrix Balancing Optimization Model**

**Step 0:** (Initialization.) Set $\nu = 0$ and choose $x^0 \in \mathbb{R}^{mn}_{++}$ and $y^0 \in \mathbb{R}^m$, $z^0 \in \mathbb{R}^n$, such that the initialization conditions of Algorithm 6.4.1 are satisfied. For example, set $x^0_{ij} = a_{ij}$, $y^0_i = 0$, $z^0_j = 0$ for all $i = 1, 2, \ldots, m$, $j = 1, 2, \ldots, n$.

**Step 1:** (Iterative step over rows of the matrix, i.e., equations (13.47).) For all $i = 1, 2, \ldots, m$, calculate:

$$\tilde{\beta} \leftarrow -\frac{(\sum_{j=1}^{n} x^{\nu}_{ij}) - y^{\nu}_i - u_i}{\frac{1}{\lambda} + \sum_{j=1}^{n} x^{\nu}_{ij}}, \tag{13.56}$$

$$x^{\nu+\frac{1}{2}}_{ij} = x^{\nu}_{ij} \exp \tilde{\beta}, \tag{13.57}$$

$$y^{\nu+1}_i = y^{\nu}_i - \frac{\tilde{\beta}}{\lambda}. \tag{13.58}$$

**Step 2:** (Iterative step over columns of the matrix, i.e., equations (13.48).) For all $j = 1, 2, \ldots, n$, calculate:

$$\tilde{\beta} \leftarrow -\frac{(\sum_{i=1}^{m} x^{\nu+\frac{1}{2}}_{ij}) - z^{\nu}_j - v_j}{\frac{1}{\lambda} + \sum_{i=1}^{m} x^{\nu+\frac{1}{2}}_{ij}}, \tag{13.59}$$

$$x^{\nu+1}_{ij} = x^{\nu+\frac{1}{2}}_{ij} \exp \tilde{\beta}, \tag{13.60}$$

$$z^{\nu+1}_j = z^{\nu}_j - \frac{\tilde{\beta}}{\lambda}. \tag{13.61}$$

**Step 3:** Replace $\nu \leftarrow \nu + 1$, and return to Step 1.

## 13.6 Stochastic Programming for Portfolio Management

Portfolio management problems can be viewed as multiperiod dynamic decision problems where transactions take place at discrete time points.

At each point in time the manager has to assess the prevailing market conditions (such as prices and interest rates) and the composition of the existing portfolio. The manager has also to assess potential fluctuations in interest rates, prices, and cashflows. This information is incorporated into a sequence of actions of buying or selling securities, and short-term borrowing or lending. Thus, at the next point in time the portfolio manager has a seasoned portfolio and, faced with a new set of possible future movements, must incorporate the new information so that transactions can be executed.

Portfolio management of equities is based on the notion of *diversification* introduced by H. Markowitz in his seminal work in the 1950s. Diversification is achieved by minimizing the variance of returns during a holding period, subject to constraints on the mean value of the returns. There is only one time interval under consideration. Therefore, future transactions are not incorporated and this is a single-period (myopic) model. The portfolio management strategy for fixed-income securities has been that of *portfolio immunization*, i.e., portfolios are developed that are hedged against small changes from the current term structure of interest rates. Such models are again single-period and ignore future transactions. Furthermore, they ignore the truly stochastic nature of interest rates, and merely hedge against (small) changes from currently observed data. The idea of immunization dates back to the actuary F.M. Reddington in the 1950s, and it has been used extensively since the mid-70s.

The increased complexity of the fixed-income securities and the increased volatility of the financial markets during the 1980s have motivated interest in mathematical models that more accurately capture the dynamic (i.e., multiperiod) and stochastic nature of the portfolio management problem. Stochastic programming models with recourse provide a versatile tool for the representation of a wide variety of portfolio management problems. This section formulates a multistage stochastic programming model for managing portfolios of fixed-income securities. We assume that readers have some familiarity with basic concepts of finance. As an introductory text we recommend Bodie, Kane, and Marcus (1989) and, for more advanced material, Elton and Gruber (1984) and Zenios (1993a).

The model specifies a sequence of investment decisions at discrete time points. Decisions are made at the beginning of each time period. The portfolio manager starts with a given portfolio and a set of scenarios about future states of the economy which she incorporates into an investment decision. The precise composition of the portfolio depends on transactions at the previous decision point and on the realized scenario. Another set of investment decisions are made that incorporate both the current status of the portfolio and new information about future scenarios.

We develop a three-stage model, with decisions made at time instances $t_0, t_1$, and $t_2$. Extension to a multistage model is straightforward. Scenarios unfold between $t_0$ and $t_1$, and then again between $t_1$ and $t_2$. A simple three-stage problem is illustrated in Figure 13.3, where it is assumed that scenarios evolve on a binomial lattice. At instance $t_0$ two scenarios are anticipated and by instance $t_1$ this uncertainty is resolved. Denote these scenarios by $s_0^0$ and $s_0^1$. At $t_1$ two more scenarios are anticipated, $s_1^0$ and $s_1^1$. A complete path is denoted by a pair of scenarios. In this example there are four paths from $t_0$ to $t_2$ denoted by the pairs $(s_0^0, s_1^0), (s_0^0, s_1^1), (s_0^1, s_1^0)$, and $(s_0^1, s_1^1)$.

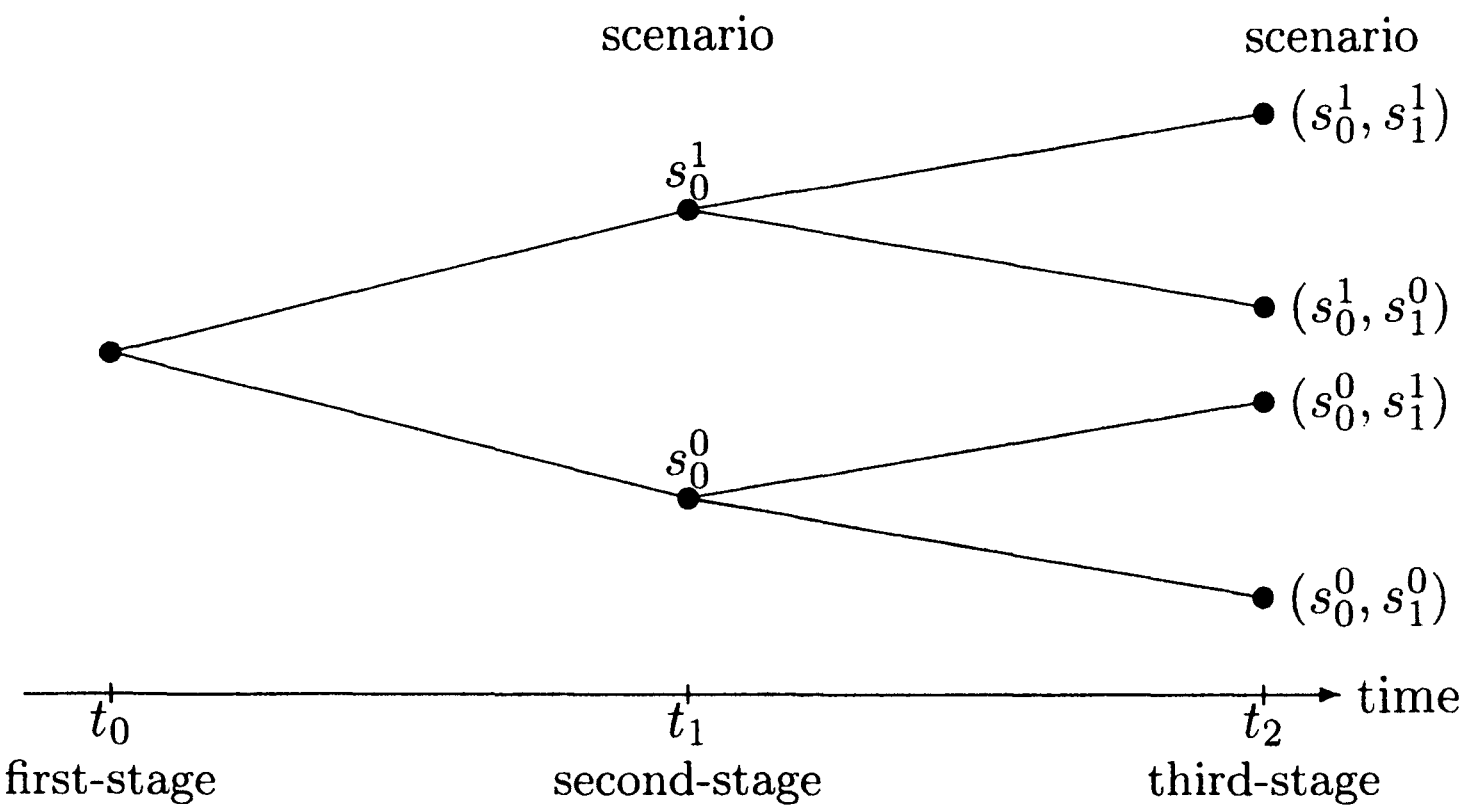

**Figure 13.3** The evolution of scenarios on a binomial lattice and the structure of the portfolio investment decisions.

The model is a nonlinear program maximizing the expected value of a utility function of terminal wealth. Expectations are computed over the set of postulated scenarios. The value of the utility function is obtained from the solution of a multistage program, which is structured in such a way that decisions at every time period anticipate future scenarios, while they adapt to new information on market conditions as it becomes available.

### 13.6.1 Notation

The model is developed using cashflow accounting equations. Investment decisions are in dollars of face value. We define first the parameters of the model.

$S_0$, $S_1$ : the sets of scenarios anticipated at $t_0$ and $t_1$ respectively. We use $s_0$ and $s_1$ to denote scenarios from $S_0$ and $S_1$, respectively. Paths are denoted by pairs of the form $(s_0, s_1)$, and with each path we associate a probability $\pi(s_0, s_1)$.

$I$ : the set of available securities. The cardinality of $I$ (i.e., number of available financial instruments) is $m$.

$c_0$ : the dollar amount of riskless asset available at $t_0$.

$b^0 \in \mathbb{R}^m$ : a vector whose components denote the composition of the initial portfolio.

$\zeta^0$, $\xi^0 \in \mathbb{R}^m$ : vectors of bid and ask prices respectively, at $t_0$. These prices are known with certainty. In order to buy an instrument the buyer has to pay the bid price, and in order to sell it the owner is asking for the ask price.

$\zeta^1(s_0)$, $\xi^1(s_0) \in \mathbb{R}^m$, for all $s_0 \in S_0$ : vectors of bid and ask prices, respectively, realized at $t_1$. These prices depend on the scenario $s_0$.

$\zeta^2(s_0, s_1)$, $\xi^2(s_0, s_1) \in \mathbb{R}^m$, for all $s_0 \in S_0$ and all $s_1 \in S_1$ : vectors of bid and ask prices, respectively, realized at $t_2$. These prices depend on the path $(s_0, s_1)$.

$\alpha^0(s_0)$, $\alpha^1(s_0, s_1) \in \mathbb{R}^m$, for all $s_0 \in S_0$ and all $s_1 \in S_1$ : vectors of *amortization factors* during the time intervals $[t_0, t_1)$ and $[t_1, t_2)$ respectively. The amortization factors indicate the fraction of outstanding face value of the securities at the end of the interval compared to the outstanding face value at the beginning of the interval. These factors capture the effects of any embedded options, such as prepayments and calls, or the effect of lapse behavior. For example, a corporate security that is called during the interval has an amortization factor of 0, and an uncalled bond has an amortization factor of 1. A mortgage security that experiences a 10 percent prepayment and that pays, through scheduled payments, an additional 5 percent of the outstanding loan has an amortization factor of 0.85. These factors depend on the scenarios.

$k^0(s_0)$, $k^1(s_0, s_1) \in \mathbb{R}^m$, for all $s_0 \in S_0$, and all $s_1 \in S_1$: vectors of *cash accrual factors* during the intervals $[t_0, t_1)$ and $[t_1, t_2)$ respectively. These factors indicate cash generated during the interval, per unit face value of the security, due to scheduled payments and exercise of the embedded options, accounting for accrued interest. For example, a corporate security that is called at the beginning of a one-year interval, in a 10 percent interest rate environment, will have a cash accrual factor of 1.10. These factors depend on the scenarios.

$\rho_0(s_0)$, $\rho_1(s_0, s_1)$ : short-term riskless reinvestment rates during the intervals $[t_0, t_1)$ and $[t_1, t_2)$ respectively. These rates depend on the scenarios.

$L_1(s_0)$, $L_2(s_0, s_1)$ : liability payments at $t_1$ and $t_2$ respectively. Liabilities may depend on the scenarios.

Now let us define decision variables. We have four distinct decisions at each point in time: how much of each security to buy, sell, or hold in the portfolio, and how much to invest in the riskless asset. All variables are constrained to be nonnegative.

First-stage variables, at $t_0$:

$x^0 \in \mathbb{R}^m$ : the components of the vector denote the face value of each security bought.

$y^0 \in \mathbb{R}^m$ : denotes componentwise the face value of each security sold.
$z^0 \in \mathbb{R}^m$ : denotes componentwise the face value of each security held in the portfolio .
$v_0$ : the dollar amount invested in the riskless asset.

Second-stage variables, at $t_1$ for each scenario $s_0$:
$x^1(s_0) \in \mathbb{R}^m$ : denotes the vector of the face values of each security bought.
$y^1(s_0) \in \mathbb{R}^m$ : denotes the vector of the face values of each security sold.
$z^1(s_0) \in \mathbb{R}^m$ : denotes the vector of the face values of each security held in the portfolio.
$v_1(s_0)$ : the dollar amount invested in the riskless asset.

Third-stage variables, at $t_2$ for each scenario $(s_0, s_1)$:
$x^2(s_0, s_1) \in \mathbb{R}^m$ : denotes the vector of the face values of each security bought.
$y^2(s_0, s_1) \in \mathbb{R}^m$ : denotes the vector of the face values of each security sold.
$z^2(s_0, s_1) \in \mathbb{R}^m$ : denotes the vector of the face values of each security held in the portfolio.
$v_2(s_0, s_1)$ : the dollar amount invested in the riskless asset.

### 13.6.2 Model formulation

The constraints of the model express cashflow accounting for the riskless asset and inventory balance for each security at all time periods.

*First-stage constraints:* At the first stage (i.e., at time $t_0$) all prices are known with certainty. The *cashflow accounting equation* specifies that the original endowment in the riskless asset, plus any proceeds from liquidating part of the existing portfolio, equal the amount invested in the purchase of new securities plus the amount invested in the riskless asset, i.e.,

$$c_0 + \sum_{i=1}^{m} \xi_i^0 y_i^0 = \sum_{i=1}^{m} \zeta_i^0 x_i^0 + v_0. \tag{13.62}$$

For each security in the portfolio we have an *inventory balance constraint*:

$$b_i^0 + x_i^0 = y_i^0 + z_i^0 \text{ for all } i \in I. \tag{13.63}$$

*Second-stage constraints:* Decisions made at the second stage (i.e., at time $t_1$) depend on the scenario $s_0$ realized during the interval $[t_0, t_1)$. Hence, we have one constraint for each scenario. These decisions also depend on the investment decisions made at the first stage.

Cashflow accounting ensures that the amount invested in the purchase of new securities and the riskless asset is equal to the income generated by the existing portfolio during the holding period, plus any cash generated from sales, less the liability payments. There is one constraint for each scenario

$$\rho_0(s_0)v_0 + \sum_{i=1}^{m} k_i^0(s_0)z_i^0 + \sum_{i=1}^{m} \xi_i^1(s_0)y_i^1(s_0)$$
$$= v_1(s_0) + \sum_{i=1}^{m} \zeta_i^1(s_0)x_i^1(s_0) + L_1(s_0), \text{ for all } s_0 \in S_0. \quad (13.64)$$

Inventory balance equations constrain the amount of each security sold or remaining in the portfolio to be equal to the outstanding amount of face value at the end of the first period, plus any amount purchased at the beginning of the second stage. There is one constraint for each security and for each scenario:

$$\alpha_i^0(s_0)z_i^0 + x_i^1(s_0) = y_i^1(s_0) + z_i^1(s_0), \text{ for all } i \in I, \ s_0 \in S_0. \quad (13.65)$$

*Third-stage constraints:* Decisions made at the third stage (i.e., at time $t_2$) depend on the path $(s_0, s_1)$ realized during the period $[t_1, t_2)$ and on the decisions made at $t_1$. The constraints are similar to those of the second stage. The cashflow accounting equation is

$$\rho_1(s_0, s_1)v_1(s_0) + \sum_{i=1}^{m} k_i^1(s_0, s_1)z_i^1(s_0) + \sum_{i=1}^{m} \xi_i^2(s_0, s_1)y_i^2(s_0, s_1)$$
$$= v_2(s_0, s_1) + \sum_{i=1}^{m} \zeta_i^2(s_0, s_1)x_i^2(s_0, s_1) + L_2(s_0, s_1), \quad (13.66)$$

for all paths $(s_0, s_1)$ such that $s_0 \in S_0$ and $s_1 \in S_1$.

The inventory balance equation is:

$$\alpha_i^1(s_0, s_1)z_i^1(s_0) + x_i^2(s_0, s_1) = y_i^2(s_0, s_1) + z_i^2(s_0, s_1), \quad (13.67)$$

for all $i \in I$, and all paths $(s_0, s_1)$ such that $s_0 \in S_0$ and $s_1 \in S_1$.

*Objective function:* The objective function maximizes the expected utility of terminal wealth. In order to measure terminal wealth all securities in the portfolio are marked-to-market, in accordance with recent U.S. Federal Accounting Standards Board (FASB) regulations that require reporting portfolio market and book values. The composition of the portfolio and its market value depend on the scenarios $(s_0, s_1)$. The objective of the portfolio optimization model is

$$\text{Maximize} \sum_{(s_0,s_1)\in S_0\times S_1} \pi(s_0, s_1)\mathcal{U}\left(W(s_0, s_1)\right),$$

where $\pi(s_0, s_1)$ is the probability associated with the path $(s_0, s_1)$; $W(s_0, s_1)$ denotes terminal wealth; and $\mathcal{U}$ denotes the utility function. Terminal wealth is given by

$$W(s_0, s_1) \doteq v_2(s_0, s_1) + \sum_{i=1}^{m} \xi_i^2(s_0, s_1) z_i^2(s_0, s_1). \tag{13.68}$$

## 13.7 Stochastic Network Models

We revisit the deterministic equivalent formulation of the two-stage model (13.14)–(13.18) and address the special case where the constraints (13.15)–(13.16) can be represented by a generalized network structure. (See Section 12.1 for the definition of networks and generalized networks.) For that structure both $\begin{pmatrix} A \\ T(\omega^s) \end{pmatrix}$ and $W(\omega^s)$ are *node-arc incidence matrices* with exactly two nonzero entries per column for each scenario $\omega^s$ in a scenario set $\Omega \doteq \{\omega^1, \omega^2, \ldots, \omega^S\}$; the complete matrix in (13.19) is not however due to the occurrence of $T(\omega^s)x$ for all scenarios $\omega^s \in \Omega$. The matrix $\left(A^{\mathrm{T}} \mid T(\omega^1)^{\mathrm{T}} \mid \cdots \mid T(\omega^s)^{\mathrm{T}}\right)^{\mathrm{T}}$ has more than two nonzero entries in every column. The recourse problem has a network structure, but the first-stage variables $x$ are *complicating variables* as they link the recourse problem constraints (via the linear equations (13.16)). The next section gives a detailed formulation of the stochastic network problem, where the network structure becomes more discernible.

The financial modeling applications of Section 13.6 can be represented using stochastic network structures. For each fixed-income security at each time period we associate a network node; and for each transaction we associate an arc. For example, an arc linking two different securities at the same time period is used to denote sale of one security and purchase of the other, while an arc linking the same security across different time periods is used to denote the inventory of that particular security. Figure 13.4 illustrates the structure of the network flow problem for two securities and three time periods. If all data are known with complete certainty the problem is the classic network flow problem. The stochastic programming model uses a different network flow problem for each scenario, but the flows on the arcs representing first-stage variables are common across scenarios.

Problems of planning hydroelectric power scheduling can also be represented as stochastic network problems. Power generation decisions made today for the coming hours depend on the current state of the system, electricity demand, water inflows, and so on. The necessary input data for modeling such a complex system consist of: the network topology specified by the geographical location of dams and the associated hydropower generation units, and their interconnections; limits on reservoir storage, level of turbine operations, pumping, and spillage; hydroelectric production coefficients for each reservoir, obtained from engineering analysis of its storage capacity and the turbine technology. Important input data are also the water level in the reservoirs and the electricity demand. Both of these quantities are uncertain, and furthermore, demand for electricity exhibits both a daily and a seasonal variation that can be estimated, at best, by a set of scenarios. The same is also true for the water level that depends on rainfall. These uncertainties are fundamental to the operation of the system. They can be incorporated in a stochastic network model.

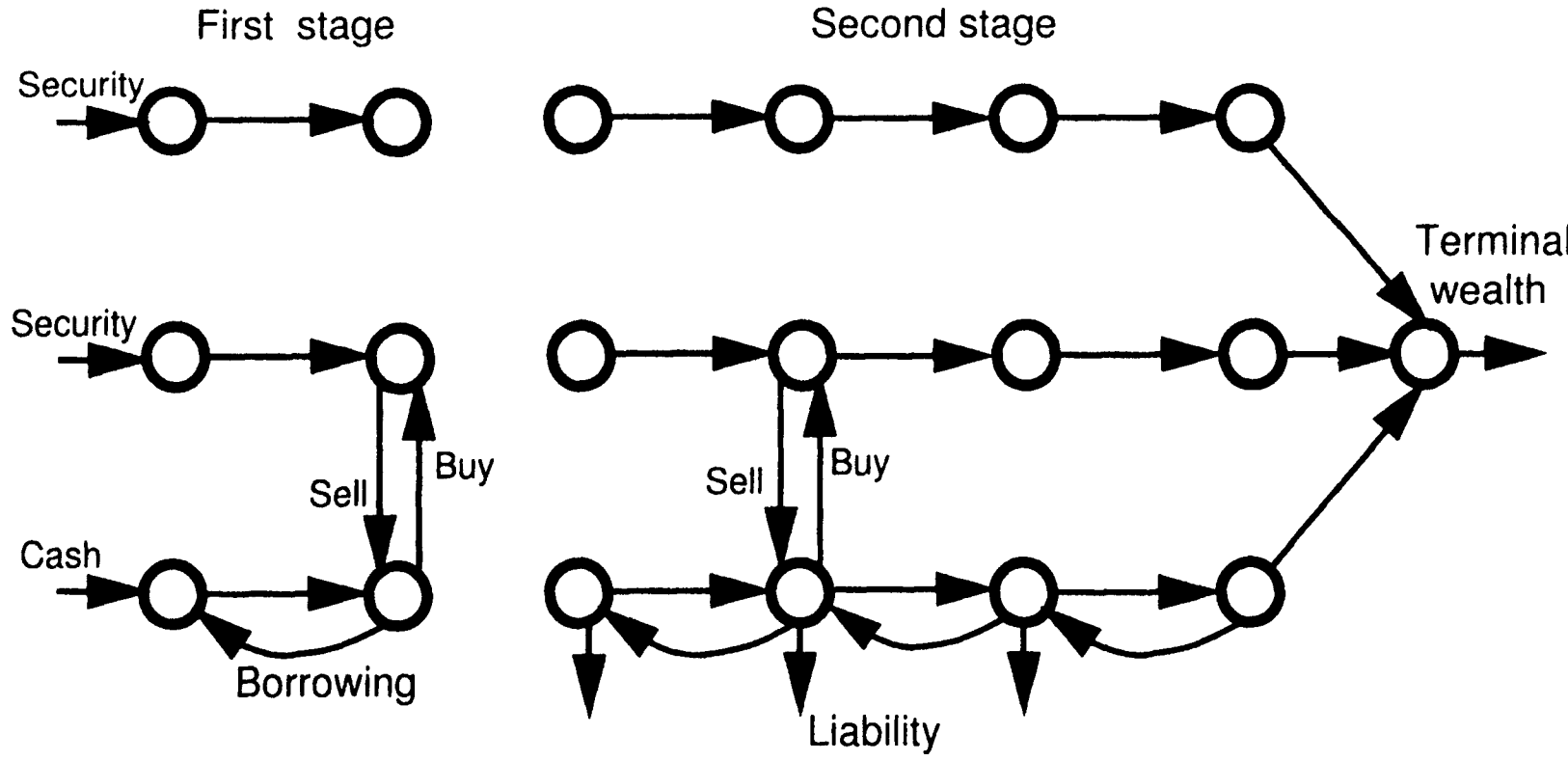

**Figure 13.4** The structure of a generalized network model for portfolio management of two securities over three time periods.

Another complex problem that has been modeled using stochastic networks is that of planning air-traffic ground holding policies (see also Section 12.3.2). The air traffic system is a complex web of airports, aircraft, and traffic controllers at all airports. In the United States a centralized flow control facility in Washington, DC, coordinates the control of this system. The complexity of the air traffic system in Europe is intensified by the need to coordinate several control centers among different countries, a situation which has become further complicated with the integration of east European countries into the system.

The systems are highly congested; traffic flow is carefully monitored and controlled so that flights proceed without risk to safety. One key control mechanism is *ground holding*, whereby a flight is delayed for departure if congestion is anticipated at the destination airport. Ground holding is a safe and relatively inexpensive solution, as opposed to holding aircraft in flight before granting landing clearance. While the air traffic control system does an excellent job monitoring traffic so that high safety standards are maintained, there is substantial room for improvement especially with regard to cost effectiveness. It is estimated that ground delays in the United States in 1986 averaged 2000 hours per day, equivalent to grounding a total of 250 airplanes (a carrier the size of Delta Airlines). A study by the West German Institute for Technology estimated the avoidable cost of air traffic delays in 1990 due to ground holding alone at 1.5 billion US dollars.

The ground holding policy problem seeks optimal holding policies, based on the number of flights scheduled for departure during the planning horizon and the travel time to the destination airport. Even for the simple case where only a single destination airport is analyzed, its capacity is uncertain due to weather conditions. The problem is complicated further by the presence of multiple airports: ground holding decisions at each one have a cascade effect on all others. The single destination airport problem has been modeled using stochastic network optimization models. Application of the model to data obtained from Logan airport (Boston, Mass., USA) proved that substantial reductions in total delay can be realized when using the stochastic programming dynamic models as opposed to more commonly used static models.

### 13.7.1 Split-variable formulation of stochastic network models

We consider here an alternative formulation of the deterministic equivalent formulation (13.14)–(13.18) that better illustrates the network structure, and is also more suitable to the development of parallel optimization algorithms (Section 13.8). The *split-variable formulation* (see also Section 13.3.3) breaks the stochastic network problem into a large number of independent deterministic network flow problems with some additional coupling constraints, by *replicating* the first-stage variables $x$ into a set of variables $x^s \in \mathbb{R}^{n_0}$ for each $\omega^s \in \Omega$. Once a different first-stage decision is allowed for each scenario the stochastic program decomposes into $S$ independent problems. Of course, the first-stage variables must be *non-anticipative*, that is, they cannot depend on as yet unobserved scenarios, a requirement that is enforced by adding the condition that $x^1 = x^2 = \cdots = x^S$. The split-variable formulation is then equivalent to the original stochastic problem (13.14)–(13.18). It can be written as follows

$$\text{Minimize} \quad \sum_{s=1}^{S} p_s(f(x^s) + q_s(y^s, \omega^s)) \tag{13.69}$$

$$\text{s.t.} \quad Ax^s = b \quad \text{for all } s \in \Omega, \tag{13.70}$$

$$T(\omega^s)x^s + W(\omega^s)y^s = h(\omega^s) \quad \text{for all } s \in \Omega, \tag{13.71}$$

$$x^1 - x^s = 0 \quad \text{for all } s \in \Omega, \tag{13.72}$$

$$x^s \in \mathbb{R}^{n_0}_{+}, \tag{13.73}$$

$$y^s \in \mathbb{R}^{n_1}_{+}. \tag{13.74}$$

The constraints (13.72) are known as *non-anticipativity constraints* and they ensure that first-stage decisions $x^s$ do not depend on future realizations. For two scenarios $s_1$ and $s_2$ that are indistinguishable when the first-stage decisions are made, we have that $x^{s_1} = x^{s_2}$.

With this reformulation the model (13.69)–(13.74) is a *network with side constraints*. In the absence of the (side) constraints (13.72), the constraint set decomposes completely into $S$ independent problems. Each problem has a network structure, since $(A^{\mathrm{T}} \mid T(\omega^s)^{\mathrm{T}})^{\mathrm{T}}$ and $W(\omega^s)$ are node-arc incidence matrices.

The constraints matrix of the split-variable formulation has a block-diagonal structure with additional (coupling) rows for the non-anticipativity constraints.

Let $M \doteq S\cdot(m_0+m_1)+(S-1)\cdot n_0$, $N \doteq S\cdot(n_0+n_1)$, and let $I$ denote the $n_0 \times n_0$ identity matrix. Recall that $m_0$ and $n_0$ are the numbers of first-stage constraints and variables respectively, and $m_1, n_1$ are the numbers of second-stage constraints and variables respectively. The constraint matrix for (13.70)–(13.72) has dimension $M \times N$. We denote this matrix by $\Phi$, (see also equation (13.20)), and it is defined as follows

$$\Phi \doteq \begin{pmatrix} A & & & & & & \\ T(\omega^1) & W(\omega^1) & & & & & \\ & & A & & & & \\ & & T(\omega^2) & W(\omega^2) & & & \\ & & & & \ddots & & \\ & & & & & A & \\ & & & & & T(\omega^S) & W(\omega^S) \\ I & & -I & & & & \\ \vdots & & & & \ddots & & \\ I & & & & & -I & \end{pmatrix}. \tag{13.75}$$

It is evident from the structure of the constraints matrix that the problem decomposes by scenario if the non-anticipativity constraints are ignored.

Let $\gamma \in \mathbb{R}^M$ denote the right-hand side of (13.70)–(13.72) and let $z \in \mathbb{R}^N$ be the vector of decision variables, i.e.,

$$z = ((x^1)^{\mathrm{T}} \mid (y^1)^{\mathrm{T}} \mid \cdots \mid (x^S)^{\mathrm{T}} \mid (y^S)^{\mathrm{T}})^{\mathrm{T}}. \tag{13.76}$$

We also introduce, for completeness, a vector $u \in \mathbb{R}^N$ whose components denote upper bounds on the variables. Finally, let $F(z)$ denote the objective function (13.69).

The split-variable formulation with bounded variables can be written in matrix form as

$$\text{Minimize} \quad F(z) \tag{13.77}$$

$$\text{s.t.} \quad \begin{pmatrix} \gamma \\ 0 \end{pmatrix} \leq \begin{pmatrix} \Phi \\ I_N \end{pmatrix} z \leq \begin{pmatrix} \gamma \\ u \end{pmatrix}, \tag{13.78}$$

$$z \in \mathbb{R}^N, \tag{13.79}$$

where $I_N$ is the $N \times N$ identity matrix.

### 13.7.2 Component-wise representation of the stochastic network problem

The previous sections illustrated the macro structure of the problem and the decomposable nature of the constraints matrix. This section gives a component-wise formulation of the stochastic network problem and illustrates its micro structure by specifying algebraically all equations. We assume for simplicity the same underlying network structure for all scenario problems, and with the notation of Section 9.2 we represent this structure by the graph $G = (\mathcal{N}, \mathcal{A})$, where $\mathcal{N} = \{1, 2, \ldots, m_0 + m_1\}$ is the set of nodes and $\mathcal{A} = \{(i,j) \mid i,j \in \mathcal{N}\} \subseteq \mathcal{N} \times \mathcal{N}$ is the set of arcs. Let $\delta_i^+ = \{j \mid (i,j) \in \mathcal{A}\}$ be the set of nodes having an arc with origin node $i$, and $\delta_j^- = \{i \mid (i,j) \in \mathcal{A}\}$ be the set of nodes having an arc with destination node $j$. We partition the set of all nodes into two disjoint sets, $\mathcal{N}_0$ and $\mathcal{N}_1$. The set $\mathcal{N}_0$ consists of the $m_0$ nodes whose incident arcs are all first stage so that their flow conservation constraints do not depend on the realization of the uncertain quantities. The resources (i.e., supply or demand) for these nodes are real numbers, denoted by $b_i$ for all $i \in \mathcal{N}_0$. The set $\mathcal{N}_1 = \mathcal{N} \setminus \mathcal{N}_0$ consists of the $m_1$ nodes with stochastic right-hand sides or incident second-stage arcs. The resources for these nodes are denoted by $r_i^s$ for all $i \in \mathcal{N}_1$, $s \in \Omega$.

We also partition the arc set $\mathcal{A}$ into two disjoint sets $\mathcal{A}_0$ and $\mathcal{A}_1$, corresponding to replicated first- and second-stage decisions respectively. The number of arcs in these sets are denoted by $n_0$ and $n_1$ respectively. Denote by $x_{ij}^s$ for $(i,j) \in \mathcal{A}_0$, and $y_{ij}^s$ for $(i,j) \in \mathcal{A}_1$ the flows on the arc with origin node $i$ and destination node $j$ under scenario index $s \in \Omega$. The upper bound of a replicated first-stage arc $x_{ij}^s$ is denoted by $u_{ij}$ and the upper bound of a second-stage arc $y_{ij}^s$ is denoted by $v_{ij}^s$. The multiplier on arc $(i,j)$ is denoted by $m_{ij}$ for $(i,j) \in \mathcal{A}_0$, and by $m_{ij}^s$ for $(i,j) \in \mathcal{A}_1$. The network optimization model for a fixed scenario index $s \in \Omega$ is given by:

$$\underset{x^s \in \mathbb{R}^{n_0}, y^s \in \mathbb{R}^{n_1}}{\text{Minimize}} \quad \sum_{(i,j)\in\mathcal{A}_0} p_s f_{ij}(x^s_{ij}) + \sum_{(i,j)\in\mathcal{A}_1} p_s q^s_{ij}(y^s_{ij}) \tag{13.80}$$

s.t.

$$\sum_{j\in\delta_i^+} x^s_{ij} - \sum_{k\in\delta_i^-} m_{ki} x^s_{ki} = b_i \quad \text{for all } i \in \mathcal{N}_0, \tag{13.81}$$

$$\sum_{j\in\delta_i^+\cap\mathcal{N}_0} x^s_{ij} - \sum_{k\in\delta_i^-\cap\mathcal{N}_0} m_{ki} x^s_{ki} + \sum_{j\in\delta_i^+\cap\mathcal{N}_1} y^s_{ij} - \sum_{k\in\delta_i^-\cap\mathcal{N}_1} m^s_{ki} y^s_{ki} = r^s_i \quad \text{for all } i \in \mathcal{N}_1, \tag{13.82}$$

$$0 \le x^s_{ij} \le u_{ij} \quad \text{for all } (i,j)\in\mathcal{A}_0, \tag{13.83}$$

$$0 \le y^s_{ij} \le v^s_{ij} \quad \text{for all } (i,j)\in\mathcal{A}_1. \tag{13.84}$$

The complete stochastic network problem (13.69)–(13.72) is obtained by replicating the network problem (13.80)–(13.84) for each scenario and including the non-anticipativity constraints

$$x^1_{ij} - x^s_{ij} = 0 \text{ for all } s \in \Omega \text{ and for all } (i,j) \in \mathcal{A}_0. \tag{13.85}$$

In this section we have been referring to quantities pertaining to an arc $(i,j) \in \mathcal{A}$ under scenario $s \in \Omega$ by using subscripts $(i,j)$ and a superscript $s$, respectively. To establish the correspondence between the matrix/vector notation of (13.77)–(13.78) and the component-wise notation of this section we impose a lexicographic order (see an example of such an order on page 362) on the arcs in $\mathcal{A}$ and let $(i_1, j_1)$ denote the first arc in $\mathcal{A}$. Then $z_1$, the first component of $z$ in (13.76), and $x^1_{i_1 j_1}$ refer to the same variable and so on.

## 13.8 Iterative Algorithm for Stochastic Network Optimization

In this section we develop an algorithm for stochastic network problems with a quadratic objective function using the split-variable formulation of the problem (Section 13.7.1). The algorithm is a specialization of the general row-action algorithm (Section 6.2) to the network structure. Since the algorithm works with the replicated problem it is convenient to use $x^s_{ij}$ to denote both first- and second-stage variables. That is, $x^s_{ij}$ for $(i,j) \in \mathcal{A}_0$ is a replicated first-stage variable, while $x^s_{ij}$ for $(i,j) \in \mathcal{A}_1$ is a second-stage variable. The objective function $F$ takes the form

$$F(x^s) = \sum_{(i,j)\in\mathcal{A},\ s\in\Omega} \frac{1}{2} w^s_{ij} (x^s_{ij})^2 + c^s_{ij} x^s_{ij}. \tag{13.86}$$

where $w_{ij}^s > 0$ and $c_{ij}^s$ are constants which are obtained by evaluating the summations in the minimand of (13.80) when the functions $f_{ij}$ and $q_{ij}^s$ are quadratic, and then adding up over all scenarios to get the expected value of the objective function values.

Let $M_1 \doteq S(m_0 + m_1)$. Then rows $1, 2, \ldots, M_1$ of the constraints matrix $\Phi$ (cf. equation (13.78)) correspond to network flow conservation constraints, and rows $M_1 + 1, \ldots, M$ correspond to the non-anticipativity constraints that take the simple form

$$x_{ij}^1 - x_{ij}^s = 0,$$

for all $(i,j) \in \mathcal{A}_0$ and all $s \in \Omega$. The dual price $\pi_\ell$, $\ell \in \{1, 2, \ldots, M_1\}$, associated with the flow conservation constraint for node $i \in \mathcal{N}$ under scenario $s \in \Omega$, is denoted by $\pi_i^s$. The dual price $\pi_\ell$, $\ell \in \{M+1, \ldots, M+N\}$, associated with the simple bound constraints for $x_{ij}^s$ (i.e., the reduced cost of $x_{ij}^s$), is denoted by $\pi_{ij}^s$. We now develop the specific projection formulae for use in the iterative steps of Algorithm 6.4.1. The row-action algorithm iterates one row at a time on the constraint matrix $\Phi$. The precise formulae for the iterative step are obtained by using the equation corresponding to the chosen matrix row, as given by equations (13.81)–(13.82). We develop the equations for only a single iterative step of the algorithm over all constraints, i.e., flow conservation equality constraints, bounds on the variables and non-anticipativity equality constraints. The complete algorithm is summarized below as Algorithm 13.8.1.

### Projection on Flow Conservation Constraints

First we derive the iterative step of the algorithm when the chosen row of the matrix $\Phi$ corresponds to the flow conservation constraint (13.81). Consider the flows on the incoming arcs, $x_{ki}^s$ for $k \in \delta_i^-$, and the flows on the outgoing arcs, $x_{ij}^s$ for $j \in \delta_i^+$ for a given node $i \in \mathcal{N}_0$ under some scenario $s \in \Omega$.

The generalized projection $\hat{z}$ of the current iterate $z$ onto the hyperplane $H(\phi^i, \gamma_i)$ where $\phi^i$ is the $i$th column of $\Phi^T$, and $\gamma_i$ is the $i$th component of $\gamma$ (see (13.77)–(13.78)), determined by the flow conservation constraint at node $i$, is obtained by solving the system (see Lemma 2.2.1):

$$\nabla F(\hat{z}) = \nabla F(z) + \beta_i^s \phi^i, \tag{13.87}$$

$$\hat{z} \in H(\phi^i, \gamma_i). \tag{13.88}$$

Of course, if $z \in H(\phi^i, \gamma_i)$, then $\beta_i^s = 0$ and $\hat{z} = z$. If the current iterate $z$ does not satisfy flow conservation at the $i$th node define the *node surplus* $\sigma_i^s$ as

$$\sigma_i^s \doteq b_i - \Big(\sum_{j \in \delta_i^+} x_{ij}^s - \sum_{k \in \delta_i^-} m_{ki} x_{ki}^s\Big). \tag{13.89}$$

Applying the iterative step (12.33) to the functional form of the objective function (13.86) and using the structure of the rows of the constraints matrix $\Phi$ that correspond to network flow constraints we get

$$\hat{x}^s_{ij} = x^s_{ij} + \beta^s_i \frac{1}{w^s_{ij}} \quad \text{for } j \in \delta^+_i, \tag{13.90}$$

$$\hat{x}^s_{ki} = x^s_{ki} - \beta^s_i \frac{m_{ki}}{w^s_{ki}} \quad \text{for } k \in \delta^-_i. \tag{13.91}$$

Substituting these expressions for $\hat{x}^s_{ij}$ and $\hat{x}^s_{ki}$ into (13.81) we get

$$\sum_{j\in\delta^+_i} (x^s_{ij} + \beta^s_i \frac{1}{w^s_{ij}}) - \sum_{k\in\delta^-_i} m_{ki}(x^s_{ki} - \beta^s_i \frac{m_{ki}}{w^s_{ki}}) = b_i. \tag{13.92}$$

From this and (13.89) we obtain

$$\beta^s_i = \frac{\sigma^s_i}{\displaystyle\sum_{j\in\delta^+_i} \frac{1}{w^s_{ij}} + \sum_{k\in\delta^-_i} \frac{m^2_{ki}}{w^s_{ki}}}. \tag{13.93}$$

Using this result in (13.90) and (13.91) gives the desired formulae for updating all primal variables incident to node $i$. The dual variable for this node $\pi^s_i$ is updated by adding $\beta^s_i$ to its current value, i.e., $\pi^s_i \leftarrow \pi^s_i + \beta^s_i$. Similar algebraic manipulations lead to the updating formulae required to apply the row-action algorithm to rows corresponding to constraints (13.82).

### Projection on Simple Bound Constraints

Now we develop the specific projections on the simple bounds (13.83)–(13.84) for the general row-action Algorithm 6.4.1. We do it in detail for (13.83) only. Denote by $x^s_{ij}$ the current value of the variable and by $\hat{x}^s_{ij}$ the projected value. If $x^s_{ij} < 0$, we get from (6.32) that

$$0 = \hat{x}^s_{ij} = x^s_{ij} + \frac{\beta}{w^s_{ij}}. \tag{13.94}$$

The primal variable is set to zero and the projection parameter is $\beta = -w^s_{ij}x^s_{ij}$. The dual price of the constraint is updated by subtracting $\beta$ from its current value, i.e., $\pi^s_{ij} \leftarrow \pi^s_{ij} - \beta$.

If $x^s_{ij} > u_{ij}$ we similarly set the primal variable $\hat{x}^s_{ij}$ to the upper bound $u_{ij}$, and from (6.32) compute the projection parameter $\beta = (u_{ij} - x^s_{ij})w^s_{ij}$, and update the dual price of the bound constraint by subtracting $\beta$.

Finally, if $0 \le x^s_{ij} \le u_{ij}$ we get from (6.61)

$$\hat{x}_{ij}^s = x_{ij}^s + \frac{\pi_{ij}^s}{w_{ij}^s} \tag{13.95}$$

and then set the dual price $\pi_{ij}^s$ to zero.

### Projections on Non-anticipativity Constraints

We now develop the iterative step of the algorithm for the equality, non-anticipativity constraints. A non-anticipativity constraint (13.85) has the form

$$x_{ij}^1 - x_{ij}^s = 0, \tag{13.96}$$

for some $(i,j) \in \mathcal{A}_0$ and some $s \in \Omega$. Let $\mu(s)$ be a row index such that $\phi^{\mu(s)}$ for $s = 2, 3, \ldots, S$ is the row of $\Phi$ that corresponds to the constraint $x_{ij}^1 - x_{ij}^s = 0$.

If $z$ is the current iterate, the generalized projection onto the hyperplane represented by this constraint is the point $\hat{z}$ which solves the system

$$\nabla F(\hat{z}) = \nabla F(z) + \beta \phi^{\mu(s)}, \tag{13.97}$$

$$x_{ij}^1 = x_{ij}^s. \tag{13.98}$$

Noting that the $\mu(s)$th row of the constraints matrix has only two nonzero entries (cf. equation (13.96)) we can write this system as

$$\begin{aligned}
\hat{x}_{ij}^1 &= x_{ij}^1 + \frac{\beta}{w_{ij}^1}, \\
\hat{x}_{ij}^s &= x_{ij}^s - \frac{\beta}{w_{ij}^s}, \\
\hat{x}_{ij}^1 &= x_{ij}^s.
\end{aligned}$$

Solving this, we get

$$\hat{x}_{ij}^1 = x_{ij}^s = \frac{w_{ij}^1 x_{ij}^1 + w_{ij}^s x_{ij}^s}{w_{ij}^1 + w_{ij}^s}, \tag{13.99}$$

i.e., the point $(x_{ij}^1, x_{ij}^s)$ is projected upon the point with coordinates equal to the weighted average of $x_{ij}^1$ and $x_{ij}^s$, with $w_{ij}^1$ and $w_{ij}^s$ being the weights.

Consider now the effect of repeated projections of the row-action algorithm on the non-anticipativity constraints (13.85). We can take advantage of the almost cyclic control of the algorithm in a way that would not have been possible with the cyclic control alone. The almost cyclic control of the row-action algorithm allows repeated projections upon these constraints alone until convergence—within some tolerance—of the variables $x_{ij}^s$, for any fixed $(i,j) \in \mathcal{A}_0$ to a limit $x_{ij}^*$. We show that $x_{ij}^*$ can be obtained analytically, rather than using the iterative scheme. This result has important implications for implementations, since the effect of repeated application of equation (13.99) for all $s \in \Omega$ can then be calculated in closed form.

The non-anticipativity constraints for the replications of a single first-stage variable $x_{ij}$ take the form

$$\begin{aligned} x_{ij}^1 - x_{ij}^2 &= 0, \\ x_{ij}^1 - x_{ij}^3 &= 0, \\ &\vdots \\ x_{ij}^1 - x_{ij}^S &= 0. \end{aligned} \tag{13.100}$$

Let $\nabla F_{ij} : \mathbb{R}^S \to \mathbb{R}^S$ denote the subvector of the gradient $\nabla F$ corresponding to the $S$ replications of the first-stage variable $x_{ij}^1, \ldots, x_{ij}^S$ and, similarly, let $\Phi_{(ij)}$ denote the submatrix of $\Phi$ consisting of the columns corresponding to $x_{ij}^1, \ldots, x_{ij}^S$.

By repeated projection onto these non-anticipativity constraints, such that the $\nu$th projection is onto the constraint $x_{ij}^1 - x_{ij}^{\ell(\nu)} = 0$, we obtain a sequence of points $x^\nu \in \mathbb{R}^S$ satisfying

$$\nabla F_{ij}(x^\nu) = \nabla F_{ij}(y) + \sum_{k=1}^{\nu} \lambda_k \phi_{(ij)}^{\mu(\ell(k))}, \tag{13.101}$$

where $\phi_{(ij)}^{\mu(\ell(k)}$ is the column of matrix $\Phi_{(ij)}$ corresponding to the $\mu(\ell(k))$th variable, $\lambda_k$ is the projection parameter corresponding to the $k$th projection, $y$ is the starting point, and $\mu(\ell)$ is the row index of the non-anticipativity constraints, see the discussion on page 408. The limit point $x^* \in \mathbb{R}^S$ satisfies

$$\nabla F_{ij}(x^*) = \nabla F_{ij}(y) + \sum_{k=1}^{\infty} \lambda_k \phi_{(ij)}^{\mu(\ell(k))} \tag{13.102}$$

and must, by the non-anticipativity constraints (13.100), have all components identical, i.e., $x^* = (x_{ij}^*, \ldots, x_{ij}^*)^T$ for some $x_{ij}^* \in \mathbb{R}$. Let now $\Lambda_s \doteq \sum_{\{k \mid \ell(k)=s\}} \lambda_k$, for $s = 2, \ldots, S$. Using the fact that $F(y)$ is the quadratic function (13.86) rewrite (13.102) as the system in $S$ variables $x_{ij}^*, \Lambda_2, \ldots, \Lambda_S$:

$$\begin{aligned} x_{ij}^* &= y_{ij}^1 + \frac{1}{w_{ij}^1} \sum_{s=2}^{S} \Lambda_s, \\ x_{ij}^* &= y_{ij}^2 - \frac{1}{w_{ij}^2} \Lambda_2, \\ &\vdots \\ x_{ij}^* &= y_{ij}^S - \frac{1}{w_{ij}^S} \Lambda_S. \end{aligned} \tag{13.103}$$

In matrix form, this is

$$Ht = y, \tag{13.104}$$

where

$$H \doteq \begin{pmatrix} 1 & \frac{-1}{w_{ij}^1} & \frac{-1}{w_{ij}^1} & \cdots & \frac{-1}{w_{ij}^1} \\ 1 & \frac{1}{w_{ij}^2} & & & \\ 1 & & \frac{1}{w_{ij}^3} & & \\ \vdots & & & \ddots & \\ 1 & & & & \frac{1}{w_{ij}^S} \end{pmatrix}, \tag{13.105}$$

and $t \doteq (x_{ij}^*, \Lambda_2, \dots, \Lambda_S)^{\mathrm{T}}$. By inverting $H$ we can solve for $t$. Since we are only interested in $x_{ij}^*$ (not in $\Lambda_2, \dots, \Lambda_S$), we need only calculate the first row of $H^{-1}$ denoted by $h = (h_s)_{s=1}^S$. Using the special structure of $H$ we easily get

$$h = \frac{1}{\det H}\left(\prod_{s=1}^{S} \frac{1}{w_{ij}^s}(w_{ij}^1, w_{ij}^2, \dots, w_{ij}^S)^{\mathrm{T}}\right),$$

where det $H$ is the determinant of $H$. The inner product of the first column of $H$, which consists of all ones, and the first row of $H^{-1}$ must equal 1. Therefore $\sum_{s=1}^S h_s = 1$. Hence

$$\det H = \left(\prod_{s=1}^{S} \frac{1}{w_{ij}^s}\right)\left(\sum_{s=1}^{S} w_{ij}^s\right).$$

Note that det $H > 0$, so that the system (13.104) has a unique solution. Solving for $x_{ij}^*$ we get

$$x_{ij}^* = \langle h, y\rangle = \frac{\sum_{s=1}^{S} w_{ij}^s y_{ij}^s}{\sum_{s=1}^{S} w_{ij}^s}. \tag{13.106}$$

Since $x$ is a first-stage variable

$$\frac{w_{ij}^1}{p_1} = \frac{w_{ij}^2}{p_2} = \dots = \frac{w_{ij}^S}{p_s}.$$

Also, $\sum_{s=1}^S p_s = 1$, so the result can be simplified to

$$x_{ij}^* = \sum_{s=1}^{S} p_s y_{ij}^s. \tag{13.107}$$

### The Row-action Algorithm for Quadratic Stochastic Networks

We have now completed all the components required to specialize the row-action algorithm to quadratic stochastic network problems. The complete algorithm proceeds as follows.

### Algorithm 13.8.1 Row-Action Algorithm for Quadratic Stochastic Networks

**Step 0:** (Initialization.) Set $\nu = 0$ and get $\pi^0$ and $z^0$ such that

$$\nabla F(z^0) = -\begin{pmatrix} \Phi \\ I_N \end{pmatrix}^{\mathrm{T}} \pi^0.$$

For example, $\pi^0 = 0$ and

$$(x_{ij}^s)^0 = -\frac{c_{ij}^s}{w_{ij}^s} \quad \text{for all } (i,j) \in \mathcal{A}_0,\ s \in \Omega, \tag{13.108}$$

$$(y_{ij}^s)^0 = -\frac{c_{ij}^s}{w_{ij}^s} \quad \text{for all } (i,j) \in \mathcal{A}_1,\ s \in \Omega. \tag{13.109}$$

**Step 1:** (Iterative step for the scenario subproblems.) For all $s \in \Omega$, do the following:

**Step 1.1:** (Iterative step for the flow conservation constraints.) Let

$$(\beta_i^s)^{\nu+\frac{1}{2}} = \frac{\sigma_i^s}{\sum_{j \in \delta_i^+} \frac{1}{w_{ij}^s} + \sum_{k \in \delta_i^-} \frac{m_{ki}^2}{w_{ki}^s}} \quad \text{for all } i \in \mathcal{N}_0, \tag{13.110}$$

$$(\beta_i^s)^{\nu+\frac{1}{2}} = \frac{\sigma_i^s}{\sum_{j \in \delta_i^+} \frac{1}{w_{ij}^s} + \sum_{k \in \delta_i^-} \frac{(m_{ki}^s)^2}{w_{ki}^s}} \quad \text{for all } i \in \mathcal{N}_1. \tag{13.111}$$

For all first-stage nodes $i \in \mathcal{N}_0$:

$$(x_{ij}^s)^{\nu+\frac{1}{2}} = (x_{ij}^s)^{\nu} + (\beta_i^s)^{\nu+\frac{1}{2}} \frac{1}{w_{ij}^s} \quad \text{for all } j \in \delta_i^+, \tag{13.112}$$

$$(x_{ki}^s)^{\nu+\frac{1}{2}} = (x_{ki}^s)^{\nu} - (\beta_i^s)^{\nu+\frac{1}{2}} \frac{m_{ki}}{w_{ki}^s} \quad \text{for all } k \in \delta_i^-, \tag{13.113}$$

$$(\pi_i^s)^{\nu+1} = (\pi_i^s)^{\nu} - (\beta_i^s)^{\nu+\frac{1}{2}}. \tag{13.114}$$

For all second-stage nodes $i \in \mathcal{N}_1$:

$$(y_{ij}^s)^{\nu+\frac{1}{2}} = (y_{ij}^s)^{\nu} + (\beta_i^s)^{\nu+\frac{1}{2}} \frac{1}{w_{ij}^s} \quad \text{for all } j \in \delta_i^+, \tag{13.115}$$

$$(y_{ki}^s)^{\nu+\frac{1}{2}} = (y_{ki}^s)^{\nu} - (\beta_i^s)^{\nu+\frac{1}{2}} \frac{m_{ki}^s}{w_{ki}^s} \quad \text{for all } k \in \delta_i^-, \tag{13.116}$$

$$(\pi_i^s)^{\nu+1} = (\pi_i^s)^{\nu} - (\beta_i^s)^{\nu+\frac{1}{2}}. \tag{13.117}$$

**Step 1.2:** (Iterative step for the simple bounds.)
For all first-stage arcs $(i,j) \in \mathcal{A}_0$:

$$(x_{ij}^s)^{\nu+1} = \begin{cases} u_{ij} & \text{if } (x_{ij}^s)^{\nu+\frac{1}{2}} \geq u_{ij}, \\ 0 & \text{if } (x_{ij}^s)^{\nu+\frac{1}{2}} \leq 0, \\ (x_{ij}^s)^{\nu+\frac{1}{2}} + \dfrac{(\pi_{ij}^s)^{\nu}}{w_{ij}^s} & \text{if } 0 < (x_{ij}^s)^{\nu+\frac{1}{2}} < u_{ij}. \end{cases} \tag{13.118}$$

and

$$(\pi_{ij}^s)^{\nu+1} = \begin{cases} (\pi_{ij}^s)^{\nu} - w_{ij}^s(u_{ij} - (x_{ij}^s)^{\nu+\frac{1}{2}}) & \text{if } (x_{ij}^s)^{\nu+\frac{1}{2}} \geq u_{ij}, \\ (\pi_{ij}^s)^{\nu} + w_{ij}^s(x_{ij}^s)^{\nu+\frac{1}{2}} & \text{if } (x_{ij}^s)^{\nu+\frac{1}{2}} \leq 0, \\ 0 & \text{if } 0 < (x_{ij}^s)^{\nu+\frac{1}{2}} < u_{ij}. \end{cases} \tag{13.119}$$

For all second-stage arcs $(i,j) \in \mathcal{A}_1$:

$$(y_{ij}^s)^{\nu+1} = \begin{cases} v_{ij}^s & \text{if } (y_{ij}^s)^{\nu+\frac{1}{2}} \geq v_{ij}^s, \\ 0 & \text{if } (y_{ij}^s)^{\nu+\frac{1}{2}} \leq 0, \\ (y_{ij}^s)^{\nu+\frac{1}{2}} + \dfrac{(\pi_{ij}^s)^{\nu}}{w_{ij}^s} & \text{if } 0 < (y_{ij}^s)^{\nu+\frac{1}{2}} < v_{ij}^s. \end{cases} \tag{13.120}$$

and

$$(\pi_{ij}^s)^{\nu+1} = \begin{cases} (\pi_{ij}^s)^{\nu} - w_{ij}^s(v_{ij}^s - (y_{ij}^s)^{\nu+\frac{1}{2}}) & \text{if } (y_{ij}^s)^{\nu+\frac{1}{2}} \geq v_{ij}^s, \\ (\pi_{ij}^s)^{\nu+\frac{1}{2}} + w_{ij}^s(y_{ij}^s)^{\nu+\frac{1}{2}} & \text{if } (y_{ij}^s)^{\nu+\frac{1}{2}} \leq 0, \\ 0 & \text{if } 0 < (y_{ij}^s)^{\nu+\frac{1}{2}} < v_{ij}^s. \end{cases} \tag{13.121}$$

**Step 2:** (Iterative step for non-anticipativity constraints.)
For all first-stage arcs $(i,j) \in \mathcal{A}_0$ set:

$$x_{ij}^* = \sum_{s=1}^{S} p_s (x_{ij}^s)^{\nu+1}, \tag{13.122}$$

$$(x_{ij}^s)^{\nu+1} = x_{ij}^* \text{ for all } s \in \Omega. \qquad (13.123)$$

**Step 3:** Let $\nu \leftarrow \nu + 1$ and return to Step 1.

**Decompositions for Parallel Computing**

The calculations in Step 1 of Algorithm 13.8.1 are repeated for multiple independent scenarios. Hence, these calculations can be executed concurrently utilizing as many processors as number of scenarios.

The calculations in Step 1.1, for a given scenario index $s$, are performed for all first- and second-stage nodes. On first examination these calculations are not independent, since nodes have arcs in common and the calculations for a given node, say $i$, cannot change the flow on arc $(i, j)$ at the same time that the calculations for node $j$ are updating the flow on this arc. However, it is possible to execute the calculations for multiple nodes concurrently if we identify nodes that do not have arcs in common. Such sets of nodes are identified by *coloring* the underlying graph, and iterating on same color nodes simultaneously. For example, in a time-staged network it is possible to iterate concurrently on all nodes corresponding to odd-order time periods, and then iterate concurrently on all nodes in even-order time periods. Another alternative is to employ a Jacobi variant of the algorithm described in Step 1.1. (see Section 13.9 for references to related literature). The calculations in Step 1.2 can be executed concurrently for all arcs and all scenarios utilizing as many processors as the number of scenarios times the number of arcs.

The calculations in Step 2 can be executed concurrently for all first-stage arcs, although the calculations involved in this step are trivial compared to the amount of work performed in Step 1. Sections 14.5 and 15.5 give details on the implementation of this algorithm, and report computational results with large-scale test problems.

## 13.9 Notes and References

Stochastic programming models were first formulated as mathematical programs in the late 1950s by Dantzig (1955) and Beale (1955). Programs with probabilistic constraints were introduced by Charnes and Cooper (1959). For general references on stochastic programming see Dempster (1980), Ermoliev and Wets (1988), Kall (1976), Kall and Wallace (1994), and Wets (1989). A textbook treatment of stochastic programming is Kall and Wallace (1994)

**13.1** For further discussion on probability theory as it applies to stochastic programming see Kall (1976), Wets (1989), and Frauendorfer's thesis (1992). For general background on probability theory refer to Billingsley (1995) or Parzen (1960).

**13.3** The problem formulations can be found in the general references cited above. See also Walkup and Wets (1966) and Wets (1966a, 1966b, 1972, 1983). Wets (1974) develops the deterministic equivalent formulation. Multistage programs are discussed in Birge (1985, 1988), Olsen (1976), Gassmann (1990), Ermoliev and Wets (1988), and Wets (1989). Dupačova (1995) compiled a bibliography.

**13.4** The robust optimization model was suggested by Mulvey, Vanderbei, and Zenios (1995). The terminology of *structural* and *control* variables is borrowed from the flexibility analysis of manufacturing systems; see Seider, Brengel, and Widagdo (1991). Applications of robust optimization are discussed in Guttierez and Kouvelis (1995), King et al. (1988), Malcolm and Zenios (1994), Paraskevopoulos, Karakitsos, and Rustem (1991) and Sengupta (1991).

**13.5.1** The diet problem was studied by Stigler (1945) and used by Dantzig (1963) as the first test problem for the simplex method. See also Dantzig (1990).

**13.5.2** Capacity planning and expansion has been a fertile ground for the application of optimization models. For a textbook treatment, in the context of manufacturing applications, see Hayes and Wheelwright (1984). The robust optimization approach to capacity planning for a multiproduct, multifacility production firm was suggested by Eppen, Martin, and Schrage (1989), who applied their model to plan car manufacturing facilities for the General Motors Company. The model described in this section is a simplified version of their application. The same reference discusses the merits of a robust optimization formulation for the capacity expansion planing model. They use a model of *expected downside risk*, and illustrate its performance with numerical results. The Markowitz criterion was introduced by Markowitz (1952); see also Perold (1984) and Dahl, Meeraus, and Zenios (1993). The downside risk function used in this section was suggested by Zenios and Kang (1993) in the context of portfolio management applications. It is the limiting case of the *mean-absolute deviation* models of Sharpe (1971) and Konno and Yamazaki (1991) with asymmetric risk functions. See Speranza (1993) for further analysis of the properties of asymmetric, piecewise linear, penalty functions.

Robust optimization models for capacity expansion planning have been developed by Guttierez and Kouvelis (1995) who consider *outsourcing* (i.e., subcontracting part of the manufacturing requirements) as the means to achieving robustness in manufacturing capacity while reducing costs. Stochastic programming models for capacity expansion for power generation firms have been proposed by Murphy, Sen, and Soyster (1982), Granville et al. (1988), and Dantzig et al. (1989). Extensions of these models, using robust optimization, are developed by Malcolm and Zenios (1994).

**13.5.3** The robust optimization model for matrix balancing was developed by Zenios and Zenios (1992). A similar model, for the more general problem of image reconstruction from projections, was suggested earlier by Elfving (1989), who also suggested the use of a single Newton step for the estimation of the Bregman parameter. Both references develop solution algorithms for the respective models, using the row-action algorithms of Chapter 6. The use of a secant approximation for the estimation of the Bregman parameter is discussed in Section 6.9; see references in Section 6.10.

**13.6** Textbook treatment of investments and portfolio management are given by Bodie, Kane, and Marcus (1989) and Elton and Gruber (1984). The classic models for portfolio management, namely Markowitz's mean/variance model and the portfolio immunization model are discussed by Markowitz (1952) and Reddington (1952) respectively. See Dahl, Meeraus, and Zenios (1993) for a modern treatment of these models and the associated mathematical programming formulations.

The applications of stochastic programming models in portfolio management are numerous. For a simple illustration of dynamic programming formulations for multiperiod portfolio optimization see Bertsekas (1987, pp. 73–77). General references are collected in Zenios (1993a) and Ziemba and Vickson (1975). The application of stochastic programming to address problems in short-term financial planning was suggested by Kallberg, White, and Ziemba (1982). Its application to bond portfolio management was suggested by Bradley and Crane (1972); for application to Bank asset/liability management see Kusy and Ziemba (1986); for application to the asset allocation see Mulvey and Vladimirou (1992) and Mulvey (1993); for application to fixed-income portfolio management see Zenios (1991c, 1993b); and for application to funding of insurance products see Nielsen and Zenios (1996b). The model discussed in this section is adapted from Golub et al. (1995) and Holmer et al. (1993).

**13.7** There exists extensive literature on stochastic programming problems with network recourse. The stochastic transportation problem was introduced by Williams (1963). See also the papers by Cooper and LeBlanc (1977), Wallace (1986, 1987), Mulvey and Vladimirou (1991), Nielsen and Zenios (1993c, 1996a), and the PhD theses by Vladimirou (1990) and Nielsen (1992).

Applications of stochastic network models to hydroelectric power scheduling are reported in Dembo et al. (1990) and Dembo, Mulvey, and Zenios (1989). Applications to airtraffic control were developed by Richetta (1991); see also Zenios (1991a). Transportation and logistics models are developed in Frantzeskakis and Powell (1989).

The split-variable formulation for the decomposition of mathematical programs is fairly standard in large-scale optimization; see, for example, Bertsekas and Tsitsiklis (1989, page 231). The use of split-variable formulations for stochastic programming problems was suggested by Rockafellar and Wets (1991). It was used by Mulvey and Vladimirou (1989, 1991) and Nielsen and Zenios (1993a, 1996a) as a device for exploiting the special structure of stochastic programs with network recourse.

**13.8** The row-action iterative algorithm for the two-stage stochastic programming problems with network recourse was developed by Nielsen and Zenios (1993a). It was further extended to the multistage problem by Nielsen and Zenios (1996a), and was used for the solution of linear stochastic networks within the context of the PMD algorithm (Chapter 3) by Nielsen and Zenios (1993d). The partitioning of network structures using graph coloring to increase the amount of parallelism was suggested in Zenios and Mulvey (1988a). They describe a graph coloring heuristic, based on earlier work of Christofides (1971), used for their test problems, and compare the idea of graph coloring to the use of a Jacobi variant. As one may expect, the parallel scheme based on graph coloring exhibits a faster (practical) convergence rate than the Jacobi algorithm. However, the Jacobi algorithm permits the use of more processors. For large-scale problems, and on computers with sufficiently many processors, the Jacobi parallel algorithm could be substantially faster, in solution time, than parallel implementations based on graph coloring. This was demonstrated in Zenios and Lasken (1988a, 1988b).

CHAPTER 14

# Decompositions for Parallel Computing

This chapter examines the parallel implementation of several of the algorithms developed earlier whose structure—often in conjunction with that of the application—makes them suitable for decomposition into independent tasks. The decomposition is sometimes facilitated by the structure of the mathematical algorithm; simultaneous and block-iterative algorithms (as characterized in Section 1.3) are natural candidates for parallel computations. Sometimes the decomposition is facilitated by the structure of the problem, such as in image reconstruction where we can partition a discretized image into domains that are reconstructed independent of each other.

In order to implement a parallel algorithm we must first partition its operations, or the problem data, into independent tasks that can be mapped on multiple processors (see also Section 1.2). That is, we must partition the operations of the algorithm into sets of independent operations that can proceed concurrently, or we must partition the problem data into blocks of data without interdependencies, so that multiple blocks can be processed in parallel. This is called *task partitioning.* Once a partition has been determined the tasks must be assigned to one or more processors for execution; this is known as the problem of *task scheduling.* Finally it may be necessary to impose some ordering to task execution and specify when information must be exchanged among the processors performing different tasks in order to guarantee correctness of the results. This is known as *task synchronization.*

Task partitioning is dependent primarily on the algorithm or the model. It can be developed with little concern for the targeted computer architecture, whereas the scheduling and synchronization of tasks depend on it. The organization of memory and the communication network are key factors for efficient task scheduling and synchronization. In spite of their importance however, we will not deal with the issues of task scheduling and synchronization, because they cannot be addressed independent of the architecture. Furthermore, it is often—but not always—the case that these issues can be resolved by the compiler, by the programming language, or by calls to a system library. Kumar et al. (1994) discuss implementation of several algorithms (such as matrix operations, sorting, dynamic pro-gramming, graph operations) in conjunction with the underlying computer architecture. Bertsekas and Tsitsiklis (1989) discuss the communication and synchronization requirements of similar algorithms in conjunction with some popular communication networks.

In this chapter we concentrate on the problem of *task partitioning.* The partitions refer to the operations of an algorithm or the data of a problem—but not to the hardware. It is assumed that enough processors are available so that each task can be assigned to a different processor. If fewer processors are available then each one will execute more than one task.

The chapter is organized as follows: Section 14.1 describes an abstract parallel machine and related operations which are then used to describe the task partitionings. Section 14.2 describes the partitioning of simple structures (dense and sparse matrices) and their mapping on multiple processors. Sections 14.3–14.6 discuss the implementation of parallel algorithms for the solution of real-world problems discussed in previous chapters. The algorithms we focus on include algorithms for matrix balancing, image reconstruction, network optimization, and interior point algorithms. Further information and references are given in Section 14.7.

## 14.1 Vector-Random Access Machine (V-RAM)

We describe now an abstract model of a data parallel computer, the *Vector-Random Access Machine* (V-RAM). (Refer to Section 1.1.3 for an introduction to the data parallel programming paradigm.) A RAM is a machine whose memory can be accessed in constant time which does not depend on the address of the word being accessed. The V-RAM is a standard RAM with the addition of *vector memory* and a parallel *vector processor* (see Figure 14.1). The vector memory is a sequence of addresses containing linearly ordered collections of scalar values, known as *simple vectors.* An important and distinctive feature of a V-RAM is that the vector lengths need not be identical (in contrast to vector computers). The vector processor executes vector instructions (in parallel) on sets of simple vectors and scalars stored in vector and scalar memories. For example, the SAXPY operation $Y \leftarrow \alpha X + Y$ corresponds to the multiplication of a vector $X = (X(j))_{j=1}^{n}$ in the vector memory by a scalar element $\alpha$ from the scalar memory, and the addition of the result to the vector $Y = (Y(j))_{j=1}^{n}$ in the vector memory. A data parallel realization of a V-RAM maps the parallel vector processor and the vector memory onto multiple processors. All addresses of the vector processor and their vector memory address(es) are mapped on separate processors. A processor, in this context, is a physical processing element (i.e., a central processing unit, CPU). However, if sufficient processing elements are unavailable, it is still possible to map multiple addresses of the vector memory and vector processor onto the same processing element. We then refer to the ensemble of processing element with each vector memory and vector address assigned to it as a *virtual processor* (VP). Each processing element thus emulates multiple VPs.

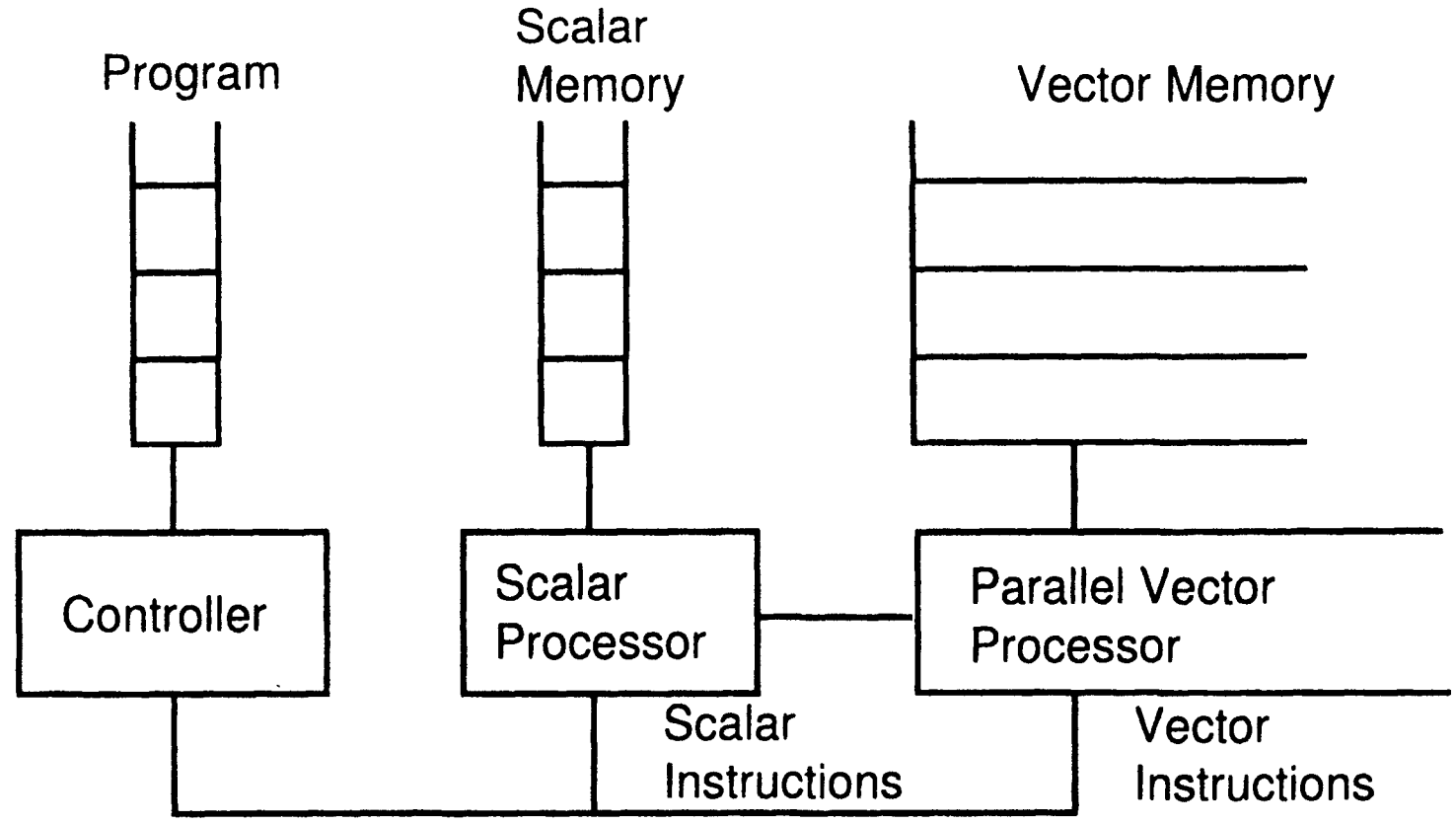

**Figure 14.1** The Vector–Random Access Machine (V-RAM) model for data parallel computing.

The set of vector instructions of a V-RAM can be used to describe the implementation of several algorithms. The level of parallelism is proportional to the amount of data in the problem—or, at least, proportionate to the length of the parallel vector processor of the specific data parallel machine. We review next the vector instructions of a V-RAM.

### 14.1.1 Parallel prefix operations

There are two classes of vector instructions: *scans*, or parallel prefix operations, and *segmented scans*, which are used for the data parallel implementation of algorithms on a V-RAM. The specific instructions relevant to this chapter are *scans* and *spreads.* The $\otimes$-*scan* primitive, for an associative binary operator $\otimes$, takes a sequence of real numbers $\{x_0, x_1, \dots, x_n\}$ and produces another sequence $\{y_0, y_1, \dots, y_n\}$ such that $y_i = x_0 \otimes x_1 \otimes \cdots \otimes x_i$. For example, *add-scan* takes as an argument a parallel variable (i.e., a variable with its $i$th element residing in the $i$th address of a simple vector) and returns at the $i$th address the value of the parallel variable summed over all $j = 0, 1, \dots, i$. A scan can only be applied to preceeding processors (e.g., sum over $j = 0, 1, \dots, i-1$) or it can be performed in reverse. The $\otimes$-*spread* primitive, for an associative binary operator $\otimes$, takes a sequence of real numbers $\{x_0, x_1, \dots, x_n\}$ and produces another sequence $\{y_0, y_1, \dots, y_n\}$ such that $y_i = x_0 \otimes x_1 \otimes \cdots \otimes x_n$. For example, *add-spread* takes as argument the elements of a parallel variable and returns at all addresses of the simple vector the sum of the parallel variable summed over $j = 0, 1, \dots, n$.

Another variation of the scan primitives allows their operation within

| | | |
|---|---|---|
| Processing Element | = | { 0 1 2 3 4 5 6 7 8 9 } |
| Parallel Variable $X$ | = | { 5 1 3 4 3 9 2 6 1 0 } |
| Segment Bits ($Sb$) | = | { 1 0 1 0 0 0 1 0 0 1 } |
| $Y = \mathit{add\text{-}scan}(X, Sb)$ | = | { 0 5 6 3 7 10 19 2 8 9 } |
| $Y = \mathit{copy\text{-}scan}(Y, Sb)$ | = | { 0 0 0 6 6 6 6 19 19 19 } |

**Figure 14.2** An example of the segmented *add-scan* and *copy-scan* parallel prefix operations.

*segments* of a parallel variable. These are known as *segmented-*$\otimes$*-scan* instructions. They take as arguments a parallel variable and a set of *segment bits*, which are specified by a partition of the vector memory into multiple simple vectors of contiguous memory address. (Segment bits are set equal to 1 at the starting location of a new segment and 0 elsewhere.) A *segmented-*$\otimes$*-scan* operation restarts at the beginning of every segment. Figure 14.2 illustrates the use of the scan primitives *segmented-add-scan* and *segmented-copy-scan*.

## 14.2 Mapping Data to Processors

We consider now the partitioning of matrices, and the mapping of these partitions onto multiple VPs using the V-RAM model. To associate a V-RAM with some given data structure (such as a matrix) we must specify the number and size of simple vectors of the V-RAM that will store the data structure. We develop here the V-RAM model for both *dense* matrices with few or no zero entries, and *sparse* matrices where many of the entries are zero.

### 14.2.1 Mapping a dense matrix

A dense $n \times n$ matrix can be uniformly partitioned into smaller submatrices, and each one distributed to a different processor. It is possible to partition the matrix into $n^2$ submatrices and use as many VPs. If only $p$ processors are available the matrix is partitioned into blocks of sizes $n/\sqrt{p} \times n/\sqrt{p}$. This partitioning of a dense matrix into submatrices is known as *checkerboard partitioning*.

An alternative partitioning is obtained by dividing the matrix into blocks of complete rows or columns. All blocks contain an equal number of rows or columns, and each block is distributed to a different processor. This is known as *striped partitioning*. The maximum number of VPs that can be effectively utilized by striped partitioning of an $n \times n$ matrix is $n$.

### The V-RAM Model of a Dense Matrix

The V-RAM model of a dense matrix specifies the number and dimensions of simple vectors required to store the matrix. The precise specifications depend on the operations that are to be performed on the dense matrix. This is analogous to the choice of data structures in uniprocessor scalar programming, where the choice of data structures also depends on the operations that are performed. For example, the appropriate V-RAM model for determining the largest entry of a matrix is one with a single simple vector of length $n^2$.

Consider now an algorithm that must compute the row and column sums of an $n \times n$ matrix. (Such computations appear in the matrix balancing algorithms of Section 9.3, and in the transportation algorithms of Section 12.4.) The V-RAM model that facilitates the computation of row sums consists of $n$ simple vectors each of length $n$. The $i$th vector holds the entries of the $i$th row of the matrix. Denoting the entries of this simple vector by $\{a_{i1}, a_{i2}, \ldots, a_{in}\}$, the row sum is computed by using the operation *add-scan*$\{a_{i1}, a_{i2}, \ldots, a_{in}\}$.

This V-RAM model can be implemented on a parallel machine using either the checkerboard partitioning of the matrix if $n^2$ VPs are available, or the (row-wise) striped partitioning if $n$ VPs are available. Likewise, the V-RAM model that facilitates the computation of column sums consists of $n$ simple vectors each of length $n$. The $j$th vector holds the entries of the $j$th column of the matrix. If we denote the entries of this column by $\{a_{1j}, a_{2j}, \ldots, a_{nj}\}$ then the column sum is computed by using the operation *add-scan* $\{a_{1j}, a_{2j}, \ldots, a_{nj}\}$. This V-RAM model can also be implemented using either the checkerboard partitioning or the (column-wise) striped partitioning.

Both V-RAM models with the checkerboard data mapping can be implemented on a parallel machine with $n^2$ physical processors, interconnected by a two-dimensional *mesh communication network* whereby each processor is connected to its four neighbors. Figure 14.3 illustrates the checkerboard mapping of the dense matrix onto a mesh network of sixteen processors, denoted by $P_0, P_1, \ldots, P_{15}$, and the V-RAM models for computing row and column sums.

### 14.2.2 Mapping a sparse matrix

Consider now an $n \times n$ sparse matrix with $q$ nonzero entries. Mapping the matrix onto more than $q$ VPs is a waste of resources. We use a format for storing only the nonzero entries of the matrix, and then map this format onto multiple VPs.

The standard format for storing a sparse matrix (without detectable sparsity structure), is the *compressed row-wise* format. It uses an array, called VAL, of dimension $q$ that stores the nonzero entries, an array $J$ of length $q$ that contains pointers to the column index of each nonzero entry,

$$\begin{pmatrix} a_{00} & a_{01} & a_{02} & a_{03} & a_{04} & a_{05} & a_{06} & a_{07} \\ a_{10} & a_{11} & a_{12} & a_{13} & a_{14} & a_{15} & a_{16} & a_{17} \\ a_{20} & a_{21} & a_{22} & a_{23} & a_{24} & a_{25} & a_{26} & a_{27} \\ a_{30} & a_{31} & a_{32} & a_{33} & a_{34} & a_{35} & a_{36} & a_{37} \\ a_{40} & a_{41} & a_{42} & a_{43} & a_{44} & a_{45} & a_{46} & a_{47} \\ a_{50} & a_{51} & a_{52} & a_{53} & a_{54} & a_{55} & a_{56} & a_{57} \\ a_{60} & a_{61} & a_{62} & a_{63} & a_{64} & a_{65} & a_{66} & a_{67} \\ a_{70} & a_{71} & a_{72} & a_{73} & a_{74} & a_{75} & a_{76} & a_{77} \end{pmatrix}$$

Dense matrix $A$

| | | | |
|---|---|---|---|
| (0,0) (0,1)<br>$P_0$<br>(1,0) (1,1) | (0,2) (0,3)<br>$P_1$<br>(1,2) (1,3) | (0,4) (0,5)<br>$P_2$<br>(1,4) (1,5) | (0,6) (0,7)<br>$P_3$<br>(1,6) (1,7) |
| (2,0) (2,1)<br>$P_4$<br>(3,0) (3,1) | (2,2) (2,3)<br>$P_5$<br>(3,2) (3,3) | (2,4) (2,5)<br>$P_6$<br>(3,4) (3,5) | (2,6) (2,7)<br>$P_7$<br>(3,6) (3,7) |
| (4,0) (4,1)<br>$P_8$<br>(5,0) (5,1) | (4,2) (4,3)<br>$P_9$<br>(5,2) (5,3) | (4,4) (4,5)<br>$P_{10}$<br>(5,4) (5,5) | (4,6) (4,7)<br>$P_{11}$<br>(5,6) (5,7) |
| (6,0) (6,1)<br>$P_{12}$<br>(7,0) (7,1) | (6,2) (6,3)<br>$P_{13}$<br>(7,2) (7,3) | (6,4) (6,5)<br>$P_{14}$<br>(7,4) (7,5) | (6,6) (6,7)<br>$P_{15}$<br>(7,6) (7,7) |

Block-checkerboard partitioning using 16 processors and an array $I$ whose $i$th entry points to the starting address of the $i$th row of the matrix in arrays $VAL$ and $J$. Figure 14.4 illustrates the compressed row-wise format of a small matrix. The matrix can also be stored using a similar *compressed column-wise* format.

| | 0 | 1 | 2 | 3 | 4 | 5 | 6 | 7 |
|---|---|---|---|---|---|---|---|---|
| 0 | $a_{00}$ | $a_{01}$ | $a_{02}$ | $a_{03}$ | $a_{04}$ | $a_{05}$ | $a_{06}$ | $a_{07}$ |
| 1 | $a_{10}$ | $a_{11}$ | $a_{12}$ | $a_{13}$ | $a_{14}$ | $a_{15}$ | $a_{16}$ | $a_{17}$ |
| 2 | $a_{20}$ | $a_{21}$ | $a_{22}$ | $a_{23}$ | $a_{24}$ | $a_{25}$ | $a_{26}$ | $a_{27}$ |
| 3 | $a_{30}$ | $a_{31}$ | $a_{32}$ | $a_{33}$ | $a_{34}$ | $a_{35}$ | $a_{36}$ | $a_{37}$ |
| 4 | $a_{40}$ | $a_{41}$ | $a_{42}$ | $a_{43}$ | $a_{44}$ | $a_{45}$ | $a_{46}$ | $a_{47}$ |
| 5 | $a_{50}$ | $a_{51}$ | $a_{52}$ | $a_{53}$ | $a_{54}$ | $a_{55}$ | $a_{56}$ | $a_{57}$ |
| 6 | $a_{60}$ | $a_{61}$ | $a_{62}$ | $a_{63}$ | $a_{64}$ | $a_{65}$ | $a_{66}$ | $a_{67}$ |
| 7 | $a_{70}$ | $a_{71}$ | $a_{72}$ | $a_{73}$ | $a_{74}$ | $a_{75}$ | $a_{76}$ | $a_{77}$ |
| 0 | $a_{00}$ | $a_{10}$ | $a_{20}$ | $a_{30}$ | $a_{40}$ | $a_{50}$ | $a_{60}$ | $a_{70}$ |
| 1 | $a_{01}$ | $a_{11}$ | $a_{21}$ | $a_{31}$ | $a_{41}$ | $a_{51}$ | $a_{61}$ | $a_{71}$ |
| 2 | $a_{02}$ | $a_{12}$ | $a_{22}$ | $a_{32}$ | $a_{42}$ | $a_{52}$ | $a_{62}$ | $a_{72}$ |
| 3 | $a_{03}$ | $a_{13}$ | $a_{23}$ | $a_{33}$ | $a_{43}$ | $a_{53}$ | $a_{63}$ | $a_{73}$ |
| 4 | $a_{04}$ | $a_{14}$ | $a_{24}$ | $a_{34}$ | $a_{44}$ | $a_{54}$ | $a_{64}$ | $a_{74}$ |
| 5 | $a_{05}$ | $a_{15}$ | $a_{25}$ | $a_{35}$ | $a_{45}$ | $a_{55}$ | $a_{65}$ | $a_{75}$ |
| 6 | $a_{06}$ | $a_{16}$ | $a_{26}$ | $a_{36}$ | $a_{46}$ | $a_{56}$ | $a_{66}$ | $a_{76}$ |
| 7 | $a_{07}$ | $a_{17}$ | $a_{27}$ | $a_{37}$ | $a_{47}$ | $a_{57}$ | $a_{67}$ | $a_{77}$ |

**Figure 14.3** Checkerboard mapping of a dense matrix on to a mesh network of processors and V-RAM model for computing row and column sums.

The natural parallel mapping of a sparse matrix stored in compressed row-wise format is the row-wise striped partitioning: the matrix is partitioned into blocks of rows, and the compressed row-wise data structure of each block is mapped onto a different processor. This partitioning of a sparse matrix may lead to load imbalance among the processors, that is, processors may receive blocks with differing numbers of nonzero entries. An efficient partitioning will select rows to add to a block in such a way that all blocks have approximately the same number of nonzero entries. For example, if a matrix has a few dense rows then those should not be placed into the same block. A simple technique to ensure the formation of blocks of homogeneous size is to select randomly the rows assigned to each block. In this way it is likely that each block will be assigned both dense and sparse columns, so that the number of nonzero entries in each block will be approximately the same.

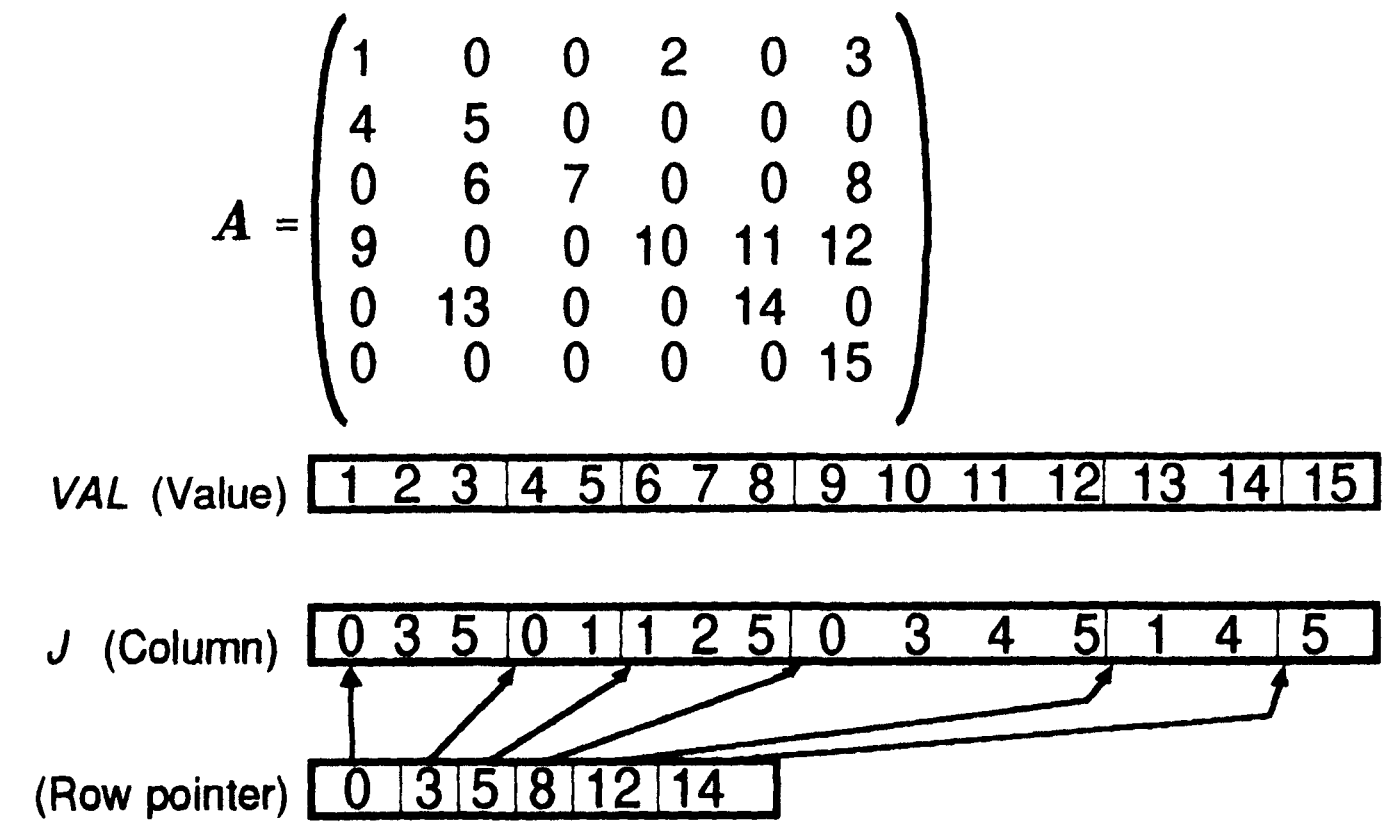

**Figure 14.4** Compressed row-wise representation of a sparse matrix $A$.

### The V-RAM Model of a Sparse Matrix

The V-RAM model of a sparse matrix specifies the number and dimensions of simple vectors needed to perform some computations on the matrix using a V-RAM. We consider the case of calculating the sums of the rows of the matrix. The appropriate V-RAM model consists of $n$ simple vectors, one for each row of the matrix. The length of each simple vector is equal to the number of nonzero entries of the respective row, and the vector stores the array $VAL$. (The calculation of the sum of the rows does not require knowledge of the column indices, stored in the array $J$. However, if these indices are needed for some calculations, then $n$ additional simple vectors will store the entries of $J$.) The array $I$ is redundant in this V-RAM model, since the starting address of each row is given implicitly by the order of the simple vectors: the first vector corresponds to the first row, the second vector corresponds to the second row, and so on.

This V-RAM model can be implemented on a linearly ordered array of $q$ VPs as follows. The array is separated into $n$ segments, with the length of the $i$th segment equal to the number of nonzero entries in the $i$th row, denoted by $q_i$. The VPs of each segment store the $VAL$ entries of the corresponding row, denoted by $\{a_{i1}, a_{i2}, \ldots, a_{iq_i}\}$. The seg-ments are indicated by a vector of length $q = \sum_{i=1}^{n} q_i$ of segment bits that has a 1 at locations $1, q_1, q_1 + q_2, \ldots, q$, and 0 elsewhere. The sums of the rows can be computed by a *segmented-add-scan* on the parallel variable $\{a_{11}, a_{12}, \ldots, a_{1q_1}, a_{21}, a_{22}, \ldots, a_{2q_2}, \ldots a_{n1}, a_{n2}, \ldots, a_{nq_n}\}$. A similar V-RAM model can be developed to map a compressed column-wise representation of the matrix on a V-RAM machine.

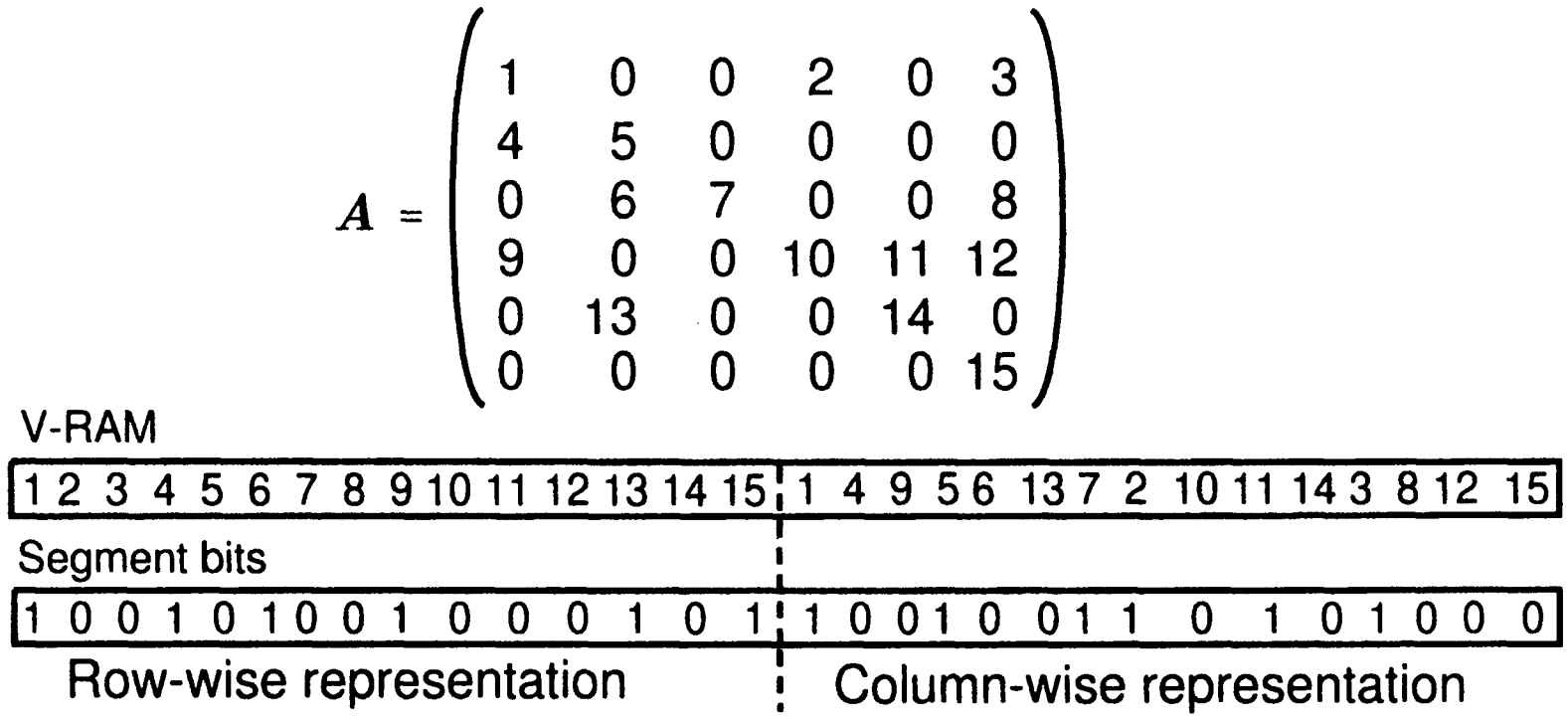

**Figure 14.5** The V-RAM model of a sparse matrix $A$ for simultaneous execution of both row and column operations.

Most applications that involve sparse matrices execute calculations on both rows and columns of the matrix (but not on just one). Such is the case for the matrix balancing algorithms of Section 9.3, and the transportation algorithms of Section 12.4. If row sums and column sums will not be computed simultaneously we need only to implement the V-RAM models for both row-wise and column-wise formats, and operate alternately between the two. A V-RAM model can also be developed to facilitate the simultaneous calculation of both row and column sums, in which case it utilizes a linearly ordered array of $2q$ processors partitioned into $2n$ segments. The length of the $i$th segment for $i = 1, 2, \ldots, n$ is equal to the number of nonzero entries in the $i$th row, and the length of the $j$th segment $j = n+1, n+2, \ldots, 2n$ is equal to the number of nonzero entries in the $(j-n)$th column. The segments of the first $q$ VPs are identical to the segments of the V-RAM model of the row-wise format, and the segments of the last $q$ VPs are identical to the segments of the column-wise format. A *segmented-add-scan* on these concatenated segments will simultaneously compute all row and column sums of the matrix. Figure 14.5 illustrates this V-RAM model.

## 14.3 Parallel Computing for Matrix Balancing

We now discuss parallel computing techniques for the implementation of algorithms for matrix balancing, focusing in particular on Problem 9.2.3, and alternative ways to parallelize the RAS algorithm (Algorithm 9.3.2). The application of RAS to this problem is well suited to both data parallel and control parallel implementations. The mapping of data to parallel processors which facilitates the parallel execution of RAS can also be used to implement in parallel the range-RAS algorithm (Algorithm 9.3.1).

### 14.3.1 Data parallel computing with RAS

We begin with a description of data parallel implementations. The implementation for solving dense problems is given first, followed by an implementation that exploits the sparsity which is typically present in large-scale problems.

#### Solving Dense Problems

The parallel execution of RAS for dense matrix balancing problems is performed using the data mapping of dense matrices (Section 14.2.1). In particular, the $m \times n$ matrix $X = (x_{ij})$ is mapped onto $mn$ VPs using a checkerboard partitioning.

The V-RAM model that we employ uses dense row-wise and column-wise models, with $2m$ simple vectors of length $n$, and $2n$ simple vectors of length $m$. The first $2m$ vectors store the entries $x_{ij}$ of the matrix in a row-wise partitioning; the remaining $2n$ vectors store the same entries in a column-wise partitioning. Recall that it is possible to implement these two models using the checkerboard data mapping on a mesh network of processors (see Figure 14.3). Hence, the representation of the matrix using both partitionings does not make redundant use of VPs.

The row totals $u_i$ are stored in $m$ vectors of length $n$ such that the $i$th vector stores $n$ copies of the value $u_i$. Although this uses redundantly the memory of the VPs, it is done to reduce the communication requirements of the implementation. Similarly, $n$ vectors of length $m$ store the column totals $v_j$, such that the $j$th vector stores $m$ copies of $v_j$

Step 1 of the RAS algorithm can be executed simultaneously for all rows of the matrix. An *add-spread* operation on the first $n$ simple vectors computes the partial sum of the entries of each row and communicates the totals to all columns in the corresponding row. Hence, all addresses of the $i$th row store the sum $\sum_{j=1}^{n} x_{ij}$. Each VP can now divide the row total $u_i$ by the computed partial sum to estimate the scaling factor $\rho_i$. Since $u_i$ is stored (redundantly) on all VPs corresponding to the $i$th row, this calculation is executed by multiple VPs simultaneously, without any need to communicate data. The scaling factor is then used to multiply the value of $x_{ij}$, which completes Step 1.

In a similar fashion Step 2 of the algorithm can be executed simultaneously for all columns of the matrix. An *add-spread* operation on the $n$ simple vectors that store the column-wise partitioning of the matrix computes the partial sum $\sum_{i=1}^{m} x_{ij}$, for each $j = 1, 2, \ldots, n$, which is then stored in all addresses of the $j$th vector. Each VP can now divide the column total $v_j$ by the computed partial sum to estimate the scaling factor $\sigma_j$. This factor is then used to multiply the value of $x_{ij}$, completing the execution of Step 2.

### Solving Sparse Problems

For execution of RAS for sparse problems we use the V-RAM models of Section 14.2.2, in particular the model with two simple vectors of length $q$, one for the compressed row-wise representation of the matrix and one for the compressed column-wise representation. An additional vector of length $q$ stores the row totals, with the segment corresponding to row $i$ storing $q_i$ copies of $u_i$. Similarly, a vector of length $q$ stores the column totals, with the segment corresponding to column $j$ storing $q_j$ copies of $v_j$. With this V-RAM model each matrix entry is stored twice: once in row-wise compressed format and once in column-wise compressed format. These entries should be, componentwise, identical. Hence, when the $(i, j)$th entry of the matrix is scaled in the simple vector that stores the matrix row-wise, its value should be communicated to the corresponding entry in the vector that stores the matrix column-wise. Two $q$-long vectors store pointers that establish the one-to-one mapping between the two formats.

The implementation of RAS for sparse problems is similar to that for dense problems. Step 1 of the algorithm is executed simultaneously for all rows of the matrix. A *segmented-add-spread* operation on the simple vector with the compressed row-wise format computes the partial sum of the entries of each row and communicates the sum to all addresses of the corresponding row. Hence, all addresses of the segment corresponding to the $i$th row have access to the sum $\sum_{j=1}^{n} x_{ij}$. Each VP divides the column total $u_i$, which is available locally, by the computed partial sum to estimate the scaling factor $\rho_i$. The scaling factor is then used to multiply the values of $x_{ij}$, completing the execution of Step 1. Before the algorithm can proceed to Step 2 the entries $x_{ij}$, just updated in the row-wise format, are communicated to the vector with the column-wise format. This operation involves the transfer of data among different VPs along the communication network.

The execution of Step 2 of the algorithm is similar to that of Step 1. A *segmented-add-spread* operation on the simple vector with the compressed column-wise format computes the partial sum of the entries of each column and communicates the total to all rows in the corresponding column. Hence, all addresses of the segment corresponding to the $j$th row have access to the sum $\sum_{i=1}^{m} x_{ij}$. Each VP can now divide the row total $v_j$, which is available locally, by the computed partial sum to estimate the scaling factor $\sigma_j$. The scaling factor is then used to multiply the value of $x_{ij}$, thus completing Step 2.

### 14.3.2 Control parallel computing with RAS

The RAS algorithm (Algorithm 9.3.2) iterates by scaling the rows of the matrix (Step 1) before it proceeds with the scaling of the columns (Step 2). The operations for multiple row updates are independent of each other. Each row can be updated given the current values of the row entries and the target row sum enabling row scaling operations to be performed simultaneously and in parallel for multiple rows. Once all row scaling operations are completed the algorithm may proceed with the simultaneous scaling of the columns. Each row or column scaling operation is a task.

The decision on how to group tasks for parallel computations leaves some room for experimentation. One possibility is to consider each row scaling operation as a task that can be scheduled for execution on any available processor. If there are more rows than processors then each processor will process the next available task. The advantage of this approach is that it achieves good load balancing; that is, all processors will terminate at approximately the same time since the slower processors process fewer rows as it takes longer to process dense rows. The disadvantage of this approach is that the processors must inquire about the availability of tasks and initiate a task. On some computers the overhead from task initiation could be substantial.

An alternative approach to task scheduling is to group rows together into as many blocks as there are processors and assign one block to each processor. Since only one task initialization is required the overhead of the implementation is reduced. To achieve load balancing the size of each block (number of rows) must be chosen so that all blocks have approximately the same number of nonzero entries.

A third approach is to work with a fully simultaneous version of RAS, which iterates concurrently on both rows and columns of the matrix. This simultaneous algorithm, in the sense of (1.9)–(1.10), is mathematically a special case of the block-iterative MART algorithm (Algorithm 6.7.1) and does not require a separate convergence analysis. It has the following form (refer to Section 9.2.3 for the notation).

#### Algorithm 14.3.1 Simultaneous RAS Algorithm

**Step 0:** (Initialization.) Set $\nu = 0$ and $x_{ij}^0 = a_{ij}$ for all $i = 1, 2, \ldots, m$, and $j = 1, 2, \ldots, n$. Choose fixed underrelaxation parameters $\lambda_\rho, \lambda_\sigma$ such that $\epsilon \leq \lambda_\rho,\ \lambda_\sigma \leq 1$ for some fixed arbitrarily small $\epsilon > 0$.

**Step 1:** (Computing row scaling factors.) For $i = 1, 2, \ldots, m$, calculate

$$\rho_i^\nu = \left( \frac{u_i}{\sum_{j \in \delta_i^+} x_{ij}^\nu} \right)^{\lambda_\rho} \tag{14.1}$$

**Step 2:** (Computing column scaling factors.) For $j = 1, 2, \ldots, n$, calculate

$$\sigma_j^\nu = \left( \frac{v_j}{\sum_{i \in \delta_j^-} x_{ij}^\nu} \right)^{\lambda_\sigma} \tag{14.2}$$

**Step 3:** (Scaling the matrix.) For $i = 1, 2, \ldots, m$, and $j = 1, 2, \ldots, n$, compute

$$x_{ij}^{\nu+1} = \rho_i^\nu x_{ij}^\nu \sigma_j^\nu. \tag{14.3}$$

**Step 4:** Replace $\nu \leftarrow \nu + 1$ and return to Step 1.

Simultaneous RAS is a parallel algorithm and can be implemented in parallel in a way similar to RAS: multiple processors will compute the scaling factors for multiple rows and multiple columns (Steps 1 and 2); and the matrix entries (Step 3) are updated in parallel for all entries. The simultaneous RAS can utilize a greater number of processors than RAS, since both row and column sums can be computed simultaneously.

## 14.4 Parallel Computing for Image Reconstruction

We discuss now the parallel implementation of row-action iterative algorithms for image reconstruction, considering in particular the block-MART algorithm (Algorithm 6.7.1) applied to the entropy optimization model for image reconstruction (see Section 10.3).

Image reconstruction problems are natural applications for data parallel programming. The pixels of a discretized image can be mapped onto a parallel machine using checkerboard partitioning (see Section 14.2.1). However, the iterative reconstruction algorithms operate on sets of pixels according to an adjacency relation dictated by the rays intersecting the image. For example, if the $i$th ray runs parallel to one of the edges of the discretized image domain in Figure 14.6, then each pixel needs to communicate with its nearest neighbor on the left or the right and the checkerboard partitioning can be mapped naturally onto a machine with a mesh communication network. However, many rays run through the image at an angle $\theta$ (see Figure 14.7), and the communication requirements are determined by the set of pixels intersected by such rays; thus there is no natural mapping of these pixels onto a mesh, for all angles $\theta$. Data parallel implementations of image reconstruction algorithms usually require special purpose parallel machines with appropriate communication networks. The use of such special purpose machines is, quite often, justified by the applications.

In this section we discuss general purpose approaches to the parallel implementation of reconstruction algorithms based on control parallelism. Parallelism is exploited by partitioning the operations of the algorithm into tasks that can be performed concurrently. We discuss three parallel variants of the block-MART algorithm (Algorithm 6.7.1), all of which are logically special cases of this algorithm such that separate mathematical analysis is not needed. However, the three schemes indicate the flexibility of block-iterative algorithms for parallel computations.

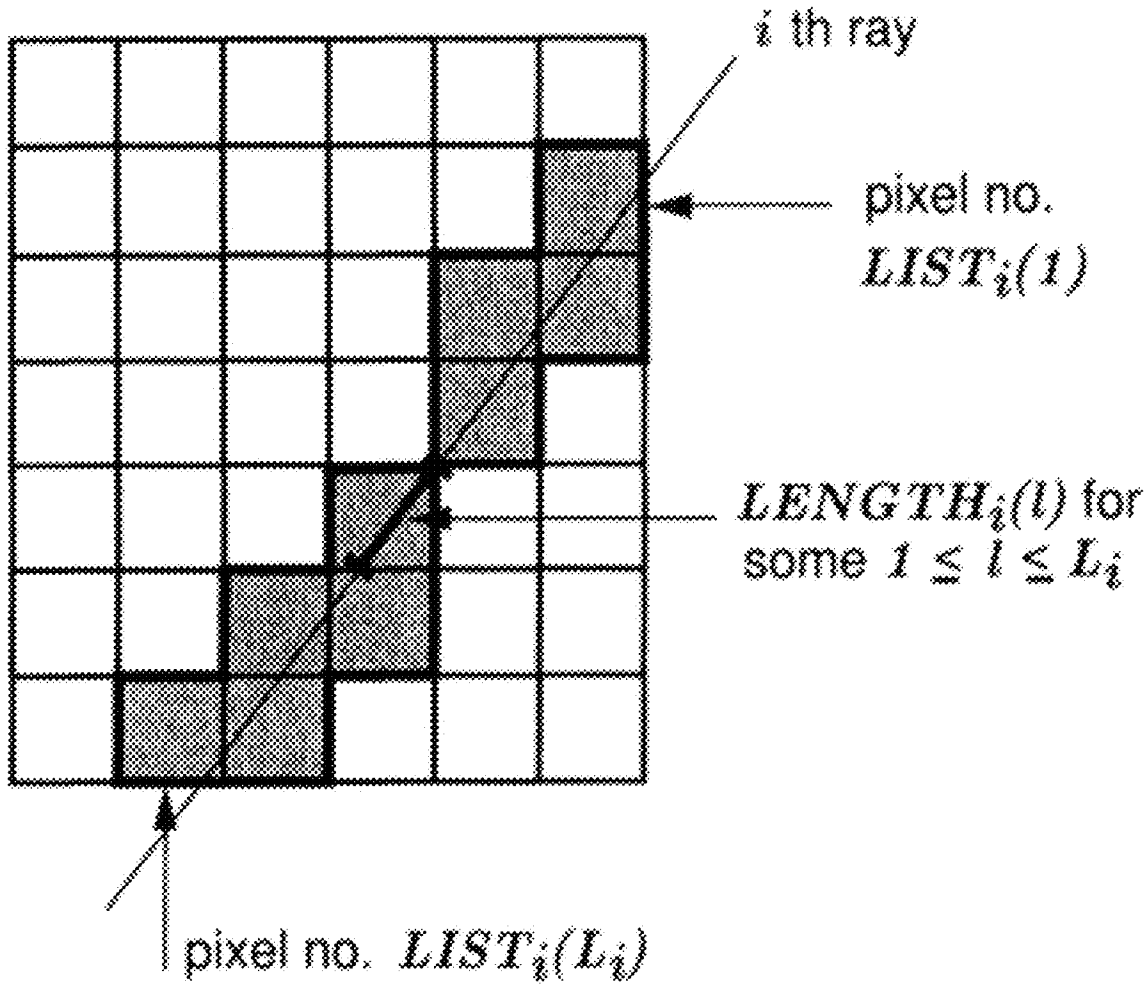

**Figure 14.6** Renumbering of pixels of a discretized image along the path of a ray.

We consider the implementation of block-MART as the reconstruction algorithm within the SNARK77 software package. SNARK77 (as well as its latest available extension SNARK93) is a programming environment designed to facilitate the implementation and testing of reconstruction algorithms. Based on user specifications of a test image and the geometry of data collection, SNARK77 creates a discretized image and simulates the data measurements. This information is then passed on to an algorithm that reconstructs the original image specified by the user. The $\nu$th iterative step of the block-MART algorithm can be written in the form

$$x_j^{\nu+1} = x_j^\nu \exp\left( \sum_{i \in I_{t(\nu)}} w_i^{t(\nu)} d_i^\nu a_j^i \right), \tag{14.4}$$

see (6.125). Instead of computing the term under the summation sign for each pixel $j$ and then updating the corresponding $x_j^\nu$, we proceed as follows. For every pixel $j$ we accumulate the terms

$$w_i^{t(\nu)} d_i^\nu a_j^i \tag{14.5}$$

in an array we call a *correction array.* Once all the terms are accumulated we update the $j$th elements of $x^\nu$. This computation is executed in two intermediate steps—i.e., accumulating the correction terms and updating the vector $x^\nu$—that are performed using vector operations.

We formalize these two steps, omitting temporarily the iteration index $\nu$ for the sake of simplicity. Given an iteration index $\nu$, a block index $t \doteq t(\nu)$ is chosen according to a cyclic control. Let the total number of pixels intersected by the $i$th ray be $L_i$, and let $LIST_i(l)$ be an array of the indices of these pixels, for $l = 1, 2, \ldots, L_i$. Together with $LIST_i(l)$ we know for each $l$ (for example, by using SNARK77), the length of intersection of the $l$th pixel $LIST_i(l)$ with the $i$th ray. This information is stored in an array $LENGTH_i(l)$ for $l = 1, 2, \ldots, L_i$ (see Figure 14.6). We denote by $J_t(\cdot)$ the array containing the indices of all pixels which are intersected by *any* ray that belongs to the block $I_t$. The total number of elements in $J_t(\cdot)$ is $S_t$ which is the number of elements in the set $\{j \mid a_j^i > 0 \text{ for all } i \in I_t\}$. $C_t(s)$ for $s = 1, 2, \ldots, S_t$ is the temporary correction array associated with the $t$th block.

For a given block, the correction array $C_t(\cdot)$ and the image vector updating are computed as follows:

**Initialization:**

$$C_t(s) \leftarrow 0\,, \qquad s = 1, 2, \ldots, S_t. \tag{14.6}$$

**Correction array calculation:** Compute, for all $i \in I_t$,

$$C_t\left(LIST_i(l)\right) \leftarrow C_t\left(LIST_i(l)\right) + w_i^t d_i LENGTH_i(l). \tag{14.7}$$

**Updating the image vector:** For every $r = J_t(s)$ and for all $s = 1, 2, \ldots, S_t$, calculate:

$$x_r \leftarrow \exp\left(\log x_r + C_t(s)\right). \tag{14.8}$$

The computations in equations (14.6)–(14.8) involve operations on vectors and can be implemented using the Basic Linear Algebra Subroutines (BLAS). Such routines are available for a wide range of computers, including vector and parallel architectures, and are coded for maximum efficiency on each machine.

Having described the implementation of the basic iterative step of the algorithm we turn now to alternative parallel implementations.

### 14.4.1 Parallelism within a block

The first parallel scheme partitions the operations within a single block into independent tasks. Recall that the following four steps are executed during one iteration over the equations of a block:

**Step 1:** Choose the rays that form the block using SNARK77 utilities. Given a user-specified size of the block (i.e., the number of rays per block) SNARK77 will determine the rays that form this particular block.

**Step 2:** For each ray $i$ in the block, determine the simulated measurement $y_i$, the list of pixels $j$ intersected by it—$LIST_i(\cdot)$—and the lengths of the intersections $a_j^i$, stored in $LENGTH_i(j)$. This information is provided by SNARK77 based on the geometry of data collection and the discretization of the image.

**Step 3:** Compute the correction terms $C_t(\cdot)$; cf. equation (14.7).

**Step 4:** Update all pixel intensities $x_r$; cf. equation (14.8).

These four steps are executed for each block and are repeated until some iteration number limit is reached or some convergence criterion is satisfied. Blocks are processed in a cyclical manner.

It is possible to utilize multiple processors during the execution of some of the steps. The execution of Step 1 requires synchronization in order to avoid situations where two or more processors assign the same ray to a block. However, this step is not particularly time-consuming, since SNARK77 just goes through the list of rays and picks a user-specified number of them according to a simple rule. For example, it can choose the first ray that has not yet been used. (SNARK77 has various options for choosing rays).

Step 2 is computationally intensive. SNARK77 uses the geometry of the phantom and the orientation of the rays to determine the list of pixels intersected by each ray and the lengths of intersections. The exact specifications of the phantom—known to SNARK77 but, of course, unknown to the reconstruction algorithm—are used to compute the simulated measurements. These operations can be executed concurrently for all rays in the block. Computations are executed in parallel for as many rays as there are processors on the target machine.

In Step 3 the correction term is computed for each ray. Again multiple rays can be processed in parallel. The updating of pixel intensities $x_r$ using equation (14.8) requires synchronization to avoid updating a single pixel by multiple processors, when such a pixel is intersected by several rays that were operated upon in parallel.

### 14.4.2 Parallelism with independent blocks

The block nature of the algorithm and the geometry of image reconstruction problems provide another mechanism for introducing parallelism with larger tasks. It is possible to instruct SNARK77 to generate *independent blocks* that can be processed in parallel. Two blocks are termed independent if their respective rays do not intersect any common pixels. With independent blocks it is known a priori that no write conflicts will arise between the tasks for a given view, so there is no need to synchronize the updating of the pixels.

Independent blocks are formed by grouping together rays that belong to the same view (Figure 14.7). With this parallel scheme, synchronization is required whenever the algorithm begins to operate on blocks in a different view, since it is possible that some processors may still be operating on blocks in the previous view. In this case write conflicts may arise in updating the intensities of pixels shared between tasks; therefore, processors are synchronized to complete execution on the blocks of a given view before operations on a new view may begin.

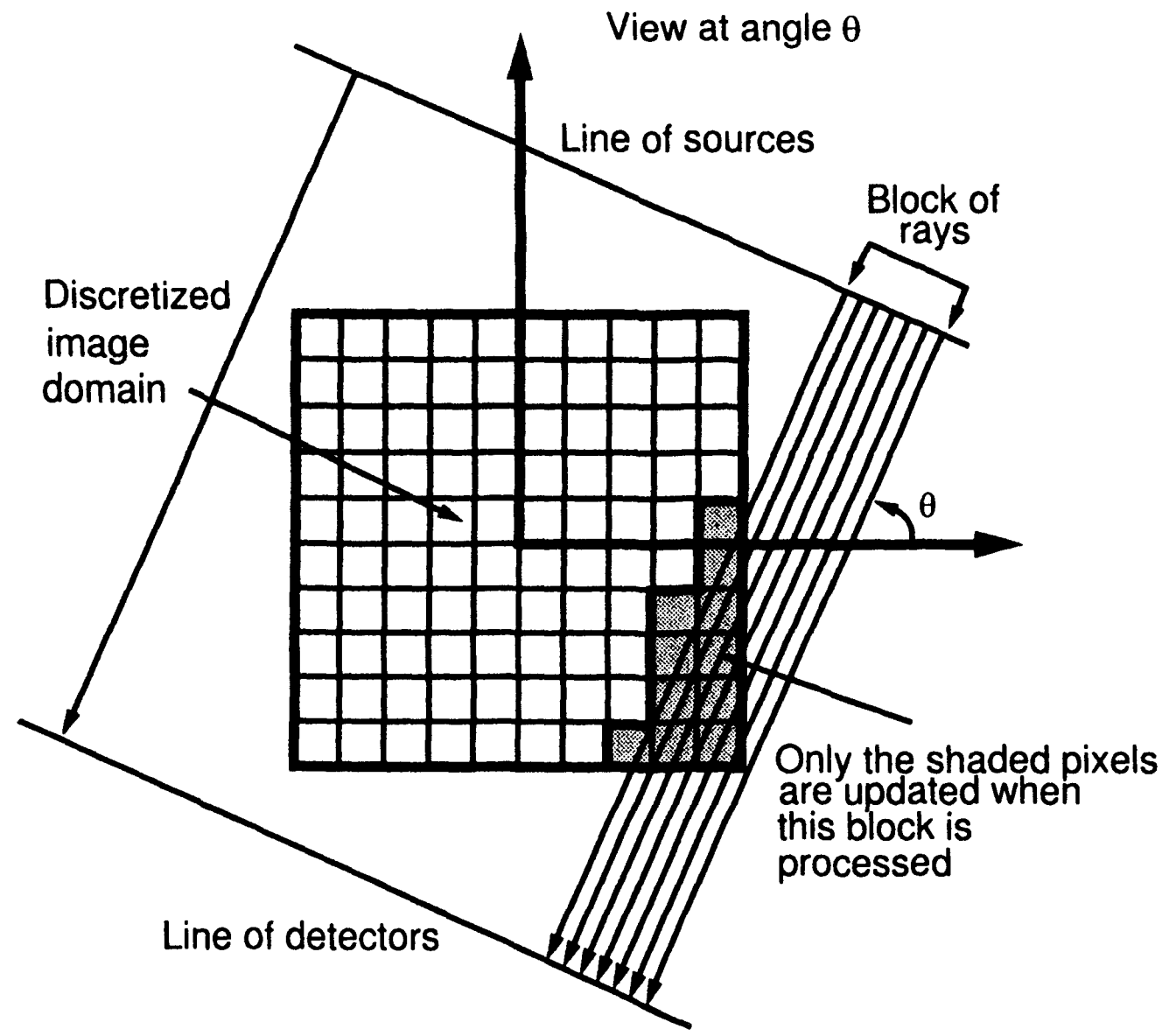

**Figure 14.7** Parallel computing with independent blocks.

### 14.4.3 Parallelism between views

A third parallel scheme is illustrated in Figure 14.8 where all rays in a single view are grouped together as a block and are operated on in parallel. Observe that the block-MART algorithm is a simultaneous block-iterative algorithm as characterized in Section 1.3. We assume $M$ different views in the parallel geometry of data collection, and let $I_t$ contain the ray indices of all rays in the $t$th view. For each $t$, $1 \leq t \leq M$ we define an *intermediate image* $x^{\nu+1,t}$ by using for all $j = 1, 2, \ldots, J$ the formula

$$x_j^{\nu+1,t} = (x_j^\nu)^{1/M} \prod_{i \in I_t} \exp\left(w_i^t d_i^\nu a_j^i\right). \tag{14.9}$$

With this as a particular representation of equation (1.11), take the operator $S$ of equation (1.12) to be given by

$$x_j^{\nu+1} = \prod_{t=1}^{M} x_j^{\nu+1,t}, \quad j = 1, 2, \ldots, J. \tag{14.10}$$

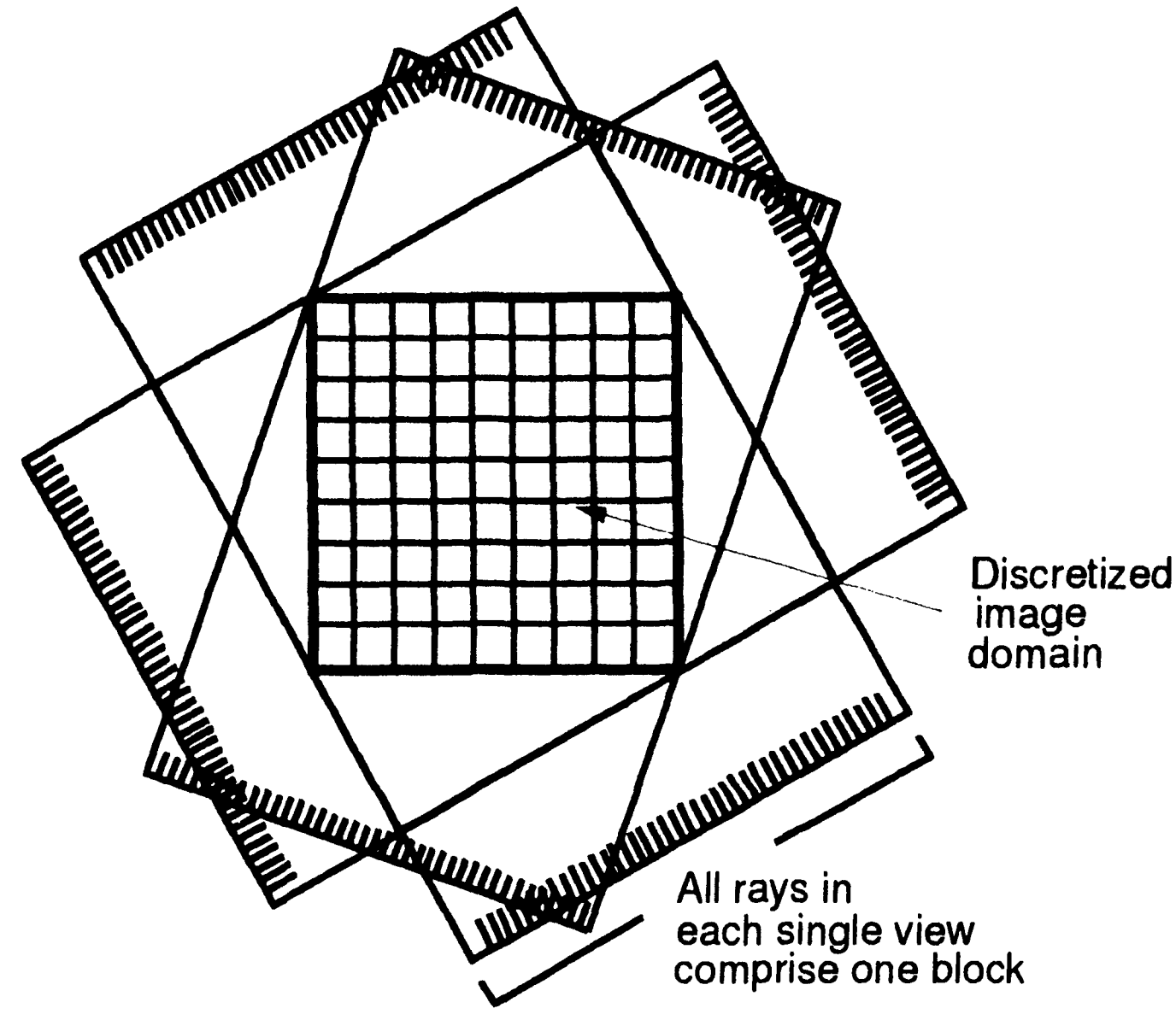

**Figure 14.8** Parallel computing between blocks.

This simultaneous block-iterative version of MART (Algorithm 6.6.2) is mathematically equivalent to a fully simultaneous algorithm with

$$x_j^{\nu+1} = x_j^\nu \prod_{i=1}^{I} \exp\left(w_i d_i^\nu a_j^i\right), \tag{14.11}$$

for all $j = 1, 2, \ldots, J$, where all rays are lumped into a single block. However, from a computational point of view it offers another, completely different, parallel scheme for implementing block-MART for image reconstruction. Since all pixel intensities must be updated by the rays of each view (and hence in parallel by multiple processors) temporary storage must be used for the correction factors computed by each view. The algorithmic operator then aggregates the intermediate correction factors and updates the image vector $x$.

## 14.5 Parallel Computing for Network-structured Problems

We discuss now the parallel implementation of the iterative algorithms of Sections 12.4 and 13.8 for problems with network structures, using a V-RAM model. We give V-RAM models for dense and sparse bipartite graphs, for sparse transshipment networks, and for problems with embedded network structures such as the multicommodity network flow problem and the stochastic programming problem with network recourse.

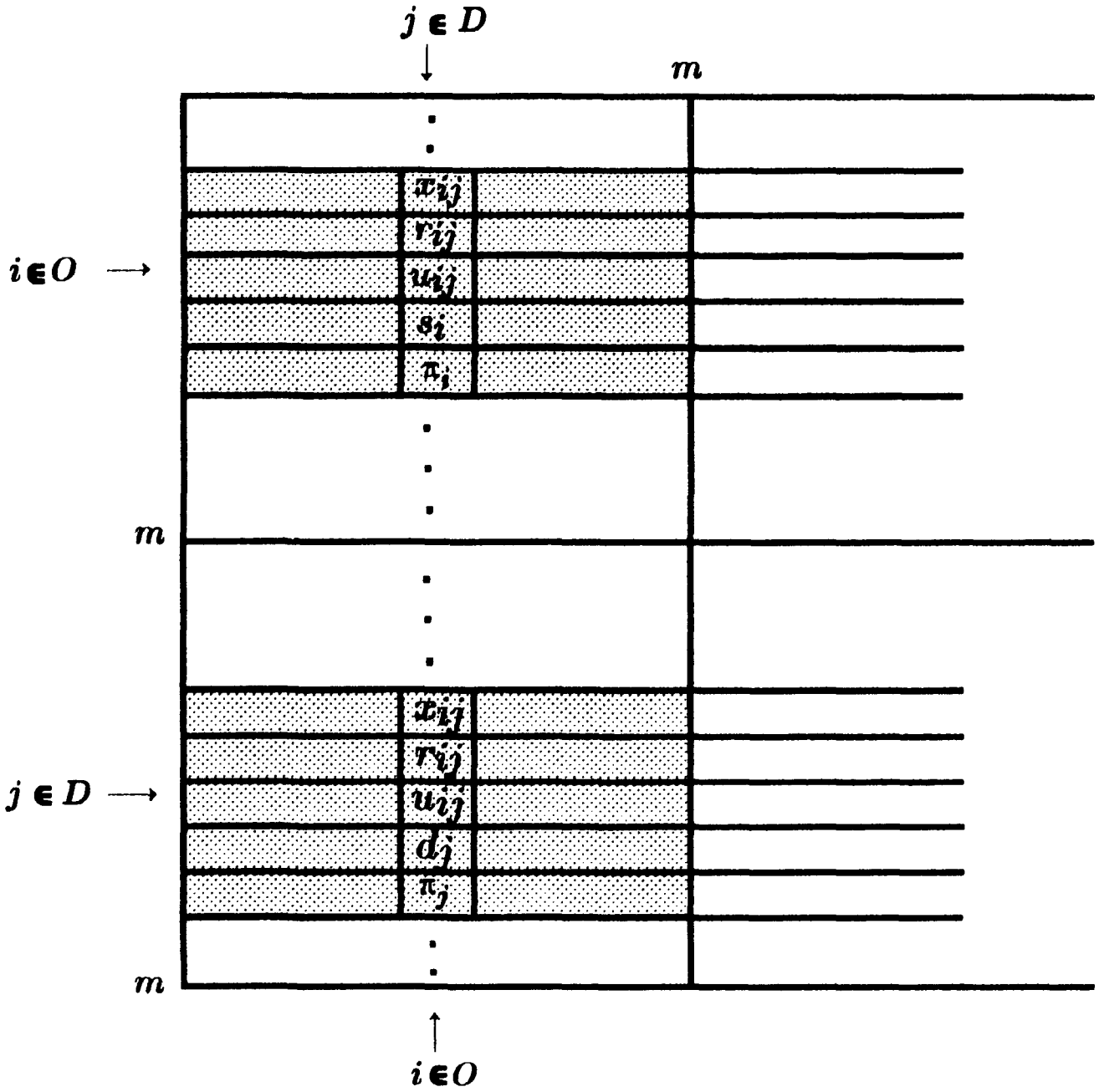

**Figure 14.9** The V-RAM model for a dense $m \times m$ bipartite graph.

### 14.5.1 Solving dense transportation problems

We develop first a V-RAM model for the dense bipartite graph of the transportation problem (Section 12.2.1). The model is similar to that for dense matrices (see Figure 14.9). Assuming that both the number of origin nodes and the number of destination nodes are equal to $m$, the model uses $2m\kappa$ simple vectors. The length of each vector is $m$, which is the number of arcs incident to each node. The constant $\kappa$ depends on the data used by the specific algorithm that is implemented. To implement Algorithm 12.4.1 for the quadratic transportation problem the following simple vectors are needed for each origin node $i$: vector of flows $x_{ij}$, vectors of upper bounds $u_{ij}$ and dual prices $r_{ij}$ for all $j = 1, 2, \ldots, m$, and vectors for the supplies and their dual prices. All addresses of the supply and dual price vectors for the $i$th node hold the same entries, i.e., $s_i$ and $\pi_i^O$ respectively. For all destination nodes similar simple vectors are needed for the arc-related information $x_{ij}, u_{ij}, r_{ij}$, and additional vectors are needed for the demands and their dual prices $d_j$ and $\pi_j^D$, (see Figure 14.9).

With these data structures, a vector instruction can be carried out simultaneously on the elements of all simple vectors corresponding to origin nodes, thus executing Step 1 of Algorithm 12.4.1. For the implementation of Step 2, vector instructions are executed on all simple vectors of desti-nation nodes. For example, an *add-spread* operation on the simple vectors that store $x_{ij}$ for each origin node computes the partial sum $\sum_{j=1}^{n} x_{ij}$ for all $i$, and communicates the total to all addresses of the $i$th node. This sum and the local copy of the supply $s_i$ are used to compute the scaling factor $\rho_i$, which is then used to update the flow $x_{ij}$. Similar operations are executed on the simple vectors corresponding to the destination nodes to complete Step 2.

In order to avoid redundant use of VPs by the implementation we note that only $m\kappa$ of the simple vectors are operated upon simultaneously: first the $m\kappa$ vectors of the origin nodes for execution of Step 1, and then the $m\kappa$ vectors of the destination nodes for execution of Step 2. Each vector is of length $m$, so that at most $m^2\kappa$ VPs are required to implement this algorithm. Because they are symmetrical, the vectors corresponding to the arc-related information $x_{ij}, u_{ij}$, and $r_{ij}$ are identical for origin and destination nodes.

### 14.5.2 Solving sparse transportation problems

The V-RAM model for a sparse bipartite graph of the transportation problem of Section 12.2.1 is similar to the V-RAM model of sparse matrices shown in Figure 14.10. It uses $(m_O + m_D)\kappa$ simple vectors for the $m_O$ origin and $m_D$ destination nodes. The length of each vector is equal to the degree of the corresponding node, and the constant $\kappa$ depends on the specific algorithm that is implemented. To implement Algorithm 12.4.1 or 12.4.2 for nonlinear transportation problems the following simple vectors are needed for each origin node $i$: vector of flows $x_{ij}$, vector of upper bounds $u_{ij}$ and vector of dual prices $r_{ij}$ for all $j = 1, 2, \ldots, m$, and vectors for the supply and dual price. Again, all addresses of the supply and dual price vectors for the $i$th node hold the same entries, i.e., $s_i$ and $\pi_i^O$ respectively. Likewise, for each destination node simple vectors are needed for the arc-related information $x_{ij}, u_{ij}, r_{ij}$, and additional vectors are needed for the demand and dual prices $d_j$ and $\pi_j^D$.

This V-RAM model can be implemented on a linearly ordered array of $n$ VPs as follows: The array is partitioned into $m_O$ segments, with the length of the $i$th segment equal to the out degree of the $i$th node, i.e., the number of arcs outgoing from node $i$. The set of these arc indices is denoted by $\delta_i^+$ with cardinality $q_i^+$. The VPs of the $i$th segment store the entries $x_{ij}, u_{ij}, r_{ij}$ for all $j \in \delta_i^+$. All $q_i^+$ addresses of the $i$th segment also store the entries $s_i$ and $\pi_i^O$. The segments are indicated by a vector of segment bits that has a 1 at locations $1, q_1^+, q_1^+ + q_2^+, \ldots, n$, and 0 elsewhere (see Figure 14.10). The sums of the flows on all arcs outgoing from a node are computed by a *segmented-add-scan* operation on the parallel variables $x_{ij}$ for $i = 1, 2, \ldots, m_O$ and all $j \in \delta_i^+$.

The same linearly ordered array can be partitioned into $m_D$ segments, with the $j$th segment of equal length to the in degree of the $j$th node, i.e., the number of arcs going into node $j$. The set of these arcs is denoted by $\delta_j^-$ with cardinality $q_j^-$. The VPs of the $j$th segment store the entries $x_{ij}, u_{ij}, r_{ij}$ for all $i \in \delta_j^-$. All $q_j^-$ addresses of the $j$th segment also store the entries $d_j$ and $\pi_j^D$. The segments are indicated by a vector of segment bits that has a 1 at locations $1, q_1^-, q_1^- + q_2^-, \ldots, n$, and 0 elsewhere (see again Figure 14.10). The sums of the flows on all arcs going into a node are computed by a *segmented-add-scan* operation on the parallel variables $x_{ij}$ for $j = 1, 2, \ldots, m_D$, and all $i \in \delta_j^-$.

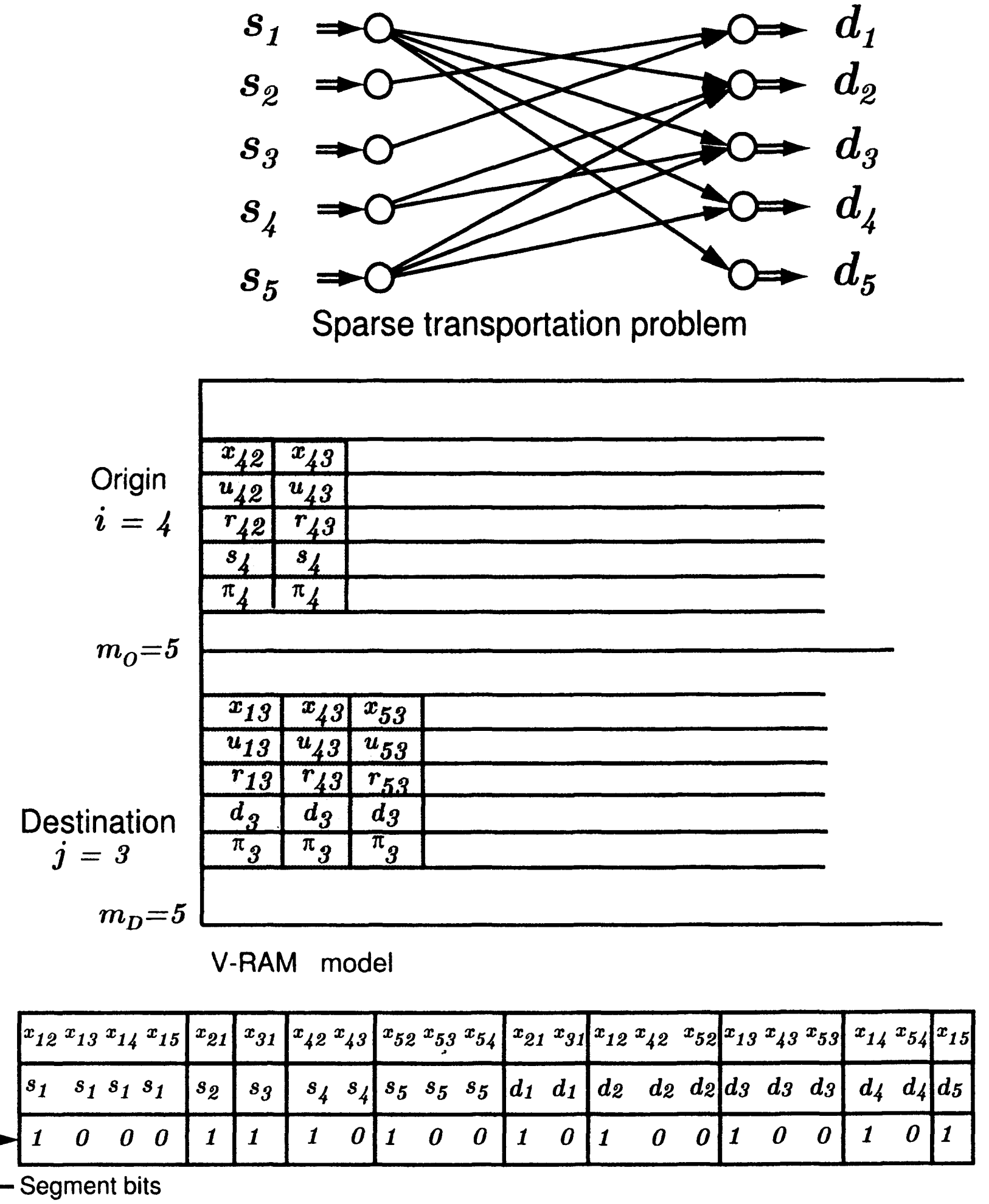

| $x_{12}$ $x_{13}$ $x_{14}$ $x_{15}$ | $x_{21}$ | $x_{31}$ | $x_{42}$ $x_{43}$ | $x_{52}$ $x_{53}$ $x_{54}$ | $x_{21}$ $x_{31}$ | $x_{12}$ $x_{42}$ $x_{52}$ | $x_{13}$ $x_{43}$ $x_{53}$ | $x_{14}$ $x_{54}$ | $x_{15}$ |
|---|---|---|---|---|---|---|---|---|---|
| $s_1$ $s_1$ $s_1$ $s_1$ | $s_2$ | $s_3$ | $s_4$ $s_4$ | $s_5$ $s_5$ $s_5$ | $d_1$ $d_1$ | $d_2$ $d_2$ $d_2$ | $d_3$ $d_3$ $d_3$ | $d_4$ $d_4$ | $d_5$ |
| 1 0 0 0 | 1 | 1 | 1 0 | 1 0 0 | 1 0 | 1 0 0 | 1 0 0 | 1 0 | 1 |

**Figure 14.10** The V-RAM model for a sparse bipartite graph.

Similar to the dense bipartite graph implementation we note that a vector instruction can be executed first on all entries of the simple vectors corresponding to origin nodes, and then on the simple vectors of destination nodes. Only $\kappa$ simple vectors are operated upon simultaneously, and this model can be implemented using at most $n\kappa$ VPs, where $n$ is the number of arcs. For a dense $m \times m$ graph we have $n = m^2$ and we need the same number of virtual processors as in the dense implementation. In the sparse implementation, however, we do not have the symmetry of the simple vectors for the arc-related quantities, thus we cannot use the same VP to store $x_{ij}, u_{ij}, r_{ij}$ for both origin and destination nodes. These parameters are stored in two copies: once in the segment of origin node $i$, and once in the segment of destination node $j$. Before the algorithm switches from iterations on origin nodes to iterations on destination nodes the values of the arc-related variables $x_{ij}$ must be updated, i.e., the two copies of the simple vectors must be brought into agreement. This is a memory-to-memory communication operation.

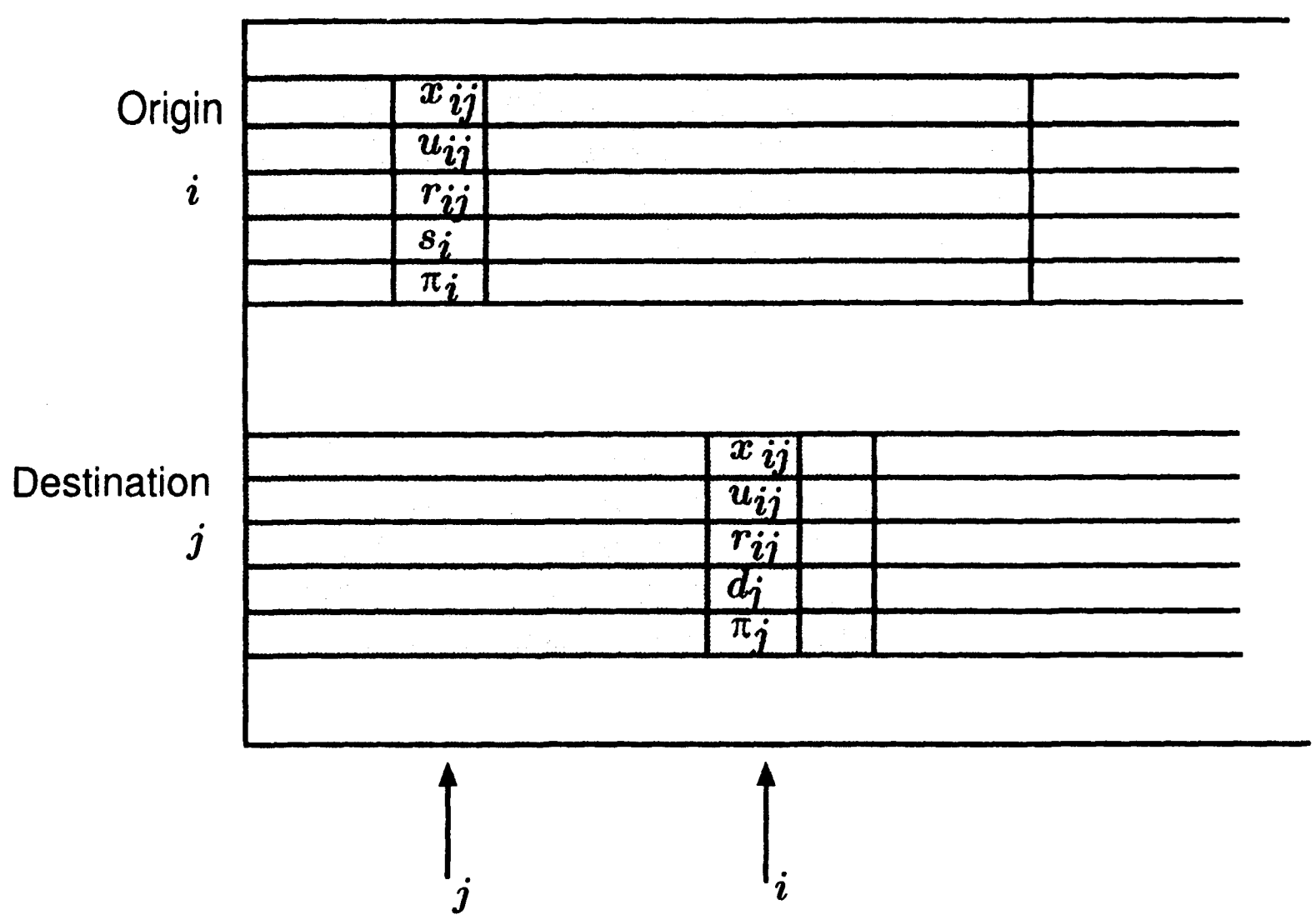

**Figure 14.11** The V-RAM model for a sparse transshipment graph.

### 14.5.3 Solving sparse transshipment graphs

The V-RAM model for sparse transshipment graphs is shown in Figure 14.11. Like the sparse bipartite graph model it uses $m\kappa$ simple vectors for all nodes. The length of each vector is equal to the sum of the in and out degrees of the node. The $i$th simple vector stores information for all arcs going into or coming out from node $i$, i.e., all arcs $(i, j)$ with $j \in \delta_i^+$ and $j \in \delta_i^-$. Hence, arc-related information is stored twice, once in the simple vector corresponding to the origin node of an arc, and once in the simple vector of the destination node. Unlike the case of the bipartite graph, it is not possible to implement this V-RAM model using only $n$ linearly ordered VPs. Instead, $2n$ VPs are needed since each arc is stored twice, and we cannot exploit the partitioning of nodes between origin and destination as we did for the bipartite graphs.

### 14.5.4 Solving network-structured problems

The V-RAM model for sparse network-structured problems, such as the multicommodity network flow problem and the stochastic program with network recourse, are easy generalizations of the model in Figure 14.11. Consider, for example, the split-variable formulation of the stochastic network problem discussed in Section 13.3.3. Multiple copies of the simple vectors of Figure 14.11 are needed, one for each scenario. With this model it is possible to iterate concurrently on the simple vectors for all scenarios.

However, these simple vectors are not totally independent. In the case of stochastic network problems, the first-stage split variables must converge to the same value (see Section 13.7). For the multicommodity network flow problem the total flow of all commodities on each arc should lie between prespecified upper and lower bound; (see Section 12.2.2). Hence, the V-RAM model for network-structured problems requires additional simple vectors to indicate the first-stage variables that were split, or the commodities flowing on the same arc.

## 14.6 Parallel Computing with Interior Point Algorithms

In this section we consider the parallel implementation of the interior point algorithm for stochastic programming problems. In particular, the computationally intensive operations involved in the dual step direction calculation can be implemented in parallel. We consider the implementation of Procedure 8.3.2 on distributed memory machines, with particular emphasis on the efficient execution of the communication steps. This procedure is well suited for parallel implementation because the various operations involving separate submatrices can be carried out independently.

The computation begins with the submatrices $T_l$, $W_l$, and $\Theta_l$ (see, e.g., Theorem 8.3.1 for the definitions of the submatrices), and the vector segment $\psi^l$ located on processor $l$, for $l = 1, 2, \ldots, N$. Processor $l$ can then compute $S_l$ and proceed independently with all computations involving these data alone. Interprocessor data communication is necessary at only three steps of this procedure. In particular, in Step 2b, the processors must communicate to form the matrix $G_1$ and the vectors $\hat{p}^1$ and $\hat{p}^2$. If Steps 2c, 2d, and 2e proceed serially on the master processor, then the processors must communicate to broadcast the computed vector $q^1$. Steps 3 and 4 require only the distributed data $S_l$, $T_l$, $q^l$, and $p^l$ on processor $l$ and so may be carried out with full parallelism. A final communication step accumulates all the partial vectors $(\Delta y)^l$ at the master processor for use in subsequent iterations. Thus, all communication steps require either broadcasting of data from one processor to all others (*one-to-all* communication in Step 2e) or gathering of data distributed among processors in Steps 2b and 4. The data gathers can be either *all-to-all* or *all-to-one*, depending on the intended use of the data and the speed of the gather routines.

We discuss the parallel implementation of this procedure on distributed-memory MIMD multiprocessors. In the next section we introduce the hypercube multiprocessor used as the target architecture. We also describe the optimal one-to-all and all-to-all communication routines that are the basis for those used in the parallel implementation.

### 14.6.1 The communication schemes on a hypercube

We consider a distributed memory MIMD message-passing parallel computer in which processors are connected according to a hypercube graph. The size of a hypercube is defined by its dimension $d$; and a hypercube graph of dimension $d$ and the multiprocessor based on it are called *d-cubes*. A $d$-cube graph has $p \doteq 2^d$ nodes, and the processors of the hypercube computer are located at the nodes of the graph.

Nodes in a $d$-cube are assigned $d$-bit binary identifiers (from 0 through $p-1$) such that all $d$ nodes connected to node $j$ have identifiers differing from $j$ in exactly one bit. A $d$-cube can be constructed by connecting corresponding nodes of two $(d-1)$-cubes in any of $d$ different ways. For example, the common three-dimensional 3-cube can be constructed by linking corresponding processors of two squares (or 2-cubes) to form its top and bottom, left and right, or front and back faces. In a $d$-cube, the $d$ neighbors of node $j$ define the $d$ nodes corresponding to node $j$ in the $d$ different $(d-1)$-cubes. Figure 14.12 illustrates a 3-cube and the mapping of a spanning binary tree onto the hypercube graph.

The most efficient codes on hypercubes are those in which messages do not cross communication wires. One way to ensure contention-free routing of messages is to ensure that only processors located on neighboring nodes of the cube communicate, and this rule is used to the extent possible in the parallel implementation of Procedure 8.3.2.

The first communication scheme broadcasts data from a master processor to all other processors via the links of a spanning tree of the hypercube graph rooted at processor 0. All communication takes place between nearest neighbors. This procedure is termed a *spanning tree broadcast* (STB). A spanning tree of height $d$ is embedded into a $d$-cube by simple bit manipulation of the node identifiers. Node 0 is at the root of the tree. The broadcast proceeds for a total of $d$ communication steps where, at step $l$ for $l = 1, 2, \ldots, d$, each processor with node identifier $j < 2^l$ pairs with the processor with a node identifier different from $j$ in the $l$th bit only. For each pair, the node with the smaller identifier is the parent node which sends its data to its child. Thus, at step $l$ of the algorithm, $2^{l-1}$ processors receive the data. We use an STB to broadcast the vector $q^1$ in the parallel implementation of Step 2e of Procedure 8.3.2.

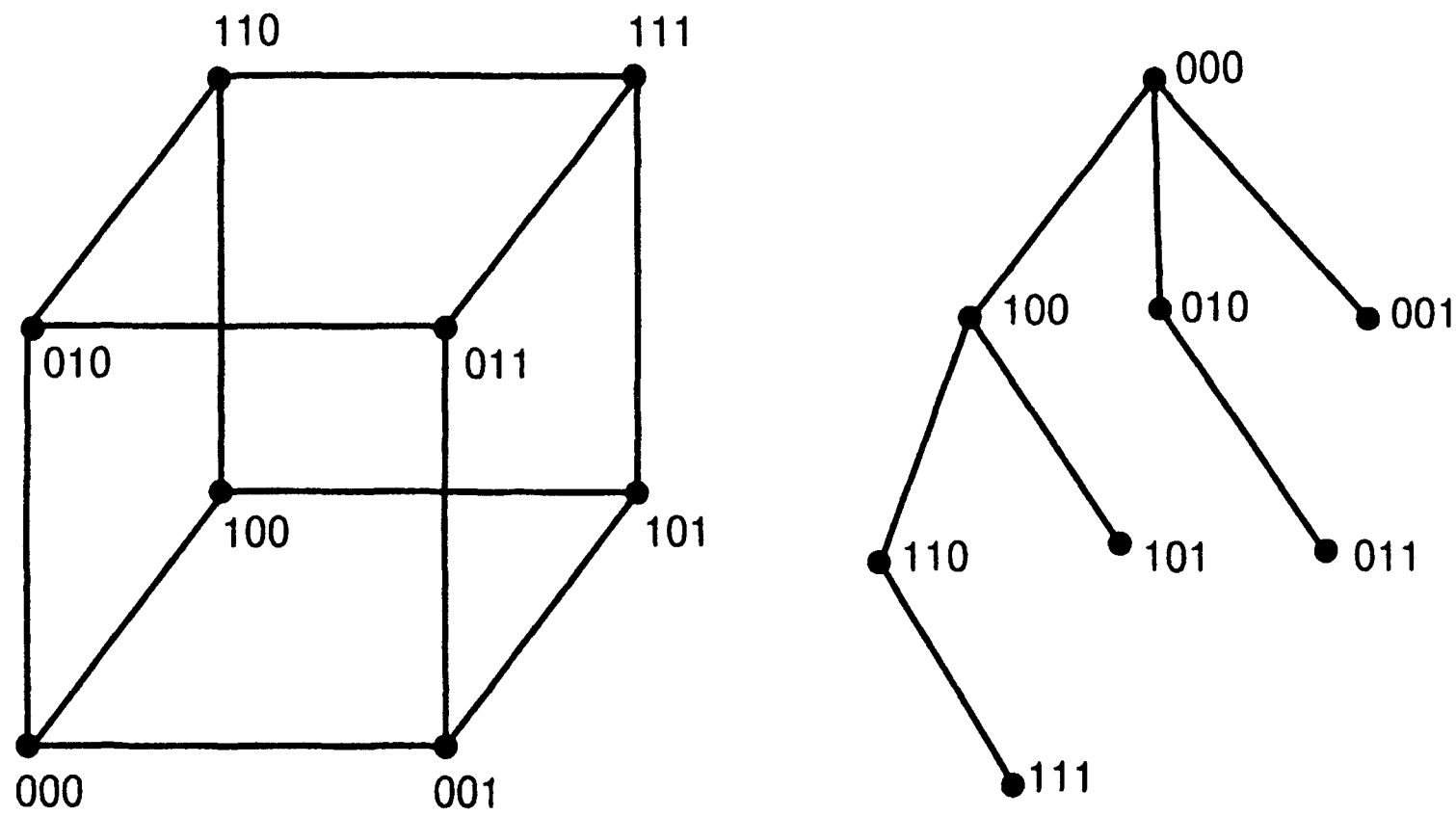

**Figure 14.12** A 3-cube and the mapping of a spanning binary tree onto the hypercube graph.

The STB moves data from one processor to all others. Suppose, in contrast, that each of the $p$ processors begins with $k$ bytes of data but that the full $pk$ bytes of data need to be accumulated at a master processor in node 0. The distributed data can be gathered by traversing the spanning tree from its leaves to its root. Again, the gather operation proceeds for $d$ steps; but every processor receives data from its child node, appends that data to its own data, and forwards the accumulated data to its parent in the next communication step. Upon completion, the processor at node 0 holds all $pk$ bytes of data. This is known as a *spanning tree gather* (STG). We use a global sum routine based on an STG in the parallel implementation of Step 2b of Procedure 8.3.2 to form the matrix $G_1$, using the distributed results computed in Steps 2a and 2b.

The second communication routine is based on an *alternate direction exchange* (ADE). An ADE can be used, for example, to accumulate in all $p = 2^d$ processors $pk$ bytes of data initially distributed with $k$ bytes per processor. In each of $d$ communication steps, the $d$-cube splits into a different pair of $(d-1)$-cubes. The $2^{d-1}$ pairs of corresponding nodes from the two cubes exchange and accumulate their data sets. Thus, as in the STG, the size of the data set received by any processor at communication step $l$ is $2^{l-1}k$ bytes, $l = 1, 2, \ldots, d$. Again, all communication takes place between processors located at neighboring nodes.

STB, STG, and ADE are optimal broadcast and gather routines in the sense that they take only the minimal number of communication steps, namely, $d$ steps on a $d$-cube. Minimizing the number of communication steps is important on distributed memory computers because the cost of communication is typically high compared to the cost of computation. In particular, the cost of communicating an $m$-byte message from one hypercube node to a neighboring one is $\beta+m\tau$, where the communication startup latency $\beta$ is generally large in comparison to the transmission time per byte $\tau$. (On the Intel iPSC/860, see Section 15.7; for example, $\beta/\tau \approx 340$ for messages of more than 100 bytes.)

### 14.6.2 The parallel implementation on a hypercube

The basic structure of our parallel implementation of Procedure 8.3.2 is as follows. We assume that the number of scenarios $N$ equals the number of processors $p$. If $N$ and $p$ differ, the scenarios are distributed equally among the processors. If $p$ does not divide $N$, then the remaining scenarios and tasks are distributed one per processor until all have been assigned.

**Steps 1, 2a, and 2b** All processors solve in parallel their systems for $p^l$ and $(u^l)^i$ for $i = 1, 2, \ldots, n_0$. (Note that the matrices $S_l$ need be factored only once to solve the $n_0 + 1$ systems.) All processors then compute their matrices $\hat{T}_l \doteq (T_l^{\mathrm{T}}(u^l)^1, \ldots, T_l^{\mathrm{T}}(u^l)^{n_0})$ in parallel.

The processors then accumulate the matrix $G_1$ on the master processor by means of a *spanning tree global sum.* This communication scheme is based on an STG, but after a processor receives the partial sums from all its children it adds them to its own partial sum and then sends the result to its parent. After the master processor receives the partial sums from its children, it adds them to $\Theta_0^{-2} + A_0^{\mathrm{T}} A_0$ to form the matrix $G_1$. At the same time, the global sum $\hat{p}^1 = A_0^{\mathrm{T}} p^0 + \sum_{l=1}^{N} \hat{T}_l p^l$ is produced at the master processor. The master processor computes $\hat{p}^2 = -p^0$ serially.

**Steps 2c, 2d, and 2e** In our hypercube implementation, all computations involving the small-order dense matrices $G_1$ and $G_2$ are executed serially on the master processor. This is more efficient than applying parallel dense matrix techniques for problems where the number of first-stage variables is low compared to second-stage variables, and for problems where many scenarios are solved at each processor (i.e., when $N/p \gg 1$). If the first-stage matrices are large we may consider alternative implementations using parallel dense linear algebra techniques to execute steps 2c, 2d, and 2e. The master processor broadcasts the computed vector $q^1$ to all other processors via an STB.

**Steps 3 and 4** Steps 3 and 4 are essentially completely parallel: the master processor determines $r^0$ and all processors compute the matrix-vector products $T_l q^1$ and solve their sparse systems. The processors then compute the vector segments $\Delta y^l$ in parallel. As a final step in the algorithm, the processors accumulate the full correction vector $\Delta y$ at the master processor. This final step can be excluded if the next iteration calls for a distributed vector.

### 14.6.3 An alternative parallel implementation

The modular structure of the parallel Procedure 8.3.2 allows it to be implemented on machines suitable for coarse-grained parallel tasks. In this section, we show how the procedure can be implemented on a different distributed memory MIMD computer—the Connection Machine CM–5. The processors of the CM–5 are connected according to a special network (known as a *fat tree*), and near-optimal communication algorithms (such as broadcast) are available. Using the communication routines of the new target machine makes the implementation straightforward. However, a simple modification to the implementation (discussed in the previous section) results in improved efficiency.

The modification concerns the global summation required to form $G_1$ and $\hat{p}^1$ in Steps 2a and 2b. In the implementation on the hypercube, $G_1$ and $\hat{p}^1$ are formed on the processor at node 0, and Steps 2c, 2d, and 2e are carried out serially at node 0. The vector $q^1$ produced at Step 2e is then broadcast to all processors before Step 3 starts. On the Connection Machine CM–5, a global vector reduction function, *add-spread*, takes vectors (or matrices) distributed in every processor as input, sums them, and leaves a single vector (or matrix) sum at every processor. Using this routine we can implement Procedure 8.3.2 such that Steps 2c, 2d, and 2e are performed in parallel, but redundantly, on all processors. Hence, every processor has a copy of vector $q^1$ (at the end of Step 2e), and the broadcast call (STB) before Step 3 is no longer necessary.

The modified parallel procedure on the CM–5, with the redundant calculations, is summarized below.

**Procedure 14.6.1 Parallel matrix factorization for dual step direction calculation on the Connection Machine CM–5**

*Begin with the following data distribution: Processor $l$ holds $S_l$, $T_l$, and $\psi^l$, $l = 1, 2, \ldots, N$, and a copy of $A_0$, $S_0$, $\Theta_0$, and $\psi^0$.*

**Step 1:** *(Solve $Sp = \psi$ in parallel.) Solve on all processors redundantly, the system $S_0 p^0 = \psi^0$. Solve in parallel on processors $l = 1, 2, \ldots,$ $N$, the system $S_l p^l = \psi^l$ for $p^l$.*

**Step 2:** *(Solve $Gq = V^{\mathrm{T}} p$ in parallel.)*

- **a.** *Solve in parallel on processors $l = 1, 2, \ldots, N$, the system $S_l (u^l)^i = (T_l)_{\cdot i}$ for $(u^l)^i$, $i = 1, 2, \ldots, n_0$.*
- **b.** *Multiply $T_l^{\mathrm{T}} (u^l)^i$ for $i = 1, 2, \ldots, n_0$, in parallel, on processors $l = 1, 2, \ldots, N$.*
  *Use a global vector reduction add-spread, to form $G_1$ and $\hat{p}^1$ on all processors $l = 1, 2, \ldots, N$.*

*Compute $\hat{p}^2$ redundantly on all processors.*

**c.** *Solve on all processors redundantly, the system $G_1 u = \hat{p}^1$ for $u$ and set $v = \hat{p}^2 + A_0 u$.*

**d.** *Form $G_2$ on all processors redundantly, by solving $(G_1)w^i = (A_0^{\mathrm{T}})_{\cdot i}$ for $w^i$ for $i = 1, 2, \ldots, m_0$, and setting $G_2 = -A_0[w^1\, w^2 \cdots w^{m_0}]$.*

**e.** *Solve on all processors redundantly, the system $G_2 q^2 = -v$ for $q^2$, and the system $G_1 q^1 = \hat{p}^1 - A_0^{\mathrm{T}} q^2$ for $q^1$.*

**Step 3:** *(Solve $Sr = Uq$ in parallel.) Set $r^0 = A_0 q^1 + q^2$ on all processors redundantly. Solve on processors $l = 1, 2, \ldots, N$, the system $S_l r^l = T_l q^1$ for $r^l$.*

**Step 4:** *(Form $\Delta y$ in parallel.) Set $\Delta y^0 = p^0 - r^0$ on all processors redundantly. Set $\Delta y^l = p^l - r^l$ on processors $l = 1, 2, \ldots, N$.*

## 14.7 Notes and References

The topic of implementation of parallel algorithms, i.e., *parallel programming*, is currently under intensive research. It is no more complex than the programming of scalar uniprocessor, machines for which many books have been written, but it is still less developed. The book by Carriero and Gelernter (1992) lays the foundation for programming techniques for a large class of parallel machines, namely general purpose asynchronous MIMD machines. A survey of programming languages for parallel machines is given in Bal, Stenier, and Tanenbaum (1989). Kumar et al. (1994) discuss the design and analysis of parallel algorithms with particular attention to implementation issues on a variety of current architectures.

**14.1** Blelloch (1990) introduced the Vector-Random Access Machine (V-RAM) as an abstract model for data parallel computing. Data parallelism was introduced in the dissertation of Hillis (1985); see also Hillis (1987). Hillis and Steele (1986) discuss data parallel algorithms. For a textbook treatment of the concepts introduced in this section see Kumar et al. (1994).

**14.2** The mapping of dense and sparse matrices on parallel processors is discussed by several authors. See, e.g., Kumar et al. (1994) or Carriero and Gelernter (1992). The V-RAM models were developed by Blelloch (1990).

**14.3.1–14.3.2** Data parallel computing for matrix balancing using the RAS algorithm is discussed in Zenios (1990), and control parallel computing with the same algorithm is discussed in Zenios and Iu (1990). The material of this section can be extended to the implementation of algorithms for interval-constrained matrix balancing, and Censor and Zenios (1991) discuss the data parallel implementation of the RRAS algorithm.

**14.4** The image reconstruction algorithms in this section were discussed in Chapters 6 and 10. The potential advantages of alternative reconstruction algorithms for parallel computations were discussed by Censor (1988). The material on the parallel implementation of block-MART is based on the work of Zenios and Censor (1991b). The SNARK77 software package is described in the SNARK77 guide; see Herman and Rowland (1978). Its latest extension, SNARK93, is available from the Medical Image Processing Group (MIPG) at the Department of Radiology of the University of Pennsylvania; see Browne, Herman, and Odhner (1993). The Basic Linear Algebra Subroutines (BLAS) library is described in Dongarra et al. (1988a, 1988b). Parallel implementation of image reconstruction algorithms from noisy data is discussed in Herman et al. (1990).

**14.5** A survey of data parallel computing for network-structured problems is given in Zenios (1994a). The data structures for the network problems of this section have their origin in the mapping of network transshipment problems on the Connection Machine CM–1 introduced by Zenios and Lasken (1988a). Similar data structures for the representation of sparse matrices were developed independently by Blelloch (1990). The data structures for dense and sparse transportation problems were developed in Zenios and Censor (1991a). The extension for the multicommodity transportation problem was done by Zenios (1991b), and the extension for stochastic network problems by Nielsen and Zenios (1993a). A comparison of alternative implementations of network data structures on the Connection Machine CM-2 is given in Nielsen and Zenios (1992a). Similar data structures for the solution of assignment and transportation problems using other algorithms (the auction algorithm, the method of multipliers, and the alternating direction method of multipliers) were used by Wein and Zenios (1991). The algorithms of Goldberg (1987) and Eckstein (1989, 1993) were implemented on the Connection Machine CM–2 using the data structures for sparse transshipment problems.

**14.6** The parallel implementation of interior point algorithms for stochastic programming problems on the Intel iPSC/860 and the Connection Machine CM–5 is discussed in Jessup, Yang, and Zenios (1994a). Details of the implementation on an Intel iPSC/860 hypercube, and a model of the performance of the parallel algorithm, are given in Jessup, Yang, and Zenios (1994b). Communication algorithms for hypercubes are discussed by Saad and Schultz (1989). See also Bertsekas and Tsitsiklis (1989) and Kumar et al. (1994) for the properties of hypercubes. The performance of the Intel iPSC/860 hypercube is discussed in Dunigan (1990). The communication algorithms for the Connection Machine CM–5 fat tree are given in Leiserson et al. (1992).

# CHAPTER 15

# Numerical Investigations

> Failures must be considered the cue for further application of effort and concentration of will power. And if substantial efforts have already been made, the failures are all the more joyous. It means that our crowbar has struck the iron box containing the treasure. Overcoming the increased difficulties is all the more valuable because in failure the *growth of the person performing the task* takes place in proportion to the difficulty encountered.
>
> Aleksandr I. Solzhenitsyn, *The First Circle*, 1968

In the first chapter of this book we argued that parallel processing can improve by several orders of magnitude the speed with which calculations are executed on a computer. Such improvements have a significant impact on various areas of application, and it is these anticipated improvements that largely motivate the investigations into methods of parallel optimization. Throughout the book we developed feasibility and optimization algorithms that are suitable for implementation on parallel machines; we described large-scale applications that are solved with these algorithms, and discussed different modes of parallel implementations. To what extent can the goal of solving large-scale applications be achieved with parallel processing? In this chapter we answer this question by reporting computational results with the parallel solution of several of the applications discussed earlier.

The main objective of this chapter is to study numerical experiments that show that the selected parallel algorithms perform well when implemented on parallel machines. To establish this assertion we rely on the performance measures introduced in Section 1.4. The reported experiments also illustrate how the combination of parallel algorithms and parallel computers can indeed solve much larger problems than have been possible with competing serial algorithms on scalar uniprocessor architectures. However, we do not discuss some important issues, such as the cost of a machine or the ease of programming an algorithm in parallel. It is left to the reader to evaluate, for his or her specific project, whether the anticipated performance of an algorithm—as illustrated by the summary results reported here—and the significance of the application at hand justify either the cost of a parallel machine or the efforts in developing the parallel code. The reader should bear in mind that the results presented here are summary, but they were chosen to be typical and illustrative. Additional results can be found in the literature discussed in Section 15.8, which can be used to complement our summarized results.

The technology of parallel computing is currently changing very rapidly, and therefore care is needed in interpreting the results of this chapter. First, absolute measures—such as solution times or computing rates in FLOPS—quickly become obsolete. In the early 1980s applications executing at MFLOPS (i.e., $10^6$ FLOPS) rates were considered very efficient. By the late 1980s the goal for efficient codes was in GFLOPS ($10^9$ FLOPS). The goal of parallel processing in the 1990s is TERAFLOPS ($10^{12}$ FLOPS) performance. Hence, absolute measures are useful for comparing two competing algorithms on the same architecture, or for comparing the implementation of an algorithm on different but contemporary architectures. Second, relative measures—such as speedup or efficiency—are indicative of the performance achieved by an algorithm on the specific parallel architecture. Hence, extrapolation to other architectures is meaningful only to the extent that the implementation on the alternative architecture does not alter the ratio of computation to communication time.

One important observation is emphasized in all experiments: we are examining the *relative* performance of a parallel algorithm implemented on a suitable parallel architecture, vis-à-vis the performance of the same, or a competing, algorithm on a contemporary high-performance workstation or mainframe. Such comparisons highlight the merits of parallel processing, and illustrate the size and complexity of models that can be solved efficiently using parallel computing.

The performance achieved today on a parallel machine may soon be possible on a high-performance workstation. But the parallel computer of the time will also exhibit a proportionate—if not higher—improvement in performance. The advantages of multiprocessor over uniprocessor computers (as we illustrate in this chapter) will be sustained, as both serial and parallel architectures are often based on the same underlying hardware technology.

We begin, in Section 15.1, with a discussion of guidelines for reporting computational experiments on parallel machines, and then proceed to report results with several algorithms and for several areas of application. The numerical results in the remaining sections of this chapter are all related to the subjects studied earlier in the book. Sections 15.2 to 15.5 report numerical results with the solution of models arising in matrix balancing, image reconstruction, network flow problems, and planning under uncertainty. Section 15.6 analyzes the performance of algorithms for linear programming problems, based on proximal minimization with $D$-functions. As an Appendix to this chapter we give in Section 15.7 a brief description of the machines used in the experiments. Notes and references in Section 15.8 conclude this chapter.

## 15.1 Reporting Computational Experiments on Parallel Machines

The task of reporting computational experiments is quite complex. Several factors affect the experimental design and influence the preparation of the report. The purpose of the computational experiment is one major factor. Is the test performed to establish the efficiency of a given machine? Or is it performed to establish the efficiency of a code? In the former case *time* may be the only relevant metric since all other resources of the machine (memory, input and output devices, etc.) are considered a fixed cost. In the latter case *work* is a more appropriate measure that incorporates not only execution time, but also time for data input and output, memory access, and so on.

Reporting the results from computational experiments becomes less ambiguous if a set of guidelines are followed. The following guidelines are based on the work of Barr and Hickman (1993), who suggested the development of reporting standards based on a survey of experts in the field of computational mathematical programming. We do not in this chapter follow the guidelines in all their detail when reporting computational results, since our goal is to provide a general view of the performance of the parallel algorithm, and we avoid details that are not essential. The articles cited in Section 15.8, from which the numerical results are obtained, follow the guidelines more faithfully.

**Guidelines for Reporting Computational Experiments on Parallel Machines**

### Document the Computational Experiment:

1. Describe the code: this includes the mathematical algorithm, any modifications to the algorithm, the design of the code, the data structures and the tuning parameters.
2. Describe the computing environment: this includes all pertinent characteristics of the machine, such as manufacturer, model, size of memories, type and number of processors.
3. Describe the communication environment: this includes the configuration of the communication network, interprocessor communication schemes, and technical characteristics of the network such as startup latency and transmission time per byte.
4. Describe the testing environment and methodology: this includes a precise definition of all performance metrics used, the data collection procedure for computing the metrics (for example, explain how time was measured), and the values of the tuning parameters.

### Use a Well-planned Experimental Design:

1. Focus on the time and cost required to solve problems that are both difficult and typical for the application or algorithm.

2. Identify those factors that affect the results, such as problem characteristics, values of tuning parameters, and amount of parallelism.
3. Provide points of reference by using well-known test problems and codes, even if these codes must be tested on a different, perhaps serial, machine.
4. Perform the reported results on a dedicated or lightly loaded system, and report the effect of the typical machine load on these results.
5. Employ techniques from statistical experimental design in order to highlight factors that contribute to the reported results.

**Provide a Comprehensive Report of the Results:**

1. Report results with all relevant performance metrics.
2. Use graphics where possible and when informative.
3. Use summary statistics such as central tendency and variability.
4. Provide as much detail as necessary.
5. Report results with failures that provide insight.

## 15.2 Matrix Balancing

In Chapter 9 we defined the matrix balancing problem and discussed some of its real-world manifestations. The size (extremely large) and particular structure of these problems make them suitable for parallel computations. Using the implementation techniques of Chapter 14 and the algorithms of Chapter 9, we show what can be achieved in practice. In particular, we report summary computational results with the parallel implementation of the RAS algorithm (Algorithm 9.3.2), the simultaneous RAS (Algorithm 14.3.1), and the range-RAS algorithm (Algorithm 9.3.1) on a variety of parallel architectures.

Test problems were derived from regional input/output accounts in the United States for 1977. Their characteristics are summarized in Table 15.1, and they are equality-constrained problems in the category of Problem 9.2.3.

**Table 15.1** Size characteristics of the matrix balancing test problems.

| Problem name | Number of rows × columns | Number of nonzero entries |
|---|---|---|
| USE537 | 504 × 473 | 57247 |
| MAKE537 | 523 × 533 | 9586 |

To provide a benchmark against which to evaluate the performance of the parallel algorithms we first run RAS on a VAX 8700 mainframe, a Floating Point Systems (FPS) M64/60 attached array processor, and a CRAY X-MP/48 vector supercomputer. The implementations on the FPS and the CRAY are coded to take advantage of the vector architecture of the machines. Results show that these large-scale problems would take several hours to solve on current mainframes, and several minutes on a vector supercomputer (Table 15.2). In the following sections we contrast this performance with the results achieved using alternative parallel computing methods.

**Table 15.2** Benchmark solution times with RAS in hrs:min:sec. Error is the largest difference between the target row or column sum, and the sum of the entries of the corresponding row or column.

| Problem name | Error | VAX 8700 | FPS M64/60 | CRAY X-MP/48 |
|---|---|---|---|---|
| USE537 | $10^{-6}$ | 0:00:10 | 0:00:1 | 0:00:00.3 |
| | $10^{-12}$ | 0:00:26 | 0:00:4 | 0:00:00.9 |
| MAKE537 | $10^{-4}$ | 0:00:30 | 0:00:11 | 0:00:06.9 |
| | $10^{-6}$ | 1:13:06 | 0:05:54 | 0:05:20.5 |
| | $10^{-8}$ | > 20:00:00 | 1:12:28 | 0:36:46.9 |

### 15.2.1 Data parallel implementations

Using the techniques developed in Section 14.3.1, the RAS algorithm was implemented on a Connection Machine CM-2 with 32K processing elements. Computational results compare the performance of both a sparse and a dense data parallel implementation of the algorithm, with the benchmark implementation on the CRAY X-MP vector computer (Table 15.3). A close look at these results tells us that, in general, the performance of the data parallel implementation is comparable to that of the serial implementation on the vector supercomputer, but it is much faster than implementation on the VAX 8700 mainframe. However, problems with a greater number of nonzero entries will be solved in the same time (per iteration) on the CM-2, while the execution time on the supercomputer and the mainframe will increase linearly with the number of nonzero entries in the matrix. We also observe from the results that these particular test problems are solved more efficiently on the CM-2, using the dense implementation. The dense code executes at a higher FLOPS rate than the sparse code, because it is implemented in a way that uses regular, and therefore very efficient, communications. Hence, the dense code solves these problems faster, even if it has to execute null operations on the zero entries of the matrix. As problems become even more sparse, however, this advantage will not be maintained. See, for example, the comparison between the dense and the sparse implementations for the row-action algorithm for transportation problems in Section 15.4.1.

**Table 15.3** Comparisons of execution times of data parallel implementations of RAS on a CM–2 with 32K processing elements (CM time in seconds) and control parallel implementation on an Alliant FX/8 with 8 processors (CPU times are in seconds), with the benchmark implementation on a single processor of the CRAY X-MP. Error is the largest difference between the target row or column sum, and the sum of the entries of the corresponding row or column.

| **Problem name** | **Error** | **Data parallel** | | **Control parallel** | **Benchmark** |
|---|---|---|---|---|---|
| | | (sparse) | (dense) | | |
| USE537 | $10^{-6}$ | 0.45 | 0.18 | 1.04 | 0.30 |
| MAKE537 | $10^{-4}$ | 11.41 | 6.13 | 7.50 | 6.90 |

**Table 15.4** Comparing the performance of range-RAS with RAS on a Connection Machine CM–2. Solution times are in seconds. Error is the largest difference between the target row or column sum, and the sum of the entries of the corresponding row or column. Iter. denotes number of iterations, where iteration is defined as the execution of an iterative step of the algorithm, once, over all rows and all columns of the matrix.

| **Problem** | **Range-RAS** | | | **RAS** | | |
|---|---|---|---|---|---|---|
| **name** | Iter. | Time | Error | Iter. | Time | Error |
| USE537 | 40 | 0.52 | $10^{-6}$ | 30 | 0.23 | $10^{-6}$ |
| MAKE537 | 1000 | 13.08 | $2.7 \times 10^{-5}$ | 1000 | 7.68 | $8.0 \times 10^{-6}$ |

There is an interesting aspect to the experimental results with the range-RAS (RRAS) algorithm (Algorithm 9.3.1) for interval-constrained problems. This algorithm can also be used to solve the equality-constrained test problems, simply by treating them as interval-constrained problems with identical values for both sides of the range constraints. Since RRAS is applicable to more general problem formulations than RAS, it is computationally more expensive. Table 15.4 compares the performance of RRAS with that of RAS when applied to the same equality-constrained test problems. RRAS is roughly twice as slow as RAS, but both algorithms solve these large-scale problems within a few seconds.

### 15.2.2 Control parallel implementations

We now examine the control parallel, as opposed to the data parallel, implementation of matrix balancing algorithms. The RAS algorithm and its simultaneous version (Algorithm 14.3.1) were implemented on an Alliant FX/8 with 8 processors. Each processor has vector functional units, and the implementations take advantage of the vector architecture.

**Table 15.5** Solution times of the parallel implementation of RAS on the Alliant FX/8 shared memory multiprocessor in CPU hrs:min:sec. Error is the largest difference between the target row or column sum, and the sum of the entries of the corresponding row or column. Iteration is defined as the execution of an iterative step of the algorithm, once, over all rows and all columns of the matrix.

| **Problem:** | **MAKE537** | **MAKE537** | **MAKE537** | **USE537** |
|---|---|---|---|---|
| **Error:** | $10^{-4}$ | $10^{-6}$ | $10^{-8}$ | $10^{-12}$ |
| **Iterations:** | 802 | 80057 | 342081 | 70 |
| **Scalar:** | 00:01:19.9 | 02:10:20.8 | 09:20:51.0 | 00:00:31.0 |
| **Vectorized:** | 00:00:51.3 | 01:23:13.4 | 05:55:15.9 | 00:00:18.5 |
| **2-CPUs:** | 00:00:25.7 | 00:41:37.6 | 02:59:33.9 | 00:00:09.5 |
| **3-CPUs:** | 00:00:17.5 | 00:28:13.0 | 02:00:46.9 | 00:00:06.4 |
| **4-CPUs:** | 00:00:13.5 | 00:21:39.1 | 01:31:54.9 | 00:00:04.9 |
| **5-CPUs:** | 00:00:11.0 | 00:17:32.9 | 01:14:57.7 | 00:00:04.1 |
| **6-CPUs:** | 00:00:09.5 | 00:14:56.1 | 01:03:37.7 | 00:00:03.5 |
| **7-CPUs:** | 00:00:08.3 | 00:13:04.0 | 00:55:39.6 | 00:00:03.1 |
| **8-CPUs:** | 00:00:07.5 | 00:11:45.0 | 00:49:45.9 | 00:00:02.7 |

Results with the control parallel implementation of RAS on the Alliant FX/8 are summarized in Table 15.5, where it is observed that the speedups achieved when using multiple processors are significant, and they increase almost linearly with the number of available processors. This is demonstrated more clearly in Figure 15.1(a) where the average relative speedup is computed for the solutions of MAKE537 (for three different levels of accuracy) and USE537, and it is plotted as a function of the number of processors.

The control parallel implementation of simultaneous RAS (with weights $\lambda_\rho = \lambda_\sigma = 0.8$) was also tested, and was compared with the uniprocessor implementation of RAS. The achieved speedup, shown in Figure 15.1(b), appears to be superlinear. A careful examination of the results indicates that simultaneous RAS is slightly more efficient than RAS even on a single processor. If we were to develop the speedup curve based on a comparison of the parallel simultaneous RAS with the uniprocessor implementation of the simultaneous RAS we would obtain a curve identical to that of Figure 15.1(a). The results of Figure 15.1(b) are of interest in their own right, since it is the experimentation with parallel computing that revealed that simultaneous RAS can be a more efficient algorithm, even when executing

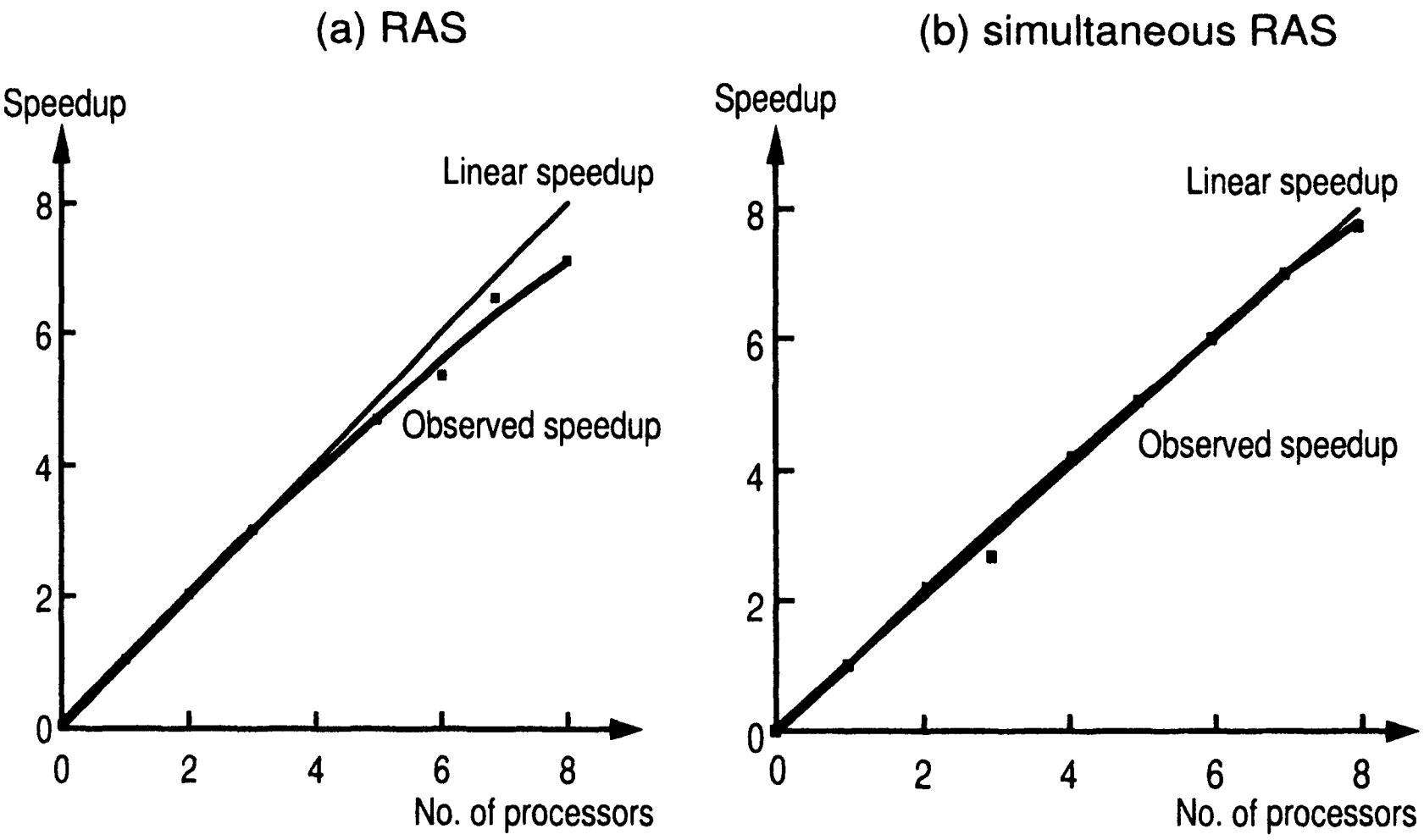

**Figure 15.1** Relative speedup of (a) RAS, and (b) simultaneous RAS, with increasing number of processors on the Alliant FX/8 shared memory multiprocessor.

serially.

## 15.3 Image Reconstruction

Problems in image reconstruction were introduced in Chapter 10, where we also mentioned the numerous real-world settings where such problems arise. Models arising in the discrete setting can be huge, as finer discretizations are introduced so as to achieve more detailed reconstructions. In this section we illustrate the performance of parallel computing in solving large-scale instances of such problems, using an example from medical imaging. The algorithm we use is block-MART (Algorithm 6.7.1), for which alternative parallel implementation schemes were developed in Section 14.4. Block-MART was implemented as the reconstruction algorithm within the SNARK77 image reconstruction testing software environment on a CRAY X-MP/48 with 4 processors. (For details on the implementation see Section 14.4.) Each processor of the CRAY has vector processing capabilities, and the algorithm was coded in a way that exploits the vector architecture.

The algorithm was tested using a series of standard problems derived from the Shepp and Logan (1976) test phantom of a head section (see Table 15.6). In Table 15.7 both the scalar and the vectorized code on a CRAY are compared to the benchmark performance of the code on a SUN-3/160c workstation. The improvement due to vectorization is approximately tenfold, while a 100-fold reduction in solution time is achieved with the use of the vector supercomputer, as opposed to a contemporary workstation.

**Table 15.6** Size characteristics of the image reconstruction test problems from the Shepp and Logan phantom of a head section.

| Problem name | Dimensions (pixels) | Approximate no. of variables | Approximate no. of equations |
|---|---|---|---|
| SHEPP64 | 64 × 64 | 4,000 | 11,000 |
| SHEPP512 | 512 × 512 | 200,000 | 92,000 |
| SHEPP1024 | 1024 × 1024 | 1,050,000 | 185,000 |

**Table 15.7** Performance of block-MART on a CRAY X-MP with and without vector computations, and benchmark performance on a SUN-3 workstation. Times are in CPU seconds. NA means not available.

| Problem name | CRAY X-MP/48 | | SUN-3/160c |
|---|---|---|---|
| | Scalar | Vectorized | |
| SHEPP64 | 496.17 | 81.19 | 7,842.0 |
| SHEPP512 | 1140.49 | 108.81 | NA |
| SHEPP1024 | NA | 312.71 | NA |

Block-MART is a block-iterative algorithm, as characterized in Section 1.3.1, and as such affords users great flexibility in grouping equations during an implementation. The parallel implementation for the specific test problems solved in this section was carried out with 45 blocks per view. This implementation is used to evaluate both the parallel scheme that exploits parallelism within a block (Sections 14.4.1), and the scheme that exploits parallelism with independent blocks (Section 14.4.2). The relative speedups achieved with the two parallel schemes are presented in Table 15.8, where it is shown that both schemes achieve more than 75 percent efficiency. Further comparison with the benchmark implementation on the SUN-3/160c reveals that a parallel vector computer yields improvements of more than two orders of magnitude over a modern workstation.

The implementation of the block-iterative algorithm block-MART allows us to address some interesting issues that arise when implementing iterative reconstruction algorithms. How does the parallel efficiency of a row-action MART compare with that of the block-iterative algorithm block-MART? Which block sizes give the best performance? Such questions can best be addressed through practical experimentation. The results shown here do not provide complete answers, but demonstrate the effect of different block sizes on the practical performance of the algorithms.

We compare the parallel efficiency of MART with that of block-MART by repeating the experiments with parallel independent blocks (see Sec-tion 14.4.2) for problem SHEPP512, and varying the block sizes. With 512 blocks per view the block-MART algorithm becomes MART, while with one block per view we have the extreme case of block-MART wherein each view is a single block. Results summarized in Figure 15.2 show that intermediate block sizes (30 blocks per view for this particular problem) are optimal. With this block size the blocks are large enough to provide large tasks and avoid excessive synchronization delays and, at the same time, there are enough tasks available to achieve load balancing between the processors. While the precise block size will vary by problem and will also be affected by the characteristics of the parallel machine, we can conclude from this experiment that the optimal block size is neither one equation per block, nor all equations of a single view in a block, but that some balance must be found. This conclusion holds true for other block-iterative algorithms as well. Experimental results reported in the literature (see Section 15.8) show similar behavior with the block versions of the EM algorithm (discussed in Section 10.2.1) and the variable block-Kaczmarz algorithm (Algorithm 10.4.1).

**Table 15.8** Relative speedups with the alternative parallel implementations of block-MART of Sections 14.4.1 and 14.4.2 on the CRAY X-MP/48 vector supercomputer with up to 4 processors.

| No. of processors | Parallelism within a block | | Parallelism with independent blocks | |
|---|---|---|---|---|
| | SHEPP512 | SHEPP1024 | SHEPP512 | SHEPP1024 |
| 2 | 1.82 | 1.80 | 1.82 | 1.79 |
| 3 | 2.49 | 1.97 | 2.32 | 1.93 |
| 4 | 2.97 | 3.15 | 3.01 | 3.11 |

## 15.4 Multicommodity Network Flows

In Chapter 14 we introduced several formulations of network flow problems and discussed examples of real-world settings where such problems arise. Algorithms for several formulations were developed, either based on the row-action framework of Chapter 6 or on the model decomposition framework of Chapter 7. In this section we examine results with numerical experiments using row-action algorithms for transportation and multicommodity network flow problems, and with a model decomposition algorithm for multicommodity network problems. Our goal is to illustrate the effectiveness with which each algorithm exploits parallelism, and the efficiency with which large-scale problems can be solved.

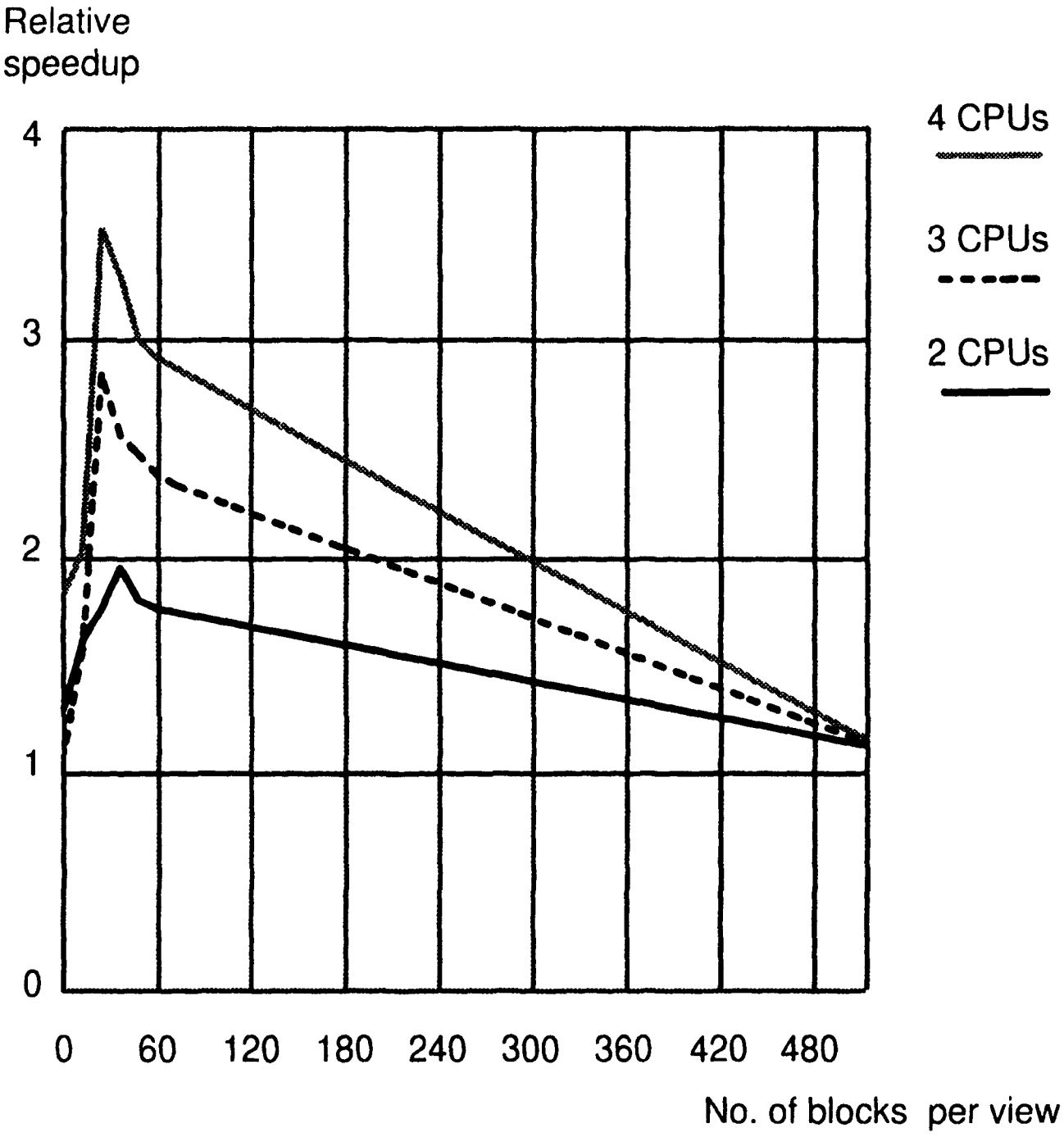

**Figure 15.2** Relative speedup with varying block sizes with the block-MART algorithm. With 512 blocks per view we recover MART. With one block per view we have the extreme case implementation of block-MART wherein all rays in a view are blocked together.

### 15.4.1 Row-action algorithm for transportation problems

The row-action algorithm (Algorithm 12.4.1) for quadratic transportation problems was implemented on a Connection Machine CM–2 using the techniques of Section 14.5. We use this implementation to illustrate the efficiency of data parallel computing for transportation problems, and to compare the performance of the data parallel row-action algorithm with that of other serial and parallel algorithms that have been published in the literature. Table 15.9 gives the size characteristics of the test problems used in these experiments.

These test problems were solved using both the dense and the sparse implementations discussed in Sections 14.5.1 and 14.5.2. An optimal implementation for dense problems executes at 1.5 GFLOPS on a Connection Machine CM–2 with 32K processing elements, while the sparse implementation executes at 275 MFLOPS on the same machine. The performance of both implementations can be seen in Table 15.10, where it is evident that, for the very sparse problems solved here, the sparse code outperforms

**Table 15.9** Size characteristics of the quadratic transportation problems.

| Problem name | Number of nodes (origin × destination) | Number of arcs |
|---|---|---|
| TSENG1 | 500 × 500 | 5063 |
| TSENG2 | 750 × 750 | 7611 |
| TSENG3 | 1000 × 1000 | 10126 |
| TSENG4 | 1250 × 1250 | 12665 |
| TSENG5 | 500 × 500 | 10051 |
| TSENG6 | 750 × 750 | 15086 |
| TSENG7 | 1000 × 1000 | 20134 |
| TSENG8 | 1250 × 1250 | 25153 |

**Table 15.10** Comparing implementation of a row-action algorithm for dense problems on a Connection Machine CM–2 with 32K processing elements with implementation for sparse problems on a CM–2 with 16K processing elements. Solution times are in CM seconds.

| Problem | Dense code | Sparse code |
|---|---|---|
| TSENG1 | 1.91 | 1.55 |
| TSENG2 | 6.52 | 2.45 |
| TSENG3 | 5.43 | 2.03 |
| TSENG4 | 18.64 | 2.06 |
| TSENG5 | 1.53 | 1.63 |
| TSENG6 | 8.71 | 3.43 |
| TSENG7 | 5.43 | 3.45 |
| TSENG8 | 15.05 | 2.85 |

the dense code, even if the latter has a much higher computing rate, as expressed in FLOPS.

The performance of the data parallel row-action algorithm on the Connection Machine CM–2 with 32K processors was compared to the control parallel implementation of a related relaxation algorithm (see Section 15.8 for references on relaxation algorithms) on an Alliant FX/8 with eight processors, and to the performance of the same algorithm on a workstation. From the results of this experiment we observe that the data parallel implementation of the row-action Algorithm 12.4.1 on the CM–2 substantially outperforms the serial implementation of the relaxation algorithm on a workstation, and it is also somewhat faster than the control parallel implementation of the relaxation algorithm on the MIMD machine (Table 15.11).

Finally, three extremely large and dense problems were solved using the data parallel row-action algorithm and two state-of-the-art serial algorithms on an IBM 3090 vector computer. This experiment aims to illustrate the efficiency with which a parallel algorithm implemented on suitable architecture can solve extremely large problems, in contrast to serial algorithms on contemporary mainframes. In particular, we compare the data parallel Algorithm 12.4.1 executing on the Connection Machine CM–2, to an equilibration algorithm and an exact Newton method (see Section 15.8 for references on these algorithms) implemented on an IBM 3090 mainframe. The results show that the massively parallel implementation of the row-action algorithm is up to two orders of magnitude faster than the implementation of state-of-the-art serial algorithms (see Table 15.12).

**Table 15.11** Solution times (in seconds) for transportation problems with the row-action Algorithm 12.4.1 on a CM–2 with 32K processing elements, and with a relaxation algorithm on an Alliant FX/8 with 8 processors and a MicroVAX-II.

| **Problem name** | **32K CM–2** | **Alliant FX/8** (8 processors) | **MicroVAX-II** |
|---|---|---|---|
| TSENG1 | 1.55 | 3.01 | 27.07 |
| TSENG2 | 1.47 | 14.47 | 42.40 |
| TSENG3 | 1.22 | 7.31 | 58.18 |
| TSENG4 | 1.23 | 9.24 | 70.67 |
| TSENG5 | 0.98 | 4.33 | 54.40 |
| TSENG6 | 2.05 | 6.61 | 78.92 |
| TSENG7 | 2.07 | 9.15 | 106.24 |
| TSENG8 | 1.71 | 11.41 | 135.62 |

**Table 15.12** Comparing the row-action Algorithm 12.4.1 on the CM–2 with an equilibration algorithm and an exact Newton algorithm on the IBM 3090 for the solution of dense quadratic optimization problems. Times are in CM and CPU hrs:min:sec, respectively.

| **Problem size** origin × destination | **Row-action method 32K CM–2** | **Equilibration method IBM 3090** | **Exact Newton method IBM 3090** |
|---|---|---|---|
| 500 × 500 | 0:00:18 | 0:01:09 | 0:30:00 |
| 1000 × 1000 | 0:00:22 | 0:08:03 | NA |
| 2000 × 2000 | 0:00:47 | 1:03:43 | NA |

**Table 15.13** Solving large quadratic multicommodity transportation problems on a Connection Machine CM–2 with 32K processing elements. (Times are in min:sec). Error denotes the largest % violation of the GUB constraints, and the largest violation of flow conservation constraints at the nodes.

| **Problem size** | | **Error** | | **Time** |
|---|---|---|---|---|
| Network formulation | Linear programming formulation | % GUB constraints | Absolute node constraints | CM–2 with 32K PEs |
| 5 commod's.<br>1,024 nodes<br>262,000 arcs | 5120 eqns.<br>26,181 GUB<br>1,310,000 vars. | .09 | $10^{-6}$ | 2:15 |
| 10 commod's.<br>1,024 nodes<br>262,144 arcs | 10,240 eqns.<br>25,957 GUB<br>2,621,140 vars. | .09 | $10^{-6}$ | 4:55 |
| 20 commod's.<br>1,024 nodes<br>262,144 arcs | 20,480 eqns.<br>26,260 GUB<br>5,424,880 vars. | .09 | $10^{-5}$ | 7:20 |
| 5 commod's.<br>2,048 nodes<br>1,048,074 arcs | 10,240 eqns.<br>104,088 GUB<br>5,240,370 vars. | .07 | $10^{-6}$ | 6:20 |
| 8 commod's.<br>2,048 nodes<br>1,048,570 arcs | 16,384 eqns.<br>104,521 GUB<br>8,384,560 vars. | .09 | $10^{-6}$ | 4:40 |

### 15.4.2 Row-action algorithm for multicommodity transportation problems

The row-action algorithm for multicommodity network flow problems (Algorithm 12.4.3) was implemented using both the techniques of Section 14.5.4 and control-level parallelism. The experiments presented in Table 15.13 illustrate the efficiency with which the data parallel implementation on a Connection Machine CM–2 with 32K processing elements solves problems with millions of variables. For example, the results show that such problems with millions of variables are solved within 4–10 minutes.

A control parallel implementation of the same row-action algorithm was tested on an Alliant FX/8. The algorithm was parallelized by distributing the multiple commodities across multiple processors. The iterative Steps 1.1–1.3 of Algorithm 12.4.3 that operate on rows of the constraint set for different commodities were executed concurrently by multiple processors. Upon completion of one iteration over all rows and all commodities, the multiple processors cooperated to compute the projection operation on the generalized upper bound (GUB) constraints (Step 2). The results illustrate that good relative speedups are achieved on machines with four and eight processors (see Figure 15.3). As expected, the relative speedups increase for problems with more commodities. Efficiency of more than 80 percent is achieved on the larger machine (eight processors) for a sufficiently large number of commodities. We also learn from Figure 15.3 that the relative speedup does not monotonically increase with the number of commodities. For example the speedup with 10 commodities is lower than that achieved with eight commodities on the eight-processor system. Clearly, the number of commodities must be a multiple of the number of processors to ensure good load balancing and achieve high efficiency.

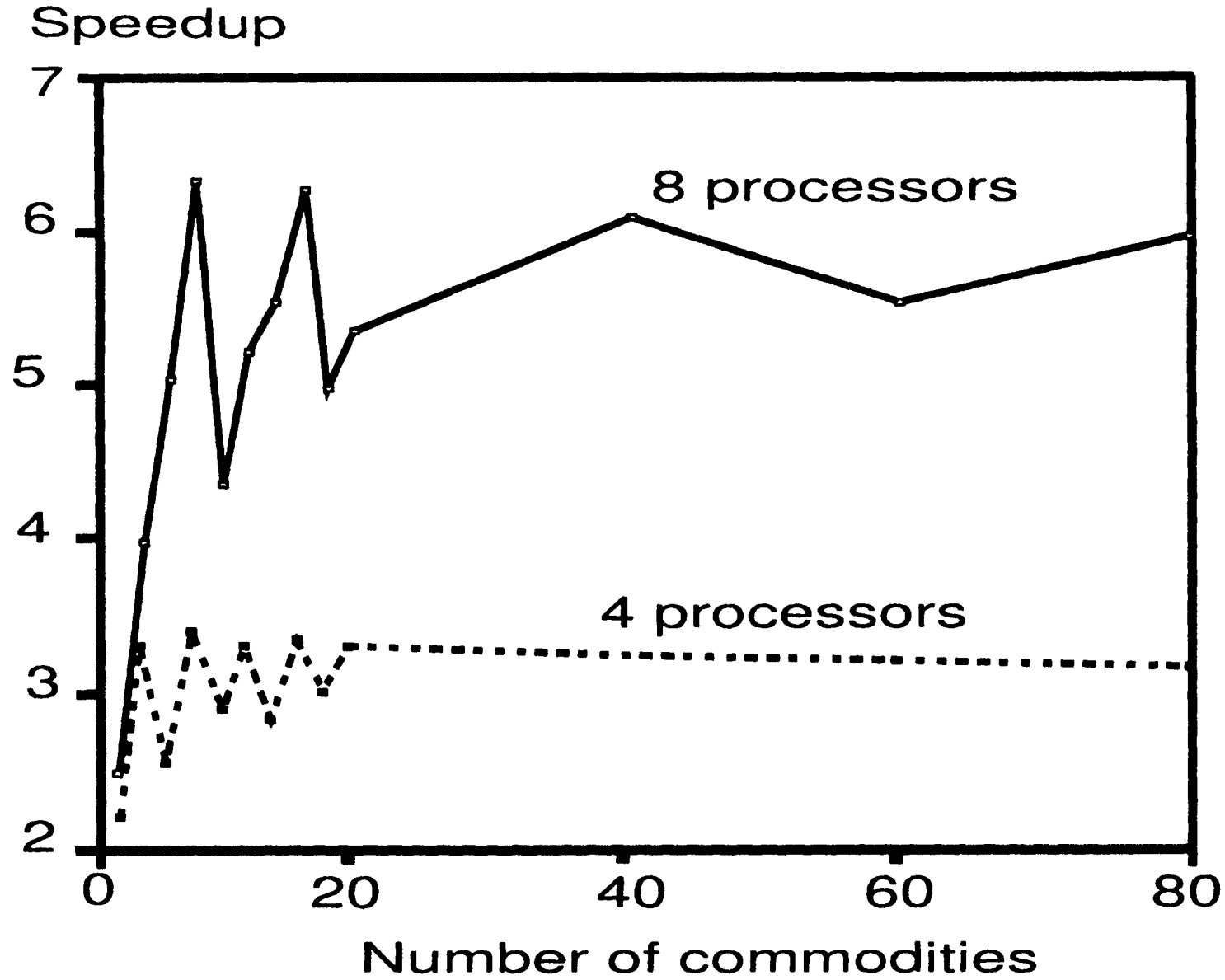

**Figure 15.3** Relative speedups with the control parallel implementation of the row-action Algorithm 12.4.3 for multicommodity network flow problems on an Alliant FX/8 shared memory multiprocessor.

### 15.4.3 Linear-quadratic penalty (LQP) algorithm for multicommodity network flow problems

A model decomposition algorithm based on the linear-quadratic penalty function was developed for the multicommodity network flow problem in Section 12.5. This algorithm was implemented using control parallelism on an MIMD machine, and was tested on a series of large-scale problems arising in military logistics applications. (Experimental results with a data parallel implementation of the same algorithm are reported in the litera-ture; see Section 15.8 for references). The test problems, referred to as PDS problems, arise in a system for the distribution of patients from war theaters to hospitals. All test problems have 11 commodities. The size of the underlying network depends on the planning horizon $\tau$ and models with different time horizons are denoted by PDS$\tau$. Table 15.14 gives the sizes of the test problems.

**Table 15.14** Size characteristics of the patient distribution system (PDS) multicommodity network flow problems. All problems have 11 commodities.

| **Problem name** | **Network formulation** | | **Linear programming formulation** |
|---|---|---|---|
| | Arcs | Nodes | |
| PDS1 | 339 | 126 | 1,473 × 3,816 |
| PDS3 | 1,117 | 390 | 4,593 × 12,590 |
| PDS5 | 2,149 | 686 | 7,546 × 23,639 |
| PDS10 | 4,433 | 1,399 | 15,389 × 48,763 |
| PDS15 | 7,203 | 2,125 | 23,375 × 79,233 |
| PDS20 | 10,116 | 2,447 | 31,427 × 105,728 |
| PDS30 | 15,126 | 4,223 | 46,453 × 154,998 |

The problems were solved on a CRAY Y-MPE246 vector supercomputer with eight processors. The implementation was vectorized using the basic linear algebra subprograms from the BLAS library (see Section 15.8). The solutions of the subproblems were parallelized by distributing single-commodity subproblems for distinct commodities to different processors for concurrent execution. The master problem mainly involves dense matrix operations and those were also implemented in a way that takes advantage of parallelism. In Table 15.15, the executions of the algorithm on a CRAY and on a VAX 6400 mainframe are reported. We notice that speedups of up to four are achieved when using the eight processors of the CRAY. The reduced efficiency, less than 50 percent, is due to the serial bottleneck in the algorithm, due to the presence of a master program and poor load balancing when solving the subproblems. As can be seen in Figure 15.4, the number of pivots used by the network simplex algorithm for solving the linear subproblems differs substantially among the subproblems since some subproblems are much more difficult than others, which contributes to load imbalance among processors. Nevertheless, the results of Table 15.15 show that the parallel algorithm on the supercomputer is substantially faster than the same algorithm executing on a contemporary mainframe.

How does the parallel LQP algorithm compare with a state-of-the-art general purpose optimizer? This is an interesting question, since the LQP algorithm exploits both the structure of the problem and the computer architecture, and we would like to know the advantage of doing so, as opposed to solving the problems using a general purpose optimizer on a general purpose workstation. The LQP algorithm was benchmarked against the state-of-the-art linear programming code OB1. The results in Figure 15.5 tell us that the LQP algorithm outperforms the general purpose linear programming code, even when both are executing on a single processor. The advantage of LQP is further increased with the use of parallelism.

**Table 15.15** Performance of the LQP algorithm (Algorithm 12.5.1) on the patient distribution system (PDS) multicommodity network flow problems. Solution times in seconds. NA means not available.

| Problem name | VAX 6400 | CRAY Y-MP | |
|---|---|---|---|
| | | 1 CPU | 8 CPU |
| PDS1 | 112.3 | 1.8 | NA |
| PDS3 | 1006.5 | 18.7 | NA |
| PDS5 | 6965.6 | 93.0 | 22.9 |
| PDS10 | NA | 408.1 | 95.9 |
| PDS15 | NA | 940.3 | 307.7 |
| PDS20 | NA | 1946.8 | 740.0 |
| PDS30 | NA | 7504.0 | 2566.0 |

## 15.5 Planning Under Uncertainty

The problem of planning under uncertainty is considered one of the major outstanding problems in management science and operations research. One reason for this lies in the significant difficulty inherent in developing appropriate models (such as the stochastic linear programming and robust optimization models of Chapter 13). These models can be extremely large, their size growing linearly with the number of scenarios and exponentially with the number of stages. Such problems, together with the algorithms introduced in Chapters 8 and 13, are ideal for demonstrating the performance of parallel machines. In this section we inspect computational results from the use of an interior point algorithm for solving stochastic linear programs, and from use of row-action algorithm for solving nonlinear stochastic network problems.

### 15.5.1 Interior point algorithm

The primal-dual path following interior point algorithm discussed in Chapter 8 was implemented using the matrix factorization techniques discussed in Section 8.3. The implementation takes advantage of the dual block-angular structure of the stochastic programming problems, and can also be implemented on a parallel machine (see Section 14.6).

We use two models, SCSD8 and SEN, increasing the number of scenarios to generate problems of extremely large size. The SCSD8 is a stochastic ver-sion of a model for finding the minimal design of a multistage truss. Problem SEN is a telecommunication network design problem. Problem SCSD8 has a first-stage constraints matrix with dimensions $m_0 = 10$, $n_0 = 70$, and a second-stage constraints matrix with dimensions $m_1 = 20$, $n_1 = 140$. SEN has dimensions $m_0 = 1$, $n_0 = 90$, $m_1 = 175$, $n_1 = 795$. Large-scale problems are generated by increasing the number of scenarios. Table 15.16 summarizes the size characteristics of the test problems. To put the size of these problems in perspective we mention that the largest version of SCSD8 solved in the published literature has 32 scenarios, while the largest SEN problem solved to date has 1000 scenarios. The test problems solved here have over 130,000 scenarios for SCSD8 and over 16,000 scenarios for SEN.

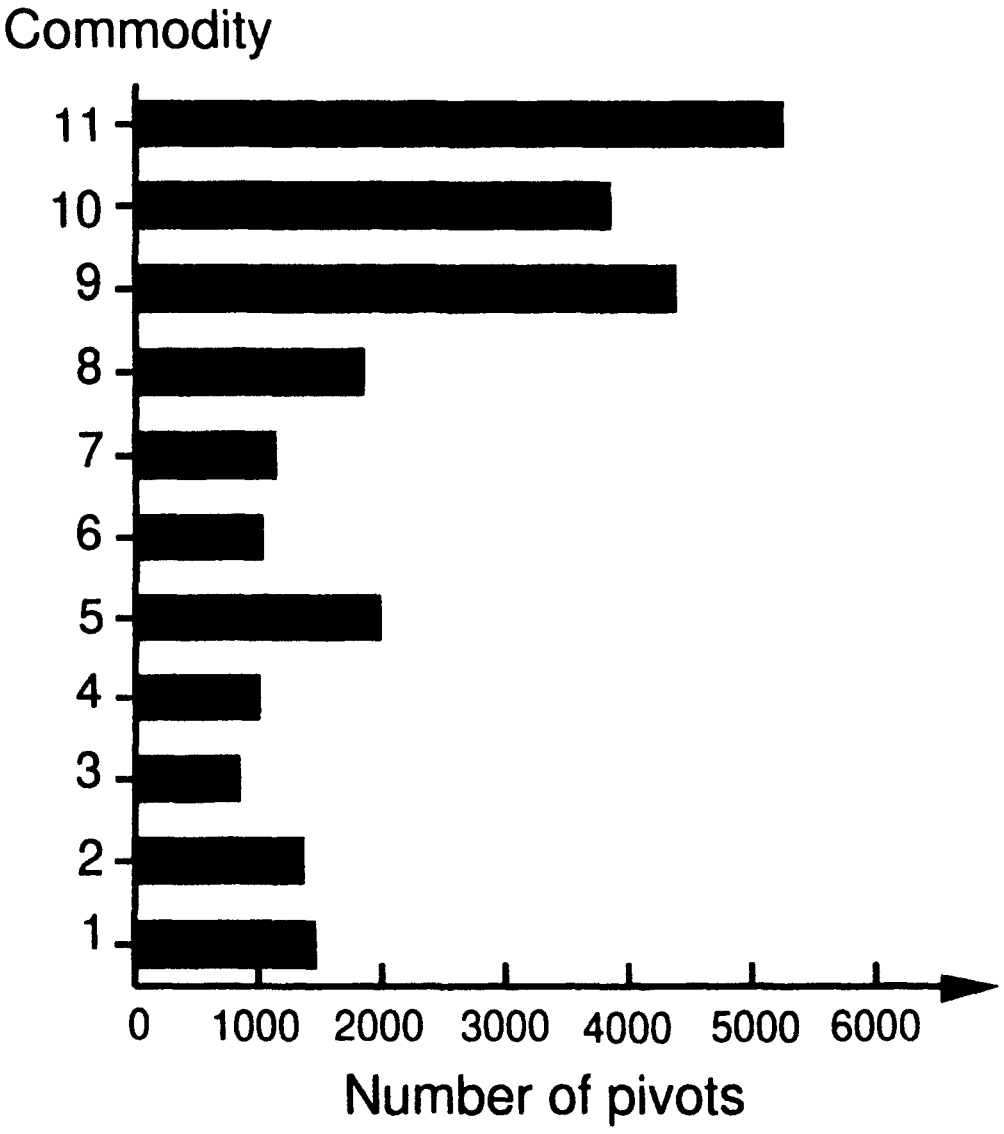

**Figure 15.4** Number of network simplex pivots for solving each commodity subproblem, at a sample iteration of the LQP algorithm during the solution of problem PDS20.

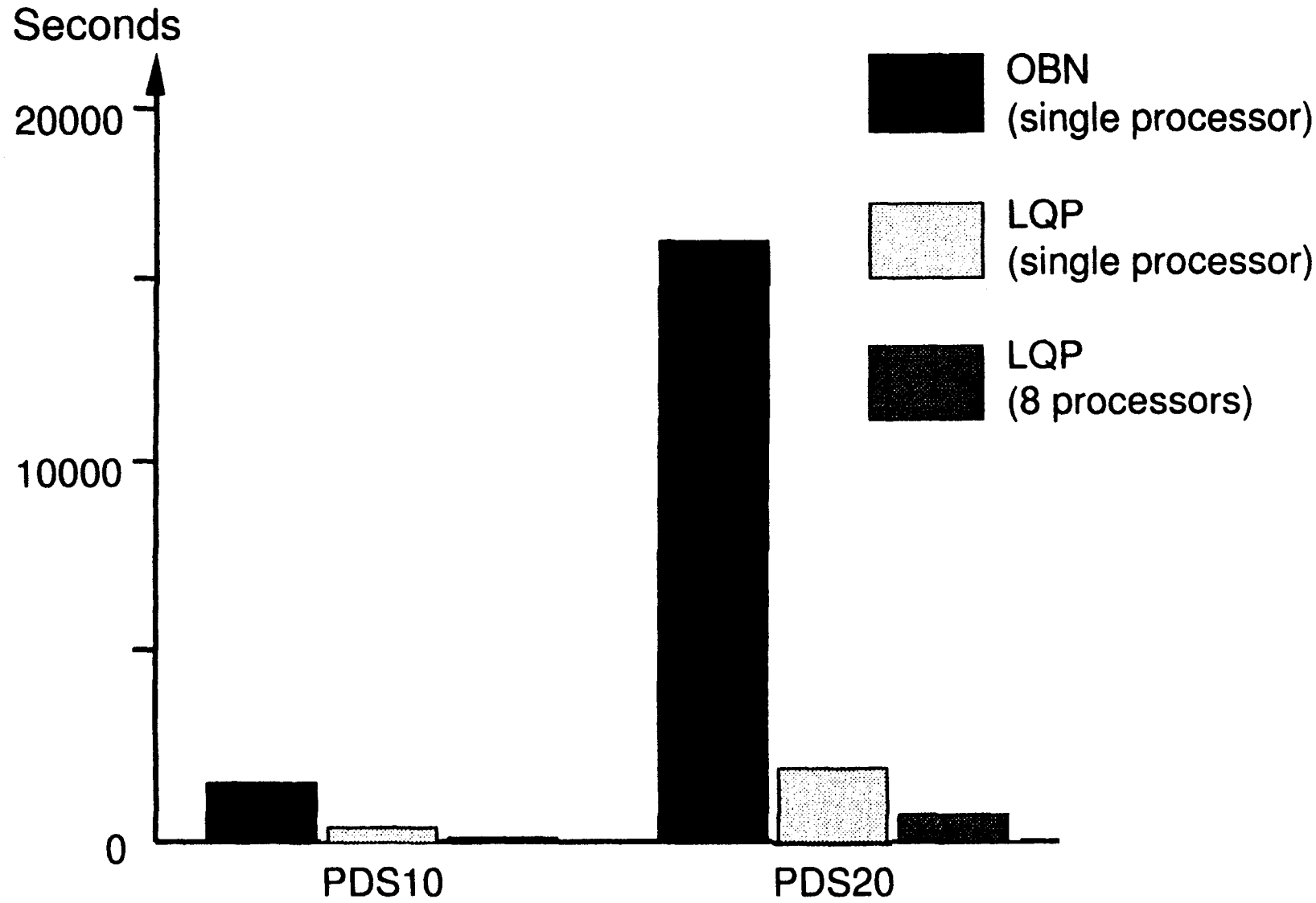

**Figure 15.5** Comparison of the linear-quadratic penalty LQP algorithm (Algorithm 12.5.1) with the general purpose linear programming code OB1 for the solution of multicommodity network flow problems. Both codes run on a CRAY Y-MPE.

**Table 15.16** Size characteristics of the stochastic programming test problems.

| Problem name | Number of scenarios | Number of constraints | Number of variables |
|---|---:|---:|---:|
| SCSD8.4 | 4 | 90 | 630 |
| SCSD8.8 | 8 | 170 | 1,190 |
| SCSD8.16 | 16 | 330 | 2,310 |
| SCSD8.32 | 32 | 650 | 4,550 |
| SCSD8.64 | 64 | 1,290 | 9,030 |
| SCSD8.128 | 128 | 2,570 | 17,990 |
| SCSD8.256 | 256 | 5,130 | 35,910 |
| SCSD8.512 | 512 | 10,250 | 71,750 |
| SCSD8.1024 | 1,024 | 20,490 | 143,360 |
| SCSD8.2048 | 2,048 | 40,970 | 286,790 |
| SCSD8.130172 | 130,172 | 2,603,440 | 18,224,080 |
| SEN.4 | 4 | 701 | 2,913 |
| SEN.8 | 8 | 1,401 | 5,737 |
| SEN.16 | 16 | 2,801 | 11,385 |
| SEN.32 | 32 | 5,601 | 22,681 |
| SEN.64 | 64 | 11,201 | 45,273 |
| SEN.128 | 128 | 22,401 | 90,457 |
| SEN.256 | 256 | 44,801 | 180,825 |
| SEN.512 | 512 | 89,601 | 361,561 |
| SEN.16384 | 16,384 | 2,867,201 | 13,025,369 |

The test problems were solved on a Connection Machine CM–5e with up to 64 processors. Figure 15.6 illustrates that very good relative speedups are achieved, especially for problems with many scenarios. The speedups are very good for machines with up to 64 processors (used in this experiment) but the upward sloping curves indicate that greater speedups can be expected with even larger machines. Figure 15.6 also reveals that superlin-ear speedup is achieved when solving SCSD8 on machines with 32 and 64 processors. This phenomenon occurs because of memory caching effects: since SCSD8 blocks are small and a few blocks can fit in the cache (fast) memory; large machines, however fit more blocks in the fast cache memory thus improving the execution time for linear algebra calculations. For the SEN problems each block is large so that it cannot fit in cache memory. Thus the speedup achieved with the SEN test problems does not benefit from caching and it is bounded by the number of processors.

**Table 15.17** Scalability of the parallel interior point algorithm on the Connection Machine CM–5e. The solution times (in seconds) remain virtually constant—when adjusting for changes in the number of iterations—as the number of processors increases to match the number of scenarios.

| Problem name | Number of processors | Number of iterations | Solution time |
|---|---|---|---|
| SCSD8.4 | 4 | 9 | 1.19 |
| SCSD8.32 | 32 | 9 | 1.22 |
| SCSD8.64 | 64 | 9 | 1.22 |
| SEN.4 | 4 | 18 | 13.6 |
| SEN.32 | 32 | 19 | 14.5 |
| SEN.64 | 64 | 21 | 16.1 |

How scalable is the parallel implementation of the interior point algorithm when implemented on machines with even more processors to solve larger problems? This is an important issue for stochastic programming as problems can have thousands of scenarios that potentially could be solved on massively parallel computers. To address this question (termed scalability, Definition 1.4.7) for the parallel interior point algorithms, we investigated solutions of problems with increasing number of scenarios on the Connection Machine CM–5, using as many processors as there are scenarios for each test problem. The results summarized in Table 15.17 show that the algorithm is perfectly scalable. For example, the solution time for problem SCSD8 remains virtually constant, as problems with increasing number of scenarios are solved using an equal number of processors. The same is true for the SEN problems, when we adjust for the fact that for larger test problems the algorithm takes more iterations to reach optimality.

Finally, we assess the potential of parallel computing to solve the extremely large-scale problems arising in stochastic programming, when compared with current state-of-the-art serial codes. The test problems were solved with both the parallel interior point algorithm and one of the best serial codes currently available for large-scale linear programming, called LOQO. This code implements the same primal-dual path following interior point algorithm used in the parallel implementation. The difference between the two codes is the way in which the linear systems of equations are solved in calculating the dual step; our code uses the parallel procedures developed in Chapter 8. The results in Table 15.18 provide a benchmark for the performance of the parallel interior point algorithm against a state-of-the-art serial code.

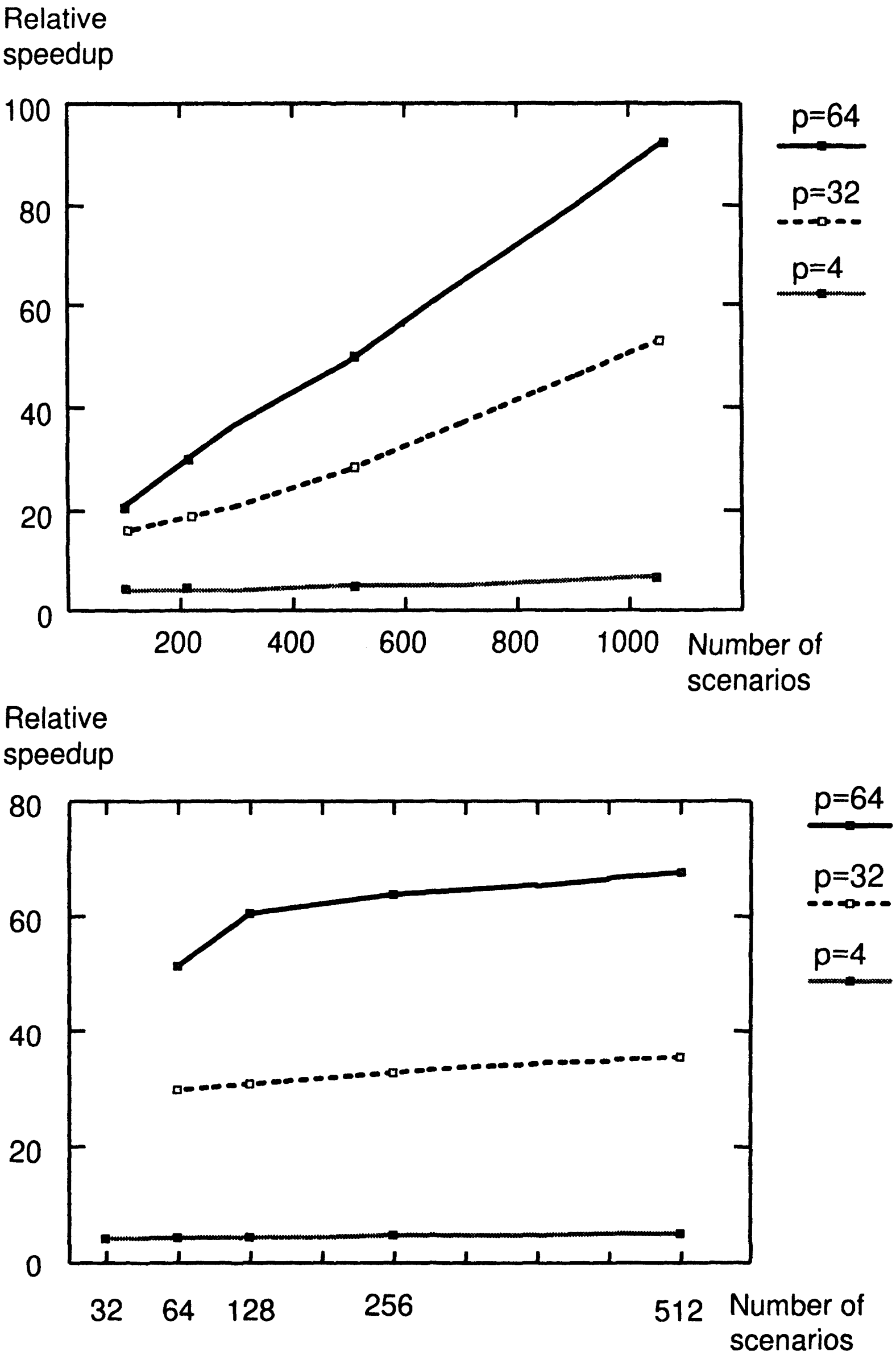

**Figure 15.6** Relative speedup of the parallel interior point code implemented on the Connection Machine CM–5e with p=4, 32 and 64 processors, for the solution of the SCSD8 (top figure) and SEN (bottom figure) test problems.

**Table 15.18** Comparing LOQO, executing on a single processor of the Connection Machine CM–5e, with the the parallel interior point algorithm. Time is given in seconds. The number of processors used by the parallel algorithm is equal to the number of scenarios up to a maximum of 64 (size of machine). NA means not available.

| Problem name | LOQO | | Parallel code | |
|---|---|---|---|---|
| | Iterations | Time | Iterations | Time |
| SCSD8.8 | 9 | 1.45 | 9 | 5.54 |
| SCSD8.16 | 9 | 3.37 | 9 | 6.03 |
| SCSD8.32 | 9 | 5.92 | 9 | 5.45 |
| SCSD8.64 | 9 | 12.8 | 9 | 5.42 |
| SCSD8.128 | NA | NA | 10 | 1.88 |
| SCSD8.256 | NA | NA | 10 | 2.39 |
| SCSD8.512 | NA | NA | 11 | 3.79 |
| SCSD8.1024 | NA | NA | 12 | 6.73 |
| SCSD8.2048 | NA | NA | 14 | 13.82 |
| SEN.8 | 14 | 37.31 | 19 | 25.1 |
| SEN.16 | 16 | 188.1 | 19 | 48.5 |
| SEN.32 | 17 | 837.2 | 20 | 14.5 |
| SEN.64 | 19 | 1702.1 | 21 | 16.1 |
| SEN.128 | NA | NA | 23 | 30.8 |
| SEN.256 | NA | NA | 31 | 78.3 |
| SEN.512 | NA | NA | 31 | 153.5 |
| SEN.16384 | NA | NA | 49 | 7638.3 |

### 15.5.2 Row-action algorithm for nonlinear stochastic networks

For stochastic programming problems with network recourse, such as the financial planning problems introduced in Section 13.7, we have developed specialized algorithms using the row-action framework. In this section we present results for the solution of stochastic network problems using Algorithm 13.8.1, for which a data parallel implementation was developed using the techniques of Section 14.5.4. As test problems we use a set of asset allocation models, whose characteristics are given in Table 15.19. All test problems used here have a quadratic objective function, thus Algorithm 13.8.1 is directly applicable. Section 15.6.2 presents further results for the solution of linear stochastic network problems, using proximal minimization with $D$-functions.

**Table 15.19** Size characteristics of the stochastic network test problems.

| Problem name | Scenarios | Nodes | Arcs | Linear prog. formulation |
|---|---|---|---|---|
| | | (per scenario) | | (rows × columns) |
| DETER0 | 18 | 121 | 334 | $1,923 \times 5,468$ |
| DETER1 | 52 | 121 | 335 | $5,527 \times 15,737$ |
| DETER2 | 80 | 91 | 249 | $6,095 \times 17,313$ |
| DETER3 | 72 | 121 | 335 | $7,647 \times 21,777$ |
| DETER4 | 70 | 61 | 163 | $3,235 \times 9,133$ |
| DETER5 | 48 | 121 | 335 | $5,103 \times 14,529$ |
| DETER6 | 40 | 121 | 335 | $4,255 \times 12,113$ |
| DETER7 | 60 | 121 | 335 | $6,375 \times 18,153$ |
| DETER8 | 36 | 121 | 335 | $3,831 \times 10,905$ |

**Table 15.20** Solution time (CM seconds) and number of major and minor iterations for quadratic stochastic network test problems on the Connection Machine CM–2 with 8K processing elements using a row-action algorithm.

| Problem name | Major iterations | Minor iterations | Time |
|---|---|---|---|
| DETER0 | 9 | 850 | 10.5 |
| DETER1 | 9 | 875 | 19.0 |
| DETER2 | 7 | 625 | 22.0 |
| DETER3 | 9 | 875 | 31.4 |
| DETER4 | 5 | 325 | 6.7 |
| DETER5 | 10 | 950 | 24.7 |
| DETER6 | 9 | 900 | 19.4 |
| DETER7 | 11 | 1025 | 22.3 |
| DETER8 | 11 | 1050 | 22.7 |

The problems were solved on a Connection Machine CM–2 with 8K processors. The algorithm was terminated when both the maximum error in the conservation of flow equation at all nodes was below a small tolerance $\epsilon > 0$, and each non-anticipativity constraint was violated by less than $\epsilon$, for $\epsilon = 10^{-3}$. Results are shown in Table 15.20, where each projection on the non-anticipativity constraints is called a *major iteration*, and an iterations for the solution of the scenario subproblems is called a *minor* iteration. The results show that these large-scale problems were solved in less than one-half minute, and required only very few major iterations. A careful inspection of the results of Table 15.20 brings up an interesting observation, namely, that solution time depends neither on the number of scenarios nor on the size of the deterministic equivalent linear program. For example, DETER4 was solved much more quickly than DETER6 although it has twice as many scenarios. This type of behavior occurs because as problems get larger the algorithm can utilize a greater number of processors to execute the iterative steps, whereas the number of iterations required to reach a solution does not necessarily depend on the size of the problem. As a result the solution time does not depend directly on the size of the problem, contrary to what is typically the case with serial algorithms.

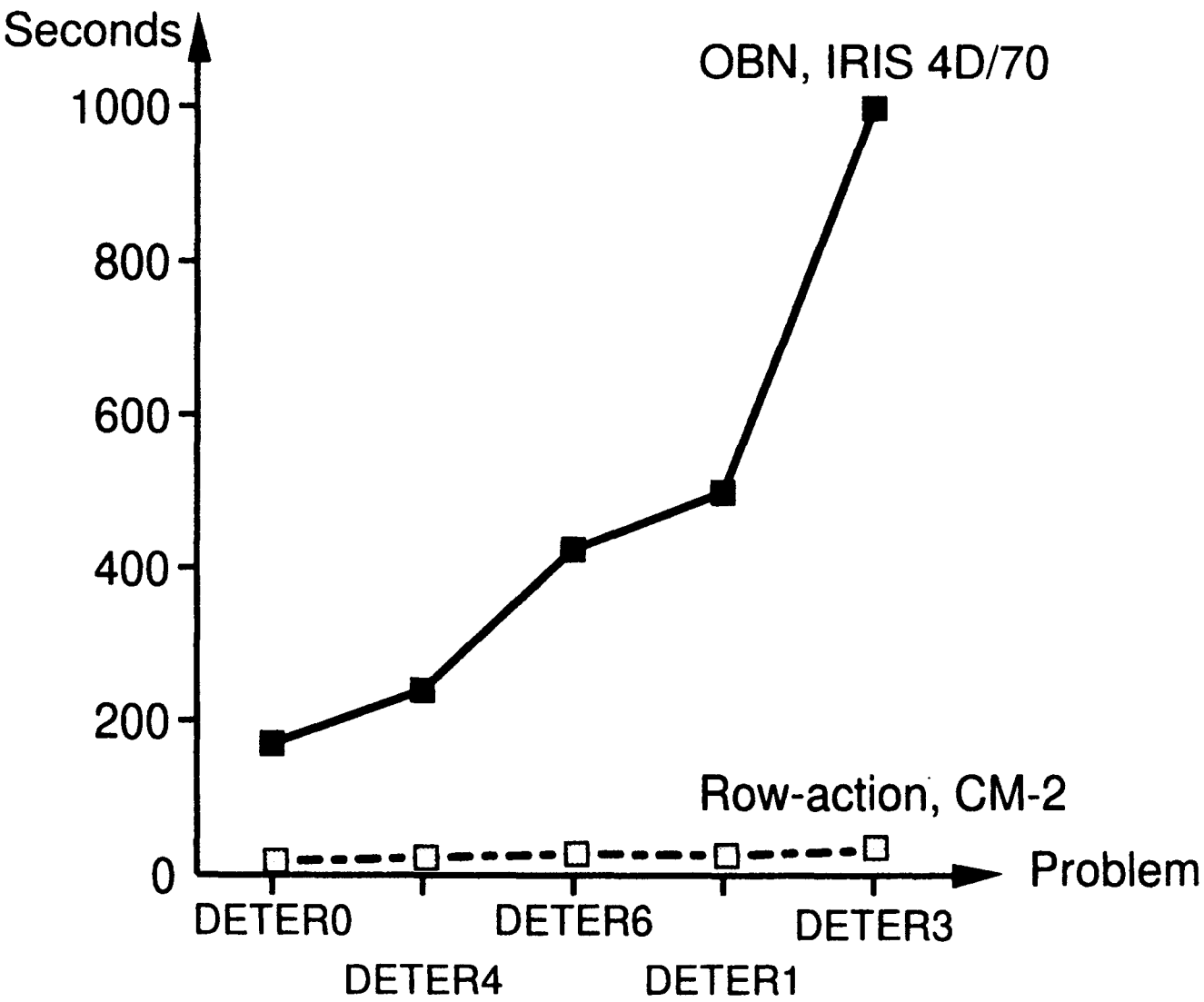

**Figure 15.7** Solving stochastic network problems with the row-action algorithm on the Connection Machine CM–2 with 8K processing elements, and with the general purpose interior point code OBN on a workstation IRIS 4D/70. Execution times are in seconds.

How do these computational results compare with those obtained with a state-of-the-art quadratic programming code, executing on a contemporary workstation? We compare the row-action Algorithm 13.8.1 executing on the CM–2 with the interior point code OBN on an IRIS 4D/70 workstation. The results show that the data parallel row-action algorithm substantially outperforms OBN (see Figure 15.7). It should be noted, however, that OBN is a general purpose optimizer while the row-action algorithm is specialized to take advantage of both the network structure and the availability of multiple processors. It is interesting to observe that the solution time for the row-action algorithm on the massively parallel machine increases only marginally with problem size, while the increase is substantial for OBN. We can thus conclude, supported by the results of Figure 15.7, that the data parallel implementation of the row-action algorithm for stochastic network problems is scalable, and very competitive with general purpose optimizers.

**Table 15.21** Solving large-scale quadratic stochastic network problems on the Connection Machine CM–2. Solution times are in CM seconds. $\epsilon$ denotes the $\ell_\infty$ norm of the error in primal feasibility upon termination. NA means not available.

| **Number of scenarios** | **Linear prog. formulation** rows × columns | **8K CM–2** | | **32K CM–2** | |
|---|---|---|---|---|---|
| | | **tolerance** $\epsilon$ | | | |
| | | $10^{-3}$ | $10^{-4}$ | $10^{-3}$ | $10^{-4}$ |
| 128 | 13583 × 38689 | 30.1 | 46.2 | 10.4 | 16.2 |
| 512 | 54287 × 154657 | 108.2 | 155.4 | 30.7 | 46.3 |
| 1024 | 108559 × 309281 | 210.8 | 326.5 | 57.3 | 86.3 |
| 2048 | 217103 × 618529 | 407.5 | 623.1 | 113.1 | 163.6 |
| 8196 | 868367 × 2474017 | NA | NA | NA | 11 min. |

To further challenge the algorithm, we tested it on extremely large-scale problems. Problem DETER3 was modified by replicating the scenarios until there were 128, 512, 1024, 2048 and 8196 scenarios, and the problems were solved on an 8K and a 32K Connection Machine CM–2. The results in Table 15.21 demonstrate the suitability of the algorithm for solving large-scale stochastic problems and prove that the algorithm scales very effectively for larger problems on larger machines. For example the 8K CM–2 solves 512 scenarios in 108.2 seconds, and a CM–2 with 32K processing elements solves four times as many scenarios in almost the same time. Figure 15.8 shows that the algorithm does not achieve very high FLOPS compared to the peak rate of CM–2, which is in GFLOPS.

## 15.6 Proximal Minimization with $D$-functions

We have seen that when implemented in parallel row-action algorithms are very efficient for various instances of nonlinear problems (matrix balancing, image reconstruction, transportation). However, they cannot be used to solve directly the important class of linear programming problems (see Chapter 3). Problems with linear objective functions can benefit from some of the row-action methods developed earlier, by employing a sequence of nonlinear perturbations using proximal minimization (PMD) algorithms. In this section we test the performance of two special cases of PMD algorithms (Algorithm 3.1.2) for solving a variety of linear programming problems. In particular, we combine the PMD algorithmic framework with parallel row-action algorithms for nonlinear network problems and nonlinear stochastic network problems. With this approach we solve large-scale linear programming problems from both problem categories.

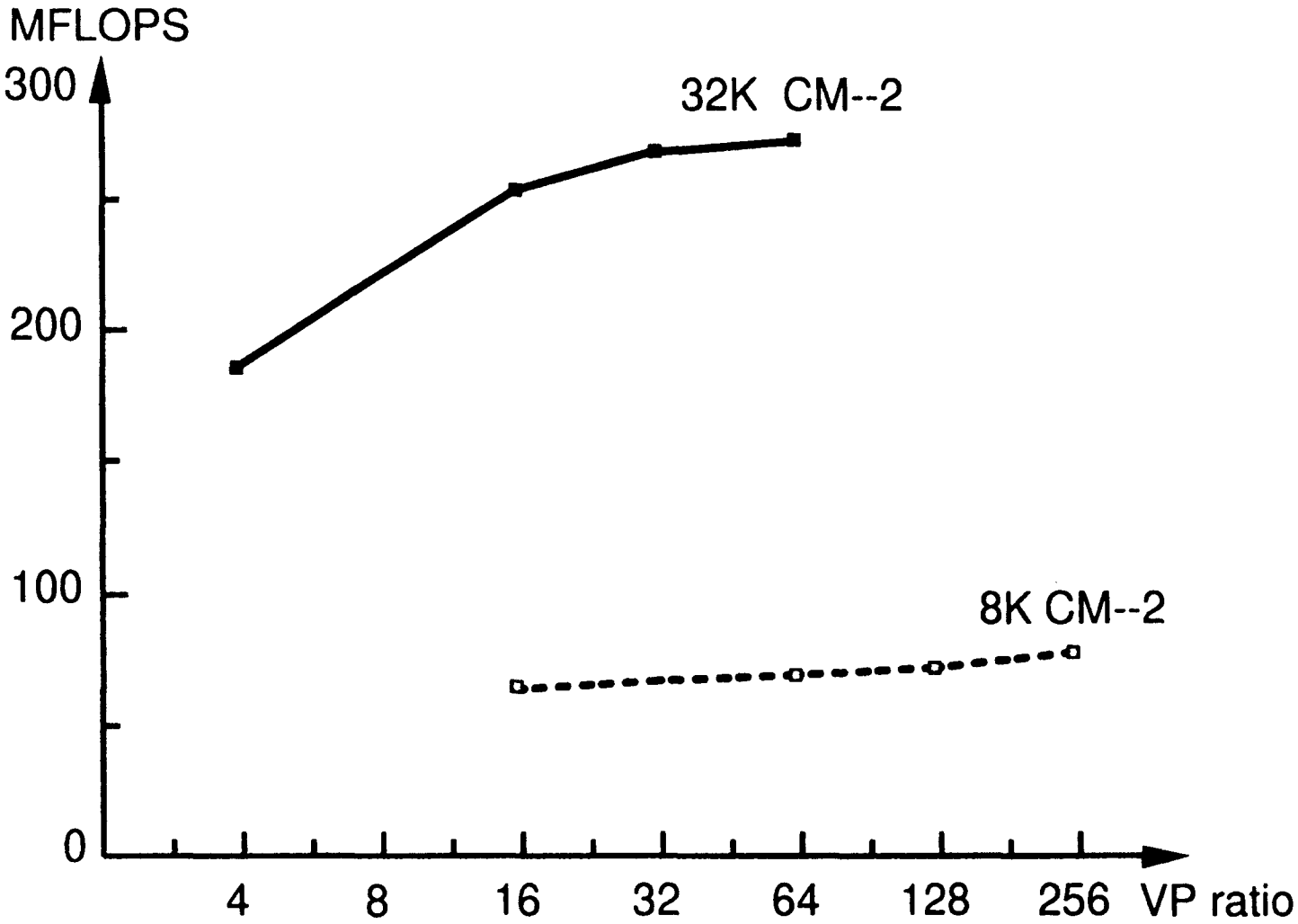

**Figure 15.8** MFLOPS rates with the row-action algorithm for stochastic network problems on a Connection Machine CM–2 with 8K and 32K processing elements.

### 15.6.1 Solving linear network problems

Consider the linear network flow problem $\min_{x \in X} \langle c, x \rangle$ where $X \subseteq \mathbb{R}^n$ is the polyhedral feasible region $X = \{x \in \mathbb{R}^n \mid Ax = b,\ 0 \leq x \leq u\}$, assumed to be nonempty, and $A$ is the $m \times n$ flow conservation constraint matrix. For some suitable choice of the Bregman function $f$, with zone $S$, and a positive sequence $\{\gamma_\nu\}_{\nu=0}^{\infty}$ with $\liminf_{\nu \to \infty} \gamma_\nu = \gamma < \infty$, the PMD algorithm proceeds from an arbitrary starting point $x^0 \in S$ according to the iteration

$$x^{\nu+1} \leftarrow \underset{x \in X \cap \bar{S}}{\operatorname{argmin}} \langle c, x \rangle + \frac{1}{\gamma_\nu} D_f(x, x^\nu). \tag{15.1}$$

As noted in Section 3.3, choosing the Bregman function $f(x)$ as in Example 2.1.1 or in Example 2.1.2 leads to the following iterative steps of two concrete algorithms.

**The Quadratic Proximal Point (QPP) Algorithm:**

$$x^{\nu+1} \leftarrow \underset{x \in X}{\operatorname{argmin}} \left( \langle c, x \rangle + \frac{1}{2\gamma_\nu} \parallel x - x^\nu \parallel^2 \right). \tag{15.2}$$

**Table 15.22** Size characteristics of the linear network test problems. Cost range refers to the range of the (cost) coefficients of the objective function.

| Problem | Nodes | Arcs | Cost range |
|---|---:|---:|---:|
| PDS1.1 | 126 | 339 | 0-99999 |
| PDS5.1 | 686 | 2,149 | 0-99999 |
| PDS10.1 | 1,399 | 4,433 | 0-99999 |
| PDS20.1 | 2,857 | 10,116 | 0-99999 |
| PDS30.1 | 4,223 | 15,526 | 0-99999 |
| Navy | 3,498 | 6,841 | 0-1 |
| Hugenavy | 30,639 | 64,542 | 0 |
| TP00 | 5,000 | 50,000 | 0-100 |
| TP01 | 5,000 | 100,000 | 0-100 |
| TP1 | 10,000 | 1,000,000 | 0-100 |
| TP2 | 20,000 | 1,000,000 | 0-100 |
| TP5 | 50,000 | 1,000,000 | 0-100 |
| TP10 | 50,000 | 2,000,000 | 0-100 |
| Gen0.25M | 32,768 | 245,760 | 0-10 |
| Gen0.5M | 65,536 | 491,520 | 0-10 |
| Gen16M | 2,097,152 | 15,728,640 | 0-10 |

**The Entropic Proximal Point (EPP) Algorithm:**

$$x^{k+1} \leftarrow \underset{x \in X \cap \bar{S}}{\text{argmin}} \left( \langle c, x \rangle + \frac{1}{\gamma_\nu} \left( \sum_{j=1}^{n} x_j \left( \log \left( \frac{x_j}{x_j^\nu} \right) - 1 \right) + \sum_{j=1}^{n} x_j^\nu \right) \right). \tag{15.3}$$

(The constant term $\sum_{j=1}^{n} x_j^\nu$ can be removed from the objective function.)

Both algorithms were implemented for the solution of linear transportation models on a Connection Machine CM-2 using the techniques developed in Section 14.5. In particular, the row-action Algorithms 12.4.1 and 12.4.2 were used to solve the resulting quadratic and entropic optimization subproblems, respectively.

The PMD algorithm was applied to solve linear programming versions of quadratic transportation problems (Table 15.9), and test problems derived from real-world applications with varying size characteristics (Table 15.22). The two algorithms are compared on medium-size problems arising from the applications (Figure 15.9), and we also examine execution times for both algorithms on very large, randomly generated test problems (Table 15.23).

How do the aforementioned results compare with the results achieved with a state-of-the-art network solver executing on a uniprocessor? When we compare the performance of the PMD algorithms executing on the Con-nection Machine CM–2 to the performance of the state-of-the-art network simplex code GENOS executing on a single processor of the CRAY Y-MP supercomputer (see Figure 15.10) we find that for smaller problems the network simplex code on the supercomputer is faster than the PMD code on the massively parallel system. However, the difference in performance is reduced as problems get larger, and the PMD code can also solve some extremely large problems with only moderate increase in solution times. Most interesting is the outcome with extremely large problems: here, the serial algorithm cannot solve the problems at all, even when using a contemporary supercomputer, while the PMD algorithm runs successfully.

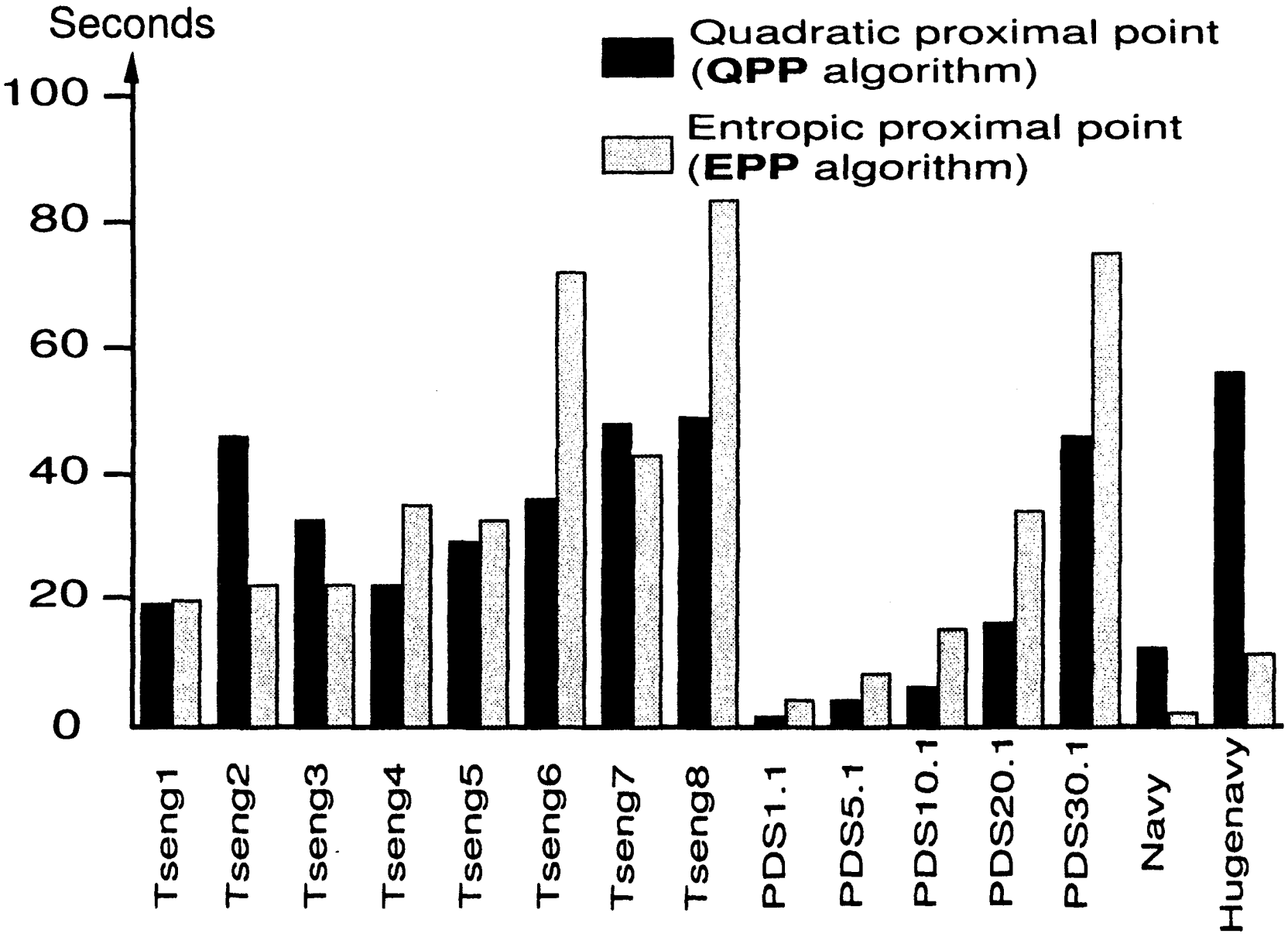

**Figure 15.9** Comparing the quadratic (QPP) and entropic (EPP) proximal point algorithms on the Connection Machine CM–2 with 16K processing elements.

### 15.6.2 Solving linear stochastic network problems

Both the QPP and EPP minimization algorithms were used to solve linear stochastic network problems. The size characteristics of the test problems are given in Table 15.24. The SPDA test problems were derived from a portfolio management application with insurance products (see Section 13.6). Computational results with the QPP algorithm for solving these test problems on a Connection Machine CM–2 are reported in Table 15.25, where we see that the PMD algorithm can reach solutions within a few percentage points of the exact values.

**Table 15.23** Solving large-scale randomly generated linear transportation problems using the quadratic (QPP) and entropic (EPP) proximal point algorithms. Test problems were solved on a Connection Machine CM–2 with 16K processing elements, with the exception of problem Gen16M that was solved on a CM–2 with 32K processing elements. Solution times are in CM minutes. NA means not available.

| Problem | QPP | | EPP | |
|---|---|---|---|---|
| | Iteration | Time | Iteration | Time |
| TP1 | 9000 | 47.5 | 1550 | 12.1 |
| TP2 | 5300 | 27.8 | 1275 | 10.0 |
| TP5 | 3200 | 17.0 | 1300 | 10.2 |
| TP10 | 4875 | 51.5 | 1725 | 27.0 |
| Gen0.25M | 2325 | 4.6 | 9900 | 25.4 |
| Gen0.5M | 2400 | 7.2 | 8350 | 44.9 |
| Gen16M | 3075 | 99 | NA | NA |

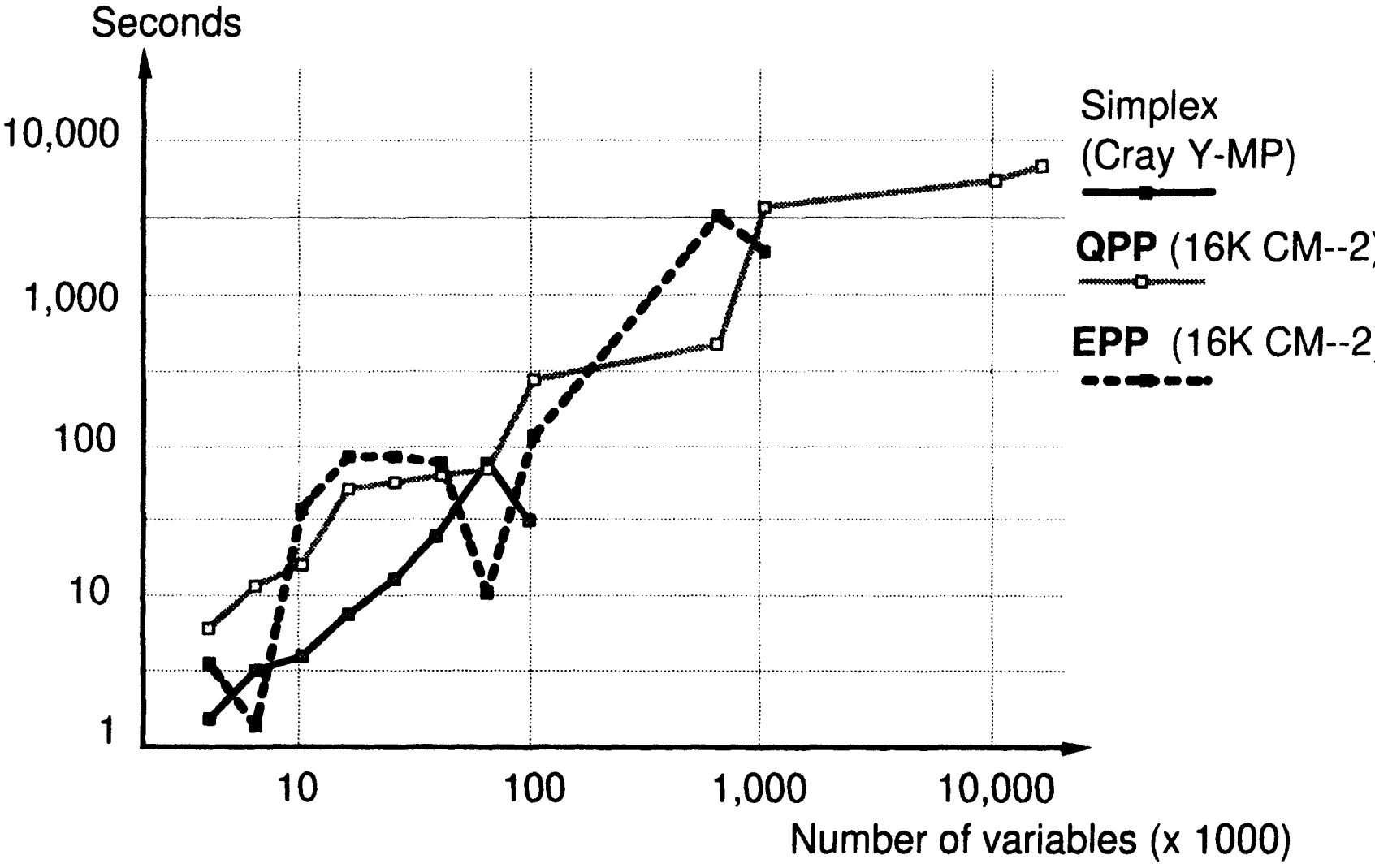

**Figure 15.10** Comparison of the quadratic (QPP) and entropic (EPP) proximal point algorithms with a network simplex algorithm running on a CRAY Y–MP vector supercomputer. The figure shows execution time (seconds) versus number of variables on a log-log scale, for problems (left to right) PDS10.1, Navy, PDS20.1, PDS30.1, TP00, Hugenavy, TP01, Gen0.5M, TP10, and Gen16M.

**Table 15.24** Size characteristics of the linear stochastic network test problems. Each scenario is a generalized network. Each scenario of the SPDA problems has 274 nodes and 697 arcs, and each scenario of the Random problems has 41 nodes and 107 arcs.

| Problem name | Number of stages | Number of scenarios | Split-variable linear program ( rows × cols.) | Objective value |
|---|---|---|---|---|
| SPDA-ext1 | 1 | 16 | $4,384 \times 11,152$ | 84.4333 |
| SPDA-ext2 | 2 | 16 | $5,524 \times 11,152$ | 68.2938 |
| SPDA-ext3 | 3 | 16 | $7,148 \times 11,152$ | 46.8079 |
| SPDA-ext4 | 4 | 16 | $8,540 \times 11,152$ | 33.2967 |
| SPDA-ext5 | 5 | 16 | $9,468 \times 11,152$ | 21.7297 |
| SPDA-ran1 | 1 | 16 | $4,384 \times 11,152$ | 50.1140 |
| SPDA-ran2 | 2 | 16 | $5,524 \times 11,152$ | 37.0613 |
| SPDA-ran3 | 3 | 16 | $7,148 \times 11,152$ | 36.5103 |
| SPDA-ran4 | 4 | 16 | $8,540 \times 11,152$ | 35.8041 |
| SPDA-ran5 | 5 | 16 | $9,468 \times 11,152$ | 34.7860 |
| Random2 | 2 | 256 | $13,301 \times 27,136$ | 122.4360 |
| Random3 | 3 | 256 | $15,941 \times 27,136$ | 108.5118 |
| Random4 | 4 | 216 | $15,511 \times 22,896$ | 98.0750 |
| Random5 | 5 | 256 | $20,825 \times 27,136$ | 91.3661 |
| Random6 | 6 | 243 | $21,997 \times 25,758$ | 89.6722 |
| Random9 | 9 | 256 | $30,219 \times 27,136$ | 91.0520 |

The QPP algorithm was also benchmarked against the state-of-the-art linear programming code OSL running on an HP 9000/730 workstation (Table 15.26). The results show that for the smaller test problems OSL executing on a workstation is uniformly faster than the massively parallel algorithm by a factor of 1.3 to 3.7. However, the massively parallel algorithm becomes increasingly competitive as the problem size increases. As Table 15.26 illustrates, when solving the Random3 problem and a similarly generated 3-stage problem with four times as many scenarios, the massively parallel algorithm took 4.2 times longer to solve the larger problem—due to a slight increase in the number of iterations—but OSL took 19.3 times longer on the HP 9000/730. (The objective value obtained with the QPP algorithm was within 0.02 percent of the exact value obtained using OSL.)

Finally, we examined the performance of the QPP algorithm on some extremely large test problems. Their impressive (by today's standards) sizes are summarized in Table 15.27 and the performance of the algorithm is summarized in Table 15.28. The algorithm reaches a high level of accuracy for these very large problems, although the execution times are substantial.

**Table 15.25** Solving the SPDA and the small random problems on the Connection Machine CM–2 using the QPP algorithm. Shown are the objective values achieved and the relative error from the exact values.

| Problem name | No. of major iterations | Time (sec.) | Objective value | Error |
|---|---|---|---|---|
| SPDA-ext1 | 1870 | 493.99 | 84.3751 | 0.07% |
| SPDA-ext2 | 2027 | 529.14 | 68.2624 | 0.05% |
| SPDA-ext3 | 1968 | 513.02 | 46.7932 | 0.03% |
| SPDA-ext4 | 4023 | 1048.60 | 33.3101 | 0.04% |
| SPDA-ext5 | 3389 | 822.64 | 21.6849 | 0.2% |
| SPDA-ran1 | 1286 | 538.93 | 50.0591 | -0.22% |
| SPDA-ran2 | 1141 | 614.80 | 37.0096 | -0.13% |
| SPDA-ran3 | 953 | 400.01 | 36.4401 | -0.2% |
| SPDA-ran4 | 883 | 468.21 | 35.7427 | -0.25% |
| SPDA-ran5 | 913 | 387.04 | 34.6951 | -0.27% |
| Random2 | 1576 | 517.18 | 122.4201 | 0.01% |
| Random3 | 1230 | 408.35 | 108.4764 | 0.03% |
| Random4 | 1271 | 422.20 | 98.0436 | 0.03% |
| Random5 | 1407 | 461.77 | 91.3439 | 0.02% |
| Random6 | 1491 | 494.03 | 89.6524 | 0.02% |
| Random9 | 954 | 313.16 | 90.9761 | 0.08% |

**Table 15.26** Comparing the quadratic proximal point (QPP) algorithm on the Connection Machine CM–2 to OSL's primal simplex algorithm on an HP 9000 workstation. Times are in seconds.

| Problem name | Number of scenarios | Proximal algorithm | | OSL |
|---|---|---|---|---|
| | | CM time | Iterations | CPU time |
| Random2 | 256 | 517.18 | 1576 | 376.83 |
| Random3 | 256 | 408.35 | 1230 | 314.13 |
| Random4 | 216 | 422.20 | 1271 | 199.50 |
| Random5 | 256 | 461.77 | 1407 | 195.47 |
| Random6 | 243 | 494.03 | 1491 | 147.23 |
| Random9 | 256 | 313.16 | 954 | 83.66 |
| Random3×4 | 1024 | 1700.96 | 1465 | 6049.29 |

**Table 15.27** Size characteristics of the large-scale, randomly generated stochastic network test problems. Large3 is a three-stage problem and problem Large5 has five stages.

| Problem name | Scenarios | Deterministic equivalent | Split-variable |
|---|---|---|---|
| Large3 | 16,392 | $492,036 \times 1,377,676$ | $1,030,773 \times 1,736,704$ |
| Large5 | 6,461 | $51,028 \times 154,178$ | $810,289 \times 695,466$ |

**Table 15.28** Solving large-scale stochastic network test problems with the QPP algorithm on a Connection Machine CM–2 with 32K processing elements. Solution times are in CM hours.

| Problem name | Number of iterations | Time | Absolute error | |
|---|---|---|---|---|
| | | | Primal error | Dual error |
| Large3 | 1871 | 4.81 | 0.0657 | 0.0010 |
| Large5 | 1782 | 1.94 | 0.01 | 0.0008 |

## 15.7 Description of Parallel Machines

We give here a brief description of the parallel machines used in the computational experiments, listing only the salient architectural features of each machine. Further details are available from the cited documentation. Reports by Dongarra (1990) and Dongarra and Duff (1985) compare several architectures available at the time of the writing of this book, including the machines discussed here, and are periodically updated to include information on new models.

### 15.7.1 Alliant FX/8

The Alliant FX/8 is a shared memory vector multiprocessor, introduced in 1985. It is configured with up to 8 processors, operating at 170 nsec (nanoseconds) clock cycle with peak rate of 94.4 MFLOPS. Memory subsystems consist of 8Mbyte units, and the FX/8 may have up to 8 units. On a fully configured system up to 2Gbytes of shared memory are available to a program. Each processor has vector functional units with hardware pieplining for each unit. Separate functional units are available for instruction fetch, address calculations, and floating point add and multiply operations. The units may be pipelined or operate in parallel. Similarily, the functional units of each processor may operate in parallel with the functional units of all other processors. Further details on the Alliant FX family can be found in the FX/FORTRAN Programmers Handbook, by Alliant Computer Systems Corporation (1985).

### 15.7.2 Connection Machine CM–2

The Connection Machine CM–2 is a fine-grain SIMD—Single Instruction, Multiple Data—system, introduced in 1988. Its basic hardware component is an integrated circuit with 16 processing elements (PEs) and a router that handles general communication. A fully configured CM has 4,096 chips for a total of 65,536 PEs. The 4,096 chips are interconnected as a 12-dimensional hypercube. Each processor is equipped with local memory of 8Kbytes (more on larger and later configurations), and for each cluster of 32 PEs a floating point accelerator handles floating point arithmetic. The machine operates at 6-7 MHz clock cycle. The peak performance achieved by a fully configured machine (i.e., 64K PEs) for optimization problems is 3GFLOPS, achieved for the solution of nonlinear dense transportation problems. Further details on the CM–2 are given in Hillis (1985, 1987) and Hillis and Steele (1986).

### 15.7.3 Connection Machine CM–5

The Connection Machine CM–5 is a distributed memory MIMD machine that supports data parallel computing. Computing nodes are based on the SPARC 2 RISC processor with up to 32Mbytes of local memory, (the experimental CM–5e model is based on the SPARC 10 processor), and each processor also supports four vector functional units. The communication network is a *fat tree*, which is reportedly scalable for configurations up to 16,384 processors. The largest configuration installed to date has 1000 processors, and the largest system used in the experiments reported in this book has 64 processors. Futher details on the CM–5 can be found in Leiserson et al. (1992).

### 15.7.4 CRAY X-MP and Y-MP

The CRAY X-MP and Y-MP are general purpose vector multiprocessors, with tightly coupled processors and shared memory. The X-MP was introduced in 1983 and through the 1980s and early 90s it and its successor, the Y-MP, were considered the yardsticks against which parallel and supercomputers were evaluated. The X-MP has up to 4 processors, and each processor has vector functional units for floating point addition, multiplication and reciprocal approximation, vector units for integer addition, logical and shift operations, and special units for integer scalar addition, logical and shift operations. The X-MP operates at 8.5 nsec clock cycle. The Y-MP is configured with up to 8 processors, and it operates at 4.2 nsec. Further details on the architecture of this family can be found in Chen (1984).

### 15.7.5 Intel iPSC/860

The Intel iPSC/860 is a distributed memory MIMD computer, with a hypercube communication network. Computing nodes are 32-bit i860 RISC

processors, with up to 16Mbytes of local memory. A fully configured system can have up to 128 processors. Further details on the iPSC/860 can be found in the System User's Guide by Intel Corporation (1992).

## 15.8 Notes and References

**15.1** Guidelines for reporting computational experiments with mathematical software have been developed through initiatives of the Committee on Algorithms (COAL) of the Mathematical Programming Society. A set of guidelines was developed by Crowder, Dembo, and Mulvey (1978), and was revised by Jackson et al. (1990); see also Greenberg (1990). Barr and Hickman (1993) discuss issues in reporting computational experiences with parallel algorithms. See also the commentaries to their article written by Gustafson (1993) and Lusk (1993). The guidelines here are based on the suggestions of Barr and Hickman (1993), and the modifications by Gustafson (1993).

**15.2** The 1977 United States national data for regional input/output were obtained from the data tape from the Interindustry Economic Division, Bureau of Economic Analysis (1977). The experiments with the data parallel implementation of RAS on the Connection Machine CM–2 are reported in Zenios (1990). Zenios and Iu (1990) report results with the control parallel implementation on an Alliant FX/8. Results with the range-RAS algorithm on the Connection Machine CM–2 are due to Censor and Zenios (1991).

**15.3** The Shepp and Logan test phantom is described in the paper by the same authors (1976). The numerical results with the block-MART algorithm for image reconstruction are due to Zenios and Censor (1991b). Experimental work with block versions of other algorithms can be found in Hudson and Larkin (1994) for the EM algorithm and in Herman and Levkowitz (1987) for the block-ART (block-Kaczmarz) algorithm; see also Section 10.7. Bramely and Sameh (1991) investigated a domain decomposition approach that allows parallelism with a row-action method such as Kaczmarz's method, see also Kamath and Sameh (1988).

**15.4.1** The experiments with the data parallel implementation of the row-action algorithms for nonlinear transportation problems are due to Zenios and Censor (1991a). McKenna and Zenios (1990) discuss an optimal implementation of the algorithm using the hypercube communication model of the Connection Machine CM–2. Computational results with a similar relaxation algorithm using control parallelism in an asynchronous environment are due to Chajakis and Zenios (1991), who based their work on earlier analysis by Bertsekas and El Baz (1987). The relaxation algorithm was developed in Bertsekas, Hossein and Tseng (1987) and Zenios and Mulvey (1986b); a textbook treatment is given in Bertsekas and Tsitsiklis (1989). See also the PhD thesis by Tseng (1986). The results for the equilibration algorithm are due to Nagurney, Kim, and Robinson (1990), and for the exact Newton algorithm due to Klincewicz (1989).

**15.4.2** The results for the data parallel implementation of the row-action algorithm for multicommodity network flow problems are reported in Zenios (1991b). Censor, Chajakis, and Zenios (1995) report numerical results with control parallel implementations of the same algorithm.

**15.4.3** The parallel performance of the LQP algorithm was tested by Pinar and Zenios (1992). The same authors discuss the data parallel implementation of the algorithm in Pinar and Zenios (1994a). See also Pinar and Zenios (1993) for a comparative study of various parallel algorithms for multicommodity network flow problems. BLAS is a set of basic linear algebra subroutines, and is implemented for peak efficiency on several machines; see Lawson et al. (1979) and Dongarra et al. (1988a, 1988b). OB1 is a linear programming code based on a primal-dual path following algorithm due to Lustig, Marsten, and Shanno (1991).

**15.5.1** The SCSD8 problems are also solved in Birge and Holmes (1992). They report solution of the largest test problem to date, with 32 scenarios on a network of workstations. The SEN problems are due to Sen, Doverspike, and Cosares (1994). The largest SEN test problem solved to date, with 1000 scenarios, was reported by G.B. Dantzig and G. Infanger at the International Institute of Applied Systems Analysis (IIASA) workshop on stochastic programming, Laxenburg, Austria, in June 1994. Results with the parallel factorization of structured matrices arising in stochastic programming are reported by Jessup, Yang and Zenios (1994a). They report results on both an Intel iPSC/860 and a Connection Machine CM–5. A report by the same authors (1994b) includes an analysis of the results based on a communication model of the hypercube. The solution of linear and quadratic stochastic programming problems using the parallel matrix factorization with interior point methods is discussed in Yang and Zenios (1996). LOQO (Linear or Quadratic Optimizer) is described in Vanderbei (1992) and Vanderbei and Carpenter (1993).

**15.5.2** The results with the row-action algorithms for stochastic network problems are from Nielsen and Zenios (1993a, 1996a) for two-stage and multistage problems, respectively. See also Nielsen and Zenios (1992b) for a discussion of the effects of problem structure on the performance of the algorithms. The test problems are based on the models from Mulvey and Vladimirou (1991, 1992). OBN is the quadratic programming optimization version of the code OB1. OB1 is a linear programming code based on a primal-dual path following algorithm due to Lustig, Marsten, and Shanno (1991).

**15.6** Results for the solution of linear network problems using PMD algorithms are reported in Nielsen and Zenios (1993b, 1993c). Nielsen and Zenios (1993d) report results for the solution of two-stage stochastic network problems, and in (1996a) they extend the algorithm and perform further experiments for multi-stage problems. The SPDA portfolio models are described in Nielsen and Zenios (1996b). The Generalized Network Optimization System, GENOS, software system is described in Mulvey and Zenios (1987b). OSL is the Optimization System Library, developed by IBM Corporation and described in the documentation by International Business Machines (IBM) Corporation (1991).

# Index

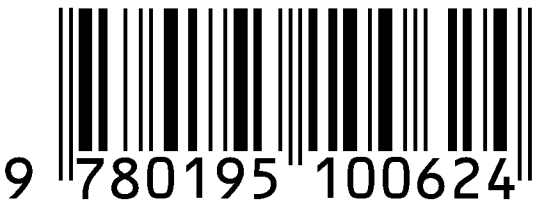
9 780195 100624